\documentclass[10pt]{amsart}
\usepackage{amssymb,amsmath,amsfonts,amsthm,graphics,mathrsfs,amscd,hyperref}
\usepackage[hmargin=1in,vmargin=1in]{geometry}
\usepackage[all,cmtip]{xy}
\usepackage{cmap}
\usepackage[T1]{fontenc}
\usepackage[utf8]{inputenc}
\usepackage{lmodern}

\usepackage{cancel}

\theoremstyle{plain}

\newtheorem{definition}[equation]{Definition}

\newtheorem{lemma}[equation]{Lemma}
\newtheorem{proposition}[equation]{Proposition}
\newtheorem{theorem}[equation]{Theorem}

\newtheorem{remark}[equation]{Remark}

\numberwithin{equation}{subsection}

\title[Obstructed and unobstructed foliated deformations]{On deformations of foliated complex analytic structures IV-1 : Obstructed and unobstructed foliated deformations}
\author{Chunghoon Kim}

\begin{document}

\maketitle
\begin{abstract}
We study obstructions of deformations of foliated complex analytic structures. We describe obstructions in terms of \v Cech resolutions. From foliated complex analytic families in Part I and unfolded complex analytic families in Part III, we deduce the integrability conditions. We identify infinitesimal deformations and obstructions in terms of Dolbeault resolutions or Dolbeault type bicomplexes from the integrability conditions. Then we describe the correspondences of infinitesimal deformations and obstructions between \v Cech resolutions, and Dolbeault resolutions or Dolbeault type bicomplexes.

\end{abstract}

\tableofcontents

\section{Infinitesimal foliated deformations and obstructions in terms of tangent sheaves by using \v Cech resolution of $\mathcal{E}_{\Theta_{\mathcal{F}_0}}^\bullet$}

\subsection{Review of deformations of foliated complex analytic structures in terms of tangent sheaves}\

We recall the definitions of a foliated analytic family in terms of tangent sheaves and associated foliated Kodaira-Spencer map in terms of tangent sheaves from Part I.
\begin{definition}\label{t3}
Suppose that given a domain $B\subset \mathbb{C}^m$, there is a set $\left\{\left(M_t, \Theta_{\mathcal{F}_t}\right)|t\in B \right\} $ of $n$-dimensional compact $($singularly$)$ foliated complex manifolds $(M_t, \Theta_{\mathcal{F}_t})$ of dimension $p$, so that for each $t\in B$, we have an exact sequence
\begin{align}\label{t1}
0\to \Theta_{\mathcal{F}_t}\to \Theta_{M_t}\to \Theta_{M_t}/\Theta_{\mathcal{F}_t}\to 0
\end{align}
where $\Theta_{M_t}/\Theta_{\mathcal{F}_t}$ is a torsion-free $\mathcal{O}_{M_t}$-module. We say that $\left\{\left(M_t, \Theta_{\mathcal{F}_t}\right)|t\in B\right\}$ is a family of $($singularly$)$ foliated complex manifold or $($singularly$)$ foliated complex analytic family in terms of tangent sheaves if there exist a complex manifold $\mathcal{M}$ and a holomorphic map $\pi:\mathcal{M}\to B$  such that
\begin{enumerate}
\item we have an exact sequence
\begin{align}\label{t2}
0\to \Theta_\mathcal{F} \to \Theta_{\mathcal{M}/B}\to \Theta_{\mathcal{M}/B}/\Theta_{\mathcal{F}}\to 0
\end{align}
where $\Theta_{\mathcal{F}}$ is a coherent subsheaf of the relative tangent sheaf $\Theta_{\mathcal{M}/B}$ over $B$.
\item for each $t\in B$, $\pi^{-1}(t)= \left(M_t, \Theta_{\mathcal{F}_t}\right)$ and the exact sequence $(\ref{t1})$ is induced from $(\ref{t2})$ by restricting to $M_t$
\item the rank of the Jacobian of $\pi$ is equal to $m$ at every point of $\mathcal{M}$.
\item $[\Theta_\mathcal{F},\Theta_\mathcal{F}]\subset \Theta_\mathcal{F}$ under the Lie bracket and $\Theta_{\mathcal{M}/ B}/\Theta_\mathcal{F}$ is flat over $B$.
\end{enumerate}
We will denote the $($singularly$)$ foliated complex analytic family by $(\mathcal{M}, \Theta_\mathcal{F}, B, \pi)$.
\end{definition}

\begin{remark}
Let $\left(\mathcal{M}, \Theta_\mathcal{F}, B, \pi \right)$ with $\Theta_\mathcal{F}$ locally free be a foliated analytic family. Let $\Delta$ be an open set of $B$. Then the restriction $\left(\mathcal{M}_\Delta= \pi^{-1}(\Delta), \Theta_{\mathcal{F}}|_\Delta, \Delta, \pi |_{\mathcal{M}_\Delta} \right)$ is also a foliated complex analytic family in terms of tangent sheaves. We will denote the family by $\left( \mathcal{M}_\Delta, \Theta_{\mathcal{F}_\Delta}, \Delta, \pi\right)$.
\end{remark}

\begin{remark}\label{o4}
Given a foliated analytic family $\left( \mathcal{M}, \Theta_\mathcal{F}, B, \pi \right)$ with $\Theta_{\mathcal{F}}$ locally free as in \textnormal{Definition \ref{t3}}, we can choose a system of local complex coordinates $\{\xi_i,...,\xi_j,...\}, \xi_j:p\to \xi_j(p)$, and coordinate polydisks $\mathcal{U}_j$ with respect to $\xi_j$, satisfying $(1)$ $\xi_j(p)=\left(\xi_j^1(p),...,\xi_j^n(p),t_1,...,t_m\right),\left(t_1,...,t_m\right)=\pi(p)$ and $(2)$ $\mathcal{U}=\left\{\mathcal{U}_j|j=1,2,...\right\}$ is locally finite. Then $\left\{p\to (\xi_j^1(p),...,\xi_j^n(p))|\mathcal{U}_j\cap M_t\ne \emptyset \right\}$ gives a system of local complex coordinates on $M_t$. In terms of these coordinates, $\pi$ is the projection given by $\pi: \left(\xi_j^1,...,\xi_j^n, t_1,...,t_m \right)\to \left(t_1,...,t_m\right)$. For $j,k$ with $\mathcal{U}_j\cap \mathcal{U}_k\ne \emptyset$, we denote that coordinate transformation from $\xi_k$ to $\xi_j$ by $f_{jk}:\left(\xi_k^1,...,\xi_k^n,t \right)\to \left(\xi_j^1,...,\xi_j^n,t \right)=f_{jk}\left(\xi_k^1,...,\xi_k^n,t \right)$. Thus $f_{jk}$ is given by $\xi_j^\alpha=f_{jk}^\alpha\left(\xi_k^1,...,\xi_k^n,t_1,...,t_m\right), \alpha=1,...,n$. On the other hand, we may assume that the locally free sheaf $\Theta_\mathcal{F}$ is trivialized in the covering $\mathcal{U}$ so that $\Gamma\left(\mathcal{U}_j,\Theta_\mathcal{F}\right)\cong \bigoplus^p \Gamma\left(\mathcal{U}_j, \mathcal{O}_\mathcal{M}\right)$ which is generated by vector fields $T_j^\alpha(\xi_j,t)\in \Gamma\left(\mathcal{U}_j, \Theta_\mathcal{F}\right),\alpha=1,...,p$ of the form
\begin{align}\label{ta3}
T_j^\alpha(\xi_j,t):=\sum_{\beta=1}^n T_j^{\alpha\beta}(\xi_j,t)\frac{\partial}{\partial \xi_j^\beta},\,\,\,\,\,\,\,\,\,\alpha=1,...,p
\end{align}
for some $T_j^{\alpha\beta}(\xi_j,t)\in \Gamma\left( \mathcal{U}_j , \mathcal{O}_\mathcal{M} \right)$, and satisfies
\begin{align}\label{ta5}
\left[T_j^\alpha(\xi_j,t), T_j^\beta(\xi_j,t)\right]=\sum_{\gamma=1}^n g_{j\alpha\beta}^\gamma(\xi_j,t) T_j^\gamma(\xi_j,t)
\end{align}
for some $g_{j\alpha\beta}^\gamma(\xi_j,t)\in \Gamma\left(\mathcal{U}_j, \mathcal{O}_\mathcal{M}\right)$ with $g_{j\alpha\beta}^\gamma(\xi_j,t)=-g_{j\beta\alpha}^\gamma(\xi_j,t)$, and for $j,k$ with $\mathcal{U}_j\cap \mathcal{U}_k \ne \emptyset$, we have
\begin{align}\label{o3}
T_j^\alpha(\xi_j,t)=\sum_{\beta=1}^p r_{jk}^{\alpha\beta}(\xi_k,t) T_k^\beta(\xi_k,t)
\end{align}
for some $r_{jk}^{\alpha\beta} (\xi_k,t)\in \Gamma\left(\mathcal{U}_j\cap \mathcal{U}_k, \mathcal{O}_\mathcal{M}\right)$, and for $i,j,k$ with $\mathcal{U}_i\cap \mathcal{U}_j\cap \mathcal{U}_k\ne \emptyset$, we have
\begin{align}\label{o1}
r_{ik}^{\alpha \gamma}(\xi_k,t)= \sum_{\beta=1}^p r_{ij}^{\alpha\beta}\left(f_{jk}(\xi_k,t),t) \right) r_{jk}^{\beta \gamma}(\xi_k,t)
\end{align}
More specifically $(\ref{o3})$ implies that
\begin{align}\label{ta4}
\sum_{\gamma=1}^n T_j^{\alpha\gamma}\left(f_{jk}(\xi_k,t),t\right)\frac{\partial}{\partial \xi_j^\gamma}=\sum_{\beta=1}^p\sum_{\eta=1}^n\sum_{\gamma=1}^n r_{jk}^{\alpha\beta}(\xi_k,t) T_k^{\beta\eta}(\xi_k,t)\frac{\partial f_{jk}^\gamma(\xi_k,t)}{\partial \xi_k^\eta}\frac{\partial}{\partial \xi_j^\gamma}
\end{align}

\end{remark}

Given a foliated complex analytic family $\left(\mathcal{M}, \Theta_\mathcal{F}, B, \pi \right)$ with $\Theta_\mathcal{F}$ locally free of deformations of $\left(M, \Theta_{\mathcal{F}_0} \right)=\pi^{-1}(0),0\in B$, the infinitesimal foliated deformations in the foliated complex analytic family  are encoded in $\mathbb{H}^1\left( M, \Theta_{\mathcal{F}_0}^\bullet \right)$ which is the first cohomology group of the following leaf complex
\begin{align}
\Theta_{\mathcal{F}_0}^\bullet : \Theta_M \to  \mathscr{H}om_{\mathcal{O}_M}\left( \Theta_{\mathcal{F}_0}, \frac{\Theta_M}{\Theta_{\mathcal{F}_0}} \right)  \to  \mathscr{H}om_{\mathcal{O}_M}\left( \bigwedge^2 \Theta_{\mathcal{F}_0}, \frac{\Theta_M}{\Theta_{\mathcal{F}_0}} \right) \to \mathscr{H}om_{\mathcal{O}_M}\left( \bigwedge^3 \Theta_{\mathcal{F}_0}, \frac{\Theta_M}{\Theta_{\mathcal{F}_0}} \right) \to \cdots
\end{align}
We will denote the $i$-th cohomology group by $\mathbb{H}^1\left( M , \Theta_{\mathcal{F}_0}^\bullet \right)$. We can compute $\mathbb{H}^i\left(M_0, \Theta_{\mathcal{F}_0}^\bullet \right)$  by the following \v Cech resolution of $\Theta_{\mathcal{F}_0}^\bullet$. Here $\delta$ is the \v Cech map and $\mathcal{U}_0= \mathcal{U}\cap M=\left\{ U_j :=\mathcal{U}_j \cap M | j=1,2,... \right\}$ is an open covering of $M$:
{\Small{\begin{equation} \label{o5}
\begin{CD}
\cdots \\
@AD_2 AA \\
C^0\left(\mathcal{U}_0, \mathscr{H}om_{\mathcal{O}_M}\left( \bigwedge^2  \Theta_{\mathcal{F}_0}, \frac{\Theta_M}{\Theta_{\mathcal{F}_0} } \right) \right) @>-\delta>> \cdots  \\
@AD_1AA @AD_1 AA \\
C^0\left(\mathcal{U}_0, \mathscr{H}om_{\mathcal{O}_M}\left(\Theta_{\mathcal{F}_0}, \frac{\Theta_{M}}{\Theta_{\mathcal{F}_0} } \right) \right)@>\delta>> C^1\left(\mathcal{U}_0, \mathscr{H}om_{\mathcal{O}_{M}}\left(\Theta_{\mathcal{F}_0}, \frac{\Theta_M}{\Theta_{\mathcal{F}_0}}  \right) \right) @>-\delta>> \cdots \\
@AD_0 AA @AD_0 AA @AD_0 AA \\
C^0(\mathcal{U}_0, \Theta_{M}) @>-\delta>> C^1(\mathcal{U}_0, \Theta_{M}) @>\delta>> C^2(\mathcal{U}_0, \Theta_{M}) @>>> \cdots 
\end{CD}
\end{equation}}}
Then we can define the foliated Kodaira-Spencer map at $0\in B$ in terms of tangent sheaves which encodes the information of infinitesimal foliated deformations of $\pi^{-1}(0)= \left( M, \Theta_{\mathcal{F}_0} \right)$ in the foliated analytic family $\left(\mathcal{M},  \Theta_\mathcal{F}, B, \pi \right)$ (see Part I).
\begin{definition}[foliated Kodaira-Spencer map in terms of tangent sheaf]\label{ta12}
Let $(\mathcal{M}, \Theta_\mathcal{F}, B, \pi)$ with $\Theta_{\mathcal{F}}$ locally free be a foliated complex analytic family of deformations of $(M, \Theta_{\mathcal{F}_0})=\pi^{-1}(0),0\in B$, where $B$ is a domain of $\mathbb{C}^m$. We keep the notations in \textnormal{Remark \ref{o4}}. For a tangent vector $\frac{\partial}{\partial t}=\sum_{\lambda=1}^m c_\lambda \frac{\partial}{\partial t_\lambda}, c_\lambda \in \mathbb{C}$, of $B$, the foliated Kodaira-Spencer map in terms of tangent sheaf at $0$ is defined to be a $\mathbb{C}$-linear map
{\small{\begin{align*}
\varphi_0:T_0(B)&\to \mathbb{H}^1\left(M, \Theta_{\mathcal{F}_0}^\bullet \right)\\
                \frac{\partial}{\partial t}&\mapsto \frac{\partial \left(M_t, \Theta_{\mathcal{F}_t}\right)}{\partial t}|_{t=0}:=\left( \left\{  \sum_{\alpha=1}^n \frac{\partial f_{jk}^\alpha(\xi_k,t)}{\partial t}|_{t=0}\frac{\partial}{\partial \xi_j^\alpha} \right\}, \left\{ \left(T_{0j}^\alpha :=T_j^\alpha(\xi_j,0)\mapsto -\overline{\frac{\partial T_j^\alpha(\xi_j,t)}{\partial t}}|_{t=0}\right)_{\alpha=1,...,p} \right\} \right) \\
                 &\,\,\,\,\,\,\,\,\,\,\in C^1\left(\mathcal{U}_0, \Theta_{M} \right) \bigoplus C^0\left( \mathcal{U}_0, \mathscr{H}om_{\mathcal{O}_{M}}\left(\Theta_{\mathcal{F}_0}, \frac{\Theta_{M}}{\Theta_{\mathcal{F}_0}} \right) \right)
\end{align*}}}
Here $\frac{\partial T_j^\alpha(\xi_j,t)}{\partial t}|_{t=0}:=\sum_{\beta=1}^n \frac{\partial T_j^{\alpha\beta}(\xi_j,t)}{\partial t}|_{t=0}\frac{\partial}{\partial \xi_j^\beta} $ is the derivative of $\left(\ref{ta3} \right)$ with respect to $t$ at $0$, and $\overline{\frac{\partial T_j^\alpha(\xi_j,t)}{\partial t}}|_{t=0}$ is the image of $\frac{\partial T_j^\alpha(\xi_j,t)}{\partial t}|_{t=0}$ under $\Theta_M \to \frac{\Theta_M}{\Theta_{\mathcal{F}_0}}$.
\end{definition}

In Part II we have shown that the complex of sheaves
\begin{align}\label{o6}
\mathcal{E}_{\Theta_{\mathcal{F}_0}}^\bullet: \mathcal{E}_{\Theta_{\mathcal{F}_0}}\to \mathscr{H}om_{\mathcal{O}_M}\left( \Theta_{\mathcal{F}_0}, \Theta_M \right) \to  \mathscr{H}om_{\mathcal{O}_M}\left( \bigwedge^2 \Theta_{\mathcal{F}_0}, \frac{\Theta_M}{\Theta_{\mathcal{F}_0}} \right) \to \mathscr{H}om_{\mathcal{O}_M}\left( \bigwedge^3 \Theta_{\mathcal{F}_0}, \frac{\Theta_M}{\Theta_{\mathcal{F}_0}} \right) \to \cdots
\end{align}
also controls foliated deformations of $\left( M , \Theta_{\mathcal{F}_0} \right)$ in terms of tangent sheaves. We can actually study obstructions of foliated deformations of $\left(M, \Theta_{\mathcal{F}_0} \right)$ in terms of tangent sheaves by using the complex of sheaves $\Theta_{\mathcal{F}_0}^\bullet$.  However obstructions of deformations of the underlying locally free sheaf $\Theta_{\mathcal{F}_0}$ are implicit in $\mathbb{H}^2\left( M, \Theta_{\mathcal{F}_0}^\bullet \right)$, i.e. in \v Cech resolution $(\ref{o5})$ of $\Theta_{\mathcal{F}_0}^\bullet$. However the complex of sheaves $\mathcal{E}_{\Theta_{\mathcal{F}_0}}^\bullet$ explicitly shows the infinitesimal deformations and obstructions of deformations of the underlying locally free sheaf $\Theta_{\mathcal{F}_0}$. Hence it is desirable to use $\mathcal{E}_{\Theta_{\mathcal{F}_0}}^\bullet$ instead of $\Theta_{\mathcal{F}_0}^\bullet$ to identify obstructions clearly. Let us identify the infinitesimal deformations in terms of the complex (\ref{o6}) of sheaves  $\mathcal{E}_{\Theta_{\mathcal{F}_0}}^\bullet$.

\subsection{Identification of infinitesimal foliated deformations by using $\mathcal{E}_{\Theta_{\mathcal{F}_t}}^\bullet$} \label{o41}\

Let $\left( \mathcal{M}, \Theta_\mathcal{F}, B, \pi \right)$ be a foliated analytic family with $\Theta_\mathcal{F}$ locally free as in Definition $\ref{t3}$. We keep the notations in Remark $\ref{o4}$. Let $\mathcal{U}_t=\mathcal{U} \cap M_t =\left\{ U_j^t:=\mathcal{U}_j\cap M_t| j=1,2,... \right\}$ be an open covering of $M_t$. We recall the definition of $\mathcal{E}_{\Theta_{\mathcal{F}_t}}$ on $M_t$ with respect to the open covering $\mathcal{U}_t=\left\{U_j^t\right\}$ of $M_t$ from Part II, and let us consider the \v Cech resolution of $\mathcal{E}_{\Theta_{\mathcal{F}_t}}^\bullet$:

{\Small{\begin{equation} \label{o7}
\begin{CD}
\cdots \\
@AD_2^t AA \\
C^0\left(\mathcal{U}_t, \mathscr{H}om_{\mathcal{O}_M}\left( \bigwedge^2  \Theta_{\mathcal{F}_t}, \frac{\Theta_M}{\Theta_{\mathcal{F}_0} } \right) \right) @>-\delta>> \cdots  \\
@AD_1'^tAA @AD_1 AA \\
C^0\left(\mathcal{U}_t, \mathscr{H}om_{\mathcal{O}_{M_t}}\left(\Theta_{\mathcal{F}_t}, \Theta_{M_t} \right) \right)@>\delta>> C^1\left(\mathcal{U}_t, \mathscr{H}om_{\mathcal{O}_{M}}\left(\Theta_{\mathcal{F}_t}, \Theta_{M_t}  \right) \right) @>-\delta>> \cdots \\
@AD_0'^t AA @A{D}_0'^tAA @AD_0'^tAA \\
C^0\left(\mathcal{U}_t,\mathcal{E}_{\Theta_{\mathcal{F}_t}}\right) @>-\delta>> C^1\left(\mathcal{U}_t, \mathcal{E}_{\Theta_{\mathcal{F}_t}}\right) @>\delta>> C^2 \left(\mathcal{U}_t, \mathcal{E}_{\Theta_{\mathcal{F}_t}} \right) @>>> \cdots 
\end{CD}
\end{equation}}}

By taking the derivative of $(\ref{ta3})$, we set $\frac{\partial T_j^{\alpha}(\xi_j,t)}{\partial t}:=\sum_{\beta=1}^n \frac{\partial T_j^{\alpha\beta}(\xi_j,t)}{\partial t}\frac{\partial}{\partial \xi_j^\beta},\alpha=1,...,p $. We define an element $\alpha_j(t) \in \Gamma\left( U_j^t, \mathscr{H}om_{\mathcal{O}_{M_t}}\left( \Theta_{\mathcal{F}_t}, \Theta_{M_t} \right) \right)$ in the following way:
\begin{align*}
\alpha_j(t):\Gamma\left( U_j^t, \Theta_{\mathcal{F}_t}   \right) &\to \Gamma\left( U_j^t, \Theta_{M_t} \right) \\
   T_j^\alpha=T_j^\alpha(\xi_j,t)&\mapsto -\frac{\partial T_j^{\alpha}(\xi_j,t)}{\partial t}= - \sum_{\beta=1}^n \frac{\partial T_j^{\alpha\beta}(\xi_j,t)}{\partial t}\frac{\partial}{\partial \xi_j^\beta},\,\,\,\,\,\,\,\,\alpha=1,...,p
\end{align*}
and linearly extend to $\Gamma\left( U_j^t , \Theta_{\mathcal{F}_t} \right)$. We define an element $A_{jk}(t) \in \Gamma\left( U_j^t \cap U_k^t, \mathscr{H}om_{\mathcal{O}_{M_t} }\left( \Theta_{\mathcal{F}_t}, \Theta_{\mathcal{F}_t} \right)  \right)$ in the following way:
\begin{align*}
A_{jk}(t):\Gamma\left( U_{jk}^t, \Theta_{\mathcal{F}_t} \right) &\to \Gamma\left( U_{jk}^t, \Theta_{\mathcal{F}_t} \right) \\
 T_j^\alpha &\mapsto \sum_{\beta=1}^p \frac{\partial r_{jk}^{\alpha\beta} (\xi_j,t)}{\partial t} T_k^\beta,\,\,\,\,\,\,\,\,\,\alpha=1,...,p
\end{align*}
and linearly extend to $\Gamma\left( U_{jk}^t, \Theta_{\mathcal{F}_t} \right)$, and we set $\theta_{jk}(t):=  \sum_{\eta=1}^n \frac{\partial f_{jk}^\eta(\xi_k,t)}{\partial t}\frac{\partial}{\partial \xi_j^\eta}\in \Gamma\left( U_{jk}^t , \Theta_{M_t} \right)$. If we consider $\left( \theta_{jk}(t), A_{jk}(t)      \right)$ as an element of $\Gamma\left(U_j^t, \Theta_{M_t}\right)|_{U_{jk}^t}\bigoplus \Gamma\left(U_j^t, \mathscr{H}om_{\mathcal{O}_{M_t}}\left( \Theta_{\mathcal{F}_t}, \Theta_{\mathcal{F}_t} \right)\right) |_{U_{jk}^t}$, then we see that $\left\{ \left( \theta_{jk}(t)  , A_{jk}(t) \right) \right\}\in C^1\left( \mathcal{U}_t, \mathcal{E}_{\Theta_{\mathcal{F}_t}} \right)$. In fact, first $\theta_{jk}(t)=- \theta_{kj}(t)$, and second, by taking the derivative of $\sum_{\beta=1}^q r_{kj}^{\alpha\beta}(f_{jk}(\xi_k, t), t) r_{jk}^{\beta \gamma}(\xi_k, t) = \delta_\alpha^\gamma$ with respect to $t_k$, we have
\begin{align*}
&\sum_{\beta=1}^p \sum_{\eta=1}^n \frac{\partial r_{kj}^{\alpha\beta}}{\partial \xi_j^\eta}\frac{\partial f_{jk}^\eta}{\partial t} r_{jk}^{\beta \gamma}+ \sum_{\beta=1}^p \frac{\partial r_{kj}^{\alpha\beta}}{\partial t_j} r_{jk}^{\beta \gamma} +\sum_{\beta=1}^p r_{kj}^{\alpha\beta} \frac{\partial r_{jk}^{\beta \gamma}}{\partial t} = 0 \Longrightarrow  \left[ \theta_{jk}, r_{kj}^{\alpha\beta}\right] + \frac{\partial r_{kj}^{\alpha\beta}}{\partial t_j} + \sum_{\gamma, \delta=1}^q r_{kj}^{\alpha\gamma} \frac{\partial r_{jk}^{\gamma \delta}}{\partial t} r_{kj}^{\delta \beta}=0
\end{align*}
We note that $\left( \theta_{jk}(t), A_{jk}(t)\right)$ on $U_j^t$ be transformed on $U_k$ into $\left(\theta_{jk}(t), A_{jk}^k(t) \right)$, where
{\small{\begin{align*}
A_{jk}^k(T_k^\alpha)= A_{jk}(t)(T_k^\alpha) + \sum_{\beta=1}^p \left[\theta_{jk}(t), r_{kj}^{\alpha\beta} \right]T_j^\beta = \sum_{\beta, \gamma=1}^p r_{kj}^{\alpha\beta} \frac{\partial r_{jk}^{\beta \gamma}}{\partial t} T_k^\gamma +  \sum_{\beta=1}^p \left[\theta_{jk}(t), r_{kj}^{\alpha\beta} \right]T_j^\beta =- \sum_{\beta=1}^p \frac{\partial r_{kj}^{\alpha\beta}}{\partial t_j} T_j^\beta = - A_{kj}(T_k^\alpha)
\end{align*}}}

We claim that $\left( \left\{\alpha_j(t)\right\}, \left\{ \left(\theta_{jk}(t), A_{jk}(t) \right) \right\}\right) $ defines a $1$-cocycle in the above \v Cech resolution $(\ref{o7})$ of $\mathcal{E}_{\Theta_{\mathcal{F}_t}}^\bullet$. As in Part I, by taking the derivative of $(\ref{ta5})$ with respect to $t$, we see that $D_1'^t\left(\alpha_j(t) \right)=0$. On the other hand, by taking the derivative of $(\ref{ta4})$ with respect to $t$, we obtain as in Part I that
\begin{align*}
\left({\alpha}_k(t)- {\alpha}_j(t)\right)\left(T_j^\alpha\right) + \left[ \theta_{jk}(t) ,  T_j^\alpha \right] =\sum_{\beta=1}^p \frac{\partial r_{jk}^{\alpha\beta}}{\partial t} T_k^\beta \Longrightarrow \left(\alpha_k- \alpha_j(t) \right)(T_j^\alpha)+\left[ \theta_{jk}(t), T_j^\alpha \right]- A_{jk}(t)\left(T_j^\alpha\right)=0
\end{align*}
This implies that $\delta\left(\{\alpha_j(t)\} \right) + D_0'\left( \left\{\left(\theta_{jk}(t), A_{jk}(t) \right) \right\} \right)=0$. It remains to show that $\delta\left(\left\{\theta_{jk}(t), A_{jk}(t) \right\} \right)=0$. In fact, first we note that $\delta\left( \left\{\theta_{jk}(t) \right\} \right)=0$. Second, by taking the derivative of $(\ref{o1})$ with respect to $t$, we obtain that
{\small{\begin{align*}
&\frac{\partial r_{ik}^{\alpha \gamma}(\xi_k, t)}{\partial t} = \sum_{\beta=1}^q \sum_{\eta=1}^n \frac{\partial r_{ij}^{\alpha\beta}}{\partial \xi_j^\eta}\left(f_{jk}(\xi_k, t) , t\right) \frac{\partial f_{jk}^{\eta}(\xi_k, t)}{\partial t} r_{jk}^{\beta \gamma}(\xi_k,t) + \sum_{\beta=1}^q \frac{\partial r_{ij}^{\alpha\beta}}{\partial t}\left(f_{jk}(\xi_k,t),t\right) r_{jk}^{\beta \gamma}(\xi_k, t) + \sum_{\beta=1}^q r_{ij}^{\alpha\beta}\left(f_{jk}(\xi_k,t),t\right)\frac{\partial r_{jk}^{\beta \gamma}}{\partial t}(\xi_k,t)
\end{align*}}}
Then we have
\begin{align}\label{o8}
 \frac{\partial r_{ik}^{\alpha \gamma}}{\partial t}= \sum_{\beta=1}^q\left[ \theta_{jk}(t), r_{ij}^{\alpha\beta}\right] r_{jk}^{\beta \gamma} + \sum_{\beta=1}^q \frac{\partial r_{ij}^{\alpha\beta}}{\partial t}  r_{jk}^{\beta \gamma} + \sum_{\beta=1}^q r_{ij}^{\alpha\beta} \frac{\partial r_{jk}^{\beta \gamma}}{\partial t}
\end{align}
Since $A_{jk}(t)$ on $U_j^t$ is identified on $U_i^t$ with $A_{jk}^i\left(T_{i}^\alpha\right)=A_{jk}\left(T_{i}^\alpha\right) + \sum_{\beta=1}^p \left[\theta_{jk}(t),  r_{ij}^{\alpha\beta} \right] T_{j}^\beta$, $(\ref{o8})$ implies
\begin{align*}
&A_{ij}(t)\left(T_{i}^\alpha\right)- A_{ik}(t)\left(T_{i}^\alpha\right)+ A_{jk}^i(t)\left(T_{i}^\alpha\right)\\
&= \sum_{\beta=1}^q \frac{\partial r_{ij}^{\alpha\beta}}{\partial t} T_{j}^\beta - \sum_{\beta=1}^q \frac{\partial r_{ik}^{\alpha\beta}}{\partial t} T_{k}^\beta + \sum_{\beta=1}^q r_{ij}^{\alpha\beta} \frac{\partial r_{jk}^{\beta \gamma}}{\partial t} T_{k}^\gamma + \sum_{\beta=1}^p \left[\theta_{jk}(t),  r_{ij}^{\alpha\beta} \right] T_{j}^\beta = 0,
\end{align*}
so that $\delta\left( \left\{ \theta_{jk}(t), A_{jk}(t) \right\} \right)=0$. Hence we have the infinitesimal foliated deformation along $\frac{\partial}{\partial t}$ at $t$:
\begin{align}\label{o10}
&\frac{\partial \left( M_t, \Theta_{\mathcal{F}_t} \right)}{\partial t}= \left( \left\{ \alpha_j (t) \right\},   \left\{ \left(\theta_{jk}(t), A_{jk}(t) \right) \right\} \right) \in  \mathbb{H}^1\left( M_t , \mathcal{E}_{\Theta_{\mathcal{F}_t } }^\bullet \right) 
\end{align}

\subsection{Description of obstructed elements in \v Cech resolution of $\mathcal{E}_{\Theta_{\mathcal{F}_0}}^\bullet$} \label{o42}

\begin{definition}
We say that a compact foliated complex manifold $\left(M, \Theta_{\mathcal{F}_0}\right)$ with $\Theta_{\mathcal{F}_0}$ locally free is unobstructed in foliated deformations in terms of tangent sheaves if there is a foliated analytic family $\left( \mathcal{M}, \Theta_\mathcal{F}, B, \pi \right)$ in terms of tangent sheaves of deformations of $\pi^{-1}(0)=\left( M, \Theta_{\mathcal{F}_0} \right), 0\in B$ such that the foliated Kodaira-Spencer map $\varphi_0:T_0B\to \mathbb{H}^1\left( M, \mathcal{E}_{\Theta_{\mathcal{F}_0}}^\bullet \right)$ in terms of tangent sheaves is an isomorphism at $0 \in B$. Otherwise, we say that $\left( M, \Theta_{\mathcal{F}_0} \right)$ is obstructed in foliated deformations in terms of tangent sheaves.
\end{definition}

\begin{remark}
We note that by theorem of existence for deformations of foliated complex analytic structures in terms of tangent sheaves $($\textnormal{Part I}$)$, if a compact foliated complex manifold $\left(M, \Theta_{\mathcal{F}_0} \right)$ with $\Theta_{\mathcal{F}_0}$ locally free satisfies $\mathbb{H}^2\left(  M,  \mathcal{E}_{\Theta_{\mathcal{F}_0}}^\bullet \right)=0$, then there is a foliated complex analytic family $\left( \mathcal{M}, \Theta_\mathcal{F}, B, \pi \right)$ such that $\pi^{-1}(0)= \left(M, \Theta_{\mathcal{F}_0}\right), 0\in B$ and the foliated Kodaira-Spencer map in terms of tangent sheaves is an isomorphism at $0\in B$, so that $\left(M, \Theta_{\mathcal{F}_0}\right)$ is unobstructed in foliated deformations in terms of tangent sheaves.
\end{remark}

As we find an obstructed element $\theta\in H^1(M, \Theta_M)$ (see \cite{Kod05} p.209-214) in deformations of complex structures in order to show that a given compact complex manifold $M$ is obstructed, analogously one way to determine obstructedness of a compact foliated complex manifold $(M, \Theta_{\mathcal{F}_0})$ with $\Theta_{\mathcal{F}_0}$ locally free in deformations of foliated complex analytic structures is to find an obstructed element as we define in the following way:

\begin{definition}
Let $\left(M, \Theta_{\mathcal{F}_0}\right)$ be a compact foliated complex manifold with $\Theta_{\mathcal{F}_0}$ locally free. We say that $a\in \mathbb{H}^1\left( M, \mathcal{E}_{\Theta_{\mathcal{F}_0}}^\bullet \right)$ is an obstructed element in foliated deformations of $\left(M, \Theta_{\mathcal{F}_0}\right)$ in terms of tangent sheaves if there is no foliated analytic family $\left( \mathcal{M}, \Theta_\mathcal{F}, B, \pi \right)$ in terms of tangent sheaves of deformations of $\left(M, \Theta_{\mathcal{F}_0}\right)= \pi^{-1}(0), 0\in B$ such that $a\in \varphi_0\left(T_0 M \right)$ where $\varphi_0: T_0M  \to \mathbb{H}^1\left( M, \mathcal{E}_{\Theta_{\mathcal{F}_0}}^\bullet \right)$ is the foliated Kodaira-Spencer map at $0\in B$ of $\left(\mathcal{M}, \Theta_\mathcal{F}, B, \pi \right)$ in terms of tangent sheaves.
\end{definition}

Hence if there is an obstructed element $a\in \mathbb{H}^1\left( M, \mathcal{E}_{\Theta_{\mathcal{F}_0}}^\bullet \right)$ in foliated deformations in terms of tangent sheaves, $\left( M, \Theta_{\mathcal{F}_0} \right)$ is obstructed. As $\theta\in H^1\left( M, \Theta_M \right)$ is an obstructed element in deformations of complex structures if $\left[ \theta, \theta \right]\ne 0 \in H^2\left( M, \Theta_M \right)$, we now describe a condition when $a\in \mathbb{H}^1\left( M, \mathcal{E}_{\Theta_{\mathcal{F}_0}}^\bullet \right)$ is an obstructed element in foliated deformations in terms of tangent sheaves. We extend the arguments in \cite{Kod05} p.210-214 in the context of a foliated analytic family in terms of tangent sheaves.

Let $\left( \mathcal{M}, \Theta_\mathcal{F}, B, \pi \right)$ with $0\in B \subset \mathbb{C}$ and $\Theta_{\mathcal{F}}$ locally free be a foliated analytic family in terms of tangent sheaves such that $\pi^{-1}(0)=\left(M, \Theta_{\mathcal{F}_0} \right)$. Take a small disk $\Delta$ with centre $0$ such that $0\in \Delta \subset B$, and represent $\left(\mathcal{M}_\Delta, \Theta_{\mathcal{F}_\Delta} \right)=\pi^{-1}(\Delta)$ in the form
\begin{align*}
\left(\mathcal{M}_\Delta, \Theta_{\mathcal{F}_\Delta} \right) = \bigcup_{j=1}^l \left( U_j\times \Delta, \Theta_{\mathcal{F}_\Delta} \right)
\end{align*}
where each $U_j$ is a polydisk, and $\Theta_{\mathcal{F}_\Delta}$ on $U_j\times \Delta$ is generated by $T_j^\alpha(\xi_j,t):=\sum_{\beta=1}^n T_j^{\alpha\beta}(\xi_j,t)\frac{\partial}{\partial \xi_j^\beta},\alpha=1,...,p$ such that $\left[ T_j^\alpha(\xi_j,t), T_j^\beta(\xi_j,t)\right]=\sum_{\gamma=1}^p g_{j\alpha\beta}^\gamma(\xi_j,t)T_j^\gamma(\xi_j,t)$ with $g_{j\alpha\beta}^\gamma(\xi_j,t)=-g_{j\alpha\beta}^\gamma(\xi_j,t)$, and $(\xi_j,t)\in U_j\times \Delta$ and $(\xi_k,t)\in U_k\times \Delta$ are the same points on $\mathcal{M}_\Delta$ if $\xi_j^\alpha= f_{jk}^\alpha(\xi_k,t)$ for $\alpha=1,...,n$. Here each $f_{jk}^\alpha(\xi_k,t)=f_{jk}^\alpha(\xi_k^1,..., \xi_k^n, t)$ is a holomorphic function of $\xi_k^1,..., \xi_k^n,t$ defined on $\left(U_k\times \Delta \right) \cap \left( U_j\times \Delta \right)\ne \emptyset$. On the other hand, $T_j^\alpha(\xi_j,t)= \sum_{\beta=1}^p r_{jk}^{\alpha\beta}\left(\xi_j,t \right) T_k^\beta(\xi_k,t)$ on $\left(U_k\times \Delta \right) \cap \left( U_j\times \Delta \right)\ne \emptyset$, and $r_{ik}^{\alpha\gamma}\left(\xi_k,t \right)= \sum_{\beta=1}^p r_{ij}^{\alpha\beta}(\xi_j,t) r_{jk}^{\beta \gamma}(\xi_k,t)$ on $\left( U_i\times \Delta \right) \cap \left( U_j\times \Delta \right) \cap \left( U_k\times \Delta \right)\ne \emptyset$. 

As in $(\ref{o10})$, the infinitesimal foliated deformation $\frac{\partial \left( M_t, \Theta_{\mathcal{F}_t} \right)}{\partial t}\in \mathbb{H}^1\left( M_t, \mathcal{E}_{\Theta_{\mathcal{F}_t} }^\bullet \right)$ is the cohomology class of the $1$-cocycle $\left( \left\{\alpha_j(t)\right\}, \left\{\left(\theta_{jk}(t), A_{jk}(t) \right) \right\} \right)\in C^0\left( \mathcal{U}_t, \mathscr{H}om_{\mathcal{O}_{M_t}}\left( \Theta_{\mathcal{F}_t}, \Theta_{M_t} \right) \right)\bigoplus C^1\left( \mathcal{U}_t, \mathcal{E}_{\Theta_{\mathcal{F}_t}} \right)$, where $\mathcal{U}_t=\left\{U_j^t=U_j\times t \right\}$ such that
\begin{align*}
\alpha_j(t):\Gamma\left( U_j^t, \Theta_{\mathcal{F}_t} \right)&\to \Gamma\left( U_j^t, \Theta_{M_t} \right) \\
 T_j^\alpha(\xi_j,t) &\mapsto  \alpha_j^\alpha(\xi_j,t):=\sum_{\beta=1}^n \alpha_j^{\alpha\beta}(\xi_j,t)\frac{\partial}{\partial \xi_j^\beta},\,\,\,\,\,\,\,\,\textnormal{where}\,\,\,\alpha_j^{\alpha\beta}(\xi_j, t)= - \frac{\partial T_j^{\alpha\beta}(\xi_j,t)}{\partial t},\alpha=1,...,p
\end{align*}
and the vector field
\begin{align*}
\theta_{jk}(t)=\sum_{\alpha=1}^n \theta_{jk}^\alpha(\xi_j,t)\frac{\partial}{\partial \xi_j^\alpha},\,\,\,\,\,\,\,\,\textnormal{where}\,\,\,\,\theta_{jk}^\alpha(\xi_j,t)=\frac{\partial f_{jk}(\xi_k,t)}{\partial t},\,\,\,\xi_k=f_{kj}(\xi_j,t)
\end{align*}
and
\begin{align*}
A_{jk}(t):\Gamma\left( U_j^t, \Theta_{\mathcal{F}_t} \right)&\to \Gamma\left( U_j^t, \Theta_{\mathcal{F}_t} \right) \\
     T_j^\alpha(\xi_j,t)&\mapsto \sum_{\gamma=1}^p A_{jk}^{\alpha\gamma}(\xi_j,t)T_j^\gamma(\xi_j,t),\,\,\textnormal{where}\,\,\, A_{jk}^{\alpha\gamma}(\xi_j,t)=\sum_{\beta=1}^p\frac{\partial r_{jk}^{\alpha\beta}(\xi_k,t)}{\partial t}r_{kj}^{\beta \gamma}(\xi_j,t),\xi_k=f_{kj}(\xi_j,t),\alpha=1,...,p
\end{align*}

On $\left(U_k\times \Delta \right)\cap \left(U_j\times \Delta \right) \ne \emptyset$, we have the equalities
\begin{align}
f_{ik}^\alpha(\xi_k,t)&=f_{ij}^\alpha\left(f_{jk}(\xi_k,t),t \right),\,\,\,\alpha=1,...,n \label{o11}\\
\sum_{\gamma=1}^n T_j^{\alpha\gamma}\left(f_{jk}(\xi_k,t),t\right)&=\sum_{\beta=1}^p\sum_{\eta,\gamma=1}^n r_{jk}^{\alpha\beta}(\xi_k,t) T_k^{\beta\eta}(\xi_k,t)\frac{\partial f_{jk}^\gamma(\xi_k,t)}{\partial \xi_k^\eta},\,\,\,\alpha=1,...,p, \gamma=1,...,n \label{o12}
\end{align}
and on $\left(U_i\cap \Delta \right)\cap \left(  U_j \times \Delta \right)\cap \left(U_k\times \Delta \right)\ne \emptyset$, we have the equalities
\begin{align}
r_{ik}^{\alpha\gamma}\left(\xi_k,t \right) &= \sum_{\beta=1}^p r_{ij}^{\alpha\beta}\left(f_{jk}(\xi_k,t),t \right) r_{jk}^{\beta \gamma}(\xi_k,t),\,\,\,\,\,\alpha,\gamma=1,...,p. \label{o13}
\end{align}
and on $U_j\times \Delta$, we have the equalities 
\begin{align}
\left[ T_j^\alpha(\xi_j,t), T_j^\beta(\xi_j,t) \right]= \sum_{\gamma=1}^q g_{j\alpha\beta}^\gamma(\xi_j,t) T_j^\gamma(\xi_j,t),\,\,\,\,\,\,\alpha,\beta=1,...,p \label{o14}
\end{align}
By differentiating both sides of $(\ref{o11}),(\ref{o12}),(\ref{o13})$ and $(\ref{o14})$ in $t$, we have the equalities
\begin{align}
&\theta_{ik}^\alpha(\xi_i,t)=\theta_{ij}^\alpha(\xi_i,t) + \sum_{\beta=1}^n \frac{\partial \xi_i^\alpha}{\partial \xi_j^\beta}\theta_{jk}^\beta(\xi_j,t),\,\,\,\alpha=1,...,n. \label{o15}\\
&\sum_{\eta=1}^n\frac{\partial T_j^{\alpha\gamma}}{\partial \xi_j^\eta}\left( f_{jk}(\xi_k,t),t\right)\frac{\partial f_{jk}^\eta}{\partial t} +\frac{\partial T_j^{\alpha\gamma}}{\partial t}\left( f_{jk}(\xi_k,t),t \right),\,\,\,\,\,\alpha  =1,...,p, \gamma=1,...,n \label{o16}\\
&= \sum_{\beta=1}^p\sum_{\eta=1}^n \frac{\partial r_{jk}^{\alpha\beta}}{\partial t} T_k^{\beta\eta}\frac{\partial f_{jk}^\gamma}{\partial \xi_k^\eta} +\sum_{\beta=1}^p\sum_{\eta=1}^n r_{jk}^{\alpha\beta}\frac{\partial T_k^{\beta\eta}}{\partial t} \frac{\partial f_{jk}^\gamma}{\partial \xi_k^\eta} +\sum_{\beta=1}^p\sum_{\eta=1}^n r_{jk}^{\alpha\beta}T_k^{\beta\eta}\frac{\partial}{\partial \xi_k^\eta}\left(\frac{\partial f_{jk}^\gamma}{\partial t} \right) \notag
\end{align}
and for $\alpha,\gamma=1,...,p$,
{\small{\begin{align}
\frac{\partial r_{ik}^{\alpha \gamma}}{\partial t} = \sum_{\beta=1}^p \sum_{\eta=1}^n \frac{\partial r_{ij}^{\alpha\beta}}{\partial \xi_j^\eta}\left(f_{jk}(\xi_k, t) , t\right) \frac{\partial f_{jk}^{\eta}}{\partial t} r_{jk}^{\beta \gamma} + \sum_{\beta=1}^p \frac{\partial r_{ij}^{\alpha\beta}}{\partial t}\left(f_{jk}(\xi_k,t),t\right) r_{jk}^{\beta \gamma} + \sum_{\beta=1}^p r_{ij}^{\alpha\beta}\left(f_{jk}(\xi_k,t),t\right)\frac{\partial r_{jk}^{\beta \gamma}}{\partial t} \label{o17}
\end{align}}}
\begin{align}
\left[\frac{\partial T_j^\alpha}{\partial t}, T_j^\beta\right] + \left[ T_j^\alpha, \frac{\partial T_j^\beta}{\partial t} \right] =\sum_{\gamma=1}^p \frac{\partial g_{j\alpha\beta}^\gamma}{\partial t} T_j^\gamma + \sum_{\gamma=1}^p g_{j\alpha\beta}^\gamma \frac{\partial T_j^\gamma}{\partial t} ,\,\,\,\alpha,\beta=1,...,p \label{o18}
\end{align}

Next we will take the derivative of $(\ref{o15}), (\ref{o16}), (\ref{o18})$ in $t$ as functions of $\xi_j^1,..., \xi_j^n,t$ and $(\ref{o17})$ in $t$ as functions of $\xi_k^1,..., \xi_k^n,t$ again. As in \cite{Kod05} p.211 we sometimes write $\left( \frac{\partial}{\partial t}\right)_j$ instead of $\left(\frac{\partial}{\partial t} \right)$ in order to make explicit that $\frac{\partial}{\partial t}$ denotes the differentiation of a function of $\xi_j^1,..., \xi_j^n,t$ with respect to $t$. On $\left( U_j \times \Delta \right)\cap \left( U_k \times \Delta \right)\ne \emptyset$, we have the following equalities: $\frac{\partial}{\partial \xi_k^\alpha}=\sum_{\beta=1}^n \frac{\partial \xi_j^\beta}{\partial \xi_k^\alpha} \frac{\partial}{\partial \xi_j^\beta}$, where $\frac{\partial \xi_j^\beta}{\partial \xi_k^\alpha}=\frac{\partial f_{jk}^\beta(\xi_k,t)}{\partial \xi_k^\alpha}$, and
\begin{align}\label{o53}
\left( \frac{\partial}{\partial t} \right)_k  = \sum_{s=1}^n \frac{\partial f_{jk}^s(\xi_k,t)}{\partial t}\frac{\partial}{\partial \xi_j^s} + \left( \frac{\partial}{\partial t} \right)_j = \sum_{s=1}^n \theta_{jk}^s(z_j,t)\frac{\partial}{\partial \xi_j^s} + \left( \frac{\partial}{\partial t} \right)_j 
\end{align}
By taking the derivative of $(\ref{o15})$ with respect to $t$, and putting
\begin{align*}
\dot{\theta}_{ij}(t)= \sum_{\alpha=1}^n \dot{\theta}_{ij}^\alpha (\xi_i,t)\frac{\partial}{\partial z_i^\alpha} , \,\,\,\,\,\textnormal{where}\,\,\,\dot{\theta}_{ij}^\alpha(\xi_i,t)=\frac{\partial \theta_{ij}^\alpha(\xi_i,t)}{\partial t}
\end{align*}
we obtain (for the detail, see \cite{Kod05} p.211-213)
\begin{align}\label{o36}
\dot{\theta}_{ij}(t) -\dot{\theta}_{ik}(t) + \dot{\theta}_{jk} (t) = \left[ \theta_{ij}(t), \theta_{jk}(t) \right]
\end{align}

By taking the derivative of $(\ref{o18})$ with respect to $t$, and putting
\begin{align*}
\dot{\alpha}_j^{\alpha}(\xi_j,t)= \frac{\partial \alpha_j^{\alpha}(\xi_j,t)}{\partial t}= \sum_{\beta=1}^n \frac{\partial \alpha_j^{\alpha\beta}(\xi_j,t)}{\partial t} \frac{\partial}{\partial \xi_j^\beta} = - \sum_{\beta=1}^n \frac{\partial^2 T_j^{\alpha\beta}(\xi_j,t)}{\partial t^2} \frac{\partial}{\partial \xi_j^\beta}, \,\,\,\alpha=1,...,p
\end{align*}
we obtain 
\begin{align*}
-\left[ \dot{\alpha}_j^\alpha(t),  T_j^\beta \right] +2\left[ \alpha_j^\alpha(t), \alpha_j^\beta(t)\right] - \left[ T_j^\alpha, \dot{\alpha}_j^\beta(t)\right]= \sum_{\gamma=1}^p \frac{\partial^2 g_{j\alpha\beta}^\gamma}{\partial t^2} T_j^\gamma - 2\sum_{\gamma=1}^p \frac{\partial g_{j\alpha\beta}^\gamma}{\partial t}\alpha_j^\gamma(t) - \sum_{\gamma=1}^p g_{j\alpha\beta}^\gamma \dot{\alpha}_j^\gamma(t)
\end{align*}
By defining $\dot{\alpha}_j(t)$ by
\begin{align*}
\dot{\alpha}_j(t):\Gamma\left( U_j^t, \Theta_{\mathcal{F}_t} \right) &\to \left( U_j^t, \Theta_{M_t} \right)  \\
T_j^\alpha &\mapsto \dot{\alpha}_j^\alpha(t)=\dot{\alpha}_j^\alpha(\xi_j,t)
\end{align*}
by linearly extending to $\Gamma\left( U_j^t, \Theta_{\mathcal{F}_t} \right)$, we obtain
\begin{align}\label{o37}
\left[ T_j^\alpha, \dot{\alpha}_j^\beta (t) \right] - \left[ T_j^\beta, \dot{\alpha}_j^\alpha(t)\right] - \sum_{\gamma=1}^p g_{j\alpha\beta}^\gamma \dot{\alpha}_j^\gamma (t)  = 2\left[ \alpha_j^\alpha(t), \alpha_j^\beta(t) \right]  + 2\alpha_j(t)\left( \sum_{\gamma=1}^p \frac{\partial g_{j\alpha\beta}^\gamma}{\partial t} T_j^\gamma \right) - \sum_{\gamma=1}^p \frac{\partial^2 g_{j\alpha\beta}^\gamma}{\partial t^2} T_j^\gamma
\end{align}

On the other hand, let us define 
\begin{align}
\dot{A}_{jk}(t):\Gamma\left( U_{jk}^t,  \Theta_{\mathcal{F}_t} \right)\to &\Gamma\left( U_{jk}^t, \Theta_{\mathcal{F}_t} \right) \label{o35} \\
 T_j^\alpha \mapsto &\sum_{\gamma=1}^p \frac{\partial}{\partial t_j}\left( A_{jk}^{\alpha \gamma}(\xi_j,t)\right) T_j^\gamma = \sum_{\beta, \gamma=1}^p \frac{\partial}{\partial t_j}\left( \frac{\partial r_{jk}^{\alpha \beta}(\xi_k,t)}{\partial t} r_{kj}^{\beta \gamma}(\xi_j,t) \right) T_j^\gamma  \notag
\end{align}
by linearly extending to $\Gamma\left( U_{jk}^t, \Theta_{\mathcal{F}_t} \right)$. 

By considering $\left(\dot{\theta}_{jk}(t), \dot{A}_{jk}(t) \right)\in \Gamma\left(U_{jk}^t, \Theta_{M_t} \right)|_{U_j^t} \bigoplus \Gamma\left(U_{jk}^t, \mathscr{H}om_{\mathcal{O}_{M_t}} \left( \Theta_{\mathcal{F}_t}, \Theta_{\mathcal{F}_t} \right) |_{U_j^t}\right)$, we claim that $\left\{\left( \dot{\theta}_{jk}(t), \dot{A}_{jk}(t) \right) \right\} \in C^1\left(\mathcal{U}_t, \mathcal{E}_{ \Theta_{\mathcal{F}_t}} \right)$. In fact, first we note that $\dot{\theta}_{jk}(t)=-\dot{\theta}_{kj}(t)$. Second, by taking the derivative of $\sum_{\beta=1}^q r_{kj}^{\alpha\beta}(f_{jk}(\xi_k, t), t) r_{jk}^{\beta \gamma}(\xi_k, t) = \delta_\alpha^\gamma$ with respect to $t_k$, we have
\begin{align}\label{o24}
&\sum_{\beta=1}^p \sum_{\eta=1}^n \frac{\partial r_{kj}^{\alpha\beta}}{\partial \xi_j^\eta}\frac{\partial f_{jk}^\eta}{\partial t} r_{jk}^{\beta \gamma}+ \sum_{\beta=1}^p \frac{\partial r_{kj}^{\alpha\beta}}{\partial t_j} r_{jk}^{\beta \gamma} +\sum_{\beta=1}^p r_{kj}^{\alpha\beta} \frac{\partial r_{jk}^{\beta \gamma}}{\partial t} = 0 \Longrightarrow  \left[ \theta_{jk}, r_{kj}^{\alpha\beta}\right] + \frac{\partial r_{kj}^{\alpha\beta}}{\partial t_j} + \sum_{\gamma, \delta=1}^q r_{kj}^{\alpha\gamma} \frac{\partial r_{jk}^{\gamma \delta}}{\partial t} r_{kj}^{\delta \beta}=0
\end{align}
By taking the derivative with respect to $t_j$ again, we have
\begin{align}\label{o20}
0=&\sum_{\beta=1}^p \sum_{\eta=1}^n  \frac{\partial}{\partial \xi_j^\eta} \left(\frac{\partial r_{kj}^{\alpha\beta}}{\partial t_j} \right) \frac{\partial f_{jk}^\eta}{\partial t} r_{jk}^{\beta \gamma} + \sum_{\beta=1}^p \sum_{\eta=1}^n \frac{\partial r_{kj}^{\alpha\beta}}{\partial \xi_j^\eta}\dot{\theta}_{jk}^\eta r_{jk}^{\beta \gamma} + \sum_{\beta=1}^p \sum_{\eta=1}^n \frac{\partial r_{kj}^{\alpha\beta}}{\partial \xi_j^\eta} \frac{\partial f_{jk}^\eta}{\partial t} \frac{\partial r_{jk}^{\beta \gamma}}{\partial t_j}\\
&+ \sum_{\beta=1}^p \frac{\partial^2 r_{kj}^{\alpha \beta}}{\partial t_j^2} r_{jk}^{\beta \gamma} + \sum_{\beta=1}^p \frac{\partial r_{kj}^{\alpha\beta}}{\partial t_j} \frac{\partial r_{jk}^{\beta\gamma}}{\partial t_j} + \sum_{\beta=1}^p \frac{\partial r_{kj}^{\alpha\beta}}{\partial t_j}\frac{\partial r_{jk}^{\beta \gamma}}{\partial t}+ \sum_{\beta=1}^p r_{kj}^{\alpha\beta} \frac{\partial}{\partial t_j}\left( \frac{\partial r_{jk}^{\beta\gamma}}{\partial t_k}\right) \notag \\
&= \sum_{\beta=1}^p \frac{\partial}{\partial t_k}\left( \frac{\partial r_{kj}^{\alpha\beta}}{\partial t_j} \right)  r_{jk}^{\beta \gamma} + \sum_{\beta=1}^p \left[\dot{\theta}_{jk}, r_{kj}^{\alpha\beta} \right] r_{jk}^{\beta \gamma}  + \sum_{\beta=1}^p \frac{\partial r_{kj}^{\alpha\beta}}{\partial t_k} \frac{\partial r_{jk}^{\beta\gamma}}{\partial t_j} + \sum_{\beta=1}^p \frac{\partial r_{kj}^{\alpha\beta}}{\partial t_j}\frac{\partial r_{jk}^{\beta \gamma}}{\partial t}+ \sum_{\beta=1}^p r_{kj}^{\alpha\beta} \frac{\partial}{\partial t_j}\left( \frac{\partial r_{jk}^{\beta\gamma}}{\partial t_k}\right)     \notag
\end{align}
We note that $\left( \dot{\theta}_{jk}(t), \dot{A}_{jk}(t) \right)$ on $U_j$ is transformed on $U_k$ into $\left( \dot{\theta}_{jk} (t) , \dot{A}_{jk}^k(t) \right)$, where
\begin{align*}
\dot{A}_{jk}^k(t)(T_k^\alpha)&= \dot{A}_{jk}(t)(T_k^\alpha) + \sum_{\beta=1}^p \left[\dot{\theta}_{jk}(t), r_{kj}^{\alpha \beta} \right] T_j^\beta = \sum_{\beta,\gamma, \delta=1}^p r_{kj}^{\alpha\beta}  \frac{\partial}{\partial t_j}\left( \frac{\partial r_{jk}^{\beta \gamma}}{\partial t} r_{kj}^{\gamma \delta} \right) T_j^\delta + \sum_{\beta=1}^p \left[\dot{\theta}_{jk}(t), r_{kj}^{\alpha \beta} \right] T_j^\beta  \\
&= \sum_{\beta, \gamma=1}^p r_{kj}^{\alpha \beta} \frac{\partial}{\partial t_j}\left(\frac{\partial r_{jk}^{\beta \gamma}}{\partial t} \right) T_k^\gamma  + \sum_{\beta, \gamma, \delta=1}^p r_{kj}^{\alpha\beta} \frac{\partial r_{jk}^{\beta \gamma}}{\partial t_k} \frac{\partial r_{kj}^{\gamma \delta}}{\partial t_j} T_j^\delta + \sum_{\beta=1}^p \left[\dot{\theta}_{jk}(t), r_{kj}^{\alpha \beta} \right] T_j^\beta  \\
&= \sum_{\beta, \gamma=1}^p r_{kj}^{\alpha \beta} \frac{\partial}{\partial t_j}\left(\frac{\partial r_{jk}^{\beta \gamma}}{\partial t} \right) T_k^\gamma  + \sum_{\beta, \gamma=1}^p  \frac{\partial r_{kj}^{\alpha \beta}}{\partial t_k} \frac{\partial r_{jk}^{\beta \gamma}}{\partial t_j} T_k^\gamma + \sum_{\beta=1}^p \left[\dot{\theta}_{jk}(t), r_{kj}^{\alpha \beta} \right] T_j^\beta 
\end{align*}
Then $(\ref{o20})$ implies that 
\begin{align*}
\dot{A}_{kj}(t)(T_k^\alpha)= \sum_{\beta, \gamma=1}^p \frac{\partial}{\partial t_k}\left( \frac{\partial r_{kj}^{\alpha\beta}}{\partial t_j} r_{jk}^{\beta \gamma} \right) T_k^\gamma = \sum_{\beta=1}^p \frac{\partial}{\partial t_k}\left(\frac{\partial r_{kj}^{\alpha\beta}}{\partial t_j} \right) T_j^\beta + \sum_{\beta,\gamma=1}^p \frac{\partial r_{kj}^{\alpha\beta}}{\partial t_j} \frac{\partial r_{jk}^{\beta \gamma}}{\partial t_k} T_k^\gamma = - \dot{A}_{jk}^k(t)(T_k^\alpha)
\end{align*}
This proves $\left\{\left( \dot{\theta}_{jk}(t), \dot{A}_{jk}(t) \right) \right\} \in C^1\left(\mathcal{U}_t, \mathcal{E}_{ \Theta_{\mathcal{F}_t}} \right)$.

On the other hand, let us take the derivative of $(\ref{o16})$ with respect to $t_j$:
{\small{\begin{align}\label{o21}
& \sum_{\eta=1}^n \frac{\partial}{\partial \xi_j^\eta}\left(\frac{\partial T_j^{\alpha\gamma}}{\partial t_j} \right) \frac{\partial f_{jk}^\eta}{\partial t_k} + \sum_{\eta=1}^n \frac{\partial T_j^{\alpha \gamma}}{\partial \xi_j^\eta}\dot{\theta}_{jk}^\eta + \frac{\partial^2 T_j^{\alpha \gamma}}{\partial t_j^2}\\
&=\left( \frac{\partial}{\partial t}\right)_j\left(\sum_{\beta=1}^p\sum_{\eta=1}^n \frac{\partial r_{jk}^{\alpha\beta}}{\partial t} T_k^{\beta\eta}\frac{\partial f_{jk}^\gamma}{\partial \xi_k^\eta} +\sum_{\beta=1}^p\sum_{\eta=1}^n r_{jk}^{\alpha\beta} \frac{\partial T_k^{\beta\eta}}{\partial t} \frac{\partial f_{jk}^\gamma}{\partial \xi_k^\eta} + \sum_{\eta=1}^nT_j^{\alpha\eta}\frac{\partial}{\partial \xi_j^\eta}\left(\frac{\partial f_{jk}^\gamma}{\partial t} \right) \right) \notag \\
&= - \sum_{\delta,\eta=1}^n\sum_{\beta=1}^p \frac{\partial f_{jk}^\delta}{\partial t} \frac{\partial}{\partial \xi_j^\delta}\left(\frac{\partial r_{jk}^{\alpha\beta}}{\partial t} T_k^{\beta\eta}\frac{\partial f_{jk}^\gamma}{\partial \xi_k^\eta}\right) -\sum_{\delta,\eta=1}^n\sum_{\beta=1}^p\frac{\partial f_{jk}^\delta}{\partial t}\frac{\partial}{\partial \xi_j^\delta} \left(r_{jk}^{\alpha\beta} \frac{\partial T_k^{\beta\eta}}{\partial t} \frac{\partial f_{jk}^\gamma}{\partial \xi_k^\eta}\right)+\sum_{\eta=1}^n \frac{\partial T_j^{\alpha\eta}}{\partial t_j}\frac{\partial}{\partial \xi_j^\eta}\left( \frac{\partial f_{jk}^\gamma}{\partial t} \right)  + \sum_{\eta=1}^n T_j^{\alpha \eta} \frac{\partial \dot{\theta}_{jk}^\gamma}{\partial \xi_j^\eta} \notag \\
&+ \sum_{\beta=1}^p \sum_{\eta=1}^n \frac{\partial^2 r_{jk}^{\alpha\beta}}{\partial t^2} T_k^{\beta\eta} \frac{\partial f_{jk}^\gamma}{\partial \xi_k^\eta} +2\sum_{\beta=1}^p\sum_{\eta=1}^n \frac{\partial r_{jk}^{\alpha\beta}}{\partial t}\frac{\partial T_k^{\beta\eta}}{\partial t}\frac{\partial f_{jk}^\gamma}{\partial \xi_k^\eta}+\sum_{\beta=1}^p\sum_{\eta=1}^n \frac{\partial r_{jk}^{\alpha\beta}}{\partial t} T_k^{\beta\eta}\frac{\partial}{\partial \xi_k^\eta}\left(\frac{\partial f_{jk}^\gamma}{\partial t} \right) \notag \\
&+\sum_{\beta=1}^p\sum_{\eta=1}^n r_{jk}^{\alpha\beta} \frac{\partial^2 T_k^{\beta\eta}}{\partial t^2}\frac{\partial f_{jk}^\gamma}{\partial \xi_k^\eta}+\sum_{\beta=1}^p\sum_{\eta=1}^n r_{jk}^{\alpha\beta} \frac{\partial T_k^{\beta\eta}}{\partial t}\frac{\partial}{\partial \xi_k^\eta}\left( \frac{\partial f_{jk}^\gamma}{\partial t}\right) \notag
\end{align}}}
Then by summing $\sum_{\gamma=1}^n \frac{\partial}{\partial \xi_j^\gamma}$ on both sides of $(\ref{o21})$, we have
{\small{\begin{align} \label{o22}
& \sum_{\eta,\gamma=1}^n  \frac{\partial f_{jk}^\eta}{\partial t_k}\frac{\partial}{\partial \xi_j^\eta}\left(\frac{\partial T_j^{\alpha\gamma}}{\partial t_j} \right) \frac{\partial}{\partial \xi_j^\gamma}+ \sum_{\eta,\gamma=1}^n \dot{\theta}_{jk}^\eta \frac{\partial T_j^{\alpha \gamma}}{\partial \xi_j^\eta} \frac{\partial}{\partial \xi_j^\gamma} + \sum_{\gamma=1}^n \frac{\partial^2 T_j^{\alpha \gamma}}{\partial t_j^2}\frac{\partial}{\partial \xi_j^\gamma}\\
&=     - \sum_{\delta,\eta,\gamma=1}^n\sum_{\beta=1}^p \frac{\partial f_{jk}^\delta}{\partial t} \frac{\partial}{\partial \xi_j^\delta}\left(\frac{\partial r_{jk}^{\alpha\beta}}{\partial t} T_k^{\beta\eta}\frac{\partial f_{jk}^\gamma}{\partial \xi_k^\eta}\right)\frac{\partial}{\partial \xi_j^\gamma} -\sum_{\delta,\eta,\gamma=1}^n\sum_{\beta=1}^p\frac{\partial f_{jk}^\delta}{\partial t}\frac{\partial}{\partial \xi_j^\delta} \left(r_{jk}^{\alpha\beta} \frac{\partial T_k^{\beta\eta}}{\partial t} \frac{\partial f_{jk}^\gamma}{\partial \xi_k^\eta}\right)\frac{\partial}{\partial \xi_j^\gamma}+\sum_{\eta,\gamma=1}^n \frac{\partial T_j^{\alpha\eta}}{\partial t_j}\frac{\partial}{\partial \xi_j^\eta}\left( \frac{\partial f_{jk}^\gamma}{\partial t} \right)\frac{\partial}{\partial \xi_j^\gamma} \notag \\
&  + \sum_{\eta,\gamma=1}^n T_j^{\alpha \eta} \frac{\partial \dot{\theta}_{jk}^\gamma}{\partial \xi_j^\eta} \frac{\partial}{\partial \xi_j^\gamma}   + \sum_{\beta=1}^p \sum_{\eta=1}^n \frac{\partial^2 r_{jk}^{\alpha\beta}}{\partial t^2} T_k^{\beta\eta} \frac{\partial }{\partial \xi_k^\eta} +2\sum_{\beta=1}^p\sum_{\eta=1}^n \frac{\partial r_{jk}^{\alpha\beta}}{\partial t}\frac{\partial T_k^{\beta\eta}}{\partial t}\frac{\partial }{\partial \xi_k^\eta}+\sum_{\beta=1}^p\sum_{\eta,\gamma=1}^n \frac{\partial r_{jk}^{\alpha\beta}}{\partial t} T_k^{\beta\eta}\frac{\partial}{\partial \xi_k^\eta}\left(\frac{\partial f_{jk}^\gamma}{\partial t} \right)  \frac{\partial}{\partial \xi_j^\gamma}      \notag \\
&+\sum_{\beta=1}^p\sum_{\eta=1}^n r_{jk}^{\alpha\beta} \frac{\partial^2 T_k^{\beta\eta}}{\partial t^2}\frac{\partial}{\partial \xi_k^\eta}+\sum_{\beta=1}^p\sum_{\eta,\gamma=1}^n r_{jk}^{\alpha\beta} \frac{\partial T_k^{\beta\eta}}{\partial t}\frac{\partial}{\partial \xi_k^\eta}\left( \frac{\partial f_{jk}^\gamma}{\partial t}\right) \frac{\partial}{\partial \xi_j^\gamma} \notag 
\end{align}
Then $(\ref{o22})$ is equivalent to
\begin{align} \label{o27}
&\left[ \theta_{jk}(t),  \frac{\partial T_j^\alpha }{\partial t_j} \right] +  \left[ \dot{\theta}_{jk}(t),  T_j^\alpha \right]+   \frac{\partial^2 T_j^\alpha }{\partial t_j^2} \\
&= -\left[ \theta_{jk}(t),  \sum_{\beta=1}^p \frac{\partial r_{jk}^{\alpha\beta}}{\partial t} T_k^\beta   \right] -\left[ \theta_{jk}(t),  \sum_{\beta=1}^p  r_{jk}^{\alpha\beta}\frac{\partial T_k^\beta}{\partial t} \right]  \notag     + \sum_{\beta=1}^p \frac{\partial^2 r_{jk}^{\alpha\beta}}{\partial t^2} T_k^{\beta} +2\sum_{\beta=1}^p\frac{\partial r_{jk}^{\alpha\beta}}{\partial t}\frac{\partial T_k^\beta}{\partial t} +\sum_{\beta=1}^p r_{jk}^{\alpha\beta} \frac{\partial^2 T_k^{\beta}}{\partial t^2} \notag 
\end{align}}}

\begin{lemma}
\begin{align}\label{o23}
&\dot{\alpha}_k(t)- \dot{\alpha}_j(t) + D_0'\left(\left(\dot{\theta}_{jk}(t), \dot{A}_{jk}(t) \right)  \right)=\dot{\alpha}_k(t)- \dot{\alpha}_j(t) + \left(\mathcal{L}_{\dot{\theta}_{jk}(t)} - \dot{A}_{jk}(t)\right)\\
& = -\alpha_j(t)\circ A_{jk}(t) - \alpha_k(t)\circ A_{jk}(t) +\mathcal{L}_{\theta_{jk}(t)}\circ \alpha_j(t) + \mathcal{L}_{\theta_{jk}(t)}\circ \alpha_k(t) \notag
\end{align}
\end{lemma}

\begin{proof}
Let us compute the left-hand side of $(\ref{o23})$. From $(\ref{o24})$, we have
{\small{\begin{align}\label{o25}
&\left(\dot{\alpha}_k(t)- \dot{\alpha}_j(t)+\mathcal{L}_{\dot{\theta}_{jk}(t)}  - \dot{A}_{jk}(t) \right)(T_j^\alpha) \\
&= - \sum_{\beta=1}^p r_{jk}^{\alpha\beta}\frac{\partial^2 T_k^\beta}{\partial t_k^2} + \frac{\partial T_j^\alpha}{\partial t_j^2} +\left[ \dot{\theta}_{jk}(t), T_j^\alpha \right]- \sum_{\beta, \gamma=1}^p \frac{\partial}{\partial t_j}\left( \frac{\partial r_{jk}^{\alpha\beta}}{\partial t} r_{kj}^{\beta \gamma} \right) T_j^\gamma \notag \\
&=  - \sum_{\beta=1}^p r_{jk}^{\alpha\beta}\frac{\partial^2 T_k^\beta}{\partial t_k^2} + \frac{\partial T_j^\alpha}{\partial t_j^2}  +\left[ \dot{\theta}_{jk}(t), T_j^\alpha \right]  -\sum_{\beta=1}^p \frac{\partial^2 r_{jk}^{\alpha\beta}}{\partial t^2} T_k^\beta + \sum_{\beta=1}^p \left[ \theta_{jk}, \frac{\partial r_{jk}^{\alpha\beta}}{\partial t}\right] T_{k}^\beta + \sum_{\beta, \gamma=1}^p \frac{\partial r_{jk}^{\alpha\beta}}{\partial t}\left[ \theta_{jk}, r_{kj}^{\beta\gamma}\right] T_{j}^\gamma + \sum_{\beta, \gamma, \delta=1}^q  \frac{\partial r_{jk}^{\alpha\beta}}{\partial t} r_{kj}^{\beta \gamma} \frac{\partial r_{jk}^{\gamma \delta}}{\partial t} T_{k}^\delta \notag
\end{align}}}

On the other hand, let us compute the right-hand side of $(\ref{o23})$. Since $\alpha_k(t)-\alpha_j(t)+ \left(\mathcal{L}_{\theta_{jk}(t)} -A_{jk}(t)\right)=0$,

{\small{\begin{align} \label{o26}
&\left(-\alpha_j(t)\circ A_{jk}(t) - \alpha_k(t)\circ A_{jk}(t) +\mathcal{L}_{\theta_{jk}(t)}\circ \alpha_j(t) +  \mathcal{L}_{\theta_{jk}(t)}\circ \alpha_k(t)\right)\left(T_j^\alpha \right) \\
&= \left( -2 \alpha_k(t)\circ A_{jk}(t) -  \left( \mathcal{L}_{\theta_{jk}(t)} - A_{jk}(t) \right) \circ A_{jk}(t) +\mathcal{L}_{\theta_{jk}(t)}\circ \alpha_j(t) + \mathcal{L}_{\theta_{jk}(t)}\circ \alpha_k(t) \right)(T_j^\alpha) \notag \\
&= 2\sum_{\beta=1}^p \frac{\partial r_{jk}^{\alpha\beta}}{\partial t}\frac{\partial T_k^\beta}{\partial t} - \sum_{\beta, \gamma =1}^p \frac{\partial r_{jk}^{\alpha\beta}}{\partial t} r_{kj}^{\beta \gamma} \left(\mathcal{L}_{\theta_{jk}(t)}- A_{jk}(t) \right) (T_j^\gamma) -\left[ \theta_{jk}(t) ,  \frac{\partial T_j^\alpha }{\partial t}  \right] - \left[ \theta_{jk}(t) ,  \sum_{\beta=1}^p r_{jk}^{\alpha\beta} \frac{\partial T_k^\beta }{\partial t}  \right] \notag \\
&=  2\sum_{\beta=1}^p \frac{\partial r_{jk}^{\alpha\beta}}{\partial t}\frac{\partial T_k^\beta}{\partial t} - \sum_{\beta=1}^p \frac{\partial r_{jk}^{\alpha\beta}}{\partial t}\left[ \theta_{jk}(t), T_k^\beta \right] + \sum_{\beta, \gamma=1}^p \frac{\partial r_{jk}^{\alpha\beta}}{\partial t} \left[ \theta_{jk}(t), r_{kj}^{\beta \gamma} \right] T_j^\gamma + \sum_{\beta, \gamma, \delta=1}^p \frac{\partial r_{jk}^{\alpha\beta}}{\partial t} r_{kj}^{\beta \gamma} \frac{\partial r_{jk}^{\gamma \delta}}{\partial t} T_k^\delta \notag \\
& - \left[ \theta_{jk}(t) ,  \frac{\partial T_j^\alpha}{\partial t} \right] - \left[ \theta_{jk}(t) ,  \sum_{\beta=1}^p  r_{jk}^{\alpha\beta} \frac{\partial T_k^\beta }{\partial t}  \right] \notag
 \end{align}}}
Then (\ref{o27}), (\ref{o25}) and (\ref{o26}) implies (\ref{o23}).

\end{proof}

Next, by taking the derivative of $(\ref{o17})$ with respect to $t_k$ again, we have
{{\begin{align}\label{o29}
\frac{\partial^2 r_{ik}^{\alpha \gamma}}{\partial t^2} &= \sum_{\beta=1}^q \sum_{\eta, \delta=1}^n \frac{\partial^2 r_{ij}^{\alpha\beta}}{\partial \xi_j^\delta \partial \xi_j^\eta} \frac{\partial f_{jk}^\delta}{\partial t}\frac{\partial f_{jk}^\eta}{\partial t} r_{jk}^{\beta \gamma}+ 2 \sum_{\beta=1}^q \sum_{\eta, \delta=1}^n \frac{\partial}{\partial \xi_j^\eta}\left( \frac{\partial r_{ij}^{\alpha\beta}}{\partial t_j} \right)\frac{\partial f_{jk}^\eta}{\partial t} r_{jk}^{\beta \gamma} + \sum_{\beta=1}^q \sum_{\eta=1}^n \frac{\partial r_{ij}^{\alpha\beta}}{\partial \xi_j^\eta}\frac{\partial^2 f_{jk}}{\partial t^2} r_{jk}^{\beta \gamma} \\
&+ 2 \sum_{\beta=1}^q \sum_{\eta=1}^n \frac{\partial r_{ij}^{\alpha \beta}}{\partial \xi_j^\eta} \frac{\partial f_{jk}^{\eta}}{\partial t} \frac{\partial r_{jk}^{\beta \gamma}}{\partial t} 
+ \sum_{\beta=1}^q \frac{\partial^2 r_{ij}^{\alpha\beta}}{\partial t_j^2} r_{jk}^{\beta \gamma} + 2 \sum_{\beta=1}^q \frac{\partial r_{ij}^{\alpha\beta}}{\partial t_j} \frac{\partial r_{jk}^{\beta \gamma}}{\partial t} + \sum_{\beta=1}^q r_{ij}^{\alpha\beta} \frac{\partial^2 r_{jk}^{\beta \gamma}}{\partial t^2} \notag
\end{align}}}

Then we show that
\begin{lemma}\label{o52}
\begin{align}\label{o32}
\dot{A}_{ij} (t) - \dot{A}_{ik} (t) + \dot{A}_{jk}^i  (t) = \left[ \theta_{ij} (t) , A_{jk}^i(t)\right] - \left[ \theta_{jk}(t), A_{ij}(t)\right] + \left[ A_{ij} (t) ,   A_{jk}^i(t)\right]
\end{align}
for the definition of the bracket $[-,-]$ in $(\ref{o32})$, see \textnormal{Appendix \ref{apd1}}.
\end{lemma}

\begin{proof}
First we note that from $(\ref{o35})$ and $(\ref{o24})$ $\left[ \theta_{jk}, r_{kj}^{\alpha\beta}\right] + \frac{\partial r_{kj}^{\alpha\beta}}{\partial t_j} + \sum_{\gamma, \delta=1}^q r_{kj}^{\alpha\gamma} \frac{\partial r_{jk}^{\gamma \delta}}{\partial t} r_{kj}^{\delta \beta}=0$, we have
{\small{\begin{align*}
\dot{A}_{jk}(T_j^\alpha)=\sum_{\beta, \gamma=1}^p \frac{\partial}{\partial t_j}\left(\frac{\partial r_{jk}^{\alpha\beta}}{\partial t_k} r_{kj}^{\beta \gamma} \right) T_j^\gamma = \sum_{\beta=1}^p \frac{\partial^2 r_{jk}^{\alpha\beta}}{\partial t^2} T_k^\beta - \sum_{\beta=1}^p \left[ \theta_{jk}, \frac{\partial r_{jk}^{\alpha\beta}}{\partial t} \right] T_k^\beta - \sum_{\beta=1}^p \frac{\partial r_{jk}^{\alpha\beta}}{\partial t}\left[ \theta_{jk}, r_{kj}^{\beta \gamma} \right] T_j^\gamma - \sum_{\beta, \gamma, \delta =1}^p \frac{\partial r_{jk}^{\alpha\beta}}{\partial t} r_{kj}^{\beta \gamma} \frac{\partial r_{jk}^{\gamma \delta}}{\partial t} T_k^\delta
\end{align*}}}
Then we have
{\Small{\begin{align*}
&\left(\dot{A}_{ij} (t) - \dot{A}_{ik} (t) + \dot{A}_{jk}^i (t) \right)(T_i^\alpha)= \dot{A}_{ij}(t)(T_i^\alpha) - \dot{A}_{ik}(t)(T_i^\alpha) + \dot{A}_{jk}(t)(T_i^\alpha) + \sum_{\gamma=1}^n \left[  \dot{\theta}_{jk}, r_{ij}^{\alpha\gamma}  \right] T_j^\gamma\\
&=\dot{A}_{ij}(T_{i}^\alpha)- \dot{A}_{ik}(T_{i}^\alpha) + \dot{A}_{jk}(T_{i}^\alpha) + \sum_{\gamma=1}^n \left[ \sum_{\beta=1}^n \frac{\partial^2 f_{jk}^\beta}{\partial t^2} \frac{\partial}{\partial \xi_j^\beta} - \sum_{\beta, \eta=1}^n\frac{\partial f_{jk}^{\eta}}{\partial t}\frac{\partial}{\partial \xi_j^\eta}\left(\frac{\partial f_{jk}^\beta}{\partial t}   \right)\frac{\partial}{\partial \xi_j^\beta} ,    r_{ij}^{\alpha \gamma}  \right] T_{j}^\gamma\\
&= \sum_{\gamma=1}^p \frac{\partial^2 r_{ij}^{\alpha \gamma}}{\partial t^2} T_j^\gamma - \sum_{\gamma=1}^p \left[ \theta_{ij},   \frac{\partial r_{ij}^{\alpha \gamma}}{\partial t}\right] T_{j}^\gamma- \sum_{\beta, \gamma=1}^p \frac{\partial r_{ij}^{\alpha\beta}}{\partial t}\left[ \theta_{ij}, r_{ji}^{\beta\gamma}\right] T_{i}^\gamma - \sum_{\beta, \gamma, \delta=1}^p  \frac{\partial r_{ij}^{\alpha\beta}}{\partial t} r_{ji}^{\beta \gamma} \frac{\partial r_{ij}^{\gamma \delta}}{\partial t} T_{j}^\delta \\
& - \sum_{\gamma=1}^p \frac{\partial^2 r_{ik}^{\alpha \gamma}}{\partial t^2} T_{k}^\gamma + \sum_{\gamma=1}^p \left[ \theta_{ik}, \frac{\partial r_{ik}^{\alpha \gamma}}{\partial t} \right] T_{k}^\gamma + \sum_{\beta, \gamma=1}^p \frac{\partial r_{ik}^{\alpha\beta}}{\partial t}\left[ \theta_{ik}, r_{ki}^{\beta\gamma}\right] T_{i}^\gamma + \sum_{\beta, \gamma, \delta=1}^p  \frac{\partial r_{ik}^{\alpha\beta}}{\partial t} r_{ki}^{\beta \gamma} \frac{\partial r_{ik}^{\gamma \delta}}{\partial t} T_{k}^\delta\\
&+ \sum_{\beta, \gamma=1}^p r_{ij}^{\alpha\beta}\frac{\partial^2 r_{jk}^{\beta \gamma}}{\partial t^2} T_{k}^\gamma - \sum_{\beta, \gamma=1}^p r_{ij}^{\alpha\beta}\left[ \theta_{jk}, \frac{\partial r_{jk}^{\beta \gamma}}{\partial t} \right] T_{k}^\gamma - \sum_{\beta, \gamma,\delta =1}^p r_{ij}^{\alpha\beta} \frac{\partial r_{jk}^{\beta\gamma}}{\partial t}\left[ \theta_{jk}, r_{kj}^{\gamma\delta}\right] T_{j}^\delta - \sum_{\beta, \gamma, \delta, \eta =1}^p  r_{ij}^{\alpha\beta} \frac{\partial r_{jk}^{\beta\gamma}}{\partial t} r_{kj}^{\gamma \delta} \frac{\partial r_{jk}^{\delta \eta}}{\partial t} T_{k}^\eta \\
&+ \sum_{\gamma,\beta=1}^n \frac{\partial^2 f_{jk}^\beta}{\partial t^2}\frac{\partial r_{ij}^\gamma}{\partial \xi_j^\beta} T_{j}^\gamma - \sum_{\gamma,  \beta, \eta=1}^n \frac{\partial f_{jk}^\eta}{\partial t}\frac{\partial}{\partial \xi_j^\eta}\left(  \frac{\partial f_{jk}^\beta}{\partial t}  \right) \frac{\partial r_{ij}^{\alpha \gamma}}{\partial \xi_j^\beta} T_{j}^\gamma
\end{align*}}}
Then from $(\ref{o29})$, we have

{\Small{\begin{align*}
&\dot{A}_{ij}(T_i^\alpha) - \dot{A}_{ik}(T_i^\alpha) + \dot{A}_{jk}^i (T_i^\alpha)\\
&= - \sum_{\gamma=1}^p \left[ \theta_{ij},   \frac{\partial r_{ij}^{\alpha \gamma}}{\partial t}\right] T_j^\gamma- \sum_{\beta, \gamma=1}^p \frac{\partial r_{ij}^{\alpha\beta}}{\partial t}\left[ \theta_{ij}, r_{ji}^{\beta\gamma}\right] T_{i}^\gamma - \sum_{\beta, \gamma, \delta=1}^p  \frac{\partial r_{ij}^{\alpha\beta}}{\partial t} r_{ji}^{\beta \gamma} \frac{\partial r_{ij}^{\gamma \delta}}{\partial t} T_{j}^\delta \\
& + \sum_{\gamma=1}^p \left[ \theta_{ik}, \frac{\partial r_{ik}^{\alpha \gamma}}{\partial t} \right] T_{k}^\gamma + \sum_{\beta, \gamma=1}^p \frac{\partial r_{ik}^{\alpha\beta}}{\partial t}\left[ \theta_{ik}, r_{ki}^{\beta\gamma}\right] T_{i}^\gamma + \sum_{\beta, \gamma, \delta=1}^p  \frac{\partial r_{ik}^{\alpha\beta}}{\partial t} r_{ki}^{\beta \gamma} \frac{\partial r_{ik}^{\gamma \delta}}{\partial t} T_{k}^\delta\\
& - \sum_{\beta, \gamma=1}^p r_{ij}^{\alpha\beta}\left[ \theta_{jk}, \frac{\partial r_{jk}^{\beta \gamma}}{\partial t} \right] T_{k}^\gamma - \sum_{\beta, \gamma,\delta =1}^p r_{ij}^{\alpha\beta} \frac{\partial r_{jk}^{\beta\gamma}}{\partial t}\left[ \theta_{jk}, r_{kj}^{\gamma\delta}\right] T_{j}^\delta - \sum_{\beta, \gamma, \delta, \eta =1}^p  r_{ij}^{\alpha\beta} \frac{\partial r_{jk}^{\beta\gamma}}{\partial t} r_{kj}^{\gamma \delta} \frac{\partial r_{jk}^{\delta \eta}}{\partial t} T_{k}^\eta \\
&+ \sum_{\gamma,\beta=1}^n \cancel{ \frac{\partial^2 f_{jk}^\beta}{\partial t^2}\frac{\partial r_{ij}^\gamma}{\partial \xi_j^\beta} T_{j}^\gamma }  - \sum_{\gamma,  \beta, \eta=1}^n \frac{\partial f_{jk}^\eta}{\partial t}\frac{\partial}{\partial \xi_j^\eta}\left(  \frac{\partial f_{jk}^\beta}{\partial t}  \right) \frac{\partial r_{ij}^{\alpha \gamma}}{\partial \xi_j^\beta} T_{j}^\gamma -\sum_{\beta=1}^p \sum_{\eta, \delta=1}^n \frac{\partial^2 r_{ij}^{\alpha\beta}}{\partial \xi_j^\delta \partial \xi_j^\eta} \frac{\partial f_{jk}^\delta}{\partial t}\frac{\partial f_{jk}^\eta}{\partial t} w_{0j}^\beta - 2 \sum_{\beta=1}^p \left[ \theta_{jk}, \frac{\partial r_{ij}^{\alpha\beta}}{\partial t}\right] T_{j}^\beta \\
&  - \sum_{\beta=1}^p \sum_{\eta=1}^n \cancel{ \frac{\partial^2 f_{jk}^\eta}{\partial t^2}\frac{\partial r_{ij}^{\alpha\beta}}{\partial \xi_j^\eta} T_{j}^\beta } - 2 \sum_{\beta, \gamma=1}^p \left[ \theta_{jk}, r_{ij}^{\alpha\beta}\right] \frac{\partial r_{jk}^{\beta \gamma}}{\partial t} T_{k}^\gamma- 2 \sum_{\beta, \gamma=1}^p \frac{\partial r_{ij}^{\alpha\beta}}{\partial t} \frac{\partial r_{jk}^{\beta \gamma}}{\partial t} T_{k}^\gamma
\end{align*}}}
We continue the computation. Then
{\Small{\begin{align*}
&= - \sum_{\gamma=1}^p \left[ \theta_{ij},   \frac{\partial r_{ij}^{\alpha \gamma}}{\partial t}\right] T_{j}^\gamma- \sum_{\beta, \gamma=1}^p \frac{\partial r_{ij}^{\alpha\beta}}{\partial t}\left[ \theta_{ij}, r_{ji}^{\beta\gamma}\right] T_{i}^\gamma - \sum_{\beta, \gamma, \delta=1}^p  \frac{\partial r_{ij}^{\alpha\beta}}{\partial t} r_{ji}^{\beta \gamma} \frac{\partial r_{ij}^{\gamma \delta}}{\partial t} T_{j}^\delta \\
& + \sum_{\gamma=1}^p \left[ \theta_{ik}, \frac{\partial r_{ik}^{\alpha \gamma}}{\partial t} \right] T_{k}^\gamma + \sum_{\beta, \gamma=1}^p \frac{\partial r_{ik}^{\alpha\beta}}{\partial t}\left[ \theta_{ik}, r_{ki}^{\beta\gamma}\right] T_{i}^\gamma + \sum_{\beta, \gamma, \delta=1}^p  \frac{\partial r_{ik}^{\alpha\beta}}{\partial t} r_{ki}^{\beta \gamma} \frac{\partial r_{ik}^{\gamma \delta}}{\partial t} T_k^\delta\\
& - \sum_{\beta, \gamma=1}^p r_{ij}^{\alpha\beta}\left[ \theta_{jk}, \frac{\partial r_{jk}^{\beta \gamma}}{\partial t} \right] T_{k}^\gamma - \sum_{\beta, \gamma,\delta =1}^p r_{ij}^{\alpha\beta} \frac{\partial r_{jk}^{\beta\gamma}}{\partial t}\left[ \theta_{jk}, r_{kj}^{\gamma\delta}\right] T_{j}^\delta - \sum_{\beta, \gamma, \delta, \eta =1}^p  r_{ij}^{\alpha\beta} \frac{\partial r_{jk}^{\beta\gamma}}{\partial t} r_{kj}^{\gamma \delta} \frac{\partial r_{jk}^{\delta \eta}}{\partial t} T_k^\eta \\
&  - \sum_{\gamma=1}^p \left[ \theta_{jk}, \left[ \theta_{jk}, r_{ij}^{\alpha \gamma} \right] \right]  T_{j}^\gamma - 2 \sum_{\beta=1}^p \left[ \theta_{jk}, \frac{\partial r_{ij}^{\alpha\beta}}{\partial t}\right] T_{j}^\beta  - 2 \sum_{\beta, \gamma=1}^p \left[ \theta_{jk}, r_{ij}^{\alpha\beta}\right] \frac{\partial r_{jk}^{\beta \gamma}}{\partial t} T_{k}^\gamma- 2 \sum_{\beta, \gamma=1}^q \frac{\partial r_{ij}^{\alpha\beta}}{\partial t} \frac{\partial r_{jk}^{\beta \gamma}}{\partial t} T_{k}^\gamma 
\end{align*}}}

We recall from $(\ref{o17})$ that $\frac{\partial r_{ik}^{\alpha \gamma}}{\partial t}= \sum_{\beta=1}^p \left[ \theta_{jk}, r_{ij}^{\alpha\beta}\right] r_{jk}^{\beta \gamma} + \sum_{\beta=1}^p \frac{\partial r_{ij}^{\alpha\beta}}{\partial t} r_{jk}^{\beta \gamma} + \sum_{\beta=1}^p r_{ij}^{\alpha\beta} \frac{\partial r_{jk}^{\beta \gamma}}{\partial t}$. Then
{\Small{\begin{align}\label{o30}
&\dot{A}_{ij}(T_i^\alpha) - \dot{A}_{ik}(T_i^\alpha) + \dot{A}_{jk}^i (T_i^\alpha)\\
 &=- \sum_{\gamma=1}^p \left[ \theta_{ij},   \frac{\partial r_{ij}^{\alpha \gamma}}{\partial t}\right] T_{j}^\gamma- \sum_{\beta, \gamma=1}^p \frac{\partial r_{ij}^{\alpha\beta}}{\partial t}\left[ \theta_{ij}, r_{ji}^{\beta\gamma}\right] T_{i}^\gamma - \sum_{\beta, \gamma, \delta=1}^p  \frac{\partial r_{ij}^{\alpha\beta}}{\partial t} r_{ji}^{\beta \gamma} \left(\sum_{\eta=1}^p \frac{\partial r_{ik}^{\gamma \eta}}{\partial t} r_{kj}^{\eta\delta} - \left[ \theta_{jk}, r_{ij}^{\gamma \delta} \right] - \sum_{\eta, \xi=1}^p r_{ij}^{\gamma \eta}\frac{\partial r_{jk}^{\eta \xi}}{\partial t} r_{kj}^{\xi \delta}\right)T_{j}^\delta \notag \\
 &+\sum_{\beta, \gamma=1}^p \left[ \theta_{ij}, \frac{\partial r_{ik}^{\alpha\beta}}{\partial t} r_{ki}^{\beta \gamma} \right] T_{i}^\gamma +\sum_{\beta, \gamma=1}^p \left[ \theta_{jk}, \frac{\partial r_{ik}^{\alpha\beta}}{\partial t} r_{ki}^{\beta \gamma} \right] T_{i}^\gamma + \sum_{\beta, \gamma, \delta=1}^p \frac{\partial r_{ik}^{\alpha\beta}}{\partial t} r_{ki}^{\beta \gamma}\left( \sum_{\eta=1}^p \left[ \theta_{jk}, r_{ij}^{\gamma \eta} \right] r_{jk}^{\eta \delta} +\sum_{\eta=1}^p \frac{\partial r_{ij}^{\gamma \eta}}{\partial t} r_{jk}^{\eta \delta} + \sum_{\eta=1}^p r_{ij}^{\gamma\eta} \frac{\partial r_{jk}^{\eta \delta}}{\partial t} \right) T_{k}^\delta \notag \\
&  - \sum_{\beta, \gamma=1}^p r_{ij}^{\alpha\beta}\left[ \theta_{jk}, \frac{\partial r_{jk}^{\beta \gamma}}{\partial t} \right] T_{k}^\gamma - \sum_{\beta, \gamma,\delta =1}^p r_{ij}^{\alpha\beta} \frac{\partial r_{jk}^{\beta\gamma}}{\partial t}\left[ \theta_{jk}, r_{kj}^{\gamma\delta}\right] T_{j}^\delta - \sum_{\gamma, \delta, \eta=1}^p \left( \frac{\partial r_{ik}^{\alpha \gamma}}{\partial t} - \sum_{\beta=1}^p \left[ \theta_{jk}, r_{ij}^{\alpha\beta}\right] r_{jk}^{\beta \gamma} -  \sum_{\beta=1}^p \frac{\partial r_{ij}^{\alpha\beta}}{\partial t} r_{jk}^{\beta \gamma} \right) r_{kj}^{\gamma \delta} \frac{\partial r_{jk}^{\delta \eta}}{\partial t} T_{k}^\eta \notag \\
&  - \sum_{\gamma=1}^p \left[ \theta_{jk}, \left[ \theta_{jk}, r_{ij}^{\alpha \gamma} \right] \right]  T_{j}^\gamma - 2 \sum_{\beta=1}^p \left[ \theta_{jk}, \frac{\partial r_{ij}^{\alpha\beta}}{\partial t}\right] T_{j}^\beta  - 2 \sum_{\beta, \gamma=1}^p \left[ \theta_{jk}, r_{ij}^{\alpha\beta}\right] \frac{\partial r_{jk}^{\beta \gamma}}{\partial t} T_{k}^\gamma - 2 \sum_{\beta, \gamma=1}^q \frac{\partial r_{ij}^{\alpha\beta}}{\partial t} \frac{\partial r_{jk}^{\beta \gamma}}{\partial t} T_{k}^\gamma  \notag
\end{align}}}
We continue the computation. Then
{\Small{\begin{align}
 &=- \sum_{\gamma=1}^p \left[ \theta_{ij},   \frac{\partial r_{ij}^{\alpha \gamma}}{\partial t}\right] T_{j}^\gamma- \sum_{\beta, \gamma=1}^p \frac{\partial r_{ij}^{\alpha\beta}}{\partial t}\left[ \theta_{ij}, r_{ji}^{\beta\gamma}\right] T_{i}^\gamma - \sum_{\beta, \gamma,  \eta =1}^p  \frac{\partial r_{ij}^{\alpha\beta}}{\partial t} r_{ji}^{\beta \gamma}  \frac{\partial r_{ik}^{\gamma \eta}}{\partial t} T_k^\eta -  \sum_{\alpha, \beta=1}^p \frac{\partial r_{ij}^{\alpha\beta}}{\partial t}\left[ \theta_{jk}, r_{ji}^{\beta  \gamma} \right] T_{i}^\gamma \notag \\
  &+\sum_{\beta, \gamma=1}^p \left[ \theta_{ij}, \frac{\partial r_{ik}^{\alpha\beta}}{\partial t} r_{ki}^{\beta \gamma} \right] T_{i}^\gamma +\sum_{\beta, \gamma, \eta=1}^p \left[ \theta_{jk}, \frac{\partial r_{ik}^{\alpha\beta}}{\partial t} r_{ki}^{\beta \gamma}  r_{ij}^{\gamma \eta} \right] T_{j}^\eta + \sum_{\beta, \gamma, \delta, \eta =1}^p \frac{\partial r_{ik}^{\alpha\beta}}{\partial t} r_{ki}^{\beta \gamma} \frac{\partial r_{ij}^{\gamma \eta}}{\partial t}  T_{j}^\eta \notag \\
  &- \sum_{\beta, \gamma, \delta=1}^p \left[ \theta_{jk}, r_{ij}^{\alpha\beta} \frac{\partial r_{jk}^{\beta \gamma}}{\partial t} r_{kj}^{\gamma \delta} \right] T_{j}^\delta    - \sum_{\gamma=1}^p \left[ \theta_{jk}, \left[ \theta_{jk}, r_{ij}^{\alpha \gamma} \right] \right]  T_{j}^\gamma  - 2 \sum_{\beta, \gamma=1}^p \left[ \theta_{jk}, r_{ij}^{\alpha\beta}\right] \frac{\partial r_{jk}^{\beta \gamma}}{\partial t} T_{k}^\gamma  - 2 \sum_{\beta=1}^p \left[ \theta_{jk}, \frac{\partial r_{ij}^{\alpha\beta}}{\partial t}\right] T_{j}^\beta \notag \\
   & = - \sum_{\gamma=1}^p \left[ \theta_{ij},   \frac{\partial r_{ij}^{\alpha \gamma}}{\partial t}\right] T_{j}^\gamma- \sum_{\beta, \gamma=1}^p \frac{\partial r_{ij}^{\alpha\beta}}{\partial t}\left[ \theta_{ij}, r_{ji}^{\beta\gamma}\right] T_{i}^\gamma - \sum_{\beta, \gamma, \eta =1}^p  \frac{\partial r_{ij}^{\alpha\beta}}{\partial t} r_{ji}^{\beta \gamma}  \frac{\partial r_{ik}^{\gamma \eta}}{\partial t} T_k^\eta -  \sum_{\alpha, \beta=1}^p \frac{\partial r_{ij}^{\alpha\beta}}{\partial t}\left[ \theta_{jk}, r_{ji}^{\beta  \gamma} \right] T_{i}^\gamma \notag \\
  &+\sum_{\beta, \gamma=1}^p \left[ \theta_{ij}, \frac{\partial r_{ik}^{\alpha\beta}}{\partial t} r_{ki}^{\beta \gamma} \right] T_{i}^\gamma  -  \sum_{\beta =1}^p \left[ \theta_{jk}, \frac{\partial r_{ij}^{\alpha\beta}}{\partial t} \right] T_{j}^\beta + \sum_{\beta, \gamma, \delta, \eta =1}^p \frac{\partial r_{ik}^{\alpha\beta}}{\partial t} r_{ki}^{\beta \gamma} \frac{\partial r_{ij}^{\gamma \eta}}{\partial t}  T_{j}^\eta   \notag
\end{align}}}

On the other hand, we recall from $(\ref{o17})$ that $\frac{\partial r_{ik}^{\alpha \gamma}}{\partial t}= \sum_{\beta=1}^p \left[ \theta_{jk}, r_{ij}^{\alpha\beta}\right] r_{jk}^{\beta \gamma} + \sum_{\beta=1}^p \frac{\partial r_{ij}^{\alpha\beta}}{\partial t} r_{jk}^{\beta \gamma} + \sum_{\beta=1}^p r_{ij}^{\alpha\beta} \frac{\partial r_{jk}^{\beta \gamma}}{\partial t}$,
\begin{align*}
A_{jk}^i\left(T_i^\alpha \right) = A_{jk}\left(T_i^\alpha\right) + \sum_{\gamma=1}^p \left[ \theta_{jk}, r_{ij}^{\alpha \gamma} \right] T_j^\gamma= \sum_{\xi=1}^p \left(\sum_{\beta=1}^p \frac{\partial r_{ik}^{\alpha \beta}}{\partial t} r_{ki}^{\beta \xi} - \sum_{\beta=1}^p \frac{\partial r_{ij}^{\alpha\beta}}{\partial t} r_{ji}^{\beta \xi} \right) T_{i}^\xi
\end{align*}
Then we have
{\Small{\begin{align}\label{o31}
 &\left[\theta_{ij}, A_{jk}^i \right]\left(T_{i}^\alpha\right) - \left[\theta_{jk}, A_{ij}\right]\left(T_{i}^\alpha\right) + \left[A_{ij}, A_{jk}^i\right]  \left( T_{i}^\alpha\right)\\
&=\sum_{\xi=1}^p \left[ \theta_{ij}, \sum_{\beta=1}^p \frac{\partial r_{ik}^{\alpha \beta}}{\partial t} r_{ki}^{\beta \xi} - \sum_{\beta=1}^p \frac{\partial r_{ij}^{\alpha\beta}}{\partial t} r_{ji}^{\beta \xi}\right] T_{i}^\xi - \sum_{\xi=1}^p \left[ \theta_{jk},  \sum_{\beta=1}^p \frac{\partial r_{ij}^{\alpha\beta}}{\partial t} r_{ji}^{\beta \xi}\right] T_{i}^\xi \notag \\
&+  \sum_{\xi,\eta=1}^r \left(\sum_{\beta=1}^p \frac{\partial r_{ik}^{\alpha \beta}}{\partial t} r_{ki}^{\beta \xi} - \sum_{\beta=1}^p \frac{\partial r_{ij}^{\alpha\beta}}{\partial t} r_{ji}^{\beta \xi} \right) \frac{\partial r_{ij}^{\xi \eta}}{\partial t} T_{j}^\eta - \sum_{\beta,\gamma=1}^p \frac{\partial r_{ij}^{\alpha\beta}}{\partial t} r_{ji}^{\beta \gamma} \left( \sum_{\eta=1}^p \frac{\partial r_{ik}^{\gamma \eta}}{\partial t} T_{k}^\eta - \sum_{\eta=1}^p \frac{\partial r_{ij}^{\gamma \eta}}{\partial t} T_{j}^\eta  \right) \notag \\
&=\sum_{\beta=1}^p \left[ \theta_{ij}, \frac{\partial r_{ik}^{\alpha \beta}}{\partial t}\right] T_{k}^\beta + \sum_{\beta, \xi=1}^p \frac{\partial r_{ik}^{\alpha \beta}}{\partial t}\left[ \theta_{ij}, r_{ki}^{\beta \xi}\right] T_{i}^\xi - \sum_{\beta=1}^p \left[ \theta_{ij}, \frac{\partial r_{ij}^{\alpha \beta}}{\partial t} \right] T_{j}^\beta - \sum_{\beta, \xi=1}^p \frac{\partial r_{ij}^{\alpha \beta}}{\partial t}\left[ \theta_{ij}, r_{ji}^{\beta \xi}\right] T_{i}^\xi \notag \\
&- \sum_{\beta=1}^p \left[ \theta_{jk}, \frac{\partial h_{ij}^{\alpha\beta}}{\partial t} \right] T_{j}^\beta - \sum_{\beta, \xi=1}^p \frac{\partial r_{ij}^{\alpha\beta}}{\partial t}\left[ \theta_{jk}, r_{ji}^{\beta \xi}\right] T_{i}^\xi +  \sum_{\xi,\eta,\beta=1}^p \frac{\partial r_{ik}^{\alpha \beta}}{\partial t} r_{ki}^{\beta \xi} \frac{\partial r_{ij}^{\xi \eta}}{\partial t} T_{j}^\eta - \sum_{\beta,\gamma,\eta=1}^p \frac{\partial r_{ij}^{\alpha\beta}}{\partial t} r_{ji}^{\beta \gamma}  \frac{\partial r_{ik}^{\gamma \eta}}{\partial t} T_{k}^\eta \notag
\end{align}}}
$(\ref{o30})$ and $(\ref{o31})$ implies $(\ref{o32})$.

\end{proof}

From $(\ref{o36}), (\ref{o37}), (\ref{o23})$ and $(\ref{o32})$, and substituting $t=0$, we obtain
{\small{\begin{align}
\dot{\theta}_{ij}(0) -\dot{\theta}_{ik}(0) + \dot{\theta}_{jk} (0) &= \left[ \theta_{ij}(0), \theta_{jk}(0) \right] \label{o38}\\
\left[ T_{0j}^\alpha, \dot{\alpha}_j^\beta (0) \right] - \left[ T_{0j}^\beta, \dot{\alpha}_j^\alpha(0)\right] - \sum_{\gamma=1}^p g_{0j\alpha\beta}^\gamma \dot{\alpha}_j^\gamma (0) & = 2\left[ \alpha_j^\alpha(0), \alpha_j^\beta(0) \right]  + 2\alpha_j(0)\left( \sum_{\gamma=1}^p \frac{\partial g_{j\alpha\beta}^\gamma}{\partial t}|_{t=0} \cdot T_{0j}^\gamma \right) - \sum_{\gamma=1}^p \frac{\partial^2 g_{j\alpha\beta}^\gamma}{\partial t^2}|_{t=0} \cdot T_{0j}^\gamma\\
\dot{\alpha}_k(0)- \dot{\alpha}_j(0) + D_0'\left(\left(\dot{\theta}_{jk}(0), \dot{A}_{jk}(0) \right)  \right) & = -\alpha_j(0)\circ A_{jk}(0) - \alpha_k(0)\circ A_{jk}(0) +\mathcal{L}_{\theta_{jk}(0)}\circ \alpha_j(0) + \mathcal{L}_{\theta_{jk}(0)}\circ \alpha_k(0)\\
\dot{A}_{ij} (0) - \dot{A}_{ik} (0) + \dot{A}_{jk}^i  (0) &= \left[ \theta_{ij} (0) , A_{jk}^i(0)\right] - \left[ \theta_{jk}(0), A_{ij}(t)\right] + \left[ A_{ij} (0) ,   A_{jk}^i(0)\right] \label{o39}
\end{align}}}

Let us define $\left[\alpha_j(0), \alpha_j(0) \right]_{\Theta_{\mathcal{F}_0}} \in \Gamma\left( U_j^0, \mathscr{H}om_{\mathcal{O}_M}\left( \bigwedge^2 \Theta_{\mathcal{F}_0}, \frac{\Theta_M}{\Theta_{\mathcal{F}_0}} \right) \right)$ as follows:
\begin{align*}
\left[ \alpha_j(0), \alpha_j(0)\right]_{\Theta_{\mathcal{F}_0}}: \Gamma\left(U_j^0, \bigwedge^2 \Theta_{\mathcal{F}_0} \right) &\to \Gamma\left( U_j^0, \frac{\Theta_M}{\Theta_{\mathcal{F}_0}} \right) \\
T_{0j}^\alpha \wedge T_{0j}^\beta &\mapsto  \overline{ 2\left[ \alpha_j^\alpha(0), \alpha_j^\beta(0) \right]  + 2\alpha_j(0)\left( \sum_{\gamma=1}^p \frac{\partial g_{j\alpha\beta}^\gamma}{\partial t}|_{t=0} \cdot T_{0j}^\gamma \right)}
\end{align*}
and linearly extends to $\Gamma\left( U_j, \bigwedge^2 \Theta_{\mathcal{F}_0} \right)$. Then by combining $(\ref{o38})$ and $(\ref{o39})$ (see Appendix \ref{apd1}), we can write $(\ref{o38})- (\ref{o39})$ as follows:
{\small{\begin{align}
&D_1'\left(\dot{\alpha}_j(0) \right)= \left[\alpha_j(0), \alpha_j(0) \right]_{\Theta_{\mathcal{F}_0}} \label{o40}  \\
\dot{\alpha}_k(0)- \dot{\alpha}_j(0) + D_0'\left(\left(\dot{\theta}_{jk}(0), \dot{A}_{jk}(0) \right)  \right)& = -\alpha_j(0)\circ A_{jk}(0) - \alpha_k(0)\circ A_{jk}(0) +\mathcal{L}_{\theta_{jk}(0)}\circ \alpha_j(0) + \mathcal{L}_{\theta_{jk}(0)}\circ \alpha_k(0)\notag \\
\left(\dot{\theta}_{ij}(0), \dot{A}_{ij}(0) \right) -& \left( \dot{\theta}_{ik} (0) , \dot{A}_{ik} (0) \right)+ \left( \dot{\theta}_{jk}(0), \dot{A}_{jk}(0) \right)  = \left[ \left( \theta_{ij} (0), A_{ij} (0)   \right) , \left( \theta_{jk} (0) , A_{jk}  (0) \right) \right] \notag
\end{align}}}

Since $(M, \Theta_{\mathcal{F}_0})=(M_0, \Theta_{\mathcal{F}_0})$ by assumption, identifying $U_j$ with $U_j^0=U_j\times 0$, we may consider $\mathcal{U}=\{U_j\}$ as a finite open covering of $M$. For a given 

\begin{lemma}\label{o54}
For any $1$-cocycle $\left( \left\{ \alpha_j \right\} , \left\{ (\theta_{jk}, A_{jk} \right\} \right)\in C^0\left( \mathcal{U}, \mathscr{H}om_{\mathcal{O}_M}\left( \Theta_{\mathcal{F}_0}, \Theta_M \right) \right) \bigoplus C^1\left( \mathcal{U}, \mathcal{E}_{\Theta_{\mathcal{F}_0}} \right)$ where $\alpha_j:\Gamma\left( U_j,\Theta_{\mathcal{F}_0} \right)\to \Gamma\left( U_j, \Theta_M \right)$, $\theta_{jk}\in \Gamma\left(U_{jk}, \Theta_M \right)$ and $A_{jk}:\Gamma\left(U_{jk}, \Theta_{\mathcal{F}_0} \right)\to \Gamma\left( U_{jk}, \Theta_{\mathcal{F}_0} \right)$ satisfying
\begin{enumerate}
\item $D_1'(\alpha_j)=0$, i.e. $\left[ T_{0j}^\alpha, \alpha_j\left(T_{0j}^\beta\right) \right] - \left[ T_{0j}^\beta, \alpha_j\left(T_{0j}^\alpha \right) \right] - \sum_{\gamma=1}^p g_{0j\alpha\beta}^\gamma \alpha_j\left(T_{0j}^\gamma\right)= Q_{j\alpha\beta}\in \Gamma\left(U_j, \Theta_{\mathcal{F}_0}\right)$ for $\alpha, \beta=1,...,p$.
\item $\alpha_k-\alpha_j+D_0'\left( \left(\theta_{jk}, A_{jk} \right) \right)=0$, i.e. $\alpha_k(T_j^\alpha)- \alpha_j(T_j^\alpha)+ \mathcal{L}_{\theta_{jk}}(T_j^\alpha) - A_{jk}(T_j^\alpha)=0$ for $\alpha=1,...,p$.
\item $\left( \theta_{ij}, A_{ij} \right)  - \left(\theta_{ik}, A_{ik} \right) + \left(\theta_{jk}, A_{jk} \right)=0$, i.e. $\theta_{ij}-\theta_{ik}+ \theta_{jk}=0$ and $A_{ij}- A_{ik}+ A_{jk}^i=0$,
\end{enumerate}
we define $R_j:=\left[\alpha_j, \alpha_j\right]_{\Theta_{\mathcal{F}_0}}\in \Gamma\left(U_j, \mathscr{H}om_{\mathcal{O}_M}\left(\bigwedge^2 \Theta_{\mathcal{F}_0}, \frac{\Theta_M}{\Theta_{\mathcal{F}_0}} \right) \right)$ by
\begin{align*}
R_j:\Gamma\left( U_j, \bigwedge^2 \Theta_{\mathcal{F}_0} \right) &\to \Gamma\left( U_j, \frac{\Theta_M}{\Theta_{\mathcal{F}_0}} \right) \\
 T_{0j}^\alpha \wedge T_{0j}^\beta &\mapsto \overline{2\left[\alpha_j\left(T_{0j}^\alpha \right), \alpha_j\left(T_{0j}^\beta \right) \right] - 2\alpha_j\left(Q_{j\alpha\beta} \right)}
\end{align*}
and linearly extends to $\Gamma\left( U_j, \bigwedge^2 \Theta_{\mathcal{F}_0} \right)$, and we define $F_{jk}\in \Gamma\left( U_{jk}, \mathscr{H}om_{\mathcal{O}_M}\left( \Theta_{\mathcal{F}_0}, \Theta_M \right) \right)$ by
\begin{align*}
F_{jk}: \Gamma\left( U_{jk}, \Theta_{\mathcal{F}_0} \right) & \to \Gamma\left( U_j, \Theta_M \right) \\
                  T_{0j}^\alpha &\mapsto -\alpha_j \left( A_{jk}(T_{0j}^\alpha) \right) - \alpha_k\left(A_{jk} (T_{0j}^\alpha) \right) +\left[ \theta_{jk}, \alpha_j \left(T_{0j}^\alpha\right) \right] + \left[ \theta_{jk}, \alpha_k \left( T_{0j}^\alpha\right) \right]
\end{align*}
and linearly extends to $\Gamma\left( U_{jk}, \Theta_{\mathcal{F}_0} \right)$. Then
{\small{\begin{align}\label{it1}
\left( \left\{ R_j\right\}, \left\{F_{jk} \right\}, \left\{ \left[  \left(\theta_{ij}, A_{ij}  \right) , \left(\theta_{jk}, A_{jk} \right)   \right]  \right\}\right) \in C^0\left( \mathcal{U}, \mathscr{H}om_{\mathcal{O}_M}\left( \bigwedge^2 \Theta_{\mathcal{F}_0}, \frac{\Theta_M}{\Theta_{\mathcal{F}_0}} \right)\right) \bigoplus C^1\left( \mathcal{U}, \mathscr{H}om_{\mathcal{O}_M}\left( \Theta_{\mathcal{F}_0}, \Theta_M \right) \right) \bigoplus C^2\left( \mathcal{U} , \mathcal{E}_{\Theta_{\mathcal{F}_0}} \right)
\end{align}}}
defines a $2$-cocycle in \v Cech resolution of $(\ref{o7})$ of $\mathcal{E}_{\Theta_{\mathcal{F}_0}}^\bullet$.
\end{lemma}

\begin{remark}
We show that $F_{jk}\in C^1\left( \mathcal{U}, \mathscr{H}om_{\mathcal{O}_M}\left( \Theta_{\mathcal{F}_0}, \Theta_M \right) \right)$. In fact, first we note that $(\theta_{jk}, A_{jk})$ on $U_j$ is transformed on $U_k$ into $\left( \theta_{jk}, A_{jk}^k \right)$ where 
\begin{align*}
A_{jk}^k\left(T_{0j}^\alpha \right)= A_{jk}\left(T_{0j}^\alpha \right) - \sum_{\beta=1}^p \left[ \theta_{jk}, r_{0jk}^{\alpha\beta} \right] T_{0k }^\beta = -A_{kj}\left(T_{0j}^\alpha \right) =- \sum_{\beta=1}^p r_{jk}^{\alpha\beta} A_{kj}\left(T_{0k}^\beta \right) 
\end{align*}
Then
{\Small{\begin{align*}
&F_{kj}\left(T_{0j}^\alpha \right) = \sum_{\beta=1}^p r_{jk}^{\alpha\beta} F_{kj}\left(T_{0k}^\beta \right) = -  \sum_{\beta=1}^p r_{0jk}^{\alpha\beta}\alpha_k\left(A_{kj}\left(T_{0k}^\beta \right) \right) - \sum_{\beta=1}^p r_{0jk}^{\alpha\beta} \alpha_j\left(A_{kj}\left(T_{0k}^\beta \right) \right) +\sum_{\beta=1}^p r_{0jk}^{\alpha\beta}\left[ \theta_{kj}, \alpha_k\left(T_{0k}^\beta \right)\right] + \sum_{\beta=1}^p r_{0jk}^{\alpha\beta}\left[ \theta_{kj}, \alpha_j\left(T_{0k}^\beta \right)  \right]\\ 
&= \alpha_k\left( A_{jk}\left( T_{0j}^\alpha \right)        \right) - \sum_{\beta=1}^p \left[ \theta_{jk},  r_{0jk}^{\alpha\beta} \right] \alpha_k\left( T_{0k}^\beta \right) + \alpha_j\left( A_{jk}\left(T_{0j}^\alpha \right) \right) - \sum_{\beta, \gamma=1}^p\left[ \theta_{jk}, r_{0jk}^{\alpha\beta} \right] r_{0kj}^{\beta \gamma} \alpha_j\left( T_{0j}^\gamma \right)   - \left[ \theta_{jk}, \alpha_k\left( T_{0j}^\alpha  \right) \right] + \sum_{\beta=1}^p \left[ \theta_{jk}, r_{0jk}^{\alpha\beta} \right] \alpha_k\left(T_{0k}^\beta \right) \\
&- \left[ \theta_{jk}, \alpha_j\left( T_{0j}^\alpha \right) \right] + \sum_{\beta=1}^p\left[ \theta_{jk}, r_{0jk}^{\alpha\beta} \right] \alpha_j\left(T_{0k}^\beta \right) =\alpha_k\left(A_{jk} \left( T_{0j}^\alpha \right) \right) + \alpha_j\left( A_{jk}\left( T_{0j}^\alpha \right) \right) - \left[  \theta_{jk}, \alpha_k\left( T_{0j}^\alpha \right)\right] - \left[ \theta_{jk}, \alpha_j \left( T_{0j}^\alpha \right)\right] = - F_{jk}\left( T_{0j}^\alpha\right)
\end{align*}}}
\end{remark}

\begin{remark}
We will deduce the equality $(\ref{o59})$ below which will be used in the proof of \textnormal{Lemma $\ref{o54}$}. We note that
{\small{\begin{align}\label{o56}
&\left[T_{0j}^\alpha, \left[ T_{0j}^\beta, T_{0j}^\gamma \right] \right] + \left[T_{0j}^\beta, \left[T_{0j}^\gamma, T_{0j}^\alpha \right] \right] + \left[T_{0j}^\gamma, \left[ T_{0j}^\alpha, T_j^\beta\right] \right] =\sum_{\xi=1}^p \left[T_{0j}^\alpha, g_{0j\beta\gamma}^\xi T_{0j}^\xi \right] + \sum_{\xi=1}^p \left[ T_{0j}^\beta, g_{0j\gamma \alpha}^\xi T_{0j}^\xi \right] + \sum_{\xi=1}^p \left[ T_{0j}^\gamma, g_{0j\alpha\beta}^\xi T_{0j}^\xi\right]\\
&= \sum_{\eta=1}^p \left[T_{0j}^\alpha,  g_{0j\beta\gamma}^\eta \right] T_{0j}^\eta + \sum_{\xi, \eta=1}^p g_{0\beta\gamma}^\xi g_{0j\alpha \xi}^\eta T_{0j}^\eta + \sum_{\eta=1}^p \left[ T_{0j}^\beta, g_{0j\gamma\alpha}^\eta\right] T_j^\eta + \sum_{\xi, \eta=1}^p g_{0j\gamma\alpha}^\xi g_{0j\beta \xi}^\eta T_{0j}^\eta + \sum_{\eta=1}^p \left[ T_{0j}^\gamma, g_{0j\alpha\beta}^\eta \right] T_j^\eta + \sum_{\xi, \eta=1}^p g_{j\alpha\beta}^\xi g_{0j\gamma \xi}^\eta T_{0j}^\eta = 0 \notag
\end{align}}} 

We compute 
{\small{\begin{align}\label{o57}
&\left[ T_{0j}^\gamma,  Q_{j\alpha\beta}\right] + \left[ T_{0j}^\alpha, Q_{j\beta\gamma} \right] + \left[ T_{0j}^\beta, Q_{j\gamma \alpha} \right]\\
&=  \left[ T_{0j}^\gamma, \left[ T_{0j}^\alpha, \alpha_j \left(T_{0j}^\beta \right) \right] \right] - \left[ T_{0j}^\gamma, \left[ T_{0j}^\beta, \alpha_j\left(T_{0j}^\alpha\right) \right] \right] - \sum_{\xi=1}^p \left[ T_{0j}^\gamma, g_{j\alpha\beta}^\xi \right] \alpha_j\left(T_{0j}^\xi \right) - \sum_{\xi=1}^p g_{0j\alpha\beta}^\xi \left[ T_{0j}^\gamma, \alpha_j\left(T_{0j}^\xi \right)\right] \notag \\
& +  \left[ T_{0j}^\alpha, \left[ T_{0j}^\beta, \alpha_j\left(T_{0j}^\gamma\right) \right] \right] - \left[ T_{0j}^\alpha, \left[ T_{0j}^\gamma, \alpha_j\left(T_{0j}^\beta\right) \right] \right] - \sum_{\xi=1}^p \left[ T_j^\alpha, g_{0j\beta\gamma}^\xi \right] \tilde{\alpha}_j(T_j^\xi ) - \sum_{\xi=1}^p g_{0j\beta\gamma}^\xi \left[ T_{0j}^\alpha, \alpha_j\left(T_{0j}^\xi \right)\right] \notag \\
& +  \left[ T_{0j}^\beta, \left[ T_{0j}^\gamma, \alpha_j\left(T_{0j}^\alpha\right) \right] \right] - \left[ T_{0j}^\beta, \left[ T_{0j}^\alpha, \alpha_j\left(T_{0j}^\gamma \right) \right] \right] - \sum_{\xi=1}^p \left[ T_{0j}^\beta, g_{0j\gamma\alpha}^\xi \right] \alpha_j \left(T_{0j}^\xi \right) - \sum_{\xi=1}^p g_{0j\gamma\alpha}^\xi \left[ T_{0j}^\beta, \alpha_j(T_{0j}^\xi)\right] \notag\\
&= \sum_{\xi=1}^p \left[g_{0j\gamma\alpha}^\xi T_{0j}^\xi, \alpha_j \left(T_{0j}^\beta \right) \right]- \sum_{\xi=1}^p \left[ T_{0j}^\gamma, g_{0j\alpha\beta}^\xi \right] \alpha_j \left(T_{0j}^\xi \right) - \sum_{\xi=1}^p g_{0j\alpha\beta}^\xi \left[ T_{0j}^\gamma, \alpha_j \left(T_{0j}^\xi \right)\right] \notag\\
&+\sum_{\xi=1}^p \left[  g_{0j\alpha\beta}^\xi T_{0j}^\xi, \alpha_j \left(T_{0j}^\gamma \right)\right] - \sum_{\xi=1}^p \left[ T_{0j}^\alpha, g_{0j\beta\gamma}^\xi \right] \alpha_j \left(T_j^\xi \right) - \sum_{\xi=1}^p g_{0j\beta\gamma}^\xi \left[ T_{0j}^\alpha, \alpha_j \left(T_{0j}^\xi \right)\right]\notag\\
&+ \sum_{\xi=1}^p \left[ g_{0j\beta\gamma}^\xi T_{0j}^\xi, \alpha_j \left(T_{0j}^\alpha \right) \right] - \sum_{\xi=1}^p \left[ T_{0j}^\beta, g_{0j\gamma\alpha}^\xi \right] \alpha_j \left(T_{0j}^\xi \right) - \sum_{\xi=1}^p g_{0j\gamma\alpha}^\xi \left[ T_{0j}^\beta, \alpha_j\left(T_{0j}^\xi \right)\right] \notag
\end{align}}}

Then $\left( \ref{o57} \right)$ is equivalent to 
{\small{\begin{align}\label{o58}
&\left[ T_{0j}^\gamma,  Q_{j\alpha\beta}\right] + \left[ T_{0j}^\alpha, Q_{j\beta\gamma} \right] + \left[ T_{0j}^\beta, Q_{j\gamma \alpha} \right]\\
&= \sum_{\xi=1}^p g_{0j\gamma\alpha}^\xi \left[ T_{0j}^\xi, \alpha_j \left(T_{0j}^\beta \right) \right] - \sum_{\xi=1}^p \left[ \alpha_j\left(T_{0j}^\beta \right), g_{0j\gamma\alpha}^\xi \right] T_{0j}^\xi - \sum_{\xi=1}^p \left[ T_{0j}^\gamma, g_{0j\alpha\beta}^\xi \right] \alpha_j \left(T_{0j}^\xi\right) - \sum_{\xi=1}^p g_{j\alpha\beta}^\xi \left[ T_{0j}^\gamma, \alpha_j\left(T_{0j}^\xi \right) \right] \notag \\
&+ \sum_{\xi=1}^p g_{0j\alpha\beta}^\xi \left[ T_{0j}^\xi, \alpha_j\left(T_{0j}^\gamma\right) \right] - \sum_{\xi=1}^p \left[ \alpha_j \left(T_{0j}^\gamma \right), g_{0j\alpha\beta}^\xi \right] T_{0j}^\xi - \sum_{\xi=1}^p \left[ T_{0j}^\alpha, g_{0j\beta\gamma}^\xi \right] \alpha_j\left(T_{0j}^\xi \right) - \sum_{\xi=1}^p g_{0j\beta\gamma}^\xi \left[ T_{0j}^\alpha, \alpha_j \left(T_{0j}^\xi \right)\right] \notag \\
&+ \sum_{\xi=1}^p g_{0j\beta\gamma}^\xi \left[ T_{0j}^\xi, \alpha_j\left(T_{0j}^\alpha \right) \right] - \sum_{\xi=1}^p   \left[\alpha_j \left(T_{0j}^\alpha\right), g_{0j\beta\gamma}^\xi \right] T_{0j}^\xi - \sum_{\xi=1}^p \left[ T_{0j}^\beta, g_{0j\gamma\alpha}^\xi \right] \alpha_j\left(T_{0j}^\xi \right) - \sum_{\xi=1}^p g_{0j\gamma\alpha}^\xi \left[ T_{0j}^\beta, \alpha_j \left(T_{0j}^\xi \right)\right] \notag\\
&= \sum_{\xi=1}^p g_{0j\gamma\alpha}^\xi Q_{\xi\beta} + \sum_{\xi, \eta=1}^p g_{ 0 j \gamma \alpha}^\xi g_{ 0 j\xi \beta}^\eta \alpha_j \left(T_{0j}^\eta \right)  - \sum_{\xi=1}^p \left[ \alpha_j\left(T_{0j}^\beta \right), g_{0j\gamma\alpha}^\xi \right] T_{0j}^\xi - \sum_{\eta=1}^p \left[ T_{0j}^\gamma, g_{0j\alpha\beta}^\eta \right] \alpha_j \left(T_{0j}^\eta \right) \notag\\
&+\sum_{\xi=1}^p g_{0j\alpha\beta}^\xi Q_{\xi\gamma} + \sum_{\xi, \eta=1}^p g_{0j \alpha \beta}^\xi g_{0j\xi \gamma}^\eta \alpha_j\left(T_{0j}^\eta \right)  - \sum_{\xi=1}^p \left[ \alpha_j\left(T_{0j}^\gamma\right), g_{0j\alpha\beta}^\xi \right] T_{0j}^\xi - \sum_{\eta=1}^p \left[ T_{0j}^\alpha, g_{0j\beta\gamma}^\eta \right] \alpha_j\left(T_{0j}^\eta\right) \notag \\
&+\sum_{\xi=1}^p g_{0j\beta\gamma}^\xi Q_{\xi\alpha} + \sum_{\xi, \eta=1}^p g_{0j \beta \gamma}^\xi g_{0j\xi \alpha}^\eta \alpha_j\left(T_{0j}^\eta\right)  - \sum_{\xi=1}^p \left[ \alpha_j\left(T_{0j}^\alpha\right), g_{0j\beta\gamma}^\xi \right] T_{0j}^\xi - \sum_{\eta=1}^p \left[ T_{0j}^\beta, g_{0j\gamma\alpha}^\eta \right] \alpha_j \left(T_{0j}^\eta\right) \notag
\end{align}}}

Then from $(\ref{o56})$ we see that $(\ref{o58})$ is equivalent to 
\begin{align}\label{o59}
&\left[ T_{0j}^\gamma, Q_{j\alpha\beta} \right] + \left[ T_{0j}^\alpha, Q_{j\beta \gamma} \right] + \left[ T_{0j}^\beta, Q_{j\gamma \alpha} \right] - \sum_{\xi=1}^p g_{0j\gamma\alpha}^\xi Q_{j\xi \beta} - \sum_{\xi=1}^p g_{0j\alpha\beta}^\xi Q_{j\xi \gamma} - \sum_{\xi=1}^p g_{0\beta \gamma}^\xi Q_{j\xi \alpha} \\
&= - \sum_{\xi=1}^p \left[ \alpha_j \left(T_{0j}^\beta \right), g_{0j\gamma\alpha}^\xi \right] T_j^\xi - \sum_{\xi=1}^p \left[ \alpha _j \left(T_j^\gamma \right), g_{0j\alpha\beta}^\xi \right] T_{0j}^\xi - \sum_{\xi=1}^p \left[ \alpha_j\left(T_{0j}^\alpha \right), g_{0j\beta \gamma }^\xi \right] T_{0j}^\xi \notag
\end{align}

\end{remark}

\begin{remark}
We will deduce the equality $(\ref{o63})$ below which will be used in the proof of \textnormal{Lemma $\ref{o54}$}. Let us denote $A_{jk}\left(T_{0j}^\alpha\right):= A_{jk}^\alpha \in \Gamma\left( U_{jk}, \Theta_{\mathcal{F}_0} \right)$, so that $\alpha_k\left(T_{0j}^\alpha \right) - \alpha_j\left( T_{0j}^\alpha \right) + \left[ \theta_{jk}, T_{0j}^\alpha \right]= A_{jk}^\alpha$. We note that
{\small{\begin{align}\label{o60}
&\left[T_{0j}^\alpha, T_{0j}^\beta \right] = \sum_{\gamma=1}^p g_{0j\alpha\beta}^\gamma T_{0j}^\gamma \iff \sum_{\delta, \eta=1}^p \left[r_{0jk}^{\alpha\delta} T_{0k}^\delta, r_{0jk}^{\beta \eta} T_{0k}^\eta \right]= \sum_{\gamma, \xi=1}^p g_{0j\alpha\beta}^\gamma r_{0jk}^{\gamma \xi} T_{0k}^\xi   \notag \\
&\iff \sum_{\delta, \eta=1}^p r_{0jk}^{\alpha \delta} r_{0jk}^{\beta \eta} \left[ T_{0k}^\delta, T_{0k}^\eta \right] + \sum_{\delta, \eta=1}^p r_{0jk}^{\alpha \delta} T_{0k}^\delta \left(r_{0jk}^{\beta \eta} \right) T_{0k}^\eta - \sum_{\delta, \eta=1}^p r_{0jk}^{\beta \eta} T_{0k}^\eta\left( r_{0jk}^{\alpha \delta} \right) T_{0k}^\delta = \sum_{\gamma, \xi=1}^p g_{0j\alpha\beta}^\gamma r_{0jk}^{\gamma \xi} T_{0k}^\xi  \\
&\iff \sum_{\delta, \eta, \xi=1}^p r_{0jk}^{\alpha \delta} r_{0jk}^{\beta \eta} g_{0k\delta \eta}^\xi T_{0k}^\xi + \sum_{\delta, \eta=1}^p r_{0jk}^{\alpha \delta} T_{0k}^\delta \left(r_{0jk}^{\beta \eta} \right) T_{0k}^\eta - \sum_{\delta, \eta=1}^p r_{0jk}^{\beta \eta} T_{0k}^\eta\left( r_{0jk}^{\alpha \delta} \right) T_{0k}^\delta = \sum_{\gamma, \xi=1}^p g_{0j\alpha\beta}^\gamma r_{0jk}^{\gamma \xi} T_{0k}^\xi \notag 
\end{align}}}
From $\left[T_{0j}^\alpha , \alpha_j \left(T_{0j}^\beta\right) \right] - \left[ T_j^\beta, \alpha_j\left(T_{0j}^\alpha \right)\right] - \sum_{\gamma=1}^p g_{0j\alpha\beta}^\gamma\alpha_j\left(T_{0j}^\gamma\right)  =Q_{j\alpha\beta}$ and $\alpha_k\left( T_{0j}^\alpha\right)- \alpha_j\left(T_{0j}^\alpha \right)+\left[ \theta_{jk}, T_{0j}^\alpha \right]=A_{jk}^\alpha$, we have
{\small{\begin{align} \label{o61}
&\left[ T_{0j}^\alpha, \alpha_k\left(T_{0j}^\beta\right)\right] +\left[ T_{0j}^\alpha, \left[ \theta_{jk}, T_{0j}^\beta \right] \right] -\left[ T_{0j}^\alpha, A_{jk}^\beta \right]-\left[ T_{0j}^\beta, \alpha_k \left(T_{0j}^\alpha\right)\right] -\left[ T_{0j}^\beta, \left[ \theta_{jk}, T_{0j}^\alpha \right] \right] +\left[ T_{0j}^\beta, A_{jk}^\alpha \right]  - \sum_{\gamma =1}^p g_{0j\alpha\beta}^\gamma \alpha_j\left( T_{0j}^\gamma\right) = Q_{j\alpha\beta}\\
&\iff \sum_{\delta, \eta=1}^p \left[ r_{0jk}^{\alpha\delta} T_{0k}^\delta, r_{jk}^{\beta \eta} \alpha_k\left(T_{0k}^\eta\right) \right] + \left[\theta_{jk}, \sum_{\gamma=1}^p g_{0j\alpha\beta}^\gamma T_{0j}^\gamma \right] - \left[ T_{0j}^\alpha, A_{jk}^\beta \right]  + \left[ T_{0j}^\beta, A_{jk}^\alpha \right] \notag\\
& - \sum_{\delta, \eta=1}^p \left[ r_{0jk}^{\beta \eta} T_{0k}^\eta,  r_{0jk}^{\alpha \delta} \alpha_k \left(T_{0k}^{\delta}\right)\right] - \sum_{\gamma=1}^p g_{0j\alpha\beta}^\gamma \alpha_k \left(T_{0j}^\gamma \right) - \sum_{\gamma=1}^p g_{0j\alpha\beta}^\gamma \left[ \theta_{jk} , T_{0j}^\gamma \right] +\sum_{\gamma=1}^p g_{0j\alpha\beta}^\gamma A_{jk}^\gamma =  Q_{j\alpha\beta} \notag \\
&\iff \sum_{\delta, \eta=1}^p r_{0jk}^{\alpha \delta} r_{0jk}^{\beta\eta}\left[ T_{0k}^\delta,\alpha_k \left(T_{0k}^\eta \right) \right] + \sum_{\delta, \eta=1}^p r_{0jk}^{\alpha \delta} T_{0k}^\delta\left( r_{0jk}^{\beta \eta} \right) \alpha_k\left(T_{0k}^\eta\right) - \sum_{\delta, \eta=1}^p r_{0jk}^{\beta \eta} \left[ \alpha_k\left(T_{0k}^\eta\right), r_{0jk}^{\alpha \delta} \right] T_{0k}^\delta + \sum_{\gamma=1}^p \left[ \theta_{jk}, g_{0j\alpha\beta}^\gamma \right] T_{0j}^\gamma \notag \\
& - \sum_{\delta, \eta=1}^p r_{0jk}^{\alpha \delta} r_{0jk}^{\beta \eta}\left[T_{0k}^\eta, \alpha_k\left(T_{0k}^\delta \right) \right]  - \sum_{\delta, \eta=1}^p r_{0jk}^{\beta \eta} T_{0k}^\eta\left(r_{0jk}^{\alpha \delta} \right) \alpha_k \left(T_{0k}^\delta \right) + \sum_{\delta, \eta=1}^p r_{jk}^{\alpha \delta} \left[\alpha_k \left(T_{0k}^\delta \right), r_{0jk}^{\beta \eta} \right] T_{0k}^\eta - \sum_{\gamma, \xi=1}^p g_{0j\alpha\beta}^\gamma r_{0jk}^{\gamma \xi} \alpha_k \left(T_{0k}^\xi \right) \notag \\
&- \left[ T_{0j}^\alpha, A_{jk}^\beta \right] + \left[ T_{0j}^\beta, A_{jk}^\alpha \right] + \sum_{\gamma=1}^p g_{j\alpha\beta}^\gamma A_{jk}^\gamma = Q_{j\alpha\beta} \notag
\end{align}}}

Then $(\ref{o60})$ implies that $\left( \ref{o61} \right)$ is equivalent to
{\small{\begin{align}
&\sum_{\delta, \eta=1}^p r_{0jk}^{\alpha \delta} r_{0jk}^{\beta\eta}\left[ T_{0k}^\delta, \alpha_k \left(T_{0k}^\eta \right) \right] -\sum_{\delta, \eta, \xi=1}^p r_{0jk}^{\alpha \delta} r_{0jk}^{\beta \eta} g_{0k\delta \eta}^\xi \alpha_k\left(T_{0k}^\xi \right) - \sum_{\delta, \eta=1}^p r_{jk}^{\beta \eta} \left[ \alpha_k\left(T_{0k}^\eta \right), r_{0jk}^{\alpha \delta} \right] T_{0k}^\delta + \sum_{\gamma =1}^p \left[ \theta_{jk}, g_{0j\alpha\beta}^\gamma \right] T_{0j}^\gamma \\
& - \sum_{\delta, \eta=1}^p r_{0jk}^{\alpha \delta} r_{0jk}^{\beta \eta}\left[T_{0k}^\eta, \alpha_k \left(T_{0k}^\delta \right) \right]  + \sum_{\delta, \eta=1}^p r_{0jk}^{\alpha \delta} \left[ \alpha_k \left(T_{0k}^\delta\right), r_{0jk}^{\beta \eta} \right] T_{0k}^\eta \notag - \left[ T_{0j}^\alpha, A_{jk}^\beta \right] + \left[ T_{0j}^\beta, A_{jk}^\alpha \right] + \sum_{\gamma=1}^p g_{0j\alpha\beta}^\gamma A_{jk}^\gamma = Q_{j\alpha\beta} \notag
\end{align}}}
which is equivalent to
\begin{align} \label{o63}
& \left[ T_{0j}^\alpha, A_{jk}^\beta \right] - \left[ T_{0j}^\beta, A_{jk}^\alpha \right] - \sum_{\gamma=1}^p g_{0j\alpha\beta}^\gamma A_{jk}^\gamma\\
&=\sum_{\delta, \eta =1}^p r_{0jk}^{\alpha \delta} r_{0jk}^{\beta\eta}Q_{k\delta \eta} - \sum_{\delta, \eta=1}^p r_{0jk}^{\beta \eta} \left[ \alpha_k \left(T_{0k}^\eta\right), r_{0jk}^{\alpha \delta} \right] T_{0k}^\delta + \sum_{\gamma=1}^p \left[ \theta_{jk}, g_{0j\alpha\beta}^\gamma \right] T_{0j}^\gamma + \sum_{\delta, \eta=1}^p r_{0jk}^{\alpha \delta} \left[\alpha_k\left(T_{0k}^\delta\right), r_{0jk}^{\beta \eta} \right] T_{0k}^\eta- Q_{j\alpha\beta} \notag
\end{align}

\end{remark}

\begin{proof}[Proof of \textnormal{Lemma \ref{o54}}]
First we note that by putting $K_{jk}:=\left(\theta_{jk}, A_{jk} \right)\in \Gamma\left( U_{jk}, \mathcal{E}_{\Theta_{\mathcal{F}_0}} \right)$, so that $K_{ij}-K_{ik}+ K_{jk}=0$, we see that $\left\{ \xi_{ijk}:=\left[ K_{ij} , K_{jk}\right]\right\}\in C^2\left(\mathcal{U}, \mathcal{E}_{\Theta_{\mathcal{F}_0}} \right)$. In fact,
\begin{align*}
\xi_{ijk}&=\left[ K_{ik}, K_{kj} \right] = - \left[K_{ij}+ K_{jk}, K_{jk}  \right] = -\left[ K_{ij}, K_{jk} \right]= -\xi_{ijk}\\
\xi_{jik}&=\left[ K_{ji}, K_{ik} \right]=- \left[ K_{ij}, K_{ij} + K_{jk} \right]= - \left[ K_{ij}, K_{jk} \right]= -\xi_{ijk}
\end{align*}
and $\xi_{jkl}- \xi_{ikl} + \xi_{ijl} - \xi_{ijk}=\left[ K_{jk}- K_{ik}, K_{kl} \right] + \left[ K_{ij}, K_{jl}- K_{jk} \right] = - \left[ K_{ij}, K_{kl}\right] - \left[ K_{ij}, K_{lk} \right]=0$. This implies that $-\delta\left( \left\{ \left[ \left(\theta_{ij}, A_{ij} \right) , \left( \theta_{jk}, A_{jk} \right)\right] \right\} \right)=0$.

We show that $D_2\left(R_j \right) =0$, equivalently $D_2\left(B_j:=\frac{1}{2} R_j \right)=0$. In fact, from $(\ref{o59})$, 
{\Small{\begin{align}
&D_2\left(B_j\right)\left(T_{0j}^\alpha \wedge T_{0j}^\beta\wedge T_{0j}^\gamma\right)\\
& =\left[T_{0j}^\alpha, B_j\left(T_{0j}^\beta, T_{0j}^\gamma\right)\right]-\left[T_{0j}^\beta, B_j\left(T_{0j}^\alpha, T_{0j}^\gamma\right)\right]+ \left[T_{0j}^\gamma, B_j\left(T_{0j}^\alpha ,  T_{0j}^\beta\right) \right] - B_j\left( \left[T_{0j}^\alpha, T_{0j}^\beta \right], T_{0j}^\gamma \right) + B_j\left(\left[T_{0j}^\alpha, T_{0j}^\gamma \right], T_{0j}^\beta \right) -B_j\left( \left[T_{0j}^\beta, T_{0j}^\gamma  \right], T_{0j}^\alpha \right) \notag \\
&= \left[ T_{0j}^\alpha, \left[ \alpha _j \left(T_{0j}^\beta \right) , \alpha_j\left(T_{0j}^\gamma \right)  \right]   \right] - \left[ T_{0j}^\alpha, \alpha_j\left(Q_{j\beta \gamma}\right)  \right] -  \left[ T_{0j}^\beta, \left[\alpha_j\left(T_{0j}^\alpha \right) , \alpha_j \left(T_{0j}^\gamma \right)  \right]   \right] +  \left[ T_{0j}^\beta, \alpha_j \left(Q_{j\alpha \gamma}\right)   \right] \notag \\
&+ \left[ T_{0j}^\gamma, \left[ \alpha_j \left(T_{0j}^\alpha \right) , \alpha_j \left( T_{0j}^\beta \right)  \right]   \right] - \left[ T_{0j}^\gamma, \alpha_j \left(Q_{j\alpha \beta} \right) \right] - \sum_{\xi=1}^p g_{0j\alpha\beta}^\xi B_j\left( T_{0j}^\xi, T_{0j}^\gamma  \right)   + \sum_{\xi=1}^p g_{0j\alpha\gamma}^\xi B_j\left( T_{0j}^\xi, T_{0j}^\beta  \right)  - \sum_{\xi=1}^p g_{0j\beta\gamma}^\xi B_j\left( T_{0j}^\xi, T_{0j}^\alpha  \right) \notag \\
&= \left[ \left[T_{0j}^\alpha, \alpha_j \left(T_{0j}^\beta \right) \right] - \left[T_{0j}^\beta , \alpha_j \left(T_{0j}^\alpha \right) \right]          , \alpha _j \left(T_{0j}^\gamma \right)   \right] + \left[ \left[T_{0j}^\beta, \alpha_j \left(T_{0j}^\gamma \right) \right] -\left[ T_{0j}^\gamma, \alpha_j\left(T_{0j}^\beta \right) \right],  \alpha_j \left(T_{0j}^\alpha \right) \right] +\left[\left[T_{0j}^\gamma, \alpha_j\left(T_{0j}^\alpha \right)\right] - \left[ T_{0j}^\alpha, \alpha_j \left(T_{0j}^\gamma \right)   \right], \alpha_j \left(T_{0j}^\beta \right)    \right] \notag \\
&-  \left[ T_{0j}^\alpha,   \alpha_j \left(Q_{j\beta\gamma} \right) \right]   + \left[ T_{0j}^\beta, \alpha_j ( Q_{j\alpha \gamma} ) \right]   -  \left[ T_{0j}^\gamma,  \alpha_j \left(Q_{j\alpha\beta} \right) \right] -\sum_{\xi=1}^p g_{0j\alpha\beta}^\xi \left[\alpha_j\left(T_{0j}^\xi \right), \alpha_j \left( T_{0j}^\gamma \right)  \right] + \sum_{\xi=1}^p g_{0j\alpha\beta}^\xi \alpha_j \left(Q_{j\xi \gamma} \right) \notag \\
&+\sum_{\xi=1}^p g_{0j\alpha \gamma}^\xi \left[\alpha_j \left(T_{0j}^\xi \right), \alpha_j \left( T_{0j}^\beta \right)  \right] - \sum_{\xi =1}^p g_{0j\alpha\gamma}^\xi \alpha_j \left( Q_{j\xi \beta} \right) -\sum_{\xi=1}^p g_{0j\beta\gamma}^\xi \left[\alpha_j \left(T_{0j}^\xi \right), \alpha_j\left( T_{0j}^\alpha\right)  \right] + \sum_{\xi=1}^p g_{0j\beta\gamma}^\xi \alpha_j\left(Q_{j\xi \alpha} \right) \notag \\
&= \sum_{\xi=1}^p \left[ g_{0j\alpha\beta}^\xi \alpha_j\left(T_{0j}^\xi\right), \alpha_j\left(T_{0j}^\gamma \right) \right] + \left[  Q_{j\alpha\beta},  \alpha_j \left(T_{0j}^\gamma \right)    \right] + \sum_{\xi=1}^p \left[ g_{0j\beta\gamma}^\xi \alpha_j \left(T_{0j}^\xi \right), \alpha_j \left(T_{0j}^\alpha \right) \right] +\left[  Q_{j\beta\gamma},  \alpha_j \left(T_{0j}^\alpha\right)    \right] +\sum_{\xi=1}^p \left[ g_{0j\gamma\alpha}^\xi \alpha_j \left(T_{0j}^\xi \right), \alpha_j\left(T_{0j}^\beta \right) \right] \notag  \\
& + \left[  Q_{j\gamma\alpha},  \alpha_j \left(T_{0j}^\beta \right)    \right] -  \left[ T_{0j}^\alpha,  \alpha_j \left(Q_{j\beta\gamma} \right) \right]   + \left[ T_{0j}^\beta, \alpha_j \left( Q_{j\alpha \gamma} \right) \right]   -  \left[ T_{0j}^\gamma,  \alpha_j \left(Q_{j\alpha\beta} \right) \right] -\sum_{\xi=1}^p g_{0j\alpha\beta}^\xi \left[ \alpha_j \left(T_{0j}^\xi \right),  \alpha_j \left( T_{0j}^\gamma \right)  \right] + \sum_{\xi=1}^p g_{0j\alpha\beta}^\xi \alpha_j \left(Q_{j\xi \gamma} \right) \notag \\
&+\sum_{\xi=1}^p g_{0j\alpha \gamma}^\xi \left[\alpha_j\left(T_{0j}^\xi \right), \alpha_j \left( T_{0j}^\beta \right)  \right] - \sum_{\xi =1}^p g_{0j\alpha\gamma}^\xi \alpha_j\left( Q_{j\xi \beta} \right) -\sum_{\xi=1}^p g_{0j\beta\gamma}^\xi \left[\alpha_j\left(T_{0j}^\xi \right), \alpha_j \left( T_{0j}^\alpha \right)  \right] + \sum_{\xi=1}^p g_{0j\beta\gamma}^\xi \alpha_j\left(Q_{j\xi \alpha} \right) \notag\\
&\equiv_{\Theta_{\mathcal{F}_0}} - \sum_{\xi=1}^p \left[\alpha_j\left(T_{0j}^\gamma \right), g_{0j\alpha\beta}^\xi \right]\alpha_j\left(T_{0j}^\xi\right)  - \alpha_j\left(\left[ T_{0j}^\gamma, Q_{j\alpha\beta}\right] \right) - \sum_{\xi=1}^p \left[ \alpha_j \left(T_{0j}^\alpha \right), g_{0j\beta\gamma}^\xi \right]\alpha_j\left(T_{0j}^\xi \right)  - \alpha_j\left(\left[ T_{0j}^\alpha, Q_{j\beta\gamma}\right] \right) \notag \\
& - \sum_{\xi=1}^p \left[\alpha_j \left(T_{0j}^\beta \right), g_{0j\gamma\alpha}^\xi \right] \alpha_j \left(T_{0j}^\xi \right)  - \alpha_j\left(\left[ T_{0j}^\beta, Q_{j\gamma\alpha}\right] \right) + \sum_{\xi=1}^p g_{0j\alpha\beta}^\xi \alpha_j (Q_{j\xi \gamma} ) - \sum_{\xi =1}^p g_{0j\alpha\gamma}^\xi \alpha_j \left( Q_{j\xi \beta}\right)  +  \sum_{\xi=1}^p g_{0j\beta\gamma}^\xi \alpha_j\left(Q_{j\xi \alpha} \right) = 0 \notag
\end{align}}}

We show that $-\delta\left( \left\{ R_j \right\} \right)+ D_1'\left( \left\{ F_{jk} \right\} \right)=0$. In fact, we note that
{\small{\begin{align}\label{o64}
R_k\left(T_{0j}^\alpha \wedge T_{0j}^\beta \right)= \sum_{\delta, \eta=1}^p r_{0jk}^{\alpha \delta} r_{0jk}^{\beta \eta} R_k\left(T_{0k}^\delta \wedge T_{0k}^\eta\right)= 2\sum_{\delta, \eta=1}^p r_{0jk}^{\alpha \delta} r_{0jk}^{\beta \eta}\left[ \alpha_k\left(T_{0k}^\delta\right), \alpha_k\left(T_{0k}^\eta\right)  \right] - 2 \sum_{\delta, \eta=1}^p \sum_{\gamma=1}^p r_{0jk}^{\alpha \delta} r_{0jk}^{\beta \eta}  \alpha_k\left(Q_{k\delta \eta} \right)
\end{align}}}
and
{\small{\begin{align}\label{o65}
&R_j \left(T_{0j}^\alpha \wedge T_{0j}^\beta\right) =2\left[ \alpha_j\left(T_{0j}^\alpha\right), \alpha_j \left(T_{0j}^\beta \right) \right]- 2 \alpha_j\left( Q_{j\alpha\beta}\right)  = \left[  \alpha_j(T_{0j}^\alpha), \alpha_j \left(T_{0j}^\beta \right) \right] -\left[ \alpha_j \left(T_{0j}^\beta \right), \alpha_j \left(T_{0j}^\alpha \right)    \right] - 2  \alpha_j\left( Q_{j\alpha\beta}\right) \\
&=  \left[ \alpha_k\left(T_{0j}^\alpha\right) + \left[\theta_{jk}, T_{0j}^\alpha \right] - A_{jk}^\alpha    ,  \alpha_j \left(T_{0j}^\beta \right)   \right] -\left[ \alpha_k \left(T_{0j}^\beta\right) + \left[ \theta_{jk}, T_{0j}^\beta \right] - A_{jk}^\beta, \alpha_j \left(T_{0j}^\alpha \right) \right] - 2  \alpha_j\left( Q_{j\alpha\beta}\right) \notag \\
&=\left[\alpha_k \left(T_{0j}^\alpha\right), \alpha_j\left(T_{0j}^\beta\right) \right] +\left[ \left[ \theta_{jk}, T_{0j}^\alpha \right], \alpha_j \left(T_{0j}^\beta\right)  \right] -\left[ A_{jk}^\alpha, \alpha_j \left(T_{0j}^\beta \right)\right] \notag \\
&-\left[ \alpha_k \left(T_{0j}^\beta \right), \alpha_j \left(T_{0j}^\alpha \right)  \right] -\left[ \left[ \theta_{jk}, T_{0j}^\beta \right], \alpha_j \left(T_{0j}^\alpha \right)   \right] +\left[ V_{jk}^\beta, \alpha_j \left(T_{0j}^\alpha \right) \right] - 2  \alpha_j\left( Q_{j\alpha\beta}\right) \notag \\
&=\left[ \alpha_k \left(T_{0j}^\alpha\right), \alpha_k \left(T_{0j}^\beta\right) +\left[\theta_{jk}, T_{0j}^\beta \right] - A_{jk}^\beta \right] + \left[ \theta_{jk}, \left[ T_{0j}^\alpha, \alpha_j \left(T_{0j}^\beta \right)\right] \right] -\left[ T_{0j}^\alpha , \left[ \theta_{jk}, \alpha_j \left(T_{0j}^\beta\right)\right] \right] -\left[A_{jk}^\alpha, \alpha_j \left(T_{0j}^\beta\right) \right] \notag \\
&-\left[ \alpha_k \left(T_{0j}^\beta\right),  \alpha_k \left(T_{0j}^\alpha\right) +\left[ \theta_{jk}, T_{0j}^\alpha   \right] - V_{jk}^\alpha   \right] -\left[ \theta_{jk}, \left[ T_j^\beta,\alpha_j\left(T_{0j}^\alpha\right)\right] \right] +\left[ T_{0j}^\beta, \left[ \theta_{jk}, \tilde{\alpha}_j(T_j^\alpha)\right] \right] +\left[ A_{jk}^\beta,  \alpha_j\left(T_{0j}^\alpha\right) \right]  - 2  \alpha_j\left( Q_{j\alpha\beta}\right) \notag \\
&=  2\left[ \alpha_k \left(T_j^\alpha\right), \alpha_k \left(T_{0j}^\beta\right)   \right] + \left[ \alpha_k \left(T_{0j}^\alpha\right), \left[ \theta_{jk}, T_{0j}^\beta \right]  \right] -\left[ \alpha_k \left(T_{0j}^\alpha\right), A_{jk}^\beta  \right] + \left[ \theta_{jk}, \left[ T_{0j}^\alpha, \alpha_j \left(T_{0j}^\beta \right)\right] \right] -\left[ T_{0j}^\alpha , \left[ \theta_{jk}, \alpha_j \left(T_{0j}^\beta\right)\right] \right] - \left[A_{jk}^\alpha, \alpha_j  \left(T_{0j}^\beta \right) \right] \notag \\
&-\left[ \alpha_k\left(T_{0j}^\beta\right), \left[ \theta_{jk}, T_{0j}^\alpha \right] \right] +\left[ \alpha_k\left(T_{0j}^\beta\right), A_{jk}^\beta \right] -\left[ \theta_{jk}, \left[ T_j^\beta, \alpha_j\left(T_{0j}^\alpha \right)\right] \right] +\left[ T_{0j}^\beta, \left[ \theta_{jk}, \alpha_j\left(T_{0j}^\alpha \right)\right] \right] +\left[ A_{jk}^\beta,  \alpha_j \left(T_{0j}^\alpha\right) \right]  - 2  \alpha_j\left( Q_{j\alpha\beta}\right) \notag
\end{align}}}
We note that
{\small{\begin{align}\label{o68}
&\left[\alpha_k(T_{0j}^\alpha), \alpha_k\left(T_{0j}^\beta\right) \right] =\sum_{\delta, \eta=1}^p \left[ r_{0jk}^{\alpha\delta} \alpha_k\left(T_k^\delta\right) , r_{jk}^{\beta \eta} \alpha_k\left(T_{0k}^\eta\right) \right]\\
& = \sum_{\delta, \eta=1}^p r_{0jk}^{\alpha \delta} r_{0jk}^{\beta \eta} \left[ \alpha_k\left(T_{0k}^\delta \right), \alpha_k \left(T_{0k}^\eta\right)\right] + \sum_{\delta, \eta=1}^p r_{0jk}^{\alpha \delta}\left[ \alpha_k \left(T_{0k}^\delta \right) , r_{0jk}^{\beta \eta}\right] \alpha_k \left(T_{0k}^\eta\right) - \sum_{\delta, \eta=1}^p r_{0jk}^{\beta \eta}\left[ \alpha_k \left(T_{0k}^\eta \right ), r_{0jk}^{\alpha \delta} \right] \alpha_k\left(T_{0k}^\delta\right) \notag
\end{align}}}

From $(\ref{o68})$, $(\ref{o65})$ is equivalent to
{\small{\begin{align}\label{o69}
&R_j \left(T_{0j}^\alpha \wedge T_{0j}^\beta\right)\\
&=  2 \sum_{\delta, \eta=1}^p r_{0jk}^{\alpha \delta} r_{0jk}^{\beta \eta} \left[ \alpha_k\left(T_{0k}^\delta \right), \alpha_k \left(T_{0k}^\eta\right)\right] + 2\sum_{\delta, \eta=1}^p r_{0jk}^{\alpha \delta}\left[ \alpha_k \left(T_{0k}^\delta \right) , r_{0jk}^{\beta \eta}\right] \alpha_k \left(T_{0k}^\eta\right) - 2\sum_{\delta, \eta=1}^p r_{0jk}^{\beta \eta}\left[ \alpha_k \left(T_{0k}^\eta \right ), r_{0jk}^{\alpha \delta} \right] \alpha_k\left(T_{0k}^\delta\right) \notag\\
&+ \left[ \alpha_k \left(T_{0j}^\alpha\right), \left[ \theta_{jk}, T_{0j}^\beta \right]  \right] -\left[ \alpha_k \left(T_{0j}^\alpha\right), A_{jk}^\beta  \right] + \left[ \theta_{jk}, \left[ T_{0j}^\alpha, \alpha_j \left(T_{0j}^\beta \right)\right] \right] -\left[ T_{0j}^\alpha , \left[ \theta_{jk}, \alpha_j \left(T_{0j}^\beta\right)\right] \right] - \left[A_{jk}^\alpha, \alpha_j  \left(T_{0j}^\beta \right) \right] \notag \\
&-\left[ \alpha_k\left(T_{0j}^\beta\right), \left[ \theta_{jk}, T_{0j}^\alpha \right] \right] +\left[ \alpha_k\left(T_{0j}^\beta\right), A_{jk}^\beta \right] -\left[ \theta_{jk}, \left[ T_j^\beta, \alpha_j\left(T_{0j}^\alpha \right)\right] \right] +\left[ T_{0j}^\beta, \left[ \theta_{jk}, \alpha_j\left(T_{0j}^\alpha \right)\right] \right] +\left[ A_{jk}^\beta,  \alpha_j \left(T_{0j}^\alpha\right) \right]  - 2  \alpha_j\left( Q_{j\alpha\beta}\right) \notag
\end{align}}}
On the other hand, we have
{\small{\begin{align}\label{o67}
&D_1'\left(F_{jk}\right)\left(T_{0j}^\alpha \wedge T_{0j}^\beta \right)= \left[ T_{0j}^\alpha, F_{jk}\left(T_{0j}^\beta \right)   \right] - \left[ T_{0j}^\beta, F_{jk}\left(T_{0j}^\alpha \right)    \right]  - \sum_{\gamma=1}^p g_{0j\alpha\beta}^\gamma F_{jk}(T_{0j}^\gamma)\\
&=\left[ T_{0j}^\alpha, -\alpha_j\left(A_{jk}^\beta\right) - \alpha_k\left(A_{jk}^\beta\right) +\left[ \theta_{jk}, \alpha_j \left(T_{0j}^\beta \right)   \right]  +\left[\theta_{jk}, \alpha_k\left( T_{0j}^\beta \right) \right] \right] + \left[ T_{0j}^\beta, \alpha_j\left(A_{jk}^\alpha \right) + \alpha_k\left(A_{jk}^\alpha\right) - \left[ \theta_{jk}, \alpha_j \left(T_{0j}^\alpha \right)   \right]  -\left[\theta_{jk}, \alpha_k\left( T_{0j}^\alpha \right) \right] \right] \notag \\
&+\sum_{\gamma=1}^p g_{0j\alpha\beta}^\gamma \alpha_j\left(A_{jk}^\gamma \right) + \sum_{\gamma=1}^p g_{0j\alpha\beta}^\gamma \alpha_k\left( A_{jk}^\gamma \right) - \sum_{\gamma=1}^p g_{0j\alpha\beta}^\gamma \left[ \theta_{jk}, \alpha_j \left( T_{0j}^\gamma \right)\right] - \sum_{\gamma=1}^p g_{0j\alpha\beta}^\gamma \left[ \theta_{jk},\alpha_k \left( T_{0j}^\gamma \right) \right] \notag \\
&=- \left[T_{0j}^\alpha, \alpha_j\left(A_{jk}^\beta \right)   \right] -\left[ T_{0j}^\alpha, \alpha_k\left(A_{jk}^\beta\right)  \right] +\left[ T_{0j}^\alpha,  \left[ \theta_{jk}, \alpha_j\left(T_{0j}^\beta\right) \right] \right] + \left[T_{0j}^\alpha, \left[  \theta_{jk},\alpha_k\left(T_{0j}^\beta \right)\right] \right] \notag \\
&+\left[ T_{0j}^\beta, \alpha_j\left(A_{jk}^\alpha\right) \right] +\left[ T_{0j}^\beta,\alpha_k \left(A_{jk}^\alpha\right)\right] -\left[ T_{0j}^\beta, \left[ \theta_{jk}, \alpha_j \left(T_{0j}^\alpha \right) \right] \right] - \left[ T_{0j}^\beta, \left[ \theta_{jk}, \alpha_k\left(T_{0j}^\alpha \right)\right] \right] \notag \\
&+\sum_{\gamma=1}^p g_{0j\alpha\beta}^\gamma \alpha_j\left(A_{jk}^\gamma \right) + \sum_{\gamma=1}^p g_{0j\alpha\beta}^\gamma \alpha_k\left( A_{jk}^\gamma \right) - \left[ \theta_{jk}, \alpha_j \left( \left[ T_{0j}^\alpha, T_{0j}^\beta\right] \right)\right] + \sum_{\gamma=1}^p \left[\theta_{jk}, g_{0j\alpha\beta}^\gamma \right] \alpha_j\left(T_{0j}^\gamma\right) \notag \\
& - \left[ \theta_{jk}, \alpha_k \left( \left[ T_{0j}^\alpha, T_{0j}^\beta\right] \right)\right] + \sum_{\gamma=1}^p \left[\theta_{jk}, g_{j\alpha\beta}^\gamma \right] \alpha_k\left(T_{0j}^\gamma\right) \notag \\
&=- \left[T_{0j}^\alpha, \alpha_j\left(A_{jk}^\beta \right)   \right] -\left[ T_{0j}^\alpha, \alpha_k\left(A_{jk}^\beta \right)  \right] +\left[ T_{0j}^\alpha,  \left[ \theta_{jk}, \alpha_j\left(T_{0j}^\beta \right) \right] \right] + \left[T_{0j}^\alpha, \left[  \theta_{jk},\alpha_k\left(T_{0j}^\beta \right)\right] \right] \notag \\
&+\left[ T_{0j}^\beta, \alpha_j \left( A_{jk}^\alpha \right) \right] +\left[ T_{0j}^\beta, \alpha_k\left(A_{jk}^\alpha \right)\right] -\left[ T_j^\beta, \left[ \theta_{jk}, \alpha_j\left(T_{0j}^\alpha \right) \right] \right] - \left[ T_j^\beta, \left[ \theta_{jk}, \alpha_k \left(T_{0j}^\alpha \right)\right] \right] \notag \\
&+\sum_{\gamma=1}^p g_{0j\alpha\beta}^\gamma \tilde{\alpha}_j\left(A_{jk}^\gamma \right) + \sum_{\gamma=1}^p g_{0j\alpha\beta}^\gamma \alpha_k\left( A_{jk}^\gamma \right)  - \left[ \theta_{jk}, \alpha_j \left( \left[ T_j^\alpha, T_j^\beta\right] \right)\right]  - \left[ \theta_{jk}, \alpha_k \left( \left[ T_{0j}^\alpha, T_{0j}^\beta\right] \right)\right]  \notag \\
&+ 2\sum_{\gamma=1}^p \left[\theta_{jk}, g_{0j\alpha\beta}^\gamma \right]\alpha_k\left(T_{0j}^\gamma\right) + \sum_{\gamma=1}^p \left[ \theta_{jk}, g_{0j\alpha\beta}^\gamma\right] \left[ \theta_{jk}, T_{0j}^\gamma \right] - \sum_{\gamma=1}^q \left[ \theta_{jk}, g_{j\alpha\beta}^\gamma \right] A_{jk}^\gamma \notag
\end{align}}}

Then from $(\ref{o64}), (\ref{o69})$ and $(\ref{o67})$, we have

{\small{\begin{align}
&\left(-R_k+ R_j+ D_1'\left(F_{jk} \right)\right)\left(T_{0j}^\alpha \wedge T_{0j}^\beta \right)\\
&=- \cancel{2\sum_{\delta, \eta=1}^p r_{0jk}^{\alpha \delta} r_{0jk}^{\beta \eta}\left[ \alpha_k\left(T_{0k}^\delta\right), \alpha_k\left(T_{0k}^\eta\right)  \right]} + 2 \sum_{\delta, \eta=1}^p \sum_{\gamma=1}^p r_{0jk}^{\alpha \delta} r_{0jk}^{\beta \eta}  \alpha_k\left(Q_{k\delta \eta} \right)  \notag \\
&+  \cancel{2 \sum_{\delta, \eta=1}^p r_{0jk}^{\alpha \delta} r_{0jk}^{\beta \eta} \left[ \alpha_k\left(T_{0k}^\delta \right), \alpha_k \left(T_{0k}^\eta\right)\right] } + 2\sum_{\delta, \eta=1}^p r_{0jk}^{\alpha \delta}\left[ \alpha_k \left(T_{0k}^\delta \right) , r_{0jk}^{\beta \eta}\right] \alpha_k \left(T_{0k}^\eta\right) - 2\sum_{\delta, \eta=1}^p r_{0jk}^{\beta \eta}\left[ \alpha_k \left(T_{0k}^\eta \right ), r_{0jk}^{\alpha \delta} \right] \alpha_k\left(T_{0k}^\delta\right) \notag\\
&+ \left[ \alpha_k \left(T_{0j}^\alpha\right), \left[ \theta_{jk}, T_{0j}^\beta \right]  \right] -\left[ \alpha_k \left(T_{0j}^\alpha\right), A_{jk}^\beta  \right] + \left[ \theta_{jk}, \left[ T_{0j}^\alpha, \alpha_j \left(T_{0j}^\beta \right)\right] \right] -\xcancel{\left[ T_{0j}^\alpha , \left[ \theta_{jk}, \alpha_j \left(T_{0j}^\beta\right)\right] \right] }- \left[A_{jk}^\alpha, \alpha_j  \left(T_{0j}^\beta \right) \right] \notag \\
&-\left[ \alpha_k\left(T_{0j}^\beta\right), \left[ \theta_{jk}, T_{0j}^\alpha \right] \right] +\left[ \alpha_k\left(T_{0j}^\beta\right), A_{jk}^\alpha \right] -\left[ \theta_{jk}, \left[ T_{0j}^\beta, \alpha_j\left(T_{0j}^\alpha \right)\right] \right] + \bcancel{\left[ T_{0j}^\beta, \left[ \theta_{jk}, \alpha_j\left(T_{0j}^\alpha \right)\right] \right] }+\left[ A_{jk}^\beta,  \alpha_j \left(T_{0j}^\alpha\right) \right]  - 2  \alpha_j\left( Q_{j\alpha\beta}\right) \notag \\
&- \left[T_{0j}^\alpha, \alpha_j\left(A_{jk}^\beta \right)   \right] -\left[ T_{0j}^\alpha, \alpha_k\left(A_{jk}^\beta \right)  \right] + \xcancel{\left[ T_{0j}^\alpha,  \left[ \theta_{jk}, \alpha_j\left(T_{0j}^\beta \right) \right] \right]} + \left[T_{0j}^\alpha, \left[  \theta_{jk},\alpha_k\left(T_{0j}^\beta \right)\right] \right] \notag \\
&+\left[ T_{0j}^\beta, \alpha_j \left( A_{jk}^\alpha \right) \right] +\left[ T_{0j}^\beta, \alpha_k\left(A_{jk}^\alpha \right)\right] - \bcancel{\left[ T_{0j}^\beta, \left[ \theta_{jk}, \alpha_j\left(T_{0j}^\alpha \right) \right] \right]} - \left[ T_j^\beta, \left[ \theta_{jk}, \alpha_k \left(T_{0j}^\alpha \right)\right] \right] \notag \\
&+\sum_{\gamma=1}^p g_{0j\alpha\beta}^\gamma \alpha_j\left(A_{jk}^\gamma \right) + \sum_{\gamma=1}^p g_{0j\alpha\beta}^\gamma \alpha_k\left( A_{jk}^\gamma \right)  - \left[ \theta_{jk}, \alpha_j \left( \left[ T_{0j}^\alpha, T_{0j}^\beta\right] \right)\right]  - \left[ \theta_{jk}, \alpha_k \left( \left[ T_{0j}^\alpha, T_{0j}^\beta\right] \right)\right]  \notag \\
&+ 2\sum_{\gamma=1}^p \left[\theta_{jk}, g_{0j\alpha\beta}^\gamma \right]\alpha_k\left(T_{0j}^\gamma\right) + \sum_{\gamma=1}^p \left[ \theta_{jk}, g_{0j\alpha\beta}^\gamma\right] \left[ \theta_{jk}, T_{0j}^\gamma \right] - \sum_{\gamma=1}^q \left[ \theta_{jk}, g_{j\alpha\beta}^\gamma \right] A_{jk}^\gamma \notag
\end{align}}}

We note that
{\small{\begin{align}\label{o70}
&\alpha_k\left(\left[ T_{0j}^\alpha, T_{0j}^\beta \right] \right) - \alpha_j\left(\left[T_{0j}^\alpha, T_{0j}^\beta \right] \right) +\left[\theta_{jk}, \left[ T_{0j}^\alpha, T_{0j}^\beta\right] \right]\\
&= \sum_{\gamma=1}^p g_{0j\alpha\beta}^\gamma \alpha_k \left(T_{0j}^\gamma\right) - \sum_{\gamma=1}^p g_{0j\alpha\beta}^\gamma \alpha_j \left(T_{0j}^\gamma \right) + \sum_{\gamma=1}^p g_{0j\alpha\beta}^\gamma\left[\theta_{jk}, T_{0j}^\gamma \right] + \sum_{\gamma=1}^p \left[ \theta_{jk}, g_{0j\alpha\beta}^\gamma \right] T_{0j}^\gamma = \sum_{\gamma =1}^p g_{0j\alpha\beta}^\gamma  A_{jk}^\gamma  + \sum_{\gamma=1}^p \left[\theta_{jk}, g_{0j\alpha\beta}^\gamma \right]T_{0j}^\gamma \notag
\end{align}}}

Then $(\ref{o63})$, $(\ref{o70})$ and $\alpha_j\left(Q_{j\alpha\beta} \right)\equiv_{\Theta_{\mathcal{F}_0}} \alpha_k\left(Q_{j\alpha\beta} \right) + \left[ \theta_{jk}, Q_{j\alpha\beta}\right]$  implies that
{\small{\begin{align*}
&\left(-R_k+R_j +D_1(F_{jk}) \right)\left(T_j^\alpha \wedge T_j^\beta\right)\\
&\equiv_{\Theta_{\mathcal{F}_0} } \cancel{2\alpha_k\left( \left[ T_{0j}^\alpha,  A_{jk}^\beta \right]\right)}- \bcancel{2\alpha_k\left(\left[ T_{0j}^\beta, A_{jk}^\alpha \right] \right)} - \xcancel{2\sum_{\gamma=1}^p g_{j\alpha\beta}^\gamma \alpha_k\left( A_{jk}^\gamma\right)}\\
&+ \left[ \alpha_k \left(T_{0j}^\alpha\right), \left[ \theta_{jk}, T_{0j}^\beta \right]  \right]^1 + \left[ \theta_{jk}, \left[ T_j^\alpha, \alpha_j \left(T_{0j}^\beta\right)\right] \right]  -\bcancel{2\left[A_{jk}^\alpha, \alpha_k \left(T_{0j}^\beta \right) \right] } - \left[ A_{jk}^\alpha, \left[ \theta_{jk}, T_{0j}^\beta \right] \right]\\
&-\left[ \alpha_k \left(T_{0j}^\beta \right), \left[ \theta_{jk}, T_{0j}^\alpha \right] \right]^2  -\left[ \theta_{jk}, \left[ T_{0j}^\beta, \alpha_j \left(T_{0j}^\alpha\right) \right] \right] +\cancel{2\left[ A_{jk}^\beta,  \alpha_k \left(T_{0j}^\alpha \right) \right]} +\left[ A_{jk}^\beta, \left[ \theta_{jk}, T_{0j}^\alpha\right] \right] - 2\left[ \theta_{jk}, Q_{j\alpha\beta} \right]\\
&- \cancel{2\left[T_{0j}^\alpha, \alpha_k \left( A_{jk}^\beta \right)   \right]} - \left[ T_{0j}^\alpha, \left[ \theta_{jk}, A_{jk}^\beta \right] \right] + \left[T_{0j}^\alpha, \left[  \theta_{jk}, \alpha_k \left(T_{0j}^\beta \right)\right] \right]^2 +\bcancel{2\left[ T_{0j}^\beta, \alpha_k \left(A_{jk}^\alpha \right) \right]} +\left[ T_{0j}^\beta, \left[\theta_{jk},  A_{jk}^\alpha  \right]\right]  - \left[ T_{0j}^\beta, \left[ \theta_{jk}, \alpha_k\left(T_{0j}^\alpha\right)\right] \right]^1\\
&+\xcancel{2\sum_{\gamma=1}^p g_{0j\alpha\beta}^\gamma \alpha_k\left(A_{jk}^\gamma \right)} + \sum_{\gamma=1}^p g_{0j\alpha\beta}^\gamma  \left[ \theta_{jk}, A_{jk}^\gamma \right]  - 2\left[ \theta_{jk}, \alpha_j \left( \left[ T_{0j}^\alpha, T_{0j}^\beta\right] \right)\right]   + \left[ \theta_{jk}, \left[ \theta_{jk}, \left[ T_{0j}^\alpha, T_{0j}^\beta \right] \right]\ \right]\\
&= -\left[ \theta_{jk}, \left[ T_{0j}^\beta, \alpha_k \left(T_{0j}^\alpha\right) \right]\right] + \left[ \theta_{jk}, \left[ T_{0j}^\alpha, \alpha_j \left(T_{0j}^\beta \right)\right] \right]   - \left[ A_{jk}^\alpha, \left[ \theta_{jk}, T_{0j}^\beta \right] \right]\\
&+\left[\theta_{jk}, \left[ T_{0j}^\alpha, \alpha_k \left( T_{0j}^\beta\right)\right]  \right]  - \left[ \theta_{jk}, \left[ T_{0j}^\beta , \alpha_j \left(T_{0j}^\alpha \right)\right] \right] +\left[ A_{jk}^\beta, \left[ \theta_{jk}, T_{0j}^\alpha\right] \right] \\
& - 2\left[ \theta_{jk},  Q_{j\alpha\beta} \right] - \left[T_{0j}^\alpha, \left[ \theta_{jk}, A_{jk}^\beta \right] \right] +\left[ T_{0j}^\beta, \left[\theta_{jk},  A_{jk}^\alpha  \right]\right]- 2\left[ \theta_{jk}, \alpha_j \left( \left[ T_{0j}^\alpha, T_{0j}^\beta\right] \right)\right] +\left[ \theta_{jk}, \left[\theta_{jk}, \left[ T_{0j}^\alpha, T_{0j}^\beta \right] \right]\right] \\
&\equiv_{\Theta_{\mathcal{F}_0}} -2 \left[ \theta_{jk}, \left[ T_{0j}^\beta, \alpha_j\left(T_{0j}^\alpha \right) \right] \right] + \left[ \theta_{jk}, \left[ T_{0j}^\beta,  \left[ \theta_{jk}, T_{0j}^\alpha \right] \right] \right] -\left[ \theta_{jk}, \left[ T_{0j}^\beta, A_{jk}^\alpha \right] \right] - \left[ A_{jk}^\alpha, \left[ \theta_{jk}, T_{0j}^\beta \right] \right]\\
&+ 2\left[ \theta_{jk}, \left[  T_{0j}^\alpha, \alpha_j\left( T_{0j}^\beta \right) \right] \right]  -\left[ \theta_{jk}, \left[ T_{0j}^\alpha, \left[ \theta_{jk}, T_{0j}^\beta \right] \right] \right] + \left[ \theta_{jk}, \left[ T_{0j}^\alpha, A_{jk}^\beta \right] \right]  +\left[ A_{jk}^\beta, \left[ \theta_{jk}, T_{0j}^\alpha\right] \right] \\
& - 2\left[ \theta_{jk},  Q_{j\alpha\beta} \right] - \left[T_{0j}^\alpha, \left[ \theta_{jk}, A_{jk}^\beta \right] \right] +\left[ T_{0j}^\beta, \left[\theta_{jk},  A_{jk}^\alpha  \right]\right]- 2\left[ \theta_{jk}, \alpha_j \left( \left[ T_{0j}^\alpha, T_{0j}^\beta\right] \right)\right] +\left[ \theta_{jk}, \left[\theta_{jk}, \left[ T_{0j}^\alpha, T_{0j}^\beta \right] \right]\right] = 0
\end{align*}}}

We show that $-\delta\left( \left\{ F_{jk} \right\} \right) +D_0'\left( \left\{\left[ \left(\theta_{ij}, A_{ij} \right) , \left( \theta_{jk} , A_{jk} \right) \right] \right\} \right)=0$. First we note that
{\small{\begin{align}\label{o71}
&\alpha_j\left(T_{0i}^\alpha\right)- \alpha_i \left(T_{0i}^\alpha\right)+ \left[ \theta_{ij}, T_i^\alpha \right] = A_{ij}^\alpha,\,\,\,\alpha_k\left(T_{0i}^\alpha\right)- \alpha_i\left(T_{0i}^\alpha\right)+ \left[ \theta_{ik}, T_{0i}^\alpha \right] = A_{ik}^\alpha, \,\,\,\alpha_k\left(T_{0j}^\beta\right)- \alpha_j\left(T_{0j}^\beta\right)+ \left[ \theta_{jk}, T_{0j}^\beta \right] = A_{jk}^\beta\\
&\Longrightarrow \alpha_k\left(T_{0i}^\alpha\right)- \alpha_j\left(T_{0i}^\alpha\right)+ \left[ \theta_{jk}, T_{0i}^\alpha \right] - \sum_{\beta=1}^p \left[ \theta_{jk}, r_{0ij}^{\alpha\beta}\right] T_{0j}^\beta= \sum_{\beta=1}^p r_{0ij}^{\alpha\beta} A_{jk}^\beta \Longrightarrow  A_{ij}^\alpha -  A_{ik}^{\alpha} + \sum_{\beta=1}^p r_{0ij}^{\alpha\beta} A_{jk}^\beta = -\sum_{\beta=1}^p \left[ \theta_{jk}, r_{0ij}^{\alpha\beta} \right] T_{0j}^\beta \notag
\end{align}}}

Then from $(\ref{o71})$, we have
{\small{\begin{align*}
&F_{ij}\left(T_{0i}^\alpha\right)- F_{ik}\left(T_{0i}^\alpha \right)+ \sum_{\beta=1}^q r_{0ij}^{\alpha\beta} F_{jk}\left(T_{0j}^\beta\right)\\
&= - \alpha_i\left(A_{ij}^\alpha\right) - \alpha_j\left(A_{ij}^\alpha\right) +\left[\theta_{ij}, \alpha_i\left(T_{0i}^\alpha\right)  \right] +\left[ \theta_{ij}, \alpha_j\left(T_{0i}^\alpha\right) \right] + \alpha_i \left(A_{ik}^\alpha \right) + \alpha_k \left(A_{ik}^\alpha\right) -\left[ \theta_{ik}, \alpha_i \left(T_{0i}^\alpha\right) \right]  -\left[\theta_{ik}, \alpha_k \left(T_{0i}^\alpha \right) \right]\\
&- \sum_{\beta=1}^q r_{0ij}^{\alpha\beta} \alpha_j \left(A_{jk}^\beta\right) - \sum_{\beta=1}^q r_{0ij}^{\alpha\beta} \alpha_k\left(A_{jk}^\beta\right) + \sum_{\beta=1}^p r_{0ij}^{\alpha\beta}\left[ \theta_{jk}, \alpha_j \left(T_{0j}^\beta\right) \right]  + \sum_{\beta=1}^p r_{0ij}^{\alpha\beta} \left[\theta_{jk}, \alpha_k \left(T_{0j}^\beta\right) \right] \\
&= - \alpha_i \left(A_{ij}^\alpha\right) - \cancel{\alpha_j \left(A_{ij}^\alpha\right)} +\left[\theta_{ij},\alpha_i \left(T_{0i}^\alpha\right)  \right] +\left[ \theta_{ij}, \alpha_j \left(T_{0i}^\alpha\right) \right]^1 + \alpha_i \left(A_{ik}^\alpha\right) + \bcancel{\alpha_k \left(A_{ik}^\alpha \right)} -\left[ \theta_{ik}, \alpha_i \left(T_{0i}^\alpha\right) \right]  -\left[\theta_{ik}, \alpha_k \left(T_{0i}^\alpha\right) \right]^2\\
&+ \cancel{\alpha_j \left(A_{ij}^\alpha\right)} - \alpha_j\left(A_{ik}^\alpha\right) + \cancel{ \sum_{\beta=1}^p \left[\theta_{jk}, r_{0ij}^{\alpha\beta} \right]\alpha_j \left(T_{0j}^\beta\right) } + \alpha_k\left(A_{ij}^\alpha\right) - \bcancel{\alpha_k \left(A_{ik}^\alpha \right) } + \xcancel{\sum_{\beta=1}^p \left[\theta_{jk}, r_{0ij}^{\alpha\beta} \right] \alpha_k \left(T_{0j}^\beta \right)}\\
& + \left[ \theta_{jk}, \alpha_j \left(T_{0i}^\alpha\right) \right]^1 - \cancel{\sum_{\beta=1}^p \left[ \theta_{jk}, r_{0ij}^{\alpha\beta} \right] \alpha_j \left(T_{0j}^\beta\right) } +\left[ \theta_{jk}, \alpha_k \left( T_{0i}^\alpha\right) \right]^2 - \xcancel{\sum_{\beta=1}^p \left[ \theta_{jk}, r_{0ij}^{\alpha\beta} \right] \alpha_k \left(T_{0j}^\beta\right)}\\
&= \left( \alpha_k - \alpha_i   \right)\left(A_{ij}^\alpha\right) - \left(\alpha_j- \alpha_i \right)\left( A_{ik}^\alpha \right) -\left[\theta_{ij}, \alpha_k \left(T_{0i}^\alpha\right)- \alpha_i \left(T_{0i}^\alpha\right)    \right] +\left[ \theta_{ik}, \alpha_j \left(T_{0i}^\alpha\right)- \alpha_i \left(T_{0i}^\alpha \right)  \right] \\
& = \left( \alpha_k - \alpha_i   \right)\left(A_{ij}^\alpha\right) - \left(\alpha_j- \alpha_i \right)\left( A_{ik}^\alpha \right) -\left[\theta_{ij}, \sum_{\beta=1}^p A_{ik}^{\alpha\beta} T_{0i}^\beta - \left[ \theta_{ik} , T_{0i}^\alpha \right]   \right] +\left[ \theta_{ik}, \sum_{\beta=1}^p A_{ij}^{\alpha \beta} T_{0i}^\beta - \left[ \theta_{ij}, T_{0i}^\alpha \right] \right] \\
& =  \sum_{\beta=1}^p A_{ij}^{\alpha\beta} \left( \alpha_k\left(T_{0i}^\beta \right) - \alpha_i \left( T_{0i}^\beta \right) + \left[ \theta_{ik}, T_{0i}^\beta \right] - A_{ik}\left( T_{0i}^\beta \right) + A_{ik}\left(T_{0i}^\beta \right) \right)  + \sum_{\beta=1}^p \left[ \theta_{ik}, A_{ij}^{\alpha\beta} \right] T_{0i}^\beta    - \left[ \theta_{ik}, \left[ \theta_{ij}, T_{0i}^\alpha \right] \right]              \\
&- \sum_{\beta=1}^p A_{ik}^{\alpha\beta} \left( \alpha_j\left( T_{0i}^\beta \right) - \alpha_i\left( T_{0i}^\beta \right) + \left[ \theta_{ij}, T_{0i}^\beta \right] - A_{ij}\left( T_{0i}^\beta \right) + A_{ij}\left( T_{0i}^\beta \right)  \right) - \sum_{\beta=1}^p \left[ \theta_{ij}, A_{ik}^{\alpha\beta} \right] T_{0i}^\beta +\left[ \theta_{ij}, \left[  \theta_{ik}, T_{0i}^\alpha  \right] \right] \\
&= A_{ik}\circ A_{ij}\left(  T_{0i}^\alpha \right) - A_{ij}\circ A_{ik}\left(T_{0i}^\alpha \right) + \sum_{\beta=1}^p \left[ \theta_{ij}+ \theta_{jk}, A_{ij}^{\alpha\beta} \right] T_{0i}^\beta - \sum_{\beta=1}^p \left[ \theta_{ij}, A_{ij}^{\alpha\beta} + A_{jk}^{i\alpha\beta} \right] T_{0i}^\beta + \left[ \left[  \theta_{ij}, \theta_{jk} \right] , T_{0i}^\alpha \right] \\
&= \left( A_{ij} + A_{jk}^i\right)\circ A_{ij}\left( T_{0i}^\alpha \right)- A_{ij}\circ \left( A_{ij} + A_{jk}^i \right) \left(T_{0i}^\alpha \right) + \sum_{\beta=1}^p \left[ \theta_{jk}, A_{ij}^{\alpha\beta} \right] T_{0i}^\beta - \sum_{\beta=1}^p \left[ \theta_{ij}, A_{jk}^{i \alpha\beta} \right] T_{0i}^\beta + \left[ \left[  \theta_{ij}, \theta_{jk} \right] , T_{0i}^\alpha \right]\\
&=\left[ \left[ \theta_{ij}, \theta_{jk} \right] , T_{0i}^\alpha \right] - \left(\left[ \theta_{ij}, A_{jk}^{i} \right] - \left[ \theta_{jk},  A_{ij} \right] + \left[ A_{ij},  A_{jk}^i \right] \right) \left(T_{0i}^\alpha \right) = D_0'\left( \left[ \left( \theta_{ij}, A_{ij} \right),  \left( \theta_{jk}, A_{jk} \right)    \right]      \right)
\end{align*}}} 
This completes Lemma \ref{o54}.

\end{proof}
Then from $(\ref{o40})$ and Lemma \ref{o54}, we obtain
\begin{theorem}
Suppose given a compact foliated complex manifold $\left( M, \Theta_{\mathcal{F}_0} \right)$ with $\Theta_{\mathcal{F}_0}$ locally free, and $\left(\left\{\alpha_j \right\}, \left\{\left(\theta_{jk}, A_{jk} \right) \right\} \right)\in \mathbb{H}^1\left( M, \mathcal{E}_{\Theta_{\mathcal{F}_0}}^\bullet \right)$. In order that there may exist a foliated analytic family $\left( \mathcal{M}, \Theta_\mathcal{F}, B, \pi \right)$ in terms of tangent sheaves such that $\pi^{-1}(0)=\left(M, \Theta_{\mathcal{F}_0}  \right)$ and that $\left(\frac{\partial \left(\mathcal{M}, \Theta_\mathcal{F} \right)}{\partial t}\right)|_{t=0}=\left( \left\{ \alpha_j \right\}, \left\{\left( \theta_{jk} , A_{jk} \right) \right\} \right)$, it is necessary that 
{\Small{\begin{align*}
\left(\left\{ R_j\right\}, \left\{ F_{jk}\right\}, \left\{\left[\left( \theta_{jk}, A_{jk} \right) , \left( \theta_{jk}, A_{jk}\right) \right] \right\} \right) \in C^0\left( \mathcal{U}, \mathscr{H}om_{\mathcal{O}_M}\left( \bigwedge^2 \Theta_{\mathcal{F}_0}, \frac{\Theta_M}{\Theta_{\mathcal{F}_0}} \right)\right) \bigoplus C^1\left( \mathcal{U}, \mathscr{H}om_{\mathcal{O}_M}\left( \Theta_{\mathcal{F}_0}, \Theta_M \right) \right) \bigoplus C^2\left( \mathcal{U} , \mathcal{E}_{\Theta_{\mathcal{F}_0}} \right)
\end{align*}}}
from $(\ref{it1})$ in \textnormal{Lemma \ref{o54}} defines the $0$ class in $\mathbb{H}^2\left( M, \mathcal{E}_{\Theta_{\mathcal{F}_0}}^\bullet  \right)$.
\end{theorem}

Given $\left(\left\{\alpha_j \right\}, \left\{\left(\theta_{jk}, A_{jk} \right) \right\} \right)\in \mathbb{H}^1\left( M, \mathcal{E}_{\Theta_{\mathcal{F}_0}}^\bullet \right)$, if the $2$-cocycle $\left(\left\{ R_j\right\}, \left\{ F_{jk}\right\}, \left\{\left[\left( \theta_{jk}, A_{jk} \right) , \left( \theta_{jk}, A_{jk}\right) \right] \right\} \right) $ is not $0$ in $\mathbb{H}^2\left( M, \mathcal{E}_{\Theta_{\mathcal{F}_0}}^\bullet \right)$, there is no foliated deformation $\left( M_t, \Theta_{\mathcal{F}_t} \right)$ in terms of tangent sheaves with $\left( M_0, \Theta_{\mathcal{F}_0} \right)=\left( M, \Theta_{\mathcal{F}_0} \right)$ and $\left( \frac{\partial\left( \mathcal{M},\Theta_\mathcal{F} \right)}{\partial t} \right)|_{t=0}=\left( \left\{ \alpha_j \right\}, \left\{ \left( \theta_{jk}, A_{jk} \right) \right\} \right)$.

\section{Infinitesimal foliated deformations and obstructions in terms of tangent sheave by using Dolbeault resolution of $\mathcal{E}_{\Theta_{\mathcal{F}_0}}^\bullet$}

In a complex analytic family $\left(\mathcal{M}, B, \pi \right)$ of deformations of a compact complex manifold $M=\pi^{-1}(0)$, the deformations near $M$ are represented by $C^\infty$ vector $(0,1)$-forms $\varphi(t)\in A^{0,1}\left( M, \Theta_M \right)$ on $M$ satisfying $\varphi(0)=0$ and the integrability condition $\bar{\partial} \varphi(t)- \frac{1}{2}\left[\varphi(t), \varphi(t) \right]=0$ where $t\in \Delta$ a sufficiently small polydisk in $B$ (see \cite{Kod05} section \S 5.3). Similarly we will describe the integrability condition in the context of foliated deformations in terms of tangent sheaves.

Let  $\left( \mathcal{M}, \Theta_\mathcal{F} , B, \pi \right)$ with $\Theta_\mathcal{F}$ locally free be a foliated analytic family  in terms of tangent sheaves of deformations of a compact foliated complex manifold $\left(M, \Theta_{\mathcal{F}_0} \right)=\pi^{-1}(0)$ with $\Theta_{\mathcal{F}_0}$ locally free. Assume that given an open covering $\{U_j\}$ of $M$ by coordinate neighborhoods, $\Gamma\left(U_j, \Theta_{\mathcal{F}_0} \right)$ is generated by $T_{0j}^1,..., T_{0j}^p \in \Gamma\left( U_j, \Theta_M \right)$ with $T_{0j}^\alpha=\sum_{\beta=1}^p r_{0jk}^{\alpha\beta}  T_{0k}^\beta$ on $U_j\cap U_k$ for some $r_{0jk}^{\alpha\beta}\in \Gamma\left(U_{jk}, \mathcal{O}_M \right),\alpha=1,...,p$ and $\left[T_{0j}^\alpha, T_{0j}^\beta \right]= \sum_{\gamma=1}^p g_{0j\alpha\beta}^\gamma T_{0j}^\gamma$ for some $g_{0j\alpha\beta}^\gamma \in \Gamma\left( U_j, \mathcal{O}_M \right),\alpha,\beta=1,...,p$. In this subsection, we show that in the foliated analytic family $\left( \mathcal{M}, \Theta_\mathcal{F} , B, \pi \right)$, the deformations near $\left( M, \Theta_{\mathcal{F}_0} \right)$ are represented by a section $\left\{\left(\varphi(t), U_j(t) \right)\right\}\in A^{0,1}\left( M,\mathcal{E}_{\Theta_{\mathcal{F}_0}} \right)$, where $\varphi(t)\in A^{0,1}\left(M, \Theta_M \right)$ and  $U_j(t)\in \Gamma\left( U_j, \mathcal{A}^{0,1}\left( \mathscr{H}om_{\mathcal{O}_M}\left( \Theta_{\mathcal{F}_0}, \Theta_{\mathcal{F}_0} \right)\right) \right)$ which is  defined by $U_j(t)\left(T_{0j}^\alpha \right)= \sum_{\beta=1}^p U_j^{\alpha\beta}(t) T_{0j}^\beta$ for some $U_j^{\alpha\beta}(t)\in \Gamma\left(U_j, \mathcal{A}^{0,1} \right),\alpha=1...,p$, and $X(t)\in  A^{0,0}\left( M, \mathscr{H}om_{\mathcal{O}_M}\left( \Theta_{\mathcal{F}_0}, \Theta_M  \right) \right)$ satisfying $\varphi(0)=0, U_j(0)=0$, and $X(0)=\textnormal{Id}$ and the integrability condition: by setting $X_j(t)\left(T_{0j}^\alpha \right)= X_j^\alpha (t) \in \Gamma\left( U_j, \mathcal{A}^{0,0}\left(\Theta_M \right) \right),\alpha=1,...,p$, 
\begin{align*}
&\bar{\partial} \varphi(t) - \frac{1}{2}\left[\varphi(t), \varphi (t)\right]=0 \\
&\left[ X_j^\alpha (t), X_j^\beta(t)\right]= \sum_{\gamma=1}^p K_{j\alpha\beta}^\gamma (t)  X_j^\gamma (t) \,\,\,\,\,\,\textnormal{for some}\,\,\,\,\, K_{j\alpha\beta}^\gamma(t)\in \Gamma\left( U_j, \mathcal{A}^{0,0} \right) ,\alpha,\beta=1,...,p\\
&\bar{\partial} X_j^\alpha (t) - \left[ \varphi (t) , X_j^\alpha (t) \right]+ \sum_{\eta=1}^p U_j^{\alpha \eta} (t) X_j^\eta (t) =0,\,\,\,\,\, \alpha=1,...,p\\
&\bar{\partial} U_j^{\alpha \delta}(t)- \left[ \varphi(t), U_j^{\alpha \delta}(t) \right]  + \sum_{\sigma=1}^p U_j^{\alpha \sigma} (t) \wedge U_j^{\sigma \delta}(t) =0,\,\,\,\,\,\alpha, \delta=1,...,p
\end{align*}
where $t\in \Delta$ a sufficiently small polydisk in $B$. To deduce the integrability condition, we extend Kodaira's approach (\cite{Kod05} section \S 5.3) in the context of a foliated analytic family in terms of tangent sheaves.

\subsection{Integrability condition} \label{oi12}\

Let $\left(\mathcal{M}, \Theta_\mathcal{F}, B, \pi \right)$ with $\Theta_\mathcal{F}$ locally free be a foliated analytic family in terms of tangent sheaves, where $B$ is a domain of $\mathbb{C}^m$ containing the origin $0$ as in Definition \ref{t3}. Define $|t|=\max_\lambda |t_\lambda|$ for $t=(t_1,...,t_m)\in \mathbb{C}^m$, and let $\Delta=\{t\in \mathbb{C}^m | |t|<r\}$ the polydisk  of radius $r>0$. If we take a sufficiently  small $\Delta \subset B$, then $\left(\mathcal{M}_\Delta, \Theta_{\mathcal{F}_\Delta} \right)=\pi^{-1}\left(\Delta\right)$ is represented in the form
\begin{align*}
\left(\mathcal{M}_\Delta, \Theta_{\mathcal{F}_\Delta}\right)=\bigcup_j \left(U_j\times \Delta, \Theta_\mathcal{F}|_{U_j\times \Delta}\right)
\end{align*}
We denote a point of $U_j$ by $\xi_j=\left(\xi_j^1,...,\xi_j^n\right)$ and $\Gamma\left(U_j\times \Delta, \Theta_\mathcal{F}\right)$ is generated by $T_j^1\left(\xi_j,t\right),...,T_j^p\left(\xi_j,t\right)$, where
\begin{align*}
T_j^\alpha\left(\xi_j,t\right)=\sum_{\beta=1}^n T_j^{\alpha\beta}\left(\xi_j, t\right) \frac{\partial}{\partial \xi_j^\beta}
\end{align*}
for some $T_j^{\alpha\beta}(z_j,t)\in \Gamma\left( U_j\times \Delta , \mathcal{O}_\mathcal{M} \right)$, and satisfies
\begin{align}\label{tai5}
\left[T_j^\alpha(z_j,t), T_j^\beta(z_j,t)\right]=\sum_{\gamma=1}^n g_{j\alpha\beta}^\gamma(z_j,t) T_j^\gamma(z_j,t)
\end{align}
for some $g_{j\alpha\beta}^\gamma(z_j,t)\in \Gamma\left(U_j\times \Delta,  \mathcal{O}_\mathcal{M}\right)$ with $g_{j\alpha\beta}^\gamma(z_j,t)=-g_{j\beta\alpha}^\gamma(z_j,t)$.

For simplicity, we assume that $U_j=\left\{\xi_j\in \mathbb{C}^m | |\xi_j|<1\right\}$ where $|\xi_j|=\max_a |\xi_j^a|$. $(\xi_j,t)\in U_j\times \Delta$ and $(\xi_k, t)\in U_k\times \Delta$ are the same point on $\mathcal{M}_\Delta$ if $\xi_j^\alpha=f_{jk}^\alpha(\xi_k,t),\alpha=1,...,n$ where $f_{jk}(\xi_k,t)$ is a holomorphic map of $\xi_k^1,..., \xi_k^n, t_1,...,t_m$, defined on $(U_k\times \Delta) \cap (U_j\times \Delta)$, and we have the following relation
\begin{align*}
T_j^\alpha(\xi_j,t)=\sum_{\beta=1}^p r_{jk}^{\alpha\beta}(\xi_k,t) T_k^\beta(\xi_k,t),\,\,\,\,\,\alpha=1,...,p
\end{align*}
for a holomorphic function $r_{jk}^{\alpha\beta}(\xi_k,t)$ on $(U_k\times \Delta) \cap ( U_j\times \Delta)$. More precisely
\begin{align}
\sum_{\alpha,\delta=1}^n T_j^{\alpha\delta}(f_{jk}(\xi_k,t),t)\frac{\partial}{\partial \xi_j^\delta } =\sum_{\beta=1}^p \sum_{\gamma,\delta=1}^n r_{jk}^{\alpha\beta}(\xi_k,t) T_k^{\beta\gamma}(\xi_k,t)\frac{\partial f_{jk}^\delta(\xi_j,t)}{\partial \xi_k^\gamma}\frac{\partial}{\partial \xi_j^\delta},\,\,\,\,\,\alpha,\beta=1,...,p
\end{align}
Equivalently
\begin{align} \label{oi5}
T_j^{\alpha\delta}(f_{jk}(\xi_k,t),t)= \sum_{\beta=1}^p \sum_{\gamma=1}^n r_{jk}^{\alpha\beta}(\xi_k,t) T_k^{\beta\gamma}(\xi_k,t)\frac{\partial f_{jk}^{\delta}(\xi_k,t)}{\partial \xi_k^\gamma}
\end{align}
and on $\left(U_i\times \Delta\right)\cap \left( U_j\times \Delta \right)\cap \left( U_k\times \Delta \right)$, we have 
\begin{align*}
r_{ik}^{\alpha \gamma}(\xi_k,t)= \sum_{\beta=1}^p r_{ij}^{\alpha\beta}\left(f_{jk}(\xi_k,t),t) \right) r_{jk}^{\beta \gamma}(\xi_k,t),\,\,\,\,\,\alpha,\gamma=1,...,p
\end{align*}

By \cite{Kod05} Theorem 2.3, $M_t$ is diffeomorphic to $M_0=\pi^{-1}(0)$ as differentiable manifolds for each $t\in \Delta$. We put $M:=M_0$. By \cite{Kod05} Theorem 2.5, if we take a sufficiently small $\Delta$, there is a diffeomorphism $\Psi$ of $M\times \Delta$ onto $\mathcal{M}_\Delta$ as differentiable manifolds such that $\pi\circ \Psi$ is the projection $M\times \Delta \to \Delta$. Let $z=\left(z^1,..., z^n \right)$ be local coordinates of $M=M_0$. Then we have $\pi\circ \Psi(z,t)=t, t\in \Delta$. For $\Psi(z,t)\in U_j\times \Delta$, put
\begin{align*}
\Psi\left( z, t\right) = \left(\xi_j^1(z,t), ..., \xi_j^n(z,t), t_1,..., t_m \right)
\end{align*}
Then each component $\xi_j^\alpha=\xi_j^\alpha(z,t),\alpha=1,...,n$ is a $C^\infty$ function. If we identify $\mathcal{M}_\Delta =\Psi\left(M\times \Delta \right)$ with $M\times \Delta$, and $\Theta_{\mathcal{F}}$ with $\Psi^* \Theta_\mathcal{F}$ via $\Psi$, $\left(\mathcal{M}_\Delta, \Theta_{\mathcal{F}_\Delta} \right)$ is considered as a complex manifold with the complex structure defined on the $C^\infty$ manifold $M\times \Delta$ by the system of local coordinates on $U_j\times \Delta$
\begin{align*}
\left\{ \left(\xi_j,t \right)| j=1 ,2, 3,...   \right\},\,\,\,\,\,(\xi_j,t)=\left( \xi_j^1(z,t), ..., \xi_j^n(z,t), t_1,..., t_m \right)
\end{align*}
and the locally free sheaf $\Psi^* \Theta_{\mathcal{F}_\Delta}$ on $M\times \Delta$, more precisely, the $C^\infty$ differentiable vector bundle on $M\times \Delta$.

We note that since $(z^1,..., z^n)$ and $\left(\xi_j^1(z,0), ..., \xi_j^n(z,0)\right)$ are local complex coordinates on $M=M_0$, $\xi_j^\alpha(z,0)$ are holomorphic functions of $z^1,...,z^n,\alpha=1,...,n$. We also note that if we take $\Delta$ sufficiently small, we have
\begin{align}\label{oi3}
\det\left( \frac{\partial \xi_j^\alpha(z,t)}{\partial z^\lambda} \right)_{\alpha, \lambda=1,...,n} \ne 0\,\,\,\,\,\, \, \textnormal{for}\,\,\,t\in \Delta
\end{align}
On the other hand, each locally free sheaf $\Psi^* \Theta_{\mathcal{F}_t}$ in a differentiable  family $\left\{ \Psi^*\Theta_{\mathcal{F}_t} |t\in \Delta \right\}$ is the restriction to $M\times t$ of the locally free sheaf $\Psi^* \Theta_\mathcal{F}$ over $M\times \Delta , \Psi^* \Theta_{\mathcal{F}_t}= \Psi^* \Theta_{\mathcal{F}}|_{M\times t}$. Then by \cite{Kod05} Lemma 7.1 the locally free sheaf $\Psi^* \Theta_\mathcal{F}$ over $M\times \Delta$ as a $C^\infty$ vector bundle is $C^\infty$ equivalent to the trivial vector bundle $\Theta_{\mathcal{F}_0}\times \Delta $. In other words, there is a bundle diffeomorphism $\sigma: \Theta_{\mathcal{F}_0}\times \Delta \to \Psi^* \Theta_{\mathcal{F}}$ over $M\times \Delta$ when we consider $\Theta_{\mathcal{F}_0}$ and $\Psi^* \Theta_{\mathcal{F}}$ as vector bundles. We note that $\Theta_\mathcal{F}$ as a vector bundle is described by the transition matrices $\left( r_{jk}(\xi_j,t)^T\right)^{-1}$ with respect to the covering $\left\{U_j\times \Delta \right\}$ of $\mathcal{M}_\Delta$, where $r_{jk}(\xi_j,t):=\left( r_{jk}^{\alpha\beta} (\xi_k,t)\right)_{\alpha, \beta=1,...,p}$. Then $\Psi^* \Theta_{\mathcal{F}}$ is described by the transition matrices $\left( r_{jk}\left(\xi_k(z,t) , t\right)^T \right)^{-1}$ with respect to the covering $\left\{\mathscr{U}_j:=\Psi^{-1}\left( U_j\times \Delta \right) \right\}$ of $M\times \Delta$, where $r_{jk}\left(\xi_k\left( z, t\right), t \right):= \left( r_{jk}^{\alpha\beta}\left(\xi_k(z,t),t \right)\right)_{\alpha, \beta=1,...,p}$. In particular $\Theta_{\mathcal{F}_0}$ is described by $\left(r_{jk}(\xi_j(z,0),0 )^T \right)^{-1}$. Then there are invertible $C^\infty$ matrix functions $\sigma_j(z,t):\mathbb{C}^p \to \mathbb{C}^p$ such that $\sigma$ is represented locally by
\begin{align*}
\sigma:\mathbb{C}^p\times\mathscr{U}_j &\to \mathbb{C}^p \times \mathscr{U}_j\\
                         \left( X=(x_1,..., x_p), z, t \right) &\mapsto \left( \sigma_j(z,t)X, z, t\right)
\end{align*}
such that
\begin{align*}
\sigma_j(z,t) \cdot \left( r_{jk}\left(\xi_k(z,0),0 \right)^T \right)^{-1} = \left( r_{jk}\left(\xi_k(z,t), t \right)^T \right)^{-1} \cdot \sigma_k(z,t)
\end{align*}

Then we have
\begin{align*}
\left(\sigma_j(z,t)^{-1} \right)^T\cdot r_{jk}\left(\xi_k(z,0),0 \right) = r_{jk}\left(\xi_k( z,t),t \right) \cdot \left( \sigma_k(z,t)^{-1} \right)^T
\end{align*}
Let $\left(\sigma_j(z,t)^{-1} \right)^T $ be represented by invertible matrix $c_j(z,t):=\left(c_j^{\alpha\beta}(z,t) \right)$ where $c_j^{\alpha\beta}(z,t)$ is a $C^\infty$ function, and set $r_{0jk}^{\alpha\beta}(z):=r_{jk}^{\alpha\beta}(\xi_j(z,0),0)$ and $r_{0jk}(z):=\left(r_{0jk}^{\alpha\beta} (z)\right)_{\alpha, \beta=1,...,p}$, so that we have
\begin{align}\label{oi14}
c_j(z,t)\cdot r_{0jk}(z)= r_{jk}\left(\xi_k(z,t),t) \right) \cdot c_k(z,t) 
\end{align}
Then we have
\begin{align}\label{oi10}
r_{jk}^{\alpha\beta}\left(\xi_k(z,t),t\right)= \sum_{\delta, \eta=1}^p c_{j}^{\alpha \delta}(z,t) r_{0jk}^{\delta \eta}(z) \left(c_k^{-1}\right)^{\eta\beta}(z,t),\,\,\,\,\,\,\,\alpha,\beta=1,...,p
\end{align}
where $c_j^{-1}(z,t)=\left( \left(  c_j^{-1}\right)^{\alpha\beta}(z,t) \right)_{\alpha, \beta=1,.., p} $ is the inverse of $c_j(z,t)$.

\subsection{Identification of deformations of complex structures with $\varphi(t)\in A^{0,1}\left(M, \Theta_M \right)$} \label{oi51}\

Recall that we put $\mathscr{U}_j=\Psi^{-1}\left(U_j\times \Delta \right)$. Then $\mathscr{U}_j \subset M\times \Delta$ is the domain of $\xi_j^{\alpha}(z,t)$. From $(\ref{oi3})$, we can define a $(0,1)$-form $\varphi_j^\lambda(z,t)= \sum_{v=1}^n \varphi_{jv}^\lambda(z,t) d\bar{z}^v$ for each $\lambda=1,...,n $ such that $\bar{\partial}\xi_j^\alpha = \sum_{\lambda=1}^n \varphi_{jv}^\lambda(z,t)\frac{\partial \xi_j^\alpha}{\partial z_\lambda}d\bar{z}^v$. Then on $\mathscr{U}_j \cap \mathscr{U}_k$, we have $\sum_{\lambda=1}^n \varphi_j^\lambda(z,t)\frac{\partial}{\partial z^\lambda} = \sum_{\lambda=1}^n \varphi_k^\lambda(z,t) \frac{\partial}{\partial z^\lambda}$ (see \cite{Kod05} p.262). Then if for $(z,t)\in \mathscr{U}_j$, we define $\varphi(z,t):=\sum_{\lambda=1}^n \varphi_j^\lambda(z,t)\frac{\partial}{\partial z^\lambda}$, then $\varphi(t):=\varphi(z,t)\in A^{0,1}\left(M, \Theta_M \right)$ for $t\in \Delta$, and satisfies
\begin{align}\label{oi23}
\varphi(0)=0,\,\,\,\,\,\,\,\bar{\partial}\varphi(t)- \frac{1}{2}\left[ \varphi(t), \varphi(t) \right]=0
\end{align}
(see \cite{Kod05} p.263, p.265), We also point out that
\begin{theorem}\label{oi13}
If we take a sufficiently small polydisk $\Delta$ as in \textnormal{subsection \ref{oi12}}, then for $t\in \Delta$, a local $C^\infty$ function $f$ on $M$ is holomorphic with respect to the complex structure $M_t$ if and only if
\begin{align*}
\left( \bar{\partial} - \bar{\partial} \varphi(t)\right)f=0
\end{align*} 
\end{theorem}
\begin{proof}
See \cite{Kod05} Theorem 5.3 p.263
\end{proof}

\subsection{Identification of deformations of complex structures and locally free sheaves with $\left\{\left( \varphi(t), U_j(t) \right)\right\}  \in A^{0,1}\left( M, \mathcal{E}_{\Theta_{\mathcal{F}_0}} \right)$} \label{oi52} \

By setting $R_{jk}^{\alpha\beta}(z,t)=r_{jk}^{\alpha\beta}\left(\xi_k(z,t),t \right)$ and $R_{jk}:=\left( r_{jk}^{\alpha\beta}\left( \xi_j(z,t),t \right) \right)_{\alpha,\beta=1,...p}$, we have from $(\ref{oi14})$,
\begin{align*}
c_j(z,t) \cdot r_{0jk}(z) = R_{jk}(z,t)\cdot c_k(z,t),
\end{align*}
so that we can write
\begin{align}\label{oi15}
\sum_{\beta=1}^p R_{jk}^{\alpha\beta} c_k^{\beta \gamma}= \sum_{\beta=1}^p c_j^{\alpha\beta} r_{0jk}^{\beta \gamma} ,\,\,\,\,\,\,\alpha,\gamma=1,...,p
\end{align}
By taking $\bar{\partial}$ on $(\ref{oi15})$ , we have
\begin{align}\label{oi16}
\sum_{\beta=1}^p \bar{\partial } R_{jk}^{\alpha\beta} \cdot c_k^{\beta \gamma} + \sum_{\beta=1}^p R_{jk}^{\alpha\beta} \cdot \bar{\partial} c_k^{\beta\gamma}= \sum_{\beta=1}^p \bar{\partial} c_j^{\alpha\beta} \cdot r_{0jk}^{\beta \gamma}
\end{align}
From Theorem \ref{oi13} and $(\ref{oi10})$, we have
\begin{align*}
 \bar{\partial} R_{jk}^{\alpha\beta}=\left[\varphi, R_{jk}^{\alpha\beta}\right]&= \sum_{\delta, \eta =1}^p \left[\varphi, c_j^{\alpha\delta}\cdot r_{0jk}^{\delta \eta}  \cdot \left(c_k^{-1}\right)^{\eta \beta} \right] \\
\Longrightarrow \sum_{\beta=1}^p \bar{\partial } R_{jk}^{\alpha\beta} \cdot c_k^{\beta \gamma} &= \sum_{\delta, \eta, \beta=1}^p \left[\varphi, c_j^{\alpha\delta}\cdot r_{0jk}^{\delta \eta}  \cdot \left(c_k^{-1}\right)^{\eta \beta} \right] c_k^{\beta \gamma} = \sum_{\delta=1}^p \left[ \varphi, c_j^{\alpha\delta} r_{0jk}^{\delta \gamma} \right] - \sum_{\delta, \eta, \beta=1}^p c_j^{\alpha \delta} r_{0jk}^{\delta \eta} \left(c_k^{-1}\right)^{\eta \beta} \left[ \varphi, c_k^{\beta \gamma}\right] \\
&= \sum_{\delta=1}^p \left[ \varphi, c_j^{\alpha \delta} \right] r_{0jk}^{\delta \gamma} + \sum_{\delta=1}^p c_j^{\alpha \delta} \left[\varphi, r_{0jk}^{\delta \gamma} \right]-   \sum_{\delta, \eta, \beta=1}^p c_j^{\alpha \delta} r_{0jk}^{\delta \eta} \left(c_k^{-1}\right)^{\eta \beta} \left[ \varphi, c_k^{\beta \gamma}\right]
\end{align*}
Then from $(\ref{oi16})$ and $(\ref{oi10})$, $ \sum_{\beta=1}^p \bar{\partial } R_{jk}^{\alpha\beta} \cdot c_k^{\beta \gamma} + \sum_{\beta=1}^p R_{jk}^{\alpha\beta} \cdot \bar{\partial} c_k^{\beta\gamma}-\sum_{\beta=1}^p \bar{\partial} c_j^{\alpha\beta} \cdot r_{0jk}^{\beta \gamma} =0$ implies that
{\small{\begin{align}\label{oi17}
&\sum_{\delta=1}^p \left[ \varphi, c_j^{\alpha \delta} \right] r_{0jk}^{\delta \gamma} + \sum_{\delta=1}^p c_j^{\alpha \delta} \left[\varphi, r_{0jk}^{\delta \gamma} \right]-   \sum_{\delta, \eta, \beta=1}^p c_j^{\alpha \delta} r_{0jk}^{\delta \eta} \left(c_k^{-1}\right)^{\eta \beta} \left[ \varphi, c_k^{\beta \gamma}\right] + \sum_{\delta, \eta, \beta=1}^p c_j^{\alpha \delta} \cdot r_{0jk}^{\delta \eta}\cdot \left(c_k^{-1}\right)^{\eta \beta} \cdot \bar{\partial} c_k^{\beta \gamma} -\sum_{\beta=1}^p \bar{\partial} c_j^{\alpha\beta} \cdot r_{0jk}^{\beta \gamma}=0
\end{align}}}
By multiplying $\sum_{\alpha=1}^p \left(c_j^{-1} \right)^{\sigma \alpha}$ on the left and multiplying $\sum_{\xi=1}^p \left(r_{0jk}^{-1}\right)^{\gamma \xi}$ on $(\ref{oi17})$, we have
\begin{align}
&\sum_{\alpha =1}^p \left(c_j^{-1}\right)^{\sigma \alpha}\left[ \varphi, c_j^{\alpha \xi} \right] + \sum_{\gamma=1}^p \left[\varphi, r_{0jk}^{\sigma \gamma} \right] \left(r_{0jk}^{-1}\right)^{\gamma \xi} -   \sum_{ \eta, \beta,\gamma=1}^p  r_{0jk}^{\sigma \eta} \left(c_k^{-1}\right)^{\eta \beta} \left[ \varphi, c_k^{\beta \gamma}\right] \left(r_{0jk}^{-1}\right)^{\gamma \xi}\\
&+ \sum_{ \eta, \beta,\gamma=1}^p r_{0jk}^{\sigma \eta}\cdot \left(c_k^{-1}\right)^{\eta \beta} \bar{\partial} c_k^{\beta \gamma}\cdot \left(r_{0jk}^{-1}\right)^{\gamma\xi} -\sum_{\sigma =1}^p \left(c_j^{-1}\right)^{\sigma \alpha}\cdot \bar{\partial} c_j^{\alpha\xi}=0\notag
\end{align}
which is equivalent to
\begin{align}\label{oi18}
\sum_{\sigma=1}^p \left(c_j^{-1}\right)^{\sigma \alpha}\left( \bar{\partial} c_j^{\alpha \xi} - \left[ \varphi, c_j^{\alpha \xi}\right] \right) =\sum_{\eta=1}^p r_{0jk}^{\sigma \eta} \left(  \sum_{\beta,\gamma=1}^p \left(c_k^{-1}\right)^{\eta \beta} \left( \bar{\partial} g_k^{\beta \gamma}- \left[ \varphi, c_k^{\beta \gamma}\right]    \right) \right) \left(r_{0jk}^{-1}\right)^{\gamma \xi} + \sum_{\gamma=1}^p \left[ \varphi, r_{0jk}^{\sigma \gamma} \right] \left(r_{0jk}^{-1}\right)^{\gamma \xi}
\end{align}

By taking $U_j^{\alpha \delta}(t):= \sum_{\eta=1}^p \left( c_j^{-1}\right)^{\alpha \eta}\left( \bar{\partial} c_j^{\eta \delta} - \left[ \varphi, c_j^{\eta \delta}\right] \right)$, let us define an element $U_j (t)\in \Gamma\left(U_j, \mathcal{A}^{0,1}\left(\mathscr{H}om_{\mathcal{O}_M}\left( \Theta_{\mathcal{F}_0}, \Theta_{\mathcal{F}_0} \right) \right)    \right)$ by
\begin{align*}
U_j(t):\Gamma\left( U_j, \Theta_{\mathcal{F}_0} \right)&\to \Gamma\left( U_j, \mathcal{A}^{0,0}\left( \Theta_{\mathcal{F}_0} \right) \right) \\
  T_{0j}^\alpha &\mapsto \sum_{\delta=1}^p U_j^{\alpha \delta} T_{0j}^\delta:=\sum_{\eta, \delta=1}^p \left(c_j^{-1} \right)^{\alpha \eta} \left(\bar{\partial} c_j^{\eta \delta}- \left[\varphi, c_j^{\eta \delta} \right] \right) T_{0j}^\delta
\end{align*}
and linearly extends to $\Gamma\left(U_j, \Theta_{\mathcal{F}_0} \right)$. Then from $(\ref{oi18})$, we see that
\begin{align}\label{oi40}
U_j\left(T_{0j}^\alpha \right) = \sum_{\beta=1}^p r_{0jk}^{\alpha\beta} \left( \sum_{\gamma, \eta=1}^p \left( c_k^{-1} \right)^{\beta \gamma} \left( \bar{\partial}  c_k^{\gamma \eta} - \left[ \varphi, c_k^{\gamma \eta}\right]\right) T_{0k}^\eta \right) + \sum_{\gamma=1}^p \left[ \varphi, r_{0jk}^{\alpha \beta} \right] T_{0k}^\beta = U_k\left( T_{0j}^\alpha \right) + \sum_{\beta=1}^p \left[ \varphi, r_{0jk}^{\alpha\beta} \right] T_{0k}^\beta
\end{align}
This implies that
\begin{align}
\left\{ \left( \varphi(t), U_j(t)  \right) \right\} \in A^{0,1}\left( M, \mathcal{E}_{\Theta_{\mathcal{F}_0}} \right)
\end{align}
and since $c_j(0)=\textnormal{Id}$, we have $U_j(0)=0$. On the other hand, we note that
\begin{align}\label{oi19}
\sum_{\beta=1}^p \left(c_j^{-1}\right)^{\alpha \beta} c_j^{\beta \eta} = \delta_\alpha^\eta \Longrightarrow \sum_{\beta=1}^p \bar{\partial} \left(c_j^{-1}\right)^{\alpha\beta}\cdot c_j^{\beta \eta} + \sum_{\beta=1}^p \left(c_j^{-1}\right)^{\alpha\beta} \cdot \bar{\partial} c_j^{\beta \eta}=0 \Longrightarrow \bar{\partial} \left(c_j^{-1}\right)^{\alpha\delta}= - \sum_{\beta, \eta=1}^p \left(c_j^{-1}\right)^{\alpha\beta} \cdot \bar{\partial} c_j^{\beta \eta}\cdot \left(c_j^{-1}\right)^{\eta \delta}
\end{align}
Then from (\ref{oi19}), we have
\begin{align*}
&\bar{\partial}\left(\sum_{\eta=1}^p \left(c_j^{-1}\right)^{\alpha\eta} \left( \bar{\partial} c_j^{\eta \delta} - \left[ \varphi, c_j^{\eta \delta}\right]   \right)  \right)= \sum_{\eta=1}^p \bar{\partial} \left(\left(c_j^{-1}\right)^{\alpha\eta}\right)\wedge   \left( \bar{\partial} c_j^{\eta \delta} - \left[ \varphi, c_j^{\eta \delta}\right]   \right)   +\sum_{\eta=1}^p \left(c_j^{-1}\right)^{\alpha\eta}\cdot \left( - \left[\bar{\partial} \varphi, c_j^{\eta \delta} \right]  + \left[\varphi, \bar{\partial} c_j^{\eta \delta} \right]     \right)\\
&=-\sum_{\eta, \xi, \sigma=1}^p \left(c_j^{-1}\right)^{\alpha  \xi} \cdot \bar{\partial } c_j^{\xi \sigma} \cdot \left(c_j^{-1}\right)^{\sigma \eta} \wedge \left( \bar{\partial } c_j^{\eta \delta} - \left[ \varphi, c_j^{\eta \delta} \right]   \right) +\sum_{\eta=1}^p \left(c_j^{-1}\right)^{\alpha \eta}\left( -\frac{1}{2}\left[\left[ \varphi, \varphi\right], c_j^{\eta \delta} \right]+ \left[ \varphi,  \bar{\partial}  c_j^{\eta \delta}  \right]    \right) \\
&=-\sum_{\eta, \xi, \sigma=1}^p \left(c_j^{-1}\right)^{\alpha  \xi} \cdot \bar{\partial }  c_j^{\xi \sigma} \cdot \left(c_j^{-1}\right)^{\sigma \eta} \wedge \left( \bar{\partial } c_j^{\eta \delta} - \left[ \varphi, c_j^{\eta \delta} \right]   \right) +\sum_{\eta=1}^p \left(c_j^{-1}\right)^{\alpha \eta}\left( -\left[ \varphi, \left[ \varphi , c_j^{\eta \delta} \right] \right]+ \left[ \varphi,  \bar{\partial}  c_j^{\eta \delta}  \right]    \right)\\
&= -\sum_{ \xi, \sigma=1}^p \left(c_j^{-1}\right)^{\alpha  \xi} \cdot \bar{\partial } c_j^{\xi \sigma} \wedge \left(\sum_{\eta=1}^p \left( c_j^{-1} \right)^{\sigma \eta}  \left( \bar{\partial } c_j^{\eta \delta} - \left[ \varphi, c_j^{\eta \delta} \right]   \right)\right) +\sum_{\eta=1}^p \left(c_j^{-1}\right)^{\alpha \eta}\left( \left[ \varphi, \bar{\partial} c_j^{\eta \delta} - \left[ \varphi , c_j^{\eta \delta} \right] \right]   \right)\\
&= -\sum_{ \xi, \sigma=1}^p \left(c_j^{-1}\right)^{\alpha  \xi} \cdot \bar{\partial } c_j^{\xi \sigma} \wedge \left(\sum_{\eta=1}^p \left(c_j^{-1}\right)^{\sigma \eta}  \left( \bar{\partial } c_j^{\eta \delta} - \left[ \varphi, c_j^{\eta \delta} \right]   \right)\right) + \left[ \varphi, \sum_{\eta=1}^p \left(c_j^{-1}\right)^{\alpha \eta}\left(  \bar{\partial} c_j^{\eta \delta} - \left[ \varphi , c_j^{\eta \delta} \right]    \right) \right]\\
& + \sum_{\xi, \sigma =1}^p \left(c_j^{-1}\right)^{\alpha \xi}\left[ \varphi, c_j^{\xi \sigma}\right] \wedge  \left(\sum_{\eta=1}^p \left(c_j^{-1}\right)^{\sigma \eta}\left(\bar{\partial} c_j^{\eta \delta}- \left[ \varphi , c_j^{\eta \delta}\right] \right)\right)\\
&= - \sum_{ \sigma=1}^p \left( \sum_{\xi=1}^q \left(c_j^{-1}\right)^{\alpha  \xi} \left( \bar{\partial } c_j^{\xi \sigma} - \left[ \varphi, c_j^{\xi\sigma} \right] \right) \right) \wedge \left(\sum_{\eta=1}^p \left(c_j^{-1} \right)^{\sigma \eta}  \left( \bar{\partial } c_j^{\eta \delta} - \left[ \varphi, g_j^{\eta \delta} \right]   \right)\right) + \left[ \varphi, \sum_{\eta=1}^p \left(c_j^{-1}\right)^{\alpha \eta}\left(  \bar{\partial} c_j^{\eta \delta} - \left[ \varphi , c_j^{\eta \delta} \right]    \right) \right]
\end{align*}
This implies that
\begin{align}\label{oi24}
\bar{\partial} U_j^{\alpha \eta} - \left[ \varphi, U_j^{\alpha \eta} \right] + \sum_{\sigma=1}^p U_j^{\alpha \sigma} \wedge U_j^{\sigma \eta}=0
\end{align}
\subsection{Identification of deformations of foliated structures $X(t)\in A^{0,0} \left( M, \mathscr{H}om_{\mathcal{O}_M}\left( \Theta_{\mathcal{F}_0}, \Theta_M \right) \right)$} \label{it20} \

From $\left( \ref{oi3} \right)$ there exists $C^\infty$ functions $V_j^{\alpha \gamma}(z,t),\gamma=1,...,p$ such that
\begin{align} \label{oi4}
\sum_{\gamma=1}^p V_j^{\alpha \gamma}(z,t)\frac{\partial \xi_j^\beta(z,t)}{\partial z^\gamma }= T_j^{\alpha\beta}\left(\xi_j(z,t),t \right) \Longrightarrow \sum_{ \beta,\gamma=1}^q V_j^{\alpha\gamma}(z,t)\frac{\partial \xi_j^\beta(z,t)}{\partial z^\gamma }\frac{\partial}{\partial \xi_j^\beta}= \sum_{ \beta=1}^p T_j^{\alpha \beta}(\xi_j(z,t), t)\frac{\partial}{\partial \xi_j^\beta}
\end{align}

Then from $(\ref{oi4})$ and $(\ref{oi5})$ we have
{\small{\begin{align*}
&\sum_{\gamma=1}^n \left( V_j^{\alpha\gamma}(z,t)- \sum_{\beta=1}^p r_{jk}^{\alpha\beta} \left(\xi_k(z,t),t\right) V_k^{\beta\gamma}(z,t) \right)\frac{\partial \xi_j^\gamma(z,t)}{\partial z^\gamma}=  \sum_{\gamma=1}^n V_j^{\alpha\gamma}(z,t)\frac{\partial \xi_j^\delta(z,t) }{\partial z^\gamma}- \sum_{\beta=1}^p\sum_{\gamma=1}^n r_{jk}^{\alpha\beta} \left(\xi_k(z,t),t \right) V_k^{\beta\gamma}(z,t) \frac{\partial f_{jk}^\delta\left(\xi_k(z,t),t \right)}{\partial z^\gamma}\\
&=  T_j^{\alpha \delta}\left(\xi_j(z,t),t \right) - \sum_{\beta=1}^p\sum_{\gamma=1}^n r_{jk}^{\alpha\beta} \left(\xi_k(z,t),t\right) V_k^{\beta\gamma}(z,t) \left( \sum_{\eta=1}^n\frac{\partial f_{jk}^{\delta}}{\partial \xi_k^\eta}\frac{\partial \xi_k^\eta}{\partial z^\gamma} + \sum_{\eta=1}^n \frac{\partial f_{jk}^\delta}{\partial \bar{\xi}_j^\eta} \frac{\partial \bar{\xi}_k^\eta}{\partial z^\gamma}    \right)\\
&= T_j^{\alpha \delta}\left(\xi_j(z,t),t \right)- \sum_{\beta=1}^p \sum_{\eta=1}^n  r_{jk}^{\alpha\beta}(\xi_k(z,t),t) T_k^{\beta \eta}(\xi_k(z,t),t)\frac{\partial f_{jk}^\delta}{\partial \xi_k^\eta} = 0
\end{align*}}}
Then from $(\ref{oi3})$, we have
\begin{align}\label{oi50}
 V_j^{\alpha\gamma}(z,t)- \sum_{\beta=1}^p r_{jk}^{\alpha\beta} \left(\xi_k(z,t),t\right) V_k^{\beta\gamma}(z,t) =0 
\end{align}

We note that
\begin{align}\label{oi6}
\left[\sum_{\gamma=1}^n V_j^{\alpha \gamma} \frac{\partial}{\partial z^\gamma} , \sum_{\delta=1}^n V_j^{\beta \eta} \frac{\partial}{\partial z^\eta}  \right] = \sum_{\gamma, \eta=1}^n \left( V_j^{\alpha\gamma} \frac{\partial V_j^{\beta \eta}}{\partial z^\gamma} - V_j^{\beta \gamma} \frac{\partial V_j^{\alpha \eta}}{\partial z^\gamma} \right) \frac{\partial}{\partial z^\eta}
\end{align}
and by taking the derivative of $(\ref{oi4})$ $\sum_{ \gamma=1}^q V_j^{\alpha\gamma}(z,t)\frac{\partial \xi_j^\beta(z,t)}{\partial z^\gamma }=T_j^{\alpha \beta}\left(\xi_j(z,t), t\right)$ with respect to $z^\eta$, we have
\begin{align}\label{oi7}
\sum_{\gamma=1}^p \frac{\partial V_j^{\alpha \gamma}}{\partial z^\eta} \frac{\partial \xi_j^\beta}{\partial z^\gamma} + \sum_{\gamma=1}^q V_j^{\alpha \gamma} \frac{\partial^2 \xi_j^\beta}{\partial z^\eta \partial z^\gamma} = \sum_{\delta=1}^n \frac{\partial T_j^{\alpha\beta}}{\partial \xi_j^\delta} \frac{\partial \xi_j^\delta}{\partial z^\eta}
\end{align}
and from $(\ref{tai5})$, we have
\begin{align}\label{oi8}
\left[\sum_{\beta=1}^n T_j^{\alpha\gamma} \frac{\partial}{\partial \xi_j^\gamma}, \sum_{\eta=1}^p T_j^{\beta \eta} \frac{\partial}{\partial \xi_j^\eta} \right]= \sum_{\gamma , \eta=1}^n \left(T_j^{\alpha \gamma}\frac{\partial T_j^{\beta \eta}}{\partial \xi_j^\gamma} - T_j^{\beta \gamma}\frac{\partial T_j^{\alpha \eta}}{\partial \xi_j^\gamma} \right)\frac{\partial}{\partial \xi_j^\eta} = \sum_{\gamma, \eta=1}^p g_{j\alpha\beta}^\gamma T_j^{\gamma \eta} \frac{\partial}{\partial \xi_j^\eta}
\end{align}
Then from $(\ref{oi4}), (\ref{oi7})$ and $(\ref{oi8})$, we have
\begin{align*}
&\sum_{\gamma=1}^n T_j^{\alpha \gamma}\frac{\partial T_j^{\beta \eta}}{\partial \xi_j^\gamma} - \sum_{\gamma=1}^n T_j^{\beta \gamma}\frac{\partial T_j^{\alpha \eta}}{\partial \xi_j^\gamma}= \sum_{\gamma=1}^p g_{j\alpha\beta}^\gamma T_j^{\gamma \eta}\iff \sum_{\gamma, \delta=1}^n V_j^{\alpha \delta} \frac{\partial \xi_j^\gamma}{\partial z^\delta} \frac{\partial T_j^{\beta \eta}}{\partial \xi_j^\gamma} - \sum_{\gamma,\delta=1}^n V_j^{\beta \delta} \frac{\partial \xi_j^\gamma}{\partial z^\delta} \frac{\partial T_j^{\alpha \eta}}{\partial \xi_j^\gamma} = \sum_{\gamma=1}^p \sum_{\delta=1}^n g_{j\alpha\beta}^\gamma V_j^{\gamma \delta}\frac{\partial \xi_j^\eta}{\partial z^\delta}\\
&\iff \sum_{\delta=1}^n V_j^{\alpha \delta}\left(   \sum_{\gamma=1}^n \frac{\partial V_j^{\alpha \gamma}}{\partial z^\delta}\frac{\partial \xi_j^\eta}{\partial z^\gamma} + \sum_{\gamma=1}^n V_j^{\beta \gamma} \frac{\partial^2 \xi_j^\eta}{\partial z^\delta \partial z^\gamma}     \right) - \sum_{\delta=1}^n V_j^{\beta \delta}\left(\sum_{\gamma=1}^p \frac{\partial V_j^{\alpha \gamma}}{\partial z^\delta} \frac{\partial \xi_j^\eta}{\partial z^\gamma} + \sum_{\gamma=1}^p V_j^{\alpha \gamma} \frac{\partial^2 \xi_j^\eta}{\partial z^\delta z^\gamma} \right)= \sum_{\gamma=1}^p \sum_{\delta=1}^n g_{j\alpha\beta}^\gamma V_j^{\gamma \delta}\frac{\partial \xi_j^\eta}{\partial z^\delta}\\
&\Longrightarrow \sum_{\delta,\gamma=1}^n \left( V_j^{\alpha \delta}\frac{\partial V_j^{\alpha \gamma}}{\partial z^\delta} - V_j^{\beta \delta} \frac{\partial V_j^{\alpha \gamma}}{\partial z^\delta} \right)\frac{\partial \xi_j^\eta}{\partial z^\gamma} = \sum_{\gamma=1}^p \sum_{\delta=1}^n g_{j\alpha\beta}^\gamma V_j^{\gamma \delta}\frac{\partial \xi_j^\eta}{\partial z^\delta}
\end{align*}
Then from $(\ref{oi3})$ and $(\ref{oi6})$ and by setting $V_j^\alpha(z,t)= \sum_{\delta=1}^n V_j^{\alpha \delta}(z,t) \frac{\partial}{\partial z^\delta}$ and $G_{j\alpha\beta}^\gamma(z,t):=g_{j\alpha\beta}^\gamma\left( \xi_j(z,t),t \right)$, we have
\begin{align}\label{oi20}
\left[ V_j^\alpha(z,t), V_j^\beta(z,t) \right]= \sum_{\gamma=1}^p g_{j\alpha\beta}^\gamma\left( \xi_j(z,t),t  \right) V_j^\gamma(z,t) = \sum_{\gamma=1}^p G_{j\alpha\beta}^\gamma(z,t)\cdot V_j^\gamma(z,t)
\end{align}

On the other hand, from $(\ref{oi10})$ and $(\ref{oi50})$, we have
\begin{align}\label{oi11}
V_j^\alpha= \sum_{\beta, \delta, \eta=1}^p c_j^{\alpha\delta} r_{0jk}^{\delta \eta} \left(c_k^{-1}\right)^{\eta \beta} V_k^\beta \iff \sum_{\delta=1}^p \left(c_j^{-1}\right)^{\alpha \delta} V_j^\delta = \sum_{\beta=1}^p r_{0jk}^{\alpha\beta} \left(\sum_{\delta=1}^p \left(c_k^{-1}\right)^{\beta \delta} V_k^\delta    \right)
\end{align}

Let us define
\begin{align*}
X_j(t):\Gamma\left(U_j, \Theta_{\mathcal{F}_0}\right)&\to \Gamma\left(U_j, \mathcal{A}^{0,0}\left(\Theta_{M}\right)\right) \\
              T_{0j}^\alpha &\mapsto X_j^\alpha :=\sum_{\delta=1}^p \left(c_j^{-1} \right)^{\alpha\delta} V_j^{\delta}
\end{align*}
and linearly extends to $\Gamma\left(U_j, \Theta_{\mathcal{F}_0} \right)$. Then $(\ref{oi11})$ implies that $X= \left\{X_j \right\}$ defines a global section in $A^{0,0}\left( M, \mathscr{H}om_{\mathcal{O}_M}\left(\Theta_{\mathcal{F}_0}, \Theta_M \right)   \right)$. We note that as in Part I, since $\bar{\partial} V_j^\delta - \left[ \varphi, V_j^\delta \right] =0$, we have from $(\ref{oi19})$
\begin{align*}
\bar{\partial}\left(\sum_{\delta=1}^p \left(c_j^{-1}\right)^{\alpha \delta} V_j^\delta    \right) -\left[ \varphi, \sum_{\delta=1}^p \left(c_j^{-1}\right)^{\alpha \delta} V_j^\delta  \right]&= -\sum_{\beta, \eta, \delta =1}^p \left(c_j^{-1} \right)^{\alpha \beta} \bar{\partial} c_j^{\beta \eta} \cdot \left(c_j^{\eta \delta}\right)^{-1}\cdot V_j^\delta + \sum_{\beta, \delta, \eta=1}^p \left(c_j^{-1}\right)^{\alpha\beta} \left[  \varphi, c_j^{\beta \eta}  \right]\left( c_j^{-1} \right)^{\eta \delta} V_j^\delta \\
&= - \sum_{\beta, \eta, \delta=1}^p \left(c_j^{-1} \right)^{\alpha\beta}\left(  \bar{\partial} c_j^{\beta \eta} - \left[  \varphi, c_j^{\beta \eta}  \right]  \right) \left(c_j^{-1} \right)^{\eta \delta} V_j^\delta
\end{align*}
This implies that
\begin{align}\label{oi25}
\bar{\partial} \left(X_j\left(T_{0j}^\alpha\right)\right) - \left[ \varphi, X_j\left(T_{0j}^\alpha\right)  \right] + X_j \left(U_j\left( T_{0j}^\alpha \right) \right)=0 \iff \bar{\partial} X_j^\alpha - \left[ \varphi, X_j^\alpha \right] + \sum_{\beta=1}^p U_j^{\alpha \beta} X_j^\beta = 0
\end{align}

We note that from $(\ref{oi20})$ and $(\ref{oi19})$, 
{\small{\begin{align}\label{oi21}
&\left[\sum_{\delta=1}^p \left(c_j^{-1}\right)^{\alpha \delta} V_j^\delta, \sum_{\eta=1}^p \left(c_j^{-1}\right)^{\beta \eta} V_j^\eta     \right] = \sum_{\delta,\eta=1}^p \left(c_j^{-1}\right)^{\alpha \delta}\left[ V_j^\delta, \left(c_j^{-1} \right)^{\beta \eta}  \right] V_j^\eta + \sum_{\delta, \eta=1}^p \left(c_j^{-1}\right)^{\beta \eta} \left[ \left(c_j^{-1}\right)^{\alpha \delta} V_j^\delta, V_j^\eta   \right]\\
&=  \sum_{\delta,\eta=1}^p \left(c_j^{-1}\right)^{\alpha \delta}\left[ V_j^\delta,\left(c_j^{-1}\right)^{\beta \eta}  \right] V_j^\eta - \sum_{\delta, \eta=1}^p \left(c_j^{-1}\right)^{\beta \eta}\left[V_j^\eta, \left(c_j^{-1}\right)^{\alpha \delta}  \right] V_j^\delta + \sum_{\delta, \eta, \gamma =1}^p \left(c_j^{-1}\right)^{\alpha \delta} \left(c_j^{-1}\right)^{\beta \eta} G_{j\delta \eta}^\gamma V_j^\gamma \notag \\
&= -\sum_{\delta, \eta, a,b=1}^n \left(c_j^{-1}\right)^{\alpha \delta} \left(c_j^{-1}\right)^{\beta a} \left[V_j^\delta, c_j^{ab} \right] \left(c_j^{-1}\right)^{b \eta} V_j^\eta + \sum_{\delta, \eta, b,d=1}^p \left(c_j^{-1}\right)^{\beta \eta} \left(c_j^{-1}\right)^{\alpha b}\left[ V_j^\eta , c_j^{b d}   \right] \left(c_j^{-1}\right)^{d\delta} V_j^\delta \notag  \\
&+ \sum_{\delta, \eta, \gamma, a, b=1}^p \left(c_j^{-1}\right)^{\alpha\delta} \left(c_j^{-1}\right)^{\beta \eta} G_{j\delta \eta}^{\gamma} c_j^{\gamma a} \left(c_j^{-1}\right)^{a b} V_j^b \notag
\end{align}}}
By setting
\begin{align*}
K_{j\alpha\beta}^\eta:= -\sum_{\delta, a=1}^n \left(c_j^{-1}\right)^{\alpha \delta} \left(c_j^{-1}\right)^{\beta a} \left[V_j^\delta, c_j^{a\eta} \right] + \sum_{\delta,d=1}^p \left(c_j^{-1}\right)^{\beta \delta} \left(c_j^{-1}\right)^{\alpha d}\left[ V_j^\delta , c_j^{d\eta}   \right] + \sum_{\delta, a, \gamma=1}^p \left(c_j^{-1}\right)^{\alpha \delta} \left(c_j^{-1}\right)^{\beta a} G_{j\delta a}^{\gamma} c_j^{\gamma \eta},
\end{align*}
$(\ref{oi21})$ implies that
\begin{align}\label{oi22}
\left[ X_j^\alpha, X_j^\beta\right]=\sum_{\gamma=1}^p K_{j\alpha\beta}^\gamma X_j^\gamma
\end{align}
We also note that $K_{j\alpha\beta}^\gamma(0)= G_{j\alpha\beta}^\gamma(z,0)= g_{j\alpha\beta}^\gamma\left(\xi_j(z,0),0 \right)$.

\subsection{Identification of Infinitesimal deformations by using Dolbeault resolution of $\mathcal{E}_{\Theta_{\mathcal{F}_0}}^\bullet$} \label{it23}\

In the previous subsection we showed from $(\ref{oi23}),(\ref{oi24}),(\ref{oi25})$ and $(\ref{oi22})$ that given a foliated analytic family $\left(\mathcal{M}, \Theta_\mathcal{F}, B, \pi \right)$ in terms of tangent sheaves of deformations of $\left(M, \Theta_{\mathcal{F}_0} \right)=\pi^{-1}(0)$, the deformations near $\left( M, \Theta_{\mathcal{F}_0} \right)$ are represented by a section $\left\{\left(\varphi(t), U_j(t) \right)\right\}\in A^{0,1}\left( M,\mathcal{E}_{\Theta_{\mathcal{F}_0}} \right)$, where $\varphi(t)\in A^{0,1}\left(M, \Theta_M \right)$ and  $U_j(t)\in \Gamma\left( U_j, \mathcal{A}^{0,1}\left( \mathscr{H}om_{\mathcal{O}_M}\left( \Theta_{\mathcal{F}_0}, \Theta_{\mathcal{F}_0} \right)\right) \right)$ which is  defined by $U_j(t)\left(T_{0j}^\alpha \right)= \sum_{\beta=1}^p U_j^{\alpha\beta}(t) T_{0j}^\beta$ for some $U_j^{\alpha\beta}(t)\in \Gamma\left(U_j, \mathcal{A}^{0,1} \right),\alpha=1...,p$, and $X(t)\in  A^{0,0}\left( M, \mathscr{H}om_{\mathcal{O}_M}\left( \Theta_{\mathcal{F}_0}, \Theta_M  \right) \right)$ satisfying $\varphi(0)=0, U_j(0)=0$, and $X(0)=\textnormal{Id}$ and the integrability condition: by setting $X_j(t)\left(T_{0j}^\alpha \right)= X_j^\alpha (t) \in \Gamma\left( U_j, \mathcal{A}^{0,0}\left(\Theta_M \right) \right),\alpha=1,...,p$, 
\begin{align}
&\bar{\partial} \varphi(t) - \frac{1}{2}\left[\varphi(t), \varphi (t)\right]=0 \label{oi26} \\
&\left[ X_j^\alpha (t), X_j^\beta(t)\right]= \sum_{\gamma=1}^p K_{j\alpha\beta}^\gamma (t)  X_j^\gamma (t) \,\,\,\,\,\,\textnormal{for some}\,\,\,\,\, K_{j\alpha\beta}^\gamma(t)\in \Gamma\left( U_j, \mathcal{A}^{0,0} \right) ,\alpha,\beta=1,...,p \label{oi27}\\
&\bar{\partial} X_j^\alpha (t) - \left[ \varphi (t) , X_j^\alpha (t) \right]+ \sum_{\eta=1}^p U_j^{\alpha \eta} (t) X_j^\eta (t) =0,\,\,\,\,\, \alpha=1,...,p \label{oi28}\\
&\bar{\partial} U_j^{\alpha \delta}(t)- \left[ \varphi(t), U_j^{\alpha \delta}(t) \right]  + \sum_{\sigma=1}^p U_j^{\alpha \sigma} (t) \wedge U_j^{\sigma \delta}(t) =0,\,\,\,\,\,\alpha, \delta=1,...,p \label{oi29}
\end{align}
where $t\in \Delta$ a sufficiently small polydisk in $B$. We identify infinitesimal deformations at $0$ by taking the derivative of $(\ref{oi26})-(\ref{oi29})$ with respect to $t$ and setting $t=0$. Let us take the derivative of $(\ref{oi26})-(\ref{oi29})$ with respect to $t$. Then we have
\begin{align}
&\bar{\partial}\left(\frac{\partial \varphi}{\partial t} \right) - \frac{1}{2}\left[ \frac{\partial \varphi}{\partial t}, \varphi(t)\right] - \frac{1}{2}\left[\varphi(t), \frac{\partial \varphi(t)}{\partial t} \right] =0 \label{oi30}\\
&\left[ \frac{\partial X_j^\alpha}{\partial t}, X_j^\beta\right] + \left[ X_j^\alpha, \frac{\partial X_j^\beta}{\partial t} \right] = \sum_{\gamma=1}^p \frac{\partial K_{j\alpha\beta}^\gamma}{\partial t} X_j^\gamma + \sum_{\gamma=1}^p K_{j\alpha\beta}^\gamma \frac{\partial X_j^\gamma}{\partial t} \label{oi31}\\
&\bar{\partial}\left(\frac{\partial X_j^\alpha}{\partial t} \right) -\left[ \frac{\partial\varphi}{\partial t}, X_j^\alpha \right] - \left[ \varphi, \frac{\partial X_j^\alpha}{\partial t} \right] + \sum_{\beta=1}^p \frac{\partial U_j^{\alpha \eta}}{\partial t} X_j^\eta + \sum_{\beta=1}^p U_j^{\alpha \eta} \frac{\partial X_j^\eta}{\partial t}\label{oi32} \\
&\bar{\partial}\left( \frac{\partial U_j^{\alpha \delta}}{\partial t} \right) - \left[ \frac{\partial \varphi}{\partial t}, U_j^{\alpha \delta}\right] - \left[ \varphi, \frac{\partial U_j^{\alpha \delta}}{\partial t} \right] + \sum_{\sigma=1}^p \frac{\partial U_j^{\alpha\sigma}}{\partial t} \wedge U_j^{\sigma \delta} + \sum_{\sigma=1}^p U_j^{\alpha \sigma} \wedge \frac{\partial U_j^{\sigma \delta}}{\partial t}=0\label{oi33}
\end{align}
By setting $t=0$ of $(\ref{oi30})-(\ref{oi33})$, we obtain 
\begin{align}
&\bar{\partial}\left(\frac{\partial \varphi}{\partial t}|_{t=0} \right)=0\label{oi42}\\
& \left[ \frac{\partial X_j^\alpha}{\partial t}|_{t=0}, T_{0j}^\beta \right] +\left[ T_{0j}^\alpha,  \frac{\partial T_j^\beta}{\partial t} |_{t=0} \right] - \sum_{\gamma=1}^p g_{0\alpha\beta}^\gamma \frac{\partial X_j^\gamma}{\partial t}|_{t=0} = \sum_{\gamma=1}^p \frac{\partial K_{j\alpha\beta}^\gamma}{\partial t}|_{t=0} T_{0j}^\gamma \label{oi43} \\
& \bar{\partial}\left(\frac{\partial X_j^\alpha}{\partial t}|_{t=0} \right) - \left[ \frac{\partial \varphi}{\partial t}|_{t=0}, T_{0j}^\alpha \right] + \sum_{\beta=1}^p \frac{\partial U_j^{\alpha \eta}}{\partial t}|_{t=0} T_{0j}^\eta  =0 \label{oi45}\\
&\bar{\partial}\left(\frac{\partial U_j^{\alpha\delta}}{\partial t}|_{t=0} \right)=0 \label{oi46}
\end{align}
We define $\frac{\partial X_j(t)}{\partial t}|_{t=0}\in A^{0,0}\left( M, \mathscr{H}om_{\mathcal{O}_M}\left(\Theta_{\mathcal{F}_0}, \Theta_M \right)   \right)$ by
\begin{align*}
\frac{\partial X_j}{\partial t}|_{t=0}:\Gamma\left( U_j, \Theta_{\mathcal{F}_0} \right) &\to \Gamma\left( U_j, \mathcal{A}^{0,0}\left(\Theta_M \right) \right) \\
 T_{0j}^\alpha &\mapsto \frac{\partial X_j^\alpha}{\partial t}
\end{align*}
and linearly extends to $\Gamma\left(U_j, \Theta_{\mathcal{F}_0} \right)$. Then $\frac{\partial X}{\partial t}|_{t=0}=\left\{ \frac{\partial X_j}{\partial t}|_{t=0} \right\}$ defines a global section in $A^{0,0}\left( M, \mathscr{H}om_{\mathcal{O}_M}\left( \Theta_{\mathcal{F}_0} , \Theta_M \right)\right)$. We define an element $\frac{\partial U_j}{\partial t}|_{t=0}\in \Gamma\left(U_j ,\mathcal{A}^{0,1}\left(\mathscr{H}om_{\mathcal{O}_M}\left(\Theta_{\mathcal{F}_0}, \Theta_{\mathcal{F}_0} \right) \right)\right)$ by
\begin{align*}
\frac{\partial U_j}{\partial t}|_{t=0}:\Gamma\left( U_j, \Theta_{\mathcal{F}_0} \right) &\to \Gamma\left( U_j,  \mathcal{A}^{0,1}\left(\Theta_{\mathcal{F}_0} \right) \right)\\
                                      T_{0j}^\alpha &\mapsto \sum_{\eta=1}^p \frac{\partial U_j^{\alpha\eta}}{\partial t}|_{t=0} T_{0j}^\eta
\end{align*}
and linearly extends to $\Gamma\left( U_j, \Theta_{\mathcal{F}_0}\right)$. We note that from $(\ref{oi40})$, we have
\begin{align}\label{oi41}
\sum_{\eta=1}^p U_j^{\alpha \eta} T_{0j}^\eta= \sum_{\beta, \eta=1}^p r_{0jk}^{\alpha\beta} U_j^{\beta \eta} T_{0k}^\eta + \sum_{\gamma=1}^p \left[ \varphi, r_{0jk}^{\alpha\beta} \right] T_{0k}^\beta
\end{align}
By taking the derivative of $(\ref{oi41})$ with respect to $t$ and setting $t=0$, we have
\begin{align*}
\sum_{\eta=1}^p \frac{\partial U_j^{\alpha \eta}}{\partial t}|_{t=0} T_{0j}^\eta = \sum_{\beta, \eta=1}^p r_{0jk}^{\alpha\beta}\frac{\partial U_j^{\beta \eta}}{\partial t}|_{t=0} T_{0k}^\eta  + \sum_{\beta=1}^p \left[ \frac{\partial \varphi}{\partial t}|_{t=0}, r_{0jk}^{\alpha\beta} \right] T_{0k}^\beta
\end{align*}
This implies that $\left\{ \left( \frac{\partial \varphi}{\partial t}|_{t=0}, \frac{\partial U_j}{\partial t}|_{t=0} \right) \right\}\in A^{0,1}\left(M, \mathcal{E}_{\Theta_{\mathcal{F}_0} }\right)$. 

Then $(\ref{oi42})-(\ref{oi46})$ implies that $\left(\frac{\partial X}{\partial t}|_{t=0}, \left\{ \left(\frac{\partial \varphi}{\partial t}|_{t=0}, \frac{\partial U_j}{\partial t}|_{t=0} \right)\right\} \right) \in A^{0,0}\left( M, \mathscr{H}om_{\mathcal{O}_M}\left(\Theta_{\mathcal{F}_0}, \Theta_M \right)\right)\bigoplus A^{0,1}\left(M, \mathcal{E}_{\Theta_{\mathcal{F}_0}} \right) $ defines a $1$-cocycle in the following Dolbeault resolution of $\mathcal{E}_{\Theta_{\mathcal{F}_0}}^\bullet$:
{\Small{\begin{equation} \label{oi47}
\begin{CD}
\cdots \\
@AD_2 AA \\
\frac{A^{0,0}\left( \mathcal{U}_0, \mathscr{H}om_{\mathcal{O}_M}\left( \bigwedge^2  \Theta_{\mathcal{F}_0}, \Theta_M  \right)\right)}{A^{0,0}\left(\mathcal{U}_0, \mathscr{H}om_{\mathcal{O}_M}\left( \bigwedge^2  \Theta_{\mathcal{F}_0}, \Theta_{\mathcal{F}_0}  \right) \right)} @>\bar{\partial}>> \cdots  \\
@AD_1'AA @AD_1 AA \\
A^{0,0}\left(\mathcal{U}_t, \mathscr{H}om_{\mathcal{O}_{M}}\left(\Theta_{\mathcal{F}_0}, \Theta_{M} \right) \right)@>-\bar{\partial}>>  A^{0,1}\left(\mathcal{U}_0, \mathscr{H}om_{\mathcal{O}_{M}}\left(\Theta_{\mathcal{F}_0}, \Theta_{M}  \right) \right) @>\bar{\partial} >> \cdots \\
@AD_0' AA @A{D}_0'AA @AD_0'AA \\
A^{0,0}\left(\mathcal{U}_0,\mathcal{E}_{\Theta_{\mathcal{F}_0}}\right) @>\bar{\partial}>> A^{0,1}\left(\mathcal{U}_0, \mathcal{E}_{\Theta_{\mathcal{F}_0}}\right) @>-\bar{\partial}>> A^{0,2} \left(\mathcal{U}_0, \mathcal{E}_{\Theta_{\mathcal{F}_0}} \right) @> \bar{\partial}>> \cdots 
\end{CD}
\end{equation}}}

Then we have the foliated Kodaira-Spencer map in terms of tangent sheaves
\begin{align}\label{a2}
\varphi_0: T_0 B &\to \mathbb{H}^1\left( M, \mathcal{E}_{\Theta_{\mathcal{F}_0}}^\bullet \right)\\
             \frac{\partial}{\partial t} &\mapsto \left(\frac{\partial X}{\partial t}|_{t=0}, \left\{\left(\frac{\partial \varphi}{\partial t}|_{t=0} , \frac{\partial U_j}{\partial t}|_{t=0}\right)\right\}   \right) \in A^{0,0}\left(M, \mathscr{H}om_{\mathcal{O}_M}\left( \Theta_{\mathcal{F}_0}, \Theta_M \right) \right) \bigoplus A^{0,1}\left( M, \mathcal{E}_{\Theta_{\mathcal{F}_0}} \right) \notag
\end{align}
In the following subsection \ref{ti1}, we will describe the explicit correspondence of infinitesimal deformations between \v Cech resolution and Dolbeault resolution.

\subsection{Identification of obstructions  by using Dolbeault resolution of $\mathcal{E}_{\Theta_{\mathcal{F}_0}}^\bullet$} \label{it30} \

We assume that the foliated analytic family $\left( \mathcal{M}, \Theta_\mathcal{F}, B, \pi \right)$ in the previous subsection is a $1$-parameter family, i.e $\Delta \subset \mathbb{C}$. Then by taking the derivative of $(\ref{oi30})-(\ref{oi33})$ with respect to $t$ again, we have
\begin{align}
&\bar{\partial}\left(\frac{\partial^2 \varphi}{\partial t^2} \right) -\left[\frac{\partial \varphi}{\partial t}, \frac{\partial \varphi}{\partial t} \right] - \frac{1}{2}\left[ \frac{\partial^2 \varphi}{\partial t^2}, \varphi \right] - \frac{1}{2} \left[ \varphi, \frac{\partial^2 \varphi}{\partial t^2}\right]=0 \label{it3}\\
&\left[\frac{\partial^2 X_j^\alpha}{\partial t^2}, X_j^\beta \right] + 2\left[ \frac{\partial X_j^\alpha}{\partial t}, \frac{\partial X_j^\beta}{\partial t} \right] + \left[X_j^\alpha, \frac{\partial^2 X_j^\beta}{\partial t^2} \right] = \sum_{\gamma=1}^p \frac{\partial^2 K_{j\alpha\beta}^\gamma}{\partial t^2} X_j^\gamma + 2 \sum_{\gamma=1}^p \frac{\partial K_{j\alpha\beta}^\gamma}{\partial t}\frac{\partial X_j^\gamma}{\partial t} + \sum_{\gamma=1}^p K_{j\alpha\beta}^\gamma \frac{\partial^2 X_j^\gamma}{\partial t^2}\\
&\bar{\partial}\left( \frac{\partial^2 X_j^\alpha}{\partial t^2} \right) - \left[ \frac{\partial^2 \varphi}{\partial t^2}|_{t=0},  T_{0j}^\alpha \right] -2 \left[\frac{\partial \varphi}{\partial t}, \frac{\partial X_j^\alpha}{\partial t} \right] - \left[\varphi, \frac{\partial^2 X_j^\alpha}{\partial t^2} \right] + \sum_{\eta=1}^p \frac{\partial^2 U_j^{\alpha \eta}}{\partial t^2} X_j^\eta + 2 \sum_{\eta=1}^p \frac{\partial U_j^{\alpha \eta}}{\partial t} \frac{\partial X_j^\eta}{\partial t} + \sum_{\eta=1}^p U_j^{\alpha \eta}\frac{\partial^2 X_j^\eta}{\partial t^2}=0 \\
&\bar{\partial}\left(\frac{\partial^2 U_j^{\alpha \delta}}{\partial t^2} \right) - 2\left[ \frac{\partial \varphi}{\partial t}, \frac{\partial U_j^{\alpha \delta}}{\partial t} \right] - \left[\frac{\partial^2 \varphi}{\partial t^2}, U_j^{\alpha \delta} \right] - \left[ \varphi, \frac{\partial^2 U_j^{\alpha \delta}}{\partial t} \right] + \sum_{\sigma=1}^p \frac{\partial^2 U_j^{\alpha \sigma}}{\partial t^2} \wedge U_j^{\sigma \delta} \label{it4}
\end{align}
Then by setting $t=0$ of $(\ref{it3})-(\ref{it4})$, we obtain
\begin{align}
&-\bar{\partial}\left( \frac{\partial^2 \varphi}{\partial t^2}|_{t=0} \right) = -\left[\frac{\partial \varphi}{\partial t}|_{t=0}, \frac{\partial \varphi}{\partial t}|_{t=0} \right] \label{it8}\\
& \left[\frac{\partial^2 X_j^\alpha}{\partial t^2}|_{t=0}, T_{0j}^\beta \right]  + \left[ T_{0j}^\alpha,  \frac{\partial^2 X_j^\beta}{\partial t^2}|_{t=0} \right] - \sum_{\gamma=1}^p g_{0j\alpha\beta}^\gamma \frac{\partial^2 X_j^\gamma}{\partial t^2}|_{t=0} - \sum_{\gamma=1}^p \frac{\partial^2 K_{j\alpha\beta}}{\partial t^2}|_{t=0}\cdot T_{0j}^\gamma \\
&= -2 \left[ \frac{\partial X_j^\alpha}{\partial t}|_{t=0}, \frac{\partial X_j^\beta}{\partial t}|_{t=0} \right] + 2\sum_{\gamma=1}^p \frac{\partial K_{j\alpha\beta}^\gamma}{\partial t}|_{t=0} \cdot \frac{\partial X_j^\gamma}{\partial t}|_{t=0} \notag \\
&- \bar{\partial}\left( \frac{\partial^2 X_j^\alpha}{\partial t^2}|_{t=0} \right) + \left[ \frac{\partial^2 \varphi}{\partial t^2}|_{t=0} , T_{0j}^\alpha \right]  - \sum_{\eta=1}^p \frac{\partial^2 U_j^{\alpha \eta}}{\partial t^2}|_{t=0}\cdot T_{0j}^\eta = - 2\left[ \frac{\partial \varphi}{\partial t}|_{t=0}, \frac{\partial X_j^\alpha}{\partial t}|_{t=0} \right] +  2 \sum_{\eta=1}^p \frac{\partial U_j^{\alpha \eta}}{\partial t}|_{t=0}\cdot \frac{\partial X_j^\eta}{\partial t}|_{t=0} \\
& - \bar{\partial}\left( \frac{\partial^2 U_j^{\alpha \delta}}{\partial t^2}|_{t=0} \right)= - 2\left[\frac{\partial \varphi}{\partial t}|_{t=0}, \frac{\partial U_j^{\alpha \delta}}{\partial t}|_{t=0} \right] +  2 \sum_{\sigma=1}^p \frac{\partial U_j^{\alpha \sigma}}{\partial t}|_{t=0} \wedge \frac{\partial U_j^{\sigma \delta}}{\partial t}|_{t=0} \label{it9}
\end{align}

\begin{lemma}\label{its30}
For any $1$-cocycle $\left(  \Phi =\left\{ \Phi_j \right\},   \left\{ \left( \psi, a_j \right)  \right\} \right) \in A^{0,0}\left( M, \mathscr{H}om_{\mathcal{O}_M}\left( \Theta_{\mathcal{F}_0}, \Theta_M  \right)  \right) \bigoplus A^{0,1}\left( M, \mathcal{E}_{\Theta_{\mathcal{F}_0}} \right)$ where $\Phi_j \in \Gamma\left( U_j, \mathcal{A}^{0,0}\left(\mathscr{H}om_{\mathcal{O}_M}\left( \Theta_{\mathcal{F}_0}, \Theta_M \right) \right) \right)$ and $\psi\in A^{0,1}\left(M, \Theta_M \right)$ and $a_j\in \Gamma\left( U_j, \mathcal{A}^{0,1}\left(\mathscr{H}om_{\mathcal{O}_M}\left( \Theta_{\mathcal{F}_0}, \Theta_{\mathcal{F}_0} \right) \right) \right)$ such that 
\begin{align*}
\Phi_j:\Gamma\left(U_j, \Theta_{\mathcal{F}_0} \right) &\to \Gamma\left( U_j, \mathcal{A}^{0,0}\left( \Theta_M \right) \right)\\
 T_{0j}^\alpha &\mapsto \Phi_j^\alpha
\end{align*}
with $\Phi_j=\Phi_k$ on $U_{jk}$ and 
\begin{align*}
a_j:\Gamma\left( U_j, \Theta_{\mathcal{F}_0} \right) &\to \Gamma\left(U_j , \mathcal{A}^{0,1}\left( \Theta_{\mathcal{F}_0} \right) \right) \\
 T_{0j}^\alpha &\mapsto a_j\left(T_{0j}^\alpha \right)= \sum_{\beta=1}^p a_j^{\alpha\beta} T_{0j}^\beta \,\,\,\,\,\,\,\,\, a_j^{\alpha\beta}\in \Gamma\left( U_j, \mathcal{A}^{0,1} \right)
\end{align*}
with $a_j\left(T_{0j}^\alpha \right)= a_k\left(T_{0j}^\alpha \right)+ \sum_{\beta=1}^p \left[ \psi, r_{0jk}^{\alpha\beta} \right] T_{0k}^\beta$, i.e. for $\alpha=1,....,p$,
\begin{align}\label{it5}
\sum_{\beta=1}^p a_j^{\alpha \beta} T_{0j}^\beta = \sum_{\beta, \eta=1}^p r_{0jk}^{\alpha\beta} a_k^{\beta \eta} T_{0k}^\eta +\sum_{\beta=1}^p \left[ \psi, r_{0jk}^{\alpha\beta} \right] T_{0k}^\beta \iff  \sum_{\beta=1}^p a_j^{\alpha\beta} r_{0jk}^{\beta \eta}= \sum_{\beta=1}^p r_{0jk}^{\alpha\beta} a_k^{\beta \eta} + \left[ \psi, r_{0jk}^{\alpha \eta} \right] 
\end{align}
 satisfying
\begin{enumerate}
\item $\left[T_{0j}^\alpha, \Phi_j\left( T_{0j}^\beta \right) \right] - \left[ T_{0j}^\beta, \Phi_j\left( T_{0j}^\alpha \right) \right] - \sum_{\gamma=1}^p g_{0j\alpha\beta}^\gamma \Phi_j\left(T_{0j}^\gamma \right) = K_{\alpha\beta}$ for some $K_{j\alpha\beta}\in \Gamma\left(U_j, \mathcal{A}^{0,0} \left(\Theta_{\mathcal{F}_0} \right)\right),\alpha,\beta=1,...,p$
\item $-\bar{\partial} \Phi_j\left(T_{0j}^\alpha \right) + \left[ \psi, T_{0j}^\alpha \right]- \sum_{\beta=1}^p a_j^{\alpha\beta} T_{0j}^\beta =0,\alpha=1,...,p$.
\item $-\bar{\partial} \psi=0$ and $-\bar{\partial} a_j^{\alpha \beta},\alpha,\beta =1,...,p$.
\end{enumerate}
we define 
\begin{align*}
\tilde{\rho}_j:\Gamma\left(U_j, \bigwedge^2 \Theta_{\mathcal{F}_0} \right) &\to \Gamma\left( U_j, \mathcal{A}^{0,0}\left(\Theta_M\right) \right)\\
 T_{0j}^\alpha \wedge T_{0j}^\beta &\mapsto -2\left[ \Phi_j\left( T_{0j}^\alpha \right), \Phi_j\left(T_{0j}^\beta \right)  \right]+ 2 \Phi_j\left(K_{j\alpha\beta} \right)
\end{align*}
Then $\tilde{\rho}=\left\{ \tilde{\rho}_i\right\}$ defines a global section in $A^{0,0}\left(M, \mathscr{H}om_{\mathcal{O}_M}\left( \bigwedge^2 \Theta_{\mathcal{F}_0}, \Theta_M  \right) \right)$, and define $\rho$ to be the image in $\frac{A^{0,0}\left(M, \mathscr{H}om_{\mathcal{O}_M}\left( \bigwedge^2 \Theta_{\mathcal{F}_0}, \Theta_M  \right) \right)}{A^{0,0}\left(M, \mathscr{H}om_{\mathcal{O}_M}\left( \bigwedge^2 \Theta_{\mathcal{F}_0}, \Theta_{\mathcal{F}_0}  \right) \right)}$, and we define $\eta_j \in \Gamma\left(U_j, \mathcal{A}^{0,1}\left( \mathscr{H}om_{\mathcal{O}_M} \left( \Theta_{\mathcal{F}_0}, \Theta_M  \right) \right) \right)$ by
\begin{align*}
\eta_j:\Gamma\left( U_j, \Theta_{\mathcal{F}_0} \right) &\to \Gamma\left( U_j, \mathcal{A}^{0,1}\left( \Theta_M \right) \right)\\
 T_{0j}^\alpha &\mapsto -2\left[ \psi, \Phi_j^\alpha \right] + 2\sum_{\eta=1}^p a_j^{\alpha \eta} \Phi_j^\eta
\end{align*}
Then $\eta=\{\eta_j\}$ defines a global section in $A^{0,0}\left( M,  \mathscr{H}om_{\mathcal{O}_M}\left( \Theta_{\mathcal{F}_0}, \Theta_M \right) \right)$. Then
{\small{\begin{align}\label{it10}
\left( \rho, \eta, \left( -\left[ \left(\psi, a_j \right), \left( \psi, a_j \right) \right] \right) \right) \in \frac{A^{0,0}\left( M,  \mathscr{H}om_{\mathcal{O}_M}\left( \bigwedge^2 \Theta_{\mathcal{F}_0}, \Theta_M \right) \right)}{A^{0,0}\left( M, \mathscr{H}om_{\mathcal{O}_M}\left( \bigwedge^2 \Theta_{\mathcal{F}_0}, \Theta_{\mathcal{F}_0} \right) \right)} \bigoplus A^{0,1}\left( M, \mathscr{H}om_{\mathcal{O}_M}\left( \Theta_{\mathcal{F}_0}, \Theta_M \right) \right) \bigoplus A^{0,2}\left( M , \mathcal{E}_{\Theta_{\mathcal{F}_0}} \right)
\end{align}}}
defines a $2$-cocycle in the Dolbeault resolution $(\ref{oi47})$ of $\mathcal{E}_{\Theta_{\mathcal{F}_0}}^\bullet$.
\end{lemma}

\begin{remark}

We show that $\left\{ \tilde{\rho}_j \right\}$ defines a global section in $A^{0,0}\left( M, \mathscr{H}om_{\mathcal{O}_M}\left( \bigwedge^2 \Theta_{\mathcal{F}_0}, \Theta_M \right) \right)$. We note that since $\Phi_j=\Phi_k$, we have
\begin{align*}
&K_{j\alpha\beta}= \left[ T_{0j}^\alpha, \Phi_j\left(T_{0j}^\beta \right) \right] -\left[ T_{0j}^\beta, \Phi_j\left( T_{0j}^\alpha \right) \right] - \sum_{\gamma=1}^p g_{0j\alpha\beta}^\gamma \Phi_j\left( T_{0j}^\gamma \right) \\
&= \left[ \sum_{\delta=1}^p r_{0jk}^{\alpha \delta} T_{0k}^\delta, \sum_{\eta=1}^p r_{0jk}^{\beta \eta} \Phi_k\left( T_{0k}^\eta\right)  \right] - \left[ \sum_{\eta=1}^p r_{0jk}^{\beta \eta } T_{0k}^\eta, \sum_{\delta=1 }^p r_{0jk}^{\alpha \delta} \Phi_k\left( T_{0k}^\delta \right)  \right] - \sum_{\xi,\gamma=1}^p g_{0j\alpha\beta}^\gamma r_{0jk}^{\gamma \xi} \Phi_k\left( T_{0k}^\xi \right)\\
&=\sum_{\delta,\eta=1}^p r_{0jk}^{\alpha \delta} T_{0k}^{\delta}\left( r_{0jk}^{\beta \eta} \right) \Phi_k\left( T_{0k}^\eta \right) + \sum_{\delta, \eta=1}^p r_{0jk}^{\alpha \delta} r_{0jk}^{\beta \eta} \left[T_{0k}^\delta, \Phi_k\left( T_{0k}^\eta \right) \right] - \sum_{\delta, \eta=1}^p r_{0jk}^{\beta \eta}\left[\Phi_k\left(T_{0k}^\eta \right) , r_{0jk}^{\alpha \delta}\right] T_{0k}^\delta\\
&- \sum_{\eta, \delta=1}^p r_{0jk}^{\beta \eta} T_{0k}^\eta\left(r_{0jk}^{\alpha \delta} \right)\Phi_k\left(T_{0k}^\delta \right) - \sum_{\eta, \delta=1}^p r_{0jk}^{\alpha\delta} r_{0jk}^{\beta \eta} \left[ T_{0k}^\eta, \Phi_k\left(T_{0k}^\delta \right) \right] + \sum_{\eta, \delta=1}^p r_{0jk}^{\alpha \delta} \left[ \Phi_k\left(T_{0k}^\delta \right),  r_{0jk}^{\beta \eta} \right] T_{0k}^\eta \\
&- \sum_{\xi=1}^p \left( \sum_{\delta, \eta=1}^p r_{0jk}^{\alpha \delta} r_{0jk}^{\beta \eta} g_{k\delta \eta}^\xi + \sum_{\delta=1}^p r_{0jk}^{\alpha \delta} T_{0j}^\delta \left( r_{0jk}^{\beta \xi } \right) - \sum_{\eta=1}^p r_{0jk}^{\beta \eta} T_{0j}^\eta \left( r_{0jk}^{\alpha \xi} \right)  \right) \Phi_k\left( T_{0k}^\xi \right)\\
&= \sum_{\delta, \eta=1}^p r_{0jk}^{\alpha \delta } r_{0jk}^{\beta \eta} K_{k  \delta \eta} - \sum_{\delta, \eta=1}^p r_{0jk}^{\beta \eta}\left[\Phi_k\left(T_{0k}^\eta \right) , r_{0jk}^{\alpha \delta}\right] T_{0k}^\delta  + \sum_{\eta, \delta=1}^p r_{0jk}^{\alpha \delta} \left[ \Phi_k\left(T_{0k}^\delta \right),  r_{0jk}^{\beta \eta} \right] T_{0k}^\eta
\end{align*}
On the other hand, we have
\begin{align*}
&\left[ \Phi_j\left(T_{0j}^\alpha \right), \Phi_j\left(T_{0j}^\alpha \right) \right] = \left[\sum_{\delta=1}^p r_{0jk}^{\alpha\delta} \Phi_k\left(T_{0k}^\delta \right) , \sum_{\eta=1}^p r_{0jk}^{\beta \eta} \Phi_k\left( T_{0k}^\eta  \right) \right]\\
&=\sum_{\delta, \eta=1}^p r_{0jk}^{\alpha \delta} \left[ \Phi_k\left(T_{0k}^\delta \right) , r_{0jk}^{\beta \eta}  \right] \Phi_k\left(T_{0k}^\eta \right)  + \sum_{\delta, \eta=1}^p r_{0jk}^{\alpha \delta} r_{0j}^{\beta \eta} \left [\Phi_k\left( T_{0k}^\delta \right), \Phi_k\left(T_{0k}^\eta \right) \right] - \sum_{\delta, \eta=1}^p r_{0jk}^{\beta \eta} \left[ \Phi_k\left( T_{0k}^\eta \right),  r_{0jk}^{\alpha \delta} \right] \Phi_k\left( T_{0k}^\delta \right)
\end{align*}

Then we have
\begin{align*}
- \left[ \Phi_j\left( T_{0j}^\alpha \right), \Phi_j \left( T_{0j}^\beta \right) \right] + \Phi_j\left( K_{j\alpha\beta} \right) =  \sum_{\delta, \eta=1}^p r_{0jk}^{\alpha \delta} r_{0j}^{\beta \eta} \left( -  \left[ \Phi_k\left( T_{0k}^\delta \right), \Phi_k \left( T_{0k}^\eta \right) \right]  + \Phi_k \left( K_{k\delta \eta}  \right) \right)
\end{align*}

This implies that $\left\{ \tilde{\rho}_j \right\}$ defines a global section.
\end{remark}

\begin{remark}
We show that $\left[ \left( \psi, a_j \right), \left( \psi, a_j \right) \right]\in A^{0,2}\left( M, \mathcal{E}_{\Theta_{\mathcal{F}_0}} \right)$. In fact, we recall from $(\ref{it5})$
\begin{align*}
\sum_{\delta=1}^p a_j^{\alpha \delta} T_{0j}^\delta= \sum_{\beta, \delta=1}^p r_{0jk}^{\alpha\beta} a_k^{\beta \delta} T_{0k}^\delta + \sum_{\beta=1}^p \left[ \psi, r_{0jk}^{\alpha\beta} \right] T_{0k}^\beta
\end{align*}
Then we have
\begin{align*}
&\sum_{\delta=1}^p \left( 2\left[ \psi, a_j^{\alpha \delta} \right]- 2\sum_{\sigma=1}^p a_j^{\alpha \sigma} \wedge a_j^{\sigma \delta}    \right) T_{0j}^\delta \\
&= \sum_{\delta=1}^p2 \left[ \psi, a_j^{\alpha \delta} T_{0j}^\delta \right] - 2 \sum_{\delta=1}^p a_j^{\alpha \delta} \wedge \left[\psi, T_{0j}^\delta \right] - 2\sum_{\sigma=1}^p a_j^{\alpha \sigma} \wedge\left( \sum_{\beta, \delta=1}^p r_{0jk}^{\sigma\beta} a_k^{\beta \delta} T_{0k}^\delta + \sum_{\beta=1}^p \left[ \psi, r_{0jk}^{\sigma\beta} \right] T_{0k}^\beta  \right)\\
&= \sum_{\delta=1}^p 2\left[ \psi, \sum_{\beta, \eta=1}^p r_{0jk}^{\alpha\beta} a_k^{\beta \eta} T_{0k}^\eta + \sum_{\beta=1}^p \left[ \psi, r_{0jk}^{\alpha \beta} \right] T_{0k}^\beta \right] - 2\sum_{\sigma=1}^p a_j^{\alpha \sigma} \wedge \left(\sum_{\beta=1}^p r_{0jk}^{\sigma \beta} \left( \sum_{\delta=1}^p a_k^{\beta \delta} T_{0k}^\delta    -\left[ \psi, T_{0k}^\beta \right]\right)   \right)\\
&= \sum_{\delta=1}^p 2\left[ \psi, \sum_{\beta, \eta=1}^p r_{0jk}^{\alpha\beta} a_k^{\beta \eta} T_{0k}^\eta + \sum_{\beta=1}^p \left[ \psi, r_{0jk}^{\alpha \beta} \right] T_{0k}^\beta \right] - 2\sum_{\beta=1}^p \left( \sum_{\sigma=1}^p r_{0jk}^{\alpha\sigma}a_k^{\sigma \beta} + \left[ \psi, r_{0jk}^{\alpha \beta} \right] \right) \wedge \left( \sum_{\delta=1}^p a_k^{\beta \delta} T_{0k}^\delta    -\left[ \psi, T_{0k}^\beta \right]   \right) \\
&= \sum_{\beta=1}^p r_{0jk}^{\alpha\beta}\left( \sum_{\eta=1}^p\left( 2\left[ \psi, a_k^{\beta \eta} \right]  -2 \sum_{\sigma=1}^p a_k^{\beta \sigma} \wedge a_k^{\sigma \eta}    \right) T_{0k}^\eta \right)  + \sum_{\beta=1}^p \left[  \left[ \psi, \psi \right], r_{0jk}^{\alpha\beta} \right] T_{0k}^\beta,
\end{align*}

\end{remark}

\begin{remark}
The following computation will be used in the proof of \textnormal{Lemma \ref{its30}}.
{\small{\begin{align}\label{it6}
\bar{\partial} K_{j\alpha\beta}&= \left[ T_{0j}^\alpha, \bar{\partial} \Phi_j^\beta \right] - \left[ T_{0j}^\beta, \bar{\partial} \Phi_j^\alpha \right] - \sum_{\gamma=1}^p g_{0j\alpha\beta}^\gamma \bar{\partial} \Phi_j^\gamma \\
&= \left[ T_{0j}^\alpha, \left[ \psi, T_{0j}^\beta \right]- \sum_{\eta=1}^p a_j^{\beta\eta} T_{0j}^\eta \right] - \left[ T_{0j}^\beta , \left[ \psi, T_{0j}^\alpha \right]- \sum_{\eta=1}^p a_j^{\alpha \eta} T_{0j}^\eta \right]- \sum_{\gamma=1}^p g_{0j\alpha\beta}^\gamma \left( \left[ \psi, T_{0j}^\gamma \right] - \sum_{\eta=1}^p a_j^{\gamma \eta} T_{0j}^\eta \right) \notag \\
&=- \sum_{\eta=1}^p \left[ T_{0j}^\alpha, a_j^{\beta \eta} \right] T_{0j}^\eta - \sum_{\eta, \xi=1}^p a_j^{\beta \eta} g_{0j\alpha\eta}^\xi T_{0j}^\xi + \sum_{\eta=1}^p \left[ T_{0j}^\beta, a_j^{\alpha\eta} \right] T_{0j}^\eta + \sum_{\eta, \xi=1}^p a_j^{\alpha \eta} g_{0j\beta\eta}^\xi T_{0j}^\xi + \sum_{\gamma=1}^p \left[ \psi, g_{0j\alpha\beta}^\gamma \right] T_{0j}^\gamma + \sum_{\gamma, \eta=1}^p g_{0j\alpha\beta}^\gamma a_j^{\gamma \eta} T_{0j}^\eta \notag \\
&=\sum_{\eta=1}^p \left( -\left[ T_{0j}^\alpha, a_j^{\beta \eta}\right] - \sum_{\xi=1}^p a_j^{\beta \xi} g_{0j\alpha \xi}^\eta + \left[ T_{0j}^\beta, a_j^{\alpha \eta} \right] + \sum_{\xi=1}^p a_j^{\alpha \xi} g_{0j\beta\xi}^\eta + \left[ \psi, g_{0j\alpha\beta}^\eta \right] + \sum_{\gamma=1}^p g_{0j\alpha\beta}^\gamma a_j^{\gamma \eta} \right) T_{0j}^\eta \notag
\end{align}}}
If we write $K_{j\alpha\beta} =\sum_{\eta=1}^p K_{j\alpha\beta}^\eta T_{0j}^\eta$ for some $K_{j\alpha\beta}^\gamma \in \Gamma\left( U_j, \mathcal{A}^{0,0} \right)$, then $\Phi_j\left( K_{j\alpha\beta} \right) = \sum_{\eta=1}^p K_{j\alpha\beta}^\eta \Phi_j\left( T_{0j}^\eta \right)$. Then from $(\ref{it6})$ we have
{\small{\begin{align}\label{it7}
&\bar{\partial} \Phi_j\left( K_{j\alpha\beta} \right)=\sum_{\eta=1}^p \bar{\partial} K_{j\alpha\beta}^\eta \cdot \Phi_j^\eta+ \sum_{\eta=1}^p K_{j\alpha\beta}^\eta \bar{\partial} \Phi_j^\eta \\
&= \sum_{\eta=1}^p \left( -\left[ T_{0j}^\alpha, a_j^{\beta \eta}\right] - \sum_{\xi=1}^p a_j^{\beta \xi} g_{0j\alpha \xi}^\eta + \left[ T_{0j}^\beta, a_j^{\alpha \eta} \right] + \sum_{\xi=1}^p a_j^{\alpha \xi} g_{0j\beta\xi}^\eta + \left[ \psi, g_{0j\alpha\beta}^\eta \right] + \sum_{\gamma=1}^p g_{0j\alpha\beta}^\gamma a_j^{\gamma \eta} \right) \Phi_{j}^\eta + \sum_{\eta=1}^p K_{j\alpha\beta}^\eta \left( \left[ \psi, T_{0j}^\eta \right]- \sum_{\xi=1}^p a_j^{\eta \xi} T_{0j}^\xi \right) \notag \\
&\equiv_{\Theta_{\mathcal{F}_0}} \sum_{\eta=1}^p \left( -\left[ T_{0j}^\alpha, a_j^{\beta \eta}\right] - \sum_{\xi=1}^p a_j^{\beta \xi} g_{0j\alpha \xi}^\eta + \left[ T_{0j}^\beta, a_j^{\alpha \eta} \right] +\sum_{\xi=1}^p a_j^{\alpha \xi} g_{0j\beta\xi}^\eta + \left[ \psi, g_{0j\alpha\beta}^\eta \right] + \sum_{\gamma=1}^p g_{0j\alpha\beta}^\gamma a_j^{\gamma \eta} \right) \Phi_{j}^\eta \notag \\
&\,\,\,\,\,\,\,\,\,\,+ \sum_{\eta=1}^p \left[ \psi,  \left[ T_{0j}^\alpha, \Phi_j^\beta \right]- \left[ T_{0j}^\beta, \Phi_j^\alpha \right] - \sum_{\gamma=1}^p g_{0j\alpha\beta}^\gamma \Phi_j^\gamma \right] \notag
\end{align}}}

\end{remark}

\begin{proof}[Proof of \textnormal{Lemma \ref{its30}}]
As in the same way of the proof of $D_2(R_j)=0$ in Lemma \ref{o54}, we can show that $\hat{D}_2\left(\rho_j \right)=0$, i.e. $\hat{D}_2(\rho)=0$. We note that from $(\ref{it7})$ 
{\small{\begin{align*} 
&\bar{\partial}\rho_j\left(T_{0j}^\alpha\wedge T_{0j}^\beta \right) + \hat{D}_1'\left( \eta_j \right)\left(T_{0j}^\alpha \wedge T_{0j}^\beta \right)\\
&= -2 \bar{\partial} \left[ \Phi_j^\alpha, \Phi_j^\beta  \right]  +2 \bar{\partial}\left(\Phi_j\left( K_{j\alpha\beta} \right) \right) + \left[ T_{0j}^\alpha, \eta_j\left(T_{0j}^\beta \right) \right] - \left[ T_{0j}^\beta, \eta_j\left( T_{0j}^\alpha \right) \right] - \sum_{\gamma=1}^p g_{0j\alpha\beta}^\gamma \eta_j\left( T_{0j}^\gamma \right)\\
&\equiv_{\Theta_{\mathcal{F}_0}} -2 \left[ \left[ \psi, T_{0j}^\alpha \right] - \sum_{\eta=1}^p a_j^{\alpha\eta} T_{0j}^\eta  , \Phi_j^\beta \right]- 2 \left[ \Phi_j^\alpha ,  \left[ \psi, T_{0j}^\beta \right] - \sum_{\eta=1}^p a_j^{\beta \eta} T_{0j}^\eta \right] + 2 \sum_{\eta=1}^p \left[ \psi,  \left[ T_{0j}^\alpha, \Phi_j^\beta \right]- \left[ T_{0j}^\beta, \Phi_j^\alpha \right] - \sum_{\gamma=1}^p g_{0j\alpha\beta}^\gamma \Phi_j^\gamma \right]\\
&+ 2\sum_{\eta=1}^p \left( -\left[ T_{0j}^\alpha, a_j^{\beta \eta}\right] - \sum_{\xi=1}^p a_j^{\beta \xi} g_{0j\alpha \xi}^\eta + \left[ T_{0j}^\beta, a_j^{\alpha \eta} \right] + \sum_{\xi=1}^p a_j^{\alpha \xi} g_{0j\beta\xi}^\eta + \left[ \psi, g_{0j\alpha\beta}^\eta \right] + \sum_{\gamma=1}^p g_{0j\alpha\beta}^\gamma a_j^{\gamma \eta} \right) \Phi_{j}^\eta\\
& - 2\left[ T_{0j}^\alpha, \left[ \psi, \Phi_j^\beta \right]  -  \sum_{\eta=1}^p a_j^{\beta \eta} \Phi_j^\eta \right] + 2 \left[ T_{0j}^\beta, \left[ \psi, \Phi_j^\alpha \right] - \sum_{\eta=1}^p a_j^{\alpha \eta} \Phi_j^\eta \right] + 2 \sum_{\gamma=1}^p g_{0j\alpha\beta}^\gamma\left(\left[ \psi, \Phi_j^\gamma \right] - \sum_{\eta=1}^p a_j^{\gamma \eta} \Phi_j^\eta \right) \\
&\equiv_{\Theta_{\mathcal{F}_0}} -2 \left[ \left[ \psi , T_{0j}^\alpha \right], \Phi_j^\beta \right] + 2 \sum_{\eta=1}^p a_j^{\alpha \eta} \left[ T_{0j}^\eta, \Phi_j^\beta \right] - 2\left[ \Phi_j^\alpha, \left[ \psi, T_{0j}^\beta \right] \right] + 2 \sum_{\eta=1}^p a_j^{\beta \eta}\left[ \Phi_j^\alpha, T_{0j}^\eta \right] + 2 \sum_{\eta=1}^p \left[ \psi, \left[ T_{0j}^\alpha, \Phi_j^\beta \right] \right]- 2 \sum_{\eta=1}^p \left[ \psi, \left[ T_{0j}^\beta, \Phi_j^\alpha \right] \right]\\
&- 2 \sum_{\xi=1}^p a_j^{\beta \xi} g_{0j\alpha\xi}^\eta \Phi_j^\eta +2  \sum_{\xi=1}^p a_j^{\alpha \xi} g_{0j\beta \xi}^\eta \Phi_j^\eta - 2 \left[ T_{0j}^\alpha, \left[ \psi, \Phi_j^\beta \right] \right] + 2 \sum_{\eta=1}^p a_j^{\beta \eta}\left[ T_{0j}^\alpha, \Phi_j^\eta \right] + 2 \left[ T_{0j}^\beta, \left[ \psi, \Phi_j^\alpha \right] \right] - 2 \sum_{\eta=1}^p a_j^{\alpha \eta} \left[ T_{0j}^\beta, \Phi_j^\eta \right] \\
&=0
\end{align*}}}

We show that $\bar{\partial} \eta_j\left(T_{0j}^\alpha \right) - \left[ \left[ \psi, \psi,\right], T_{0j}^\alpha \right] + 2\sum_{\delta=1}^p \left(\left[ \psi, a_j^{\alpha \delta} \right] - \sum_{\sigma=1}^p a_j^{\alpha \sigma} \wedge a_j^{\sigma \delta} \right) T_{0j}^\delta=0$. In fact,
\begin{align*}
\bar{\partial} \eta_j\left(T_{0j}^\alpha \right) &=\bar{\partial}\left(-2\left[\psi, \Phi_j^\alpha \right]  +  2\sum_{\eta=1}^p a_j^{\alpha \eta} \Phi_j^\eta \right) = 2 \left[ \psi, \bar{\partial} \Phi_j^\alpha \right]  - 2 \sum_{\eta=1}^p a_j^{\alpha \eta} \wedge \bar{\partial} \Phi_j^\eta \\
&= 2 \left[ \psi, \left[ \psi, T_{0j}^\alpha \right] - \sum_{\beta=1}^p a_j^{\alpha\beta} T_{0j}^\beta \right] - 2 \sum_{\eta=1}^p  a_j^{\alpha \eta} \wedge \left( \left[ \psi, T_{0j}^\eta \right] - \sum_{\xi=1}^p a_j^{\eta \xi} T_{0j}^\xi \right)\\
&= 2\left[ \left[ \psi, \psi\right], T_{0j}^\alpha \right] -2 \sum_{\beta=1}^p \left[ \psi, a_j^{\alpha\beta} \right] T_{0j}^\beta +2 \sum_{\eta, \xi=1}^p  a_j^{\alpha \eta} \wedge a_j^{\eta \xi} T_{0j}^\xi \\
&=2 \left[ \left[ \psi, \psi \right], T_{0j}^\alpha \right] - 2 \left( \sum_{\beta=1}^p \left[ \psi, a_j^{\alpha \beta} \right]  - \sum_{\eta=1}^p a_j^{\alpha \eta} \wedge a_j^{\eta \beta} \right) T_{0j}^\beta
\end{align*} 

Since $\bar{\partial}\psi=0, \bar{\partial} a_j^{\alpha \sigma}=0$, we have
\begin{align*}
\bar{\partial}\left[\psi, \psi \right]=0,\,\,\,\,\,\,\,\,\,\,\,\bar{\partial}\left(\left[ \psi, a_j^{\alpha \delta} \right] - \sum_{\sigma=1}^p a_j^{\alpha \sigma } \wedge a_j^{\sigma \delta}    \right) =0  \Longrightarrow \bar{\partial}\left( -\left[ \left( \psi, a_j \right),\left( \psi, a_j \right)      \right] \right) = 0
\end{align*}
This completes the proof of Lemma \ref{its30}.
\end{proof}

From $(\ref{it8})-(\ref{it9})$ and Lemma \ref{its30}, we obtain
\begin{theorem}
Suppose given a compact foliated complex manifold $\left( M, \Theta_{\mathcal{F}_0} \right)$ with $\Theta_{\mathcal{F}_0}$ locally free, and $\left( \Phi, \left\{\left( \psi, a_j\right) \right\} \right)\in \mathbb{H}^1\left( M, \mathcal{E}_{\Theta_{\mathcal{F}_0}}^\bullet \right)$. In order that there may exist a foliated analytic family $\left( \mathcal{M} , \Theta_{\mathcal{F}}, B, \pi \right)$ in terms of tangent sheaves such that $\pi^{-1}(0)=\left( M, \Theta_{\mathcal{F}_0} \right)$ and that $\left( \frac{\partial \left( \mathcal{M}, \Theta_\mathcal{F} \right)}{\partial t} \right)|_{t=0}=\left( \Phi, \left\{\left( \psi, a_j \right) \right\} \right)$, it is necessary that
{\small{\begin{align*}
\left( \rho, \eta, \left( -\left[ \left(\psi, a_j \right), \left( \psi, a_j \right) \right] \right) \right) \in \frac{A^{0,0}\left( M,  \mathscr{H}om_{\mathcal{O}_M}\left( \bigwedge^2 \Theta_{\mathcal{F}_0}, \Theta_M \right) \right)}{A^{0,0}\left( M, \mathscr{H}om_{\mathcal{O}_M}\left( \bigwedge^2 \Theta_{\mathcal{F}_0}, \Theta_{\mathcal{F}_0} \right) \right)} \bigoplus A^{0,1}\left( M, \mathscr{H}om_{\mathcal{O}_M}\left( \Theta_{\mathcal{F}_0}, \Theta_M \right) \right) \bigoplus A^{0,2}\left( M , \mathcal{E}_{\Theta_{\mathcal{F}_0}} \right)
\end{align*}}}
from $(\ref{it10})$ in \textnormal{Lemma \ref{its30}} defines the $0$ class in $\mathbb{H}^2\left( M, \mathcal{E}_{\Theta_{\mathcal{F}_0}}^\bullet \right)$.
\end{theorem}

\section{Isomorphisms of infinitesimal foliated deformations and obstructions in terms of tangent sheaves between \v Cech resolution and Dolbeault resolution of $\mathcal{E}_{\Theta_{\mathcal{F}_0}}^\bullet$} \label{ti2}

In this section we explicitly construct isomorphisms of first and second cohomology groups between \v Cech resolution of $\mathcal{E}_{\Theta_{\mathcal{F}_0}}^\bullet$ and \v Dolbeault resolution of $\mathcal{E}_{\Theta_{\mathcal{F}_0}}^\bullet$. Let $\mathcal{U}=\left\{U_j \right\}$ be an Stein open covering of $M$ by coordinate neighborhoods with coordinates $z_j=\left(z_j^1,..., z_j^n \right)$ on $U_j$, and $\Theta_{\mathcal{F}_0}$ is generated by $T_{0j}^1,..., T_{0j}^p \in \Gamma\left( U_j , \Theta_M \right)$ with the relation $T_{0j}^\alpha= \sum_{\beta=1}^p r_{0jk}^{\alpha\beta} T_{0j}^\beta$ on $U_{jk}$ for some $r_{0jk}^{\alpha\beta}\in \Gamma\left( U_{jk}, \mathcal{O}_M \right)$ and satisfies $\left[ T_{0j}^\alpha, T_{0j}^\beta \right]= \sum_{\gamma=1}^q g_{0j\alpha\beta}^\gamma T_{0j}^\gamma$ for some $g_{0j\alpha\beta}^\gamma \in \Gamma\left( U_j , \mathcal{O}_M \right),\alpha, \beta=1,...,p$.

\subsection{Isomorphisms of infinitesimal deformations between \v Cech and Dolbeault resolutions of $\mathcal{E}_{\Theta_{\mathcal{F}_0}}^\bullet$ }\label{ti1} \

We will explicitly define an isomorphism between the first cohomology groups $\mathbb{H}^1\left( M, \mathcal{E}_{\Theta_{\mathcal{F}_0}}^\bullet \right)$ from \v Cech resolution and the Dolbeault resolution of $\mathcal{E}_{\Theta_{\mathcal{F}_0}}^\bullet$.

Let $\left( \left\{ \alpha_j\right\}, \left\{ \left( \theta_{jk}, A_{jk} \right)      \right\}\right) \in C^0\left( \mathcal{U}, \mathscr{H}om_{\mathcal{O}_M}\left( \Theta_{\mathcal{F}_0}, \Theta_M \right) \right) \bigoplus C^1\left( \mathcal{U}, \mathcal{E}_{\Theta_{\mathcal{F}_0}} \right)$ be an $1$-cocycle in the \v Cech resolution of $\mathcal{E}_{\Theta_{\mathcal{F}_0}}^\bullet$ with respect to the covering $\mathcal{U}$, so that $\left(\theta_{ij}, A_{ij}\right)- \left(\theta_{ik}, A_{ik}\right)+ \left(\theta_{jk}, A_{jk}^i\right)=0$ and $\alpha_k\left( T_{0j}^\alpha \right) - \alpha_j\left(T_{0j}^\alpha \right)+ \mathcal{L}_{\theta_{jk}}\left(T_{0j}^\alpha \right)- A_{jk}\left( T_{0j}^\alpha \right)=0$, where $A_{jk}^\alpha=A_{jk}\left(T_{0j}^\alpha \right) = \sum_{\beta=1}^p A_{jk}^{\alpha \delta} T_{0j}^\delta$ with $A_{jk}^{\alpha\delta}\in \Gamma\left(U_{jk}, \mathcal{O}_M \right)$, and $\left[ T_{0j}^\alpha,\alpha_j\left( T_{0j}^\alpha \right) \right] - \left[ T_{0j}^\beta, \alpha_j\left( T_{0j}^\alpha\right) \right]- \sum_{\gamma=1}^p g_{0j\alpha\beta}^\gamma \alpha_j\left( T_{0j}^\gamma \right)= Q_{j\alpha\beta}\in \Gamma\left(U_j, \Theta_{\mathcal{F}_0} \right),\alpha,\beta=1,...,p$. Since $\delta\left( \left\{\theta_{jk}, A_{jk} \right\}\right)=0$, there exists $-\delta\left(\left\{ \left(\xi_j, A_j\right) \right\} \right)= \left\{ \left( \theta_{jk}, A_{jk}\right) \right\}$ where $\left\{\left( \xi_j, A_j\right) \right\}\in C^0\left(\mathcal{U}, \mathcal{A}^{0,0}\left( \mathcal{E}_{\Theta_{\mathcal{F}_0}} \right) \right)$ such that $\xi_j\in \Gamma\left( U_j, \mathcal{A}^{0,0}\left(\Theta_M \right) \right)$ and $A_j\in \Gamma\left( U_j, \mathcal{A}^{0,0}\left( \mathscr{H}om_{\mathcal{O}_M}\left(\Theta_{\mathcal{F}_0}, \Theta_{\mathcal{F}_0} \right)\right) \right)$ with $A_j\left(T_{0j}^\alpha \right) = \sum_{\beta=1}^p A_j^{\alpha \beta} T_{0j}^\beta$ satisfying $\xi_j-\xi_k= \theta_{jk}$ and $A_j- A_k^j= A_{jk}$ with $A_k^j\left(T_{0j}^\alpha \right)= A_k \left( T_{0j}^\alpha \right)+ \sum_{\beta=1}^p \left[ \xi_k, r_{0jk}^{\alpha\beta} \right] T_{0k}^\beta$. Then $\left\{\left( \psi =\bar{\partial} \xi_j, \bar{\partial} A_j \right) \right\}\in A^{0,0}\left(M, \mathcal{E}_{\Theta_{\mathcal{F}_0}} \right)$ where $\phi\in A^{0,1}\left( M, \Theta_M \right)$ and $\bar{\partial} A_j\in \Gamma\left(U_j, \mathcal{A}^{0,1}\left( \mathscr{H}om_{\mathcal{O}_M}\left( \Theta_{\mathcal{F}_0}, \Theta_{\mathcal{F}_0} \right) \right) \right)$. On the other hand, let us define $\Phi_j\in \Gamma\left(U_j, \mathcal{A}^{0,0}\left(\mathscr{H}om_{\mathcal{O}_M}\left( \Theta_{\mathcal{F}_0}, \Theta_M \right) \right) \right)$ by
\begin{align*}
\Phi_j:\Gamma\left( U_j, \Theta_{\mathcal{F}_0} \right) &\to \Gamma\left( U_j, \mathcal{A}^{0,0} \left( \Theta_M \right) \right)\\
 T_{0j}^\alpha &\mapsto \left[ \xi_j, T_{0j}^\alpha \right]- A_j\left(T_{0j}^\alpha \right) - \alpha_j\left(T_{0j}^\alpha \right)
\end{align*}
Then we have
\begin{align*}
&\Phi_j\left( T_{0j}^\alpha \right) - \sum_{\beta=1}^q r_{0jk}^{\alpha\beta} \Phi_k\left( T_{0k}^\beta \right) = \left[ \xi_j, T_{0j}^\alpha \right]- A_j\left(T_{0j}^\alpha \right) - \alpha_j\left(T_{0j}^\alpha \right) - \sum_{\beta=1}^p r_{0jk}^{\alpha\beta}\left(\left[ \xi_k, T_{0k}^\beta \right]- A_k\left(T_{0k}^\beta \right) - \alpha_k\left(T_{0k}^\beta \right) \right)\\
&= \left[ \xi_j, T_{0j}^\alpha \right]- A_j\left(T_{0j}^\alpha \right) - \alpha_j\left(T_{0j}^\alpha \right) - \left[ \xi_k, T_{0j}^\alpha \right] + \sum_{\beta=1}^p \left[ \xi_k , r_{0jk}^{\alpha\beta}  \right] T_{0k}^\beta + A_k\left( w_{0j}^\alpha \right) + \alpha_k\left(w_{0j}^\alpha \right)\\
&= \left[ \theta_{jk}, T_{0j}^\alpha \right] - A_j\left(T_{0j}^\alpha \right) + A_k^j\left( T_{0j}^\alpha \right) + \alpha_k\left( T_{0j}^\alpha\right)- \alpha_j\left( T_{0j}^\alpha \right) =\left[ \theta_{jk}, T_{0j}^\alpha \right] - A_{jk}\left(w_{0j}^\alpha \right) + \alpha_k\left( T_{0j}^\alpha \right)- \alpha_j\left( T_{0j}^\alpha \right)=0
\end{align*}
This implies that $\Phi:=\left\{ \Phi_j\right\}\in A^{0,0}\left(M, \mathscr{H}om_{\mathcal{O}_M}\left(\Theta_{\mathcal{F}_0} , \Theta_M \right) \right)$. Then we claim that $\left( \Phi, \left\{ \left(  \psi ,  \bar{\partial } A_j  \right) \right\} \right)$ defines a $1$-cocycle in the Dolbeault resolution. In fact, it is enough to check that
\begin{align*}
-\bar{\partial}\left( \Phi_j\left(T_{0j}^\alpha \right) \right) + \mathcal{L}_{\psi}\left(T_{0j}^\alpha \right) - \bar{\partial} A_j\left( T_{0j}^\alpha \right) = - \bar{\partial}\left(\left[ \xi_j, T_{0j}^\alpha \right] - A_j\left(T_{0j}^\alpha \right) \right) + \left[ \bar{\partial} \xi_j, T_{0j}^\alpha \right] - \bar{\partial} A_j\left(T_{0j}^\alpha \right) = 0
\end{align*}

Then we define the correspondence of $1$-cocycle between \v Cech resolution and Dolbeault resolution of $\mathcal{E}_{\Theta_{\mathcal{F}_0}}^\bullet$ by 
\begin{align}\label{a1}
\left(\left\{\alpha_j \right\}, \left\{ \left( \theta_{jk}, A_{jk} \right) \right\} \right) \mapsto \left( \Phi, \left\{ \left(\psi, \bar{\partial} A_j\right) \right\}\right).
\end{align}

\begin{proposition} \label{b18}
We recall that given a foliated analytic family $\left(\mathcal{M}, \Theta_\mathcal{F}, B, \pi \right)$ of deformations of $\left( M, \Theta_{\mathcal{F}_0} \right)= \pi^{-1}(0)$ as in \textnormal{Remark \ref{o4}}, the infinitesimal deformation at $t=0$ is given by $\left( \left\{ \alpha_j(0)\right\}, \left\{ \left( \theta_{jk}(0), A_{jk}(0) \right) \right\} \right)$ from $(\ref{o10})$ in the first cohomology group of the \v Cech resolution of $\mathcal{E}_{\Theta_{\mathcal{F}_0}}^\bullet$. On the other hand, the infinitesimal deformation at $t=0$ is also given by $\left( \frac{\partial X}{\partial t}|_{t=0}, \left\{  \left( \frac{\partial \varphi}{\partial t}|_{t=0}    , \frac{\partial U_j}{\partial t}|_{t=0}  \right)  \right\} \right)$ from $(\ref{a2})$in the first cohomology of the Dolbeault resolution of $\mathcal{E}_{\Theta_{\mathcal{F}_0}}^\bullet$. Under the isomorphism $(\ref{a1})$,
\begin{align*}
\left( \left\{ \alpha_j(0)\right\}, \left\{ \theta_{jk}(0), A_{jk}(0) \right\} \right)\,\,\,\,\,\,\,\,\,\textnormal{corresponds to}\,\,\,\,\,\,\left( \frac{\partial X}{\partial t}|_{t=0}, \left\{ \left(  \frac{\partial \varphi}{\partial t}|_{t=0}, \frac{\partial U_j}{\partial t}|_{t=0}   \right)  \right\}  \right)
\end{align*}
\end{proposition}

\begin{proof}
We fix a tangent vector $\frac{\partial}{\partial t}\in T_0\left( \Delta\right)$, denote $\left(\frac{\partial f(t)}{\partial t} \right)_{t=0}$ by $\dot{f}$ for a $C^\infty$ function $f(t),t\in \Delta$. With this notation, we put
\begin{align*}
\xi_j=\sum_{\alpha=1}^n \dot{\xi}_j^\alpha \frac{\partial}{\partial \xi_j^\alpha},\,\,\,\,\,\textnormal{where}\,\,\,\dot{\xi}_j^\alpha = \left(  \frac{\partial \xi_j(z,t)}{\partial t} \right)_{t=0}
\end{align*}
for each $j$. Then we have
\begin{align*}
\xi_j- \xi_k = \theta_{jk},\,\,\,\,\,\, \textnormal{and}\,\,\,\,\,\bar{\partial} \xi_j= \frac{\partial \varphi}{\partial t}|_{t=0}
\end{align*}
(for the detail, see \cite{Kod05} Theorem 5.4 p.266). Differentiating $(\ref{oi14})$, we have
\begin{align*}
\sum_{\beta=1}^p \frac{\partial c_j^{\alpha\beta}}{\partial t}\cdot r_{0jk}^{\beta \eta}= \sum_{\beta,\gamma=1}^p \frac{\partial r_{jk}^{\alpha\beta}}{\partial \xi_k^\gamma} \frac{\partial \xi_k^\gamma}{\partial t} \cdot c_k^{\beta \eta} +\sum_{\beta=1}^p \frac{\partial r_{jk}^{\alpha\beta}}{\partial t}\cdot c_k^{\beta \eta} + \sum_{\beta=1}^p r_{jk}^{\alpha\beta}\left( \xi_k(z,t,),t \right) \cdot \frac{\partial c_k^{\beta \eta}}{\partial t}
\end{align*}
and setting $t=0$, we obtain
\begin{align*}
&\sum_{\beta=1}^p \frac{\partial c_j^{\alpha\beta}}{\partial t}|_{t=0}\cdot r_{0jk}^{\beta \eta} = \sum_{\beta=1}^p \frac{\partial r_{jk}^{\alpha\beta}}{\partial \xi_k^\eta}|_{t=0}\cdot \frac{\partial \xi_k^\eta}{\partial t}|_{t=0} +  \frac{\partial r_{jk}^{\alpha\eta}}{\partial t}|_{t=0} + \sum_{\beta=1}^p r_{0jk}^{\alpha\beta}\cdot \frac{\partial c_k^{\beta \eta}}{\partial t}|_{t=0} \\
&\Longrightarrow \frac{\partial c_j^{\alpha\beta}}{\partial t}|_{t=0} = \sum_{\eta=1}^p\left[ \xi_k, r_{0jk}^{\alpha\eta} \right] r_{0kj}^{\eta \beta} + \sum_{\eta=1}^p \frac{\partial r_{jk}^{\alpha \eta}}{\partial t}|_{t=0}\cdot r_{0kj}^{\eta \beta} + \sum_{\eta, \gamma=1}^p r_{0jk}^{\alpha\eta} \cdot \frac{\partial c_k^{\eta \gamma}}{\partial t}|_{t=0} \cdot r_{0kj}^{\gamma \beta} \\
&\Longrightarrow \sum_{\beta=1}^p \frac{\partial c_j^{\alpha\beta}}{\partial t}|_{t=0} T_{0j}^\beta= \sum_{\eta=}^p \left[ \xi_k, r_{0jk}^{\alpha\eta} \right] T_{0k}^\eta +\sum_{\eta=1}^p \frac{\partial r_{jk}^{\alpha\eta}}{\partial t} |_{t=0} T_{0k}^\eta +  \sum_{\eta, \gamma=1}^p r_{0jk}^{\alpha \eta}\cdot \frac{\partial c_k^{\eta \gamma}}{\partial t}|_{t=0} T_{0k}^\gamma
\end{align*}
If we define $A_j\in \Gamma\left( \mathcal{U}, \mathscr{H}om_{\mathcal{O}_M}\left( \Theta_{\mathcal{F}_0}, \Theta_{\mathcal{F}_0} \right) \right)$ by
\begin{align*}
A_j:\Gamma\left( U_j,  \Theta_{\mathcal{F}_0} \right) &\to \Gamma\left( U_j, \mathcal{A}^{0,0}\left( \Theta_{\mathcal{F}_0} \right) \right) \\
              T_{0j}^\alpha &\mapsto \sum_{\beta=1}^p \frac{\partial c_j^{\alpha\beta}}{\partial t}|_{t=0} T_{0j}^\beta
\end{align*}
Then we see that $-\delta\left(  \left\{ \left(  \xi_j, A_j \right) \right\} \right)=\left( \theta_{jk}(0), A_{jk}(0) \right)$ and $\bar{\partial} A_j= \frac{\partial U_j}{\partial t}|_{t=0}$.

On the other hand, differentiating $(\ref{oi4})$, we have
\begin{align*}
\sum_{\gamma=1}^n \frac{\partial V_j^{\alpha \gamma}}{\partial t} \frac{\partial \xi_j^\beta}{\partial z^\gamma} + \sum_{\gamma=1}^n V_j^{\alpha \gamma} \frac{\partial}{\partial z^\gamma}\left( \frac{\partial \xi_j^\beta}{\partial t} \right) = \sum_{\gamma=1}^n \frac{\partial T_j^{\alpha\beta}}{\partial \xi_j^\gamma} \frac{\partial \xi_j^\gamma}{\partial t} + \frac{\partial T_j^{\alpha\beta}}{\partial t}
\end{align*}
and setting $t=0$, we have
\begin{align*}
 &\frac{\partial V_j^{\alpha \beta}}{\partial t}|_{t=0} + \sum_{\gamma=1}^n T_{0j}^{\alpha \gamma} \frac{\partial}{\partial z^\gamma}\left( \frac{\partial \xi_j^\beta}{\partial t}|_{t=0} \right) = \sum_{\gamma=1}^n \frac{\partial T_{0j}^{\alpha\beta}}{\partial \xi_j^\gamma}\cdot \frac{\partial \xi_j^\gamma}{\partial t}|_{t=0} + \frac{\partial T_j^{\alpha\beta}}{\partial t}|_{t=0} \\
 & \Longrightarrow \frac{\partial V_j^\alpha}{\partial t}|_{t=0} =  \left[ \xi_j, T_{0j}^\alpha\right] + \frac{\partial T_j^\alpha}{\partial t}|_{t=0} \Longrightarrow  \frac{\partial V_j^\alpha}{\partial t}|_{t=0} - \sum_{\beta=1}^p \frac{\partial c_j^{\alpha\beta}}{\partial t}|_{t=0} T_{0j}^\beta = \left[ \dot{\xi}_j, T_{0j}^\alpha \right] + \frac{\partial T_j^\alpha}{\partial t}|_{t=0} - \sum_{\beta=1}^p \frac{\partial c_j^{\alpha\beta}}{\partial t} T_{0j}^\beta \\
 &\Longrightarrow \frac{\partial X_j^\alpha(t)}{ \partial t}|_{t=0}= \left[ \dot{\xi}_j, T_{0j}^\alpha \right] - A_j \left( T_{0j}^\alpha \right) - \left( - \frac{\partial T_j^\alpha}{\partial t}|_{t=0} \right)
\end{align*}
This completes the proof of Proposition \ref{b18}.
\end{proof}

\subsection{Isomorphism of obstructions between \v Cech and Dolbeault resolutions of $\mathcal{E}_{\Theta_{\mathcal{F}_0}}^\bullet$ } \label{a15} \

We will explicitly define an isomorphism between the second cohomology groups $\mathbb{H}^2\left( M, \Theta_{\mathcal{F}_0}^\bullet \right)$ between \v Cech resolution and Dolbeault resolution of $\mathcal{E}_{\Theta_{\mathcal{F}_0}}^\bullet$.

Let $\left(\left\{ C_i \right\}, \left\{ B_{ij} \right\}, \left\{ a_{ijk}, G_{ijk} \right\} \right)\in C^0\left( \mathcal{U}, \mathscr{H}om_{\mathcal{O}_M}\left( \bigwedge^2 \Theta_{\mathcal{F}_0}, \frac{\Theta_M}{\Theta_{\mathcal{F}_0}} \right) \right) \bigoplus  C^1\left( \mathcal{U} , \mathscr{H}om_{\mathcal{O}_M}\left( \Theta_{\mathcal{F}_0}, \Theta_M \right) \right) \bigoplus C^2\left( \mathcal{U}, \mathcal{E}_{\Theta_{\mathcal{F}_0}} \right) $ be a $2$-cocycle. Since $-\delta\left( \left\{ \left( a_{ijk}, G_{ijk}\right) \right\} \right)=0$, we have $\left\{\left(d_{ij}, E_{ij} \right) \right\}\in C^1\left(\mathcal{U}, \mathcal{A}^{0,0}\left(\mathcal{E}_{\Theta_{\mathcal{F}_0}} \right) \right)$ such that $\delta\left( \left\{\left( d_{ij}, E_{ij} \right) \right\}\right)=\left\{ \left( a_{ijk}, G_{ijk}\right) \right\}$ where $d_{ij}\in \Gamma\left( U_{ij}, \mathcal{A}^{0,0}\left(\Theta_M \right) \right)$ and $E_{ij}\in \Gamma\left(U_{ij},  \mathcal{A}^{0,0}\left(\mathscr{H}om_{\mathcal{O}_M}\left( \Theta_{\mathcal{F}_0}, \Theta_{\mathcal{F}_0} \right) \right) \right)$ defined by
\begin{align*}
E_{ij}:\Gamma\left( U_{ij},  \Theta_{\mathcal{F}_0}\right) &\to \Gamma\left( U_{ij}, \mathcal{A}^{0,0}\left( \Theta_{\mathcal{F}_0} \right)\right)\\
 T_{0i}^\alpha &\mapsto E_{ij}\left(T_{0i}^\alpha \right)= E_{ij}^\alpha
\end{align*}
satisfying $d_{ij}-d_{ik}+ d_{jk}= a_{ijk}$ and $E_{ij}\left( T_{0i}^\alpha\right) - E_{ik}\left(T_{0i}^\alpha\right) + E_{jk}\left(T_{0i}^\alpha\right) + \sum_{\beta=1}^p \left[d_{jk}, r_{0ij}^{\alpha\beta} \right]T_{0j}^\beta= G_{ijk}\left(T_{0i}^\alpha \right)$. Since $\delta\left( \left\{\left(\bar{\partial} d_{jk} , \bar{\partial} E_{jk} \right) \right\}\right)=0$, there exists $\left\{\left( p_j, U_j \right) \right\} \in C^0\left( \mathcal{U}, \mathcal{A}^{0,1}\left(\mathcal{E}_{\Theta_{\mathcal{F}_0}}\right) \right)$ with $p_j \in \Gamma\left( U_j,  \mathcal{A}^{0,1}\left(\Theta_M \right) \right)$  and $U_j\in \Gamma\left( U_j, \mathcal{A}^{0,1}\left( \mathscr{H}om_{\mathcal{O}_M}\left( \Theta_{\mathcal{F}_0}, \Theta_{\mathcal{F}_0} \right) \right) \right)$ defined by
\begin{align*}
U_j:\Gamma\left( U_j, \Theta_{\mathcal{F}_0}\right) &\to \Gamma\left(U_j, \mathcal{A}^{0,1}\left(\Theta_{\mathcal{F}_0} \right) \right)\\
 T_{0j}^\alpha &\mapsto U_j^\alpha
\end{align*}
 such that $p_k- p_j=\bar{\partial} d_{jk}$, and $U_k\left( T_{0j}^\alpha \right)+ \sum_{\beta=1}^p \left[ p_k, r_{0jk}^{\alpha\beta} \right] T_{0k}^\beta - U_j\left(T_{0j}^\alpha \right)= \bar{\partial} E_{jk}^\alpha$. Then $\bar{\partial} p_j- \bar{\partial} p_k=0$ and $\bar{\partial} U_k\left( T_{0j}^\alpha \right)+ \sum_{\beta=1}^p \left[ \bar{\partial} p_k, r_{0jk}^{\alpha\beta} \right] T_{0k}^\beta - \bar{\partial}U_j\left(T_{0j}^\alpha \right)= 0$. This implies that $\left\{\left(\bar{\partial} p_j, \bar{\partial} U_j \right)\right\}\in A^{0,2}\left( M, \mathcal{E}_{\Theta_{\mathcal{F}_0}} \right)$.

On the second floor in the diagram in subsection \ref{ti1}, we note that
\begin{align*}
&\mathcal{L}_{d_{ij}}\left(T_{0i}^\alpha\right) - E_{ij}\left(T_{0i}^\alpha\right) - \mathcal{L}_{d_{ik}}\left(T_{0i}^\alpha\right) + E_{ik}\left(T_{0i}^\alpha\right) + \sum_{\beta=1}^p r_{0ij}^{\alpha\beta}\left(\mathcal{L}_{d_{jk}}\left(T_{0j}^\beta\right)- E_{jk}\left(T_{0j}^\beta\right)\right) -\left( B_{ij}^\alpha - B_{ik}^\alpha + \sum_{\beta=1}^p r_{0ij}^{\alpha\beta} B_{jk}^\beta \right)\\
&= \mathcal{L}_{a_{ijk}}\left(T_{0i}^\alpha\right) - \sum_{\beta=1}^p \left[ d_{jk}, r_{0ij}^{\alpha\beta}\right] T_{0j}^\beta - G_{ijk}\left(T_{0i}^\alpha\right) + \sum_{\beta=1}^p \left[ d_{ijk}, r_{0ij}^{\alpha\beta} \right] T_{0j}^\beta - \left( B_{ij}^\alpha- B_{ik}^\alpha + \sum_{\beta=1}^p r_{0ij}^{\alpha\beta} B_{jk}^\beta \right) = 0
\end{align*}

This implies that there exists $F_j\in \Gamma\left( U_j, \mathcal{A}^{0,0} \left( \mathscr{H}om_{\mathcal{O}_M}\left( \Theta_{\mathcal{F}_0}, \Theta_M \right) \right) \right)$ defined by
\begin{align*}
F_j:\Gamma\left( U_j, \Theta_{\mathcal{F}_0} \right)&\to \Gamma\left( U_j,  \mathcal{A}^{0,0}\left( \Theta_M \right) \right)\\
 T_{0j}^\alpha &\mapsto F_j^\alpha
\end{align*}
such that $F_j\left( T_{0j}^\alpha \right)- F_k\left( T_{0j}^\alpha \right)= B_{jk}^\alpha-\mathcal{L}_{d_{jk}}\left(T_{0j}^\alpha \right)+ E_{jk}^\alpha$. Then $\bar{\partial} F_j\left(T_{0j}^\alpha \right)-\bar{\partial} F_k\left(T_{0j}^\alpha \right)= - \mathcal{L}_{\bar{\partial} d_{jk}}\left(T_{0j}^\alpha \right) + \bar{\partial} E_{jk}\left(T_{0j}^\alpha \right)$ and
\begin{align}\label{a3}
&\sum_{\beta=1}^p r_{0jk}^{\alpha\beta} \left( \mathcal{L}_{p_k}\left( T_{0k}^\beta \right) - U_k^\beta - \bar{\partial} F_k^\beta \right)- \left( \mathcal{L}_{p_j}\left(T_{0j}^\alpha \right) - U_j^\alpha - \bar{\partial} F_j^\alpha \right)\\
&=\mathcal{L}_{p_k}\left( T_{0j}^\alpha \right) - \sum_{\beta=1}^p \left[ p_k, r_{0jk}^{\alpha\beta} \right] T_{0k}^\beta  - U_k\left(T_{0j}^\alpha \right) - \bar{\partial} F_k\left(T_{0j}^\alpha \right) -  \mathcal{L}_{p_j}\left(T_{0j}^\alpha \right) + U_j^\alpha + \bar{\partial} F_j^\alpha \notag \\
&=\mathcal{L}_{\bar{\partial} d_{jk}}\left( T_{0j}^\alpha \right) - \bar{\partial} E_{jk}\left(T_{0j}^\alpha \right) - \mathcal{L}_{\bar{\partial} d_{jk}}\left(T_{0j}^\alpha \right) + \bar{\partial} E_{jk}\left(T_{0j}^\alpha \right) = 0  \notag
\end{align}
Let us define $L_j\in \Gamma\left( U_j, \mathcal{A}^{0,1} \left( \mathscr{H}om_{\mathcal{O}_M}\left( \Theta_{\mathcal{F}_0}, \Theta_M \right) \right) \right)$ by
\begin{align*}
L_j:\Gamma\left( U_j, \Theta_{\mathcal{F}_0} \right) &\to \Gamma\left( U_j, \mathcal{A}^{0,1}\left(\Theta_M \right) \right)\\
T_{0j}^\alpha &\mapsto  - \mathcal{L}_{p_j}\left( T_{0j}^\alpha \right) + U_j^\alpha + \bar{\partial} F_j^\alpha
\end{align*}
Then we see from $(\ref{a3})$ that $L=\left\{L_j \right\} \in A^{0,1}\left( M, \mathscr{H}om_{\mathcal{O}_M}\left(. \Theta_{\mathcal{F}_0}, \Theta_M \right) \right)$.

Let $\tilde{C}_j\in \Gamma\left( \mathcal{U}, \mathscr{H}om_{\mathcal{O}_M}\left( \bigwedge^2 \Theta_{\mathcal{F}_0}, \Theta_M \right)\right)$ be a lifting of $C_j\in \Gamma\left(U_j, \mathscr{H}om_{\mathcal{O}_M}\left( \bigwedge^2 \Theta_{\mathcal{F}_0}, \frac{\Theta_M}{\Theta_{\mathcal{F}_0}} \right) \right)$. Then let us consider
{\Small{\begin{align*}
&C_{jk}^{\alpha\beta}:=\tilde{C}_j\left( T_{0j}^\alpha \wedge T_{0j}^\beta \right) + \left[ T_{0j}^\alpha, F_j^\beta \right] - \left[ T_{0j}^\beta , F_j^\alpha \right] - \sum_{\gamma=1}^p g_{0j\alpha\beta}^\gamma F_j^\gamma - \sum_{\delta, \eta=1}^p r_{0jk}^{\alpha\delta} r_{0jk}^{\beta \eta} \left(\tilde{C}_j\left( T_{0k}^\delta \wedge T_{0k}^\eta \right) + \left[ T_{0k}^\delta, F_k^\eta\right] - \left[ T_{0k}^\eta, F_k^\delta \right] - \sum_{\gamma=1}^p g_{0k\delta \eta}^\gamma F_k^\gamma \right)\\
&= \left( \tilde{C}_j- \tilde{C}_k \right)\left(T_{0j}^\alpha \wedge T_{0j}^\beta \right) + \left[ \sum_{\delta=1}^p r_{0jk}^{\alpha \delta} T_{0k}^\delta , \sum_{\eta=1}^p r_{0jk}^{\beta \eta} F_k^\eta - \mathcal{L}_{d_{jk}}\left( T_{0j}^\beta \right) + E_{jk}^\beta + B_{ij}^\beta  \right] \\
&- \left[\sum_{\eta=1}^p r_{0jk}^{\beta \eta} T_{0k}^\eta ,  \sum_{\delta=1}^p r_{0jk}^{\alpha \delta} F_k^\delta - \mathcal{L}_{d_{jk}}\left( T_{0j}^\alpha \right) + E_{jk}^\alpha + B_{ij}^\alpha      \right] - \sum_{\gamma=1}^p g_{0j\alpha\beta}^\gamma\left(\sum_{\eta=1}^p r_{0jk}^{\gamma \eta} F_k^\eta - \mathcal{L}_{d_{jk}}\left( T_{0j}^\gamma \right) + E_{jk}^\gamma + B_{ij}^\gamma \right) \\
&- \sum_{\delta, \eta=1}^p r_{0jk}^{\alpha\delta} r_{0jk}^{\beta \eta} \left( \left[ T_{0k}^\delta, F_k^\eta\right] - \left[ T_{0k}^\eta, F_k^\delta \right] - \sum_{\gamma=1}^p g_{0k\delta \eta}^\gamma F_k^\gamma \right)\\
&\equiv_{\Theta_{\mathcal{F}_0}}  \left[ \sum_{\delta=1}^p r_{0jk}^{\alpha \delta} T_{0k}^\delta , \sum_{\eta=1}^p r_{0jk}^{\beta \eta} F_k^\eta - \mathcal{L}_{d_{jk}}\left( T_{0j}^\beta \right) + E_{jk}^\beta  \right] - \left[\sum_{\eta=1}^p r_{0jk}^{\beta \eta} T_{0k}^\eta ,  \sum_{\delta=1}^p r_{0jk}^{\alpha \delta} F_k^\delta - \mathcal{L}_{d_{jk}}\left( T_{0j}^\alpha \right) + E_{jk}^\alpha      \right]\\
& - \sum_{\gamma=1}^p g_{0j\alpha\beta}^\gamma\left(\sum_{\eta=1}^p r_{0jk}^{\gamma \eta} F_k^\eta - \mathcal{L}_{d_{jk}}\left( T_{0j}^\gamma \right) + E_{jk}^\gamma  \right) - \sum_{\delta, \eta=1}^p r_{0jk}^{\alpha\delta} r_{0jk}^{\beta \eta} \left( \left[ T_{0k}^\delta, F_k^\eta\right] - \left[ T_{0k}^\eta, F_k^\delta \right] - \sum_{\gamma=1}^p g_{0k\delta \eta}^\gamma F_k^\gamma \right) \equiv_{ \Theta_{\mathcal{F}_0} } 0
\end{align*}}}
Then we define
\begin{align*}
C_{jk}:\Gamma\left( U_j, \bigwedge^2 \Theta_{\mathcal{F}_0} \right)&\to \Gamma\left( U_j, \mathcal{A}^{0,0}\left( \Theta_{\mathcal{F}_0} \right) \right) \\
       T_{0j}^\alpha \wedge T_{0j}^\beta &\mapsto C_{jk}^{\alpha\beta}
\end{align*}
Then $C_{ij}- C_{ik}+ C_{jk} =0$, so that there exists $Q_j\in \Gamma\left( U_j, \mathcal{A}^{0,0}\left( \mathscr{H}om_{\mathcal{O}_M}\left(\bigwedge^2 \Theta_{\mathcal{F}_0}, \Theta_{\mathcal{F}_0} \right) \right)\right)$ such that $Q_k- Q_j= C_{jk} $. Then by setting
\begin{align*}
P_j:\Gamma\left( U_j, \bigwedge^2 \Theta_{\mathcal{F}_0} \right)&\to \Gamma\left( U_j, \mathcal{A}^{0,0}\left(\Theta_M \right)\right)\\
 T_{0j}^\alpha \wedge T_{0j}^\beta &\mapsto  - \tilde{C}_j\left( T_{0j}^\alpha \wedge T_{0j}^\beta \right) - \left[ T_{0j}^\alpha, F_j^\beta \right] + \left[ T_{0j}^\beta , F_j^\alpha \right] + \sum_{\gamma=1}^p g_{0j\alpha\beta}^\gamma F_j^\gamma  - Q_j\left(T_{0j}^\alpha \wedge T_{0j}^\beta \right)
\end{align*}
we see that $P:=\left\{ P_j \right\}$ defines a global section in $A^{0,0}\left(M, \mathscr{H}om_{\mathcal{O}_M}\left( \bigwedge^2 \Theta_{\mathcal{F}_0}, \Theta_M \right) \right)$. We denote by  $\bar{P}$ the image of $P$ in $\frac{A^{0,0}\left( M, \mathscr{H}om_{\mathcal{O}_M}\left( \bigwedge^2 \Theta_{\mathcal{F}_0}, \Theta_M \right) \right)}{A^{0,0}\left( M, \mathscr{H}om_{\mathcal{O}_M}\left( \bigwedge^2 \Theta_{\mathcal{F}_0}, \Theta_{\mathcal{F}_0}  \right) \right)}$. Then we claim that 
\begin{align*}
\left(\bar{P}, L,  \left\{ \left(\bar{\partial} p_j, \bar{\partial} U_j  \right) \right\} \right) \in \frac{A^{0,0}\left( M, \mathscr{H}om_{\mathcal{O}_M}\left( \bigwedge^2 \Theta_{\mathcal{F}_0}, \Theta_M \right) \right)}{A^{0,0}\left( M, \mathscr{H}om_{\mathcal{O}_M}\left( \bigwedge^2   \Theta_{\mathcal{F}_0}, \Theta_{\mathcal{F}_0} \right) \right)} \bigoplus A^{0,1}\left( M, \mathscr{H}om_{\mathcal{O}_M}\left( \Theta_{\mathcal{F}_0} , \Theta_M \right) \right) \bigoplus A^{0,2}\left( M, \mathcal{E}_{\Theta_{\mathcal{F}_0}} \right)
\end{align*}
defines a $2$-cocycle. In fact, 
\begin{align*}
&\bar{\partial}\left( \bar{\partial} p_j, \bar{\partial} U_j \right)=0\\
&\bar{\partial}\left( -\mathcal{L}_{p_j}\left( T_{0j}^\alpha\right) + U_j\left(T_{0j}^\alpha \right) + \bar{\partial} F_j^\alpha    \right) + \mathcal{L}_{\bar{\partial} p_j}\left( T_{0j}^\alpha \right) - \bar{\partial} U_j\left( T_{0j}^\alpha \right) =0 \\
&\bar{\partial}\left(- \tilde{C}_j\left(T_{0j}^\alpha \wedge T_{0j}^\beta \right) - \left[ T_{0j}^\alpha, F_j^\beta \right]  + \left[ T_{0j}^\beta, F_j^\alpha \right]  + \sum_{\gamma=1}^p g_{0j\alpha\beta}^\gamma F_j^\gamma  \right) + D_1\left(- \mathcal{L}_{p_j}+ U_j +\bar{\partial} F_j \right)\left(T_{0j}^\alpha \wedge T_{0j}^\beta \right)\equiv 0\\
&D_2\left(- C_j + D_1\left( F_j \right) \right)\equiv 0
\end{align*}

Then we define the correspondence of $2$-cocycle between \v Cech resolution and Dolbeault resolution by 
\begin{align}\label{a6}
\left(\left\{C_j \right\}, \left\{ B_{jk}  \right\},  \left\{ \left( a_{ijk}, G_{ijk} \right) \right\} \right)\mapsto \left(\bar{P}, L,  \left\{ \left(\bar{\partial} p_j, \bar{\partial} U_j  \right) \right\} \right)
\end{align}

Conversely, let $\left(\bar{P}, L, \left( \phi, \left\{ M_j \right\} \right)      \right)$ be a $2$-cocycle. Since $\bar{\partial} \left(\phi, M_j\right)=0$, there exists $\left(p_j, U_j \right) \in \Gamma\left( U_j, \mathcal{A}^{0,1}\left(\mathcal{E}_{\Theta_{\mathcal{F}_0}} \right) \right)$ such that $\bar{\partial} \left( p_j, U_j \right)= \left(\phi, M_j \right)$. Then since $\bar{\partial}\left(\left(p_k, U_k \right) - \left( p_j, U_j \right) \right)=0$, there exists $\left(d_{jk}, E_{jk} \right)\in \Gamma\left( U_{jk}, \mathcal{A}^{0,0}\left( \mathcal{E}_{\Theta_{\mathcal{F}_0}} \right) \right)$ such that $\bar{\partial}\left( d_{jk}, E_{jk} \right)=\left( p_k, U_k \right)-\left(p_j, U_j \right)$.  Then we set
$\left(a_{ijk}, G_{ijk} \right):=\left(d_{ij}, E_{ij} \right)- \left(d_{ik}, E_{ik} \right) + \left( d_{jk}, E_{jk}\right)$. Then we have $\left\{ \left( a_{ijk}, G_{ijk}\right)\right\}\in C^2\left(\mathcal{U}, \mathcal{E}_{\Theta_{\mathcal{F}_0}} \right)$. On the other hand,  since $\bar{\partial} L\left(T_{0j}^\alpha \right)+ \mathcal{L}_{\phi}\left(T_{0j}^\alpha \right) - M_j\left(T_{0j}^\alpha \right)=0$, we have
\begin{align*}
\bar{\partial} L_j\left(T_{0j}^\alpha \right) + \mathcal{L}_{\bar{\partial}p_j}\left(T_{0j}^\alpha \right) - \bar{\partial} U_j\left(T_{0j}^\alpha \right)=0
\end{align*}
Then there exist $F_j\in \Gamma\left( U_j, \mathcal{A}^{0,0}\left( \mathscr{H}om_{\mathcal{O}_M}\left( \Theta_{\mathcal{F}_0}, \Theta_M \right) \right) \right) $ such that $\bar{\partial} F_j\left( T_{0j}^\alpha \right)= L_j^\alpha + \mathcal{L}_{p_j}\left(T_{0j}^\alpha \right) - U_j\left(T_{0j}^\alpha \right)$. Then we note that
\begin{align*}
&\bar{\partial} F_j\left( T_{0j}^\alpha \right)- \bar{\partial} F_k\left(T_{0j}^\alpha \right) =  \mathcal{L}_{p_j}\left(T_{0j}^\alpha \right) - U_j\left(T_{0j}^\alpha \right) - \sum_{\beta=1}^p r_{0jk}^{\alpha\beta} \left(\mathcal{L}_{p_k}\left( T_{0k}^\beta \right)  - U_k\left(T_{0k}^\beta \right)   \right)\\
&= \mathcal{L}_{p_j- p_k}\left(T_{0j}^\alpha \right) - U_j\left(T_{0j}^\alpha \right) + \sum_{\beta=1}^p \left[ p_k, r_{0jk}^{\alpha\beta} \right] T_{0k}^\beta + U_k\left( T_{0j}^\alpha \right) = \mathcal{L}_{-\bar{\partial}d_{jk}}\left(T_{0j}^\alpha \right) + \bar{\partial}E_{jk}\left( T_{0j}^\alpha \right)
\end{align*}
Then we define 
\begin{align*}
B_{jk}:\Gamma\left( U_{jk}, \Theta_{\mathcal{F}_0} \right) &\to \Gamma\left( U_{jk}, \Theta_M \right) \\
 T_{0j}^\alpha &\mapsto  F_j\left(T_{0j}^\alpha \right) - F_k\left( T_{0j}^\alpha \right) + \mathcal{L}_{d_{jk}}\left( T_{0j}^\alpha \right) - E_{jk}\left(T_{0j}^\alpha \right)
\end{align*}
Then  $\left\{ B_{jk}\right\} \in C^1 \left( \mathcal{U}, \mathscr{H}om_{\mathcal{O}_M}\left( \Theta_{\mathcal{F}_0}, \Theta_M \right) \right)$.

Since $\bar{\partial} P\left(T_{0j}^\alpha \wedge T_{0j}^\beta \right)+ \left[ T_{0j}^\alpha, L_j^\alpha \right] - \left[ T_{0j}^\beta, L_j^\alpha\right]- \sum_{\gamma=1}^p g_{0j\alpha\beta}^\gamma L_j^\gamma = \sum_{\eta=1}^p R_j^{\alpha\eta} T_{0j}^\eta \in \Gamma\left( U_j, \mathcal{A}^{0,1}\left(\Theta_{\mathcal{F}_0} \right) \right)$ with some $R_j^{\alpha\beta} \in \Gamma\left( U_j, \mathcal{A}^{0,1} \right)$, we note that
{\Small{\begin{align*}
&\bar{\partial} \left(\sum_{\eta=1}^p R_j^{\alpha \eta} T_{0j}^\eta \right) = \left[ T_{0j}^\alpha,  \mathcal{L}_{-\bar{\partial} p_j}\left(T_{0j}^\beta \right) + \sum_{\eta=1}^p \bar{\partial} U_j^{\beta \eta} T_{0j}^\eta \right]  - \left[ T_{0j}^\beta, - \mathcal{L}_{\bar{\partial} p_j}\left(T_{0j}^\alpha \right) + \sum_{\eta=1}^p \bar{\partial} U_j^{\alpha \eta} T_{0j}^\eta \right]  - \sum_{\gamma=1}^p g_{0j\alpha\beta}^\gamma \left( - \mathcal{L}_{\bar{\partial} p_j}\left(T_{0j}^\alpha \right) + \sum_{\eta=1}^p \bar{\partial} U_j^{\alpha \eta} T_{0j}^\eta \right)
\end{align*}}}
implies that
{\Small{\begin{align*}
& \bar{\partial}\left(\left[ T_{0j}^\alpha,  -\mathcal{L}_{p_j}\left(T_{0j}^\beta \right) + \sum_{\eta=1}^p U_j^{\beta \eta} T_{0j}^\eta \right]  - \left[ T_{0j}^\beta, - \mathcal{L}_{ p_j}\left(T_{0j}^\alpha \right) + \sum_{\eta=1}^p U_j^{\alpha \eta} T_{0j}^\eta \right]  - \sum_{\gamma=1}^p g_{0j\alpha\beta}^\gamma \left( - \mathcal{L}_{  p_j}\left(T_{0j}^\alpha \right) + \sum_{\eta=1}^p  U_j^{\alpha \eta} T_{0j}^\eta \right) - \left(\sum_{\eta=1}^p R_j^{\alpha \eta} T_{0j}^\eta \right) \right)= 0
\end{align*}}}
Then there exists $D_j^{\alpha\beta} \in \Gamma\left(U_j, \mathcal{A}^{0,0}\left( \Theta_{\mathcal{F}_0} \right) \right)$ such that
{\Small{\begin{align*}
\bar{\partial} D_j^{\alpha\beta} = \left[ T_{0j}^\alpha,  -\mathcal{L}_{p_j}\left(T_{0j}^\beta \right) + \sum_{\eta=1}^p U_j^{\beta \eta} T_{0j}^\eta \right]  - \left[ T_{0j}^\beta, - \mathcal{L}_{ p_j}\left(T_{0j}^\alpha \right) + \sum_{\eta=1}^p U_j^{\alpha \eta} T_{0j}^\eta \right]  - \sum_{\gamma=1}^p g_{0j\alpha\beta}^\gamma \left( - \mathcal{L}_{  p_j}\left(T_{0j}^\alpha \right) + \sum_{\eta=1}^p  U_j^{\alpha \eta} T_{0j}^\eta \right) - \left(\sum_{\eta=1}^p R_j^{\alpha \eta} T_{0j}^\eta \right)
\end{align*}}}
and we note that
\begin{align*}
0&=\bar{\partial} P\left(T_{0j}^\alpha \wedge T_{0j}^\beta \right)+ \left[ T_{0j}^\alpha, L_j^\alpha \right] - \left[ T_{0j}^\beta, L_j^\alpha\right]- \sum_{\beta=1}^p g_{0j\alpha\beta}^\gamma L_j^\gamma -\sum_{\eta=1}^p R_j^{\alpha\eta} T_{0j}^\eta\\
&= \bar{\partial} P\left(T_{0j}^\alpha \wedge T_{0j}^\beta \right) + \left[ T_{0j}^\alpha, \bar{\partial} F_j^\beta - \left[ p_j, T_{0j}^\beta \right]  + \sum_{\eta=1}^p U_j^{\beta \eta} T_{0j}^\eta \right] -\left[ T_{0j}^\beta, \bar{\partial} F_j^\alpha - \left[ p_j, T_{0j}^\alpha \right]  + \sum_{\eta=1}^p U_j^{\alpha \eta} T_{0j}^\eta \right]\\
&- \sum_{\gamma=1}^p g_{0j\alpha\beta}^\gamma \left(  \bar{\partial} F_j^\gamma - \left[ p_j, T_{0j}^\gamma \right]  +  \sum_{\eta=1}^p U_j^{\gamma \eta} T_{0j}^\eta \right) -\sum_{\eta=1}^p R_j^{\alpha\eta} T_{0j}^\eta\\
&= \bar{\partial} P\left( T_{0j}^\alpha \wedge T_{0j}^\beta \right) + \bar{\partial} \left(\left[ T_{0j}^\alpha, F_j^\beta \right]  - \left[ T_{0j}^\beta, F_j^\alpha \right] - \sum_{\gamma=1}^p g_{0j\alpha\beta}^\gamma F_j^\gamma  \right) - \bar{\partial} D_j^{\alpha\beta}
\end{align*}

Then we define $C_j\in \Gamma\left( U_j, \mathscr{H}om_{\mathcal{O}_M}\left( \bigwedge^2 \Theta_{\mathcal{F}_0}, \frac{\Theta_M}{\Theta_{\mathcal{F}_0}} \right) \right)$
\begin{align*}
C_j:\Gamma\left( U_j, \bigwedge^2 \Theta_{\mathcal{F}_0} \right) &\to \Gamma\left( U_j, \frac{\Theta_M}{\Theta_{\mathcal{F}_0}} \right) \\
          T_{0j}^\alpha\wedge T_{0j}^\beta &\mapsto \overline{- P\left(T_{0j}^\alpha\wedge T_{0j}^\beta \right) - \left[ T_{0j}^\alpha, F_j^\beta \right] +\left[ T_{0j}^\beta , F_j^\alpha \right] + \sum_{\gamma=1}^p g_{0j\alpha\beta}^\gamma F_j^\gamma + D_j^{\alpha\beta}}
\end{align*}
Then 
\begin{align*}
\left( \left\{ C_j \right\} , \left\{ B_{jk} \right\}, \left\{ \left( a_{ijk}, G_{ijk} \right) \right\} \right) \in  C^0\left(  \mathcal{U}, \mathscr{H}om_{\mathcal{O}_M}\left( \bigwedge^2 \Theta_{\mathcal{F}_0}, \frac{\Theta_M}{\Theta_{\mathcal{F}_0}} \right) \right) \bigoplus C^1\left( \mathcal{U}, \mathscr{H}om_{\mathcal{O}_M} \left(  \Theta_{\mathcal{F}_0}, \Theta_M  \right) \right) \bigoplus C^2\left( \mathcal{U}, \mathcal{E}_{\Theta_{\mathcal{F}_0}} \right)
\end{align*}
defines a $2$-cocycle and we define the inverse correspondence by $\left( P, L, \left(\left\{ \phi, M_j \right\} \right) \right)\mapsto \left(\left\{C_j\right\}, \left\{ B_{jk}\right\}, \left\{\left( a_{ijk}, G_{ijk} \right) \right\} \right)$.

\begin{proposition}\label{a13}

We recall from \textnormal{Lemma \ref{o54}} that given a $1$-cocycle $\left(\left\{ \alpha_j \right\}, \left\{\left(\theta_{jk}, A_{jk} \right) \right\} \right)$ in $\mathbb{H}^1 \left( M, \mathcal{E}_{\Theta_{\mathcal{F}_0}}^\bullet \right)$, we have the associated obstruction class from $(\ref{it1})$
\begin{align}\label{a4}
\left( \left\{ R_j\right\}, \left\{F_{jk} \right\}, \left\{ \left[  \left(\theta_{ij}, A_{ij}  \right) , \left(\theta_{jk}, A_{jk} \right)   \right]  \right\}\right)  \in \mathbb{H}^2\left( M, \mathcal{E}_{\Theta_{\mathcal{F}_0}}^\bullet \right)
\end{align}
from the \v Cech resolution of $ \mathcal{E}_{ \Theta_{\mathcal{F}_0}}^\bullet$. Assume that the $1$-cocycle $\left( \left\{ \alpha_j \right\}, \left\{ \left( \theta_{jk}, A_{jk} \right) \right\} \right)$ from the \v Cech resolution of $\mathcal{E}_{\Theta_{\mathcal{F}_0}}^\bullet$ corresponds to the $1$-cocycle $\left( \Phi=\left\{ \Phi_j\right\}, \left\{\left( \psi, a_j \right) \right\} \right)$ from the Dolbeault reslolution of $\mathcal{E}_{\Theta_{\mathcal{F}_0}}^\bullet$ under the isomorphism $(\ref{a1})$.  We recall from \textnormal{Lemma \ref{its30}} that we have the associated obstruction class from $(\ref{it10})$
\begin{align}\label{a5}
\left( \rho, \eta, \left( -\left[ \left(\psi, a_j \right), \left( \psi, a_j \right) \right] \right) \right) \in \mathbb{H}^2\left( M, \mathcal{E}_{\Theta_{\mathcal{F}_0}}^\bullet \right)
\end{align}
from the Dolbeault resolution of $\mathcal{E}_{\Theta_{\mathcal{F}_0}}^\bullet$. Then $(\ref{a4})$ corresponds to $(\ref{a5})$ under the isomorphism $(\ref{a6})$.
\end{proposition}

\begin{proof}
Let $\left( \left\{ \alpha_j \right\}, \left\{ \left( \theta_{jk}, A_{jk} \right) \right\} \right)$ be a $1$-cocycle in Lemma \ref{o54}. Then $\left( \left\{ R_j\right\}, \left\{F_{jk} \right\}, \left\{ \left[  \left(\theta_{ij}, A_{ij}  \right) , \left(\theta_{jk}, A_{jk} \right)   \right]  \right\}\right) $ defines a $2$-cocycle by Lemma \ref{o54}. We have the corresponding $1$-cocycle $(\Phi, \left\{ \left( \psi, a_j:=\bar{\partial} A_j \right) \right\})$ from the Dolbeault resolution of $\mathcal{E}_{\Theta_{\mathcal{F}_0}}^\bullet$ from $(\ref{a1})$. Since $-\delta\left( \left\{ \left(\xi_j, A_j \right) \right\}\right)= \left\{ \left(\theta_{jk}, A_{jk} \right) \right\}$, i.e $\left(\xi_i, A_i \right)-\left( \xi_j, A_j\right)=\left( \theta_{jk}, A_{jk} \right)$, we have
{\small{\begin{align*}
&\left[\left( \xi_i, A_i \right), \left( \xi_j, A_j \right) \right] - \left[\left( \xi_i, A_i\right), \left( \xi_k, A_k \right) \right] + \left[\left( \xi_j, A_j \right), \left( \xi_k , A_k \right) \right]= \left[ \left(\xi_i, A_i \right), \left(\theta_{jk}, A_{jk} \right) \right] + \left[ \left(\xi_j, A_j \right), \left( \xi_k, A_k \right)- \left( \xi_j, A_j \right) \right]\\
&=\left[\left( \xi_i, B_i \right), \left( \theta_{jk}, B_{jk} \right) \right] - \left[\left( \xi_j, B_j \right), \left( \theta_{jk}, B_{jk} \right) \right] = \left[ \left( \theta_{ij}, B_{ij} \right), \left( \theta_{jk} , B_{jk} \right) \right]
\end{align*}}}
This implies that
\begin{align*}
\delta\left( \left\{\left[ \left(\xi_i, A_i\right), \left(\xi_j, A_j \right)\right]  \right\} \right) = \left\{ \left[ \left(\theta_{ij}, A_{ij}\right) , \left(\theta_{jk}, A_{jk} \right) \right]   \right]
\end{align*}
Then we set $\left( d_{jk}, E_{jk} \right):=\left[ \left( \xi_j, A_j \right), \left(\xi_k, A_k \right) \right]$ and we note that
\begin{align*}
\bar{\partial}\left[ \left(\xi_j, A_j\right), \left(\xi_k, A_k \right)\right]  =\left[ \bar{\partial}\left(\xi_j, A_j \right), \left( \xi_k, A_k \right) \right] + \left[ \left(\xi_j, A_j \right) , \bar{\partial} \left(\xi_k, A_k \right) \right] = \left[\left(\psi, a_j \right), \left(\xi_k, A_k \right) \right] - \left[\left(\psi , a_k \right), \left( \xi_j, A_j\right) \right]
\end{align*}
Then we take $\left( p_j, U_j\right):=\left[\left( \psi, a_j \right), \left(\xi_j, A_j \right) \right]$, so that $\bar{\partial}\left( p_j, U_j\right)= -\left[ \left( \psi, a_j \right), \left(\psi, a_j \right)  \right]$.

On the other hand, we note that
\begin{align*}
A_j\left( T_{0j}^\alpha \right)- A_k^j\left(T_{0j}^\alpha \right) = A_j\left(T_{0j}^\alpha \right) - \left( A_k\left( T_{0j}^\alpha \right) + \sum_{\gamma=1}^p \left[ \xi_k, r_{0jk}^{\alpha \gamma}\right] T_{0k}^\gamma  \right)= A_{jk} \left( T_{0j}^\alpha \right)= A_{jk}^\alpha = \sum_{\delta=1}^p A_{jk}^{\alpha \delta} T_{0j}^\delta
\end{align*}
This implies that 
\begin{align}\label{a7}
A_{jk}^{\alpha \delta} = A_j^{\alpha \delta} - \sum_{\beta, \eta=1}^p r_{0jk}^{\alpha \beta} A_k^{\beta\eta} r_{0kj}^{\eta \delta} - \sum_{\gamma, \eta=1}^p \left[ \xi_k, r_{0jk}^{\alpha \gamma}\right] r_{0kj}^{\gamma \delta}
\end{align}

We note that $\left( d_{jk}, E_{jk} \right)=\left[\left( \xi_j, A_j\right), \left(\xi_k, A_k \right) \right]=\left( \left[ \xi_j, \xi_k \right], \left[ \xi_j, A_k^j\right] - \left[ \xi_k, A_j\right] + A_j \circ A_k^j - A_k^j \circ A_j  \right)$. Then we have
\begin{align}\label{a9}
&- \mathcal{L}_{d_{jk}}\left( T_{0j}^\alpha \right)+ E_{jk}\left( T_{0j}^\alpha \right) =-\mathcal{L}_{\left[ \xi_j, \xi_k \right]}\left( T_{0j}^\alpha \right) +\left(\left[ \xi_j, A_k^j\right] - \left[ \xi_k, A_j\right] + A_j \circ A_k^j - A_k^j \circ A_j \right) \left( T_{0j}^\alpha \right)\\
&=-\mathcal{L}_{\xi_j}\mathcal{L}_{\xi_k}\left( T_{0j}^\alpha \right) + \mathcal{L}_{\xi_k} \mathcal{L}_{\xi_j}\left( T_{0j}^\alpha \right) + \sum_{\delta=1}^p \left[ \xi_j,  \sum_{\beta, \eta=1}^p r_{0jk}^{\alpha\beta} A_k^{\beta \eta} r_{0kj}^{\eta \delta} + \sum_{\eta=1}^p \left[\xi_k, r_{0jk}^{\alpha \eta} \right]r_{0kj}^{\eta \delta} \right] T_{0j}^\delta - \sum_{\delta=1}^p \left[ \xi_k, A_j^{\alpha \delta} \right] T_{0j}^\delta \notag \\
&+ \sum_{\xi=1}^p \left(\sum_{\beta, \eta=1}^p r_{0jk}^{\alpha\beta} A_k^{\beta \eta} r_{0kj}^{\eta \delta} + \sum_{\eta=1}^p \left[\xi_k, r_{0jk}^{\alpha \eta} \right]r_{0kj}^{\eta \delta}\right)  A_j^{\delta \xi} T_{0j}^\xi - \sum_{\delta, \xi =1}^p A_j^{\alpha \xi} \left( \sum_{\beta, \eta=1}^q r_{0jk}^{\xi\beta} A_k^{\beta \eta} r_{0kj}^{\eta \delta} + \sum_{\eta=1}^p \left[\xi_k, r_{0jk}^{\xi \eta} \right]r_{0kj}^{\eta \delta} \right) T_{0j}^\delta \notag \\
&= -\mathcal{L}_{\xi_j}\mathcal{L}_{\xi_k}\left( T_{0j}^\alpha \right) + \mathcal{L}_{\xi_k} \mathcal{L}_{\xi_j}\left( T_{0j}^\alpha \right) + \sum_{\beta, \eta=1}^p \mathcal{L}_{\xi_j} \left( r_{0jk}^{\alpha \beta} A_k^{\beta \eta} T_{0k}^\eta \right)  - \sum_{\delta, \beta, \eta=1}^p r_{0jk}^{\alpha \beta} A_k^{\beta \eta} r_{0kj}^{\eta \delta} \mathcal{L}_{\xi_j}\left( T_{0j}^\delta \right) + \sum_{\eta=1}^p \mathcal{L}_{\xi_j}\left( \left[ \xi_k, r_{0jk}^{\alpha \eta} \right] T_{0k}^\eta \right) \notag \\
&-\sum_{\eta, \delta=1}^p \left[ \xi_k, r_{0jk}^{\alpha \eta}\right] r_{0kj}^{\eta \delta}\mathcal{L}_{\xi_j}\left(T_{0j}^\delta \right)- \sum_{\delta=1}^p \left[ \xi_k, A_j^{\alpha \delta} \right] T_{0j}^\delta + \sum_{\xi, \beta, \eta=1}^p r_{0jk}^{\alpha\beta} A_k^{\beta \eta} r_{0kj}^{\eta \delta} A_j^{\delta \xi} T_{0j}^\xi + \sum_{\xi, \eta=1}^p \left[ \xi_k, r_{0jk}^{\alpha \eta} \right] r_{0kj}^{\eta \delta} A_j^{\delta \xi} T_{0j}^\xi \notag \\
&-\sum_{\beta, \eta, \xi=1}^p A_j^{\alpha \xi} r_{0jk}^{\xi\beta} A_k^{\beta \eta} T_{0k}^\eta - \sum_{\delta, \xi, \eta=1}^p A_j^{\alpha \xi} \left[ \xi_k, r_{0jk}^{\xi\eta}\right] T_{0k}^\eta \notag
\end{align}

On the other hand, we note that $F_{jk}\left(T_{0j}^\alpha \right)= - \alpha_j\left( A_{jk}^\alpha \right) - \alpha_k \left( A_{jk}^\alpha \right) + \mathcal{L}_{\theta_{jk}}\left( \alpha_j\left(T_{0j}^\alpha \right) \right) + \mathcal{L}_{\theta_{jk}}\left( \alpha_k \left( T_{0j}^\alpha \right)  \right)$. Then from $(\ref{a7})$

{\Small{\begin{align*}
&F_{jk}\left( T_{0j}^\alpha \right)=- \alpha_j\left( A_{jk}^\alpha \right) - \alpha_k \left( A_{jk}^\alpha \right) + \mathcal{L}_{\theta_{jk}}\left( \alpha_j\left(T_{0j}^\alpha \right) \right) + \mathcal{L}_{\theta_{jk}}\left( \alpha_k \left( T_{0j}^\alpha \right)  \right)\\
&=-\sum_{\delta=1}^p \left(A_j^{\alpha \delta} - \sum_{\beta, \eta=1}^q r_{0jk}^{\alpha\beta } A_k^{\beta\eta} r_{0kj}^{\eta \delta} - \sum_{\gamma, \eta=1}^p \left[ \xi_k, r_{0jk}^{\alpha \gamma}\right] r_{0kj}^{\gamma \delta}\right) \alpha_j\left( T_{0j}^\delta \right) - \sum_{\delta=1}^p \left( A_j^{\alpha \delta} - \sum_{\beta, \eta=1}^p r_{0jk}^{\alpha \beta} A_k^{\beta \eta} r_{0kj}^{\eta \delta} - \sum_{\gamma, \eta=1}^p \left[ \xi_k, r_{0jk}^{\alpha \gamma}\right] r_{0kj}^{\gamma \delta}\right) \alpha_k\left(T_{0j}^\delta \right)  \\
&+ \mathcal{L}_{\xi_j}\left(\alpha_j \left( T_{0j}^\alpha \right) \right) - \mathcal{L}_{\xi_k}\left(\alpha_j \left( T_{0j}^\alpha \right) \right) + \mathcal{L}_{\xi_j}\left(\alpha_k \left( T_{0j}^\alpha \right) \right) - \mathcal{L}_{\xi_k}\left(\alpha_k \left( T_{0j}^\alpha \right) \right) \\
&= - \alpha_j\left( A_j^\alpha \right)  + \sum_{\delta, \beta, \eta=1}^p r_{0jk}^{\alpha\beta } A_k^{\beta \eta} r_{0kj}^{\eta \delta}\left(\alpha_k\left( T_{0j}^\delta \right) + \mathcal{L}_{\theta_{jk}}\left(T_{0j}^\delta \right) - A_{jk}^\delta \right) + \sum_{\delta, \gamma, \eta=1}^p \left[ \xi_k, r_{0jk}^{\alpha \gamma } \right] r_{0kj}^{\gamma \delta} \left( \alpha_k\left( T_{0j}^\delta \right) + \mathcal{L}_{\theta_{jk}}\left(T_{0j}^\delta \right) - A_{jk}^\delta \right)  \\
&- \sum_{\delta=1}^p A_j^{\alpha \delta}\left(\alpha_j \left( T_{0j}^\delta \right) - \mathcal{L}_{\theta_{jk}} \left( T_{0j}^\delta \right) + A_{jk}^\delta \right) + \sum_{\beta=1}^p r_{0jk}^{\alpha \beta}\alpha_k\left(A_k^\beta \right) + \sum_{\gamma, \eta=1}^p \left[ \xi_k , r_{0jk}^{\alpha \gamma} \right] \alpha_k \left( T_{0k}^\gamma \right)  \\
&+ \mathcal{L}_{\xi_j}\left( \alpha_j\left( T_{0j}^\alpha \right) \right) - \mathcal{L}_{\xi_k}\left( \alpha_k\left( T_{0j}^\alpha \right) + \mathcal{L}_{\theta_{jk}}\left(T_{0j}^\alpha \right) - A_{jk}^\alpha \right) + \mathcal{L}_{\xi_j}\left( \alpha_j\left(T_{0j}^\alpha \right) - \mathcal{L}_{\theta_{jk}}\left(T_{0j}^\alpha \right) + A_{jk}^\alpha \right)  - \mathcal{L}_{\xi_k}\left(\alpha_k \left( T_{0j}^\alpha \right) \right)
\end{align*}}}

We continue the computation. Then we have
{\tiny{\begin{align}\label{a8}
&F_{jk}\left( T_{0j}^\alpha \right)=- 2\alpha_j\left( A_j^\alpha \right) + 2\sum_{\beta=1}^p r_{0jk}^{\alpha\beta} \alpha_k\left( A_k^\beta \right) + \sum_{\delta, \beta, \eta=1}^p r_{0jk}^{\alpha\beta} A_k^{\beta \eta} r_{0kj}^{\eta \delta} \mathcal{L}_{\xi_j- \xi_k}\left( T_{0j}^\delta \right)  \\
&- \sum_{\delta, \beta, \eta,\xi=1}^p r_{0jk}^{\alpha\beta} A_k^{\beta \eta} r_{0kj}^{\eta \delta}\left(A_j^{\delta \xi} - \sum_{ \sigma,\gamma=1}^p r_{0jk}^{\delta \sigma} A_k^{\sigma \gamma } r_{0kj}^{\gamma  \xi} - \sum_{\gamma =1}^p \left[ \xi_k, r_{0jk}^{\delta \gamma}\right] r_{0kj}^{\gamma \xi } \right) T_{0j}^\xi +2 \sum_{\gamma, \eta=1}^p \left[ \xi_k, r_{0jk}^{\alpha \gamma} \right] \alpha_k \left(T_{0k}^\gamma \right) \notag \\
&+ \sum_{\delta, \gamma, \eta=1}^p \left[ \xi_k, r_{0jk}^{\alpha \gamma}\right] r_{0kj}^{\gamma \delta} \mathcal{L}_{\xi_j-\xi_k}\left( T_{0j}^\delta \right) - \sum_{\delta, \eta,\xi=1}^p \left[ \xi_k, r_{0jk}^{\alpha \eta} \right] r_{0kj}^{\eta \delta} \left(A_j^{\delta \xi} - \sum_{ \sigma,\gamma=1}^q r_{0jk}^{\delta \sigma} A_k^{\sigma \gamma }r_{0kj}^{\gamma  \xi} - \sum_{\gamma =1}^p \left[ \xi_k, r_{0jk}^{\delta \gamma}\right] r_{0kj}^{\gamma \xi } \right) T_{0j}^\xi \notag \\
&+ \sum_{\delta=1}^p A_j^{\alpha \delta} \mathcal{L}_{\xi_j-\xi_k}\left( T_{0j}^\delta \right) - \sum_{\delta,\xi=1}^p A_j^{\alpha \delta}\left(A_j^{\delta \xi} - \sum_{ \sigma,\gamma,\xi=1}^q r_{0jk}^{\delta \sigma} A_k^{\sigma \gamma }r_{0kj}^{\gamma  \xi} - \sum_{\gamma =1}^p \left[ \xi_k, r_{0jk}^{\delta \gamma}\right] r_{0kj}^{\gamma \xi } \right) T_{0j}^\xi \notag \\
&+ 2 \mathcal{L}_{\xi_j}\left(\alpha_j\left( T_{0j}^\alpha \right) \right) - 2 \mathcal{L}_{\xi_k}\left(\alpha_k \left(T_{0j}^\alpha \right) \right) - \mathcal{L}_{\xi_k}\mathcal{L}_{\xi_j-\xi_k}\left(T_{0j}^\alpha \right) + \mathcal{L}_{\xi_k}\left( \sum_{\delta=1}^p A_j^{\alpha\delta }T_{0j}^\delta - \sum_{\beta, \eta=1}^p r_{0jk}^{\alpha \beta} A_k^{\beta\eta} T_{0k}^\eta - \sum_{\gamma, \eta=1}^p \left[ \xi_k, r_{0jk}^{\alpha \gamma}\right] T_{0k}^\eta \right) \notag \\
&- \mathcal{L}_{\xi_j}\mathcal{L}_{\xi_j-\xi_k}\left( T_{0j}^\alpha \right) + \mathcal{L}_{\xi_j}\left( \sum_{\delta=1}^p A_j^{\alpha \delta } T_{0j}^\delta - \sum_{\beta, \eta=1}^p r_{0jk}^{\alpha \beta} A_k^{\beta\eta} T_{0k}^\eta - \sum_{\gamma, \eta=1}^p \left[ \xi_k, r_{0jk}^{\alpha \gamma}\right] T_{0k}^\eta  \right) \notag
\end{align}}} 

Then from $(\ref{a9})$ and $(\ref{a8})$ we have
{\tiny{\begin{align*}
& -\mathcal{L}_{d_{jk}}\left( T_{0j}^\alpha \right) + E_{jk}\left( T_{0j}^\alpha \right) + F_{jk}\left( T_{0j}^\alpha \right)  \\
&=-\mathcal{L}_{\left[ \xi_j, \xi_k \right]}\left( T_{0j}^\alpha \right) +\left(\left[ \xi_j, A_k^j\right] - \left[ \xi_k, A_j\right] + A_j \circ A_k^j - A_k^j \circ A_j \right) \left( T_{0j}^\alpha \right) - \alpha_j \left( A_{jk}^\alpha \right) - \alpha_k \left( A_{jk}^\alpha \right) + \mathcal{L}_{\theta_{jk}}\left( \alpha_j \left(T_{0j}^\alpha \right) \right) + \mathcal{L}_{\theta_{jk}}\left( \alpha_k \left( T_{0j}^\alpha \right)  \right)\\
&= -\cancel{\mathcal{L}_{\xi_j}\mathcal{L}_{\xi_k}\left( w_{0j}^\alpha \right) }+ \cancel{\mathcal{L}_{\xi_k} \mathcal{L}_{\xi_j}\left( T_{0j}^\alpha \right) }+ \cancel{\sum_{\beta, \eta=1}^p \mathcal{L}_{\xi_j} \left( r_{0jk}^{\alpha \beta} A_k^{\beta \eta} T_{0k}^\eta \right) }  - \cancel{\sum_{\delta, \beta, \eta=1}^p r_{0jk}^{\alpha \beta} A_k^{\beta \eta} r_{0kj}^{\eta \delta} \mathcal{L}_{\xi_j}\left( T_{0j}^\delta \right)} + \cancel{\sum_{\eta=1}^p \mathcal{L}_{\xi_j}\left( \left[ \xi_k, r_{0jk}^{\alpha \eta} \right] T_{0k}^\eta \right) }\\
&-\bcancel{\sum_{\eta, \delta=1}^p \left[ \xi_k, r_{0jk}^{\alpha \eta}\right] r_{0kj}^{\eta \delta}\mathcal{L}_{\xi_j}\left(T_{0j}^\delta \right)}- \bcancel{\sum_{\delta=1}^p \left[ \xi_k, A_j^{\alpha \delta} \right] T_{0j}^\delta }+ \bcancel{ \sum_{\xi, \beta, \eta=1}^p r_{0jk}^{\alpha\beta} A_k^{\beta \eta} r_{0kj}^{\eta \delta} A_j^{\delta \xi} T_{0j}^\xi} + \bcancel{ \sum_{\xi, \eta=1}^p \left[ \xi_k, r_{0jk}^{\alpha \eta} \right] r_{0kj}^{\eta \delta} A_j^{\delta \xi} T_{0j}^\xi }\\
&- \bcancel{\sum_{\beta, \eta, \xi=1}^p A_j^{\alpha \xi} r_{0jk}^{\xi\beta} A_k^{\beta \eta} T_{0k}^\eta} - \bcancel{ \sum_{\delta, \xi, \eta=1}^p A_j^{\alpha \xi} \left[ \xi_k, r_{0jk}^{\xi\eta}\right] T_{0k}^\eta}\\
&- 2\alpha_j \left( A_j^\alpha \right) + 2\sum_{\beta=1}^p r_{0jk}^{\alpha\beta} \alpha_k \left( A_k^\beta \right) + \cancel{\sum_{\delta, \beta, \eta=1}^p r_{0jk}^{\alpha\beta} A_k^{\beta \eta} r_{0kj}^{\eta \delta} \mathcal{L}_{\xi_j}\left( T_{0j}^\delta \right)} - \xcancel{\sum_{\delta, \beta, \eta=1}^p r_{0jk}^{\alpha \beta} A_k^{\beta \eta} r_{0kj}^{\eta \delta} \mathcal{L}_{\xi_k} \left( T_{0j}^\delta \right)}\\
&-\bcancel{ \sum_{\delta, \beta, \eta, \xi=1}^p r_{0jk}^{\alpha\beta} A_k^{\beta \eta} r_{0kj}^{\eta \delta} A_j^{\delta \xi} T_{0j}^\xi }+ \sum_{\beta,  \eta=1}^p r_{0jk}^{\alpha\beta} A_k^{\beta \eta} A_k^{\eta \gamma} T_{0k}^\gamma + \xcancel{\sum_{\beta, \eta, \delta =1}^p  r_{0jk}^{\alpha\beta} A_k^{\beta \eta} r_{0kj}^{\eta \delta} \mathcal{L}_{\xi_k}\left(T_{0j}^\delta \right) }- \sum_{\beta, \eta =1}^p r_{0jk}^{\alpha\beta} A_k^{\beta \eta} \mathcal{L}_{\xi_k}\left(T_{0k}^\eta \right)\\
& - 2\sum_{\gamma ,\eta=1}^p r_{0jk}^{\alpha \gamma} \left[ \xi_k, \alpha_k\left(T_{0k}^\gamma \right) \right]  + \bcancel{ \sum_{\delta, \gamma, \eta=1}^p \left[ \xi_k, r_{0jk}^{\alpha \gamma} \right] r_{0kj}^{\gamma \delta} \mathcal{L}_{\xi_j}\left(T_{0j}^\delta \right)} -  \xcancel{\sum_{\delta, \gamma, \eta=1}^p \left[ \xi_k, r_{0jk}^{\alpha \gamma} \right] r_{0kj}^{\gamma \delta} \mathcal{L}_{\xi_k}\left(T_{0j}^\delta \right) } - \bcancel{ \sum_{\delta, \eta, \xi=1}^p \left[ \xi_k, r_{0jk}^{\alpha \eta} \right] r_{0kj}^{\eta \delta} A_j^{\delta\xi} T_{0j}^\xi } \\
& + \xcancel{\sum_{ \eta,\gamma=1}^p \mathcal{L}_{\xi_k}\left( r_{0jk}^{\alpha \eta} A_k^{\eta \gamma} T_{0k}^\gamma \right) } - \sum_{\eta, \gamma=1}^p r_{0jk}^{\alpha \eta} \mathcal{L}_{\xi_k}\left( A_k^{\eta \gamma} T_{0k}^\gamma  \right) - \xcancel{ \sum_{\delta, \eta=1}^p \left[ \xi_k, r_{0jk}^{\alpha \eta} \right]\left[ \xi_k, r_{0kj}^{\eta \delta} \right] T_{0j}^\delta} + \sum_{\delta=1}^p A_j^{\alpha \delta} \mathcal{L}_{\xi_j}\left(T_{0j}^\delta \right) -  \cancel{ \sum_{\delta=1}^p A_j^{\alpha \delta} \mathcal{L}_{\xi_k}\left(T_{0j}^\delta \right)} \\
&-\sum_{\delta, \xi=1}^p A_j^{\alpha \delta} A_j^{\delta \xi} T_{0j}^\xi + \bcancel{\sum_{\delta, \sigma, \gamma=1}^p A_j^{\alpha \delta} r_{0jk}^{\delta \sigma} A_k^{\sigma \gamma} T_{0k}^\gamma } + \bcancel{\sum_{\gamma, \delta=1}^p A_j^{\alpha \delta} \left[ \xi_k, r_{0jk}^{\delta \gamma} \right] T_{0k}^\gamma  }+ 2\mathcal{L}_{\xi_j}\left(\alpha_j \left( T_{0j}^\alpha \right) \right) - \cancel{\mathcal{L}_{\xi_k} \mathcal{L}_{\xi_j} \left( T_{0j}^\alpha \right) } \\
& + \cancel{ \sum_{\beta=1}^p \mathcal{L}_{\xi_k}\left( \left[ \xi_k, r_{0jk}^{\alpha\beta}\right] T_{0k}^\beta \right) } +\xcancel{ \sum_{\beta,\eta=1}^p \left[ \xi_k, r_{0jk}^{\alpha\beta}\right]  \left[ \xi_k,  r_{0kj}^{\beta \eta}\right] T_{0j}^\eta} +\xcancel{ \sum_{\beta, \eta=1}^p \left[ \xi_k, r_{0jk}^{\alpha\beta}\right] r_{0kj}^{\beta \eta} \mathcal{L}_{\xi_k}\left(T_{0j}^\eta \right) } +\sum_{\beta=1}^p r_{0jk}^{\alpha\beta} \mathcal{L}_{\xi_k}\mathcal{L}_{\xi_k} \left( T_{0k}^\beta \right)  \\
&+ \bcancel{\sum_{\delta=1}^p \left[ \xi_k, A_j^{\alpha \delta} \right] T_{0j}^\delta }+ \cancel{\sum_{\delta=1}^p A_j^{\alpha \delta} \mathcal{L}_{\xi_k}\left(T_{0j}^\delta \right)} - \xcancel{\sum_{\beta, \eta=1}^p \mathcal{L}_{\xi_k}\left( r_{0jk}^{\alpha \beta} A_k^{\beta \eta} T_{0k}^\eta  \right) } - \cancel{ \sum_{\gamma, \eta=1}^p \mathcal{L}_{\xi_k}\left( \left[ \xi_k, r_{0jk}^{\alpha \gamma} \right] T_{0k}^\eta \right) } - \mathcal{L}_{\xi_j} \mathcal{L}_{\xi_j}\left( T_{0j}^\alpha \right) + \cancel{\mathcal{L}_{\xi_j} \mathcal{L}_{\xi_k} \left( T_{0j}^\alpha \right)} \\
& + \sum_{\delta=1}^p \left[ \xi_j, A_j^{\alpha \delta} \right] T_{0j}^\delta + \sum_{\delta=1}^p A_j^{\alpha \delta} \mathcal{L}_{\xi_j}\left( T_{0j}^\delta \right) - \cancel{\sum_{\beta, \eta=1}^p \mathcal{L}_{\xi_j}\left( r_{0jk}^{\alpha \beta} A_k^{\beta\eta} T_{0k}^\eta\right)} - \cancel{\sum_{\gamma, \eta=1}^p \mathcal{L}_{\xi_j}\left(\left[ \xi_k, r_{0jk}^{\alpha \gamma}\right] T_{0k}^\gamma \right)}
\end{align*}}}

Then we have
{\small{\begin{align*}
& -\mathcal{L}_{d_{jk}}\left( T_{0j}^\alpha \right) + E_{jk}\left( T_{0j}^\alpha \right) + F_{jk}\left( T_{0j}^\alpha \right) \\
&= -2 \alpha_j \left(A_j^\alpha \right) + 2 \sum_{\beta=1}^p r_{0jk}^{\alpha\beta} \alpha_k \left( A_k^\beta  \right) + \sum_{\beta, \eta=1}^p r_{0jk}^{\alpha\beta} A_k^{\beta \eta} A_k^{\eta \gamma} T_{0k}^\gamma - 2\sum_{\gamma, \eta=1}^p r_{0jk}^{\alpha\gamma} \mathcal{L}_{\xi_k}\left( \alpha_k \left( T_{0k}^\gamma  \right) \right)  - 2\sum_{\beta, \eta=1}^p r_{0jk}^{\alpha \beta} A_k^{\beta \eta} \mathcal{L}_{\xi_k}\left( T_{0k}^\eta \right) \\
&- \sum_{\eta, \gamma=1}^p r_{0jk}^{\alpha \eta} \left[ \xi_k, A_k^{\eta \gamma} \right] T_{0k}^\gamma +2 \sum_{\delta=1}^p A_j^{\alpha \delta} \mathcal{L}_{\xi_j}\left( T_{0j}^\alpha \right) - \sum_{\delta, \xi=1}^p A_j^{\alpha \delta} A_j^{\delta \xi} T_{0j}^\xi + 2 \mathcal{L}_{\xi_j}\left( \alpha_j \left( T_{0j}^\alpha \right) \right) + \sum_{\beta=1}^p r_{0jk}^{\alpha\beta} \mathcal{L}_{\xi_k}\mathcal{L}_{\xi_k} \left(T_{0k}^\beta \right) + \sum_{\delta=1}^p \left[ \xi_j, A_j^{\alpha \delta} \right] T_{0j}^\delta \\
&- \mathcal{L}_{\xi_j} \mathcal{L}_{\xi_j}\left( T_{0j}^\alpha \right) \\
&= - 2 \sum_{\delta=1}^p A_j^{\alpha \delta} \alpha_j \left( T_{0j}^\delta \right) +2 \mathcal{L}_{\xi_j}\left(\alpha_j \left( T_{0j}^\alpha \right) \right)+ 2\sum_{\delta=1}^p A_j^{\alpha \delta} \mathcal{L}_{\xi_j}\left(T_{0j}^\delta \right) + \sum_{\delta=1}^p \left[ \xi_j, A_j^{\alpha \delta} \right] T_{0j}^\delta -\sum_{\delta, \xi=1}^p A_j^{\alpha \delta} A_j^{\delta \xi} T_{0j}^\xi - \mathcal{L}_{\xi_j} \mathcal{L}_{\xi_j} \left( T_{0j}^\alpha \right)\\
&- \sum_{\beta=1}^p r_{0jk}^{\alpha\beta}\left( - 2 \sum_{\delta=1}^p A_k^{\beta \delta} \alpha_k \left( T_{0k}^\delta \right) +2 \mathcal{L}_{\xi_k}\left(\alpha_k \left( T_{0k}^\beta \right) \right)+ 2\sum_{\delta=1}^p A_k^{\beta \delta} \mathcal{L}_{\xi_k}\left(T_{0k}^\delta \right) + \sum_{\delta=1}^p \left[ \xi_k, A_k^{\beta \delta} \right] T_{0k}^\delta -\sum_{\delta, \xi=1}^p A_k^{\beta \delta} A_k^{\delta \xi} T_{0k}^\xi - \mathcal{L}_{\xi_k} \mathcal{L}_{\xi_k} \left( T_{0k}^\beta \right) \right)
\end{align*}}}

Then we set 
\begin{align*}
F_j\left(T_{0j}^\alpha \right) =  - 2 \sum_{\delta=1}^p A_j^{\alpha \delta} \alpha_j \left( T_{0j}^\delta \right) +2 \mathcal{L}_{\xi_j}\left(\alpha_j \left( T_{0j}^\alpha \right) \right)+ 2\sum_{\delta=1}^p A_j^{\alpha \delta} \mathcal{L}_{\xi_j}\left(T_{0j}^\delta \right) + \sum_{\delta=1}^p \left[ \xi_j, A_j^{\alpha \delta} \right] T_{0j}^\delta -\sum_{\delta, \xi=1}^p A_j^{\alpha \delta} A_j^{\delta \xi} T_{0j}^\xi - \mathcal{L}_{\xi_j} \mathcal{L}_{\xi_j} \left( T_{0j}^\alpha \right)
\end{align*}
Then we have
\begin{align*}
\bar{\partial} F_j^\alpha&= - 2 \sum_{\delta=1}^p a_j^{\alpha \delta} \wedge \alpha_j\left( T_{0j}^\delta \right) + 2 \mathcal{L}_{\psi}\left( \alpha_j\left( T_{0j}^\alpha \right) \right) + 2 \sum_{\delta=1}^p a_j^{\alpha \delta}\wedge \mathcal{L}_\xi\left( T_{0j}^\delta \right) + 2 \sum_{\delta=1}^p A_j^{\alpha \delta} \mathcal{L}_{\psi}\left( T_{0j}^\delta \right) + \sum_{\delta=1}^p \left[ \psi, A_j^{\alpha \delta }\right] T_{0j}^\delta\\
&+ \sum_{\delta=1}^p \left[ \xi_j, a_j^{\alpha \delta}\right] T_{0j}^\delta - \sum_{\delta, \xi=1}^p a_j^{\alpha \delta} \wedge A_j^{\delta \xi} T_{0j}^\xi - \sum_{\delta , \xi=1}^p A_j^{\alpha \delta} a_j^{\delta \xi} \wedge T_{0j}^\xi - \mathcal{L}_\psi \mathcal{L}_\xi\left(T_{0j}^\alpha \right) - \mathcal{L}_{\xi_j}\mathcal{L}_\psi \left(T_{0j}^\alpha \right)
\end{align*}

Then since $\left( p_j, U_j\right)=\left[\left(\psi, a_j \right),\left( \xi_j, A_j\right) \right]$, we have
{\small{\begin{align*}
&-\mathcal{L}_{p_j}\left( T_{0j}^\alpha \right) + U_j^\alpha + \bar{\partial} F_j^\alpha\\
& =  - \mathcal{L}_{\left[ \psi, \xi_j \right]}\left(T_{0j}^\alpha \right) + \left( \sum_{\beta=1}^p \left[ \psi, A_j^{\alpha \delta} \right]T_{0j}^\delta - \cancel{ \sum_{\beta=1}^p \left[ \xi_j, a_j^{\alpha \delta} \right] T_{0j}^\delta} +  \bcancel{ \sum_{\delta, \xi=1}^p A_j^{\alpha \delta} a_j^{\delta \xi} \wedge T_{0j}^\xi } - \sum_{\delta, \xi=1}^p a_j^{\alpha \delta} \wedge A_j^{\delta \xi} T_{0j}^\xi\right) \\
& - 2 \sum_{\delta=1}^p a_j^{\alpha \delta} \wedge \alpha_j\left(T_{0j}^\delta \right) + 2 \mathcal{L}_{\psi}\left(\alpha_j\left(T_{0j}^\alpha \right) \right) + 2 \sum_{\delta=1}^p a_j^{\alpha \delta} \wedge \mathcal{L}_{\xi_j}\left( T_{0j}^\delta\right) + 2 \sum_{\delta=1}^p A_j^{\alpha \delta} \mathcal{L}_\psi \left( T_{0j}^\delta \right) + \sum_{\delta=1}^p \left[ \psi, A_j^{\alpha \delta} \right] T_{0j}^\delta\\
&+ \cancel{\sum_{\delta=1}^p \left[ \xi_j, a_j^{\alpha \delta}\right] T_{0j}^\delta} - \sum_{\delta, \xi=1}^p a_j^{\alpha \delta} \wedge A_j^{\delta \xi} T_{0j}^\xi - \bcancel{\sum_{\delta , \xi=1}^p A_j^{\alpha \delta} a_j^{\delta \xi} \wedge T_{0j}^\xi }- \mathcal{L}_\psi \mathcal{L}_{\xi_j}\left(T_{0j}^\alpha \right)- \mathcal{L}_{\xi_j}\mathcal{L}_\psi \left(T_{0j}^\alpha \right)\\
&=- 2\mathcal{L}_\psi\mathcal{L}_{\xi_j}\left(T_{0j}^\alpha \right) +2 \mathcal{L}_\psi\left( \sum_{\delta=1}^p A_j^{\alpha \delta} T_{0j}^\delta \right) + 2 \mathcal{L}_\psi\left(\alpha_j\left(T_{0j}^\alpha \right) \right) + 2 \sum_{\delta=1}^p a_j^{\alpha \delta} \wedge \mathcal{L}_{\xi_j}\left(T_{0j}^\delta \right) - 2 \sum_{\delta, \xi=1}^q a_j^{\alpha \delta} \wedge A_j^{\delta \xi} T_{0j}^\xi - 2\sum_{\delta=1}^p a_j^{\alpha \delta}\wedge \alpha_j\left(T_{0j}^\delta \right)\\
&=-2 \mathcal{L}_\psi\left( \mathcal{L}_{\xi_j} \left( T_{0j}^\alpha \right) - A_j^\alpha - \alpha_j\left( T_{0j}^\alpha\right)        \right) + 2\sum_{\delta=1}^p a_j^{\alpha \delta} \wedge \left( \mathcal{L}_{\xi_j}\left(T_{0j}^\delta \right)    - A_j^\delta - \alpha_j\left(T_{0j}^\delta \right) \right) = -2 \mathcal{L}_\psi \left(\Phi_j^\alpha  \right) + 2 \sum_{\delta=1}^p a_j^{\alpha \delta} \wedge \Phi_j^\delta
\end{align*}}}

Let us consider
{\Small{\begin{align}\label{a10}
&\left[ T_{0j}^\alpha, F_j\left( T_{0j}^\beta \right) \right] - \left[ T_{0j}^\beta, F_j\left(T_{0j}^\alpha \right) \right]- \sum_{\gamma=1}^p g_{0j\alpha\beta}F_j\left(T_{0j}^\gamma \right)\\
&\equiv_{\Theta_{\mathcal{F}_0}} \left[ T_{0j}^\alpha, - 2 \sum_{\delta=1}^p A_j^{\beta \delta} \alpha_j \left( T_{0j}^\delta \right) +2 \left[ \xi_j, \alpha_j \left( T_{0j}^\beta \right) \right]+ 2\sum_{\delta=1}^p A_j^{\beta \delta} \left[ \xi_j, T_{0j}^\delta \right]  - \left[ \xi_j,  \left[ \xi_j,  T_{0j}^\beta \right]\right] \right] \notag \\
&- \left[ T_{0j}^\beta, - 2 \sum_{\delta=1}^p A_j^{\alpha \delta} \alpha_j \left( T_{0j}^\delta \right) +2 \left[ \xi_j,  \alpha_j \left( T_{0j}^\alpha \right) \right]+ 2\sum_{\delta=1}^p A_j^{\alpha \delta} \left[ \xi_j, T_{0j}^\delta  \right] - \left[ \xi_j, \left[\xi_j ,  T_{0j}^\alpha \right] \right]\right] \notag\\
&- \sum_{\gamma=1}^p g_{0j\alpha\beta}^\gamma \left(- 2 \sum_{\delta=1}^p A_j^{\gamma \delta} \alpha_j \left( T_{0j}^\delta \right) +2 \left[ \xi_j,  \alpha_j \left( T_{0j}^\gamma \right) \right]+ 2\sum_{\delta=1}^p A_j^{\gamma \delta} \left[ \xi_j, T_{0j}^\delta  \right] - \left[ \xi_j, \left[\xi_j ,  T_{0j}^\gamma \right] \right] \right) \notag\\
&=-2 \sum_{\delta=1}^p \left[ T_{0j}^\alpha , A_j^{\beta \delta} \right] \alpha_j\left(T_{0j}^\delta \right) -2 \sum_{\delta=1}^p A_j^{\beta \delta} \left[ T_{0j}^\alpha , \alpha_j\left( T_{0j}^\delta \right) \right] +2 \left[ T_{0j}^\alpha, \left[ \xi_j, \alpha_j\left( T_{0j}^\beta \right) \right] \right] +2 \sum_{\delta=1}^p \left[ T_{0j}^\alpha , A_j^{\beta \delta} \right] \left[ \xi_j, T_{0j}^\delta \right] + 2 \sum_{\delta=1}^p A_j^{\beta \delta} \left[ T_{0j}^\alpha, \left[ \xi_j, T_{0j}^\delta \right] \right] \notag \\
&- \left[ T_{0j}^\alpha, \left[ \xi_j , \left[ \xi_j, T_{0j}^\beta \right] \right] \right] +2 \sum_{\delta=1}^p \left[ T_{0j}^\beta, A_j^{\alpha \delta} \right] \alpha_j\left( T_{0j}^\delta  \right) + 2\sum_{\delta=1}^p A_j^{\alpha \delta} \left[ T_{0j}^\beta, \alpha_j\left( T_{0j}^\delta \right) \right] - 2 \left[ T_{0j}^\beta, \left[ \xi_j, \alpha_j\left( T_{0j}^\alpha \right) \right] \right] - 2\sum_{\delta=1}^p \left[ T_{0j}^\beta, A_j^{\alpha\delta}\right]\left[\xi_j, T_{0j}^\delta  \right] \notag \\
&-2 \sum_{\delta=1}^p A_j^{\alpha \delta} \left[ T_{0j}^\beta, \left[ \xi_j, T_{0j}^\delta\right] \right] + \left[ T_{0j}^\beta, \left[ \xi_j , \left[ \xi_j, T_{0j}^\alpha \right] \right] \right] + 2\sum_{\gamma, \delta=1}^p g_{0j\alpha\beta}^\gamma A_j^{\gamma \delta} \alpha_j\left(T_{0j}^\delta \right) - 2\sum_{\gamma=1}^p g_{0j\alpha\beta}^\gamma \left[ \xi_j, \alpha_j\left( T_{0j}^\gamma \right) \right] - 2\sum_{\gamma, \delta=1}^p g_{0j\alpha\beta}^\gamma A_j^{\gamma \delta}\left[ \xi_j, T_{0j}^\delta \right] \notag \\
&+ \sum_{\gamma=1}^p g_{0j\alpha\beta}^\gamma \left[ \xi_j, \left[ \xi_j, T_{0j}^\gamma \right] \right] \notag
\end{align}}}

We note that 
\begin{align}\label{a11}
\left[ T_{0j}^\beta, \left[ \xi_j, \left[ \xi_j, T_{0j}^\alpha \right] \right] \right]&= \left[ \left[  \xi_j, T_{0j}^\alpha \right], \left[ \xi_j, T_{0j}^\beta \right] \right] + \left[ \xi_j, \left[ T_{0j}^\beta , \left[ \xi_j, T_{0j}^\alpha \right] \right]\right] \\
- \left[ T_{0j}^\alpha, \left[ \xi_j , \left[ \xi_j, T_{0j}^\beta \right] \right] \right] &= \left[ \left[ \xi_j, T_{0j}^\alpha \right], \left[ \xi_j, T_{0j}^\beta \right] \right]  - \left[ \xi_j, \left[ T_{0j}^\alpha, \left[ \xi_j,  T_{0j}^\beta \right] \right]\right]  \notag\\
\sum_{\gamma=1}^p g_{0j\alpha\beta}^\gamma \left[ \xi_j, \left[ \xi_j, T_{0j}^\gamma \right] \right] & = \sum_{\gamma=1}^p \left[ \xi_j, g_{0j\alpha\beta}^\gamma \left[ \xi_j, T_{0j}^\gamma \right] \right] - \sum_{\gamma=1}^p \left[ \xi_j, g_{0j\alpha\beta}^\gamma \right]\left[ \xi_j, T_{0j}^\gamma \right] \notag \\
&= \left[ \xi_j, \left[ \xi_j,  \left[ T_{0j}^\alpha, T_{0j}^\beta \right] \right]\right] - \sum_{\gamma=1}^p \left[ \xi_j, \left[ \xi_j, g_{0j\alpha\beta}^\gamma \right] T_{0j}^\gamma \right]  - \sum_{\gamma=1}^p \left[ \xi_j, g_{0j\alpha\beta}^\gamma \right]\left[ \xi_j, T_{0j}^\gamma \right] \notag \\
&\equiv_{\Theta_{\mathcal{F}_0}} \left[ \xi_j, \left[ \xi_j,  \left[ T_{0j}^\alpha, T_{0j}^\beta \right] \right]\right]  - 2 \sum_{\gamma=1}^p \left[ \xi_j, g_{0j\alpha\beta}^\gamma \right]\left[ \xi_j, T_{0j}^\gamma \right] \notag
\end{align}

We note that since $\left[ T_{0j}^\alpha, \alpha_j\left( T_{0j}^\beta\right) \right]- \left[ T_{0j}^\beta, \alpha_j\left( T_{0j}^\alpha \right) \right] - \sum_{\gamma=1}^p g_{0j\alpha\beta}^\gamma \alpha_j\left( T_{0j}^\gamma \right)= Q_{j\alpha\beta}$, we have
{\Small{\begin{align*}
K_{j\alpha\beta}:&= \left[ T_{0j}^\alpha, \Phi_j\left( T_{0j}^\beta \right) \right] - \left[ T_{0j}^\beta,  \Phi_j\left(T_{0j}^\alpha \right) \right] - \sum_{\gamma=1}^p g_{0j\alpha\beta}^\gamma \Phi_j\left(T_{0j}^\gamma \right) \\
&= \left[ T_{0j}^\alpha, \left[ \xi_j, T_{0j}^\beta \right]   - A_j\left(T_{0j}^\beta \right) - \alpha_j\left(T_{0j}^\beta \right) \right] - \left[ T_{0j}^\beta, \left[ \xi_j, T_{0j}^\alpha \right] - A_j\left( T_{0j}^\alpha \right) - \alpha_j\left(T_{0j}^\alpha \right) \right] - \sum_{\gamma=1}^p g_{0j\alpha\beta}^\gamma\left(\left[ \xi_j, T_{0j}^\gamma \right] - A_j\left( T_{0j}^\gamma \right)  - \alpha_j\left(T_{0j}^\gamma \right) \right) \\
&=- Q_{j\alpha\beta} - \left[ T_{0j}^\alpha, A_j\left(T_{0j}^\beta \right) \right] + \left[ T_{0j}^\beta, A_j\left( T_{0j}^\alpha \right) \right] +  \sum_{\gamma=1}^p g_{0j\alpha\beta}^\gamma  A_j\left( T_{0j}^\gamma \right) + \sum_{\gamma=1}^p \left[ \xi_j, g_{0j\alpha\beta}^\gamma \right] T_{0j}^\gamma\\
&=- Q_{j\alpha\beta}- \sum_{\eta=1}^p \left[ T_{0j}^\alpha, A_j^{\beta \eta} T_{0j}^\eta \right] + \sum_{\delta=1}^p \left[ T_{0j}^\beta, A_j^{\alpha \delta} T_{0j}^\delta \right] + \sum_{\gamma, \xi=1}^p g_{0j\alpha\beta}^\gamma A_j^{\gamma \xi } T_{0j}^\xi + \sum_{\gamma=1}^p  \left[ \xi_j, g_{0j\alpha\beta}^\gamma \right] T_{0j}^\gamma \\
&= - Q_{j\alpha\beta}- \sum_{\eta=1}^p \left[T_{0j}^\alpha, A_j^{\beta \eta} \right] T_{0j}^\eta - \sum_{\eta, \gamma=1}^p A_j^{\beta \eta} g_{0j\alpha \eta}^\gamma T_{0j}^\gamma + \sum_{\delta=1}^p \left[ T_{0j}^\beta, A_j^{\alpha \delta} \right] T_{0j}^\delta + \sum_{\delta, \gamma=1}^p A_j^{\alpha \delta} g_{0j\beta \delta}^\gamma T_{0j}^\gamma + \sum_{\gamma, \xi=1}^p g_{0j\alpha\beta}^\gamma A_j^{\gamma \xi } T_{0j}^\xi + \sum_{\gamma=1}^p  \left[ \xi_j, g_{0j\alpha\beta}^\gamma \right] T_{0j}^\gamma
\end{align*}}}

Let us compute
{\Small{\begin{align}\label{a12}
&-2\left[ \Phi_j\left(T_{0j}^\alpha \right), \Phi_j\left(T_{0j}^\alpha \right) \right] + 2 \Phi_j\left( K_{j\alpha\beta} \right)\\
&=-2\left[ \left[ \xi_j, T_{0j}^\alpha \right] - \sum_{\delta=1}^p A_j^{\alpha \delta} T_{0j}^\delta - \alpha_j\left(T_{0j}^\alpha \right) , \left[ \xi_j, T_{0j}^\beta \right]- \sum_{\eta=1}^p A_j^{\beta \eta} T_{0j}^\eta -\alpha_j\left(T_{0j}^\beta \right) \right] +2\Phi_j\left( - Q_{j\alpha\beta}\right) \notag \\
& + 2\Phi_j\left( - \sum_{\eta=1}^p \left[T_{0j}^\alpha, A_j^{\beta \eta} \right] T_{0j}^\eta - \sum_{\eta, \gamma=1}^p A_j^{\beta \eta} g_{0j\alpha \eta}^\gamma T_{0j}^\gamma + \sum_{\delta=1}^p \left[ T_{0j}^\beta, A_j^{\alpha \delta} \right] T_{0j}^\delta + \sum_{\delta, \gamma=1}^p A_j^{\alpha \delta} g_{0j\beta \delta}^\gamma T_{0j}^\gamma + \sum_{\gamma, \xi=1}^p g_{0j\alpha\beta}^\gamma A_j^{\gamma \xi } T_{0j}^\xi + \sum_{\gamma=1}^p  \left[ \xi_j, g_{0j\alpha\beta}^\gamma \right] T_{0j}^\gamma\right) \notag \\
&\equiv_{\Theta_{\mathcal{F}_0}} -2 \left[ \left[ \xi_j, T_{0j}^\alpha \right], \left[ \xi_j, T_{0j}^\beta \right] \right] + 2\sum_{\eta=1}^p A_j^{\beta \eta} \left[ \left[ \xi_j, T_{0j}^\alpha\right], T_{0j}^\eta \right] + 2 \left[ \left[ \xi_j, T_{0j}^\alpha \right], \alpha_j\left( T_{0j}^\beta \right) \right] + 2\sum_{\delta=1}^p A_j^{\alpha \delta} \left[ T_{0j}^\delta , \left[ \xi_j, T_{0j}^\beta \right] \right] - 2\sum_{\delta=1}^p A_j^{\alpha \delta} \left[ T_{0j}^\delta, \alpha_j\left( T_{0j}^\beta \right) \right] \notag \\
&+ 2\left[ \alpha_j\left( T_{0j}^\alpha \right), \left[ \xi_j, T_{0j}^\beta \right] \right] - 2\sum_{\eta=1}^p A_j^{\beta \eta}\left[ \alpha_j\left( T_{0j}^\alpha \right), T_{0j}^\eta \right] - 2 \left[ \alpha_j\left( T_{0j}^\alpha  \right) , \alpha_j\left( T_{0j}^\beta \right) \right] - 2\left[ \xi_j, Q_{j\alpha\beta} \right] + 2\alpha_j\left( Q_{j\alpha\beta} \right) \notag \\
&- 2\sum_{\eta=1}^p \left[ T_{0j}^\alpha, A_j^{\beta \eta} \right] \left[ \xi_j, T_{0j}^\eta \right] + 2\sum_{\eta=1}^p \left[ T_{0j}^\alpha, A_j^{\beta \eta} \right] \alpha_j\left(  T_{0j}^\eta\right) - 2 \sum_{\eta, \gamma=1}^p A_j^{\beta \eta} g_{0j\alpha \eta}^\gamma \left[ \xi_j, T_{0j}^\gamma \right] + 2 \sum_{\eta, \gamma=1}^p A_j^{\beta \eta} g_{0j\alpha\eta}^\gamma \alpha_j\left( T_{0j}^\gamma \right) \notag \\
&+2 \sum_{\delta=1}^p \left[ T_{0j}^\beta, A_j^{\alpha \delta} \right] \left[ \xi_j, T_{0j}^\delta \right] - 2\sum_{\delta=1}^p \left[ T_{0j}^\beta, A_j^{\alpha \delta } \right]\alpha_j\left(T_{0j}^\delta  \right) + 2 \sum_{\delta, \gamma=1}^p A_j^{\alpha \delta} g_{0j\beta \delta}^\gamma \left[ \xi_j, T_{0j}^\gamma \right] -  2 \sum_{\delta, \gamma=1}^p A_j^{\alpha \delta} g_{0j\beta \delta}^\gamma \alpha_j\left( T_{0j}^\gamma  \right)\notag \\
&+2 \sum_{\gamma , \xi=1}^p g_{0j\alpha\beta}^\gamma A_j^{\gamma \xi} \left[ \xi_j, T_{0j}^\xi \right] - 2 \sum_{\gamma , \xi=1}^p g_{0j\alpha\beta}^\gamma A_j^{\gamma \xi} \alpha_j\left( T_{0j}^\xi \right) + 2 \sum_{\gamma=1}^p \left[ \xi_j, g_{0j\alpha\beta}^\gamma \right] \left[ \xi_j , T_{0j}^\gamma \right] - 2 \sum_{\gamma=1}^p \left[ \xi_j, g_{0j\alpha\beta}^\gamma \right] \alpha_j\left( T_{0j}^\gamma \right) \notag
\end{align}}}

Then  from $(\ref{a10})$ and $(\ref{a11})$ and $(\ref{a12})$, we have
{\tiny{\begin{align*}
&\left[ T_{0j}^\alpha, F_j\left( T_{0j}^\beta \right) \right] - \left[ T_{0j}^\beta, F_j\left(T_{0j}^\alpha \right) \right]- \sum_{\gamma=1}^p g_{0j\alpha\beta}F_j\left(T_{0j}^\gamma \right) -2\left[ \Phi_j\left(T_{0j}^\alpha \right), \Phi_j\left(T_{0j}^\alpha \right) \right] + 2 \Phi_j\left( K_{j\alpha\beta} \right) \\
&\equiv_{\Theta_{\mathcal{F}_0}} \cancel{-2 \sum_{\delta=1}^p \left[ T_{0j}^\alpha , A_j^{\beta \delta} \right] \alpha_j\left(T_{0j}^\delta \right)}^1 -\cancel{2 \sum_{\delta=1}^p A_j^{\beta \delta} \left[ T_{0j}^\alpha , \alpha_j\left( T_{0j}^\delta \right) \right] }^2 +2 \left[ T_{0j}^\alpha, \left[ \xi_j, \alpha_j\left( T_{0j}^\beta \right) \right] \right] + \cancel{2 \sum_{\delta=1}^p \left[ T_{0j}^\alpha , A_j^{\beta \delta} \right] \left[ \xi_j, T_{0j}^\delta \right] }^3 + \cancel{2 \sum_{\delta=1}^p A_j^{\beta \delta} \left[ T_{0j}^\alpha, \left[ \xi_j, T_{0j}^\delta \right] \right]}^4 \\
&+\bcancel{2 \sum_{\delta=1}^p \left[ T_{0j}^\beta, A_j^{\alpha \delta} \right] \alpha_j\left( T_{0j}^\delta  \right)}^1 + \bcancel{2\sum_{\delta=1}^p A_j^{\alpha \delta} \left[ T_{0j}^\beta, \alpha_j\left( T_{0j}^\delta \right) \right] }^2 - 2 \left[ T_{0j}^\beta, \left[ \xi_j, \alpha_j\left( T_{0j}^\alpha \right) \right] \right] - \bcancel{2\sum_{\delta=1}^p \left[ T_{0j}^\beta, A_j^{\alpha\delta}\right]\left[\xi_j, T_{0j}^\delta  \right] }^3 - \bcancel{ 2 \sum_{\delta=1}^p A_j^{\alpha \delta} \left[ T_{0j}^\beta, \left[ \xi_j, T_{0j}^\delta\right] \right] }^4 \\
& + \xcancel{2\sum_{\gamma, \delta=1}^p g_{0j\alpha\beta}^\gamma A_j^{\gamma \delta} \alpha_j\left(T_{0j}^\delta \right)}^1 - 2\sum_{\gamma=1}^p g_{0j\alpha\beta}^\gamma \left[ \xi_j, \alpha_j\left( T_{0j}^\gamma \right) \right] - \xcancel{2\sum_{\gamma, \delta=1}^p g_{0j\alpha\beta}^\gamma A_j^{\gamma \delta}\left[ \xi_j, T_{0j}^\delta \right]}^2 +\xcancel{2 \left[ \left[ \xi_j , T_{0j}^\alpha \right], \left[ \xi_j, T_{0j}^\beta \right] \right] }^3 - \xcancel{ 2 \sum_{\gamma=1}^p \left[ \xi_j, g_{0j\alpha\beta}^\gamma \right] \left[\xi_j, T_{0j}^\gamma \right] }^4 \\
&-\xcancel{2 \left[ \left[ \xi_j, T_{0j}^\alpha \right], \left[ \xi_j, T_{0j}^\beta \right] \right]}^3 + \cancel{2\sum_{\eta=1}^p A_j^{\beta \eta} \left[ \left[ \xi_j, T_{0j}^\alpha\right], T_{0j}^\eta \right]}^4 + 2 \left[ \left[ \xi_j, T_{0j}^\alpha \right], \alpha_j\left( T_{0j}^\beta \right) \right] +  \bcancel{2\sum_{\delta=1}^p A_j^{\alpha \delta} \left[ T_{0j}^\delta , \left[ \xi_j, T_{0j}^\beta \right] \right]}^4 - \bcancel{2\sum_{\delta=1}^p A_j^{\alpha \delta} \left[ T_{0j}^\delta, \alpha_j\left( T_{0j}^\beta \right) \right]}^2 \\
&+ 2\left[ \alpha_j\left( T_{0j}^\alpha \right), \left[ \xi_j, T_{0j}^\beta \right] \right] -\cancel{ 2\sum_{\eta=1}^p A_j^{\beta \eta}\left[ \alpha_j\left( T_{0j}^\alpha \right), T_{0j}^\eta \right]}^2 - 2 \left[ \alpha_j\left( T_{0j}^\alpha  \right) , \alpha_j\left( T_{0j}^\beta \right) \right] - 2\left[ \xi_j, Q_{j\alpha\beta} \right] + 2\alpha_j\left( Q_{j\alpha\beta} \right) \\
&- \cancel{2\sum_{\eta=1}^p \left[ T_{0j}^\alpha, A_j^{\beta \eta} \right] \left[ \xi_j, T_{0j}^\eta \right]}^3 + \cancel{2\sum_{\eta=1}^p \left[ T_{0j}^\alpha, A_j^{\beta \eta} \right] \alpha_j\left(  T_{0j}^\eta\right) }^1 - \cancel{2 \sum_{\eta, \gamma=1}^p A_j^{\beta \eta} g_{0j\alpha \eta}^\gamma \left[ \xi_j, T_{0j}^\gamma \right] }^4 + \cancel{2 \sum_{\eta, \gamma=1}^p A_j^{\beta \eta} g_{0j\alpha\eta}^\gamma \alpha_j\left( T_{0j}^\gamma \right) }^2\\
&+\bcancel{2 \sum_{\delta=1}^p \left[ T_{0j}^\beta, A_j^{\alpha \delta} \right] \left[ \xi_j, T_{0j}^\delta \right]}^3 - \bcancel{2\sum_{\delta=1}^p \left[ T_{0j}^\beta, A_j^{\alpha \delta } \right]\alpha_j\left(T_{0j}^\delta  \right) }^1+ \bcancel{ 2 \sum_{\delta, \gamma=1}^p A_j^{\alpha \delta} g_{0j\beta \delta}^\gamma \left[ \xi_j, T_{0j}^\gamma \right] }^4 -  \bcancel{ 2 \sum_{\delta, \gamma=1}^p A_j^{\alpha \delta} g_{0j\beta \delta}^\gamma \alpha_j\left( T_{0j}^\gamma  \right) }^2 \\
&+\xcancel{2 \sum_{\gamma , \xi=1}^p g_{0j\alpha\beta}^\gamma A_j^{\gamma \xi} \left[ \xi_j, T_{0j}^\xi \right]}^2 - \xcancel{ 2 \sum_{\gamma , \xi=1}^p g_{0j\alpha\beta}^\gamma A_j^{\gamma \xi} \alpha_j\left( T_{0j}^\xi \right)}^1 + \xcancel{2 \sum_{\gamma=1}^p \left[ \xi_j, g_{0j\alpha\beta}^\gamma \right] \left[ \xi_j , T_{0j}^\gamma \right] }^4 - 2 \sum_{\gamma=1}^p \left[ \xi_j, g_{0j\alpha\beta}^\gamma \right] \alpha_j\left( T_{0j}^\gamma \right) \\
&=2 \left[ \xi_j, \left[ T_{0j}^\alpha, \alpha_j\left( T_{0j}^\beta \right) \right] \right] - 2 \left[ \xi_j, \left[ T_{0j}^\beta, \alpha_j\left( T_{0j}^\alpha \right) \right] \right] - 2\sum_{\gamma=1}^p \left[ \xi_j, \alpha_j\left( \left[ T_{0j}^\alpha, T_{0j}^\beta \right] \right) \right] - 2\left[ \alpha_j \left( T_{0j}^\alpha \right), \alpha_j\left( T_{0j}^\beta \right) \right] -2 \left[ \xi_j, Q_{j\alpha\beta}\right] + 2\alpha_j\left( Q_{j\alpha\beta} \right) \\
&= - 2\left[ \alpha_j \left( T_{0j}^\alpha \right), \alpha_j\left( T_{0j}^\beta \right) \right] +  2\alpha_j\left( Q_{j\alpha\beta} \right) 
\end{align*}}}

Then we have
{\Small{\begin{align*}
-2\left[ \Phi_j\left(T_{0j}^\alpha \right), \Phi_j\left(T_{0j}^\alpha \right) \right] + 2 \Phi_j\left( K_{j\alpha\beta} \right) \equiv_{\Theta_{\mathcal{F}_0}}  - 2\left[ \alpha_j \left( T_{0j}^\alpha \right), \alpha_j\left( T_{0j}^\beta \right) \right] +  2\alpha_j\left( Q_{j\alpha\beta} \right) - \left[ T_{0j}^\alpha, F_j\left( T_{0j}^\beta \right) \right] + \left[ T_{0j}^\beta, F_j\left(T_{0j}^\alpha \right) \right] + \sum_{\gamma=1}^p g_{0j\alpha\beta}F_j\left(T_{0j}^\gamma \right) 
\end{align*}}}
This completes the proof of Proposition \ref{a13}.
\end{proof}

\section{Infinitesimal foliated deformations and obstructions in terms of cotangent sheaves by using \v Cech resolution of $\mathcal{E}_{\mathcal{N}_{\mathcal{F}_0}^* }^\bullet$}

\subsection{Review of deformations of foliated complex analytic structures in terms of cotangent sheaves}\

We recall the definitions of a foliated analytic family in terms of cotangent sheaves and associated foliated Kodaira-Spencer map in terms of cotangent sheaves from Part I.

\begin{definition}\label{n3}
Suppose that given a domain $B\subset \mathbb{C}^m$, there is a set $\left\{\left(M_t, \mathcal{N}_{\mathcal{F}_t}^*\right)|t\in B \right\}$ of $n$-dimensional compact $($singularly$)$ foliated complex manifold $\left(M_t, \mathcal{N}_{\mathcal{F}_t}^*\right)$ of codimension $q$, so that for each $t\in B$, we have an exact sequence
\begin{align}\label{n1}
0\to \mathcal{N}_{\mathcal{F}_t}^*\to \Omega_{M_t}^1\to \Omega_{M_t}^1/\mathcal{N}_{\mathcal{F}_t}^*\to 0
\end{align}
where $\Omega_{M_t}/\mathcal{N}_{\mathcal{F}_t}^*$ is a torsion-free $\mathcal{O}_{M_t}$-module. We say that $\left\{\left(M_t, \mathcal{N}_{\mathcal{F}_t}^*\right)|t\in B\right\}$ is a family of $($singularly$)$ foliated complex manifold or $($singularly$)$ foliated complex analytic family in terms of cotangent sheaves if there exist a complex manifold $\mathcal{M}$ and a holomorphic map $\pi:\mathcal{M}\to B$ such that
\begin{enumerate}
\item we have an exact sequence
\begin{align}\label{n2}
0 \to \mathcal{N}_{\mathcal{F}}^*\to \Omega_{\mathcal{M}/B}^1\to \Omega_{\mathcal{M}/B}^1/\mathcal{N}_{\mathcal{F}}^*\to 0
\end{align}
where $\mathcal{N}_\mathcal{F}^*$ is a coherent subsheaf of the relative cotangent sheaf $\Omega_{\mathcal{M}/B}^1$ over $B$.
\item for each $t\in B$, $\pi^{-1}(t)= \left(M_t, \mathcal{N}_{\mathcal{F}_t }^*\right)$ and the exact sequence $(\ref{n2})$ is induced from $(\ref{n1})$ by restricting to $M_t$.
\item the rank of the Jacobian of $\pi$ is equal to $m$ at every point of $\mathcal{M}$.
\item let $S_t=\textnormal{Sing}(\mathcal{F}_t)$ for each $t\in B$. Then we have $d\mathcal{N}_{\mathcal{F}}^*\subset \mathcal{N}_{\mathcal{F}}^*\bigwedge \Omega_{\mathcal{M}/B}^1$ at every point of $\mathcal{M}-\bigcup_{t\in B} S_t$.
\item $\Omega_{\mathcal{M}/B}^1/\mathcal{N}_{\mathcal{F}}^*$ is flat over $B$.
\end{enumerate}
We will denote the $($singularly$)$ foliated complex analytic family by $(\mathcal{M}, \mathcal{N}_\mathcal{F}^*, B, \pi)$.

\end{definition}

\begin{remark}
Let $\left(\mathcal{M}, \mathcal{N}_\mathcal{F}^*, B, \pi \right)$ with $\mathcal{N}_\mathcal{F}^*$ locally free be a foliated analytic family. Let $\Delta$ be an open set of $B$. Then the restriction $\left(\mathcal{M}_\Delta= \pi^{-1}(\Delta), \mathcal{N}_{\mathcal{F}}^*|_\Delta, \Delta, \pi |_{\mathcal{M}_\Delta} \right)$ is also a foliated complex analytic family in terms of cotangent sheaves. We will denote the family by $\left( \mathcal{M}_\Delta, \mathcal{N}_{\mathcal{F}_\Delta}^*, \Delta, \pi\right)$.
\end{remark}

\begin{remark}\label{oc2}
Given a foliated analytic family $\left( \mathcal{M}, \mathcal{N}_\mathcal{F}^*, B, \pi \right)$ with $\mathcal{N}_{\mathcal{F}}^*$ locally free as in \textnormal{Definition \ref{n3}}, we can choose a system of local complex coordinates $\{\xi_i,...,\xi_j,...\}, \xi_j:p\to \xi_j(p)$, and coordinate polydisks $\mathcal{U}_j$ with respect to $\xi_j$, satisfying $(1)$ $\xi_j(p)=\left(\xi_j^1(p),...,\xi_j^n(p),t_1,...,t_m\right),\left(t_1,...,t_m\right)=\pi(p)$ and $(2)$ $\mathcal{U}=\left\{\mathcal{U}_j|j=1,2,...\right\}$ is locally finite. Then $\left\{p\to (\xi_j^1(p),...,\xi_j^n(p))|\mathcal{U}_j\cap M_t\ne \emptyset \right\}$ gives a system of local complex coordinates on $M_t$. In terms of these coordinates, $\pi$ is the projection given by $\pi: \left(\xi_j^1,...,\xi_j^n, t_1,...,t_m \right)\to \left(t_1,...,t_m\right)$. For $j,k$ with $\mathcal{U}_j\cap \mathcal{U}_k\ne \emptyset$, we denote that coordinate transformation from $\xi_k$ to $\xi_j$ by $f_{jk}:\left(\xi_k^1,...,\xi_k^n,t \right)\to \left(\xi_j^1,...,\xi_j^n,t \right)=f_{jk}\left(\xi_k^1,...,\xi_k^n,t \right)$. Thus $f_{jk}$ is given by $\xi_j^\alpha=f_{jk}^\alpha\left(\xi_k^1,...,\xi_k^n,t_1,...,t_m\right), \alpha=1,...,n$. On the other hand, we may assume that the locally free sheaf $\mathcal{N}_\mathcal{F}^*$ is trivialized in terms of the covering $\mathcal{U}$ so that $\Gamma\left(\mathcal{U}_j, \mathcal{N}_\mathcal{F}^* \right)\cong \bigoplus^q \Gamma\left(\mathcal{U}_j, \mathcal{O}_\mathcal{M} \right)$ which is generated by holomorphic $1$-forms $w_j^\alpha(\xi_j,t)\in \Gamma\left(\mathcal{U}_j,\mathcal{N}_\mathcal{F}^* \right),\alpha=1,...,q$ of the form
\begin{align}\label{n11}
w_j^\alpha(\xi_j,t):=\sum_{\beta=1}^n w_j^{\alpha\beta}(\xi_j,t)d\xi_j^\beta, \,\,\,\,\,\,\alpha=1,...,q
\end{align}
for some $w_j^{\alpha\beta}(\xi_j,t)\in \Gamma\left(\mathcal{U}_j, \mathcal{O}_\mathcal{M} \right)$, and satisfies
\begin{align}\label{n6}
w_j^1(\xi_j,t)\wedge \cdots \wedge w_j^q(\xi_j,t)\wedge dw_j^\alpha(\xi_j,t)=0,\,\,\,\,\,\,\,\alpha=1,...,q,\,\,\,\,\,\textnormal{on}\,\,\,\,\,\mathcal{U}_j.
\end{align}
For $j,k$ with $\mathcal{U}_j\cap \mathcal{U}_k\ne \emptyset$, we have
\begin{align}\label{oc1}
w_j^\alpha(\xi_j,t)= \sum_{\beta=1}^q h_{jk}^{\alpha\beta}(\xi_k, t)w_j^\beta(\xi_k,t)
\end{align}
for some $h_{jk}^{\alpha\beta}(\xi_k,t)\in \Gamma\left( \mathcal{U}_j\cap \mathcal{U}_j, \mathcal{O}_\mathcal{M} \right)$,  and for $i,j,k$ with $\mathcal{U}_i\cap \mathcal{U}_j\cap \mathcal{U}_k\ne \emptyset$, we have
\begin{align}
h_{ik}^{\alpha \gamma}(\xi_k,t)= \sum_{\beta=1}^q h_{ij}^{\alpha\beta}\left(f_{jk}(\xi_k,t),t\right) h_{jk}^{\beta \gamma}(\xi_k,t)
\end{align}
More specifically, $(\ref{oc1})$ implies that
\begin{align}\label{n10}
\sum_{\gamma,\eta=1}^n w_j^{\alpha\gamma}(f_{jk}(\xi_k,t),t)\frac{\partial f_{jk}^\gamma(\xi_k,t)}{\partial \xi_k^\eta} d\xi_k^\eta=\sum_{\beta=1}^q\sum_{\eta=1}^n h_{jk}^{\alpha\beta}(\xi_k,t)w_k^{\beta\eta}(\xi_k,t)d\xi_k^\eta
\end{align}

\end{remark}

Given a foliated complex analytic family $\left(\mathcal{M}, \mathcal{N}_\mathcal{F}^*, B, \pi \right)$ with $\mathcal{N}_\mathcal{F}^*$ locally free of deformations of $\left(M, \mathcal{N}_{\mathcal{F}_0}^* \right)=\pi^{-1}(0),0\in B$, the infinitesimal foliated deformations in the foliated complex analytic family  are encoded in $\mathbb{H}^1\left( M, \mathcal{N}_{\mathcal{F}_0}^{*\bullet} \right)$ which is the first cohomology group of the following dual leaf complex
\begin{align}
\mathcal{N}_{\mathcal{F}_0}^{*\bullet} : \Theta_M \to  \mathscr{H}om_{\mathcal{O}_M}\left( \mathcal{N}_{\mathcal{F}_0}^*, \frac{\Omega_M^1}{\mathcal{N}_{\mathcal{F}_0}^*} \right)  \to  \mathscr{H}om_{\mathcal{O}_M}\left(\mathcal{N}_{\mathcal{F}_0}^*, \tilde{\mathcal{S}}^2 \right) \to \mathscr{H}om_{\mathcal{O}_M}\left( \mathcal{N}_{\mathcal{F}_0}^* , \tilde{\mathcal{S}}^3 \right) \to \cdots
\end{align}
We will denote the $i$-th cohomology group by $\mathbb{H}^1\left( M , \mathcal{N}_{\mathcal{F}_0}^{*\bullet} \right)$. We can compute $\mathbb{H}^i\left(M_0, \Theta_{\mathcal{F}_0}^\bullet \right)$  by the following \v Cech resolution of $\mathcal{N}_{\mathcal{F}_0}^{*\bullet}$. Here $\delta$ is the \v Cech map and $\mathcal{U}_0= \mathcal{U}\cap M=\left\{ U_j :=\mathcal{U}_j \cap M | j=1,2,... \right\}$ is an open covering of $M$:
{\Small{\begin{equation} \label{oc3}
\begin{CD}
\cdots \\
@AE_2 AA \\
C^0\left(\mathcal{U}_0, \mathscr{H}om_{\mathcal{O}_M}\left( \mathcal{N}_{\mathcal{F}_0}^* , \tilde{\mathcal{S}}^2\right) \right) @>-\delta>> \cdots  \\
@AE_1AA @AE_1 AA \\
C^0\left(\mathcal{U}_0, \mathscr{H}om_{\mathcal{O}_M}\left(\mathcal{N}_{\mathcal{F}_0}^*, \frac{\Omega_M^1}{\mathcal{N}_{\mathcal{F}_0}^*} \right) \right)@>\delta>> C^1\left(\mathcal{U}_0, \mathscr{H}om_{\mathcal{O}_{M}}\left(\mathcal{N}_{\mathcal{F}_0}^*, \frac{\Omega_M^1}{\mathcal{N}_{\mathcal{F}_0}^*}  \right) \right) @>-\delta>> \cdots \\
@AE_0 AA @AE_0 AA @AE_0 AA \\
C^0(\mathcal{U}_0, \Theta_{M}) @>-\delta>> C^1(\mathcal{U}_0, \Theta_{M}) @>\delta>> C^2(\mathcal{U}_0, \Theta_{M}) @>>> \cdots 
\end{CD}
\end{equation}}}
Then we can define the foliated Kodaira-Spencer map at $0\in B$ in terms of cotangent sheaves which encodes the information of infinitesimal foliated deformations of $\pi^{-1}(0)= \left( M, \mathcal{N}_{\mathcal{F}_0}^* \right)$ in the foliated analytic family $\left(\mathcal{M},  \mathcal{N}_\mathcal{F}^*, B, \pi \right)$ (see Part I).
\begin{definition}[foliated Kodaira-Spencer map in terms of cotangent sheaf]\label{n15}
Let $\left(\mathcal{M}, \mathcal{N}_\mathcal{F}^*, B, \pi \right)$ with $\mathcal{N}_{\mathcal{F}}^*$ locally free be a foliated complex analytic family of deformations of $\left(M_t, \mathcal{N}_{\mathcal{F}_t}^*\right)=\pi^{-1}(t),t\in B$, where $B$ is a domain of $\mathbb{C}^m$. We keep the notations in \textnormal{Remark \ref{oc2}}. For a tangent vector $\frac{\partial}{\partial t}=\sum_{\lambda=1}^m c_\lambda \frac{\partial}{\partial t_\lambda}, c_\lambda \in \mathbb{C}$, of $B$, the foliated Kodaira-Spencer map in terms of cotangent sheaf at $t$ is defined to be a $\mathbb{C}$-linear map
{\small{\begin{align*}
\varphi_t:T_t(B)&\to \mathbb{H}^1\left(M_t, \mathcal{N}_{\mathcal{F}_t}^{*\bullet} \right)\\
                \frac{\partial}{\partial t}&\mapsto \frac{\partial \left(M_t, \mathcal{N}_{\mathcal{F}_t}^* \right)}{\partial t}|_{t=0}:=\left(\left\{  \sum_{\alpha=1}^n \frac{\partial f_{jk}^\alpha(\xi_k,t)}{\partial t}\frac{\partial}{\partial \xi_j^\alpha} \right\}|_{t=0}, \left\{ \left( w_{0j}^\alpha:=w_j^\alpha(\xi_j,0)\mapsto -\overline{\frac{\partial w_j^\alpha(\xi_j,t)}{\partial t}}|_{t=0}\right)_{\alpha=1,...,q} \right\} \right) \\
                 &\,\,\,\,\,\,\,\,\,\,\in C^1\left(\mathcal{U}_0, \Theta_{M} \right) \bigoplus C^0\left( \mathcal{U}_0, \mathscr{H}om_{\mathcal{O}_{M}}\left(\mathcal{N}_{\mathcal{F}_0}^*, \frac{\Omega_{M}^1}{\mathcal{N}_{\mathcal{F}_0}^*} \right) \right)
\end{align*}}}
Here $\frac{\partial w_j^\alpha(\xi_j,t)}{\partial t}|_{t=0}:=\sum_{\beta=1}^n \frac{\partial w_j^{\alpha\beta}(\xi_j,t)}{\partial t}|_{t=0}d\xi_j^\beta $ is the derivative of $\left(\ref{n11} \right)$ with respect to $t$ at $0$, and $\overline{\frac{\partial w_j^\alpha(\xi_j,t)}{\partial t}}|_{t=0}$ is the image of $\frac{\partial w_j^\alpha(\xi_j,t)}{\partial t}|_{t=0}$ under $\Omega_M^1 \to \frac{\Omega_M^1}{\mathcal{N}_{\mathcal{F}_0}^*}$.
\end{definition}
In Part II we have shown that the complex of sheaves
\begin{align}\label{oc4}
\mathcal{E}_{\mathcal{N}_{\mathcal{F}_0}^*}^\bullet: \mathcal{E}_{\mathcal{N}_{\mathcal{F}_0}^*}\to \mathscr{H}om_{\mathcal{O}_M}\left( \mathcal{N}_{\mathcal{F}_0}^*,  \Omega_M^1 \right) \to  \mathscr{H}om_{\mathcal{O}_M}\left( \mathcal{N}_{\mathcal{F}_0}^* , \tilde{\mathcal{S}}^2 \right) \to \mathscr{H}om_{\mathcal{O}_M}\left( \mathcal{N}_{\mathcal{F}_0}^* , \tilde{\mathcal{S}}^3 \right) \to \cdots
\end{align}
also controls foliated deformations of $\left( M , \mathcal{N}_{\mathcal{F}_0}^* \right)$ in terms of tangent sheaves. We can actually study obstructions of foliated deformations of $\left(M, \mathcal{N}_{\mathcal{F}_0}^* \right)$ in terms of cotangent sheaves by using the complex of sheaves $\mathcal{N}_{\mathcal{F}_0}^{*\bullet}$.  However obstructions of deformations of the underlying locally free sheaf $\mathcal{N}_{\mathcal{F}_0}^*$ are implicit in $\mathbb{H}^2\left( M, \mathcal{N}_{\mathcal{F}_0}^{*\bullet} \right)$, i.e. in \v Cech resolution $(\ref{oc3})$ of $\mathcal{N}_{\mathcal{F}_0}^{*\bullet}$. However the complex of sheaves $\mathcal{E}_{\mathcal{N}_{\mathcal{F}_0}^*}^\bullet$ explicitly shows the infinitesimal deformations and obstructions of deformations of the underlying locally free sheaf $\mathcal{N}_{\mathcal{F}_0}^*$. Hence it is desirable to use $\mathcal{E}_{\mathcal{N}_{\mathcal{F}_0}^*}^\bullet$ instead of $\mathcal{N}_{\mathcal{F}_0}^{*\bullet}$ to identify obstructions clearly. Let us identify the infinitesimal deformations in terms of the complex (\ref{oc4}) of sheaves  $\mathcal{E}_{\mathcal{N}_{\mathcal{F}_0}^*}^\bullet$.

\subsection{Identification of infinitesimal foliated deformations by using $\mathcal{E}_{\mathcal{N}_{\mathcal{F}_0}^*}^\bullet$} \label{oc8} \

Let $\left( \mathcal{M}, \mathcal{N}_\mathcal{F}^*, B, \pi \right)$ be a foliated analytic family with $\mathcal{N}_\mathcal{F}^*$ locally free as in Definition $\ref{n3}$. We keep the notations in Remark $\ref{oc2}$. Let $\mathcal{U}_t=\mathcal{U} \cap M_t =\left\{ U_j^t:=\mathcal{U}_j\cap M_t| j=1,2,... \right\}$ be an open covering of $M_t$. We recall the definition of $\mathcal{E}_{\mathcal{N}_{\mathcal{F}_t}^*}$ on $M_t$ with respect to the open covering $\mathcal{U}_t=\left\{U_j^t\right\}$ of $M_t$ from Part II, and let us consider the \v Cech resolution of $\mathcal{E}_{\mathcal{N}_{\mathcal{F}_t}}^\bullet$:

{\Small{\begin{equation} \label{oc5}
\begin{CD}
\cdots \\
@AE_2^t AA \\
C^0\left(\mathcal{U}_t, \mathscr{H}om_{\mathcal{O}_M}\left( \mathcal{N}_{\mathcal{F}_t}^* , \tilde{\mathcal{S}}_t^2 \right) \right) @>-\delta>> \cdots  \\
@AE_1'^tAA @AE_1 AA \\
C^0\left(\mathcal{U}_t, \mathscr{H}om_{\mathcal{O}_{M_t}}\left(\mathcal{N}_{\mathcal{F}_t}^*, \Omega_{M_t}^1 \right) \right)@>\delta>> C^1\left(\mathcal{U}_t, \mathscr{H}om_{\mathcal{O}_{M}}\left(\mathcal{N}_{\mathcal{F}_t}^*, \Omega_{M_t}^1  \right) \right) @>-\delta>> \cdots \\
@AE_0'^t AA @AE_0'^tAA @AE_0'^tAA \\
C^0\left(\mathcal{U}_t,\mathcal{E}_{\mathcal{N}_{\mathcal{F}_t}^*}\right) @>-\delta>> C^1\left(\mathcal{U}_t, \mathcal{E}_{\mathcal{N}_{\mathcal{F}_t}^*}\right) @>\delta>> C^2 \left(\mathcal{U}_t, \mathcal{E}_{\mathcal{N}_{\mathcal{F}_t}^*} \right) @>>> \cdots 
\end{CD}
\end{equation}}}

By taking the derivative of $(\ref{n11})$, we set $\frac{\partial w_j^{\alpha}(\xi_j,t)}{\partial t}:=\sum_{\beta=1}^n \frac{\partial w_j^{\alpha\beta}(\xi_j,t)}{\partial t}\frac{\partial}{\partial \xi_j^\beta},\alpha=1,...,q $. We define an element $\beta_j(t) \in \Gamma\left( U_j^t, \mathscr{H}om_{\mathcal{O}_{M_t}}\left( \mathcal{N}_{\mathcal{F}_t}^*, \Omega_{M_t}^1 \right) \right)$ in the following way:
\begin{align*}
\beta_j(t):\Gamma\left( U_j^t, \mathcal{N}_{\mathcal{F}_t}^*   \right) &\to \Gamma\left( U_j^t, \Omega_{M_t}^1 \right) \\
   w_j^\alpha=w_j^\alpha(\xi_j,t)&\mapsto -\frac{\partial w_j^{\alpha}(\xi_j,t)}{\partial t}= - \sum_{\beta=1}^n \frac{\partial w_j^{\alpha\beta}(\xi_j,t)}{\partial t}\frac{\partial}{\partial \xi_j^\beta},\,\,\,\,\,\,\,\,\alpha=1,...,q
\end{align*}
and linearly extend to $\Gamma\left( U_j^t , \mathcal{N}_{\mathcal{F}_t}^* \right)$. We define an element $B_{jk}(t) \in \Gamma\left( U_j^t \cap U_k^t, \mathscr{H}om_{\mathcal{O}_{M_t} }\left( \mathcal{N}_{\mathcal{F}_t}^*, \mathcal{N}_{\mathcal{F}_t}^* \right)  \right)$ in the following way:
\begin{align*}
B_{jk}(t):\Gamma\left( U_{jk}^t, \mathcal{N}_{\mathcal{F}_t}^* \right) &\to \Gamma\left( U_{jk}^t, \mathcal{N}_{\mathcal{F}_t}^* \right) \\
 w_j^\alpha &\mapsto \sum_{\beta=1}^q \frac{\partial h_{jk}^{\alpha\beta} (\xi_j,t)}{\partial t} w_k^\beta,\,\,\,\,\,\,\,\,\,\alpha=1,...,q
\end{align*}
and linearly extend to $\Gamma\left( U_{jk}^t, \mathcal{N}_{\mathcal{F}_t}^* \right)$, and we set $\theta_{jk}(t):=  \sum_{\eta=1}^n \frac{\partial f_{jk}^\eta(\xi_k,t)}{\partial t}\frac{\partial}{\partial \xi_j^\eta}\in \Gamma\left( U_{jk}^t , \Theta_{M_t} \right)$. If we consider $\left( \theta_{jk}(t), B_{jk}(t)      \right)$ as an element of $\Gamma\left(U_j^t, \Theta_{M_t}\right)|_{U_{jk}^t}\bigoplus \Gamma\left(U_j^t, \mathscr{H}om_{\mathcal{O}_{M_t}}\left( \mathcal{N}_{\mathcal{F}_t}^*, \mathcal{N}_{\mathcal{F}_t}^* \right)\right) |_{U_{jk}^t}$, then we see that $\left\{ \left(\theta_{jk}(t)  , B_{jk}(t)\right)\right\}\in C^1\left( \mathcal{U}_t, \mathcal{E}_{\mathcal{N}_{\mathcal{F}_t}^*} \right)$ as in the same arguments with $\left\{\left(\theta_{jk}(t), A_{jk}(t)\right) \right\}\in C^1\left(\mathcal{U}_t, \mathcal{E}_{\Theta_{\mathcal{F}_t}} \right)$ in subsection \ref{o41}. Then we claim that $\left( \left\{\beta_j(t)\right\}, \left\{ \left(\theta_{jk}(t), B_{jk}(t) \right) \right\}\right) $ defines a $1$-cocycle in the above \v Cech resolution $(\ref{oc5})$ of $\mathcal{E}_{\mathcal{N}_{\mathcal{F}_t}^*}^\bullet$. As in Part I, by taking the derivative of $(\ref{n6})$ with respect to $t$, we see that $E_1'^t\left(\beta_j(t) \right)=0$. On the other hand, by taking the derivative of $(\ref{n10})$ with respect to $t$, we obtain as in Part I that
\begin{align*}
\left({\beta}_k(t)- {\beta}_j(t)\right)\left(w_j^\alpha\right) + \mathcal{L}_{ \theta_{jk}(t)} \left( w_j^\alpha \right) =\sum_{\beta=1}^p \frac{\partial h_{jk}^{\alpha\beta}}{\partial t} w_k^\beta \Longrightarrow \left(\beta_k(t)- \beta_j(t) \right)(w_j^\alpha)+\mathcal{L}_{\theta_{jk}(t)}\left( w_j^\alpha \right) - B_{jk}(t)(w_j^\alpha)=0
\end{align*}
This implies that $\delta\left(\{\beta_j(t)\} \right) + E_0'\left( \left\{\left(\theta_{jk}(t), B_{jk}(t) \right) \right\} \right)=0$. It remains to show that $\delta\left(\left\{\theta_{jk}(t), B_{jk}(t) \right\} \right)=0$ which follows from the same arguments with $\delta\left(\{\alpha_j(t),  A_{jk}(t)\} \right)=0$ in subsection \ref{o41}. Hence we have the infinitesimal foliated deformation along $\frac{\partial}{\partial t}$ at $t$:
\begin{align}\label{oc7}
&\frac{\partial \left( M_t, \mathcal{N}_{\mathcal{F}_t}^* \right)}{\partial t}= \left( \left\{ \beta_j (t) \right\},   \left\{ \left(\theta_{jk}(t), B_{jk}(t) \right) \right\} \right) \in  \mathbb{H}^1\left( M_t , \mathcal{E}_{\mathcal{N}_{\mathcal{F}_t }^* }^\bullet \right) 
\end{align}

\subsection{Description of obstructed elements in \v Cech resolution of $\mathcal{E}_{\mathcal{N}_{\mathcal{F}_0}^*}^\bullet$}\label{oc9}

\begin{definition}
We say that a compact foliated complex manifold $\left(M, \mathcal{N}_{\mathcal{F}_0}^*\right)$ with $\mathcal{N}_{\mathcal{F}_0}^*$ locally free is unobstructed in foliated deformations in terms of cotangent sheaves if there is a foliated analytic family $\left( \mathcal{M}, \mathcal{N}_\mathcal{F}^*, B, \pi \right)$ in terms of cotangent sheaves of deformations of $\pi^{-1}(0)=\left( M, \mathcal{N}_{\mathcal{F}_0}^* \right), 0\in B$ such that the foliated Kodaira-Spencer map $\varphi_0:T_0B\to \mathbb{H}^1\left( M, \mathcal{E}_{\mathcal{N}_{\mathcal{F}_0}^*}^\bullet \right)$ in terms of cotangent sheaves is an isomorphism at $0 \in B$. Otherwise, we say that $\left( M, \mathcal{N}_{\mathcal{F}_0}^* \right)$ is obstructed in foliated deformations in terms of cotangent sheaves.
\end{definition}

\begin{remark}
We note that by theorem of existence for deformations of foliated complex analytic structures in terms of cotangent sheaves $($\textnormal{Part I}$)$, if a compact foliated complex manifold $\left(M, \mathcal{N}_{\mathcal{F}_0}^* \right)$ with $\mathcal{N}_{\mathcal{F}_0}^*$ locally free satisfies $\mathbb{H}^2\left(  M,  \mathcal{E}_{\mathcal{N}_{\mathcal{F}_0}^*}^\bullet \right)=0$, then there is a foliated complex analytic family $\left( \mathcal{M}, \mathcal{N}_\mathcal{F}^* , B, \pi \right)$ such that $\pi^{-1}(0)=\left(M, \mathcal{N}_{\mathcal{F}_0}^*\right), 0\in B$ and the foliated Kodaira-Spencer map in terms of cotangent sheaves is an isomorphism at $0\in B$, so that $\left(M, \mathcal{N}_{\mathcal{F}_0}^* \right)$ is unobstructed in foliated deformations in terms of cotangent sheaves.
\end{remark}

\begin{definition}
Let $\left(M, \mathcal{N}_{\mathcal{F}_0}^* \right)$ be a compact foliated complex manifold with $\mathcal{N}_{\mathcal{F}_0}^*$ locally free. We say that $a\in \mathbb{H}^1\left( M, \mathcal{E}_{\mathcal{N}_{\mathcal{F}_0}^*}^\bullet \right)$ is an obstructed element in foliated deformations of $\left(M, \mathcal{N}_{\mathcal{F}_0}^* \right)$ in terms of cotangent sheaves if there is no foliated analytic family $\left( \mathcal{M}, \mathcal{N}_\mathcal{F}^* , B, \pi \right)$ in terms of cotangent sheaves of deformations of $\left(M, \mathcal{N}_{\mathcal{F}_0}^* \right)= \pi^{-1}(0), 0\in B$ such that $a\in \varphi_0\left(T_0 M \right)$ where $\varphi_0: T_0 M  \to \mathbb{H}^1\left( M, \mathcal{E}_{\mathcal{N}_{\mathcal{F}_0}^*}^\bullet \right)$ is the foliated Kodaira-Spencer map at $0\in B$ of $\left(\mathcal{M}, \mathcal{N}_\mathcal{F}^* , B, \pi \right)$ in terms of cotangent sheaves.
\end{definition}

Hence if there is an obstructed element $a\in \mathbb{H}^1\left( M, \mathcal{E}_{\mathcal{N}_{\mathcal{F}_0}^*}^{\bullet} \right)$ in foliated deformations in terms of cotangent sheaves, $\left( M, \mathcal{N}_{\mathcal{F}_0}^* \right)$ is obstructed. We now describe a condition when $a\in \mathbb{H}^1\left( M, \mathcal{E}_{\mathcal{N}_{\mathcal{F}_0}^*}^\bullet \right)$ is an obstructed element in foliated deformations in terms of cotangent sheaves. We extend the arguments in \cite{Kod05} p.210-214 in the context of a foliated analytic family in terms of cotangent sheaves.

Let $\left( \mathcal{M}, \mathcal{N}_\mathcal{F}^*, B, \pi \right)$ with $0\in B \subset \mathbb{C}$ and $\mathcal{N}_{\mathcal{F}}^*$ locally free be a foliated analytic family in terms of cotangent sheaves such that $\pi^{-1}(0)=\left(M, \mathcal{N}_{\mathcal{F}_0}^* \right)$. Take a small disk $\Delta$ with centre $0$ such that $0\in \Delta \subset B$, and represent $\left(\mathcal{M}_\Delta, \mathcal{N}_{\mathcal{F}_\Delta}^* \right)=\pi^{-1}(\Delta)$ in the form
\begin{align*}
\left(\mathcal{M}_\Delta, \mathcal{N}_{\mathcal{F}_\Delta}^* \right) = \bigcup_{j=1}^l \left( U_j\times \Delta, \mathcal{N}_{\mathcal{F}_\Delta}^* \right)
\end{align*}
where each $U_j$ is a polydisk, and $\mathcal{N}_{\mathcal{F}_\Delta}^*$ on $U_j\times \Delta$ is generated by $w_j^\alpha(\xi_j,t):=\sum_{\beta=1}^n w_j^{\alpha\beta}(\xi_j,t) d\xi_j^\beta,\alpha=1,...,q$ such that $w_j^1(\xi_j,t)\wedge \cdots \wedge w_j^q(\xi_j,t)\wedge dw_j^\alpha(\xi_j,t)=0,\alpha=1,...,q$, and $(\xi_j,t)\in U_j\times \Delta$ and $(\xi_k,t)\in U_k\times \Delta$ are the same points on $\mathcal{M}_\Delta$ if $\xi_j^\alpha= f_{jk}^\alpha(\xi_k,t)$ for $\alpha=1,...,n$. Here each $f_{jk}^\alpha(\xi_k,t)=f_{jk}^\alpha(\xi_k^1,..., \xi_k^n, t)$ is a holomorphic function of $\xi_k^1,..., \xi_k^n,t$ defined on $\left(U_k\times \Delta \right) \cap \left( U_j\times \Delta \right)\ne \emptyset$. On the other hand, $w_j^\alpha(\xi_j,t)= \sum_{\beta=1}^q h_{jk}^{\alpha\beta}\left(\xi_j,t \right) w_k^\beta(\xi_k,t),\alpha=1,...,q$ on $\left(U_k\times \Delta \right) \cap \left( U_j\times \Delta \right)\ne \emptyset$, and $h_{ik}^{\alpha\gamma}\left(\xi_k,t \right)= \sum_{\beta=1}^p h_{ij}^{\alpha\beta}(\xi_j,t) h_{jk}^{\beta \gamma}(\xi_k,t)$ on $\left( U_i\times \Delta \right) \cap \left( U_j\times \Delta \right) \cap \left( U_k\times \Delta \right)\ne \emptyset$.

As in $(\ref{oc7})$, the infinitesimal foliated deformation $\frac{\partial \left( M_t, \Theta_{\mathcal{F}_t} \right)}{\partial t}\in \mathbb{H}^1\left( M_t, \mathcal{E}_{\Theta_{\mathcal{F}_t} }^\bullet \right)$ is the cohomology class of the $1$-cocycle $\left( \left\{\beta_j(t)\right\}, \left\{\left(\theta_{jk}(t), B_{jk}(t) \right) \right\} \right)\in C^0\left( \mathcal{U}_t, \mathscr{H}om_{\mathcal{O}_{M_t}}\left( \mathcal{N}_{\mathcal{F}_t}^*, \Omega_{M_t}^1 \right) \right)\bigoplus C^1\left( \mathcal{U}_t, \mathcal{E}_{\mathcal{N}_{\mathcal{F}_t}^*} \right)$, where $\mathcal{U}_t=\left\{U_j^t=U_j\times t \right\}$ such that
\begin{align*}
\beta_j(t):\Gamma\left( U_j^t, \mathcal{N}_{\mathcal{F}_t}^* \right)&\to \Gamma\left( U_j^t, \Omega_{M_t}^1 \right) \\
 w_j^\alpha(\xi_j,t) &\mapsto  \beta_j^\alpha(\xi_j,t):=\sum_{\beta=1}^n \beta_j^{\alpha\beta}(\xi_j,t)d\xi_j^\beta ,\,\,\,\,\,\,\,\,\textnormal{where}\,\,\,\beta_j^{\alpha\beta}(\xi_j, t)= - \frac{\partial w_j^{\alpha\beta}(\xi_j,t)}{\partial t},\alpha=1,...,q
\end{align*}
and the vector field
\begin{align*}
\theta_{jk}(t)=\sum_{\alpha=1}^n \theta_{jk}^\alpha(\xi_j,t)\frac{\partial}{\partial \xi_j^\alpha},\,\,\,\,\,\,\,\,\textnormal{where}\,\,\,\,\theta_{jk}^\alpha(\xi_j,t)=\frac{\partial f_{jk}(\xi_k,t)}{\partial t},\,\,\,\xi_k=f_{kj}(\xi_j,t)
\end{align*}
and
\begin{align*}
B_{jk}(t):\Gamma\left( U_j^t, \mathcal{N}_{\mathcal{F}_t}^* \right)&\to \Gamma\left( U_j^t, \mathcal{N}_{\mathcal{F}_t}^* \right) \\
     w_j^\alpha(\xi_j,t)&\mapsto \sum_{\gamma=1}^p B_{jk}^{\alpha\gamma}(\xi_j,t)w_j^\gamma(\xi_j,t),\,\,\textnormal{where}\,\,\, B_{jk}^{\alpha\gamma}(\xi_j,t)=\sum_{\beta=1}^q\frac{\partial h_{jk}^{\alpha\beta}(\xi_k,t)}{\partial t}h_{kj}^{\beta \gamma}(\xi_j,t),\xi_k=f_{kj}(\xi_j,t),\alpha=1,...,q
\end{align*}

On $\left(U_k\times \Delta \right)\cap \left(U_j\times \Delta \right) \ne \emptyset$, we have the equalities
\begin{align}
f_{ik}^\alpha(\xi_k,t)&=f_{ij}^\alpha\left(f_{jk}(\xi_k,t),t \right),\,\,\,\alpha=1,...,n \label{oc13}\\
\sum_{\gamma=1}^n w_j^{\alpha\gamma}(f_{jk}(\xi_k,t),t)\frac{\partial f_{jk}^\gamma(\xi_k,t)}{\partial \xi_k^\eta}&=\sum_{\beta=1}^qh_{jk}^{\alpha\beta}(\xi_k,t)w_k^{\beta\eta}(\xi_k,t),\,\,\,\alpha=1,...,q,\,\,\, \eta=1,...,n \label{oc10}
\end{align}
and on $\left(U_i\cap \Delta \right)\cap \left(  U_j \times \Delta \right)\cap \left(U_k\times \Delta \right)\ne \emptyset$, we have the equalities
\begin{align}
h_{ik}^{\alpha\gamma}\left(\xi_k,t \right) &= \sum_{\beta=1}^q h_{ij}^{\alpha\beta}\left(f_{jk}(\xi_k,t),t \right) h_{jk}^{\beta \gamma}(\xi_k,t),\,\,\,\,\,\alpha,\gamma=1,...,q. \label{oc11}
\end{align}
and on $U_j\times \Delta$, we have the equalities 
\begin{align}
w_j^1(\xi_j,t)\wedge \cdots \wedge w_j^q(\xi_j,t)\wedge dw_j^\alpha(\xi_j,t)=0,\,\,\,\,\,\,\alpha=1,...,q \label{oc12}
\end{align}
By differentiating both sides of $(\ref{oc13}),(\ref{oc10}),(\ref{oc11})$ and $(\ref{oc12})$ in $t$, we have the equalities
{\small{\begin{align}
&\theta_{ik}^\alpha(\xi_i,t)=\theta_{ij}^\alpha(\xi_i,t) + \sum_{\beta=1}^n \frac{\partial \xi_i^\alpha}{\partial \xi_j^\beta}\theta_{jk}^\beta(\xi_j,t),\,\,\,\alpha=1,...,n. \label{oc14}\\
&\sum_{\gamma,\delta=1}^n \frac{\partial w_j^{\alpha\gamma}}{\partial \xi_j^\delta}\left(f_{jk}(\xi_k,t),t) \right)\frac{\partial f_{jk}^\delta}{\partial t}\frac{\partial f_{jk}^\gamma}{\partial \xi_k^\eta}+\sum_{\gamma=1}^n \frac{\partial w_j^{\alpha\gamma}}{\partial t}\left(f_{jk}(\xi_k,t),t \right)\frac{\partial f_{jk}^\gamma}{\partial \xi_k^\eta} +\sum_{\gamma=1}^n w_j^{\alpha\gamma}\left(f_{jk}(\xi_k,t),t \right)\frac{\partial}{\partial \xi_k^\eta}\left(\frac{\partial f_{jk}^\gamma}{\partial t} \right) \label{oc15}\\
&= \sum_{\beta=1}^q \frac{\partial h_{jk}^{\alpha\beta}}{\partial t}w_k^{\beta\eta} +\sum_{\beta=1}^q h_{jk}^{\alpha\beta}\frac{\partial w_k^{\beta\eta}}{\partial t} ,\,\,\,\,\,\,\,\,\alpha=1,...,q, \,\,\, \eta=1,...,n\notag
\end{align}}}
and for $\alpha,\gamma=1,...,q$,
{\small{\begin{align}
\frac{\partial h_{ik}^{\alpha \gamma}}{\partial t} = \sum_{\beta=1}^q \sum_{\eta=1}^n \frac{\partial h_{ij}^{\alpha\beta}}{\partial \xi_j^\eta}\left(f_{jk}(\xi_k, t) , t\right) \frac{\partial f_{jk}^{\eta}}{\partial t} h_{jk}^{\beta \gamma} + \sum_{\beta=1}^q \frac{\partial h_{ij}^{\alpha\beta}}{\partial t}\left(f_{jk}(\xi_k,t),t\right) h_{jk}^{\beta \gamma} + \sum_{\beta=1}^q h_{ij}^{\alpha\beta}\left(f_{jk}(\xi_k,t),t\right)\frac{\partial h_{jk}^{\beta \gamma}}{\partial t} \label{oc16}
\end{align}}}
\begin{align}
w_j^1\wedge \cdots \wedge w_j^q \wedge d\left(\frac{\partial w_j^\alpha}{\partial t} \right)+\sum_{\gamma=1}^q w_j^1\wedge \cdots \wedge \overbrace{\frac{\partial w_j^\gamma}{\partial t}}^{\gamma-\textnormal{th}} \wedge \cdots \wedge w_j^q\wedge dw_j^\alpha=0 ,\,\,\,\alpha=1,...,q\label{oc17}
\end{align}

Next we will take the derivative of $(\ref{oc14}), (\ref{oc15}), (\ref{oc17})$ in $t$ as functions of $\xi_j^1,..., \xi_j^n,t$ and $(\ref{oc16})$ in $t$ as functions of $\xi_k^1,..., \xi_k^n,t$ again. We recall $(\ref{o53})$ in subsection \ref{o42}.
By taking the derivative of $(\ref{oc14})$ with respect to $t$, and putting
\begin{align*}
\dot{\theta}_{ij}(t)= \sum_{\alpha=1}^n \dot{\theta}_{ij}^\alpha (z_i,t)\frac{\partial}{\partial z_i^\alpha} , \,\,\,\,\,\textnormal{where}\,\,\,\dot{\theta}_{ij}^\alpha(z_i,t)=\frac{\partial \theta_{ij}^\alpha(z_i,t)}{\partial t}
\end{align*}
we obtain (for the detail, see \cite{Kod05} p.211-213)
\begin{align}\label{oc18}
\dot{\theta}_{ij}(t) -\dot{\theta}_{ik}(t) + \dot{\theta}_{jk} (t) = \left[ \theta_{ij}(t), \theta_{jk}(t) \right]
\end{align}
By taking the derivative of $(\ref{oc17})$ with respect to $t$, and putting
\begin{align*}
\dot{\beta}_j^{\alpha}(\xi_j,t)= \frac{\partial \beta_j^{\alpha}(\xi_j,t)}{\partial t}= \sum_{\beta=1}^n \frac{\partial \beta_j^{\alpha\beta}(\xi_j,t)}{\partial t} d\xi_j^\beta = - \sum_{\beta=1}^n \frac{\partial^2 w_j^{\alpha\beta}(\xi_j,t)}{\partial t^2}d\xi_j^\beta, \,\,\,\alpha=1,...,q
\end{align*}
we obtain 
\begin{align} \label{oc28}
&w_j^1\wedge \cdots \wedge w_j^q\wedge d\left(\dot{\beta}_j^\alpha(t) \right)+\sum_{\gamma=1}^q w_j^1\wedge \cdots \wedge\overbrace{\dot{\beta}_j^\gamma(t)}^{\gamma-\textnormal{th}}\wedge \cdots \wedge w_j^q\wedge dw_j^\alpha\\
&=2\sum_{\gamma=1}^q w_j^1\wedge \cdots \wedge \overbrace{\beta_j^\gamma(t)}^{\gamma-\textnormal{th}} \wedge \cdots \wedge w_j^q\wedge d\left(\beta_j^\alpha(t)\right)+2\sum_{1\leq \gamma<\eta\leq q} w_j^1\wedge \cdots \wedge \overbrace{\beta_j^\gamma(t) }^{\gamma-\textnormal{th}} \wedge \cdots \wedge \overbrace{\beta_j^\eta(t)}^{\eta-\textnormal{th}} \wedge \cdots \wedge w_j^q \wedge dw_j^\alpha \notag
\end{align}

Let us define $\dot{\beta}_j(t) \in \Gamma\left( U_j^t, \mathscr{H}om_{\mathcal{O}_{M_t}}\left( \mathcal{N}_{\mathcal{F}_t}^* , \Omega_{M_t}^1\right) \right)$ by 
\begin{align*}
\dot{\beta}_j(t):\Gamma\left(U_j^t, \mathcal{N}_{\mathcal{F}_t}^* \right) &\to \Gamma\left( U_j^t, \Omega_{M_t}^1 \right) \\
 w_j^\alpha &\mapsto \dot{\beta}_j(t)= \dot{\beta}_j(\xi_j,t)
\end{align*}
and linearly extend to $\Gamma\left( U_j^t, \mathcal{N}_{\mathcal{F}_t}^* \right)$. Let us define 
\begin{align}
\dot{B}_{jk}(t):\Gamma\left( U_{jk}^t,  \mathcal{N}_{\mathcal{F}_t}^* \right)\to &\Gamma\left( U_{jk}^t, \mathcal{N}_{\mathcal{F}_t}^* \right) \label{oc19} \\
 w_j^\alpha \mapsto &\sum_{\gamma=1}^q \frac{\partial}{\partial t_j}\left( B_{jk}^{\alpha \gamma}(\xi_j,t)\right) w_j^\gamma = \sum_{\beta, \gamma=1}^q \frac{\partial}{\partial t_j}\left( \frac{\partial h_{jk}^{\alpha \beta}(\xi_k,t)}{\partial t} h_{kj}^{\beta \gamma}(\xi_j,t) \right) w_j^\gamma  \notag
\end{align}
and linearly extends to $\Gamma\left( U_{jk}^t, \mathcal{N}_{\mathcal{F}_t}^* \right)$. By considering {\small{$\left(\dot{\theta}_{jk}(t), \dot{B}_{jk}(t) \right)\in \Gamma\left(U_{jk}^t, \Theta_{M_t} \right)|_{U_j^t} \bigoplus \Gamma\left(U_{jk}^t, \mathscr{H}om_{\mathcal{O}_{M_t}} \left( \mathcal{N}_{\mathcal{F}_t}^*, \mathcal{N}_{\mathcal{F}_t}^* \right) |_{U_j^t}\right)$}}, we see that $\left\{\left( \dot{\theta}_{jk}(t), \dot{B}_{jk}(t) \right) \right\} \in C^1\left(\mathcal{U}_t, \mathcal{E}_{ \mathcal{N}_{\mathcal{F}_t}^*} \right)$ as in the same arguments with the proof of $\left\{\left( \dot{\theta}_{jk}(t), \dot{A}_{jk}(t) \right) \right\} \in C^1\left(\mathcal{U}_t, \mathcal{E}_{ \Theta_{\mathcal{F}_t}} \right)$ in subsection \ref{o42}. On the other hand, let us take the derivative of $(\ref{oc15})$ with respect to $t_k$:
\begin{align}\label{oc20}
&\sum_{\gamma,\sigma,\delta=1}^n\frac{\partial^2 w_j^{\alpha\gamma}}{\partial \xi_j^\sigma \xi_j^\delta}\frac{\partial f_{jk}^\sigma}{\partial t} \frac{\partial f_{jk}^\delta}{\partial t} \frac{\partial f_{jk}^\gamma}{\partial \xi_k^\eta}+2\sum_{\gamma,\delta=1}^n \frac{\partial}{\partial \xi_j^\delta}\left(\frac{\partial w_j^{\alpha\gamma}}{\partial t_j} \right) \frac{\partial f_{jk}^\delta}{\partial t}\frac{\partial f_{jk}^\gamma}{\partial \xi_k^\eta} +\sum_{\gamma,\delta=1}^n \frac{\partial w_j^{\alpha\gamma}}{\partial \xi_j^\delta}\frac{\partial^2 f_{jk}^\delta}{\partial t^2}\frac{\partial f_{jk}^\gamma}{\partial \xi_k^\eta}+ 2 \sum_{\gamma,\delta=1}^n \frac{\partial w_j^{\alpha\gamma}}{\partial \xi_j^\delta}\frac{\partial f_{jk}^\delta }{\partial t}\frac{\partial}{\partial \xi_k^\eta}\left(\frac{\partial f_{jk}^\gamma}{\partial t} \right)\\
&+\sum_{\gamma=1}^n \frac{\partial^2 w_j^{\alpha\gamma}}{\partial t_j^2}\frac{\partial f_{jk}^\gamma}{\partial \xi_k^\eta}+2\sum_{\gamma=1}^n \frac{\partial w_j^{\alpha\gamma}}{\partial t_j}\frac{\partial}{\partial \xi_k^\eta}\left(\frac{\partial f_{jk}^\gamma}{\partial t} \right) +\sum_{\gamma=1}^n w_j^{\alpha\gamma}\frac{\partial}{\partial \xi_k^\eta}\left(\frac{\partial^2 f_{jk}^\gamma}{\partial t^2} \right) \notag \\
&= \sum_{\beta=1}^q \frac{\partial^2 h_{jk}^{\alpha\beta}}{\partial t^2}w_k^{\beta\eta} + \sum_{\beta=1}^q \frac{\partial h_{jk}^{\alpha\beta}}{\partial t} \frac{\partial w_k^{\beta\eta}}{\partial t} +\sum_{\beta=1}^q \frac{\partial h_{jk}^{\alpha\beta}}{\partial t} \frac{\partial w_k^{\beta\eta}}{\partial t}+\sum_{\beta=1}^q h_{jk}^{\alpha\beta} \frac{\partial^2 w_k^{\beta\eta}}{\partial t^2} \notag
\end{align}
By summing $\sum_{\eta=1}^n d\xi_k^\eta$ on both sides of $(\ref{oc20})$ and by noting $\dot{\theta}_{jk}^\beta=\frac{\partial f_{jk}^\beta}{\partial t^2} - \sum_{\sigma=1}^n \frac{\partial f_{jk}^\sigma}{\partial t} \frac{\partial}{\partial \xi_j^\sigma}\left(\frac{\partial f_{jk}^\beta}{\partial t} \right)$, we have
\begin{align} \label{oc21}
&\sum_{\gamma,\sigma,\delta=1}^n\frac{\partial^2 w_j^{\alpha\gamma}}{\partial \xi_j^\sigma \xi_j^\delta}\frac{\partial f_{jk}^\sigma}{\partial t} \frac{\partial f_{jk}^\delta}{\partial t}d\xi_j^\gamma +2\sum_{\gamma,\delta=1}^n \frac{\partial}{\partial \xi_j^\delta}\left(\frac{\partial w_j^{\alpha\gamma}}{\partial t_j} \right) \frac{\partial f_{jk}^\delta}{\partial t}d\xi_j^\gamma +\sum_{\gamma,\delta=1}^n \frac{\partial w_j^{\alpha\gamma}}{\partial \xi_j^\delta}\left( \dot{\theta}_{jk}^\delta+ \sum_{\sigma=1}^n \frac{\partial f_{jk}^\sigma}{\partial t} \frac{\partial}{\partial \xi_j^\sigma}\left(\frac{\partial f_{jk}^\delta}{\partial t} \right) \right) d\xi_j^\gamma\\
&+ 2 \sum_{\gamma,\delta=1}^n \frac{\partial w_j^{\alpha\gamma}}{\partial \xi_j^\delta}\frac{\partial f_{jk}^\delta }{\partial t} d\left(\frac{\partial f_{jk}^\gamma}{\partial t} \right)+\sum_{\gamma=1}^n \frac{\partial^2 w_j^{\alpha\gamma}}{\partial t_j^2}d\xi_j^\gamma +2\sum_{\gamma=1}^n \frac{\partial w_j^{\alpha\gamma}}{\partial t_j} d \left(\frac{\partial f_{jk}^\gamma}{\partial t} \right) +\sum_{\gamma=1}^n w_j^{\alpha\gamma}d \left( \dot{\theta}_{jk}^\gamma + \sum_{\sigma=1}^n \frac{\partial f_{jk}^\sigma}{\partial t} \frac{\partial}{\partial \xi_j^\sigma}\left(\frac{\partial f_{jk}^\gamma}{\partial t} \right) \right) \notag \\
&= \sum_{\beta=1}^q \frac{\partial^2 h_{jk}^{\alpha\beta}}{\partial t^2}w_k^\beta  + 2\sum_{\beta,\eta=1}^q \frac{\partial h_{jk}^{\alpha\beta}}{\partial t} \frac{\partial w_k^{\beta\eta}}{\partial t} dz_k^\eta + \sum_{\beta , \eta =1}^q h_{jk}^{\alpha\beta} \frac{\partial^2 w_k^{\beta\eta}}{\partial t^2} dz_k^\eta \notag
\end{align}
Then $\left( \ref{oc21} \right)$ is equivalent to
\begin{align}\label{oc22}
&\sum_{\gamma,\sigma,\delta=1}^n\frac{\partial^2 w_j^{\alpha\gamma}}{\partial \xi_j^\sigma \xi_j^\delta}\frac{\partial f_{jk}^\sigma}{\partial t} \frac{\partial f_{jk}^\delta}{\partial t}d\xi_j^\gamma + 2 \mathcal{L}_{\theta_{jk}(t)}\left(  \frac{\partial w_j^\alpha}{\partial t}\right) +\mathcal{L}_{\dot{\theta}_{jk}(t)}\left(w_j^\alpha  \right) + \sum_{\gamma, \delta, \sigma=1}^n \frac{\partial w_j^{\alpha \gamma}}{\partial \xi_j^\delta} \frac{\partial f_{jk}^\sigma}{\partial t} \frac{\partial}{\partial \xi_j^\sigma}\left( \frac{\partial f_{jk}^\delta}{\partial t} \right) d\xi_j^\gamma \\
&+ 2 \sum_{\gamma,\delta=1}^n \frac{\partial w_j^{\alpha\gamma}}{\partial \xi_j^\delta}\frac{\partial f_{jk}^\delta }{\partial t} d\left(\frac{\partial f_{jk}^\gamma}{\partial t} \right)+\sum_{\gamma=1}^n \frac{\partial^2 w_j^{\alpha\gamma}}{\partial t_j^2}d\xi_j^\gamma + \sum_{\gamma, \sigma=1}^n w_j^{\alpha \gamma} d\left( \frac{\partial f_{jk}^{\sigma}}{\partial t} \frac{\partial}{\partial \xi_j^\sigma}\left(\frac{\partial f_{jk}^\gamma}{\partial t} \right)   \right) \notag \\
&= \sum_{\beta=1}^q \frac{\partial^2 h_{jk}^{\alpha\beta}}{\partial t^2}w_k^\beta  + 2\sum_{\beta,\eta=1}^q \frac{\partial h_{jk}^{\alpha\beta}}{\partial t} \frac{\partial w_k^{\beta\eta}}{\partial t} dz_k^\eta + \sum_{\beta , \eta =1}^q h_{jk}^{\alpha\beta} \frac{\partial^2 w_k^{\beta\eta}}{\partial t^2} dz_k^\eta \notag
\end{align}

We also note that
{\Small{\begin{align*}
&\mathcal{L}_{\theta_{jk}}\mathcal{L}_{\theta_{jk}}\left(w_j^\alpha \right) =\mathcal{L}_{\sum_{\sigma=1}^n \frac{\partial f_{jk}^\sigma}{\partial t} \frac{\partial}{\partial \xi_j^\sigma}} \left( \mathcal{L}_{\sum_{\delta=1}^n \frac{\partial f_{jk}^\delta}{\partial t} \frac{\partial}{\partial \xi_j^\delta}} \left(\sum_{\gamma=1}^n w_j^{\alpha \gamma} d\xi_j^\gamma \right)   \right) =\mathcal{L}_{\sum_{\sigma=1}^n \frac{\partial f_{jk}^\sigma}{\partial t} \frac{\partial}{\partial \xi_j^\sigma}} \left( \sum_{\delta, \gamma=1}^n \frac{\partial f_{jk}^\delta}{\partial t} \frac{\partial w_j^{\alpha \gamma}}{\partial \xi_j^\delta} d\xi_j^\gamma + \sum_{\gamma=1}^n w_j^{\alpha \gamma} d\left(\frac{\partial f_{jk}^\gamma}{\partial t} \right)   \right) \\
&= \sum_{\sigma,\delta, \gamma=1}^n \frac{\partial f_{jk}^\sigma}{\partial t}\frac{\partial }{\partial \xi_j^\sigma}\left( \frac{\partial f_{jk}^\delta}{\partial t}\right)\frac{\partial w_j^{\alpha \gamma}}{\partial \xi_j^\delta} d\xi_j^\gamma + \sum_{\sigma, \delta, \gamma=1}^n \frac{\partial f_{jk}^\sigma}{\partial t} \frac{\partial f_{jk}^\delta}{\partial t} \frac{\partial^2 w_j^{\alpha \gamma}}{\partial \xi_j^\sigma \partial \xi_j^\delta} d\xi_j^\gamma + 2 \sum_{\delta, \gamma=1}^n \frac{\partial f_{jk}^\delta}{\partial t}\frac{\partial w_j^{\alpha \gamma}}{\partial \xi_j^\delta}d\left(\frac{\partial f_{jk}^\gamma}{\partial t} \right) + \sum_{\gamma,\sigma=1}^n w_j^{\alpha \gamma}d\left(\frac{\partial f_{jk}^\sigma}{\partial t}\frac{\partial }{\partial \xi_j^\sigma}\left(\frac{\partial f_{jk}^\gamma}{\partial t} \right) \right)
\end{align*}}}

This implies that $(\ref{oc22})$ is equivalent to
\begin{align}\label{oc23}
 2 \mathcal{L}_{\theta_{jk}(t)}\left(  \frac{\partial w_j^\alpha}{\partial t}\right) +\mathcal{L}_{\dot{\theta}_{jk}(t)}\left(w_j^\alpha  \right)  + \frac{\partial^2 w_j^{\alpha}}{\partial t_j^2} + \mathcal{L}_{\theta_{jk}}\mathcal{L}_{\theta_{jk}}\left(w_j^\alpha \right) = \sum_{\beta=1}^q \frac{\partial^2 h_{jk}^{\alpha\beta}}{\partial t^2}w_k^\beta  + 2\sum_{\beta=1}^q \frac{\partial h_{jk}^{\alpha\beta}}{\partial t} \frac{\partial w_k^\beta}{\partial t} + \sum_{\beta =1}^q h_{jk}^{\alpha\beta} \frac{\partial^2 w_k^\beta }{\partial t^2}  
\end{align}

\begin{lemma}
\begin{align}
&\dot{\beta}_k(t)-\dot{\beta}_j(t)+ E_0'\left( \left( \dot{\theta}_{jk} (t) , \dot{B}_{jk} (t) \right) \right) \label{oc24} \\
&= - \beta_j(t)\circ B_{jk}(t) - \beta_k (t) \circ B_{jk}(t)+ \mathcal{L}_{\theta_{jk}(t)}\circ \beta_j(t)+ \mathcal{L}_{\theta_{jk}(t)}\circ \beta_k(t) \notag
\end{align}
\end{lemma}

\begin{proof}
Let us compute the left-hand side of $(\ref{oc24})$. We note that we have the similar equality with $(\ref{o24})$ in terms of $h_{jk}^{\alpha\beta}$. Then
{\Small{\begin{align} \label{oc25}
&\left(\dot{\beta}_k(t)- \dot{\beta}_j(t)+ \mathcal{L}_{\dot{\theta}_{jk}(t)}- \dot{B}_{jk} \right)(w_j^\alpha)\\
&= - \sum_{\beta=1}^q h_{jk}^{\alpha\beta} \frac{\partial^2 w_k^\beta}{\partial t_k^2} + \frac{\partial^2 w_j^\beta}{\partial t_j^2} + \mathcal{L}_{\dot{\theta}_{jk}(t)}(w_j^\alpha) - \sum_{\beta, \gamma=1}^q \frac{\partial}{\partial t_j}\left( \frac{\partial h_{jk}^{\alpha\beta}}{\partial t} h_{kj}^{\beta \gamma} \right) w_j^\gamma \notag  \\
&= - \sum_{\beta=1}^q h_{jk}^{\alpha\beta} \frac{\partial^2 w_k^\beta}{\partial t_k^2} + \frac{\partial^2 w_j^\beta}{\partial t_j^2} + \mathcal{L}_{\dot{\theta}_{jk}(t)}(w_j^\alpha)  -\sum_{\beta=1}^q \frac{\partial^2 h_{jk}^{\alpha\beta}}{\partial t^2} w_k^\beta + \sum_{\beta=1}^q \left[ \theta_{jk}, \frac{\partial h_{jk}^{\alpha\beta}}{\partial t}\right] w_{k}^\beta + \sum_{\beta, \gamma=1}^q \frac{\partial h_{jk}^{\alpha\beta}}{\partial t}\left[ \theta_{jk}, h_{kj}^{\beta\gamma}\right] w_{j}^\gamma + \sum_{\beta, \gamma, \delta=1}^q  \frac{\partial h_{jk}^{\alpha\beta}}{\partial t} h_{kj}^{\beta \gamma} \frac{\partial h_{jk}^{\gamma \delta}}{\partial t} w_{k}^\delta \notag
\end{align}}}

On the other hand, let us compute right-hand side of $(\ref{oc24})$. Since $\beta_k(t)- \beta_j(t)+\left(\mathcal{L}_{\theta_{jk}(t)}- B_{jk}(t) \right)=0$, 
{\small{\begin{align} \label{oc26}
&\left( - \beta_j(t)\circ B_{jk}(t) - \beta_k(t) \circ B_{jk}(t)+ \mathcal{L}_{\theta_{jk}(t)}\circ \beta_j(t)+ \mathcal{L}_{\theta_{jk}(t)}\circ \beta_k(t) \right)(w_j^\alpha)\\
&=\left( - 2\beta_k(t)\circ B_{jk}(t) -\left( \mathcal{L}_{\theta_{jk}(t)}- B_{jk}(t) \right)\circ B_{jk}(t) + 2 \mathcal{L}_{\theta_{jk}(t)}\circ \beta_j(t) -  \mathcal{L}_{\theta_{jk}(t)} \circ \left(\mathcal{L}_{\theta_{jk}(t)}- B_{jk}(t) \right)\right)(w_j^\alpha) \notag \\
&= -2 \beta_k(t)\circ B_{jk}(t)(w_j^\alpha) - \sum_{\beta, \gamma=1}^q \frac{\partial h_{jk}^{\alpha\beta}}{\partial t} h_{kj}^{\beta \gamma}\left( \mathcal{L}_{\theta_{jk}(t)}- B_{jk}(t) \right)(w_j^\gamma) + 2 \mathcal{L}_{\theta_{jk}(t)}\circ \beta_j(t) (w_j^\alpha) - \mathcal{L}_{\theta_{jk}(t)} \circ \mathcal{L}_{\theta_{jk}}(t)(w_j^\alpha) + \mathcal{L}_{\theta_{jk}(t)}\circ B_{jk}(t)(w_j^\alpha) \notag \\
&= 2 \sum_{\beta=1}^q \frac{\partial h_{jk}^{\alpha\beta}}{\partial t}\frac{\partial w_k^\beta}{\partial t}  - \sum_{\beta=1}^q \frac{\partial h_{jk}^{\alpha\beta}}{\partial t} \mathcal{L}_{\theta_{jk}(t)} \left(w_k^\beta\right) + \sum_{\beta, \gamma=1}^q \frac{\partial h_{jk}^{\alpha\beta}}{\partial t} \left[ \theta_{jk}(t), h_{kj}^{\beta \gamma} \right] w_j^\gamma + \sum_{\beta, \gamma, \delta=1}^q  \frac{\partial h_{jk}^{\alpha\beta}}{\partial t} h_{kj}^{\beta \gamma} \frac{\partial h_{jk}^{\gamma \delta}}{\partial t} w_{k}^\delta - 2 \mathcal{L}_{\theta_{jk}(t)}\left(\frac{\partial w_j^\alpha}{\partial t}\right) \notag \\
&- \mathcal{L}_{\theta_{jk}(t)}\circ \mathcal{L}_{\theta_{jk}(t)}(w_j^\alpha) +  \sum_{\beta=1}^q\mathcal{L}_{\theta_{jk}(t)}\left( \frac{\partial h_{jk}^{\alpha\beta}}{\partial t} w_k^\beta \right) \notag
\end{align}}}
Then $(\ref{oc23}),(\ref{oc25})$ and $(\ref{oc26})$ implies $(\ref{oc24})$.

\end{proof}

Next, by taking the derivative of $(\ref{oc16})$ with respect to $t_k$ again and by the same arguments of Lemma \ref{o52}, we obtain
\begin{align} \label{oc27}
\dot{B}_{ij} (t) - \dot{B}_{ik} (t) + \dot{B}_{jk}^i  (t) = \left[ \theta_{ij} (t) , B_{jk}^i(t)\right] - \left[ \theta_{jk}(t), B_{ij}(t)\right] + \left[ B_{ij} (t) ,   B_{jk}^i(t)\right]
\end{align}
for the definition of the bracket $[-,-]$ in $(\ref{oc27})$, see \textnormal{Appendix \ref{apd2}}.

From $(\ref{oc18}),(\ref{oc28}),(\ref{oc24})$ and $(\ref{oc27})$, substituting $t=0$, we obtain
\begin{align}\label{oc30}
\dot{\theta}_{ij}(0)- \dot{\theta}_{ik}(0) + \dot{\theta}_{jk}(0)= \left[ \theta_{ij}(0),\theta_{jk}(0) \right]
\end{align}
{\small{\begin{align} \label{oc31}
&w_{0j}^1\wedge \cdots \wedge w_{0j}^q\wedge d\left(\dot{\beta}_j^\alpha(0) \right)+\sum_{\gamma=1}^q w_{0j}^1\wedge \cdots \wedge\overbrace{\dot{\beta}_j^\gamma(0)}^{\gamma-\textnormal{th}}\wedge \cdots \wedge w_{0j}^q\wedge dw_{0j}^\alpha\\
&=2\sum_{\gamma=1}^q w_{0j}^1\wedge \cdots \wedge \overbrace{\beta_j^\gamma(0)}^{\gamma-\textnormal{th}} \wedge \cdots \wedge w_{0j}^q\wedge d\left(\beta_j^\alpha(0)\right)+2\sum_{1\leq \gamma<\eta\leq q} w_{0j}^1\wedge \cdots \wedge \overbrace{\beta_j^\gamma(0) }^{\gamma-\textnormal{th}} \wedge \cdots \wedge \overbrace{\beta_j^\eta(0)}^{\eta-\textnormal{th}} \wedge \cdots \wedge w_{0j}^q \wedge dw_{0j}^\alpha \notag
\end{align}}}
\begin{align}\label{oc32}
\dot{\beta}_k(0)-\dot{\beta}_j(0)+ E_0'\left( \left( \dot{\theta}_{jk} (0) , \dot{B}_{jk} (0) \right) \right) = - \beta_j(0)\circ B_{jk}(0) - \beta_k (0) \circ B_{jk}(0)+ \mathcal{L}_{\theta_{jk}(0)}\circ \beta_j(0)+ \mathcal{L}_{\theta_{jk}(0)}\circ \beta_k(0)
\end{align}
\begin{align} \label{oc33}
\dot{B}_{ij} (0) - \dot{B}_{ik} (0) + \dot{B}_{jk}^i  (0) = \left[ \theta_{ij} (0) , B_{jk}^i(0)\right] - \left[ \theta_{jk}(0), B_{ij}(0)\right] + \left[ B_{ij} (0) ,   B_{jk}^i(0)\right]
\end{align}

Let us define $\left[ \beta_j(0), \beta_j(0)\right]_{\mathcal{N}_{\mathcal{F}_0}^*} \in \Gamma\left( U_j^0, \mathscr{H}om_{\mathcal{O}_M}\left( \mathcal{N}_{\mathcal{F}_0}^*, \tilde{\mathcal{S}}^2 \right) \right)$ as follows: 
{\small{\begin{align*}
\left[ \beta_j(0), \beta_j(0)\right]&_{\mathcal{N}_{\mathcal{F}_0}^*}: \Gamma\left( U_j^0, \mathcal{N}_{\mathcal{F}_0}^*\right) \to \Gamma\left( U_j^0, \tilde{\mathcal{S}}^2 \right) \\
 w_{0j}^\alpha &\mapsto 2\sum_{\gamma=1}^q w_{0j}^1\wedge \cdots \wedge \overbrace{\beta_j^\gamma(0)}^{\gamma-\textnormal{th}} \wedge \cdots \wedge w_{0j}^q\wedge d\left(\beta_j^\alpha(0)\right)+2\sum_{1\leq \gamma<\eta\leq q} w_{0j}^1\wedge \cdots \wedge \overbrace{\beta_j^\gamma(0) }^{\gamma-\textnormal{th}} \wedge \cdots \wedge \overbrace{\beta_j^\eta(0)}^{\eta-\textnormal{th}} \wedge \cdots \wedge w_{0j}^q \wedge dw_{0j}^\alpha
\end{align*}}}
and linearly extends to $\Gamma\left( U_j^0, \mathcal{N}_{\mathcal{F}_0}^* \right)$. Then by combining $(\ref{oc30})$ and $(\ref{oc33})$ (see Appendix \ref{apd2}), we can write $(\ref{oc30})- (\ref{oc33})$ as follows:
\begin{align}\label{oc34}
&E_1'\left( \dot{\beta}_j(0)\right)= \left[ \beta_j(0), \beta_j(0)\right]_{\mathcal{N}_{\mathcal{F}_0}^*}\\
\dot{\beta}_k(0)-\dot{\beta}_j(0)+ E_0'\left( \left( \dot{\theta}_{jk} (0) , \dot{B}_{jk} (0) \right) \right) & = - \beta_j(0)\circ B_{jk}(0) - \beta_k (0) \circ B_{jk}(0)+ \mathcal{L}_{\theta_{jk}(0)}\circ \beta_j(0)+ \mathcal{L}_{\theta_{jk}(0)}\circ \beta_k(0)  \notag \\
\left( \dot{\theta}_{ij}(0), \dot{B}_{ij}(0)\right)- &\left( \dot{\theta}_{ik}(0), \dot{B}_{ik}(0) \right) + \left( \dot{\theta}_{jk}(0), \dot{B}_{jk}(0) \right)= \left[  \left( \theta_{jk} (0) , B_{jk}(0) \right) , \left( \theta_{jk}(0), B_{jk}(0) \right)   \right] \notag
\end{align}

\begin{lemma}\label{oc29}
For any $1$-cocycle $\left( \left\{\beta_j \right\} , \left\{\left( \theta_{jk}, B_{jk} \right) \right\} \right)\in C^0\left(\mathcal{U},  \mathscr{H}om_{\mathcal{O}_M}\left( \mathcal{N}_{\mathcal{F}_0}^*, \Omega_M^1 \right) \right) \bigoplus C^1\left(  \mathcal{U}, \mathcal{E}_{\mathcal{N}_{\mathcal{F}_0}^*} \right)$ where $\beta_j:\Gamma\left( U_j, \mathcal{N}_{\mathcal{F}_0}^* \right)\to \Gamma\left( U_j, \Omega_M^1 \right) , \theta_{jk}\in\Gamma\left( U_{jk} , \Theta_M \right)$ and $B_{jk}:\Gamma\left(U_{jk}, \mathcal{N}_{\mathcal{F}_0}^* \right) \to \Gamma\left( U_{jk}, \mathcal{N}_{\mathcal{F}_0}^* \right)$ satisfying
\begin{enumerate}
\item $E_1'(\beta_j)= 0$, i.e. $w_{0j}^1\wedge \cdots \wedge w_{0j}^q \wedge d \beta_j\left(w_{0j}^\alpha\right) + \sum_{\gamma=1}^q w_{0j}^1 \wedge \cdots \wedge \beta_j\left(w_{0j}^\gamma\right)\wedge \cdots \wedge w_{0j}^q \wedge dw_{0j}^\alpha,\alpha=1,...,q$
\item $\beta_k- \beta_j + E_0'\left( (\theta_{jk}, B_{jk})\right)=0$, i.e. $\beta_k\left(w_{0j}^\alpha\right) - \beta_j\left(w_{0j}^\alpha \right) +\mathcal{L}_{\theta_{jk}}\left(w_{0j}^\alpha \right) - B_{jk}\left(w_{0j}^\alpha \right)=0$ for $\alpha=1,...,q$.
\item $\left( \theta_{ij}, B_{ij} \right)- \left( \theta_{ik}, B_{ik} \right) + \left( \theta_{jk}, B_{jk} \right)=0$, i.e. $\theta_{ij}-\theta_{ik}+ \theta_{jk}=0$ and $B_{ij}- B_{ik} + B_{jk}^i=0$
\end{enumerate}
we define $S_j=\left[\beta_j, \beta_j \right]_{\mathcal{N}_{\mathcal{F}_0}^*} \in \Gamma\left( U_j, \mathscr{H}om_{\mathcal{O}_M}\left(  \mathcal{N}_{\mathcal{F}_0}^*, \tilde{\mathcal{S}}^2\right) \right)$ by
{\small{\begin{align*}
S_j:&\Gamma\left( U_j, \mathcal{N}_{\mathcal{F}_0}^* \right) \to \Gamma\left( U_j, \tilde{\mathcal{S}}^2\right) \\
w_{0j}^\alpha &\mapsto 2\sum_{\gamma=1}^q w_{0j}^1\wedge \cdots \wedge \beta_j \left(w_{0j}^\gamma\right)\wedge \cdots \wedge w_{0j}^q \wedge d \beta_j\left(w_{0j}^\alpha\right)+2\sum_{1\leq \gamma <\eta\leq q} w_{0j}^1\wedge \cdots \wedge \beta_j\left(w_{0j}^\gamma\right)\wedge \cdots  \wedge \beta_j\left(w_{0j}^\eta\right)\wedge \cdots \wedge w_{0j}^q\wedge dw_{0j}^\alpha
\end{align*}}}
and linearly extends to $\Gamma\left( U_j, \mathcal{N}_{\mathcal{F}_0}^* \right)$, and we define $G_{jk}\in \Gamma\left( U_{jk}, \mathscr{H}om_{\mathcal{O}_M} \left( \mathcal{N}_{\mathcal{F}_0}^*, \Omega_M^1 \right) \right)$ by 
\begin{align*}
G_{jk}: \Gamma\left( U_{jk}, \mathcal{N}_{\mathcal{F}_0}^* \right) &\to \Gamma\left(  U_{jk}, \Omega_M^1  \right) \\
w_{0j}^\alpha &\mapsto - \beta_j\left( B_{jk}\left(w_{0j}^\alpha \right) \right) - \beta_k\left( B_{jk}\left( w_{0j}^\alpha \right) \right) + \mathcal{L}_{\theta_{jk}}\left(\beta\left(w_{0j}^\alpha \right) \right)  +  \mathcal{L}_{\theta_{jk}}\left( \beta_k\left(w_{0j}^\alpha \right)\right)
\end{align*}
and linearly extends to $\Gamma\left( U_{jk}, \mathcal{N}_{\mathcal{F}_0}^* \right)$. Then
\begin{align}\label{it2}
\left( \left\{ S_j \right\}, \left\{ G_{jk} \right\}, \left\{\left[\left( \theta_{ij}, B_{ij} \right) , \left(\theta_{jk}, B_{jk} \right) \right]  \right\} \right)\in C^0\left( \mathcal{U}, \mathscr{H}om_{\mathcal{O}_M}\left( \mathcal{N}_{\mathcal{F}_0}^*, \tilde{\mathcal{S}}^2 \right) \right)\bigoplus C^1\left( \mathcal{U} , \mathscr{H}om_{\mathcal{O}_M}\left( \mathcal{N}_{\mathcal{F}_0}^*, \Omega_M^1 \right) \right) \bigoplus C^2\left( \mathcal{U} , \mathcal{E}_{\mathcal{N}_{\mathcal{F}_0}^*}\right)
\end{align}
defines a $2$-cocycle in \v Cech resolution of $(\ref{oc5})$ of $\mathcal{E}_{\mathcal{N}_{\mathcal{F}_0}^*}^\bullet$.
\end{lemma}

\begin{remark}\label{oc80}
We show that $G_{jk}\in C^1\left( \mathcal{U}, \mathscr{H}om_{\mathcal{O}_M}\left( \mathcal{N}_{\mathcal{F}_0}^*, \Omega_M^1 \right) \right)$.  In fact, first we note that $(\theta_{jk}, B_{jk})$ on $U_j$ is transformed on $U_k$ into $\left( \theta_{jk}, B_{jk}^k \right)$ where 
\begin{align}\label{oc75}
B_{jk}^k\left(T_{0j}^\alpha \right)= B_{jk}\left(w_{0j}^\alpha \right) - \sum_{\beta=1}^q \left[ \theta_{jk}, h_{0jk}^{\alpha\beta} \right] w_{0k }^\beta = -B_{kj}\left(w_{0j}^\alpha \right) =- \sum_{\beta=1}^q h_{0jk}^{\alpha\beta} B_{kj}\left(w_{0k}^\beta \right) 
\end{align}
Then
{\Small{\begin{align*}
&G_{kj}\left(w_{0j}^\alpha \right) = \sum_{\beta=1}^q h_{0jk}^{\alpha\beta} G_{kj}\left(w_{0k}^\beta \right)\\
& = -  \sum_{\beta=1}^q h_{0jk}^{\alpha\beta}\beta_k\left(B_{kj}\left(w_{0k}^\beta \right) \right) - \sum_{\beta=1}^q h_{0jk}^{\alpha\beta} \beta_j \left(B_{kj}\left(w_{0k}^\beta \right) \right) +\sum_{\beta=1}^q h_{0jk}^{\alpha\beta} \mathcal{L}_{\theta_{kj}}\left(\beta_k\left(w_{0k}^\beta \right)\right) + \sum_{\beta=1}^q h_{0jk}^{\alpha\beta}\mathcal{L}_{\theta_{kj}}\left( \beta_j \left(w_{0k}^\beta \right)\right)\\ 
&= \beta_k\left( B_{jk}\left( w_{0j}^\alpha \right)        \right) - \sum_{\beta=1}^q \left[ \theta_{jk},  h_{0jk}^{\alpha\beta} \right] \beta_k\left( w_{0k}^\beta \right) + \beta_j\left( B_{jk}\left(w_{0j}^\alpha \right) \right) - \sum_{\beta=1}^q\left[ \theta_{jk}, h_{0jk}^{\alpha\beta} \right]\beta_k\left( w_{0j}^\beta \right)   - \mathcal{L}_{ \theta_{jk}}\left( \beta_k\left( w_{0j}^\alpha  \right) \right) + \sum_{\beta=1}^q \left[ \theta_{jk}, h_{0jk}^{\alpha\beta} \right] \beta_k \left(w_{0k}^\beta \right) \\
&- \mathcal{L}_{ \theta_{jk}}\left( \beta_j\left( w_{0j}^\alpha \right) \right) + \sum_{\beta=1}^q\left[ \theta_{jk}, h_{0jk}^{\alpha\beta} \right] \beta_j\left(w_{0k}^\beta \right) =\beta_k\left(B_{jk} \left( w_{0j}^\alpha \right) \right) + \beta_j\left( B_{jk}\left( w_{0j}^\alpha \right) \right) - \mathcal{L}_{  \theta_{jk}}\left( \beta_k \left( w_{0j}^\alpha \right)\right) - \mathcal{L}_{\theta_{jk}}\left( \beta_j \left( w_{0j}^\alpha \right)\right) = - G_{jk}\left( w_{0j}^\alpha\right)
\end{align*}}}

\end{remark}

\begin{remark}\label{oc81}

We show that $S_j:=\left[ \beta_j, \beta_j\right]_{\mathcal{N}_{\mathcal{F}_0}^*}$ is well-defined. In other words, $S_j\left(w_{0j}^\alpha \right)$ actually lies in $\Gamma\left( U_j, \tilde{\mathcal{S}}^2 \right)$. In fact, for $x\in  U_j-S$, there is a neighborhood $U_x\subset U_j-S$ such that $dw_{0j}^\alpha- \sum_{\delta=1}^q a_{j_x}^{\alpha \delta} \wedge w_{0j}^\delta=0$ on $U_x$ for some $a_{j_x}^{\alpha\delta}\in \Gamma\left( U_x, \Omega_M^1\right), \alpha=1,...,q$. Then we have 
\begin{align*}
w_{0j}^1\wedge \cdots \wedge w_{0j}^q \wedge \left(d \beta_j\left(w_{0j}^\alpha\right)- \sum_{\delta=1}^q a_{j_x}^{\alpha \delta}\wedge \beta_j \left(w_{0j}^\delta\right) \right)=0\,\,\,\,\,\,\,\,\,\text{on}\,\,\,\,\, U_{x} , 
\end{align*}
so that we can write
\begin{align} \label{oc37}
d\beta_j\left(w_{0j}^\alpha \right) = \sum_{\delta=1}^q a_{j_x}^{\alpha\delta}\wedge \beta_j\left(w_{0j}^\delta \right) + \sum_{\gamma=1}^q b_{j_x}^{\alpha \gamma} \wedge w_{0j}^\gamma \,\,\,\,\,\,\textnormal{on}\,\,\,U_x\,\,\,\,\textnormal{for some $b_{j_x}^{\alpha \gamma}\in \Gamma\left( U_x, \Omega_M^1 \right)$}
\end{align}
Then we have
{\Small{\begin{align*}
&2\sum_{\gamma=1}^q w_{0j}^1\wedge \cdots \wedge \beta_j \left(w_{0j}^\gamma \right)\wedge \cdots \wedge w_{0j}^q \wedge d \beta_j \left(w_{0j}^\alpha\right)+2\sum_{1\leq \gamma <\eta\leq q} w_{0j}^1\wedge \cdots \wedge \beta_j\left(w_{0j}^\gamma\right)\wedge \cdots  \wedge  \beta_j\left(w_{0j}^\eta\right)\wedge \cdots \wedge w_{0j}^q\wedge dw_{0j}^\alpha\\
&=2\sum_{\gamma=1}^q w_{0j}^1\wedge \cdots \wedge \beta_j\left(w_{0j}^\gamma\right)\wedge \cdots \wedge w_{0j}^q \wedge d \beta_j \left(w_{0j}^\alpha\right)+2\sum_{1\leq \gamma <\eta\leq q} w_{0j}^1\wedge \cdots \wedge \beta_j\left(w_{0j}^\gamma\right)\wedge \cdots  \wedge \beta_j\left(w_{0j}^\eta\right)\wedge \cdots \wedge w_{0j}^q\wedge \left(\sum_{\delta=1}^q a_{j_x}^{\alpha \delta}\wedge w_{0j}^\delta \right)\\
&=2\sum_{\gamma=1}^q w_{0j}^1\wedge \cdots \wedge \beta_j\left(w_{0j}^\gamma\right)\wedge \cdots \wedge w_{0j}^q \wedge d \beta_j\left(w_{0j}^\alpha\right)+2\sum_{1\leq \gamma <\eta\leq q} w_{0j}^1\wedge \cdots \wedge \beta_j\left(w_{0j}^\gamma\right)\wedge \cdots  \wedge \beta_j\left(w_{0j}^\eta\right)\wedge \cdots \wedge w_{0j}^q\wedge \left(a_{j_x}^{\alpha \gamma }\wedge w_{0j}^\gamma+ a_{j_x}^{\alpha \eta}\wedge w_{0j}^\eta \right)\\
&=2\sum_{\gamma=1}^q w_{0j}^1\wedge \cdots \wedge \beta_j\left(w_{0j}^\gamma\right)\wedge \cdots \wedge w_{0j}^q \wedge d \beta_j\left(w_{0j}^\alpha\right)-2\sum_{1\leq \gamma <\eta\leq q} w_{0j}^1\wedge  \cdots  \wedge \beta_j\left(w_{0j}^\eta\right)\wedge \cdots \wedge w_{0j}^q\wedge \left(a_{j_x}^{\alpha \gamma }\wedge \beta_j\left(w_{0j}^\gamma \right)\right)\\
&-2\sum_{1\leq \gamma <\eta\leq q} w_{0j}^1\wedge \cdots \wedge \beta_j\left(w_{0j}^\gamma\right)\wedge  \cdots \wedge w_{0j}^q\wedge \left(a_{j_x}^{\alpha \eta}\wedge \beta_j\left(w_{0j}^\eta\right) \right)\\
&=2\sum_{\gamma=1}^q w_{0j}^1\wedge \cdots \wedge \beta_j\left(w_{0j}^\gamma\right)\wedge \cdots \wedge w_{0j}^q \wedge d \beta_j\left(w_{0j}^\alpha\right)-2\sum_{1\leq \gamma \ne \eta\leq q} w_{0j}^1\wedge  \cdots  \wedge \beta_j\left(w_{0j}^\eta\right)\wedge \cdots \wedge w_{0j}^q\wedge \left(a_{j_x}^{\alpha \gamma}\wedge \beta_j\left(w_{0j}^\gamma\right)\right)\\
&=2\sum_{\gamma=1}^q w_{0j}^1\wedge \cdots \wedge \beta_j\left(w_{0j}^\gamma\right)\wedge \cdots \wedge w_{0j}^q \wedge d \beta_j\left(w_{0j}^\alpha\right)-2\sum_{1\leq  \eta\leq q} w_{0j}^1\wedge  \cdots  \wedge \beta_j\left(w_{0j}^\eta\right)\wedge \cdots \wedge w_j^q\wedge \left(\sum_{\gamma=1,\gamma\ne \eta}^q a_{j_x}^{\alpha \gamma}\wedge \beta_j\left(w_{0j}^\gamma\right)\right)\\
&=2\sum_{\gamma=1}^q w_{0j}^1\wedge \cdots \wedge \beta_j \left(w_{0j}^\gamma\right)\wedge \cdots \wedge w_{0j}^q \wedge d \beta_j\left(w_{0j}^\alpha\right)-2\sum_{1\leq  \eta\leq q} w_{0j}^1\wedge  \cdots  \wedge \beta_j\left(w_{0j}^\eta\right)\wedge \cdots \wedge w_{0j}^q\wedge \left(\sum_{\gamma=1}^q a_{j_x}^{\alpha \gamma}\wedge \beta_j\left(w_{0j}^\gamma \right) - a_{j_x}^{\alpha \eta} \wedge \beta_j\left(w_{0j}^\eta\right)\right)\\
&=2\sum_{1\leq  \eta\leq q} w_{0j}^1\wedge  \cdots  \wedge \beta_j\left(w_{0j}^\eta\right)\wedge \cdots \wedge w_{0j}^q\wedge \left(\sum_{\gamma=1}^q b_{j_x}^{\alpha \gamma}\wedge w_{0j}^\gamma  \right) = -2\sum_{\eta=1}^q w_{0j}^1 \wedge \cdots \wedge w_{0j}^q \wedge \left( b_{j_x}^{\alpha \eta} \wedge \beta_j\left(w_{0j}^\eta\right) \right)\\
&= -2 w_{0j}^1 \wedge \cdots \wedge w_{0j}^q\wedge \left(\sum_{\gamma=1}^q b_{j_x}^{\alpha \gamma} \wedge \beta_j \left(w_{0j}^\gamma\right)  \right) \in \Gamma\left(U_x, \tilde{\mathcal{S}}^2 \right)
\end{align*}}}

\end{remark}

\begin{remark}
We will deduce the equalities $(\ref{oc60}),(\ref{oc61}), (\ref{oc62})$ and $(\ref{oc67})$ below which will be used in the proof of  \textnormal{Lemma \ref{oc29}}. Let $B_{jk}^\alpha:=B_{jk}\left(w_{0j}^\alpha \right)= \sum_{\beta=1}^q B_{jk}^{\alpha\beta} w_{0j}^\beta$ with $B_{jk}^{\alpha\beta}\in \Gamma\left( U_{jk}, \mathcal{O}_M \right)$. Then we have
\begin{align*}
&w_{0j}^1 \wedge \cdots \wedge w_{0j}^q \wedge d \beta_k\left(B_{jk}^\alpha \right)=w_{0j}^1 \wedge \cdots \wedge w_{0j}^q \wedge \left( \sum_{\beta=1}^q d B_{jk}^{\alpha\beta} \wedge \beta_k\left(w_{0j}^\beta \right)  + \sum_{\beta=1}^q B_{jk}^{\alpha\beta} d\beta_k\left( w_{0j}^\beta \right)   \right)\\
&= \sum_{\beta=1}^q w_{0j}^1 \wedge \cdots \wedge w_{0j}^q \wedge d B_{jk}^{\alpha\beta} \wedge \beta_k \left(w_{0j}^\beta \right) - \sum_{\beta=1}^q B_{jk}^{\alpha\beta}\left( \sum_{\gamma=1}^q w_{0j}^1 \wedge \cdots \wedge \beta_k \left(w_{0j}^\gamma \right)\wedge \cdots \wedge w_{0j}^q \wedge dw_{0j}^\beta \right) \\
&=\sum_{\beta=1}^q w_{0j}^1 \wedge \cdots \wedge w_{0j}^q \wedge d B_{jk}^{\alpha\beta} \wedge \beta_k \left(w_{0j}^\beta \right) - \sum_{\gamma=1}^q w_{0j}^1 \wedge \cdots \wedge \beta_k\left( w_{0j}^\gamma \right) \wedge \cdots \wedge w_{0j}^q \wedge\left( d B_{jk}^\alpha - \sum_{\beta=1}^q d B_{jk}^{\alpha\beta}   \wedge w_{0j}^\beta  \right) \\
&=- \sum_{\gamma=1}^q w_{0j}^1 \wedge \cdots \wedge \beta_k\left( w_{0j}^\gamma \right) \wedge \cdots \wedge w_{0j}^q \wedge d B_{jk}^\alpha
\end{align*}

Hence we have
\begin{align}\label{oc60}
w_{0j}^1 \wedge \cdots \wedge w_{0j}^q \wedge d \beta_k\left(B_{jk}^\alpha \right) + \sum_{\gamma=1}^q w_{0j}^1 \wedge \cdots \wedge \beta_k\left( w_{0j}^\gamma \right) \wedge \cdots \wedge w_{0j}^q \wedge d B_{jk}^\alpha =0
\end{align}
Similarly we have
\begin{align}\label{oc68}
w_{0j}^1 \wedge \cdots \wedge w_{0j}^q \wedge d \mathcal{L}_{\theta_{jk}}\left(B_{jk}^\alpha \right) + \sum_{\gamma=1}^q w_{0j}^1 \wedge \cdots \wedge \mathcal{L}_{\theta_{jk}}\left( w_{0j}^\gamma \right) \wedge \cdots \wedge w_{0j}^q \wedge d B_{jk}^\alpha =0
\end{align}

We note that by taking $2\mathcal{L}_{\theta_{jk}}$ on $w_{0j}^1\wedge \cdots w_{0j}^q\wedge d \beta_k \left(w_{0j}^\alpha\right)+\sum_{\gamma=1}^q w_{0j}^1\wedge \cdots \wedge \beta_k\left(w_{0j}^\gamma\right)\wedge \cdots \wedge w_{0j}^q\wedge dw_{0j}^\alpha=0$, we have
{\Small{\begin{align}\label{oc61}
&2\sum_{\gamma=1}^q w_{0j}^1\wedge \cdots \wedge \mathcal{L}_{\theta_{jk}}\left(w_{0j}^\gamma\right)\wedge \cdots \wedge w_{0j}^q\wedge d\beta_k\left(w_{0j}^\alpha \right) + 2w_{0j}^1\wedge \cdots \wedge w_{0j}^q \wedge dL_{\theta_{jk}}\left(\beta_k\left(w_{0j}^\alpha\right) \right)\\
&+2\sum_{1\leq \eta<\gamma \leq q} w_{0j}^1\wedge \cdots \wedge \mathcal{L}_{\theta_{jk}}\left(w_{0j}^\eta\right) \wedge \cdots  \wedge \beta_k\left(w_{0j}^\gamma\right)\wedge \cdots \wedge w_{0j}^q\wedge dw_{0j}^\alpha \notag \\
&+2\sum_{1\leq \gamma <\eta \leq q} w_{0j}^1\wedge \cdots \wedge \beta_k\left(w_{0j}^\gamma\right)\wedge \cdots \wedge L_{\theta_{jk}}\left(w_{0j}^\eta\right)\wedge \cdots \wedge w_{0j}^q\wedge dw_{0j}^\alpha \notag \\
&+2\sum_{1\leq \gamma \leq q} w_{0j}^1\wedge \cdots \wedge \mathcal{L}_{\theta_{jk}}\left(\beta_k \left(w_{0j}^\gamma\right)\right)\wedge \cdots \wedge w_{0j}^q\wedge dw_{0j}^\alpha +2\sum_{1\leq \gamma \leq q} w_{0j}^1\wedge \cdots \wedge \beta_k\left(w_{0j}^\gamma\right)\wedge \cdots \wedge w_{0j}^q \wedge d \mathcal{L}_{\theta_{jk}}\left(w_{0j}^\alpha\right)=0 \notag
\end{align}}}

We note that by taking $\mathcal{L}_{\theta_{jk}}$ on $\sum_{\gamma=1}^q w_{0j}^1\wedge \cdots \wedge \mathcal{L}_{\theta_{jk}}\left(w_{0j}^\gamma\right)\wedge \cdots \wedge w_{0j}^q\wedge dw_{0j}^\alpha+ w_{0j}^1\wedge \cdots \wedge w_{0j}^q\wedge d\mathcal{L}_{\theta_{jk}}\left(w_{0j}^\alpha\right)=0$, we have
{\small{\begin{align} \label{oc62}
&2\sum_{1\leq \gamma <\eta\leq q} w_{0j}^1\wedge \cdots \wedge \mathcal{L}_{\theta_{jk}}\left(w_{0j}^\gamma\right)\wedge \cdots \wedge \mathcal{L}_{\theta_{jk}}\left(w_{0j}^\eta\right)\wedge \cdots \wedge w_{0j}^q \wedge dw_{0j}^\alpha+2\sum_{\gamma=1}^q w_{0j}^1\wedge \cdots \wedge \mathcal{L}_{\theta_{jk}}\left(w_{0j}^\gamma \right) \wedge \cdots \wedge w_{0j}^q\wedge d \mathcal{L}_{\theta_{jk}}\left(w_{0j}^\alpha\right)\\
&+ \sum_{\gamma=1}^q w_{0j}^1\wedge \cdots \wedge \mathcal{L}_{\theta_{jk}}\left(\mathcal{L}_{\theta_{jk}}\left(w_{0j}^\gamma\right)\right) \wedge \cdots \wedge w_{0j}^q \wedge dw_{0j}^\alpha +w_{0j}^1\wedge \cdots \wedge w_{0j}^q \wedge d\mathcal{L}_{\theta_{jk}}\left(\mathcal{L}_{\theta_{jk}}\left(w_{0j}^\alpha\right)\right)=0 \notag
\end{align}}}

We note that by setting $H_{0jk} = \det \left( h_{0jk}^{\alpha\beta} \right)$, we have
{\Small{\begin{align}\label{oc63}
&\sum_{\gamma=1}^q w_{0j}^1\wedge \cdots \wedge \beta_k\left(w_{0j}^\gamma\right)\wedge \cdots \wedge w_{0j}^q \wedge d \beta_k \left(w_{0j}^\alpha\right)\\
&=\sum_{\gamma=1}^q\left(\sum_{\beta_1=1}^q h_{0jk}^{1\beta_1} w_{0k}^{\beta_1}  \right)\wedge \cdots  \wedge \left(\sum_{\beta_\gamma=1}^q h_{0jk}^{\gamma \beta_\gamma} \beta_k \left(w_{0k}^{\beta_\gamma}\right)    \right) \wedge \cdots \wedge  \left( \sum_{\beta_q=1}^q h_{0jk}^{q \beta_q} w_{0k}^{\beta_q} \right) \wedge d\left(\sum_{\beta=1}^q h_{0jk}^{\alpha \beta} \beta_k \left(w_{0k}^{\beta}\right) \right) \notag \\
&=\sum_{\beta=1}^q h_{0jk}^{\alpha\beta}\left(\sum_{\gamma=1}^q H_{0jk} w_{0k}^1\wedge \cdots \wedge \beta_k \left(w_{0k}^\gamma\right)\wedge \cdots \wedge w_{0k}^q \wedge d \beta_k\left(w_{0k}^\beta\right) \right)-\sum_{\beta=1}^q\sum_{\gamma=1}^q H_{0jk} w_{0k}^1\wedge \cdots \wedge \beta_k\left(w_{0k}^\gamma\right)\wedge \cdots \wedge w_{0k}^q\wedge \beta_k\left(w_{0k}^\beta\right) \wedge dh_{0jk}^{\alpha\beta} \notag
\end{align}}}
On the other hand,
{\Small{\begin{align}\label{oc64}
&\sum_{\gamma\ne \eta} w_{0j}^1 \wedge \cdots \wedge \beta_k \left(w_{0j}^\gamma\right)\wedge \cdots \wedge \beta_k \left(w_{0j}^\eta\right)\wedge \cdots \wedge w_{0j}^q \wedge dw_{0j}^\alpha=\sum_{\gamma \ne \eta} H_{0jk} w_{0k}^1 \wedge \cdots \wedge \beta_k \left(w_{0k}^\gamma\right) \wedge \cdots  \wedge \beta_k\left(w_{0k}^\eta\right) \wedge \cdots \wedge w_{0k}^q \wedge d\left( \sum_{\beta=1}^q h_{0jk}^{\alpha\beta} w_{0k}^\beta \right)\\
&=H_{0jk}\sum_{\beta=1}^q h_{0jk}^{\alpha\beta}\left( \sum_{\gamma\ne \eta} w_{0k}^1 \wedge \cdots \wedge \beta_k\left(w_{0k}^\gamma\right)\wedge \cdots \wedge \beta_k\left(w_{0k}^\eta\right) \wedge \cdots \wedge w_{0k}^q \wedge dw_{0k}^\beta      \right) \notag \\
& - H_{0jk}\sum_{\beta=1}^q\sum_{\gamma\ne \eta} w_{0k}^1 \wedge \cdots \wedge \beta_k\left(w_{0k}^\gamma\right) \wedge \cdots \wedge \beta_k \left(w_{0k}^\eta\right) \wedge \cdots \wedge w_{0k}^q \wedge w_{0k}^\beta \wedge dh_{0jk}^{\alpha\beta} \notag \\
&=H_{0jk}\sum_{\beta=1}^q h_{0jk}^{\alpha\beta}\left( \sum_{\gamma\ne \eta} w_{0k}^1 \wedge \cdots \wedge \beta_k\left(w_{0k}^\gamma\right)\wedge \cdots \wedge \beta_k\left(w_{0k}^\eta\right) \wedge \cdots \wedge w_{0k}^q \wedge dw_{0k}^\beta      \right) \notag \\
& - H_{0jk}\sum_{\gamma\ne \beta}  w_{0k}^1 \wedge \cdots \wedge \beta_k\left(w_{0k}^\gamma\right) \wedge \cdots \wedge \beta_k\left(w_{0k}^\beta\right) \wedge \cdots \wedge w_{0k}^q \wedge w_{0k}^\beta \wedge dh_{0jk}^{\alpha\beta} \notag
\end{align}}}
Since we have 
{\Small{\begin{align}\label{oc65}
&\sum_{\beta\ne \gamma=1}^q H_{0jk} w_{0k}^1\wedge \cdots \wedge \beta_k\left(w_{0k}^\gamma\right)\wedge \cdots \wedge w_{0k}^q\wedge \beta_k\left(w_{0k}^\beta\right) \wedge dh_{0jk}^{\alpha\beta}\\
&=- \sum_{\beta\ne \gamma=1}^q H_{0jk} w_{0k}^1\wedge \cdots \wedge \beta_k\left(w_{0k}^\gamma\right)\wedge \cdots \wedge \beta_k\left(w_{0k}^\beta\right) \wedge \cdots \wedge w_{0k}^q\wedge w_{0k}^\beta \wedge d h_{0jk}^{\alpha\beta}  \notag
\end{align}}}
Then $(\ref{oc63}),(\ref{oc64})$ and $(\ref{oc65})$ implies that
{\small{\begin{align}\label{oc67}
&\sum_{\gamma=1}^q w_{0j}^1\wedge \cdots \wedge \beta_k\left(w_{0j}^\gamma\right)\wedge \cdots \wedge w_{0j}^q \wedge d \beta_k \left(w_{0j}^\alpha\right) + \sum_{\gamma\ne \eta} w_{0j}^1 \wedge \cdots \wedge \beta_k \left(w_{0j}^\gamma\right)\wedge \cdots \wedge \beta_k \left(w_{0j}^\eta\right)\wedge \cdots \wedge w_{0j}^q \wedge dw_{0j}^\alpha\\
&= H_{0jk}\sum_{\beta=1}^q h_{0jk}^{\alpha\beta}\left(\sum_{\gamma=1}^q w_{0k}^1\wedge \cdots \wedge \beta_k\left(w_{0k}^\gamma\right)\wedge \cdots \wedge w_{0k}^q \wedge d \beta_k \left(w_{0k}^\alpha\right) + \sum_{\gamma\ne \eta} w_{0k}^1 \wedge \cdots \wedge \beta_k \left(w_{0k}^\gamma\right)\wedge \cdots \wedge \beta_k \left(w_{0k}^\eta\right)\wedge \cdots \wedge w_{0k}^q \wedge dw_{0k}^\alpha \right) \notag
\end{align}}}

\end{remark}

\begin{remark}\label{oc43}
The following computation will be used in the proof of \textnormal{Lemma \ref{oc9}} and \textnormal{Lemma \ref{ou50}}. For $x\in U_j-S$, choose a neighborhood $U_x\subset U_j-S$ such that $dw_{0j}^\alpha- \sum_{\delta=1}^q a_{j_x}^{\alpha\delta} \wedge w_{0j}^\delta$ for some $a_{j_x}^{\alpha\delta}\in \Gamma\left( U_x, \Omega_M^1 \right),\alpha=1,...,q$. Then by taking $d$ on $dw_{0j}^\alpha- \sum_{\delta=1}^q a_{j_x}^{\alpha\delta} \wedge w_{0j}^\delta=0$, we have
{\small{\begin{align*}
\sum_{\delta=1}^q da_{j_x}^{\alpha \delta} \wedge w_{0j}^\delta - \sum_{\delta=1}^q a_{j_x}^{\alpha \delta} \wedge dw_{0j}^\delta =0 \iff \sum_{\delta=1}^q d a_{j_x}^{\alpha\delta} \wedge w_{0j}^\delta - \sum_{\delta, \gamma=1}^q a_{j_x}^{\alpha \delta}\wedge a_{j_x}^{\delta \gamma }\wedge w_{0j}^\gamma=0 \iff \sum_{\delta=1}^q \left( da_{j_x}^{\alpha \delta} - \sum_{\gamma=1}^q a_{j_x}^{\alpha \gamma} \wedge a_{j_x}^{\gamma \delta}     \right) \wedge w_{0j}^\delta =0
\end{align*}}}
Then we can write  
\begin{align}\label{oc40}
da_{j_x}^{\alpha \delta} - \sum_{\gamma=1}^q a_{j_x}^{\alpha \gamma} \wedge a_{j_x}^{\gamma \delta}= \sum_{\eta=1}^q g_{j_x}^{\alpha \delta \eta } \wedge w_{0j}^\eta\,\,\,\,\,\,\textnormal{on}\,\,\,U_x\,\,\,\textnormal{for some}\,\,\, g_{j_x}^{\alpha\delta \eta}\in \Gamma\left(U_x, \Omega_M^1 \right)
\end{align}
Then we have
\begin{align*}
\sum_{\delta,\eta=1}^q g_{j_x}^{\alpha \delta \eta } \wedge w_{0j}^\eta \wedge w_{0j}^\delta=0 \iff \sum_{  \eta < \delta } \left(g_{j_x}^{\alpha \delta \eta}- g_{j_x}^{\alpha \eta \delta} \right) w_{0j}^\eta \wedge w_{0j}^\delta =0
\end{align*}
By multiplying $w_{0j}^1\wedge \cdots \wedge \widehat{w}_{0j}^\eta \wedge \cdots \wedge \widehat{w}_{0j}^\delta \wedge \cdots \wedge w_{0j}^q$, we see that 
\begin{align}\label{oc41}
g_{j_x}^{\alpha \delta \eta}- g_{j_x}^{\alpha \eta\delta} \in \Gamma\left(U_x ,  \mathcal{N}_{\mathcal{F}_0}^*\right)\,\,\,\,\,\textnormal{for}\,\,\, \delta\ne \eta.
\end{align}

On the other hand, from $(\ref{oc37})$, we have $d\beta_j\left(w_{0j}^\alpha\right)=\sum_{\delta=1}^q a_{j_x}^{\alpha \delta} \wedge \beta_j \left(w_{0j}^\delta\right)+ \sum_{\gamma=1}^q b_{j_x}^{\alpha \gamma} \wedge w_{0j}^\gamma$ for some $b_{j_x}^{\alpha\gamma}\in \Gamma\left( U_x, \Omega_M^1 \right)$. Then by taking $d$, we have
\begin{align*}
\sum_{\delta=1}^q d a_{j_x}^{\alpha \delta} \wedge \beta_j \left(w_{0j}^\delta\right) -\sum_{\delta=1}^q a_{j_x}^{\alpha\delta} \wedge d \beta_j \left(w_{0j}^\delta\right) + \sum_{\gamma=1}^q db_{j_x}^{\alpha \gamma} \wedge w_{0j}^\gamma - \sum_{\gamma=1}^q b_{j_x}^{\alpha \gamma} \wedge dw_{0j}^\gamma = 0
\end{align*}
Then from $(\ref{oc37})$ and $(\ref{oc40})$, we obtain
\begin{align*}
&\sum_{\delta=1}^q \left( \sum_{\eta=1}^q g_{j_x}^{\alpha \delta \eta} \wedge w_j^\eta     \right) \wedge \beta_j \left(w_{0j}^\delta\right) - \sum_{\delta, \eta=1}^q a_{j_x}^{\alpha \delta} \wedge b_{j_x}^{\delta \eta} \wedge w_{0j}^\eta + \sum_{\eta=1}^q d b_{j_x}^{\alpha \eta} \wedge w_{0j}^\eta - \sum_{\gamma,\eta=1}^q b_{j_x}^{\alpha \gamma}\wedge a_{j_x}^{\gamma \eta}\wedge w_{0j}^\eta\\
&= \sum_{\eta=1}^q \left(  - \sum_{\delta=1}^q g_{j_x}^{\alpha \delta \eta} \wedge \beta_j\left(w_{0j}^\delta\right) - \sum_{\delta=1}^q a_{j_x}^{\alpha \delta} \wedge b_{j_x}^{\delta \eta} + db_{j_x}^{\alpha \eta} -\sum_{\gamma=1}^q b_{j_x}^{\alpha \gamma}\wedge a_{j_x}^{\gamma \eta}  \right) \wedge w_{0j}^\eta = 0
\end{align*}
This implies that
\begin{align}\label{oc42}
- \sum_{\delta=1}^q a_{j_x}^{\alpha \delta} \wedge b_{j_x}^{\delta \eta} + db_{j_x}^{\alpha \eta} -\sum_{\gamma=1}^q b_{j_x}^{\alpha \gamma}\wedge a_{j_x}^{\gamma \eta} \equiv_{\mathcal{N}_{\mathcal{F}_0}^*} - \sum_{\delta=1}^q g_{j_x}^{\alpha \delta \eta} \wedge \beta_j\left(w_{0j}^\delta\right)
\end{align}

\end{remark}

\begin{proof}[Proof of \textnormal{Lemma \ref{oc29}}]
First we note that $-\delta\left( \left\{\left[(\theta_{ij}, B_{ij}), \left(\theta_{jk}, B_{jk} \right) \right] \right\}\right)$ as in the same arguments of the proof of  $-\delta\left(\left\{ \left[\left(\theta_{ij}, A_{ij} \right), \left(\theta_{jk}, A_{jk} \right) \right] \right\} \right)=0$  in Lemma \ref{o54}.

We show that $E_2\left(S_j \right)=0$. In fact, for $x\in U_j-S$, choose a neighborhood $U_x\subset U_j-S$ such that $dw_{0j}^\alpha- \sum_{\delta=1}^q a_{j_x}^{\alpha\delta} \wedge w_{0j}^\delta$ for some $a_{j_x}^{\alpha\delta}\in \Gamma\left( U_x, \Omega_M^1 \right)$. We can check $E_2\left( S_j\right)=0$ locally on $U_x$. We keep the notations in Remark \ref{oc81}, so that $S_j\left(w_{0j}^\alpha  \right)= -2 w_{0j}^1\wedge \cdots \wedge w_{0j}^q \wedge \left(\sum_{\gamma=1}^q b_{j_x}^{\alpha \gamma} \wedge \beta_j\left(w_{0j}^\gamma \right) \right)\in \Gamma\left(U_x, \tilde{\mathcal{S}}^2 \right)$.

We also keep the notations in Remark \ref{oc43}. Then from $(\ref{oc37}),(\ref{oc42}), (\ref{oc41})$, we have
\begin{align}\label{its46}
&d\left(\sum_{\gamma=1}^q b_{j_x}^{\alpha \gamma}\wedge \beta \left(w_{0j}^\gamma\right)\right)- \sum_{\delta=1}^q a_{j_x}^{\alpha\delta} \wedge \left( \sum_{\gamma=1}^q b_{j_x}^{\delta \gamma}\wedge \beta \left(w_{0j}^\gamma \right)    \right)\\
&=\sum_{\gamma=1}^q db_{j_x}^{\alpha \gamma} \wedge \beta\left(w_{0j}^\gamma \right) - \sum_{\gamma=1}^q b_{j_x}^{\alpha \gamma} \wedge d\beta\left(w_{0j}^\gamma \right) - \sum_{\gamma, \delta=1}^q a_{j_x}^{\alpha \delta} \wedge b_{j_x}^{\delta \gamma} \wedge \beta\left(w_{0j}^\gamma\right) \notag \\
&= \sum_{\gamma=1}^q d b_{j_x}^{\alpha\gamma} \wedge \beta\left(w_{0j}^\gamma\right) - \sum_{\gamma=1}^q b_{j_x}^{\alpha \gamma}\wedge \left( \sum_{\delta=1}^q a_{j_x}^{\gamma \delta} \wedge  \beta_j \left(w_{0j}^\delta\right)  + \sum_{\eta=1}^q b_{j_x}^{\gamma \eta} \wedge w_{0j}^\eta  \right) - \sum_{\gamma, \delta=1}^q a_{j_x}^{\alpha \delta}\wedge b_{j_x}^{\delta \gamma} \wedge \beta\left(w_{0j}^\gamma\right) \notag \\
&\equiv_{\mathcal{N}_{\mathcal{F}_0}^*}  \sum_{\gamma=1}^q \left(db_{j_x}^{\alpha \gamma} - \sum_{\delta=1}^q b_{j_x}^{\alpha \delta} \wedge a_{j_x}^{\delta \gamma} - \sum_{\delta=1}^q a_{j_x}^{\alpha \delta} \wedge b_{j_x}^{\delta \gamma}      \right) \wedge \beta_j \left(w_{0j}^\gamma\right) \notag \\
& \equiv_{\mathcal{N}_{\mathcal{F}_0}^*} \sum_{\gamma=1}^q \left( \sum_{\delta=1}^q g_{j_x}^{\alpha \delta \gamma} \wedge \beta_j \left(w_{0j}^\delta\right)   \right) \wedge \beta_j \left(w_{0j}^\gamma \right) \equiv_{\mathcal{N}_{\mathcal{F}_0}^*} 0 \notag
\end{align}
This implies that $E_2(S_j)=0$ on $U_x$, and so $U_j$.

We show that $-\delta\left( \left\{ S_j \right\} \right) + E_1'\left(  \left\{ G_{jk} \right\} \right)=0$. From $(\ref{oc60}),(\ref{oc68}), (\ref{oc61}),(\ref{oc62})$ and $(\ref{oc67})$, we have

{\Small{\begin{align*}
&2\sum_{\gamma=1}^q w_{0j}^1\wedge \cdots \wedge \beta_j\left(w_{0j}^\gamma\right)\wedge \cdots \wedge w_{0j}^q \wedge d \beta_j \left(w_{0j}^\alpha\right)+2\sum_{1\leq \gamma < \eta\leq q} w_j^1\wedge \cdots \wedge \beta_j \left(w_{0j}^\gamma\right)\wedge \cdots  \wedge \beta_j\left(w_{0j}^\eta\right)\wedge \cdots \wedge w_{0j}^q\wedge dw_{0j}^\alpha\\
&=2\sum_{\gamma=1}^qw_{0j}^1\wedge \cdots \wedge\left(\beta_k\left(w_{0j}^\gamma\right)+L_{\theta_{jk}}\left(w_{0j}^\gamma\right)-B_{jk}^\gamma \right)\wedge \cdots \wedge w_{0j}^q\wedge d\beta_j\left( w_{0j}^\alpha \right)\\
&+2\sum_{1\leq \gamma < \eta\leq q} w_{0j}^1\wedge \cdots \wedge \left( \beta_k\left(w_{0j}^\gamma\right) + L_{\theta_{jk}}\left(w_{0j}^\gamma\right) - B_{jk}^\gamma \right)\wedge \cdots \wedge\left(\beta_k\left(w_{0j}^\eta\right) + \mathcal{L}_{\theta_{jk}}\left(w_{0j}^\eta \right) - B_{jk}^\eta    \right)\wedge \cdots \wedge w_j^q\wedge dw_j^\alpha\\
&=2\sum_{\gamma=1}^qw_{0j}^1\wedge \cdots \wedge\left(\beta_k\left(w_{0j}^\gamma\right)+L_{\theta_{jk}}\left(w_{0j}^\gamma\right) \right)\wedge \cdots \wedge w_{0j}^q\wedge d\beta_j\left( w_{0j}^\alpha \right)-2\sum_{\gamma=1}^qw_{0j}^1\wedge \cdots \wedge B_{jk}^\gamma \wedge \cdots \wedge w_{0j}^q\wedge d\beta_j\left( w_{0j}^\alpha \right)\\
&-2\sum_{1\leq \gamma < \eta\leq q} w_j^1\wedge \cdots \wedge B_{jk}^\gamma \wedge \cdots \wedge\beta_j\left( w_{0j}^\eta \right)\wedge \cdots \wedge w_j^q\wedge dw_j^\alpha - 2 \sum_{1\leq \eta < \gamma\leq q} w_j^1\wedge \cdots \wedge \beta_j\left( w_{0j}^\eta \right) \wedge \cdots \wedge B_{jk}^\gamma\wedge \cdots \wedge w_j^q\wedge dw_j^\alpha\\
&+2\sum_{1\leq \gamma < \eta\leq q} w_j^1\wedge \cdots \wedge \left( \beta_k\left(w_{0j}^\gamma\right) + L_{\theta_{jk}}(w_j^\gamma) \right)\wedge \cdots \wedge\left(\beta_k\left(w_{0j}^\eta\right) + \mathcal{L}_{\theta_{jk}}\left(w_{0j}^\eta \right)    \right)\wedge \cdots \wedge w_j^q\wedge dw_j^\alpha\\
&=2\sum_{\gamma=1}^qw_{0j}^1\wedge \cdots \wedge\left(\beta_k\left(w_{0j}^\gamma\right)+L_{\theta_{jk}}\left(w_{0j}^\gamma\right) \right)\wedge \cdots \wedge w_{0j}^q\wedge d\beta_j\left( w_{0j}^\alpha \right)-2\sum_{\gamma=1}^q B_{jk}^{\gamma\gamma}w_{0j}^1\wedge \cdots \wedge w_{0j}^\gamma \wedge \cdots \wedge w_{0j}^q\wedge d\beta_j\left( w_{0j}^\alpha \right)\\
&-2\sum_{1\leq \gamma < \eta\leq q} B_{jk}^{\gamma \gamma}w_j^1\wedge \cdots \wedge w_{0j}^\gamma  \wedge \cdots \wedge\beta_j\left( w_{0j}^\eta \right)\wedge \cdots \wedge w_j^q\wedge dw_j^\alpha - 2 \sum_{1\leq \eta < \gamma\leq q} B_{jk}^{\gamma\gamma}w_j^1\wedge \cdots \wedge \beta_j\left( w_{0j}^\eta \right) \wedge \cdots \wedge w_{0j}^\gamma \wedge \cdots \wedge w_j^q\wedge dw_j^\alpha\\
&-2\sum_{1\leq \gamma < \eta\leq q} w_j^1\wedge \cdots \wedge B_{jk}^{\gamma \eta} w_{0j}^\eta\wedge \cdots \wedge\beta_j\left( w_{0j}^\eta \right)\wedge \cdots \wedge w_j^q\wedge dw_j^\alpha - 2 \sum_{1\leq \eta < \gamma\leq q} w_j^1\wedge \cdots \wedge \beta_j\left( w_{0j}^\eta \right) \wedge \cdots \wedge B_{jk}^{\gamma\eta } w_{0j}^\eta \wedge \cdots \wedge w_j^q\wedge dw_j^\alpha\\
&+2\sum_{1\leq \gamma < \eta\leq q} w_j^1\wedge \cdots \wedge \left( \beta_k\left(w_{0j}^\gamma\right) + L_{\theta_{jk}}(w_j^\gamma) \right)\wedge \cdots \wedge\left(\beta_k\left(w_{0j}^\eta\right) + \mathcal{L}_{\theta_{jk}}\left(w_{0j}^\eta \right)    \right)\wedge \cdots \wedge w_j^q\wedge dw_j^\alpha
\end{align*}}}
We continue the computation. Then
{\Small{\begin{align*}
&= 2\sum_{\gamma=1}^qw_{0j}^1\wedge \cdots \wedge\left(\beta_k\left(w_{0j}^\gamma\right)+L_{\theta_{jk}}\left(w_{0j}^\gamma\right) \right)\wedge \cdots \wedge w_{0j}^q\wedge d\beta_j\left( w_{0j}^\alpha \right) + 2 \sum_{\gamma=1}^q B_{jk}^{\gamma \gamma} w_j^1 \wedge \cdots \wedge \beta_j\left(w_{0j}^\gamma \right) \wedge \cdots \wedge w_{0j}^q \wedge dw_{0j}^\alpha \\
&+2\sum_{1\leq \gamma < \eta\leq q} w_j^1\wedge \cdots \wedge B_{jk}^{\gamma \eta} \beta_j\left(w_{0j}^\eta\right)\wedge \cdots \wedge  w_{0j}^\eta \wedge \cdots \wedge w_j^q\wedge dw_j^\alpha + 2 \sum_{1\leq \eta < \gamma\leq q} w_j^1\wedge \cdots \wedge w_{0j}^\eta  \wedge \cdots \wedge B_{jk}^{\gamma\eta } \beta_j\left(w_{0j}^\eta\right) \wedge \cdots \wedge w_j^q\wedge dw_j^\alpha\\
&+2\sum_{1\leq \gamma < \eta\leq q} w_j^1\wedge \cdots \wedge \left( \beta_k\left(w_{0j}^\gamma\right) + L_{\theta_{jk}}(w_j^\gamma) \right)\wedge \cdots \wedge\left(\beta_k\left(w_{0j}^\eta\right) + \mathcal{L}_{\theta_{jk}}\left(w_{0j}^\eta \right)    \right)\wedge \cdots \wedge w_j^q\wedge dw_j^\alpha\\
&= 2\sum_{\gamma=1}^qw_{0j}^1\wedge \cdots \wedge\left(\beta_k\left(w_{0j}^\gamma\right)+L_{\theta_{jk}}\left(w_{0j}^\gamma\right) \right)\wedge \cdots \wedge w_{0j}^q\wedge d\left(\beta_k\left( w_{0j}^\alpha \right) +\mathcal{L}_{\theta_{jk}}\left( w_{0j}^\alpha \right) - B_{jk}^\alpha \right)+ 2\sum_{\gamma=1}^q w_{0j}^1 \wedge \cdots \wedge \beta_j \left( B_{jk}^\gamma \right) \wedge \cdots \wedge w_{0j}^q \wedge dw_{0j}^\alpha \\
&+2\sum_{1\leq \gamma < \eta\leq q} w_{0j}^1\wedge \cdots \wedge \left( \beta_k\left(w_{0j}^\gamma\right) + L_{\theta_{jk}}(w_{0j}^\gamma) \right)\wedge \cdots \wedge\left(\beta_k\left(w_{0j}^\eta\right) + \mathcal{L}_{\theta_{jk}}\left(w_{0j}^\eta \right)    \right)\wedge \cdots \wedge w_j^q\wedge dw_j^\alpha\\
&=2 H_{0jk} \sum_{\beta=1}^q h_{0jk}^{\alpha\beta}\left( \sum_{\gamma=1}^q w_{0k}^1 \wedge \cdots \wedge \beta_k\left(w_{0k}^\gamma \right)\wedge \cdots \wedge w_{0k}^q \wedge d \beta_k \left( w_{0k}^\beta \right) +2 \sum_{1\leq \gamma < \eta \leq q} w_{0k}^1 \wedge \cdots \wedge \beta_k\left(w_{0k}^\gamma \right) \wedge \cdots \wedge \beta_k\left( w_{0k}^\eta\right) \wedge \cdots \wedge w_{0k}^q \wedge dw_{0k}^\beta \right) \\
&+2 w_{0j}^1 \wedge \cdots \wedge w_{0j}^q \wedge d\beta_k\left(B_{jk}^\alpha \right) + 2 w_{0j}^1 \wedge \cdots  \wedge w_{0j}^q \wedge d \mathcal{L}_{\theta_{jk}}\left(B_{jk}^\alpha \right) +2 \sum_{\gamma=1}^q w_{0j}^1 \wedge \cdots \wedge \left( \beta_k \left(  B_{jk}^\gamma \right) + \mathcal{L}_{\theta_{jk}}\left( B_{jk}^\gamma \right)  \right) \wedge \cdots \wedge w_{0j}^q \wedge dw_{0j}^q\\
&-2 \sum_{\gamma=1}^q w_{0j}^1 \wedge \cdots \wedge \mathcal{L}_{\theta_{jk}}\left(\beta_k\left(w_{0j}^\gamma \right) \right) \wedge \cdots \wedge w_{0j}^q \wedge dw_{0j}^\alpha- 2w_{0j}^1 \wedge \cdots \wedge w_{0j}^q \wedge d \mathcal{L}_{\theta_{jk}}\left(\beta_k\left( w_{0j}^\alpha \right) \right)\\
&- \sum_{\gamma=1}^q w_{0j}^1 \wedge \cdots \wedge \mathcal{L}_{\theta_{jk}}\mathcal{L}_{\theta_{jk}}\left( w_{0j}^\gamma \right) \wedge \cdots \wedge w_{0j}^q \wedge dw_{0j}^\alpha - w_{0j}^1 \wedge \cdots \wedge w_{0j}^q \wedge d\mathcal{L}_{\theta_{jk}}\mathcal{L}_{\theta_{jk}}\left(w_{0j}^\alpha \right)
\end{align*}}}

We note that
{\small{\begin{align*}
&-\beta_j\left(B_{jk}^\alpha \right)-\beta_k\left(B_{jk}^\alpha \right) +\mathcal{L}_{\theta_{jk}}\left( \beta_j\left( w_{0j}^\alpha \right) \right) +\mathcal{L}_{\theta_{jk}}\left( \beta_k\left(w_{0j}^\alpha \right) \right)\\
&= -2 \beta_k\left( B_{jk}^\alpha \right) -\sum_{\beta=1}^q B_{jk}^{\alpha\beta} \left( \mathcal{L}_{\theta_{jk}}- B_{jk} \right)\left(w_{0j}^\beta \right) +2 \mathcal{L}_{\theta_{jk}}\left( \beta_k\left(w_{0j}^\alpha \right) \right) + \mathcal{L}_{\theta_{jk}}\left(\mathcal{L}_{\theta_{jk}} \left(w_{0j}^\alpha \right)\right)- \mathcal{L}_{\theta_{jk}}\left(B_{jk}^\alpha \right)\\
&= -2 \beta_k\left( B_{jk}^\alpha \right) -2  \mathcal{L}_{\theta_{jk}}\left(B_{jk}^\alpha \right) +2 \mathcal{L}_{\theta_{jk}}\left( \beta_k\left(w_{0j}^\alpha \right) \right) + \mathcal{L}_{\theta_{jk}}\left(\mathcal{L}_{\theta_{jk}} \left(w_{0j}^\alpha \right)\right)  + \sum_{\beta=1}^q\left[ \theta_{jk}, B_{jk}^{\alpha\beta} \right] w_{0j}^\beta + \sum_{\beta, \gamma=1}^q B_{jk}^{\alpha\beta} B_{jk}^{\beta \gamma} w_{0j}^\gamma
\end{align*}}}
Hence $-\delta\left( \left\{ S_j \right\} \right) + E_1'\left(  \left\{ G_{jk} \right\} \right)=0$.

We show that $-\delta\left(\left\{ G_{jk} \right\} \right) +E_0'\left( \left\{\left[ \left(\theta_{ij}, B_{ij} \right),\left(\theta_{jk}, B_{jk} \right)  \right]  \right\}    \right)=0$. First we note that
{\small{\begin{align}\label{oc69}
&\beta_j\left(w_{0i}^\alpha\right)- \beta_i\left(w_{0i}^\alpha\right)+ \mathcal{L}_{\theta_{ij}}\left(w_{0i}^\alpha\right)= B_{ij}^\alpha,\,\,\,\beta_k \left(w_{0i}^\alpha\right)-\beta_i\left(w_{0i}^\alpha\right)+\mathcal{L}_{\theta_{ik}}\left(w_{0i}^\alpha\right)=B_{ik}^\alpha, \,\,\,\beta_k\left(w_{0j}^\beta\right)- \beta_j\left(w_{0j}^\beta\right)+\mathcal{L}_{\theta_{jk}}\left(w_{0j}^\beta\right)=B_{jk}^\beta \\
&\Longrightarrow  \beta_k\left( w_{0i}^\alpha \right)  - \beta_j\left( w_{0i}^\alpha \right) + \mathcal{L}_{\theta_{jk}}\left( w_{0i}^\alpha \right)  - \sum_{\beta=1}^q \left[ \theta_{jk}, h_{0ij}^{\alpha\beta} \right] w_{0j}^\beta= \sum_{\beta=1}^q h_{0ij}^{\alpha\beta} B_{jk}^\beta \Longrightarrow  B_{ij}^\alpha - B_{ik}^\alpha + \sum_{\beta=1}^q h_{0ij}^{\alpha\beta} B_{jk}^\beta= - \sum_{\beta=1}^w \left[ \theta_{jk}, h_{0ij}^{\alpha\beta} \right] w_{0j}^\beta \notag
\end{align}}}
Then from $(\ref{oc69})$, we have
{\small{\begin{align*}
&G_{ij}\left( w_{0i}^\alpha\right)- G_{ik}\left( w_{0i}^\alpha \right) + \sum_{\beta=1}^q h_{0ij}^{\alpha\beta} G_{jk}\left(w_{0j}^\beta \right)\\
&=- \beta_i\left(B_{ij}^\alpha\right) - \beta_j\left(B_{ij}^\alpha\right)+ \mathcal{L}_{\theta_{ij}}\left(\beta_i\left(w_{0i}^\alpha\right)\right)+\mathcal{L}_{\theta_{ij}}\left(\beta_j\left(w_{0i}^\alpha\right)\right) + \beta_i\left(B_{ik}^\alpha\right) + \beta_k\left(B_{ik}^\alpha\right)- \mathcal{L}_{\theta_{ik}} \left(\beta_i\left(w_{0i}^\alpha\right)\right)-\mathcal{L}_{\theta_{ik}}\left(\beta_k\left(w_{0i}^\alpha\right)\right)  \\
&-\sum_{\beta=1}^q h_{0ij}^{\alpha\beta}\beta_j\left(B_{jk}^\beta\right) - \sum_{\beta=1}^q h_{0ij}^{\alpha\beta} \beta_k\left(B_{jk}^\beta\right) + \sum_{\beta=1}^q h_{0ij}^{\alpha\beta} \mathcal{L}_{\theta_{jk}}\left(\beta_j\left(w_{0j}^\beta\right)\right) + \sum_{\beta=1}^q h_{ij}^{\alpha\beta}  \mathcal{L}_{\theta_{jk}} \left(\beta_k\left(w_{0j}^\beta\right)\right)\\
&=- \beta_i\left(B_{ij}^\alpha\right) - \cancel{\beta_j\left(B_{ij}^\alpha\right)}+ \mathcal{L}_{\theta_{ij}} \left(\beta_i\left(w_{0i}^\alpha\right)\right)+\mathcal{L}_{\theta_{ij}}\left(\beta_j\left(w_{0i}^\alpha\right)\right) + \beta_i\left(B_{ik}^\alpha\right) + \bcancel{\beta_k\left(B_{ik}^\alpha\right)}- \mathcal{L}_{\theta_{ik}} \left(\beta_i\left(w_{0i}^\alpha\right)\right)-\mathcal{L}_{\theta_{ik}}\left(\beta_k\left(w_{0i}^\alpha\right) \right) \\
&+\cancel{\beta_j\left(B_{ij}^\alpha\right) }- \beta_j\left(B_{ik}^\alpha \right) +\cancel{\sum_{\beta=1}^q \left[\theta_{jk}, h_{0ij}^{\alpha\beta}\right] \beta_j\left(w_{0j}^\beta\right) } +\beta_k\left(B_{ij}^\alpha\right)- \bcancel{\beta_k\left(B_{ik}^\alpha\right) }+ \xcancel{  \sum_{\beta=1}^q \left[\theta_{jk}, h_{0ij}^{\alpha\beta}\right] \beta_k\left(w_{0j}^\beta\right)} \\
&+\mathcal{L}_{\theta_{jk}}\left(\beta_j\left(w_{0i}^\alpha\right) \right)- \cancel{\sum_{\beta=1}^q \left[ \theta_{jk}, h_{0ij}^{\alpha\beta}\right]\beta_j\left(w_{0j}^\beta\right)} +\mathcal{L}_{\theta_{jk}}\left(   \beta_k\left(w_{0i}^\alpha\right) \right) - \xcancel{\sum_{\beta=1}^q \left[ \theta_{jk}, h_{0ij}^{\alpha\beta} \right] \beta_k\left(w_{0j}^\beta\right) }\\
&=\left( \beta_k- \beta_i \right)\left( B_{ij}^\alpha  \right) -\left(\beta_j - \beta_i \right)\left( B_{ik}^\alpha \right) - \mathcal{L}_{\theta_{ij}}\left( \beta_k\left(w_{0i}^\alpha \right) - \beta_i\left( w_{0i}^\alpha  \right)\right) + \mathcal{L}_{\theta_{ik}}\left(\beta_j\left( w_{0i}^\alpha  \right) - \beta_i \left( w_{0i}^\alpha \right) \right)\\
&= \left( \beta_k- \beta_i \right)\left( B_{ij}^\alpha  \right) -\left(\beta_j - \beta_i \right)\left( B_{ik}^\alpha \right) - \mathcal{L}_{\theta_{ij}}\left( \sum_{\beta=1}^q B_{ik}^{\alpha\beta} w_{0i}^\beta - \mathcal{L}_{\theta_{ik}}\left(w_{0i}^\alpha \right)\right) + \mathcal{L}_{\theta_{ik}}\left( \sum_{\beta=1}^q B_{ij}^{\alpha\beta} w_{0i}^\beta - \mathcal{L}_{\theta_{ij}}\left(  w_{0i}^\alpha \right) \right) \\
&= \sum_{\beta=1}^q B_{ij}^{\alpha\beta}\left( \beta_k\left(w_{0i}^\beta \right)  - \beta_i\left( w_{0i}^\beta \right) + \mathcal{L}_{\theta_{ik}}\left( w_{0i}^\beta\right) - B_{ik}\left( w_{0i}^\beta \right) + B_{ik}\left( w_{0i}^\beta \right)    \right) + \sum_{\beta=1}^q \left[ \theta_{ik}, B_{ij}^{\alpha\beta} \right] w_{0i}^\beta - \mathcal{L}_{\theta_{ik}}\mathcal{L}_{\theta_{ij}}\left( w_{0i}^\alpha \right)\\
&- \sum_{\beta=1}^q B_{ik}^{\alpha\beta}\left( \beta_j\left(w_{0i}^\beta \right)  - \beta_i\left( w_{0i}^\beta \right) + \mathcal{L}_{\theta_{ij}}\left( w_{0i}^\beta\right) - B_{ij}\left( w_{0i}^\beta \right) + B_{ij}\left( w_{0i}^\beta \right)    \right) - \sum_{\beta=1}^q \left[ \theta_{ij}, B_{ik}^{\alpha\beta} \right] w_{0i}^\beta + \mathcal{L}_{\theta_{ij}}\mathcal{L}_{\theta_{ik}}\left( w_{0i}^\alpha \right)\\
&= B_{ik}\circ B_{ij}\left(  w_{0i}^\alpha \right) - B_{ij}\circ B_{ik}\left(w_{0i}^\alpha \right) + \sum_{\beta=1}^q \left[ \theta_{ij}+ \theta_{jk}, B_{ij}^{\alpha\beta} \right] w_{0i}^\beta - \sum_{\beta=1}^q \left[ \theta_{ij}, B_{ij}^{\alpha\beta} + B_{jk}^{i\alpha\beta} \right] T_{0i}^\beta + \mathcal{L}_{ \left[  \theta_{ij}, \theta_{jk} \right]}\left( w_{0i}^\alpha \right) \\
&= \left( B_{ij} + B_{jk}^i\right)\circ B_{ij}\left( w_{0i}^\alpha \right)- B_{ij}\circ \left( B_{ij} + B_{jk}^i \right) \left(w_{0i}^\alpha \right) + \sum_{\beta=1}^q \left[ \theta_{jk}, B_{ij}^{\alpha\beta} \right] w_{0i}^\beta - \sum_{\beta=1}^q \left[ \theta_{ij}, B_{jk}^{i \alpha\beta} \right] w_{0i}^\beta + \mathcal{L}_{\left[  \theta_{ij}, \theta_{jk} \right]}\left( w_{0i}^\alpha\right)\\
&= \mathcal{L}_{ \left[ \theta_{ij}, \theta_{jk} \right]} \left( w_{0i}^\alpha \right)  - \left(\left[ \theta_{ij}, B_{jk}^{i} \right] - \left[ \theta_{jk},  B_{ij} \right] + \left[ B_{ij},  B_{jk}^i \right] \right) \left(w_{0i}^\alpha \right) = E_0'\left( \left[ \left( \theta_{ij}, B_{ij} \right),  \left( \theta_{jk}, B_{jk} \right)    \right]      \right)
\end{align*}}}

This completes Lemma \ref{oc29}.

\end{proof}

Then from $(\ref{oc34})$ and Lemma \ref{oc29}, we obtain
\begin{theorem}
Suppose given a compact foliated complex manifold $\left( M, \mathcal{N}_{\mathcal{F}_0}^* \right)$ with $\mathcal{N}_{\mathcal{F}_0}^*$ locally free, and $\left( \left\{ \beta_j\right\}, \left\{\left( \theta_{jk}, B_{jk}\right) \right\} \right)\in \mathbb{H}^1\left( M, \mathcal{E}_{\mathcal{N}_{\mathcal{F}_0}^*}^\bullet \right)$. In order that there exist a foliated complex analytic family $\left( \mathcal{M}, \mathcal{N}_\mathcal{F}^*, \pi, B \right)$ in terms of cotangent sheaves such that $\pi^{-1}(0)=\left( M, \mathcal{N}_{\mathcal{F}_0}^* \right)$, and that $\left( \frac{\partial\left( \mathcal{M}, \mathcal{N}_\mathcal{F}^* \right)}{\partial t} \right)|_{t=0}= \left( \left\{ \beta_j \right\} , \left\{ \left( \theta_{jk}, B_{jk}\right) \right\}  \right)$, it is necessary that
{\Small{\begin{align*}
\left( \left\{ S_j \right\}, \left\{ G_{jk} \right\}, \left\{\left[\left( \theta_{ij}, B_{ij} \right) , \left(\theta_{jk}, B_{jk} \right) \right]  \right\} \right)\in C^0\left( \mathcal{U}, \mathscr{H}om_{\mathcal{O}_M}\left( \mathcal{N}_{\mathcal{F}_0}^*, \tilde{\mathcal{S}}^2 \right) \right)\bigoplus C^1\left( \mathcal{U} , \mathscr{H}om_{\mathcal{O}_M}\left( \mathcal{N}_{\mathcal{F}_0}^*, \Omega_M^1 \right) \right) \bigoplus C^2\left( \mathcal{U} , \mathcal{E}_{\mathcal{N}_{\mathcal{F}_0}^*}\right)
\end{align*}}}
from $(\ref{it2})$ in \textnormal{Lemma \ref{oc29}} defines the $0$ class in $\mathbb{H}^2\left( M, \mathcal{E}_{\mathcal{N}_{\mathcal{F}_0}^*}^\bullet \right)$.
\end{theorem}

Given $\left( \left\{ \beta_j\right\}, \left\{ \left( \theta_{jk}, B_{jk} \right) \right\} \right)\in \mathbb{H}^1\left( M, \mathcal{E}_{\mathcal{N}_{\mathcal{F}_0}^*}^\bullet \right)$, if the $2$-cocycle $\left( \left\{ S_j \right\}, \left\{ G_{jk} \right\}, \left\{\left[\left( \theta_{ij}, B_{ij} \right) , \left(\theta_{jk}, B_{jk} \right) \right]  \right\} \right)$ is not $0$ in $\mathbb{H}^2\left( M, \mathcal{E}_{\mathcal{N}_{\mathcal{F}_0}^*}^\bullet \right)$, there is no foliated deformation $\left( M_t, \mathcal{N}_{\mathcal{F}_t}^* \right)$ in terms of cotangent sheaves with $\left(M_0, \mathcal{N}_{\mathcal{F}_0}^* \right)=\left( M, \mathcal{N}_{\mathcal{F}_0}^* \right)$ and $\left( \frac{\partial\left( \mathcal{M}, \mathcal{N}_\mathcal{F}^* \right)}{\partial t} \right)|_{t=0}=\left(\left\{ \beta_j \right\}, \left\{ \left( \theta_{jk}, B_{jk}\right) \right\} \right)$. 

\section{Infinitesimal foliated deformations and obstructions in terms of cotangent sheave by using Dolbeault resolution of $\mathcal{E}_{\mathcal{N}_{\mathcal{F}_0}^*}^\bullet$}

Let $\left( \mathcal{M}, \mathcal{N}_{\mathcal{F}}^*, B, \pi \right)$ with $\mathcal{N}_{\mathcal{F}}^*$ locally free  be a foliated analytic family in terms of cotangent sheaves of deformations of a compact foliated complex manifold $\left( M, \mathcal{N}_{\mathcal{F}_0}^* \right)=\pi^{-1}(0)$ with $\mathcal{N}_{\mathcal{F}_0}^*$ locally free. Assume that given an open covering $\{U_j\}$ of $M$ by coordinate neighborhoods, $\Gamma\left( U_j, \mathcal{N}_{\mathcal{F}_0}^* \right)$ is generated by $w_{0j}^1,..., w_{0j}^q\in \Gamma\left( U_j, \Omega_M^1 \right)$ with $w_{0j}^\alpha= \sum_{\beta=1}^q h_{0jk}^{\alpha\beta} w_{0k}^\beta$ on $U_j\cap U_k$ for some $h_{0jk}^{\alpha\beta} \in \Gamma\left( U_{jk}, \mathcal{O}_M \right),\alpha=1,...,q$ and $w_{0j}^1 \wedge \cdots \wedge w_{0j}^q \wedge dw_{0j}^\alpha=0,\alpha=1,...,q$. In this subsection, we show that in the foliated analytic family $\left( \mathcal{M}, \mathcal{N}_\mathcal{F}^*, B, \pi \right)$, the deformations near $\left( M, \mathcal{N}_{\mathcal{F}_0}^* \right)$ are represented by a section $\left\{\left(\varphi(t), Z_j(t) \right) \right\} \in A^{0,1}\left( M, \mathcal{E}_{\mathcal{N}_{\mathcal{F}_0}^*} \right)$ where $\varphi(t)\in A^{0,1}\left( M, \Theta_M \right)$ and $Z_j(t)\in \Gamma\left( U_j, \mathcal{A}^{0,1}\left( \mathscr{H}om_{\mathcal{O}_M}\left( \mathcal{N}_{\mathcal{F}_0}^*, \mathcal{N}_{\mathcal{F}_0}^* \right) \right) \right)$ which is defined by $Z_j(t)\left( w_{0j}^\alpha \right) = \sum_{beta=1}^q Z_j^{\alpha\beta} w_{0j}^\beta$ for some $Z_j^{\alpha\beta}(t) \in \Gamma\left( U_j, \mathcal{A}^{0,1} \right),\alpha = 1,...,q$, and $Y(t)= \left\{Y_j(t) \right\}\in A^{0,0}\left( M, \mathscr{H}om_{\mathcal{O}_M}\left( \mathcal{N}_{\mathcal{F}_0}^*, \Omega_M^1 \right) \right)$ with $Y_j(t) \in \Gamma\left( U_j, \mathcal{A}^{0,0}\left(\mathscr{H}om_{\mathcal{O}_M} \left( \mathcal{N}_{\mathcal{F}_0}^*, \Omega_M^1 \right) \right) \right)$ satisfying $\varphi(0)=0$ and $Z_j(0)=0$, and $Y(0)=\textnormal{Id}$ and the integrability conditions : by setting $ Y_j(t)\left(w_{0j}^\alpha\right)=Y_j^\alpha(t)\in \Gamma\left(U_j, \mathcal{A}^{0,0}\left(\Omega_M^1 \right) \right)$ and $Z_j(t)\left( w_{0j}^\alpha \right)= Z_j^\alpha(t) \in \Gamma\left( U_j, \mathcal{A}^{0,0}\left( \Omega_M^1 \right) \right),\alpha=1,...,q$,
\begin{align*}
&\bar{\partial} \varphi(t) - \frac{1}{2} \left[ \varphi(t), \varphi(t) \right]= 0\\
&Y_j^1(t)\wedge \cdots \wedge Y_j^q(t) \wedge \partial Y_j^\alpha(t)=0,\,\,\,\alpha=1,...,q\\
&\bar{\partial} Y_j^\alpha(t) - \mathcal{L}_{\varphi(t)} \left( Y_j^\alpha(t) \right) + \sum_{\eta=1}^q Z_j^{\alpha \eta}(t) Y_j^\eta (t)=0,\,\,\,\alpha=1,...,q \\
&\bar{\partial} Z_j^{\alpha \delta}(t)- \left[ \varphi(t), Z_j^{\alpha \delta}(t)\right] + \sum_{\sigma=1}^q Z_j^{\alpha \sigma}(t) \wedge Z_j^{\sigma \delta}(t) = 0 , \,\,\,\,\,\alpha, \delta=1,...,q
\end{align*}
where $t\in \Delta$ a sufficiently small polydisk in $B$. To deduce the integrability condition, we extend Kodaira's approach (\cite{Kod05} section \S 5.3) in the context of a foliated analytic family in terms of cotangent sheaves.

\subsection{Integrability conditions}\label{oci80} \

Let $\left(\mathcal{M}, \mathcal{N}_\mathcal{F}^* , B, \pi \right)$ with $\mathcal{N}_\mathcal{F}^*$ locally free be a foliated analytic family in terms of cotangent sheaves, where $B$ is a domain of $\mathbb{C}^m$ containing the origin $0$ as in Definition \ref{n3}. Define $|t|=\max_\lambda |t_\lambda|$ for $t=(t_1,...,t_m)\in \mathbb{C}^m$, and let $\Delta=\{t\in \mathbb{C}^m | |t|<r\}$ the polydisk  of radius $r>0$. If we take a sufficiently  small $\Delta \subset B$, then $\left(\mathcal{M}_\Delta, \mathcal{N}_{\mathcal{F}_\Delta}^* \right)=\pi^{-1}\left(\Delta\right)$ is represented in the form
\begin{align*}
\left(\mathcal{M}_\Delta, \mathcal{N}_{\mathcal{F}_\Delta}^* \right)=\bigcup_j \left(U_j\times \Delta, \mathcal{N}_\mathcal{F}^*|_{U_j\times \Delta}\right)
\end{align*}
We denote a point of $U_j$ by $\xi_j=\left(\xi_j^1,...,\xi_j^n\right)$ and $\Gamma\left(U_j\times \Delta, \mathcal{N}_\mathcal{F}^* \right)$ is generated by $w_j^1\left(\xi_j,t\right),...,w_j^q\left(\xi_j,t\right)$, where
\begin{align*}
w_j^\alpha\left(\xi_j,t\right)=\sum_{\beta=1}^n w_j^{\alpha\beta}\left(\xi_j, t\right) d\xi_j^\beta
\end{align*}
for some $w_j^{\alpha\beta}(\xi_j,t)\in \Gamma \left( U_j \times \Delta , \mathcal{O}_\mathcal{M} \right)$, and satisfies
\begin{align}\label{ni6}
w_j^1(\xi_j,t)\wedge \cdots \wedge w_j^q(\xi_j,t)\wedge dw_j^\alpha(\xi_j,t)=0,\,\,\,\,\,\,\,\alpha=1,...,q,\,\,\,\,\,\textnormal{on}\,\,\,\,\, U_j \times \Delta.
\end{align}

For simplicity, we assume that $U_j=\left\{\xi_j\in \mathbb{C}^m | |\xi_j|<1\right\}$ where $|\xi_j|=\max_a |\xi_j^a|$. $(\xi_j,t)\in U_j\times \Delta$ and $(\xi_k, t)\in U_k\times \Delta$ are the same point on $\mathcal{M}_\Delta$ if $\xi_j^\alpha=f_{jk}^\alpha(\xi_k,t),\alpha=1,...,n$ where $f_{jk}(\xi_k,t)$ is a holomorphic map of $\xi_k^1,..., \xi_k^n, t_1,...,t_m$, defined on $(U_k\times \Delta) \cap (U_j\times \Delta)$, and we have the following relation
\begin{align*}
w_j^\alpha(\xi_j,t)=\sum_{\beta=1}^q h_{jk}^{\alpha\beta}(\xi_k,t) w_k^\beta(\xi_k,t),\,\,\,\,\,\alpha=1,...,q
\end{align*}
for a holomorphic function $h_{jk}^{\alpha\beta}(\xi_k,t)$ on $(U_k\times \Delta) \cap ( U_j\times \Delta)$. More precisely
\begin{align*}
\sum_{\gamma,\eta=1}^n w_j^{\alpha\gamma}(f_{jk}(\xi_k,t),t)\frac{\partial f_{jk}^\gamma(\xi_k,t)}{\partial \xi_k^\eta} d\xi_k^\eta=\sum_{\beta=1}^q\sum_{\eta=1}^n h_{jk}^{\alpha\beta}(\xi_k,t)w_k^{\beta\eta}(\xi_k,t)d\xi_k^\eta
\end{align*}
Equivalently
\begin{align}\label{oci1}
\sum_{\gamma=1}^n w_j^{\alpha\gamma}(f_{jk}(\xi_k,t),t)\frac{\partial f_{jk}^\gamma(\xi_k,t)}{\partial \xi_k^\eta}=\sum_{\beta=1}^q h_{jk}^{\alpha\beta}(\xi_k,t)w_k^{\beta\eta}(\xi_k,t)
\end{align}
and on $\left(U_i\times \Delta \right)\cap \left(U_j \times \Delta \right) \cap \left( U_k \times \Delta \right)$, we have
\begin{align}
h_{ik}^{\alpha \gamma}(\xi_k,t)= \sum_{\beta=1}^q h_{ij}^{\alpha\beta}\left(f_{jk}(\xi_k,t),t\right) h_{jk}^{\beta \gamma}(\xi_k,t),\,\,\,\,\,\alpha,\gamma=1,...,q
\end{align}

By \cite{Kod05} Theorem 2.3, $M_t$ is diffeomorphic to $M_0=\pi^{-1}(0)$ as differentiable manifolds for each $t\in \Delta$. We put $M:=M_0$. By \cite{Kod05} Theorem 2.5, if we take a sufficiently small $\Delta$, there is a diffeomorphism $\Psi$ of $M\times \Delta$ onto $\mathcal{M}_\Delta$ as differentiable manifolds such that $\pi\circ \Psi$ is the projection $M\times \Delta \to \Delta$. Let $z=\left(z^1,..., z^n \right)$ be local coordinates of $M=M_0$. Then we have $\pi\circ \Psi(z,t)=t, t\in \Delta$. For $\Psi(z,t)\in U_j\times \Delta$, put
\begin{align*}
\Psi\left( z, t\right) = \left(\xi_j^1(z,t), ..., \xi_j^n(z,t), t_1,..., t_m \right)
\end{align*}
Then each component $\xi_j^\alpha=\xi_j^\alpha(z,t),\alpha=1,...,n$ is a $C^\infty$ function. If we identify $\mathcal{M}_\Delta =\Psi\left(M\times \Delta \right)$ with $M\times \Delta$, and $\mathcal{N}_{\mathcal{F}}^*$ with $\Psi^* \mathcal{N}_\mathcal{F}^*$ via $\Psi$, $\left(\mathcal{M}_\Delta, \mathcal{N}_{\mathcal{F}_\Delta}^* \right)$ is considered as a complex manifold with the complex structure defined on the $C^\infty$ manifold $M\times \Delta$ by the system of local coordinates on $U_j\times \Delta$
\begin{align*}
\left\{ \left(\xi_j,t \right)| j=1 ,2, 3,...   \right\},\,\,\,\,\,(\xi_j,t)=\left( \xi_j^1(z,t), ..., \xi_j^n(z,t), t_1,..., t_m \right)
\end{align*}
and the locally free sheaf $\Psi^* \mathcal{N}_{\mathcal{F}_\Delta}^*$ on $M\times \Delta$, more precisely, the $C^\infty$ differentiable vector bundle on $M\times \Delta$.

We note that since $(z^1,..., z^n)$ and $\left(\xi_j^1(z,0), ..., \xi_j^n(z,0)\right)$ are local complex coordinates on $M=M_0$, $\xi_j^\alpha(z,0)$ are holomorphic functions of $z^1,...,z^n,\alpha=1,...,n$. We also note that if we take $\Delta$ sufficiently small, we have
\begin{align}
\det\left( \frac{\partial \xi_j^\alpha(z,t)}{\partial z^\lambda} \right)_{\alpha, \lambda=1,...,n} \ne 0\,\,\,\,\,\, \, \textnormal{for}\,\,\,t\in \Delta
\end{align}
On the other hand, each locally free sheaf $\Psi^* \mathcal{N}_{\mathcal{F}_t}^*$ in a differentiable  family $\left\{ \Psi^*\mathcal{N}_{\mathcal{F}_t}^* |t\in \Delta \right\}$ is the restriction to $M\times t$ of the locally free sheaf $\Psi^* \mathcal{N}_\mathcal{F}^*$ over $M\times \Delta , \Psi^* \mathcal{N}_{\mathcal{F}_t}^*= \Psi^* \mathcal{N}_{\mathcal{F}}^*|_{M\times t}$. Then by \cite{Kod05} Lemma 7.1 the locally free sheaf $\Psi^* \Theta_\mathcal{F}$ over $M\times \Delta$ as a $C^\infty$ vector bundle is $C^\infty$ equivalent to the trivial vector bundle $\mathcal{N}_{\mathcal{F}_0}^* \times \Delta $. In other words, there is a bundle diffeomorphism $\zeta: \mathcal{N}_{\mathcal{F}_0}^* \times \Delta \to \Psi^* \mathcal{N}_{\mathcal{F}}^*$ over $M\times \Delta$ when we consider $\mathcal{N}_{\mathcal{F}_0}^*$ and $\Psi^* \mathcal{N}_{\mathcal{F}}^*$ as vector bundles. We note that $\mathcal{N}_\mathcal{F}^*$ as a vector bundle is described by the transition matrices $\left( h_{jk}(\xi_j,t)^T\right)^{-1}$ with respect to the covering $\left\{U_j\times \Delta \right\}$ of $\mathcal{M}_\Delta$, where $h_{jk}(\xi_j,t):=\left( h_{jk}^{\alpha\beta} (\xi_k,t)\right)_{\alpha, \beta=1,...,q}$. Then $\Psi^* \mathcal{N}_{\mathcal{F}}^*$ is described by the transition matrices $\left( h_{jk}\left(\xi_k(z,t) , t\right)^T \right)^{-1}$ with respect to the covering $\left\{\mathscr{U}_j:=\Psi^{-1}\left( U_j\times \Delta \right) \right\}$ of $M\times \Delta$, where $h_{jk}\left(\xi_k\left( z, t\right), t \right):= \left( h_{jk}^{\alpha\beta}\left(\xi_k(z,t),t \right)\right)_{\alpha, \beta=1,...,q}$. In particular $\mathcal{N}_{\mathcal{F}_0}^*$ is described by $\left(h_{jk}(\xi_j(z,0),0 )^T \right)^{-1}$. Then there are invertible $C^\infty$ matrix functions $\zeta_j(z,t):\mathbb{C}^q \to \mathbb{C}^q$ such that $\sigma$ is represented locally by
\begin{align*}
\zeta :\mathbb{C}^q\times\mathscr{U}_j &\to \mathbb{C}^q \times \mathscr{U}_j\\
                         \left( X=(x_1,..., x_p), z, t \right) &\mapsto \left( \zeta_j (z,t)X, z, t\right)
\end{align*}
such that
\begin{align*}
\zeta_j(z,t) \cdot \left( h_{jk}\left(\xi_k(z,0),0 \right)^T \right)^{-1} = \left( h_{jk}\left(\xi_k(z,t), t \right)^T \right)^{-1} \cdot \zeta_k (z,t)
\end{align*}

Then we have
\begin{align*}
\left(\zeta_j(z,t)^{-1} \right)^T\cdot h_{jk}\left(\xi_k(z,0),0 \right) = h_{jk}\left(\xi_k( z,t),t \right) \cdot \left( \zeta_k (z,t)^{-1} \right)^T
\end{align*}
Let $\left(\zeta_j (z,t)^{-1} \right)^T $ be represented by invertible matrix $g_j(z,t):=\left(g_j^{\alpha\beta}(z,t) \right)$ where $g_j^{\alpha\beta}(z,t)$ is a $C^\infty$ function, and set $h_{0jk}^{\alpha\beta}(z):=h_{jk}^{\alpha\beta}(\xi_j(z,0),0)$ and $h_{0jk}(z):=\left(h_{0jk}^{\alpha\beta} (z)\right)_{\alpha, \beta=1,...,q}$, so that we have
\begin{align}\label{b19}
g_j(z,t)\cdot h_{0jk}(z)= h_{jk}\left(\xi_k(z,t),t) \right) \cdot g_k(z,t) 
\end{align}
Then we have
\begin{align}\label{oci11}
h_{jk}^{\alpha\beta}\left(\xi_k(z,t),t\right)= \sum_{\delta, \eta=1}^q g_{j}^{\alpha \delta}(z,t) h_{0jk}^{\delta \eta}(z) \left(g_k^{-1}\right)^{\eta\beta}(z,t),\,\,\,\,\,\,\,\alpha,\beta=1,...,q
\end{align}
where $g_j^{-1}(z,t)=\left( \left(  g_j^{-1}\right)^{\alpha\beta}(z,t) \right)_{\alpha, \beta=1,.., q} $ is the inverse of $g_j(z,t)$. 

\subsection{Identification of deformations of complex structures with $\varphi(t)\in A^{0,1}\left(M, \Theta_M \right)$}\

We keep the notations in subsection \ref{oi51}, so that we have $\varphi(t)\in A^{0,1}\left( M, \Theta_M \right)$ such that
\begin{align}\label{oci9}
\varphi(0)=0,\,\,\,\,\,\,\,\bar{\partial}\varphi(t)- \frac{1}{2}\left[ \varphi(t), \varphi(t) \right]=0
\end{align}

\subsection{Identification of deformations of complex structures and locally free sheaves with $\left\{ \left( \varphi(t), Z_j(t)\right) \right\} \in A^{0,1}\left( M, \mathcal{E}_{\mathcal{N}_{\mathcal{F}_0}^*} \right)$} \

As in the same way in subsection \ref{oi52}, by taking $Z_j^{\alpha \delta}(t):= \sum_{\eta=1}^q \left( g_j^{-1}\right)^{\alpha \eta}\left( \bar{\partial} g_j^{\eta \delta} - \left[ \varphi, g_j^{\eta \delta}\right] \right)$, let us define an element $Z_j (t)\in \Gamma\left(U_j, \mathcal{A}^{0,1}\left(\mathscr{H}om_{\mathcal{O}_M}\left( \mathcal{N}_{\mathcal{F}_0}^*, \mathcal{N}_{\mathcal{F}_0}^* \right) \right)    \right)$ by
\begin{align*}
Z_j(t):\Gamma\left( U_j, \mathcal{N}_{\mathcal{F}_0}^* \right)&\to \Gamma\left( U_j, \mathcal{A}^{0,0}\left( \mathcal{N}_{\mathcal{F}_0}^* \right) \right) \\
  w_{0j}^\alpha &\mapsto \sum_{\delta=1}^q Z_j^{\alpha \delta} w_{0j}^\delta:=\sum_{\eta, \delta=1}^q \left(g_j^{-1} \right)^{\alpha \eta} \left(\bar{\partial} g_j^{\eta \delta}- \left[\varphi, g_j^{\eta \delta} \right] \right) w_{0j}^\delta
\end{align*}
and linearly extends to $\Gamma\left(U_j, \mathcal{N}_{\mathcal{F}_0}^* \right)$. Then we have
\begin{align}\label{ni42}
Z_j\left( w_{0j}^\alpha \right) = Z_k\left( w_{0j}^\alpha \right) + \sum_{\beta=1}^q \left[ \varphi, h_{0jk}^{\alpha \beta} \right] w_{0j}^\beta,
\end{align}
 so that $\left\{\left(\varphi(t), Z_j(t) \right)    \right\}\in A^{0,1}\left( M, \mathcal{E}_{\mathcal{N}_{\mathcal{F}_0}^*} \right)$ such that $Z_j(0)=0$ and
\begin{align} \label{oci10}
\bar{\partial} Z_j^{\alpha \eta}- \left[\varphi, Z_j^{\alpha \eta} \right] + \sum_{\sigma=1}^q Z_j^{\alpha \sigma} \wedge Z_j^{\sigma \eta} = 0
\end{align}

\subsection{Identification of deformations of foliated structures with $Y(t)\in A^{0,0}\left( M, \mathscr{H}om_{\mathcal{O}_M}\left( \mathcal{N}_{\mathcal{F}_0}^*, \Omega_M^1 \right) \right)$ }\label{oci81}

Let us define 
\begin{align} \label{b20}
W_j^{\alpha\gamma}(z,t):= \sum_{\beta=1}^n w_j^{\alpha\beta}\left( \xi_j(z,t), t\right) \frac{\partial \xi_j^\beta(z,t)}{\partial z^\gamma},\alpha=1,...,q
\end{align}
and set $W_j^\alpha(z,t)= \sum_{\gamma=1}^n W_j^{\alpha \gamma}(z,t) dz^\gamma$. Then we have
\begin{align*}
&\sum_{\gamma=1}^n W_j^{\alpha\gamma}(z,t) dz^\gamma - \sum_{\beta=1}^q h_{jk}^{\alpha\beta}\left(\xi_k(z,t),t \right) \left(\sum_{\gamma=1}^n W_k^{\beta \gamma}(z,t)dz^\gamma   \right)\\
&= \sum_{\gamma=1}^n \left( W_j^{\alpha \gamma}(z,t)- \sum_{\beta=1}^q h_{jk}^{\alpha\beta}\left(\xi_k(z,t),t \right) W_k^{\beta \gamma}(z,t)   \right)dz^\gamma\\
&=\sum_{\gamma=1}^n \left( \sum_{\beta=1}^n w_j^{\alpha\beta}\left(\xi_j(z,t),t \right)\frac{\partial \xi_j^\beta}{\partial z^\gamma} - \sum_{\beta=1}^q \sum_{\eta=1}^n h_{jk}^{\alpha \beta}\left(\xi_k(z,t),t \right) w_k^{\beta \eta} \left(\xi_k(z,t),t \right) \frac{\partial \xi_k^\eta}{\partial z^\gamma}       \right)dz^\gamma\\
&=\sum_{\gamma=1}^n \left( \sum_{\beta, \eta=1}^n w_j^{\alpha\beta}\left(\xi_j(z,t),t \right)\frac{\partial \xi_j^\beta}{\partial  \xi_k^\eta } \frac{\partial \xi_k^\eta}{\partial z^\gamma} - \sum_{\beta=1}^q \sum_{\eta=1}^n h_{jk}^{\alpha \beta} \left(\xi_k(z,t),t \right) w_k^{\beta \eta} \left(\xi_k(z,t),t \right) \frac{\partial \xi_k^\eta}{\partial z^\gamma}       \right)dz^\gamma\\
&=\sum_{\gamma,\eta=1}^n \left(  \sum_{\beta=1}^n w_j^{\alpha\beta}\left(\xi_j(z,t),t\right)\frac{\partial \xi_j^\beta}{\partial \xi_k^\eta}  - \sum_{\beta=1}^q h_{jk}^{\alpha\beta}\left(\xi_k(z,t),t \right) w_k^{\beta \eta}\left(\xi_k(z,t),t\right)  \right)\frac{\partial \xi_k^\eta}{\partial z^\gamma} dz^\gamma
\end{align*}
Then from $(\ref{oci1})$, we obtain
\begin{align}\label{oci3}
W_j^\alpha(z,t)= \sum_{\beta=1}^q h_{jk}^{\alpha\beta} \left(\xi_k(z,t),t \right) W_k^{\beta}(z,t)
\end{align}

Let us consider
\begin{align*}
w_j^\alpha\left(\xi_j,t\right)=  W_j^{\alpha} (z,t)+ W_j'^{\alpha}(z,t)= \sum_{\beta=1}^n w_j^{\alpha\beta}\left(\xi_j(z,t),t \right)\frac{\partial \xi_j^\beta}{\partial z^\gamma}dz^\gamma + \sum_{\beta=1}^n w_j^{\alpha\beta}\left(\xi_j(z,t),t\right) \frac{\partial \xi_j^\beta}{\partial \bar{z}^\gamma} d\bar{z}^\gamma
\end{align*}
and we set $\partial W_j^\alpha(z,t):=\sum_{\gamma, \delta=1}^n \frac{\partial W_j^{\alpha \gamma}(z,t)}{\partial z^\delta} dz^\delta \wedge dz^\gamma $. Then from $(\ref{ni6})$ we have
\begin{align*}
0=&w_j^1 \left(\xi_j,t\right) \wedge \cdots \wedge w_j^q\left(\xi_j,t\right) \wedge dw_j^\alpha(\xi_j,t)\\
=& \left(W_j^1(z,t)+ W_j'^1(z,t)\right) \wedge \cdots \wedge \left(W_j^q(z,t)+ W_j'^q(z,t)\right)\wedge d\left(W_j^\alpha (z,t) + W_j'^\alpha(z,t) \right)\\
=& W_j^1(z,t) \wedge \cdots \wedge W_j^q (z,t) \wedge \partial W_j^\alpha(z,t) + \cdots
\end{align*}
Here $W_j^1\wedge \cdots \wedge W_j^q \wedge \partial W_j^\alpha$ only involves $dz^\xi$ and no $d\bar{z}^\xi$, and other terms contains at least one $d\bar{z}^\xi$.  Hence we have
\begin{align}\label{oci4}
 W_j^1 \left(z,t \right) \wedge \cdots \wedge W_j^q \left(z,t \right) \wedge \partial W_j^\alpha \left(z,t \right)=0
\end{align}

On the other hand, from $(\ref{oci3})$ and $(\ref{oci11})$, we have
\begin{align}\label{oci2}
W_j^\alpha= \sum_{\beta, \delta, \eta=1}^q g_j^{\alpha\delta} h_{0jk}^{\delta \eta} \left(g_k^{-1}\right)^{\eta \beta} W_k^\beta \iff \sum_{\delta=1}^q \left(g_j^{-1}\right)^{\alpha \delta} W_j^\delta = \sum_{\beta=1}^q h_{0jk}^{\alpha\beta} \left(\sum_{\delta=1}^q \left(g_k^{-1}\right)^{\beta \delta} W_k^\delta    \right)
\end{align}

Let us define
\begin{align*}
Y_j(t):\Gamma\left(U_j, \mathcal{N}_{\mathcal{F}_0}^*\right)&\to \Gamma \left(U_j, \mathcal{A}^{0,0}\left(\Omega_{M}^1\right)\right) \\
              w_{0j}^\alpha &\mapsto Y_j^\alpha :=\sum_{\delta=1}^q \left(g_j^{-1} \right)^{\alpha\delta} W_j^{\delta}
\end{align*}
and linearly extends to $\Gamma\left(U_j, \mathcal{N}_{\mathcal{F}_0}^* \right)$. Then $(\ref{oci2})$ implies that $Y= \left\{Y_j \right\}$ defines a global section in $A^{0,0}\left( M, \mathscr{H}om_{\mathcal{O}_M}\left(\mathcal{N}_{\mathcal{F}_0}^* , \Omega_M^1 \right)   \right)$. We note that as in Part I, since $\bar{\partial} W_j^\delta - \left[ \varphi, W_j^\delta \right] =0$, we have from $(\ref{oi19})$
\begin{align*}
\bar{\partial}\left( \sum_{\delta=1}^q \left(g_j^{-1}\right)^{\alpha \delta} W_j^\delta     \right)  - \mathcal{L}_{\varphi}\left(\sum_{\delta=1}^q \left(g_j^{-1} \right)^{\alpha \delta} W_j^\delta  \right) & =-\sum_{\beta , \eta , \delta=1}^q \left(g_j^{-1}\right)^{\alpha\beta}\bar{\partial } g_j^{\beta \eta} \cdot \left(g_j^{-1}\right)^{\eta \delta} \cdot W_j^\delta + \sum_{\delta, \eta , \delta=1}^q \left(g_j^{-1}\right)^{\alpha\beta}\left[\varphi, g_j^{\beta \eta}  \right] \left(g_j^{-1}\right)^{\eta \delta}  W_j^\delta \\
&=  - \sum_{\beta, \eta, \delta=1}^q \left(g_j^{-1}\right)^{\alpha\beta}\left(\bar{\partial} g_j^{\beta \eta}- \left[ \varphi, g_j^{\beta \eta} \right] \right) \left(g_j^{-1}\right)^{\eta \delta} W_j^\delta
\end{align*}
This implies that
\begin{align}\label{oci5}
\bar{\partial} \left(Y_j\left(w_{0j}^\alpha\right)\right) -  \mathcal{L}_{\varphi}\left( Y_j\left(w_{0j}^\alpha\right) \right) + Y_j \left(Z_j\left( w_{0j}^\alpha \right) \right)=0 \iff \bar{\partial} Y_j^\alpha - \mathcal{L}_{ \varphi}\left( Y_j^\alpha \right) + \sum_{\beta=1}^q Z_j^{\alpha \beta} Y_j^\beta = 0
\end{align}

We note that since $\partial Y_j^\alpha= \sum_{\beta=1}^q \partial \left(g_j^{-1}\right)^{\alpha\beta} \wedge W_j^\beta + \sum_{\beta=1}^q \left(g_j^{-1} \right)^{\alpha\beta} \partial W_j^\beta$, $(\ref{oci4})$ implies that
\begin{align}\label{oci8}
Y_j^1(z,t)\wedge \cdots \wedge Y_j^q(z,t) \wedge  \partial Y_j^\alpha(z,t)= 0,\,\,\,\,\,\alpha=1,...,q
\end{align}

\subsection{Identification of Infinitesimal deformations by using Dolbeault resolution of $\mathcal{E}_{\mathcal{N}_{\mathcal{F}_0}^*}^\bullet$} \label{oci83}\

In the previous subsection we showed from $(\ref{oci9}),(\ref{oci10}),(\ref{oci8})$ and $(\ref{oci5})$ that given a foliated analytic family $\left(\mathcal{M}, \mathcal{N}_\mathcal{F}^*, B, \pi \right)$ in terms of cotangent sheaves of deformations of $\left(M, \mathcal{N}_{\mathcal{F}_0}^* \right)=\pi^{-1}(0)$, the deformations near $\left( M, \mathcal{N}_{\mathcal{F}_0}^* \right)$ are represented by a section $\left\{\left(\varphi(t), Z_j(t) \right)\right\}\in A^{0,1}\left( M,\mathcal{E}_{\mathcal{N}_{\mathcal{F}_0}^*} \right)$, where $\varphi(t)\in A^{0,1}\left(M, \Theta_M \right)$ and  $Z_j(t)\in \Gamma\left( U_j, \mathcal{A}^{0,1}\left( \mathscr{H}om_{\mathcal{O}_M}\left( \mathcal{N}_{\mathcal{F}_0}^* , \mathcal{N}_{\mathcal{F}_0}^* \right)\right) \right)$ which is  defined by $Z_j(t)\left(w_{0j}^\alpha \right)= \sum_{\beta=1}^q Z_j^{\alpha\beta}(t) w_{0j}^\beta$ for some $U_j^{\alpha\beta}(t)\in \Gamma\left(U_j, \mathcal{A}^{0,1} \right),\alpha=1...,q$, and $Y(t)\in  A^{0,0}\left( M, \mathscr{H}om_{\mathcal{O}_M}\left( \mathcal{N}_{\mathcal{F}_0}^* , \Omega_M^1  \right) \right)$ satisfying $\varphi(0)=0, Z_j(0)=0$, and $Y(0)=\textnormal{Id}$ and the integrability condition: by setting $Y_j(t)\left(w_{0j}^\alpha \right)= Y_j^\alpha (t) \in \Gamma\left( U_j, \mathcal{A}^{0,0}\left(\Omega_M^1 \right) \right),\alpha=1,...,q$, 
\begin{align}
&\bar{\partial} \varphi(t) - \frac{1}{2}\left[\varphi(t), \varphi (t)\right]=0 \label{oci6}  \\
&  Y_j^1 \left(t \right) \wedge \cdots \wedge Y_j^q \left(t \right) \wedge \partial Y_j^\alpha \left(t \right)=0 , \,\,\,\,\, \alpha=1,...,q \\
&\bar{\partial} Y_j^\alpha (t) - \mathcal{L}_{ \varphi(t)}\left( Y_j^\alpha(t) \right) + \sum_{\beta=1}^q Z_j^{\alpha \beta}(t) Y_j^\beta(t) = 0,\,\,\,\,\,\alpha=1,...,q \\
&\bar{\partial} Z_j^{\alpha \delta}(t)- \left[ \varphi(t), Z_j^{\alpha \delta}(t) \right]  + \sum_{\sigma=1}^p Z_j^{\alpha \sigma} (t) \wedge Z_j^{\sigma \delta}(t) =0,\,\,\,\,\,\alpha, \delta=1,...,q \label{oci7}
\end{align}
where $t\in \Delta$ a sufficiently small polydisk in $B$.  We identify infinitesimal deformations at $0$ by taking the derivative of $(\ref{oci6})-(\ref{oci7})$ with respect to $t$ and setting $t=0$. Let us take the derivative of $(\ref{oci6})-(\ref{oci7})$ with respect to $t$. Then we have
\begin{align}
&\bar{\partial}\left(\frac{\partial \varphi}{\partial t} \right) - \frac{1}{2}\left[ \frac{\partial \varphi}{\partial t}, \varphi\right] - \frac{1}{2}\left[\varphi, \frac{\partial \varphi}{\partial t} \right] =0 \label{oci13} \\
&\sum_{\beta=1}^q Y_j^1\wedge \cdots \wedge \overbrace{ \frac{\partial Y_j^\beta}{\partial t}}^{\beta-\textnormal{th}} \wedge \cdots \wedge Y_j^q\wedge \partial Y_j^\alpha+ Y_j^1\wedge \cdots Y_j^q \wedge \partial \left(\frac{\partial Y_j^\alpha}{\partial t}   \right) =0 \\
&\bar{\partial}\left(\frac{\partial Y_j^\alpha}{\partial t} \right) - \mathcal{L}_{\frac{\partial\varphi}{\partial t}}\left( Y_j^\alpha \right) - \mathcal{L}_{ \varphi}\left(\frac{\partial Y_j^\alpha}{\partial t} \right) + \sum_{\beta=1}^q \frac{\partial Z_j^{\alpha \eta}}{\partial t} Y_j^\eta + \sum_{\beta=1}^q Z_j^{\alpha \eta} \frac{\partial Y_j^\eta}{\partial t} = 0\\
&\bar{\partial}\left( \frac{\partial Z_j^{\alpha \delta}}{\partial t} \right) - \left[ \frac{\partial \varphi}{\partial t}, Z_j^{\alpha \delta}\right] - \left[ \varphi, \frac{\partial Z_j^{\alpha \delta}}{\partial t} \right] + \sum_{\sigma=1}^p \frac{\partial Z_j^{\alpha\sigma}}{\partial t} \wedge Z_j^{\sigma \delta} + \sum_{\sigma=1}^q Z_j^{\alpha \sigma} \wedge \frac{\partial Z_j^{\sigma \delta}}{\partial t}=0 \label{oci14}
\end{align}
By setting $t=0$ of $(\ref{oci13})-(\ref{oci14})$, we obtain
\begin{align}
&\bar{\partial}\left(\frac{\partial \varphi}{\partial t}|_{t=0} \right)=0 \label{it21}\\
&\sum_{\beta=1}^q w_{0j}^1\wedge \cdots \wedge \overbrace{\frac{\partial Y_j^\beta}{\partial t}|_{t=0} }^{\beta- \textnormal{th}} \wedge \cdots \wedge w_{0j}^q \wedge dw_{0j}^q + w_{0j}^1\wedge \cdots \wedge w_{0j}^q \wedge \partial \left( \frac{\partial Y_j^\alpha}{\partial t}|_{t=0} \right) = 0\\
&\bar{\partial}\left(\frac{\partial Y_j^\alpha}{\partial t}|_{t=0} \right) - \mathcal{L}_{\frac{\partial\varphi}{\partial t}|_{t=0}}\left( w_{0j}^\alpha \right)  + \sum_{\beta=1}^q \frac{\partial Z_j^{\alpha \eta}}{\partial t}|_{t=0} w_{0j}^\eta = 0 \label{itt7}\\
&\bar{\partial}\left( \frac{\partial Z_j^{\alpha \delta}}{\partial t}|_{t=0} \right) = 0 \label{it22}
\end{align}
We define $\frac{\partial Y_j(t)}{\partial t}|_{t=0}\in A^{0,0}\left( M, \mathscr{H}om_{\mathcal{O}_M}\left(\mathcal{N}_{\mathcal{F}_0}^*, \Omega_M^1 \right)   \right)$ by
\begin{align*}
\frac{\partial Y_j}{\partial t}|_{t=0}:\Gamma\left( U_j, \mathcal{N}_{\mathcal{F}_0}^* \right) &\to \Gamma\left( U_j, \mathcal{A}^{0,0}\left(\Omega_M^1 \right) \right) \\
 w_{0j}^\alpha &\mapsto \frac{\partial Y_j^\alpha}{\partial t}|_{t=0}
\end{align*}
and linearly extends to $\Gamma\left(U_j, \mathcal{N}_{\mathcal{F}_0}^* \right)$. Then $\frac{\partial Y}{\partial t}|_{t=0} : =\left\{ \frac{\partial Y_j}{\partial t}|_{t=0} \right\}$ defines a global section in $A^{0,0}\left( M, \mathscr{H}om_{\mathcal{O}_M}\left( \mathcal{N}_{\mathcal{F}_0}^* , \Omega_M^1 \right)\right)$. We define an element $\frac{\partial Z_j}{\partial t}|_{t=0}\in \Gamma\left(U_j ,\mathcal{A}^{0,1}\left(\mathscr{H}om_{\mathcal{O}_M}\left(\mathcal{N}_{\mathcal{F}_0}^* , \mathcal{N}_{\mathcal{F}_0}^* \right) \right)\right)$ by
\begin{align*}
\frac{\partial Z_j}{\partial t}|_{t=0}:\Gamma\left( U_j, \mathcal{N}_{\mathcal{F}_0}^* \right) &\to \Gamma\left( U_j, \mathcal{A}^{0,1}\left(\mathcal{N}_{\mathcal{F}_0}^* \right) \right)\\
                                      w_{0j}^\alpha &\mapsto \sum_{\eta=1}^q \frac{\partial Z_j^{\alpha\eta}}{\partial t}|_{t=0} w_{0j}^\eta
\end{align*}
and linearly extends to $\Gamma\left( U_j, \mathcal{N}_{\mathcal{F}_0}^* \right)$. We note that from $(\ref{ni42})$, we have
\begin{align}\label{ni41}
\sum_{\eta=1}^q Z_j^{\alpha \eta} w_{0j}^\eta= \sum_{\beta, \eta=1}^q h_{0jk}^{\alpha\beta} Z_j^{\beta \eta} w_{0k}^\eta + \sum_{\gamma=1}^q \left[ \varphi, h_{0jk}^{\alpha\beta} \right] w_{0k}^\beta
\end{align}
By taking the derivative of $(\ref{ni41})$ with respect to $t$ and setting $t=0$, we have
\begin{align*}
\sum_{\eta=1}^q \frac{\partial Z_j^{\alpha \eta}}{\partial t}|_{t=0} w_{0j}^\eta = \sum_{\beta, \eta=1}^q h_{0jk}^{\alpha\beta}\frac{\partial Z_j^{\beta \eta}}{\partial t}|_{t=0} w_{0k}^\eta  + \sum_{\beta=1}^q \left[ \frac{\partial \varphi}{\partial t}|_{t=0}, h_{0jk}^{\alpha\beta} \right] w_{0k}^\beta
\end{align*}
This implies that $\left\{ \left( \frac{\partial \varphi}{\partial t}|_{t=0}, \frac{\partial Z_j}{\partial t}|_{t=0} \right) \right\}\in A^{0,1}\left(M, \mathcal{E}_{\mathcal{N}_{\mathcal{F}_0}^* }\right)$. 

Then $(\ref{oi42})-(\ref{oi46})$ implies that $\left(\frac{\partial Y}{\partial t}|_{t=0}, \left\{ \left(\frac{\partial \varphi}{\partial t}|_{t=0}, \frac{\partial Z_j}{\partial t}|_{t=0} \right)\right\} \right) \in A^{0,0}\left( M, \mathscr{H}om_{\mathcal{O}_M}\left(\Theta_{\mathcal{F}_0}, \Omega_M^1 \right)\right)\bigoplus A^{0,1}\left(M, \mathcal{E}_{\mathcal{N}_{\mathcal{F}_0}^*} \right) $ defines a $1$-cocycle in the following double complex:
{\Small{\begin{equation} \label{oci15}
\begin{CD}
\cdots @. \cdots \\
@AE_2|_{M-S} AA @AE_2|_{M-S}AA \\
A^{0,0}\left( M-S, \mathscr{H}om_{\mathcal{O}_M}\left( \mathcal{N}_{\mathcal{F}_0}^* \tilde{\mathcal{S}}^3 \right)\right) @>-(-1)^q\bar{\partial}>> A^{0,1}\left( M-S, \tilde{\mathcal{S}}^3\right) @> (-1)^q \bar{\partial} >> \cdots \\
@AE_2|_{M-S} AA @AE|_2|_{M-S} AA \\
A^{0,0}\left(M-S, \mathscr{H}om_{\mathcal{O}_M}\left( \mathcal{N}_{\mathcal{F}_0}^*,  \tilde{\mathcal{S}}^2 \right) \right) @> (-1)^q\bar{\partial}>> A^{0,1}\left(M-S, \tilde{\mathcal{S}}^2 \right) @>-(-1)^q\bar{\partial}>> \cdots  \\
@AE_1'AA @AD_1 AA \\
A^{0,0}\left(\mathcal{U}_t, \mathscr{H}om_{\mathcal{O}_{M}}\left(\Theta_{\mathcal{F}_0}, \Theta_{M} \right) \right)@>-\bar{\partial}>>  A^{0,1}\left(M, \mathscr{H}om_{\mathcal{O}_{M}}\left(\mathcal{N}_{\mathcal{F}_0}^*, \Omega_{M}^1  \right) \right) @>\bar{\partial} >> \cdots \\
@AE_0' AA @A{E}_0'AA @AE_0'AA \\
A^{0,0}\left(M,\mathcal{E}_{\mathcal{N}_{\mathcal{F}_0}}^*\right) @>\bar{\partial}>> A^{0,1}\left( M, \mathcal{E}_{\mathcal{N}_{\mathcal{F}_0}^*}\right) @>-\bar{\partial}>> A^{0,2} \left( M , \mathcal{E}_{\mathcal{N}_{\mathcal{F}_0}^*} \right) @> \bar{\partial}>> \cdots 
\end{CD}
\end{equation}}}
We will denote the $i$-th cohomology group of the above double complex $(\ref{oci15})$ by $\textnormal{B}^i $. In the next section \ref{ni2}, we will show that we have 
\begin{align}\label{ii1}
\textnormal{B}^i  \cong \mathbb{H}^i\left( M, \mathcal{E}_{\mathcal{N}_{\mathcal{F}_0}^*}^\bullet \right)\,\,\,\,\,\textnormal{for}\,\,\,i=0,1,2
\end{align}

Then we have the foliated Kodaira-Spencer map in terms of cotangent sheaves
\begin{align}\label{ii2}
\varphi_0: T_0 B &\to \mathbb{H}^1\left( M, \mathcal{E}_{\mathcal{N}_{\mathcal{F}_0}^*}^\bullet \right)\\
             \frac{\partial}{\partial t} &\mapsto \left(\frac{\partial Y}{\partial t}|_{t=0}, \left\{\left(\frac{\partial \varphi}{\partial t}|_{t=0} , \frac{\partial Z_j}{\partial t}|_{t=0}\right)\right\}   \right) \in A^{0,0}\left(M, \mathscr{H}om_{\mathcal{O}_M}\left( \Theta_{\mathcal{F}_0}, \Theta_M \right) \right) \bigoplus A^{0,1}\left( M, \mathcal{E}_{\Theta_{\mathcal{F}_0}} \right) \notag
\end{align}
In the next section \ref{ni1}, we will describe the explicit correspondence of infinitesimal deformations between \v Cech resolution and Dolbeault resolution.

\subsection{Identification of obstructions by using Dolbeault resolution of $\mathcal{E}_{\mathcal{N}_{\mathcal{F}_0}^*}^\bullet$}\label{oci84}\

We assume that the foliated analytic family $\left( \mathcal{M}, \mathcal{N}_\mathcal{F}^*, B, \pi \right)$ in the previous subsection is a $1$-parameter family, i.e. $\Delta \subset \mathbb{C}$. Then by taking the derivative of $(\ref{oci13})-(\ref{oci14})$ with respect to $t$ again, we have
\begin{align}
&\bar{\partial}\left(\frac{\partial^2 \varphi}{\partial t^2} \right) -\left[\frac{\partial \varphi}{\partial t}, \frac{\partial \varphi}{\partial t} \right] - \frac{1}{2}\left[ \frac{\partial^2 \varphi}{\partial t^2}, \varphi \right] - \frac{1}{2} \left[ \varphi, \frac{\partial^2 \varphi}{\partial t^2}\right]=0 \label{it11}\\
& \sum_{\beta=1}^q Y_j^1\wedge \cdots \wedge \frac{\partial^2 Y_j^\beta}{\partial t^2} \wedge \cdots \wedge Y_j^q \wedge \partial Y_j^\alpha +2 \sum_{1\leq \beta <\gamma \leq q} Y_j^1\wedge \cdots \wedge \frac{\partial Y_j^\beta}{\partial t} \wedge \cdots \wedge \frac{\partial Y_j^\gamma}{\partial t} \wedge \cdots Y_j^q \wedge \partial Y_j^\alpha \\
&+Y_j^1 \wedge \cdots \wedge Y_j^q \wedge \partial\left( \frac{\partial^2 Y_j^\alpha}{\partial t^2} \right) +2 \sum_{\beta=1}^q Y_j^1 \wedge \cdots \wedge \frac{\partial Y_j^\beta}{\partial t} \wedge \cdots \wedge Y_j^q \wedge \partial \left(\frac{\partial Y_j^\alpha}{\partial t} \right)=0   \notag \\
&\bar{\partial}\left(\frac{\partial^2 Y_j^{\alpha}}{\partial t^2} \right) - \mathcal{L}_{\frac{\partial^2 \varphi}{\partial t^2}}\left( Y_j^\alpha \right) -2 \mathcal{L}_{\frac{\partial \varphi}{\partial t}}\left(\frac{\partial Y_j^\alpha}{\partial t} \right) - \mathcal{L}_{\varphi}\left( \frac{\partial^2 Y_j^\alpha}{\partial t^2} \right) + \sum_{\beta=1}^q \frac{\partial^2 Z_j^{\alpha \eta}}{\partial t^2} Y_j^\eta + 2 \sum_{\beta=1}^q \frac{\partial Z_j^{\alpha \eta}}{\partial t}\frac{\partial Y_j^\eta}{\partial t} + \sum_{\beta=1}^q Z_j^{\alpha \eta} \frac{\partial^2 Y_j^\eta}{\partial t^2}=0 \\
&\bar{\partial}\left( \frac{\partial^2 Z_j^{\alpha \delta}}{\partial t^2} \right) - 2\left[ \frac{\partial \varphi}{\partial t}, \frac{\partial Z_j^{\alpha \delta}}{\partial t} \right] -\left[ \varphi, \frac{\partial^2 Z_j^{\alpha \delta}}{\partial t^2} \right] + \sum_{\sigma=1}^q \frac{\partial^2 Z_j^{\alpha \sigma}}{\partial t^2} \wedge Z_j^{\alpha \delta} + 2\sum_{\sigma=1}^q \frac{\partial Z_j^{\alpha \delta}}{\partial t} \wedge \frac{\partial Z_j^{\sigma \delta}}{\partial t} + \sum_{\sigma=1}^q Z_j^{\alpha \sigma} \wedge \frac{\partial^2 Z_j^{\alpha \delta}}{\partial t^2} = 0 \label{it12}
\end{align}

By setting $t=0$ of $(\ref{it11})-(\ref{it12})$, we obtain
{\small{\begin{align}
& -\bar{\partial}\left( \frac{\partial^2 \varphi}{\partial t^2}|_{t=0} \right)= -\left[ \frac{\partial \varphi}{\partial t}|_{t=0}, \frac{\partial \varphi}{\partial t}|_{t=0} \right]\label{it31} \\
& \sum_{\beta=1}^q w_{0j}^1 \wedge \cdots \wedge \frac{\partial^2 Y_j^\beta}{\partial t^2}|_{t=0}\wedge \cdots \wedge w_{0j}^q \wedge dw_{0j}^\alpha +  w_{0j}^1 \wedge \cdots \wedge w_{0j}^q \wedge \partial \left( \frac{\partial^2 Y_j^\alpha}{\partial t^2}|_{t=0}  \right) \\
&= -2 \sum_{1\leq \beta <\gamma \leq q} w_{0j}^1 \wedge \cdots \wedge \frac{\partial Y_j^\beta}{\partial t}|_{t=0} \wedge \cdots \wedge \frac{\partial Y_j^\gamma}{\partial t}|_{t=0} \wedge \cdots \wedge w_{0j}^q \wedge d w_{0j}^\alpha - 2 \sum_{\beta=1}^q w_{0j}^1 \wedge \cdots \wedge \frac{\partial Y_j^\beta}{\partial t}|_{t=0} \wedge \cdots \wedge w_{0j}^q \wedge \partial\left( \frac{\partial Y_j^\alpha}{\partial t}|_{t=0} \right) \notag \\
&-\bar{\partial}\left(\frac{\partial^2 Y_j^{\alpha}}{\partial t^2}|_{t=0} \right) + \mathcal{L}_{\frac{\partial^2 \varphi}{\partial t^2}|_{t=0} }\left( w_{0j}^\alpha \right) - \sum_{\beta=1}^q \frac{\partial^2 Z_j^{\alpha \eta}}{\partial t^2}|_{t=0} w_{0j}^\eta  = -2 \mathcal{L}_{\frac{\partial \varphi}{\partial t}|_{t=0}}\left(\frac{\partial Y_j^\alpha}{\partial t} |_{t=0} \right) + 2 \sum_{\beta=1}^q \frac{\partial Z_j^{\alpha \eta}}{\partial t}|_{t=0}\frac{\partial Y_j^\eta}{\partial t}|_{t=0} \label{itt13} \\
&- \bar{\partial}\left( \frac{\partial^2 Z_j^{\alpha \delta}}{\partial t^2} |_{t=0} \right) = - 2 \left[ \frac{\partial \varphi}{\partial t}|_{t=0},  \frac{\partial Z_j^{\alpha\delta}}{\partial t}|_{t=0}  \right] + 2\sum_{\sigma=1}^q \frac{\partial Z_j^{\alpha \sigma}}{\partial t}|_{t=0} \wedge \frac{\partial Z_j^{\sigma \delta}}{\partial t}|_{t=0} \label{it32}
\end{align}}}

\begin{lemma}\label{oci85}
For any $1$-cocycle $\left( \Psi=\left\{ \Psi_j \right\}, \left\{  \left( \psi, k_j \right) \right\}\right)\in A^{0,0}\left(M, \mathscr{H}om_{\mathcal{O}_M}\left( \mathcal{N}_{\mathcal{F}_0}^*, \Omega_M^1 \right) \right) \bigoplus A^{0,1}\left( M, \mathcal{E}_{\mathcal{N}_{\mathcal{F}_0}^*} \right)$ where $\Psi_j \in \Gamma\left( U_j, \mathcal{A}^{0,0}\left( \mathscr{H}om_{\mathcal{O}_M}\left( \mathcal{N}_{\mathcal{F}_0}^*, \Omega_M^1 \right) \right) \right)$ and $\psi\in A^{0,1}\left( M, \Theta_M \right)$ and $k_j \in \Gamma\left( U_j, \mathscr{H}om_{\mathcal{O}_M}\left( \mathcal{N}_{\mathcal{F}_0}^*, \mathcal{N}_{\mathcal{F}_0}^* \right) \right)$
such that
\begin{align*}
\Psi_j:\Gamma\left( U_j, \mathcal{N}_{\mathcal{F}_0}^* \right) &\to \Gamma\left( U_j, \mathcal{A}^{0,0}\left( \Omega_M^1 \right) \right)\\
 w_{0j}^\alpha &\mapsto \Psi_j^\alpha
\end{align*}
with $\Psi_j=\Psi_k$ on $U_{jk}$ and
\begin{align*}
k_j:\Gamma\left( U_j, \mathcal{N}_{\mathcal{F}_0}^* \right) &\to \Gamma\left( U_j, \mathcal{A}^{0,1}\left( \mathcal{N}_{\mathcal{F}_0}^* \right)    \right) \\
 w_{0j}^\alpha &\mapsto \sum_{\beta=1}^q k_j^{\alpha\beta} w_{0j}^\beta \,\,\,\,\,\,\,\,\,k_j^{\alpha \beta}\in \Gamma\left( U_j, \mathcal{A}^{0,1} \right)
\end{align*}
with $k_j\left(T_{0j}^\alpha \right)= k_k\left(T_{0j}^\alpha \right)+ \sum_{\beta=1}^p \left[ \psi, h_{0jk}^{\alpha\beta} \right] w_{0k}^\beta$, i.e. for $\alpha=1,....,q$,
\begin{align*}
\sum_{\beta=1}^q k_j^{\alpha \beta} w_{0j}^\beta = \sum_{\beta, \eta=1}^q h_{0jk}^{\alpha\beta} k_k^{\beta \eta} w_{0k}^\eta +\sum_{\beta=1}^q \left[ \psi, h_{0jk}^{\alpha\beta} \right] w_{0k}^\beta \iff  \sum_{\beta=1}^q k_j^{\alpha\beta} h_{0jk}^{\beta \eta}= \sum_{\beta=1}^q h_{0jk}^{\alpha\beta} k_k^{\beta \eta} + \left[ \psi, h_{0jk}^{\alpha \eta} \right] 
\end{align*}
satisfying
\begin{enumerate}
\item $\sum_{\gamma=1}^q w_{0j}^1 \wedge \cdots \wedge \Psi_j^\gamma \wedge \cdots \wedge w_{0j}^q \wedge dw_{0j}^\alpha + w_{0j}^1 \wedge \cdots \wedge w_{0j}^q \wedge \partial \Psi_j^\alpha = 0 , \alpha=1,...,q$.
\item $-\bar{\partial} \Psi_j\left( w_{0j}^\alpha \right) +\mathcal{L}_\psi\left( w_{0j}^\alpha \right)- \sum_{\beta=1}^q k_j^{\alpha\beta} w_{0j}^\beta, \alpha=1,...,q$
\item $-\bar{\partial} \psi=0$ and $-\bar{\partial} b_j^{\alpha\beta}=0, \alpha,\beta=1,..., q$.
\end{enumerate}
we define
\begin{align*}
l_j:\Gamma\left( U_j, \mathcal{N}_{\mathcal{F}_0}^* \right) &\to \Gamma\left( U_j, \mathcal{A}^{0,0}\left( \bigwedge^{q+2} \Omega_M^1 \otimes \left(\bigwedge^q \mathcal{N}_{\mathcal{F}_0}^* \right)^* \right) \right) \\
 w_{0j}^\alpha &\mapsto - 2\sum_{\gamma=1}^q w_{0j}^1 \wedge \cdots \wedge \Psi_j^\gamma \wedge \cdots \wedge w_{0j}^q \wedge \partial \Psi_j^\gamma - 2\sum_{1\leq \gamma < \eta \leq q} w_{0j}^1 \wedge \cdots \wedge \Psi_j^\gamma \wedge \cdots \wedge \Psi_j^\eta \wedge \cdots \wedge w_{0j}^q \wedge dw_{0j}^\alpha
\end{align*}
Then $\tilde{l}=\{l_j\}$ defines a global section in $A^{0,0}\left( M, \mathscr{H}om_{\mathcal{O}_M}\left(\mathcal{N}_{\mathcal{F}_0}^*,  \bigwedge^{q+2} \Omega_M^1 \otimes \left(\bigwedge^q \mathcal{N}_{\mathcal{F}_0}^* \right)^* \right)   \right)$ and the restriction $l$ of $\tilde{l}$ on $M-S$ defines a global section $A^{0,0}\left( M-S, \mathscr{H}om_{\mathcal{O}_M}\left( \mathcal{N}_{\mathcal{F}_0}^*, \tilde{\mathcal{S}}^2 \right) \right)$, and we define $\kappa_i \in \Gamma\left( U_j, \mathcal{A}^{0,1}\left( \mathscr{H}om_{\mathcal{O}_M}\left( \mathcal{N}_{\mathcal{F}_0}^*, \Omega_M^1 \right) \right) \right)$ by
\begin{align*}
\kappa_i :  \Gamma\left( U_j , \mathcal{N}_{\mathcal{F}_0}^* \right) &\to  \Gamma\left( U_j, \mathcal{A}^{0,1} \left( \Omega_M^1 \right) \right) \\
 w_{0j}^\alpha &\mapsto - 2 \mathcal{L}_\psi \left( \Psi_j^\alpha \right) + 2 \sum_{\eta=1}^q k_j^{\alpha\eta} \wedge \Psi_j^\eta
\end{align*}
Then $\kappa=\{\kappa_i\}$ defines a global section in $A^{0,1}\left( M, \mathscr{H}om_{\mathcal{O}_M}\left( \mathcal{N}_{\mathcal{F}_0}^*, \Omega_M^1 \right) \right)$, and we have $\left[\left(\psi, k_j \right), \left( \psi, k_j \right) \right] \in A^{0,2}\left(M, \mathcal{E}_{\mathcal{N}_{\mathcal{F}_0}^*} \right)$. Then
{\small{\begin{align}\label{it13}
\left( l, \kappa, \left\{ - \left[\left(\psi, k_j \right), \left(\psi, k_j \right) \right]   \right\} \right) \in A^{0,0}\left( M-S, \mathscr{H}om_{\mathcal{O}_M}\left(\mathcal{N}_{\mathcal{F}_0}^* , \tilde{\mathcal{S}}^2 \right) \right) \bigoplus A^{0,1}\left( M, \mathscr{H}om_{\mathcal{O}_M}\left( \mathcal{N}_{\mathcal{F}_0}^*, \Omega_M^1 \right) \right) \bigoplus A^{0,2}\left( M, \mathcal{E}_{\mathcal{N}_{\mathcal{F}_0}^*} \right)
\end{align}}}
defines a $2$-cocycle in the double complex $(\ref{oci15})$.

\end{lemma}

\begin{remark}
As in the same way in \textnormal{Remark \ref{oc81}}, we can show that $l_j\left(w_{0j}^\alpha \right)$ on $U_j-S$ actually lies in $\Gamma\left( U_j-S, \mathcal{A}^{0,0}\left(\tilde{\mathcal{S}}^2 \right) \right)$, i.e. $l_j|_{U_j-S} \in \Gamma\left( U_j-S, \mathcal{A}^{0,0}\left( \tilde{\mathcal{S}}^2 \right) \right)$. More precisely, for $x\in U_j-S$, choose a neighborhood $U_x\subset U_j-S$ such that $\partial \Psi_j^\alpha= \sum_{\delta=1}^q a_{j_x}^{\alpha\beta} \wedge \Psi_j^\beta= \sum_{\beta=1}^q b_{j_x}^{\alpha\beta}\in \Gamma\left(U_x, \mathcal{A}^{0,0}\left( \Omega_M^1 \right) \right)$. Then we can write
\begin{align} \label{its43}
l_j\left( w_{0j}^\alpha \right) &= -2 \sum_{\gamma=1}^q w_{0j}^1 \wedge \cdots \wedge \Psi_j^\gamma \wedge \cdots \wedge w_{0j}^q \wedge dw_{0j}^\alpha - 2 \sum_{1\leq \gamma < \leq \eta \leq q} w_{0j}^1 \wedge \cdots \wedge \Psi_j^\gamma \wedge \cdots \wedge \Psi_j^\eta \wedge \cdots \wedge w_{0j}^q \wedge dw_{0j}^\alpha \\
 &= 2 w_{0j}^1 \wedge \cdots \wedge w_{0j}^q \wedge \left( \sum_{\gamma=1}^q b_{j_x}^{\alpha \gamma} \wedge \Psi_j^\gamma \right) \in \Gamma\left( U_x, \mathcal{A}^{0,0}\left( \tilde{\mathcal{S}}^2 \right) \right) \notag
\end{align}
\end{remark}

\begin{remark}
We show that $\tilde{l}$ defines a global section in $A^{0,0}\left( M, \mathscr{H}om_{\mathcal{O}_M}\left( \mathcal{N}_{\mathcal{F}_0}^*, \bigwedge^{q+2}\Omega_M^1\otimes \mathcal{L}_0 \right) \right)$. In fact,
{\small{\begin{align*}
&2\sum_{\gamma=1}^q w_{0j}^1 \wedge \cdots \wedge \Psi_j^\gamma \wedge \cdots \wedge w_{0j}^q \wedge \partial \Psi_j^\gamma + 2\sum_{1\leq \gamma < \eta \leq q} w_{0j}^1 \wedge \cdots \wedge \Psi_j^\gamma \wedge \cdots \wedge \Psi_j^\eta \wedge \cdots \wedge w_{0j}^q \wedge dw_{0j}^\alpha\\
&=2 H_{0jk}\sum_{\gamma=1}^q w_{0k}^1 \wedge \cdots \wedge \Psi_k^\gamma \wedge \cdots \wedge w_{0k}^q \wedge \partial\left( \sum_{\beta=1}^q h_{0jk}^{\alpha\beta} \Psi_k^\beta \right) + 2H_{0jk}\sum_{1\leq \gamma < \eta \leq q} w_{0k}^1 \wedge \cdots \wedge \Psi_k^\gamma \wedge \cdots \wedge \Psi_k^\eta \wedge \cdots \wedge w_{0k}^q \wedge d\left( \sum_{\beta=1}^q h_{0jk}^{\alpha\beta} w_{0k}^\beta \right)\\
&= H_{0jk} \sum_{\beta=1}^q h_{0jk}^{\alpha\beta} \left(2\sum_{\gamma=1}^q  w_{0k}^1 \wedge \cdots \wedge \Psi_k^\gamma \wedge \cdots \wedge w_{0j}^q \wedge  \partial \Psi_k^\beta +2 \sum_{1\leq \gamma < \eta \leq q} w_{0k}^1 \wedge \cdots \wedge \Psi_k^\gamma \wedge \cdots \wedge \Psi_k^\eta \wedge \cdots \wedge w_{0k}^q \wedge dw_{0k}^\beta \right)
\end{align*}}}

\end{remark}

\begin{remark}
The following computations will be used in the proof of \textnormal{Lemma \ref{oci85}}. We note that for $x\in U_j-S$, choose a neighborhood $U_x\subset U_j-S$ such that we have
\begin{align}\label{oci20}
\partial \Psi_j^\alpha- \sum_{\beta=1}^q a_{j_x}^{\alpha\beta} \wedge \Psi_j^\beta = \sum_{\beta=1}^q b_{j_x}^{\alpha\beta} \wedge w_{0j}^\beta,\,\,\,\,\,\,\textnormal{for some}\,\,\,b_{j_x}^{\alpha\beta} \in \Gamma\left( U_x, \mathcal{A}^{0,0}\left(\Omega_M^1 \right) \right)
\end{align}

By taking $\bar{\partial}$ on $(\ref{oci20})$, we have
{\small{\begin{align*}
&-\partial \bar{\partial} \Psi_j^\alpha + \sum_{\beta=1}^q a_{j_x}^{\alpha\beta} \wedge \bar{\partial} \Psi_j^\beta = \sum_{\beta=1}^q \bar{\partial} b_{j_x}^{\alpha\beta} \wedge w_{0j}^\beta\\
&\Longrightarrow -\partial\left(\mathcal{L}_{\psi}\left(w_{0j}^\alpha \right) - \sum_{\beta=1}^q k_j^{\alpha\beta} \wedge w_{0j}^\beta  \right) + \sum_{\beta=1}^q a_{j_x}^{\alpha\beta}\wedge \left( \mathcal{L}_{\psi}\left(w_{0j}^\beta \right)  - \sum_{\gamma=1}^q k_j^{\beta \gamma} \wedge w_{0j}^\beta  \right) = \sum_{\beta=1}^q \bar{\partial} b_{j_x}^{\alpha\beta} \wedge w_{0j}^\beta\\
&\Longrightarrow  \mathcal{L}_{\psi}\left(\sum_{\beta=1}^q a_{j_x}^{\alpha\beta} \wedge w_{0j}^\beta \right) + \sum_{\beta=1}^q \partial k_j^{\alpha\beta} \wedge w_{0j}^\beta - \sum_{\beta,\gamma=1}^q k_j^{\alpha\beta} \wedge a_{j_x}^{\beta \gamma} \wedge w_{0j}^\gamma + \sum_{\beta=1}^q a_{j_x}^{\alpha\beta} \wedge \mathcal{L}_{\psi}\left( w_{0j}^\beta \right) - \sum_{\beta, \gamma=1}^q a_{j_x}^{\alpha\beta} \wedge k_j^{\beta \gamma} \wedge w_{0j}^\gamma = \sum_{\beta=1}^q \bar{\partial} b_{j_x}^{\alpha\beta} \wedge w_{0j}^\beta \\
&\Longrightarrow \sum_{\beta=1}^q \left( \mathcal{L}_{\psi}\left(a_{j_x}^{\alpha\beta} \right)  + \partial k_j^{\alpha\beta} - \sum_{\gamma=1}^q k_j^{\alpha \gamma} \wedge a_{j_x}^{\gamma \beta} - \sum_{\gamma=1}^q a_{j_x}^{\alpha \gamma} \wedge k_j^{\gamma \beta} - \bar{\partial} b_{j_x}^{\alpha\beta}      \right) \wedge w_{0j}^\beta = 0
\end{align*}}}

This implies that 
\begin{align}\label{its40}
\bar{\partial} b_{j_x}^{\alpha\beta} \equiv_{\mathcal{N}_{\mathcal{F}_0}^*} \mathcal{L}_{\psi}\left( a_{j_x}^{\alpha\beta} \right)  + \partial k_j^{\alpha\beta} - \sum_{\gamma=1}^q k_j^{\alpha \gamma} \wedge a_{j_x}^{\gamma \beta} - \sum_{\gamma=1}^q a_{j_x}^{\alpha \gamma} \wedge k_j^{\gamma \beta}
\end{align}

\end{remark}

\begin{proof}[Proof of \textnormal{Lemma \ref{oci85}}]
Since $\bar{\partial}\left( \psi, k_j \right)=0$, we have $\bar{\partial}\left( -\left[ \left(\psi, k_j \right),\left(\psi, k_j \right) \right]   \right)=0$. Since $-\bar{\partial}\Psi_j^\alpha + \mathcal{L}_{\psi}\left( w_{0j}^\alpha \right) - \sum_{\beta=1}^q k_j^{\alpha\beta} \wedge w_{0j}^\beta = 0$, we have from $(\ref{its40})$
\begin{align}\label{its41}
\bar{\partial}\left(\sum_{\beta=1}^q b_{j_x}^{\alpha\beta}\wedge \Psi_j^\beta \right) &=  \sum_{\beta=1}^q \bar{\partial} b_{j_x}^{\alpha \beta}  \Psi_j^\beta  -\sum_{\beta=1}^q b_{j_x}^{\alpha\beta} \wedge \bar{\partial}\Psi_j^\beta  \\ &\equiv_{\mathcal{N}_{\mathcal{F}_0}^*}  \sum_{\beta=1}^q\left(  \mathcal{L}_{\psi}\left( a_{j_x}^{\alpha\beta} \right)  + \partial k_j^{\alpha\beta} - \sum_{\gamma=1}^q k_j^{\alpha \gamma} \wedge a_{j_x}^{\gamma \beta} - \sum_{\gamma=1}^q a_{j_x}^{\alpha \gamma} \wedge k_j^{\gamma \beta} \right)  \Psi_j^\beta - \sum_{\beta=1}^q b_{j_x}^{\alpha\beta} \wedge \mathcal{L}_{\psi}\left( w_{0j}^\beta \right) \notag
\end{align}
On the other hand, we have
\begin{align}\label{its42}
&-\partial\left(-\mathcal{L}_{\psi}\left(\Psi_j^\alpha \right) +\sum_{\eta=1}^q k_j^{\alpha \eta} \Psi_j^\eta\right)+ \sum_{\beta=1}^q a_{j_x}^{\alpha\beta} \wedge\left(-\mathcal{L}_{\psi} \left( \Psi_j^\beta \right)  + \sum_{\eta=1}^q k_j^{\beta \eta}  \Psi_j^\eta    \right)\\
&= -\mathcal{L}_{\psi}\left(\partial \Psi_j^\alpha \right) - \sum_{\eta=1}^q \partial k_j^{\alpha \eta}  \Psi_j^\eta  + \sum_{\eta=1}^q k_j^{\alpha \eta} \wedge \partial \Psi_j^\eta - \sum_{\beta=1}^q a_{j_x}^{\alpha \beta} \wedge \mathcal{L}_{\psi}\left( \Psi_j^\beta \right) + \sum_{\beta, \eta=1}^q a_{j_x}^{\alpha \beta} \wedge k_j^{\beta \eta} \Psi_j^\eta \notag\\
&= - \mathcal{L}_{\psi}\left( \sum_{\beta=1}^q a_{j_x}^{\alpha\beta} \wedge \Psi_j^\beta + \sum_{\beta=1}^q b_{j_x}^{\alpha\beta} \wedge w_{0j}^\beta \right) - \sum_{\eta=1}^q \partial k_j^{\alpha \eta} \Psi_j^\eta + \sum_{\eta=1}^q k_j^{\alpha \eta } \wedge \left( \sum_{\beta=1}^q a_{j_x}^{\eta\beta} \wedge \Psi_j^\beta + \sum_{\beta=1}^q b_{j_x}^{\eta \beta} \wedge w_{0j}^\beta  \right) \notag \\
&- \sum_{\beta=1}^q a_{j_x}^{\alpha \beta} \wedge \mathcal{L}_{\psi}\left( \Psi_j^\beta \right) + \sum_{\beta=1}^q a_{j_x}^{\alpha \beta} \wedge k_j^{\beta \eta} \wedge \Psi_j^\eta \notag \\
&\equiv_{\mathcal{N}_{\mathcal{F}_0}^*} -\sum_{\beta=1}^q \mathcal{L}_\psi\left(a_{j_x}^{\alpha\beta} \right) \wedge \Psi_j^\beta + \sum_{\beta=1}^q b_{j_x}^{\alpha\beta}  \wedge \mathcal{L}_\psi\left( w_{0j}^\beta \right) - \sum_{\eta=1}^q \partial k_j^{\alpha\eta}\Psi_j^\eta + \sum_{\eta, \beta=1}^q k_j^{\alpha \eta} \wedge a_{j_x}^{\eta \beta} \wedge \Psi_j^\beta + \sum_{\beta,\eta=1}^q a_{j_x}^{\alpha\beta}\wedge k_j^{\beta \eta} \wedge \Psi_j^\eta \notag
\end{align}
Then $(\ref{its43}), (\ref{its41})$ and $(\ref{its42})$ imply that $(-1)^q\bar{\partial} l_j + \hat{E}_1\left( \Psi \right)=0$.

We show that $\bar{\partial} \left( \kappa \right) + E_0'\left( -\left[ (\psi, k_j), (\psi ,k_j) \right] \right)=0$. In fact,

\begin{align*}
\bar{\partial} \kappa_j\left(w_{0j}^\alpha \right) & = \bar{\partial} \left( -2 \mathcal{L}_{\psi}\left( \Psi_j^\alpha  \right) + 2 \sum_{\eta=1}^q k_j^{\alpha \eta} \wedge \Psi_j^\eta \right) = 2 \mathcal{L}_{\psi}\left( \bar{\partial} \Psi_j^\alpha  \right) +2 \sum_{\eta=1}^q \bar{\partial} k_j^{\alpha \eta} \Psi_j^\eta -2 \sum_{\eta=1}^q k_j^{\alpha \eta } \wedge \bar{\partial} \Psi_j^\eta \\
&= 2\mathcal{L}_{\psi}\left(  \mathcal{L}_{\psi} \left( w_{0j}^\alpha \right) - \sum_{\beta=1}^q k_j^{\alpha \beta} \wedge w_{0j}^\beta  \right) - 2\sum_{\eta=1}^q k_j^{\alpha \eta} \wedge \left( \mathcal{L}_\psi \left(w_{0j}^\eta \right) -  \sum_{\beta=1}^q  k_j^{\eta \beta} \wedge w_{0j}^\beta  \right) \\
&=2 \mathcal{L}_\psi \mathcal{L}_\psi \left(w_{0j}^\alpha \right) - 2 \sum_{\beta=1}^q \left[ \psi , k_j^{\alpha \beta}\right] \wedge w_{0j}^\beta + 2\sum_{\eta, \beta =1}^q k_j^{\alpha \eta} \wedge k_j^{\eta \beta} \wedge w_{0j}^\beta\\
&= \mathcal{L}_{\left[ \psi, \psi \right]} \left( w_{0j}^\alpha \right)  + 2\sum_{\beta=1}^q \left(-\left[ \psi, k_j^{\alpha\beta} \right] + \sum_{\eta=1}^q k_j^{\alpha \eta } \wedge k_j^{\eta \beta} \right) \wedge w_{0j}^\beta
\end{align*}
From $(\ref{its43})$, as in the same way with $(\ref{its46})$, we can show that $\hat{E}_2\left(l_j \right)=0$ on $U_j-S$, and so $M-S$. This completes the proof of Lemma \ref{oci85}.

\end{proof}

From $(\ref{it31})-(\ref{it32})$ and Lemma \ref{oci85}, we obtain
\begin{theorem}
Suppose given a compact foliated complex manifold $\left( M, \mathcal{N}_{\mathcal{F}_0}^* \right)$ with $\mathcal{N}_{\mathcal{F}_0}^*$ locally free, and $\left(\Psi, \left\{\left( \psi, k_j \right) \right\} \right) \in \mathbb{H}^1\left( M, \mathcal{E}_{\mathcal{N}_{\mathcal{F}_0}^*}^\bullet \right)$. In order that there may exist a foliated analytic family $\left( \mathcal{M}, \mathcal{N}_\mathcal{F}^*, B, \pi \right)$ in terms of cotangent sheaves such that $\pi^{-1}(0)= \left( M, \mathcal{N}_{\mathcal{F}_0}^* \right)$ and that $\left( \frac{\partial\left( \mathcal{M}, \mathcal{N}_\mathcal{F}^* \right)}{\partial t}\right)|_{t=0}= \left( \Psi, \left\{ \left( \psi, k_j \right) \right\} \right)$, it is necessary that
{\small{\begin{align*}
\left( l, \kappa, \left\{ - \left[\left(\psi, k_j \right), \left(\psi, k_j \right) \right]   \right\} \right) \in A^{0,0}\left( M-S, \mathscr{H}om_{\mathcal{O}_M}\left(\mathcal{N}_{\mathcal{F}_0}^* , \tilde{\mathcal{S}}^2 \right) \right) \bigoplus A^{0,1}\left( M, \mathscr{H}om_{\mathcal{O}_M}\left( \mathcal{N}_{\mathcal{F}_0}^*, \Omega_M^1 \right) \right) \bigoplus A^{0,2}\left( M, \mathcal{E}_{\mathcal{N}_{\mathcal{F}_0}^*} \right)
\end{align*}}}
from $(\ref{it13})$ in \textnormal{Lemma \ref{oci85}} defines the $0$ class in $\mathbb{H}^2\left( M, \mathcal{E}_{\mathcal{N}_{\mathcal{F}_0}^*}^\bullet \right)$.
\end{theorem}

\section{Isomorphisms of infinitesimal foliated deformations and obstructions in terms of cotangent sheaves between \v Cech resolution and Dolbeault type bicomplex associated to $\mathcal{E}_{\mathcal{N}_{\mathcal{F}_0}^*}^\bullet$}\label{ni2}

In this section we explicitly construct isomorphisms of first and second cohomology groups between \v Cech resolution of $\mathcal{E}_{\mathcal{N}_{\mathcal{F}_0}^*}^\bullet$ and \v Dolbeault type bicomplex associated to $\mathcal{E}_{\mathcal{N}_{\mathcal{F}_0}^*}^\bullet$. Let $\mathcal{U}=\left\{U_j \right\}$ be an Stein open covering of $M$ by coordinate neighborhoods with coordinates $z_j=\left(z_j^1,..., z_j^n \right)$ on $U_j$, and $\mathcal{N}_{\mathcal{F}_0}^*$ is generated by $w_{0j}^1,..., w_{0j}^p \in \Gamma\left( U_j , \Omega_M^1 \right)$ with the relation $w_{0j}^\alpha= \sum_{\beta=1}^q h_{0jk}^{\alpha\beta} w_{0j}^\beta$ on $U_{jk}$ for some $h_{0jk}^{\alpha\beta}\in \Gamma\left( U_{jk}, \mathcal{O}_M \right)$ and satisfies $w_{0j}^1 \wedge \cdots \wedge w_{0j}^q \wedge dw_{0j}^\alpha=0,\alpha=1,...,q$.  

\subsection{Isomorphisms of infinitesimal deformations between \v Cech resolution and Dolbeault type bicomplex associated to $\mathcal{E}_{\mathcal{N}_{\mathcal{F}_0}^*}^\bullet$ }\label{ni1}\

We will explicitly define an isomorphism between the first cohomology group $\mathbb{H}^1\left( M, \mathcal{E}_{\mathcal{N}_{\mathcal{F}_0}^*}^\bullet \right)$ from \v Cech resolution of $\mathcal{E}_{\mathcal{N}_{\mathcal{F}_0}^*}^\bullet$ and the first cohomology group $B^1$ from $(\ref{ii1})$, and will show that the infinitesimal foliated deformations $\left(\left\{ \beta_j(0)\right\}, \left\{\left( \theta_{jk}(0), B_{jk}(0) \right) \right\} \right)$ $(\ref{oc7})$ from \v Cech resolution of $\mathcal{E}_{\mathcal{N}_{\mathcal{F}_0}^*}^\bullet$ and  $\left( \frac{\partial Y}{\partial t}|_{t=0}, \left\{ \left( \frac{\partial \varphi}{\partial t}|_{t=0}, \frac{\partial Z_j}{\partial t}|_{t=0} \right) \right\} \right)$ $(\ref{ii2})$ from the Dolbeaul type bicomplex under the isomorphism. 

Let $\left( \left\{ \beta_j\right\}, \left\{ \left( \theta_{jk}, B_{jk} \right)      \right\}\right) \in C^0\left( \mathcal{U}, \mathscr{H}om_{\mathcal{O}_M}\left( \mathcal{N}_{\mathcal{F}_0}^*, \Omega_M^1 \right) \right) \bigoplus C^1\left( \mathcal{U}, \mathcal{E}_{\mathcal{N}_{\mathcal{F}_0}^*} \right)$ be an $1$-cocycle in the \v Cech resolution of $\mathcal{E}_{\mathcal{N}_{\mathcal{F}_0}^*}^\bullet$ with respect to the covering $\mathcal{U}$, so that $\left(\theta_{ij}, B_{ij}\right)- \left(\theta_{ik}, B_{ik}\right)+ \left(\theta_{jk}, B_{jk}^i\right)=0$ and $\beta_k\left( w_{0j}^\alpha \right) - \beta_j\left(w_{0j}^\alpha \right)+ \mathcal{L}_{\theta_{jk}}\left(w_{0j}^\alpha \right)- B_{jk}\left( w_{0j}^\alpha \right)=0$, where $B_{jk}^\alpha=B_{jk}\left(w_{0j}^\alpha \right) = \sum_{\beta=1}^q B_{jk}^{\alpha \delta} w_{0j}^\delta$ with $B_{jk}^{\alpha\delta}\in \Gamma\left(U_{jk}, \mathcal{O}_M \right)$, and $\sum_{\beta=1}^q w_{0j}^1 \wedge \cdots \wedge \beta_j\left( w_{0j}^\beta\right)\wedge \cdots \wedge w_{0j}^q \wedge dw_{0j}^\alpha + w_{0j}^1 \wedge \cdots \wedge w_{0j}^q \wedge d\beta_j\left(w_{0j}^\alpha \right)=0, \alpha=1,...,q$. Since $\delta\left( \left\{\theta_{jk}, B_{jk} \right\}\right)=0$, there exists $-\delta\left(\left\{ \left(\xi_j, B_j\right) \right\} \right)= \left\{ \left( \theta_{jk}, B_{jk}\right) \right\}$ where $\left\{\left( \xi_j, B_j\right) \right\}\in C^0\left(\mathcal{U}, \mathcal{A}^{0,0}\left( \mathcal{E}_{\mathcal{N}_{\mathcal{F}_0}^*} \right) \right)$ such that $\xi_j\in \Gamma\left( U_j, \mathcal{A}^{0,0}\left(\Theta_M \right) \right)$ and $B_j\in \Gamma\left( U_j, \mathcal{A}^{0,0}\left( \mathscr{H}om_{\mathcal{O}_M}\left(\mathcal{N}_{\mathcal{F}_0}^*, \mathcal{N}_{\mathcal{F}_0}^* \right)\right) \right)$ with $B_j\left(w_{0j}^\alpha \right) = \sum_{\beta=1}^q B_j^{\alpha \beta} w_{0j}^\beta$ satisfying $\xi_j-\xi_k= \theta_{jk}$ and $B_j-B_k^j= B_{jk}$ with $B_k^j\left(w_{0j}^\alpha \right)= B_k \left( w_{0j}^\alpha \right)+ \sum_{\beta=1}^q \left[ \xi_k, h_{0jk}^{\alpha\beta} \right] w_{0k}^\beta$. Then $\left\{\left( \psi =\bar{\partial} \xi_j, \bar{\partial} B_j \right) \right\}\in A^{0,0}\left(M, \mathcal{E}_{\mathcal{N}_{\mathcal{F}_0}^*} \right)$ where $\phi\in A^{0,1}\left( M, \Theta_M \right)$ and $\bar{\partial} B_j\in \Gamma\left(U_j, \mathcal{A}^{0,1}\left( \mathscr{H}om_{\mathcal{O}_M}\left( \mathcal{N}_{\mathcal{F}_0}^*, \mathcal{N}_{\mathcal{F}_0}^* \right) \right) \right)$. On the other hand, let us define $\Psi_j \in \Gamma\left(U_j, \mathcal{A}^{0,0}\left(\mathscr{H}om_{\mathcal{O}_M}\left( \Theta_{\mathcal{F}_0}, \Omega_M^1 \right) \right) \right)$ by
\begin{align*}
\Psi_j:\Gamma\left(U_j,  \mathcal{N}_{\mathcal{F}_0}^* \right) &\to \Gamma\left( U_j, \Omega_M^1 \right) \\
 w_{0j}^\alpha &\mapsto \mathcal{L}_{\xi_j}\left(w_{0j}^\alpha \right) - B_j\left( w_{0j}^\alpha \right) - \beta_j\left(w_{0j}^\alpha \right)
\end{align*}
Then we have
\begin{align*}
&\Psi_j\left(w_{0j}^\alpha \right) - \sum_{\beta=1}^q h_{0jk}^{\alpha\beta} \Phi_k\left(w_{0k}^\beta \right) = \mathcal{L}_{\xi_j}\left(w_{0j}^\alpha \right) - B_j\left( w_{0j}^\alpha \right) - \beta_j\left(w_{0j}^\alpha \right) - \sum_{\beta=1}^q h_{0jk}^{\alpha\beta} \left(\mathcal{L}_{\xi_k}\left(w_{0k}^\beta \right) - B_k\left( w_{0j}^\beta \right) - \beta_k\left(w_{0k}^\beta \right)\right)\\
&= \mathcal{L}_{\xi_j}\left( w_{0j}^\alpha \right)  - B_j\left( w_{0j}^\alpha \right) - \beta_j\left(w_{0j}^\alpha \right) - \mathcal{L}_{\xi_k}\left(  w_{0j}^\alpha \right)+ \sum_{\beta=1}^q \left[ \xi_k, h_{0jk}^{\alpha\beta} \right] w_{0k}^\beta  + B_k\left(w_{0j}^\alpha \right)  + \beta_k\left( w_{0j}^\alpha \right) \\
&= \mathcal{L}_{\theta_{jk}}\left( w_{0j}^\alpha \right) - B_{jk}\left(w_{0j}^\alpha \right) - \beta_j\left( w_{0j}^\alpha \right) + \beta_k\left( w_{0j}^\alpha \right)=0
\end{align*}
This implies that $\Psi:=\left\{ \Psi_j \right\}\in A^{0,0}\left(M, \mathscr{H}om_{\mathcal{O}_M}\left(\mathcal{N}_{\mathcal{F}_0}^* , \Omega_M^1 \right) \right)$. Then we claim that $\left( \Psi, \left\{ \left(  \psi ,  \bar{\partial } B_j  \right) \right\} \right)$ defines a $1$-cocycle in the Dolbeault type bicomplex. In fact, it is enough to check that
\begin{align*}
-\bar{\partial}\left( \Psi_j \left( w_{0j}^\alpha \right) \right) + \mathcal{L}_{\psi}\left(w_{0j}^\alpha \right) - \bar{\partial} B_j\left( w_{0j}^\alpha \right) = - \bar{\partial}\left(\mathcal{L}_{ \xi_j}\left( w_{0j}^\alpha \right) - B_j\left(w_{0j}^\alpha \right) \right) + \mathcal{L}_{ \bar{\partial} \xi_j}\left( w_{0j}^\alpha \right) - \bar{\partial} B_j\left(w_{0j}^\alpha \right) = 0
\end{align*}

Then we define the correspondence of $1$-cocycle between \v Cech resolution and Dolbeault type bicomplex by 
\begin{align}\label{ii5}
\left(\left\{\beta_j \right\}, \left\{ \left( \theta_{jk}, B_{jk} \right) \right\} \right) \mapsto \left( \Phi, \left\{ \left(\psi, \bar{\partial} B_j\right) \right\}\right).
\end{align}

\begin{proposition}\label{c1}

We recall that given a foliated analytic family $\left( \mathcal{M}, \mathcal{N}_\mathcal{F}^*, B, \pi \right)$ of deformations of $\left( M, \mathcal{N}_{\mathcal{F}_0}^* \right)=\pi^{-1}(0)$ as in \textnormal{Remark \ref{oc2}}, the infinitesimal deformation at $t=0$ is given by $\left( \left\{ \beta_j(0) \right\}, \left\{\left( \theta_{jk}(0), B_{jk}(0) \right) \right\} \right)$ from $(\ref{oc7})$ in the first cohomology group of the \v Cech resolution of $\mathcal{E}_{\mathcal{N}_{\mathcal{F}_0}^*}^\bullet$. On other hand, the infinitesimal deformation at $t=0$ is also given by $\left( \frac{\partial Y}{\partial t}|_{t=0}, \left\{ \left(\frac{\partial \varphi}{\partial t}|_{t=0}, \frac{\partial Z_j}{\partial t}|_{t=0} \right) \right\} \right)$ from $(\ref{ii2})$ in the first cohomology group $B^1$ of the Dolbeault type bicomplex associated to $\mathcal{E}_{{N}_{\mathcal{F}_0}^*}^\bullet$. Under the isomorphism $(\ref{ii5})$ from $\mathbb{H}^1\left(M, \mathcal{E}_{\mathcal{N}_{\mathcal{F}_0}^*}^\bullet \right)$ to $\textnormal{B}^1$, 
\begin{align*}
\left( \left\{ \beta_j(0) \right\}, \left\{\left( \theta_{jk}(0), B_{jk}(0) \right) \right\} \right)\,\,\,\,\,\,\,\textnormal{corresponds to}\,\,\,\,\,\, \left( \frac{\partial Y}{\partial t}|_{t=0}, \left\{ \left(\frac{\partial \varphi}{\partial t}|_{t=0}, \frac{\partial Z_j}{\partial t}|_{t=0} \right) \right\} \right)
\end{align*}

\end{proposition}

\begin{proof}
We fix a tangent vector $\frac{\partial}{\partial t}\in T_0\left( \Delta\right)$, denote $\left(\frac{\partial f(t)}{\partial t} \right)_{t=0}$ by $\dot{f}$ for a $C^\infty$ function $f(t),t\in \Delta$. With this notation, we put
\begin{align*}
\xi_j=\sum_{\alpha=1}^n \dot{\xi}_j^\alpha \frac{\partial}{\partial \xi_j^\alpha},\,\,\,\,\,\textnormal{where}\,\,\,\dot{\xi}_j^\alpha = \left(  \frac{\partial \xi_j(z,t)}{\partial t} \right)_{t=0}
\end{align*}
for each $j$. Then we have
\begin{align*}
\xi_j- \xi_k = \theta_{jk},\,\,\,\,\,\, \textnormal{and}\,\,\,\,\,\bar{\partial} \xi_j= \frac{\partial \varphi}{\partial t}|_{t=0}
\end{align*}
(for the detail, see \cite{Kod05} Theorem 5.4 p.266). Differentiating $(\ref{b19})$
\begin{align*}
\sum_{\beta=1}^q \frac{\partial g_j^{\alpha\beta}}{\partial t}\cdot h_{0jk}^{\beta \eta}= \sum_{\beta,\gamma=1}^q \frac{\partial h_{jk}^{\alpha\beta}}{\partial \xi_k^\gamma} \frac{\partial \xi_k^\gamma}{\partial t} \cdot g_k^{\beta \eta} +\sum_{\beta=1}^q \frac{\partial h_{jk}^{\alpha\beta}}{\partial t}\cdot g_k^{\beta \eta} + \sum_{\beta=1}^q h_{jk}^{\alpha\beta}\left( \xi_k(z,t,),t \right) \cdot \frac{\partial g_k^{\beta \eta}}{\partial t}
\end{align*}
and setting $t=0$, we obtain
\begin{align*}
&\sum_{\beta=1}^q \frac{\partial g_j^{\alpha\beta}}{\partial t}|_{t=0}\cdot h_{0jk}^{\beta \eta} = \sum_{\beta=1}^q \frac{\partial h_{jk}^{\alpha\beta}}{\partial \xi_k^\eta}|_{t=0}\cdot \frac{\partial \xi_k^\eta}{\partial t}|_{t=0} +  \frac{\partial h_{jk}^{\alpha\eta}}{\partial t}|_{t=0} + \sum_{\beta=1}^q h_{0jk}^{\alpha\beta}\cdot \frac{\partial g_k^{\beta \eta}}{\partial t}|_{t=0} \\
&\Longrightarrow \frac{\partial g_j^{\alpha\beta}}{\partial t}|_{t=0} = \sum_{\eta=1}^q \left[ \xi_k, h_{0jk}^{\alpha\eta} \right] h_{0kj}^{\eta \beta} + \sum_{\eta=1}^q \frac{\partial h_{jk}^{\alpha \eta}}{\partial t}|_{t=0}\cdot h_{0kj}^{\eta \beta} + \sum_{\eta, \gamma=1}^q h_{0jk}^{\alpha\eta} \cdot \frac{\partial g_k^{\eta \gamma}}{\partial t}|_{t=0} \cdot h_{0kj}^{\gamma \beta} \\
&\Longrightarrow \sum_{\beta=1}^q \frac{\partial g_j^{\alpha\beta}}{\partial t}|_{t=0} w_{0j}^\beta= \sum_{\eta=}^q \left[ \xi_k, h_{0jk}^{\alpha\eta} \right] w_{0k}^\eta +\sum_{\eta=1}^q \frac{\partial h_{jk}^{\alpha\eta}}{\partial t} |_{t=0} w_{0k}^\eta +  \sum_{\eta, \gamma=1}^q h_{0jk}^{\alpha \eta}\cdot \frac{\partial g_k^{\eta \gamma}}{\partial t}|_{t=0} w_{0k}^\gamma
\end{align*}
If we define $B_j\in \Gamma\left( \mathcal{U}, \mathscr{H}om_{\mathcal{O}_M}\left( \mathcal{N}_{\mathcal{F}_0}^* , \mathcal{N}_{\mathcal{F}_0}^* \right) \right)$ by
\begin{align*}
B_j:\Gamma\left( U_j,  \mathcal{N}_{\mathcal{F}_0}^* \right) &\to \Gamma\left( U_j, \mathcal{A}^{0,0}\left( \mathcal{N}_{\mathcal{F}_0}^* \right) \right) \\
              w_{0j}^\alpha &\mapsto \sum_{\beta=1}^q \frac{\partial g_j^{\alpha\beta}}{\partial t}|_{t=0} w_{0j}^\beta
\end{align*}
Then we see that $-\delta\left(  \left\{ \left(  \xi_j, B_j \right) \right\} \right)=\left( \theta_{jk}(0), B_{jk}(0) \right)$ and $\bar{\partial} B_j= \frac{\partial Z_j}{\partial t}|_{t=0}$.

On the other hand, differentiating $(\ref{b20})$, we have
\begin{align*}
\frac{\partial W_j^{\alpha \gamma}}{\partial t}= \sum_{\beta,\eta=1}^n \frac{\partial w_j^{\alpha \beta}}{\partial \xi_j^\eta} \frac{ \partial \xi_j^\eta}{\partial t} \frac{\partial \xi_j^\beta}{\partial z^\gamma} + \sum_{\beta=1}^n \frac{ \partial w_j^{\alpha\beta}}{\partial t} \frac{\partial \xi_j^\beta}{\partial z^\gamma} + \sum_{\beta=1}^n w_j^{\alpha\beta} \frac{\partial}{\partial z^\gamma}\left( \frac{\partial \xi_j^\beta}{\partial t} \right)
\end{align*}
and setting $t=0$, we have
\begin{align*}
&\frac{\partial W_j^{\alpha \gamma}}{\partial t}|_{t=0} = \sum_{\beta=1}^n \frac{\partial w_{0j}^{\alpha\gamma}}{\partial \xi_j^\eta} \frac{\partial \xi_j^\eta}{\partial t}|_{t=0} + \frac{\partial w_j^{\alpha \gamma}}{\partial t}|_{t=0} + \sum_{\beta=1}^n w_{0j}^{\alpha\beta}\frac{\partial}{\partial z^\gamma}\left( \frac{\partial \xi_j^\beta}{\partial t}|_{t=0} \right) \\
&\Longrightarrow \frac{\partial W_j^\alpha}{\partial t}|_{t=0} = \mathcal{L}_{\xi_j} \left( w_{0j}^\alpha \right) + \frac{\partial w_j^\alpha }{\partial t}|_{t=0} \Longrightarrow \frac{\partial W_j^\alpha}{\partial t}|_{t=0} - \sum_{\beta=1}^q \frac{\partial g_j^{\alpha \beta}}{\partial t}|_{t=0} w_{0j}^\beta = \mathcal{L}_{\xi_j}\left( w_{0j}^\alpha \right) + \frac{\partial w_j^\alpha}{\partial t}|_{t=0} - \sum_{\beta=1}^q \frac{\partial g_j^{\alpha\beta}}{\partial t}|_{t=0} w_{0j}^\beta \\
&\Longrightarrow \frac{\partial Y_j^\alpha (t)}{\partial t}|_{t=0} = \mathcal{L}_{\xi_j}\left( w_{0j}^\alpha \right)- B_j \left( w_{0j}^\alpha \right) - \left( - \frac{\partial w_j^\alpha}{\partial t}|_{t=0} \right)
\end{align*}
This completes the proof of Proposition \ref{c1}.

\end{proof}

\subsection{Isomorphism of obstructions between \v Cech resolution and Dolbeault type bicomplex associated to $\mathcal{E}_{\mathcal{N}_{F_0}^*}^\bullet$ } \label{a16} \

We will explicitly define an isomorphism between the second cohomology group $\mathbb{H}^2\left( M, \mathcal{E}_{\mathcal{N}_{\mathcal{F}_0}^*}^\bullet  \right)$ from \v Cech resolution of $\mathcal{E}_{\mathcal{N}_{\mathcal{F}_0}^*}^\bullet$ and the second cohomology group $\textnormal{B}^2$ from the Dolbeault type bicomplex associated to $\mathcal{E}_{\mathcal{N}_{\mathcal{F}_0}^*}^\bullet$ $(\ref{ii1})$.

Let $\left(\left\{ c_i \right\}, \left\{ b_{ij} \right\}, \left\{ a_{ijk}, g_{ijk} \right\} \right)\in C^0\left( \mathcal{U}, \mathscr{H}om_{\mathcal{O}_M}\left( \mathcal{N}_{\mathcal{F}_0}^* ,  \tilde{\mathcal{S}}^2 \right) \right) \bigoplus  C^1\left( \mathcal{U} , \mathscr{H}om_{\mathcal{O}_M}\left( \mathcal{N}_{\mathcal{F}_0}^*, \Omega_M^1 \right) \right) \bigoplus C^2\left( \mathcal{U}, \mathcal{E}_{\mathcal{N}_{\mathcal{F}_0}^*} \right) $ be a $2$-cocycle. Since $-\delta\left( \left\{ \left( a_{ijk}, g_{ijk}\right) \right\} \right)=0$, we have $\left\{\left(d_{ij}, e_{ij} \right) \right\}\in C^1\left(\mathcal{U}, \mathcal{A}^{0,0}\left(\mathcal{E}_{\mathcal{N}_{\mathcal{F}_0}^*} \right) \right)$ such that $\delta\left( \left\{\left( d_{ij}, e_{ij} \right) \right\}\right)=\left\{ \left( a_{ijk}, g_{ijk}\right) \right\}$ where $d_{ij}\in \Gamma\left( U_{ij}, \mathcal{A}^{0,0}\left(\Theta_M \right) \right)$ and $e_{ij}\in \Gamma\left(U_{ij},  \mathcal{A}^{0,0}\left(\mathscr{H}om_{\mathcal{O}_M}\left( \mathcal{N}_{\mathcal{F}_0}^*, \mathcal{N}_{\mathcal{F}_0}^* \right) \right) \right)$ defined by
\begin{align*}
e_{ij}:\Gamma\left( U_{ij},  \mathcal{N}_{\mathcal{F}_0}^* \right) &\to \Gamma\left( U_{ij}, \mathcal{A}^{0,0}\left( \mathcal{N}_{\mathcal{F}_0}^* \right)\right)\\
 w_{0i}^\alpha &\mapsto e_{ij}\left(w_{0i}^\alpha \right)= e_{ij}^\alpha
\end{align*}
satisfying $d_{ij}-d_{ik}+ d_{jk}= a_{ijk}$ and $e_{ij}\left( w_{0i}^\alpha\right) - e_{ik}\left(w_{0i}^\alpha\right) + e_{jk}\left(w_{0i}^\alpha\right) + \sum_{\beta=1}^q \left[d_{jk}, h_{0ij}^{\alpha\beta} \right]w_{0j}^\beta= g_{ijk}\left(w_{0i}^\alpha \right)$. Since $\delta\left( \left\{\left(\bar{\partial} d_{jk} , \bar{\partial} e_{jk} \right) \right\}\right)=0$, there exists $\left\{\left( p_j, u_j \right) \right\} \in C^0\left( \mathcal{U}, \mathcal{A}^{0,1}\left(\mathcal{E}_{\mathcal{N}_{\mathcal{F}_0}^*}\right) \right)$ with $p_j \in \Gamma\left( U_j,  \mathcal{A}^{0,1}\left(\Theta_M \right) \right)$  and $u_j\in \Gamma\left( U_j, \mathcal{A}^{0,1}\left( \mathscr{H}om_{\mathcal{O}_M}\left( \mathcal{N}_{\mathcal{F}_0}^* , \mathcal{N}_{\mathcal{F}_0}^* \right) \right) \right)$ defined by
\begin{align*}
u_j:\Gamma\left( U_j, \mathcal{N}_{\mathcal{F}_0}^* \right) &\to \Gamma\left(U_j, \mathcal{A}^{0,1}\left(\mathcal{N}_{\mathcal{F}_0}^* \right) \right)\\
 w_{0j}^\alpha &\mapsto u_j^\alpha
\end{align*}
 such that $p_k- p_j=\bar{\partial} d_{jk}$, and $u_k\left( w_{0j}^\alpha \right)+ \sum_{\beta=1}^q \left[ p_k, h_{0jk}^{\alpha\beta} \right] w_{0k}^\beta - u_j\left( w_{0j}^\alpha \right)= \bar{\partial} e_{jk}^\alpha$. Then $\bar{\partial} p_j- \bar{\partial} p_k=0$ and $\bar{\partial} u_k\left( w_{0j}^\alpha \right)+ \sum_{\beta=1}^q \left[ \bar{\partial} p_k, h_{0jk}^{\alpha\beta} \right] w_{0k}^\beta - \bar{\partial}u_j\left(w_{0j}^\alpha \right)= 0$. This implies that $\left\{\left(\bar{\partial} p_j, \bar{\partial} u_j \right)\right\}\in A^{0,2}\left( M, \mathcal{E}_{\mathcal{N}_{\mathcal{F}_0}^* } \right)$.

On the second floor, we note that
\begin{align*}
&\mathcal{L}_{d_{ij}}\left(w_{0i}^\alpha\right) - e_{ij}\left(w_{0i}^\alpha\right) - \mathcal{L}_{d_{ik}}\left(w_{0i}^\alpha\right) + e_{ik}\left(w_{0i}^\alpha\right) + \sum_{\beta=1}^q h_{0ij}^{\alpha\beta}\left(\mathcal{L}_{d_{jk}}\left(w_{0j}^\beta\right)- e_{jk}\left(w_{0j}^\beta\right)\right) -\left( b_{ij}^\alpha - b_{ik}^\alpha + \sum_{\beta=1}^q h_{0ij}^{\alpha\beta} b_{jk}^\beta \right)\\
&= \mathcal{L}_{a_{ijk}}(w_{0i}^\alpha) - \sum_{\beta=1}^q \left[ d_{jk}, h_{0ij}^{\alpha\beta}\right] w_{0j}^\beta - g_{ijk}\left(w_{0i}^\alpha\right) + \sum_{\beta=1}^q \left[ d_{ijk}, h_{0ij}^{\alpha\beta} \right] w_{0j}^\beta - \left( b_{ij}^\alpha- b_{ik}^\alpha + \sum_{\beta=1}^q h_{0ij}^{\alpha\beta} b_{jk}^\beta \right) = 0
\end{align*}

This implies that there exists $f_j\in \Gamma\left( U_j, \mathcal{A}^{0,0} \left( \mathscr{H}om_{\mathcal{O}_M}\left( \mathcal{N}_{\mathcal{F}_0}^*, \Omega_M^1 \right) \right) \right)$ defined by
\begin{align*}
f_j:\Gamma\left( U_j, \mathcal{N}_{\mathcal{F}_0}^* \right)&\to \Gamma\left( U_j,  \mathcal{A}^{0,0}\left( \Omega_M^1 \right) \right)\\
 w_{0j}^\alpha &\mapsto f_j^\alpha
\end{align*}
such that $f_j\left( w_{0j}^\alpha \right)- f_k\left( w_{0j}^\alpha \right)= b_{jk}^\alpha-\mathcal{L}_{d_{jk}}\left(w_{0j}^\alpha \right)+ e_{jk}^\alpha$. Then $\bar{\partial}f_j\left(w_{0j}^\alpha \right)-\bar{\partial} f_k\left(w_{0j}^\alpha \right)= - \mathcal{L}_{\bar{\partial} d_{jk}}\left(w_{0j}^\alpha \right) + \bar{\partial} e_{jk}\left(w_{0j}^\alpha \right)$ and
\begin{align*}
&\sum_{\beta=1}^q h_{0jk}^{\alpha\beta} \left( \mathcal{L}_{p_k}\left( w_{0k}^\beta \right) - u_k^\beta - \bar{\partial} f_k^\beta \right)- \left( \mathcal{L}_{p_j}\left(w_{0j}^\alpha \right) - u_j^\alpha - \bar{\partial} f_j^\alpha \right)\\
&=\mathcal{L}_{p_k}\left( w_{0j}^\alpha \right) - \sum_{\beta=1}^q \left[ p_k, h_{0jk}^{\alpha\beta} \right] w_{0k}^\beta  - u_k\left(w_{0j}^\alpha \right) - \bar{\partial} f_k\left(w_{0j}^\alpha \right) -  \mathcal{L}_{p_j}\left(w_{0j}^\alpha \right) + u_j^\alpha + \bar{\partial} f_j^\alpha\\
&=\mathcal{L}_{\bar{\partial} d_{jk}}\left( w_{0j}^\alpha \right) - \bar{\partial} e_{jk}\left(w_{0j}^\alpha \right) - \mathcal{L}_{\bar{\partial} d_{jk}}\left(w_{0j}^\alpha \right) + \bar{\partial} e_{jk}\left(w_{0j}^\alpha \right) = 0 
\end{align*}
Let us define $Q_j\in \Gamma\left( U_j, \mathcal{A}^{0,1} \left( \mathscr{H}om_{\mathcal{O}_M}\left( \mathcal{N}_{\mathcal{F}_0}^*, \Omega_M^1 \right) \right) \right)$ by
\begin{align*}
Q_j:\Gamma\left( U_j, \mathcal{N}_{\mathcal{F}_0}^* \right) &\to \Gamma\left( U_j, \mathcal{A}^{0,1}\left(\Omega_M^1 \right) \right)\\
T_{0j}^\alpha &\mapsto  - \mathcal{L}_{p_j}\left( w_{0j}^\alpha \right) + u_j^\alpha + \bar{\partial} f_j^\alpha
\end{align*}
Then we see that $Q=\{Q_j\} \in A^{0,1}\left( M, \mathscr{H}om_{\mathcal{O}_M}\left( \mathcal{N}_{\mathcal{F}_0}^*, \Omega_M^1 \right) \right)$.

On the third floor, for $x\in U_{jk}-S$, choose a neighborhood $U_x\subset U_{jk}-S$ such that $dw_{0j}^\alpha - \sum_{\beta=1}^q a_{j_x}^{\alpha\beta} \wedge w_{0j}^\beta, \alpha=1,...q$ for some $a_{j_x}^{\alpha\beta}\in \Gamma\left( U_x, \Omega_M^1 \right)$. Then we have 
\begin{align*}
dw_j^\alpha &= \sum_{\beta=1}^q a_{j_x}^{\alpha\beta} \wedge w_j^\beta= \sum_{\beta, \eta=1}^q a_{j_x}^{\alpha\beta} h_{0jk}^{\beta \eta} \wedge w_{0k}^\eta\\
dw_j^\alpha& = d\left(\sum_{\beta=1}^q h_{0jk}^{\alpha\beta} w_k^\beta \right) = \sum_{\beta=1}^q d h_{0jk}^{\alpha\beta}\wedge w_k^\beta + \sum_{\beta=1}^q h_{0jk}^{\alpha\beta}\left( \sum_{\eta=1}^q a_{k_x}^{\beta \eta }\wedge w_k^\eta  \right)
\end{align*}
This implies that
\begin{align*}
\sum_{\beta=1}^q a_{j_x}^{\alpha\beta} h_{0jk}^{\beta \eta} -  d h_{0jk}^{\alpha \eta}- \sum_{\beta=1}^q h_{0jk}^{\alpha\beta} a_{k_x}^{\beta \eta}\in \Gamma\left(U_{x}, \mathcal{N}_{\mathcal{F}_0}^* \right)
\end{align*}
Then we have
{\small{\begin{align*}
&c_j^\alpha + w_{0j}^1\wedge \cdots \wedge w_{0j}^q\wedge \left( df_j^\alpha - \sum_{\beta=1}^q a_{j_x}^{\alpha\beta} \wedge f_j^\beta    \right) -\sum_{\beta=1}^q h_{0jk}^{\alpha\beta} \det\left(h_{0jk}\right)\left( c_k^\beta  + w_{0k}^1\wedge \cdots \wedge w_{0k}^q\wedge \left(df_k^\beta - \sum_{\eta=1}^q a_{k_x}^{\beta \eta} \wedge f_k^\eta \right)      \right)\\
&= c_j^\alpha  - \sum_{\beta=1}^q h_{0jk}^{\alpha\beta} \det\left(h_{0jk}\right) c_k^\beta\\
&+\det\left(h_{0jk}\right) w_{0k}^1 \wedge \cdots \wedge w_{0k}^q \wedge \left(  \partial \left(\sum_{\beta=1}^q h_{0jk}^{\alpha\beta} f_k^\beta - \mathcal{L}_{d_{jk}}(w_{0j}^\alpha)+ e_{jk}^\alpha + b_{jk}^\alpha \right) - \sum_{\beta=1}^q a_{j_x}^{\alpha\beta} \wedge \left( \sum_{\eta=1}^q h_{0jk}^{\beta \eta} f_k^\eta  - \mathcal{L}_{d_{jk}}\left(w_{0j}^\beta\right)+ e_{jk}^\beta + b_{jk}^\beta \right) \right)\\
& -\sum_{\beta=1}^q h_{0jk}^{\alpha\beta} \det(h_{0jk}) w_{0k}^1\wedge \cdots \wedge w_{0k}^q\wedge \left(\partial f_k^\beta - \sum_{\eta=1}^q a_{k_x}^{\beta \eta} \wedge f_k^\eta \right)      \\
&= c_j^\alpha - \sum_{\beta=1}^q h_{0jk}^{\alpha\beta} \det\left( h_{0jk}\right) c_k^\beta+ w_{0j}^1\wedge \cdots \wedge w_{0j}^q \wedge \left( \partial b_{jk}^\alpha - \sum_{\beta=1}^q a_{j_x}^{\alpha\beta} \wedge b_{jk}^\beta    \right) = 0
\end{align*}}}

Let us define $R_j\in \Gamma\left( U_j- S,  \mathcal{A}^{0,0}\left(\mathscr{H}om_{\mathcal{O}_M}\left( \mathcal{N}_{\mathcal{F}_0}^*, \tilde{\mathcal{S}}^2 \right) \right) \right)$ by
\begin{align*}
R_j:\Gamma\left( U_j-S, \mathcal{N}_{\mathcal{F}_0}^*  \right) &\to \Gamma\left( U_j- S, \tilde{\mathcal{S}}^2 \right) \\
                                     w_{0j}^\alpha &\mapsto  -c_j^\alpha - \sum_{\beta=1}^q w_{0i}^1\wedge \cdots \wedge f_j^\beta \wedge \cdots \wedge w_{0j}^q \wedge dw_{0j}^\alpha - w_{0j}^1\wedge \cdots \wedge w_{0j}^q \wedge \partial f_j^\alpha
\end{align*}
Then $R:=\{R_j\}$ defines a global section in $A^{0,0}\left( M-S, \mathscr{H}om_{\mathcal{O}_M}\left(  \mathcal{N}_{\mathcal{F}_0}^*,  \tilde{\mathcal{S}}^2 \right) \right)$. Then we claim that
\begin{align*}
 \left( R, M, \left\{\left(\bar{\partial} p_j, \bar{\partial} u_j \right) \right\} \right) \in A^{0,0}\left(M-S , \mathscr{H}om_{\mathcal{O}_M}\left( \mathcal{N}_{\mathcal{F}_0}^*, \tilde{\mathcal{S}}^2 \right) \right) \bigoplus  A^{0,1}\left( M, \mathscr{H}om_{\mathcal{O}_M}\left( \mathcal{N}_{\mathcal{F}_0}^*, \Omega_M^1 \right)\right) \bigoplus A^{0,2}\left( M, \mathcal{E}_{\mathcal{N}_{\mathcal{F}_0}^*} \right)
\end{align*}
defines a $2$-cocycle. In fact, it is clear that $\hat{E}_2|_{M-S}\left( R_j\right)=0$ and $\bar{\partial}\left( \left(\bar{\partial} p_j, \bar{\partial} u_j \right) \right)=0$. On the other hand, on $U_x$
\begin{align*}
(-1)^q\bar{\partial} R_j\left(w_{0j}^\alpha \right)& = -(-1)^q \bar{\partial}\left( w_{0j}^1 \wedge \cdots \wedge w_{0j}^q \wedge \left( \partial f_j^\alpha - \sum_{\beta=1}^q a_{j_x}^{\alpha\beta} \wedge f_j^\beta \right) \right) = w_{0j}^1 \wedge \cdots \wedge w_{0j}^q \wedge\left( \partial \left( \bar{\partial} f_j^\alpha  \right) - \sum_{\beta=1}^q a_{j_x}^{\alpha\beta} \wedge \bar{\partial} f_j^\beta  \right) \\
&= - \hat{E}_1'\left(Q_j \right)\left(w_{0j}^\alpha \right) \Longrightarrow (-1)^q \bar{\partial} R_j + \hat{E}_1'\left( Q_j\right)
\end{align*}
and
\begin{align*}
\bar{\partial} Q_j\left(w_{0j}^\alpha \right)= \bar{\partial}\left( -\mathcal{L}_{p_j}\left(w_{0j}^\alpha \right) +u_j^\alpha + \bar{\partial} f_j^\alpha \right)= - \mathcal{L}_{\bar{\partial} p_j}\left(w_{0j}^\alpha \right) + \bar{\partial} u_j^\alpha = - \hat{E}_0'\left( \bar{\partial} p_j, \bar{\partial} u_j \right)\left(w_{0j}^\alpha \right) \Longrightarrow \bar{\partial } Q_j+ \hat{E}_0'\left(\bar{\partial} p_j, \bar{\partial} u_j \right)
\end{align*}

Then we define the correspondence of $2$-cocycles between \v Cech resolution and Dolbeault type bicomplex by 
\begin{align}\label{ii6}
\left(\left\{c_j \right\}, \left\{ b_{jk}  \right\},  \left\{ \left( a_{ijk}, g_{ijk} \right) \right\} \right)\mapsto \left(R, Q ,  \left\{ \left(\bar{\partial} p_j, \bar{\partial} u_j  \right) \right\} \right)
\end{align}
We show that the correspondence is well-defined. In other words, if $\left(\left\{c_j'\right\},\left\{b_{jk}'\right\}, \left\{ \left( a_{ijk}', g_{ijk}'\right) \right\} \right)$ is cohomologous to $\left(\left\{c_j\right\},\left\{b_{jk}\right\}, \left\{ \left( a_{ijk}, g_{ijk}\right) \right\} \right)$, then the the associated $2$-cocycle $\left( R' ,  Q', \left\{\left( \bar{\partial} p_j', \bar{\partial} u_j' \right) \right\} \right)$ is cohomologous to $\left( R, Q, \left\{ \bar{\partial} p_j, \bar{\partial} u_j \right\} \right)$. In fact, there exists $\left(\left\{ Z_j \right\} ,\left\{\left( X_{jk}, Y_{jk} \right) \right\} \right)\in C^0\left( \mathcal{U}, \mathscr{H}om_{\mathcal{O}_M}\left( \mathcal{N}_{\mathcal{F}_0}^*, \Omega_M^1 \right) \right) \bigoplus C^1\left(\mathcal{U}, \mathcal{E}_{\mathcal{N}_{\mathcal{F}_0}^*} \right)$ such that 
\begin{align*}
c_j'^\alpha - c_j^\alpha&= \sum_{\beta=1}^q w_{0j}^1 \wedge \cdots \wedge Z_j^\beta \wedge \cdots \wedge w_{0j}^q \wedge dw_{0j}^\alpha + w_{0j}^1 \wedge \cdots \wedge w_{0j}^q \wedge \partial Z_j^\alpha\\
b_{jk}'^\alpha- b_{jk}^\alpha&=\sum_{\beta=1}^q h_{0jk}^{\alpha\beta} Z_k^\beta - Z_j^\alpha + \mathcal{L}_{X_{jk}}\left(w_{0j}^\alpha \right) - Y_{jk}\left( w_{0j}^\alpha \right)\\
a_{ijk}'- a_{ijk}&= X_{ij}- X_{ik} + X_{jk}\\
g_{ijk}'\left(w_{0i}^\alpha \right) - g_{ijk}\left(w_{0i}^\alpha \right) &= Y_{ij}^\alpha - Y_{ik}^\alpha + \sum_{\beta=1}^q h_{0ij}^{\alpha\beta} Y_{jk}^\beta + \sum_{\beta=1}^q \left[ X_{jk}, h_{0ij}^{\alpha\beta} \right] w_{0j}^\beta
\end{align*}
Then $\left( d_{ij}'-d_{ij}\right)-\left( d_{ik}'-d_{ik}\right)+\left( d_{jk}'- d_{jk}\right)= a_{ijk}'- a_{ijk}= X_{ij}- X_{ik}+ X_{jk}$, so that there exists $\left\{U_j\right\}$ such that $U_k- U_j= d_{jk}'-d_{jk}- X_{jk}$. Then $\bar{\partial} U_k- \bar{\partial} U_j= \bar{\partial} d_{jk}'- \bar{\partial} d_{jk}=\left(p_k' - p_k \right) -\left(p_j'- p_j \right)$. Then $\left\{ -p_j' +p_j + \bar{\partial} U_j \right\}\in A^{0,1}\left(M, \Theta_M \right)$ and $-\bar{\partial}\left( - p_j' + p_j + \bar{\partial} U_j \right)= \bar{\partial} p_j'- \bar{\partial} p_j$. We note that $\left( e_{ij}'^\alpha - e_{ij}^\alpha \right)- \left( e_{ik}'^\alpha - e_{ik}^\alpha \right) + \sum_{\beta=1}^q h_{0ij}^{\alpha\beta}\left(e_{jk}'^\beta - e_{jk}^\beta \right) + \sum_{\beta=1}^q \left[d_{jk}' - d_{jk}, h_{0ij}^{\alpha\beta} \right] w_{0j}^\beta = g_{ijk}'\left(w_{0i}^\alpha \right)- g_{ijk}\left( w_{0i}^\alpha \right)= Y_{ij}^\alpha - Y_{ik}^\alpha + \sum_{\beta=1}^q h_{0ij}^{\alpha\beta} Y_k^\beta + \sum_{\beta=1}^q \left[ X_{jk}, h_{0ij}^{\alpha\beta} \right] w_{0j}^\beta$. Then there exists $\left\{ W_j\right\}\in C^0\left( \mathcal{U},\mathcal{A}^{0,0}\left( \mathscr{H}om_{\mathcal{O}_M}\left(\mathcal{N}_{\mathcal{F}_0}^*, \mathcal{N}_{\mathcal{F}_0}^* \right) \right) \right)$ such that $\sum_{\beta=1}^q h_{0jk}^{\alpha\beta} W_k^\beta + \sum_{\beta=1}^q \left[ U_k, h_{0jk}^{\alpha\beta} \right]w_{0k}^\beta - W_j^\alpha = e_{jk}'^\alpha- e_{jk}^\alpha- Y_{jk}^\alpha $. Then we have
\begin{align*}
&\sum_{\beta=1}^q h_{0jk}^{\alpha\beta} \bar{\partial} W_k^\beta + \sum_{\beta=1}^q \left[ \bar{\partial} U_k, h_{0jk}^{\alpha\beta} \right] w_{0k}^\beta - \bar{\partial} W_j^\alpha = \bar{\partial} e_{jk}'^\alpha - \bar{\partial} e_{jk}^\alpha\\
&=\left( u_k'\left( w_{0j}^\alpha \right)+ \sum_{\beta=1}^q \left[ p_k', h_{0jk}^{\alpha\beta} \right] w_{0k}^\beta - u_j'\left( w_{0j}^\alpha \right) \right)-\left( u_k\left( w_{0j}^\alpha \right)+ \sum_{\beta=1}^q \left[ p_k, h_{0jk}^{\alpha\beta} \right] w_{0k}^\beta - u_j\left( w_{0j}^\alpha \right)  \right)
\end{align*}
This implies that $\left( \left\{ \left(- p_j'+ p_j+ \bar{\partial} U_j, - u_j'+ u_j+\bar{\partial} W_j \right) \right\} \right) \in A^{0,1}\left(M, \mathcal{E}_{\mathcal{N}_{\mathcal{F}_0}^*} \right)$. On the other hand, we note that
{\small{\begin{align*}
&\sum_{\beta=1}^q h_{0jk}^{\alpha\beta} \left(f_k'^\beta- f_k^\beta - \mathcal{L}_{ U_k}\left( w_{0k}^\beta \right)+ W_k\left( w_{0k}^\beta \right)+ Z_k\left( w_{0k}^\beta \right) \right) -\left( f_j'^\alpha- f_j^\alpha - \mathcal{L}_{ U_j}\left( w_{0j}^\alpha \right)+ W_j\left( w_{0j}^\alpha \right)+ Z_j\left( w_{0j}^\alpha \right)  \right)\\
&= - b_{jk}'^\alpha + \mathcal{L}_{d_{jk}'}\left( w_{0j}^\alpha \right) - e_{jk}'^\alpha + b_{jk}^\alpha- \mathcal{L}_{d_{jk}}\left( w_{0j}^\alpha \right) + e_{jk}^\alpha - \mathcal{L}_{U_k- U_j}\left( w_{0j}^\alpha \right)+ \sum_{\beta=1}^q \left[ U_k, h_{0jk}^{\alpha\beta} \right] w_{0k}^\beta + \sum_{\beta=1}^q h_{0jk}^{\alpha\beta} W_k^\beta - W_j^\alpha \\
&+ b_{jk}'^\alpha - b_{jk}^\alpha - \mathcal{L}_{X_{jk}}\left( w_{0j}^\alpha \right) + Y_{jk}^\alpha = 0
\end{align*}}}
This implies that $S=\left\{ S_j\right\} :=\left\{ w_{0j}^\alpha \mapsto f_j'^\alpha- f_j^\alpha - \mathcal{L}_{U_j}\left( w_{0j}^\alpha \right) + W_j\left( w_{0j}^\alpha \right) + Z_j\left( w_{0j}^\alpha \right)\right\}$ defines a global section in $A^{0,0}\left( M, \mathscr{H}om_{\mathcal{O}_M}\left( \mathcal{N}_{\mathcal{F}_0}^*, \Omega_M^1 \right) \right)$. Then we claim that $\left( -S, \left\{ \left( -p_j'+ p_j+ \bar{\partial} U_j, - u_j'+ u_j+ \bar{\partial} W_j \right) \right\}\right)$ is mapped to $\left( R', Q', \left\{\left( \bar{\partial}p_j', \bar{\partial} u_j'  \right) \right\}\right)-\left( R, Q ,\left\{ \left( \bar{\partial} p_j,  \bar{\partial} u_j \right) \right\}\right)$. In fact,
\begin{align*}
&-\bar{\partial}\left(- p_j'+ p_j+ \bar{\partial} U_j, - u_j'+ u_j+ \bar{\partial} W_j \right) =\left( \bar{\partial} p_j' - \bar{\partial} p_j  , \bar{\partial} u_j' - \bar{\partial} u_j\right)\\
&-\bar{\partial}\left( -f_j'^\alpha+ f_j^\alpha + \mathcal{L}_{U_j}\left( w_{0j}^\alpha \right)- W_j\left(w_{0j}^\alpha \right) - Z_j\left(w_{0j}^\alpha \right)\right)    + \mathcal{L}_{-p_j'+ p_j+ \bar{\partial} U_j}\left( w_{0j}^\alpha \right) -\left(-u_j'+ u_j+ \bar{\partial} W_j \right)\left(w_{0j}^\alpha \right)\\
&=\bar{\partial} f_j'^\alpha- \bar{\partial } f_j^\alpha - \mathcal{L}_{p_j'}\left(w_{0j}^\alpha \right) + \mathcal{L}_{p_j}\left(w_{0j}^\alpha \right) + u_j'\left( w_{0j}^\alpha\right)- u_j\left(w_{0j}^\alpha \right) = \left( Q_j'- Q_j \right)\left( w_{0j}^\alpha \right)\\
&\sum_{\beta=1}^q w_{0j}^1 \wedge \cdots \wedge \left(- f_j'^\beta+ f_j^\alpha+\mathcal{L}_{U_j}\left(w_{0j}^\beta \right) - W_j\left(w_{0j}^\beta \right)- Z_j\left( w_{0j}^\beta \right) \right)\wedge \cdots \wedge w_{0j}^q \wedge dw_{0j}^\alpha \\
&+ w_{0j}^1\wedge \cdots \wedge w_{0j}^q\wedge \partial\left( - f_j'^\alpha + f_j^\alpha+ \mathcal{L}_{U_j}\left( w_{0j}^\alpha \right)- W_j\left( w_{0j}^\alpha \right)- Z_j\left( w_{0j}^\alpha \right)    \right) = \left(R'-R\right)\left( w_{0j}^\alpha \right)
\end{align*}

Conversely, we define the inverse map. Let $\left(R, Q , \left( \phi, \left\{ N_j \right\} \right)      \right)$ be a $2$-cocycle in the Dolbeault type bicomplex. Since $\bar{\partial} \left(\phi, N_j\right)=0$, there exists $\left(p_j, u_j \right) \in \Gamma\left( U_j, \mathcal{A}^{0,1}\left(\mathcal{E}_{\mathcal{N}_{\mathcal{F}_0}^*} \right) \right)$ such that $\bar{\partial} \left( p_j, u_j \right)= \left(\phi, N_j \right)$. Then since $\bar{\partial}\left(\left(p_k, u_k \right) - \left( p_j, u_j \right) \right)=0$, there exists $\left(d_{jk}, e_{jk} \right)\in \Gamma\left( U_{jk}, \mathcal{A}^{0,0}\left( \mathcal{E}_{\mathcal{N}_{\mathcal{F}_0}^*} \right) \right)$ such that $\bar{\partial}\left( d_{jk}, e_{jk} \right)=\left( p_k, u_k \right)-\left(p_j, u_j \right)$.  Then we set
$\left(a_{ijk}, g_{ijk} \right):=\left(d_{ij}, e_{ij} \right)- \left(d_{ik}, e_{ik} \right) + \left( d_{jk}, e_{jk}\right)$. Then we have $\left\{ \left( a_{ijk}, g_{ijk}\right)\right\}\in C^2\left(\mathcal{U}, \mathcal{E}_{\mathcal{N}_{\mathcal{F}_0}^*} \right)$. On the other hand,  since $\bar{\partial} Q \left(w_{0j}^\alpha \right)+ \mathcal{L}_{\phi}\left(w_{0j}^\alpha \right) - N_j\left(w_{0j}^\alpha \right)=0$, we have
\begin{align*}
\bar{\partial} Q_j\left(w_{0j}^\alpha \right) + \mathcal{L}_{\bar{\partial}p_j}\left(w_{0j}^\alpha \right) - \bar{\partial} u_j\left(w_{0j}^\alpha \right)=0 \Longrightarrow \bar{\partial}\left(Q_j\left(w_{0j}^\alpha \right)  + \mathcal{L}_{p_j}\left( w_{0j}^\alpha \right) - u_j\left(w_{0j}^\alpha \right) \right) = 0
\end{align*}
Then there exist $f_j\in \Gamma\left( U_j, \mathcal{A}^{0,0}\left( \mathscr{H}om_{\mathcal{O}_M}\left( \mathcal{N}_{\mathcal{F}_0}^*, \Omega_M^1 \right) \right) \right) $ such that $\bar{\partial} f_j\left( w_{0j}^\alpha \right)= Q_j^\alpha + \mathcal{L}_{p_j}\left(w_{0j}^\alpha \right) - u_j\left(w_{0j}^\alpha \right)$. Then we note that
\begin{align*}
&\bar{\partial} f_j\left( w_{0j}^\alpha \right)- \bar{\partial} f_k\left(w_{0j}^\alpha \right) =  \mathcal{L}_{p_j}\left(w_{0j}^\alpha \right) - u_j\left(w_{0j}^\alpha \right) - \sum_{\beta=1}^q h_{0jk}^{\alpha\beta} \left(\mathcal{L}_{p_k}\left( w_{0k}^\beta \right)  - u_k\left(w_{0k}^\beta \right)   \right)\\
&= \mathcal{L}_{p_j- p_k}\left(w_{0j}^\alpha \right) - u_j\left(w_{0j}^\alpha \right) + \sum_{\beta=1}^q \left[ p_k, h_{0jk}^{\alpha\beta} \right] w_{0k}^\beta + u_k\left( w_{0j}^\alpha \right) = \mathcal{L}_{-\bar{\partial}d_{jk}}\left(w_{0j}^\alpha \right) + \bar{\partial}e_{jk}\left( w_{0j}^\alpha \right)
\end{align*}
Then we define 
\begin{align*}
b_{jk}:\Gamma\left( U_{jk}, \mathcal{N}_{\mathcal{F}_0}^* \right) &\to \Gamma\left( U_{jk}, \Omega_M^1 \right) \\
 w_{0j}^\alpha &\mapsto  f_j\left(w_{0j}^\alpha \right) - f_k\left( w_{0j}^\alpha \right) + \mathcal{L}_{d_{jk}}\left( w_{0j}^\alpha \right) - e_{jk}\left(w_{0j}^\alpha \right)
\end{align*}
Then  $\left\{ b_{jk}\right\} \in C^1 \left( \mathcal{U}, \mathscr{H}om_{\mathcal{O}_M}\left( \mathcal{N}_{\mathcal{F}_0}^* , \Omega_M^1 \right) \right)$. 

For $x\in U_j-S$, choose neighborhood $U_x\subset U_j-S$ such that $dw_{0j}^\alpha - \sum_{\beta=1}^q a_{j_x}^{\alpha\beta} \wedge w_{0j}^\beta$ for some $a_{j_x}^{\alpha\beta}\in \Gamma\left( U_x, \Omega_M^1 \right),\alpha=1,...,q$. We can write $R\left(w_{0j}^\alpha \right) = w_{0j}^1 \wedge \cdots \wedge w_{0j}^q \wedge R_{j_x}^\alpha$ for some $R_{j_x}\in \Gamma\left( U_x, \mathcal{A}^{0,0}\left( \bigwedge^2 \Omega_M^1 \right) \right)$. Then since $w_{0j}^1 \wedge \cdots \wedge w_{0j}^q \wedge \left(\bar{\partial} R_{j_x}^\alpha - \partial Q_j^\alpha + \sum_{\beta=1}^q a_{j_x}^{\alpha\beta} \wedge Q_j^\beta \right) =0 $. Then we have
\begin{align*}
&w_{0j}^1 \wedge \cdots \wedge w_{0j}^q \wedge \left( \bar{\partial} R_{j_x}^\alpha - \partial \left( \bar{\partial} f_j^\alpha - \mathcal{L}_{p_j}\left(w_{0j}^\alpha \right) + u_j^\alpha \right) + \sum_{\beta=1}^q a_{j_x}^{\alpha\beta} \wedge \left( \bar{\partial} f_j^\beta - \mathcal{L}_{p_j}\left(w_{0j}^\beta \right) + u_j^\beta\right) \right) \\
&= w_{0j}^1 \wedge \cdots \wedge w_{0j}^q \wedge \left(\bar{\partial} R_{j_x}^\alpha +\bar{\partial} \left( \partial f_j^\alpha - \sum_{\beta=1}^q a_{j_x}^{\alpha\beta} \wedge f_j^\beta     \right) \right)  = 0
\end{align*}

Then
\begin{align*}
c_j^\alpha:=-R_j\left(w_{0j}^\alpha \right) - \sum_{\beta=1}^q w_{0j}^1 \wedge \cdots \wedge f_j^\beta \wedge \cdots \wedge w_{0j}^q \wedge dw_{0j}^\alpha - w_{0j}^1 \wedge \cdots \wedge w_{0j}^q \wedge \partial f_j^\alpha \in \Gamma\left( U_j-S, \tilde{\mathcal{S}}^2 \right)
\end{align*}
Hence $c_j^\alpha \in \Gamma\left( U_j, \tilde{\mathcal{S}}^2 \right)$. We define $c_j\in \Gamma\left( U_j, \mathscr{H}om_{\mathcal{O}_M}\left( \mathcal{N}_{\mathcal{F}_0}^*, \tilde{\mathcal{S}}^2 \right) \right)$ by 
\begin{align*}
c_j:\Gamma\left(U_j, \mathcal{N}_{\mathcal{F}_0}^* \right) &\to \Gamma\left( U_j,  \tilde{\mathcal{S}}^2 \right)\\
w_{0j}^\alpha &\mapsto c_j^\alpha
\end{align*}
We claim that $\left(\left\{ c_j \right\}, \left\{  b_{jk} \right\}, \left\{\left( a_{ijk}, g_{ijk} \right)  \right\} \right)$ defines a $2$-cocycle in the \v Cech resolution of $\mathcal{E}_{\mathcal{N}_{\mathcal{F}_0}^*}^\bullet$. In fact, it is clear that $E_2\left( c_j \right) =0$ and $\delta\left(\left\{\left(a_{ijk},  g_{ijk} \right) \right\}\right)=0$. On the other hand, we have
{\tiny{\begin{align*}
&b_{ij}^\alpha - b_{ik}^\alpha + \sum_{\beta=1}^q h_{0ij}^{\alpha\beta} b_{jk}^\beta  \\
&=  f_i\left(w_{0i}^\alpha \right) - f_j\left( w_{0i}^\alpha \right) + \mathcal{L}_{d_{ij}}\left( w_{0i}^\alpha \right) - e_{ij}\left( w_{0i}^\alpha \right)- f_i\left(w_{0i}^\alpha \right)+ f_k\left(w_{0i}^\alpha \right) - \mathcal{L}_{d_{ik}}\left(w_{0i}^\alpha \right) + e_{ik}\left(w_{0i}^\alpha \right) + \sum_{\beta=1}^q h_{0ij}^{\alpha\beta} \left(f_j\left(w_{0j}^\beta \right)- f_k\left(w_{0j}^\beta \right) + \mathcal{L}_{d_{jk}}\left(w_{0j}^\beta \right) - e_{jk}\left(w_{0j}^\beta \right)\right)\\
&=\left( \mathcal{L}_{d_{ij}}\left( w_{0i}^\alpha\right) - e_{ij}\left(w_{0i}^\alpha \right) \right) - \left( \mathcal{L}_{d_{ik}}\left(w_{0i}^\alpha \right)- e_{ik}\left(w_{0i}^\alpha \right)\right) + \sum_{\beta=1}^q h_{0ij}^{\alpha\beta} \left( \mathcal{L}_{d_{jk}}\left(w_{0j}^\beta \right)- e_{jk}\left(w_{0j}^\beta \right) \right)
\end{align*}}}
and we have
{\tiny{\begin{align*}
&c_j^\alpha - \det\left(h_{0jk} \right)\sum_{\beta=1}^q h_{0jk}^{\alpha\beta} c_k^\beta \\
&= - \sum_{\beta=1}^q w_{0j}^1 \wedge \cdots \wedge f_j^\beta \wedge\cdots \wedge w_{0j}^q \wedge dw_{0j}^\alpha - w_{0j}^1 \wedge \cdots \wedge w_{0j}^q \wedge \partial f_j^\alpha + \det\left(h_{0jk} \right)\sum_{\gamma=1}^q h_{0jk}^{\alpha\gamma}\left(  \sum_{\beta=1}^q w_{0k}^1 \wedge \cdots \wedge f_k^\beta \wedge\cdots \wedge w_{0k}^q \wedge dw_{0k}^\gamma + w_{0k}^1 \wedge \cdots \wedge w_{0k}^q \wedge \partial f_k^\gamma  \right)\\
&= - \sum_{\beta=1}^q w_{0j}^1 \wedge \cdots \wedge f_j^\beta \wedge\cdots \wedge w_{0j}^q \wedge dw_{0j}^\alpha - w_{0j}^1 \wedge \cdots \wedge w_{0j}^q \wedge \partial f_j^\alpha + \sum_{\gamma=1}^q h_{0jk}^{\alpha \gamma} \left(\sum_{\beta=1}^q w_{0j}^1 \wedge \cdots \wedge f_k\left(w_{0j}^\beta \right) \wedge \cdots \wedge w_{0j}^q \wedge d w_{0k}^\gamma + w_{0j}^1 \wedge \cdots \wedge w_{0j}^q \wedge \partial f_k^\gamma \right)\\
&= - \sum_{\beta=1}^q w_{0j}^1 \wedge \cdots \wedge \left( f_j\left(w_{0j}^\beta \right) - f_k\left(w_{0j}^\beta \right) \right) \wedge \cdots \wedge w_{0j}^q \wedge dw_{0j}^\alpha- w_{0j}^1 \wedge \cdots \wedge w_{0j}^q \wedge \partial\left( f_j\left( w_{0j}^\alpha \right) - f_k\left(w_{0j}^\alpha \right) \right)\\
& - \sum_{\beta, \gamma=1}^q w_{0j}^1 \wedge \cdots \wedge \left( \sum_{\eta=1}^q h_{0jk}^{\beta \eta} f_k^\eta \right) \wedge \cdots \wedge w_{0j}^q \wedge d h_{0jk}^{\alpha \gamma} \wedge w_{0k}^\gamma - \sum_{\beta=1}^q w_{0j}^1 \wedge \cdots \wedge w_{0j}^q \wedge d h_{0jk}^{\alpha \beta} \wedge  f_k^\beta\\
&= - \sum_{\beta=1}^q w_{0j}^1 \wedge \cdots \wedge \left( f_j\left(w_{0j}^\beta \right) - f_k\left(w_{0j}^\beta \right) \right) \wedge \cdots \wedge w_{0j}^q \wedge dw_{0j}^\alpha- w_{0j}^1 \wedge \cdots \wedge w_{0j}^q \wedge \partial\left( f_j\left( w_{0j}^\alpha \right) - f_k\left(w_{0j}^\alpha \right) \right)\\
& - \sum_{\beta, \gamma=1}^q w_{0j}^1 \wedge \cdots \wedge \left(  h_{0jk}^{\beta \gamma} f_k^\gamma \right) \wedge \cdots \wedge w_{0j}^q \wedge d h_{0jk}^{\alpha \gamma} \wedge w_{0k}^\gamma - \sum_{\beta=1}^q w_{0j}^1 \wedge \cdots \wedge w_{0j}^q \wedge d h_{0jk}^{\alpha \beta} \wedge  f_k^\beta\\
&= - \sum_{\beta=1}^q w_{0j}^1 \wedge \cdots \wedge \left( f_j\left(w_{0j}^\beta \right) - f_k\left(w_{0j}^\beta \right) \right) \wedge \cdots \wedge w_{0j}^q \wedge dw_{0j}^\alpha- w_{0j}^1 \wedge \cdots \wedge w_{0j}^q \wedge \partial\left( f_j\left( w_{0j}^\alpha \right) - f_k\left(w_{0j}^\alpha \right) \right)\\
\end{align*}}}
This prove the claim. Then we define the inverse correspondence by 
\begin{align}
\left( R, Q, \left(\left\{ \phi, N_j \right\} \right) \right)\mapsto \left(\left\{c_j\right\}, \left\{ b_{jk}\right\}, \left\{\left( a_{ijk}, g_{ijk} \right) \right\} \right)
\end{align}
We show that the correspondence is well-defined. In other words, if $\left(R', Q', \left\{\left( \phi', N_j'\right) \right\} \right)$ is cohomologous to $\left( R, Q, \left\{ \left( \phi, N_j \right) \right\} \right)$, then the associated $2$-cocycle $\left(\left\{ c_j' \right\}, \left\{ b_{jk}' \right\}, \left\{ \left( a_{ijk}', g_{ijk}' \right) \right\} \right)$ is cohomologous to $\left( \left\{ c_j \right\}, \left\{ b_{jk} \right\}, \left\{ \left( a_{ijk}, g_{ijk} \right) \right\} \right)$. In fact, there exist $\left( Z, \left( X, \left\{ Y_j \right\} \right)\right) \in A^{0,0}\left( M, \mathscr{H}om_{\mathcal{O}_M}\left( \mathcal{N}_{\mathcal{F}_0}^*, \Omega_M^1 \right) \right) \bigoplus A^{0,1}\left( M, \mathcal{E}_{\mathcal{N}_{\mathcal{F}_0}^*} \right) $ such that
\begin{align*}
R' \left(w_{0j}^\alpha \right)- R\left( w_{0j}^\alpha \right)&= \sum_{\beta=1}^q w_{0j}^1 \wedge \cdots \wedge Z\left( w_{0j}^\beta \right) \wedge \cdots \wedge w_{0j}^q \wedge d w_{0j}^\alpha + w_{0j}^1\wedge \cdots \wedge w_{0j}^\alpha \wedge \partial Z\left( w_{0j}^\alpha \right)  \\
Q_j'^\alpha - Q_j^\alpha&= - \bar{\partial} Z\left(w_{0j}^\alpha \right) + \mathcal{L}_{X}\left( w_{0j}^\alpha \right) - Y_j\left( w_{0j}^\alpha \right) \\
\phi'- \phi&= - \bar{\partial} X\\
N_j'- N_j &= - \bar{\partial} Y_j
\end{align*}
Then we have $\bar{\partial}\left(X + p_j' - p_j , Y_j + u_j'- u_j \right)=0$, so that there exists $\left( U_j, W_j\right)\in \Gamma\left(U_j, \mathcal{A}^{0,0}\left( \mathcal{E}_{\mathcal{N}_{\mathcal{F}_0}^*} \right) \right)$ such that $\bar{\partial} U_j=X+ p_j'- p_j$ and $\bar{\partial} W_j= Y_j+ u_j'- u_j$. Then $\bar{\partial}\left(\left( U_j, W_j\right)- \left( U_k, W_k\right) -\left( d_{jk}, e_{jk}\right) + \left( d_{jk}', e_{jk}'\right)\right)=0$, so that we have $\left\{ \left( U_j, W_j\right)- \left( U_k, W_k\right) -\left( d_{jk}, e_{jk}\right) + \left( d_{jk}', e_{jk}'\right)  \right\}\in C^1\left( \mathcal{U}, \mathcal{E}_{\mathcal{N}_{\mathcal{F}_0}^*} \right)$. On the other hand, we note that $\bar{\partial}\left( Z\left( w_{0j}^\alpha \right) - \mathcal{L}_{U_j}\left( w_{0j}^\alpha \right) + W_j\left( w_{0j}^\alpha \right) - f_j^\alpha + f_j'^\alpha \right)= Q_j^\alpha- Q_j'^\alpha -\mathcal{L}_{p_j'-p_j}\left( w_{0j}^\alpha \right) + \left( u_j'- u_j\right)\left( w_{0j}^\alpha \right)- f_j^\alpha + f_j'^\alpha=0$, so that $\left\{ Z - D_0'\left( U_j, W_j\right)- f_j+ f_j'\right\}\in C^0\left( \mathcal{U}, \mathscr{H}om_{\mathcal{O}_M}\left( \mathcal{N}_{\mathcal{F}_0}, \Omega_M^1 \right)  \right)$. 

Then we see that $\left( \left\{ -Z + D_0'\left( U_j, W_j\right)+ f_j- f_j'\right\}  , \left\{ \left( U_j, W_j\right)- \left( U_k, W_k\right) -\left( d_{jk}, e_{jk}\right) + \left( d_{jk}', e_{jk}'\right)  \right\} \right)$ is mapped to $\left(\left\{ c_j'\right\},\left\{ b_{jk}'\right\}, \left\{ \left(a_{ijk}', g_{ijk}'\right) \right\}  \right) - \left( \left\{ c_j \right\}, \left\{ b_{jk} \right\}, \left\{ \left( a_{ijk}, g_{ijk}\right) \right\} \right)$.

\begin{proposition}\label{a24}
We recall from \textnormal{Lemma \ref{oc29}} that given a $1$-cocycle $\left( \left\{ \beta_j \right\}, \left\{ \left( \theta_{jk}, B_{jk} \right)\right\} \right)$ in $\mathbb{H}^1\left(M, \mathcal{E}_{\mathcal{N}_{\mathcal{F}_0}^*}^\bullet \right)$, we have the associated obstruction class from $(\ref{it2})$
\begin{align}
\left(\left\{ S_j \right\}, \left\{ G_{jk} \right\}, \left\{\left[\left( \theta_{jk}, B_{jk}\right), \left( \theta_{jk}, B_{jk} \right) \right] \right\} \right) \in \mathbb{H}^2\left( M, \mathcal{E}_{\mathcal{N}_{\mathcal{F}_0}^*}^\bullet \right) \label{ii7}
\end{align}
from the \v Cech resolution of $\mathcal{E}_{\mathcal{N}_{\mathcal{F}_0}^*}^\bullet$. Assume that the $1$-cocycle $\left( \left\{ \beta_j \right\}, \left\{ \left( \theta_{jk}, B_{jk} \right) \right\} \right)$ from the \v Cech resolution of $\mathcal{E}_{\mathcal{N}_{\mathcal{F}_0}^*}^\bullet$ corresponds to the $1$-cocycle $\left( \Psi, \left\{ \left( \psi, k_j \right) \right\} \right)$ from the Dolbeault type bicomplex associated to $\mathcal{E}_{\mathcal{N}_{\mathcal{F}_0}^*}^\bullet$ under the isomorphism $(\ref{ii5})$ from $\mathbb{H}^1\left( M, \mathcal{E}_{\mathcal{N}_{\mathcal{F}_0}^*}^\bullet \right)$ to $\textnormal{B}^1$. We recall from \textnormal{Lemma \ref{oci85}} that we have the associated obstruction class from $(\ref{it13})$
\begin{align} \label{ii8}
\left( l, \kappa, \left\{ -\left[\left( \psi, k_j \right), \left( \psi, k_j \right) \right] \right\} \right) \in \textnormal{B}^2
\end{align}
from the Dolbeault type bicomplex associated to $\mathcal{E}_{\mathcal{N}_{\mathcal{F}_0}^*}^\bullet$. Then $(\ref{ii7})$ corresponds to $(\ref{ii8})$ under the isomorphism $(\ref{ii6})$ from $\mathbb{H}^2\left( M, \mathcal{E}_{\mathcal{N}_{\mathcal{F}_0}^*}^\bullet \right)$ to $\textnormal{B}^2$.

\end{proposition}

\begin{proof}
Let $\left(\left\{ \beta_j \right\}, \left\{\left( \theta_{jk}, B_{jk} \right) \right\} \right)$ be a $1$-cocycle in Lemma \ref{oc29}. Then $\left( \{S_j\}, \{G_{jk} \}, \left\{\left[ \left(\theta_{jk}, B_{jk}\right), \left( \theta_{jk}, B_{jk} \right) \right] \right\} \right)$ defines a $2$-cocycle by Lemma \ref{oc29}. Then we have the corresponding $1$-cocycle $\left(\Psi, \left(\psi, k_j:=\bar{\partial} B_j \right) \right)$ from Doubault type bicomplex from $(\ref{ii5})$. Since $-\delta\left( \left\{\left(\xi_j, B_j \right) \right\}\right)=\left\{\left( \theta_{jk}, B_{jk} \right)\right\}$, i.e. $\left( \xi_i, B_i\right)-\left( \xi_j , A_j \right)= \left( \theta_{jk}, A_{jk}\right)$, we have
\begin{align*}
&\left[\left( \xi_i, B_i \right), \left( \xi_j, B_j \right) \right] - \left[\left( \xi_i, B_i\right), \left( \xi_k, B_k \right) \right] + \left[\left( \xi_j, B_j \right), \left( \xi_k , B_k \right) \right]= \left[ \left(\xi_i, B_i \right), \left(\theta_{jk}, B_{jk} \right) \right] + \left[ \left(\xi_j, B_j \right), \left( \xi_k, B_k \right)- \left( \xi_j, B_j \right) \right]\\
&=\left[\left( \xi_i, B_i \right), \left( \theta_{jk}, B_{jk} \right) \right] - \left[\left( \xi_j, B_j \right), \left( \theta_{jk}, B_{jk} \right) \right] = \left[ \left( \theta_{ij}, B_{ij} \right), \left( \theta_{jk} , B_{jk} \right) \right]
\end{align*}

This implies that
\begin{align*}
\delta\left( \left\{\left[ \left(\xi_i, B_i\right), \left(\xi_j, B_j \right)\right]  \right\} \right) = \left\{ \left[ \left(\theta_{ij}, B_{ij}\right) , \left(\theta_{jk}, B_{jk} \right) \right]   \right]
\end{align*}

Then we set $\left(d_{jk}, e_{jk} \right):=\left[\left( \xi_j, B_j \right), \left( \xi_k, B_k \right) \right]$ and we note that
\begin{align*}
\bar{\partial}\left[ \left(\xi_j, B_j\right), \left(\xi_k, B_k \right)\right]  =\left[ \bar{\partial}\left(\xi_j, B_j \right), \left( \xi_k, B_k \right) \right] + \left[ \left(\xi_j, B_j \right) , \bar{\partial} \left(\xi_k, B_k \right) \right] = \left[\left(\psi, k_j \right), \left(\xi_k, B_k \right) \right] - \left[\left(\psi , k_k \right), \left( \xi_j, B_j\right) \right]
\end{align*}
Then we set $\left(p_j, u_j \right) :=  \left[ \left(\psi, k_j\right), \left(\xi_j, B_j \right) \right]$. Then $\bar{\partial}\left(p_j, u_j \right)= - \left[\left(\psi, k_j \right), \left( \psi, k_j \right) \right]$.

On the other hand, we note that
\begin{align*}
B_j\left(w_{0j}^\alpha \right)- B_k^j\left( w_{0j}^\alpha \right)=B_j\left(w_{0j}^\alpha \right) - \left( B_k\left(w_{0j}^\alpha  \right) + \sum_{\gamma=1}^q \left[ \xi_k, h_{0jk}^{\alpha \gamma } \right] w_{0k}^\gamma \right) = B_{jk}\left(w_{0j}^\alpha \right) = B_{jk}^\alpha
\end{align*}

Then we have
\begin{align}\label{b15}
B_{jk}^{\alpha \delta} = B_j^{\alpha \delta} - \sum_{\beta, \eta=1}^q h_{0jk}^{\alpha \beta} B_k^{\beta\eta} h_{0kj}^{\eta \delta} - \sum_{\gamma, \eta=1}^q \left[ \xi_k, h_{0jk}^{\alpha \gamma}\right] h_{0kj}^{\gamma \delta}
\end{align}

Since $\left( d_{jk}, e_{jk}\right)=\left[\left(\xi_j, B_j \right), \left( \xi_k ,B_k\right)    \right] =\left( \left[ \xi_j, \xi_k \right], \left[ \xi_j, B_k^j\right] - \left[ \xi_k, B_j\right] + B_j \circ B_k^j - B_k^j \circ B_j \right)$
and
\begin{align*}
B_k^j \left(w_{0j}^\alpha \right)= \sum_{\delta=1}^q  \left( \sum_{\beta, \eta=1}^q h_{0jk}^{\alpha\beta} B_k^{\beta \eta} h_{0kj}^{\eta \delta} + \sum_{\eta=1}^q \left[\xi_k, h_{0jk}^{\alpha \eta} \right]h_{0kj}^{\eta \delta} \right) w_{0j}^\delta,
\end{align*}
we see that
\begin{align}
& -\mathcal{L}_{\left[ \xi_j, \xi_k\right]}\left( w_{0j}^\alpha\right) + e_{jk}\left( w_{0j}^\alpha \right) =-\mathcal{L}_{\left[ \xi_j, \xi_k \right]}\left( w_{0j}^\alpha \right) +\left(\left[ \xi_j, B_k^j\right] - \left[ \xi_k, B_j\right] + B_j \circ B_k^j - B_k^j \circ B_j \right) \left( w_{0j}^\alpha \right)\label{b17}\\
&=-\mathcal{L}_{\xi_j}\mathcal{L}_{\xi_k}\left( w_{0j}^\alpha \right) + \mathcal{L}_{\xi_k} \mathcal{L}_{\xi_j}\left( w_{0j}^\alpha \right) + \sum_{\delta=1}^q \left[ \xi_j,  \sum_{\beta, \eta=1}^q h_{0jk}^{\alpha\beta} B_k^{\beta \eta} h_{0kj}^{\eta \delta} + \sum_{\eta=1}^q \left[\xi_k, h_{0jk}^{\alpha \eta} \right]h_{0kj}^{\eta \delta} \right] w_{0j}^\delta - \sum_{\delta=1}^q \left[ \xi_k, B_j^{\alpha \delta} \right] w_{0j}^\delta \notag \\
&+ \sum_{\xi=1}^q \left(\sum_{\beta, \eta=1}^q h_{0jk}^{\alpha\beta} B_k^{\beta \eta} h_{0kj}^{\eta \delta} + \sum_{\eta=1}^q \left[\xi_k, h_{0jk}^{\alpha \eta} \right]h_{0kj}^{\eta \delta}\right)  B_j^{\delta \xi} w_{0j}^\xi - \sum_{\delta, \xi =1}^q B_j^{\alpha \xi} \left( \sum_{\beta, \eta=1}^q h_{0jk}^{\xi\beta} B_k^{\beta \eta} h_{0kj}^{\eta \delta} + \sum_{\eta=1}^q \left[\xi_k, h_{0jk}^{\xi \eta} \right]h_{0kj}^{\eta \delta} \right) w_{0j}^\delta \notag \\
&= -\mathcal{L}_{\xi_j}\mathcal{L}_{\xi_k}\left( w_{0j}^\alpha \right) + \mathcal{L}_{\xi_k} \mathcal{L}_{\xi_j}\left( w_{0j}^\alpha \right) + \sum_{\beta, \eta=1}^q \mathcal{L}_{\xi_j} \left( h_{0jk}^{\alpha \beta} B_k^{\beta \eta} w_{0k}^\eta \right)  - \sum_{\delta, \beta, \eta=1}^q h_{0jk}^{\alpha \beta} B_k^{\beta \eta} h_{0kj}^{\eta \delta} \mathcal{L}_{\xi_j}\left( w_{0j}^\delta \right) + \sum_{\eta=1}^q \mathcal{L}_{\xi_j}\left( \left[ \xi_k, h_{0jk}^{\alpha \eta} \right] w_{0k}^\eta \right) \notag \\
&-\sum_{\eta, \delta=1}^q \left[ \xi_k, h_{0jk}^{\alpha \eta}\right] h_{0kj}^{\eta \delta}\mathcal{L}_{\xi_j}\left(w_{0j}^\delta \right)- \sum_{\delta=1}^q \left[ \xi_k, B_j^{\alpha \delta} \right] w_{0j}^\delta + \sum_{\xi, \beta, \eta=1}^q h_{0jk}^{\alpha\beta} B_k^{\beta \eta} h_{0kj}^{\eta \delta} B_j^{\delta \xi} w_{0j}^\xi + \sum_{\xi, \eta=1}^q \left[ \xi_k, h_{0jk}^{\alpha \eta} \right] h_{0kj}^{\eta \delta} B_j^{\delta \xi} w_{0j}^\xi \notag  \\
&-\sum_{\beta, \eta, \xi=1}^q B_j^{\alpha \xi} h_{0jk}^{\xi\beta} B_k^{\beta \eta} w_{0k}^\eta - \sum_{\delta, \xi, \eta=1}^q B_j^{\alpha \xi} \left[ \xi_k, h_{0jk}^{\xi\eta}\right] w_{0k}^\eta  \notag
\end{align}
On the other hand, we recall that $G_{jk}\left( w_{0j}^\alpha \right)= - \beta_j\left( B_{jk}^\alpha \right) - \beta_k \left( B_{jk}^\alpha \right) + \mathcal{L}_{\theta_{jk}}\left( \beta_j\left(w_{0j}^\alpha \right) \right) + \mathcal{L}_{\theta_{jk}}\left( \beta_k\left( w_{0j}^\alpha \right)  \right)$. Then from $(\ref{b15})$, we have
{\Small{\begin{align}
&G_{jk}\left( w_{0j}^\alpha \right)- \beta_j\left( B_{jk}^\alpha \right) - \beta_k \left( B_{jk}^\alpha \right) + \mathcal{L}_{\theta_{jk}}\left( \beta_j\left(w_{0j}^\alpha \right) \right) + \mathcal{L}_{\theta_{jk}}\left( \beta_k\left( w_{0j}^\alpha \right)  \right) \label{b16}\\
&=-\sum_{\delta=1}^q \left(B_j^{\alpha \delta} - \sum_{\beta, \eta=1}^q h_{0jk}^{\alpha\beta } B_k^{\beta\eta} h_{0kj}^{\eta \delta} - \sum_{\gamma, \eta=1}^q \left[ \xi_k, h_{0jk}^{\alpha \gamma}\right] h_{0kj}^{\gamma \delta}\right) \beta_j\left( w_{0j}^\delta \right) - \sum_{\delta=1}^q \left( B_j^{\alpha \delta} - \sum_{\beta, \eta=1}^q h_{0jk}^{\alpha \beta} B_k^{\beta \eta} h_{0kj}^{\eta \delta} - \sum_{\gamma, \eta=1}^q \left[ \xi_k, h_{0jk}^{\alpha \gamma}\right] h_{0kj}^{\gamma \delta}\right) \beta_k\left(w_{0j}^\delta \right) \notag \\
&+ \mathcal{L}_{\xi_j}\left(\beta_j \left( w_{0j}^\alpha \right) \right) - \mathcal{L}_{\xi_k}\left(\beta_j \left( w_{0j}^\alpha \right) \right) + \mathcal{L}_{\xi_j}\left(\beta_k \left( w_{0j}^\alpha \right) \right) - \mathcal{L}_{\xi_k}\left(\beta_k \left( w_{0j}^\alpha \right) \right) \notag\\
&= - \beta_j\left( B_j^\alpha \right)  + \sum_{\delta, \beta, \eta=1}^q h_{0jk}^{\alpha\beta } B_k^{\beta \eta} h_{0kj}^{\eta \delta}\left(\beta_k\left( w_{0j}^\delta \right) + \mathcal{L}_{\theta_{jk}}\left(w_{0j}^\delta \right) - B_{jk}^\delta \right) + \sum_{\delta, \gamma, \eta=1}^q \left[ \xi_k, h_{0jk}^{\alpha \gamma } \right] h_{0kj}^{\gamma \delta} \left( \beta_k\left( w_{0j}^\delta \right) + \mathcal{L}_{\theta_{jk}}\left(w_{0j}^\delta \right) - B_{jk}^\delta \right) \notag\\
&- \sum_{\delta=1}^q B_j^{\alpha \delta}\left(\beta_j\left( w_{0j}^\delta \right) - \mathcal{L}_{\theta_{jk}} \left( w_{0j}^\delta \right) + B_{jk}^\delta \right) + \sum_{\beta=1}^q h_{0jk}^{\alpha \beta}\beta_k\left(B_k^\beta \right) + \sum_{\gamma, \eta=1}^q \left[ \xi_k , h_{0jk}^{\alpha \gamma} \right] \beta_k\left( w_{0k}^\gamma \right) \notag\\
&+ \mathcal{L}_{\xi_j}\left( \beta_j\left( w_{0j}^\alpha \right) \right) - \mathcal{L}_{\xi_k}\left( \beta_k\left( w_{0j}^\alpha \right) + \mathcal{L}_{\theta_{jk}}\left(w_{0j}^\alpha \right) - B_{jk}^\alpha \right) + \mathcal{L}_{\xi_j}\left( \beta_j\left(w_{0j}^\alpha \right) - \mathcal{L}_{\theta_{jk}}\left(w_{0j}^\alpha \right) + B_{jk}^\alpha \right)  - \mathcal{L}_{\xi_k}\left(\beta_k \left( w_{0j}^\alpha \right) \right) \notag\\
&=- 2\beta_j\left( B_j^\alpha \right) + 2\sum_{\beta=1}^q h_{0jk}^{\alpha\beta} \beta_k\left( B_k^\beta \right) + \sum_{\delta, \beta, \eta=1}^q h_{0jk}^{\alpha\beta} B_k^{\beta \eta} h_{0kj}^{\eta \delta} \mathcal{L}_{\xi_j- \xi_k}\left( w_{0j}^\delta \right) \notag\\
&- \sum_{\delta, \beta, \eta,\xi=1}^q h_{0jk}^{\alpha\beta} B_k^{\beta \eta} h_{0kj}^{\eta \delta}\left(B_j^{\delta \xi} - \sum_{ \sigma,\gamma=1}^q h_{0jk}^{\delta \sigma} B_k^{\sigma \gamma }h_{0kj}^{\gamma  \xi} - \sum_{\gamma =1}^q \left[ \xi_k, h_{0jk}^{\delta \gamma}\right] h_{0kj}^{\gamma \xi } \right) w_{0j}^\xi +2 \sum_{\gamma, \eta=1}^q \left[ \xi_k, h_{0jk}^{\alpha \gamma} \right] \beta_k\left(w_{0k}^\gamma \right) \notag \\
&+ \sum_{\delta, \gamma, \eta=1}^q \left[ \xi_k, h_{0jk}^{\alpha \gamma}\right] h_{0kj}^{\gamma \delta} \mathcal{L}_{\xi_j-\xi_k}\left( w_{0j}^\delta \right) - \sum_{\delta, \eta,\xi=1}^q \left[ \xi_k, h_{0jk}^{\alpha \eta} \right] h_{0kj}^{\eta \delta} \left(B_j^{\delta \xi} - \sum_{ \sigma,\gamma=1}^q h_{0jk}^{\delta \sigma} B_k^{\sigma \gamma }h_{0kj}^{\gamma  \xi} - \sum_{\gamma =1}^q \left[ \xi_k, h_{0jk}^{\delta \gamma}\right] h_{0kj}^{\gamma \xi } \right) w_{0j}^\xi  \notag \\
&+ \sum_{\delta=1}^q B_j^{\alpha \delta} \mathcal{L}_{\xi_j-\xi_k}\left( w_{0j}^\delta \right) - \sum_{\delta,\xi=1}^q B_j^{\alpha \delta}\left(B_j^{\delta \xi} - \sum_{ \sigma,\gamma,\xi=1}^q h_{0jk}^{\delta \sigma} B_k^{\sigma \gamma }h_{0kj}^{\gamma  \xi} - \sum_{\gamma =1}^q \left[ \xi_k, h_{0jk}^{\delta \gamma}\right] h_{0kj}^{\gamma \xi } \right) w_{0j}^\xi \notag \\
&+ 2 \mathcal{L}_{\xi_j}\left(\beta_j\left( w_{0j}^\alpha \right) \right) - 2 \mathcal{L}_{\xi_k}\left(\beta_k\left(w_{0j}^\alpha \right) \right) - \mathcal{L}_{\xi_k}\mathcal{L}_{\xi_j-\xi_k}\left(w_{0j}^\alpha \right) + \mathcal{L}_{\xi_k}\left( \sum_{\delta=1}^q B_j^{\alpha\delta }w_{0j}^\delta - \sum_{\beta, \eta=1}^q h_{0jk}^{\alpha \beta} B_k^{\beta\eta} w_{0k}^\eta - \sum_{\gamma, \eta=1}^q \left[ \xi_k, h_{0jk}^{\alpha \gamma}\right] w_{0k}^\eta \right) \notag \\
&- \mathcal{L}_{\xi_j}\mathcal{L}_{\xi_j-\xi_k}\left( w_{0j}^\alpha \right) + \mathcal{L}_{\xi_j}\left( \sum_{\delta=1}^q B_j^{\alpha \delta } w_{0j}^\delta - \sum_{\beta, \eta=1}^q h_{0jk}^{\alpha \beta} B_k^{\beta\eta} w_{0k}^\eta - \sum_{\gamma, \eta=1}^q \left[ \xi_k, h_{0jk}^{\alpha \gamma}\right] w_{0k}^\eta  \right) \notag
\end{align}}}

Then from $(\ref{b17})$ and $(\ref{b16})$, we have
{\Small{\begin{align*}
& -\mathcal{L}_{d_{jk}}\left( w_{0j}^\alpha \right) + e_{jk}\left( w_{0j}^\alpha \right) + G_{jk}\left( w_{0j}^\alpha \right) \\
&-\mathcal{L}_{\left[ \xi_j, \xi_k \right]}\left( w_{0j}^\alpha \right) +\left(\left[ \xi_j, B_k^j\right] - \left[ \xi_k, B_j\right] + B_j \circ B_k^j - B_k^j \circ B_j \right) \left( w_{0j}^\alpha \right) - \beta_j\left( B_{jk}^\alpha \right) - \beta_k \left( B_{jk}^\alpha \right) + \mathcal{L}_{\theta_{jk}}\left( \beta_j\left(w_{0j}^\alpha \right) \right) + \mathcal{L}_{\theta_{jk}}\left( \beta_k\left( w_{0j}^\alpha \right)  \right)\\
&= -\cancel{\mathcal{L}_{\xi_j}\mathcal{L}_{\xi_k}\left( w_{0j}^\alpha \right) }+ \cancel{\mathcal{L}_{\xi_k} \mathcal{L}_{\xi_j}\left( w_{0j}^\alpha \right) }+ \cancel{\sum_{\beta, \eta=1}^q \mathcal{L}_{\xi_j} \left( h_{0jk}^{\alpha \beta} B_k^{\beta \eta} w_{0k}^\eta \right) }  - \cancel{\sum_{\delta, \beta, \eta=1}^q h_{0jk}^{\alpha \beta} B_k^{\beta \eta} h_{0kj}^{\eta \delta} \mathcal{L}_{\xi_j}\left( w_{0j}^\delta \right)} + \cancel{\sum_{\eta=1}^q \mathcal{L}_{\xi_j}\left( \left[ \xi_k, h_{0jk}^{\alpha \eta} \right] w_{0k}^\eta \right) }\\
&-\bcancel{\sum_{\eta, \delta=1}^q \left[ \xi_k, h_{0jk}^{\alpha \eta}\right] h_{0kj}^{\eta \delta}\mathcal{L}_{\xi_j}\left(w_{0j}^\delta \right)}- \bcancel{\sum_{\delta=1}^q \left[ \xi_k, B_j^{\alpha \delta} \right] w_{0j}^\delta }+ \bcancel{ \sum_{\xi, \beta, \eta=1}^q h_{0jk}^{\alpha\beta} B_k^{\beta \eta} h_{0kj}^{\eta \delta} B_j^{\delta \xi} w_{0j}^\xi} + \bcancel{ \sum_{\xi, \eta=1}^q \left[ \xi_k, h_{0jk}^{\alpha \eta} \right] h_{0kj}^{\eta \delta} B_j^{\delta \xi} w_{0j}^\xi }\\
&- \bcancel{\sum_{\beta, \eta, \xi=1}^q B_j^{\alpha \xi} h_{0jk}^{\xi\beta} B_k^{\beta \eta} w_{0k}^\eta} - \bcancel{ \sum_{\delta, \xi, \eta=1}^q B_j^{\alpha \xi} \left[ \xi_k, h_{0jk}^{\xi\eta}\right] w_{0k}^\eta}\\
&- 2\beta_j\left( B_j^\alpha \right) + 2\sum_{\beta=1}^q h_{0jk}^{\alpha\beta} \beta_k\left( B_k^\beta \right) + \cancel{\sum_{\delta, \beta, \eta=1}^q h_{0jk}^{\alpha\beta} B_k^{\beta \eta} h_{0kj}^{\eta \delta} \mathcal{L}_{\xi_j}\left( w_{0j}^\delta \right)} - \xcancel{\sum_{\delta, \beta, \eta=1}^q h_{0jk}^{\alpha \beta} B_k^{\beta \eta} h_{0kj}^{\eta \delta} \mathcal{L}_{\xi_k} \left( w_{0j}^\delta \right)}\\
&-\bcancel{ \sum_{\delta, \beta, \eta, \xi=1}^q h_{0jk}^{\alpha\beta} B_k^{\beta \eta} h_{0kj}^{\eta \delta} B_j^{\delta \xi} w_{0j}^\xi }+ \sum_{\beta,  \eta=1}^q h_{0jk}^{\alpha\beta} B_k^{\beta \eta} B_k^{\eta \gamma} w_{0k}^\gamma + \xcancel{\sum_{\beta, \eta, \delta =1}^q  h_{0jk}^{\alpha\beta} B_k^{\beta \eta} h_{0kj}^{\eta \delta} \mathcal{L}_{\xi_k}\left(w_{0j}^\delta \right) }- \sum_{\beta, \eta =1}^q h_{0jk}^{\alpha\beta} B_k^{\beta \eta} \mathcal{L}_{\xi_k}\left(w_{0k}^\eta \right)\\
& - 2\sum_{\gamma ,\eta=1}^q h_{0jk}^{\alpha \gamma} \mathcal{L}_{ \xi_k}\left( \beta_k\left(w_{0k}^\gamma \right) \right) + \bcancel{ \sum_{\delta, \gamma, \eta=1}^q \left[ \xi_k, h_{0jk}^{\alpha \gamma} \right] h_{0kj}^{\gamma \delta} \mathcal{L}_{\xi_j}\left(w_{0j}^\delta \right)} -  \xcancel{\sum_{\delta, \gamma, \eta=1}^q \left[ \xi_k, h_{0jk}^{\alpha \gamma} \right] h_{0kj}^{\gamma \delta} \mathcal{L}_{\xi_k}\left(w_{0j}^\delta \right) } - \bcancel{ \sum_{\delta, \eta, \xi=1}^q \left[ \xi_k, h_{0jk}^{\alpha \eta} \right] h_{0kj}^{\eta \delta} B_j^{\delta\xi} w_{0j}^\xi } \\
& + \xcancel{\sum_{ \eta,\gamma=1}^q \mathcal{L}_{\xi_k}\left( h_{0jk}^{\alpha \eta} B_k^{\eta \gamma} w_{0k}^\gamma \right) } - \sum_{\eta, \gamma=1}^q h_{0jk}^{\alpha \eta} \mathcal{L}_{\xi_k}\left( B_k^{\eta \gamma} w_{0k}^\gamma  \right) - \xcancel{ \sum_{\delta, \eta=1}^q \left[ \xi_k, h_{0jk}^{\alpha \eta} \right]\left[ \xi_k, h_{0kj}^{\eta \delta} \right] w_{0j}^\delta} + \sum_{\delta=1}^q B_j^{\alpha \delta} \mathcal{L}_{\xi_j}\left(w_{0j}^\delta \right) -  \cancel{ \sum_{\delta=1}^q B_j^{\alpha \delta} \mathcal{L}_{\xi_k}\left(w_{0j}^\delta \right)} \\
&-\sum_{\delta, \xi=1}^q B_j^{\alpha \delta} B_j^{\delta \xi} w_{0j}^\xi + \bcancel{\sum_{\delta, \sigma, \gamma=1}^q B_j^{\alpha \delta} h_{0jk}^{\delta \sigma} B_k^{\sigma \gamma} w_{0k}^\gamma } + \bcancel{\sum_{\gamma, \delta=1}^q B_j^{\alpha \delta} \left[ \xi_k, h_{0jk}^{\delta \gamma} \right] w_{0k}^\gamma  }+ 2\mathcal{L}_{\xi_j}\left(\beta_j\left( w_{0j}^\alpha \right) \right) - \cancel{\mathcal{L}_{\xi_k} \mathcal{L}_{\xi_j} \left( w_{0j}^\alpha \right) } \\
& + \cancel{ \sum_{\beta=1}^q \mathcal{L}_{\xi_k}\left( \left[ \xi_k, h_{0jk}^{\alpha\beta}\right] w_{0k}^\beta \right) } +\xcancel{ \sum_{\beta,\eta=1}^q \left[ \xi_k, h_{0jk}^{\alpha\beta}\right]  \left[ \xi_k,  h_{0kj}^{\beta \eta}\right]w_{0j}^\eta} +\xcancel{ \sum_{\beta, \eta=1}^q \left[ \xi_k, h_{0kj}^{\alpha\beta}\right] h_{0kj}^{\beta \eta} \mathcal{L}_{\xi_k}\left(w_{0j}^\eta \right) } +\sum_{\beta=1}^q h_{0jk}^{\alpha\beta} \mathcal{L}_{\xi_k}\mathcal{L}_{\xi_k} \left( w_{0k}^\beta \right)  \\
&+ \bcancel{\sum_{\delta=1}^q \left[ \xi_k, B_j^{\alpha \delta} \right] w_{0j}^\delta }+ \cancel{\sum_{\delta=1}^q B_j^{\alpha \delta} \mathcal{L}_{\xi_k}\left(w_{0j}^\delta \right)} - \xcancel{\sum_{\beta, \eta=1}^q \mathcal{L}_{\xi_k}\left( h_{0jk}^{\alpha \beta} B_k^{\beta \eta} w_{0k}^\eta  \right) } - \cancel{ \sum_{\gamma, \eta=1}^q \mathcal{L}_{\xi_k}\left( \left[ \xi_k, h_{0jk}^{\alpha \gamma} \right] w_{0k}^\eta \right) } - \mathcal{L}_{\xi_j} \mathcal{L}_{\xi_j}\left( w_{0j}^\alpha \right) + \cancel{\mathcal{L}_{\xi_j} \mathcal{L}_{\xi_k} \left( w_{0j}^\alpha \right)} \\
& + \sum_{\delta=1}^q \left[ \xi_j, B_j^{\alpha \delta} \right] w_{0j}^\delta + \sum_{\delta=1}^q B_j^{\alpha \delta} \mathcal{L}_{\xi_j}\left( w_{0j}^\delta \right) - \cancel{\sum_{\beta, \eta=1}^q \mathcal{L}_{\xi_j}\left( h_{0jk}^{\alpha \beta} B_k^{\beta\eta} w_{0k}^\eta\right)} - \cancel{\sum_{\gamma, \eta=1}^q \mathcal{L}_{\xi_j}\left(\left[ \xi_k, h_{0jk}^{\alpha \gamma}\right] w_{0k}^\gamma \right)}
\end{align*}}}

Then we have
{\small{\begin{align*}
&-\mathcal{L}_{d_{jk}}\left( w_{0j}^\alpha \right) + e_{jk}\left( w_{0j}^\alpha \right) + G_{jk}\left( w_{0j}^\alpha \right)\\
&= -2 \beta_j\left(B_j^\alpha \right) + 2 \sum_{\beta=1}^q h_{0jk}^{\alpha\beta} \beta_k\left( B_k^\beta  \right) + \sum_{\beta, \eta=1}^q h_{0jk}^{\alpha\beta} B_k^{\beta \eta} B_k^{\eta \gamma} w_{0k}^\gamma - 2\sum_{\gamma, \eta=1}^q h_{0jk}^{\alpha\gamma} \mathcal{L}_{\xi_k}\left( \beta_k\left( w_{0k}^\gamma  \right) \right)  - 2\sum_{\beta, \eta=1}^q h_{0jk}^{\alpha \beta} B_k^{\beta \eta} \mathcal{L}_{\xi_k}\left( w_{0k}^\eta \right) \\
&- \sum_{\eta, \gamma=1}^q h_{0jk}^{\alpha \eta} \left[ \xi_k, B_k^{\eta \gamma} \right] w_{0k}^\gamma +2 \sum_{\delta=1}^q B_j^{\alpha \delta} \mathcal{L}_{\xi_j}\left( w_{0j}^\alpha \right) - \sum_{\delta, \xi=1}^q B_j^{\alpha \delta} B_j^{\delta \xi} w_{0j}^\xi + 2 \mathcal{L}_{\xi_j}\left(\beta_j\left(w_{0j}^\alpha \right) \right) + \sum_{\beta=1}^q h_{0jk}^{\alpha\beta} \mathcal{L}_{\xi_k}\mathcal{L}_{\xi_k} \left(w_{0k}^\beta \right) + \sum_{\delta=1}^q \left[ \xi_j, B_j^{\alpha \delta} \right] w_{0j}^\delta \\
&- \mathcal{L}_{\xi_j} \mathcal{L}_{\xi_j}\left( w_{0j}^\alpha \right) \\
&= - 2\beta_j\left( B_j^\alpha \right) +2 \mathcal{L}_{\xi_j}\left(\beta_j\left( w_{0j}^\alpha \right) \right)+ 2\sum_{\delta=1}^q B_j^{\alpha \delta} \mathcal{L}_{\xi_j}\left(w_{0j}^\delta \right) + \sum_{\delta=1}^q \left[ \xi_j, B_j^{\alpha \delta} \right] w_{0j}^\delta -\sum_{\delta, \xi=1}^q B_j^{\alpha \delta} B_j^{\delta \xi} w_{0j}^\xi - \mathcal{L}_{\xi_j} \mathcal{L}_{\xi_j} \left(w_{0j}^\alpha \right) \\
& - \sum_{\beta=1}^q h_{0jk}^{\alpha\beta}\left(  -2 \beta_k\left( w_{0k}^\beta \right) + 2\mathcal{L}_{\xi_k}\left( \beta_k\left( w_{0k}^\beta\right)  \right) +2 \sum_{\eta=1}^q B_k^{\beta \eta}\mathcal{L}_{\xi_k}\left( w_{0k}^\eta \right) + \sum_{\gamma=1}^q \left[ \xi_k, B_k^{\eta \gamma}\right] w_{0k}^\gamma - \sum_{\eta, \gamma=1}^q B_k^{\beta \eta}B_k^{\eta \gamma}  w_{0k}^{\gamma } - \mathcal{L}_{\xi_k}\mathcal{L}_{\xi_k}\left( w_{0k}^\beta \right) \right)
\end{align*}}}

Then we define
\begin{align*}
f_j^\alpha : = -2 \sum_{\delta=1}^q B_j^{\alpha \delta} \beta_j\left(w_{0j}^\delta \right) +2 \mathcal{L}_{\xi_j}\left(\beta_j\left( w_{0j}^\alpha \right) \right)+ 2\sum_{\delta=1}^q B_j^{\alpha \delta} \mathcal{L}_{\xi_j}\left(w_{0j}^\delta \right) + \sum_{\delta=1}^q \left[ \xi_j, B_j^{\alpha \delta} \right] w_{0j}^\delta -\sum_{\delta, \xi=1}^q B_j^{\alpha \delta} B_j^{\delta \xi} w_{0j}^\xi - \mathcal{L}_{\xi_j} \mathcal{L}_{\xi_j} \left(w_{0j}^\alpha \right)
\end{align*}

Then we have
\begin{align*}
\bar{\partial} f_j^\alpha&= - 2 \sum_{\delta=1}^q k_j^{\alpha \delta} \wedge \beta_j\left(w_{0j}^\delta \right) + 2 \mathcal{L}_{\psi}\left(\beta_j\left(w_{0j}^\alpha \right) \right) + 2 \sum_{\delta=1}^q k_j^{\alpha \delta} \wedge \mathcal{L}_{\xi_j}\left( w_{0j}^\delta\right) + 2 \sum_{\delta=1}^q B_j^{\alpha \delta} \mathcal{L}_\psi \left( w_{0j}^\delta \right) + \sum_{\delta=1}^q \left[ \psi, B_j^{\alpha \delta} \right] w_{0j}^\delta\\
&+ \sum_{\delta=1}^q \left[ \xi_j, k_j^{\alpha \delta}\right] w_{0j}^\delta - \sum_{\delta, \xi=1}^q k_j^{\alpha \delta} \wedge B_j^{\delta \xi} w_{0j}^\xi - \sum_{\delta , \xi=1}^q B_j^{\alpha \delta} k_j^{\delta \xi} \wedge w_{0j}^\xi - \mathcal{L}_\psi \mathcal{L}_\xi\left(w_{0j}^\alpha \right) - \mathcal{L}_{\xi_j}\mathcal{L}_\psi \left(w_{0j}^\alpha \right)
\end{align*}
Then 
{\small{\begin{align*}
& - \mathcal{L}_{p_j}\left(w_{0j}^\alpha \right) + u_j^\alpha + \bar{\partial} f_j^\alpha \\
& =  - \mathcal{L}_{\left[ \psi, \xi_j \right]}\left(w_{0j}^\alpha \right) + \left( \sum_{\beta=1}^q \left[ \psi, B_j^{\alpha \delta} \right]w_{0j}^\delta - \cancel{ \sum_{\beta=1}^q \left[ \xi_j, k_j^{\alpha \delta} \right] w_{0j}^\delta} +  \bcancel{ \sum_{\delta, \xi=1}^q B_j^{\alpha \delta} k_j^{\delta \xi} \wedge w_{0j}^\xi } - \sum_{\delta, \xi=1}^q k_j^{\alpha \delta} \wedge B_j^{\delta \xi} w_{0j}^\xi\right) \\
& - 2 \sum_{\delta=1}^q k_j^{\alpha \delta} \wedge \beta_j\left(w_{0j}^\delta \right) + 2 \mathcal{L}_{\psi}\left(\beta_j\left(w_{0j}^\alpha \right) \right) + 2 \sum_{\delta=1}^q k_j^{\alpha \delta} \wedge \mathcal{L}_{\xi_j}\left( w_{0j}^\delta\right) + 2 \sum_{\delta=1}^q B_j^{\alpha \delta} \mathcal{L}_\psi \left( w_{0j}^\delta \right) + \sum_{\delta=1}^q \left[ \psi, B_j^{\alpha \delta} \right] w_{0j}^\delta\\
&+ \cancel{\sum_{\delta=1}^q \left[ \xi_j, k_j^{\alpha \delta}\right] w_{0j}^\delta} - \sum_{\delta, \xi=1}^q k_j^{\alpha \delta} \wedge B_j^{\delta \xi} w_{0j}^\xi - \bcancel{\sum_{\delta , \xi=1}^q B_j^{\alpha \delta} k_j^{\delta \xi} \wedge w_{0j}^\xi }- \mathcal{L}_\psi \mathcal{L}_{\xi_j}\left(w_{0j}^\alpha \right)- \mathcal{L}_{\xi_j}\mathcal{L}_\psi \left(w_{0j}^\alpha \right)\\
&=- 2\mathcal{L}_\psi\mathcal{L}_{\xi_j}\left(w_{0j}^\alpha \right) +2 \mathcal{L}_\psi\left( \sum_{\delta=1}^q B_j^{\alpha \delta} w_{0j}^\delta \right) + 2 \mathcal{L}_\psi\left(\beta_j\left(w_{0j}^\alpha \right) \right) + 2 \sum_{\delta=1}^q k_j^{\alpha \delta} \wedge \mathcal{L}_{\xi_j}\left(w_{0j}^\delta \right) - 2 \sum_{\delta, \xi=1}^q k_j^{\alpha \delta} \wedge B_j^{\delta \xi} w_{0j}^\xi - 2\sum_{\delta=1}^q k_j^{\alpha \delta}\wedge \beta_j\left(w_{0j}^\delta \right)\\
&=-2 \mathcal{L}_\psi\left( \mathcal{L}_{\xi_j}- B_j^\alpha - \beta_j\left( w_{0j}^\alpha\right)        \right) + 2\sum_{\delta=1}^q k_j^{\alpha \delta} \wedge \left( \mathcal{L}_{\xi_j}\left(w_{0j}^\delta \right)    - B_j^\delta - \beta_j\left(w_{0j}^\delta \right) \right) = -2 \mathcal{L}_\psi \left(\Psi_j^\alpha  \right) + 2 \sum_{\delta=1}^q k_j^{\alpha \delta} \wedge \Psi_j^\delta
\end{align*}}}

We also note that we have
\begin{align}
&c_j^\alpha + \sum_{\beta=1}^q w_{0j}^1\wedge \cdots \wedge f_j^\beta \wedge \cdots \wedge w_{0j}^q \wedge dw_{0j}^\alpha + w_{0j}^1 \wedge \cdots \wedge w_{0j}^q \wedge \partial f_j^\alpha \label{aa1} \\
&=  2\sum_{\beta=1}^q w_{0j}^1\wedge \cdots \wedge \beta_j\left( w_{0j}^\beta \right) \wedge \cdots \wedge w_{0j}^q \wedge d\beta_j\left(w_{0j}^\alpha \right) + 2\sum_{1\leq \gamma <\eta \leq q} w_{0j}^1 \wedge \cdots \wedge \beta_j\left(w_{0j}^\gamma \right) \wedge \cdots \wedge \beta_j\left( w_{0j}^\eta \right) \wedge \cdots \wedge w_{0j}^q \wedge dw_{0j}^\alpha \notag  \\
&+ \sum_{\beta=1}^q w_{0j}^1 \wedge \cdots \wedge \left( -2 \beta_j\left( B_j^\beta \right)  +2 \mathcal{L}_{\xi_j}\left(\beta_j\left(w_{0j}^\beta \right) \right) +2 \mathcal{L}_{\xi_j}\left(B_j^\beta \right)  -\mathcal{L}_{\xi_j}\mathcal{L}_{\xi_j}\left( w_{0j}^\beta \right) \right) \wedge \cdots \wedge w_{0j}^q \wedge dw_{0j}^\alpha  \notag \\
&+ w_{0j}^1 \wedge \cdots \wedge w_{0j}^q \wedge \partial\left( -2 \beta_j\left( B_j^\alpha  \right)  +2 \mathcal{L}_{\xi_j}\left(\beta_j\left(w_{0j}^\alpha \right) \right) +2 \mathcal{L}_{\xi_j}\left( B_j^\alpha \right)  -\mathcal{L}_{\xi_j}\mathcal{L}_{\xi_j}\left( w_{0j}^\alpha \right) \right) \notag
\end{align}

On the other hand, in a similar way to $(\ref{oc60})$ and $(\ref{oc68})$ and $(\ref{oc61})$ and $(\ref{oc62})$, and the following equality,  
{\small{\begin{align*}
&\sum_{\beta=1}^q w_{0j}^1 \wedge \cdots \wedge B_j^\beta \wedge \cdots \wedge w_{0j}^q \wedge d w_{0j}^\alpha =0 \\
& \Longrightarrow \sum_{1 \leq \gamma <\eta \leq q}^q w_{0j}^1 \wedge \cdots \wedge \mathcal{L}_{\xi_j}\left( w_{0j}^\gamma \right) \wedge \cdots \wedge B_j^\eta \wedge \cdots \wedge w_{0j}^q \wedge dw_{0j}^\alpha + \sum_{1\leq \gamma  < \eta \leq q} w_{0j}^1 \wedge \cdots \wedge B_j^\gamma \wedge \cdots \wedge \mathcal{L}_{\xi_j}\left( w_{0j}^\eta \right)\wedge \cdots \wedge w_{0j}^q \wedge dw_{0j}^\alpha \\
&+ \sum_{\beta=1}^q w_{0j}^1 \wedge \cdots \wedge B_j^\beta \wedge \cdots \wedge w_{0j}^q \wedge \partial \mathcal{L}_{\xi_j}\left(w_{0j}^\alpha \right) + \sum_{\beta=1}^q w_{0j}^1 \wedge \cdots \wedge \mathcal{L}_{\xi_j}\left( B_j^\beta \right) \wedge \cdots \wedge w_{0j}^q \wedge dw_{0j}^\alpha = 0
\end{align*}}}
and the equality
\begin{align*}
&\sum_{\beta=1}^q w_{0j}^1 \wedge \cdots \wedge \beta_j\left( B_j^\beta \right) \wedge \cdots \wedge w_{0j}^q \wedge dw_{0j}^\alpha\\
&= \sum_{\beta=1}^q w_{0j}^1 \wedge \cdots \wedge \left(\sum_{\gamma=1}^q B_j^{\beta \gamma} \beta_j\left( w_{0j}^\gamma \right) \right) \wedge \cdots \wedge w_{0j}^q \wedge d w_{0j}^\alpha\\
&= \sum_{\beta\ne\gamma} B_j^{\beta \gamma} w_{0j}^1 \wedge \cdots \wedge \overbrace{\beta_j\left(w_{0j}^\gamma \right)}^{\beta-\textnormal{th}} \wedge \cdots \wedge w_{0j}^\gamma \wedge \cdots \wedge w_{0j}^q \wedge dw_{0j}^\alpha  + \sum_{\beta=1}^q B_j^{\beta \beta} w_{0j}^1 \wedge \cdots \wedge \beta_j\left( w_{0j}^\beta \right) \wedge \cdots \wedge w_{0j}^q \wedge d w_{0j}^\alpha \\
&=  \sum_{\beta \ne \gamma}  w_{0j}^1 \cdots \wedge \beta_j\left(  w_{0j}^\gamma \right) \wedge \cdots \wedge B_j^{\beta \gamma} w_{0j}^\gamma \wedge \cdots \wedge w_{0j}^q \wedge dw_{0j}^\alpha + \sum_{\beta=1}^q B_j^{\beta \beta} w_{0j}^1 \wedge \cdots \wedge \beta_j\left( w_{0j}^\beta \right) \wedge \cdots \wedge w_{0j}^q \wedge d w_{0j}^\alpha \\
&=  \sum_{\beta\ne \gamma} w_{0j}^1 \wedge \cdots \wedge \beta_j\left( w_{0j}^\gamma \right) \wedge \cdots \wedge \left( B_j^\beta - B_j^{\beta \beta} w_{0j}^\beta \right) \wedge \cdots \wedge w_{0j}^q \wedge dw_{0j}^\alpha + \sum_{\beta=1}^q B_j^{\beta \beta} w_{0j}^1 \wedge \cdots \wedge \beta_j\left( w_{0j}^\beta \right) \wedge \cdots \wedge w_{0j}^q \wedge d w_{0j}^\alpha \\
&=  \sum_{\beta \ne \gamma} w_{0j}^1 \wedge \cdots \wedge \beta_j \left(  w_{0j}^\gamma \right) \wedge \cdots \wedge B_j^\beta \wedge \cdots \wedge w_{0j}^q \wedge dw_{0j}^\alpha + \sum_{\beta , \gamma=1}^q B_j^{\beta \beta} w_{0j}^1 \wedge \cdots \wedge \beta_j\left(w_{0j}^\gamma \right) \wedge \cdots \wedge w_{0j}^q \wedge dw_{0j}^\alpha \\
&= \sum_{\beta \ne \gamma} w_{0j}^1 \wedge \cdots \wedge \beta_j \left(  w_{0j}^\gamma \right) \wedge \cdots \wedge B_j^\beta \wedge \cdots \wedge w_{0j}^q \wedge dw_{0j}^\alpha  - \sum_{\beta=1}^q B_j^{\beta \beta} w_{0j}^1 \wedge \cdots \wedge w_{0j}^q \wedge d \beta_j\left( w_{0j}^\alpha \right)\\
&= - \sum_{\beta \ne \gamma} w_{0j}^1 \wedge \cdots \wedge  B_j^\beta \wedge \cdots \wedge \beta_j \left(  w_{0j}^\gamma \right)  \wedge \cdots \wedge w_{0j}^q \wedge dw_{0j}^\alpha  - \sum_{\beta=1}^q w_{0j}^1 \wedge \cdots \wedge B_j^\beta \wedge \cdots \wedge w_{0j}^q \wedge d \beta_j\left( w_{0j}^\alpha \right)
\end{align*}
we have
\begin{align}
&2\sum_{\beta=1}^q w_{0j}^1 \wedge \cdots \wedge \Psi_j^\beta \wedge \cdots \wedge w_{0j}^q \wedge \partial \Psi_j^\alpha +2 \sum_{1\leq \gamma < \eta \leq q} w_{0j}^1 \wedge \cdots \wedge \Psi_j^\gamma \wedge \cdots \wedge \Psi_j^\eta \wedge \cdots \wedge w_{0j}^q \wedge dw_{0j}^\alpha \label{aa2} \\
&=2 \sum_{\beta=1}^q w_{0j}^1 \wedge \cdots \wedge \left( \mathcal{L}_{\xi_j}\left(w_{0j}^\beta \right) - B_j^\beta - \beta_j\left( w_{0j}^\beta \right) \right) \wedge \cdots \wedge w_{0j}^q \wedge \partial\left( \mathcal{L}_{\xi_j}\left(w_{0j}^\alpha \right) - B_j^\alpha - \beta_j\left( w_{0j}^\alpha \right) \right) \notag \\
&+ 2\sum_{1\leq \gamma < \eta \leq q} w_{0j}^1 \wedge \cdots \wedge \left( \mathcal{L}_{\xi_j}\left(w_{0j}^\gamma \right) - B_j^\gamma - \beta_j\left( w_{0j}^\gamma \right) \right) \wedge \cdots \wedge \left(\mathcal{L}_{\xi_j}\left(w_{0j}^\eta \right) - B_j^\eta - \beta_j\left( w_{0j}^\eta \right) \right) \wedge w_{0j}^q \wedge dw_{0j}^\alpha \notag \\
&=  2\sum_{\beta=1}^q w_{0j}^1\wedge \cdots \wedge \beta_j\left( w_{0j}^\beta \right) \wedge \cdots \wedge w_{0j}^q \wedge d\beta_j\left(w_{0j}^\alpha \right) + 2\sum_{1\leq \gamma <\eta \leq q} w_{0j}^1 \wedge \cdots \wedge \beta_j\left(w_{0j}^\gamma \right) \wedge \cdots \wedge \beta_j\left( w_{0j}^\eta \right) \wedge \cdots \wedge w_{0j}^q \wedge dw_{0j}^\alpha \notag \\
&- \sum_{\beta=1}^q w_{0j}^1 \wedge \cdots \wedge \mathcal{L}_{\xi_j}\mathcal{L}_{\xi_j}\left( w_{0j}^\beta \right) \wedge \cdots \wedge w_{0j}^q \wedge dw_{0j}^\alpha - w_{0j}^1 \wedge \cdots \wedge w_{0j}^q \wedge \partial\left( \mathcal{L}_{\xi_j}\mathcal{L}_{\xi_j}\left( w_{0j}^\alpha \right) \right) + 2w_{0j}^1 \wedge \cdots \wedge w_{0j}^q \wedge \partial \mathcal{L}_{\xi_j}\left(B_j^\alpha \right) \notag \\
&+ 2 w_{0j}^1 \wedge \cdots \wedge w_{0j}^q \wedge \partial \mathcal{L}_{\xi_j}\left( \beta_j\left(w_{0j}^\alpha \right) \right) + 2 \sum_{\beta=1}^q w_{0j}^1 \wedge \cdots \wedge \mathcal{L}_{\xi_j}\left( \beta_j\left( w_{0j}^\beta \right) \right) \wedge \cdots \wedge w_{0j}^q \wedge d w_{0j}^\alpha - 2 w_{0j}^1 \wedge \cdots \wedge w_{0j}^q \wedge \partial \beta_j\left( B_j^\alpha \right) \notag \\
&+ 2 \sum_{\beta=1}^q w_{0j}^1 \wedge \cdots \wedge \mathcal{L}_{\xi_j}\left( B_j^\beta \right) \wedge \cdots \wedge w_{0j}^q \wedge dw_{0j}^\alpha - 2 \sum_{\beta=1}^q w_{0j}^1 \wedge \cdots \wedge \beta_j \left( B_j^\beta \right) \wedge \cdots \wedge w_{0j}^q \wedge dw_{0j}^\alpha \notag
\end{align}
Then $(\ref{aa1})$ and $(\ref{aa2})$ imply that $ \left\{ S_j \right\}$ corresponds to $l$. This completes the proof of Proposition \ref{ii8}.
\end{proof}

\section{Infinitesimal foliated deformations and obstructions in terms of both tangent and cotangent sheaves by using \v Cech resolution of $\mathcal{E}_{\mathcal{F}_0}^\bullet$}

\subsection{Review of deformations of foliated complex analytic structures in terms of both tangent and cotangent sheaves}\

We recall the definitions of a foliated analytic family in terms of both tangent and cotangent sheaves and associated foliated Kodaira-Spencer map in terms of both tangent and cotangent sheaves from Part I.

\begin{definition}\label{s1}
Suppose that given a domain $B\subset \mathbb{C}^m$, there is a set $\left\{\left(M_t, \Theta_{\mathcal{F}_t}, \mathcal{N}_{\mathcal{F}_t}^*\right)|t\in B\right\}$ of $n$-dimensional compact $($singularly$)$ foliated complex manifolds $(M_t,\mathcal{F}_t)$ of dimension $p$ and codimension $q=n-p$, so that for each $t\in B$, we have exact sequences
\begin{align}
0\to \Theta_{\mathcal{F}_t}\to \Theta_{M_t}\to \Theta_{M_t}/\Theta_{\mathcal{F}_t}\to 0\label{s2}\\
0\to \mathcal{N}_{\mathcal{F}_t}^*\to \Omega_{M_t}^1\to \Omega_{M_t}^1/\mathcal{N}_{\mathcal{F}_t}^*\to 0 \label{s3}
\end{align}
where $\Theta_{\mathcal{F}_t}=\left(\Omega_{M_t}^1/\mathcal{N}_{\mathcal{F}}^*\right)^*$ and $\mathcal{N}_{\mathcal{F}_t}^*= \left(\Theta_{M_t}/ \Theta_{\mathcal{F}_t}\right)^*$, and $\Theta_{M_t}/\Theta_{\mathcal{F}_t}$ and $\Omega_{M_t}^1/\mathcal{N}_{\mathcal{F}_t}^*$ are torsion-free $\mathcal{O}_{M_t}$-modules. We say that $\left\{\left(M_t, \Theta_{\mathcal{F}_t}, \mathcal{N}_{\mathcal{F}_t}\right)|t\in B\right\}$ is a family of $($singularly$)$ foliated complex manifold or $($singularly$)$ foliated complex analytic family in terms of both tangent and cotangent sheaves, or foliated analytic family for simultaneous deformations if there exists a complex manifold $\mathcal{M}$ and a holomorphic map $\pi:\mathcal{M}\to B$ such that
\begin{enumerate}
\item we have exact sequences
\begin{align}
0 \to \Theta_\mathcal{F}\to \Theta_{\mathcal{M}/B}\to \Theta_{\mathcal{M}/B}/\Theta_\mathcal{F} \to 0 \label{s4}\\
0\to \mathcal{N}_\mathcal{F}^*\to \Omega_{\mathcal{M}/B}^1\to \Omega_{\mathcal{M}/B}^1/ \mathcal{N}_\mathcal{F}^*\to 0 \label{s5}
\end{align}
where $\Theta_\mathcal{F}$ is a coherent subsheaf of the relative tangent sheaf $\Theta_{\mathcal{M}/B}$ over $B$, and $\mathcal{N}_\mathcal{F}^*$ is a coherent subsheaf of the relative cotangent sheaf $\Omega_{\mathcal{M}/B}^1$ over $B$
\item $\Theta_{\mathcal{F}}=(\Omega_{\mathcal{M}/B}^1/\mathcal{N}_\mathcal{F}^*)^*$ and $\mathcal{N}_\mathcal{F}^*=(\Theta_{\mathcal{M}/B}/\Theta_\mathcal{F})^*$.
\item for each $t\in B$, $\pi^{-1}(t)=\left(M_t, \Theta_{\mathcal{F}_t},\mathcal{N}_{\mathcal{F}_t}^*\right)$, and the exact sequence $(\ref{s2})$ is induced from $(\ref{s4})$ and $(\ref{s3})$ is induced from $(\ref{s5})$ by restricting to $M_t$, respectively.
\item the rank of the Jacobian of $\pi$ is equal to $m$ at every point of $\mathcal{M}$.
\item $[\Theta_\mathcal{F},\Theta_\mathcal{F}]\subset \Theta_\mathcal{F}$ under the Lie bracket, and $\Theta_\mathcal{M}/\Theta_\mathcal{F}$ is flat over $B$.
\item $d \mathcal{N}_\mathcal{F}^* \subset \mathcal{N}_\mathcal{F}^*\bigwedge \Omega_{\mathcal{M}/B}^1$ at every point of $\mathcal{M}-\bigcup_{t\in B} S_t$, where $S_t=\textnormal{Sing}(\mathcal{F}_t)$ for $t\in B$. 
\item $\Omega_{\mathcal{M}/B}^1/\mathcal{N}_\mathcal{F}^*$ is flat over $B$.
\end{enumerate}
We will denote the $($singularly$)$ foliated complex analytic family by $( \mathcal{M} , \Theta_\mathcal{F}, \mathcal{N}_\mathcal{F}^*, B, \pi)$ or simply by $(\mathcal{M}, \mathcal{F}, B, \pi)$.
\end{definition}

\begin{remark}
Let $\left(\mathcal{M}, \Theta_\mathcal{F}, \mathcal{N}_\mathcal{F}^* ,  B, \pi \right)$ with $\Theta_\mathcal{F}$ and $\mathcal{N}_\mathcal{F}^*$ locally free be a foliated analytic family. Let $\Delta$ be an open set of $B$. Then the restriction $\left(\mathcal{M}_\Delta= \pi^{-1}(\Delta), \Theta_{\mathcal{F}}|_\Delta , \mathcal{N}_\mathcal{F}^*|_\Delta, \pi |_{\mathcal{M}_\Delta} \right)$ is also a foliated complex analytic family in terms of both tangent sheaves and cotangent sheaves. We will denote the family by $\left( \mathcal{M}_\Delta, \Theta_{\mathcal{F}_\Delta}, \mathcal{N}_{\mathcal{F}_\Delta}^*,  \Delta, \pi\right)$.
\end{remark}

\begin{remark}\label{os3}
We shall keep the notations in \textnormal{Remark \ref{o4}} and \textnormal{Remark \ref{oc2}}. Then we only need to discuss the condition $(2)$  $\Theta_{\mathcal{F}}=\left(\frac{\Omega_{\mathcal{M}/B}^1}{\mathcal{N}_\mathcal{F}^*}\right)^*$ and $\mathcal{N}_\mathcal{F}^*= \left( \frac{\Theta_{\mathcal{M}/B}}{\Theta_\mathcal{F} } \right)^*$ in \textnormal{Definition \ref{s1}}. The condition $(2)$ implies that  for $\alpha=1,...,p$ and $\beta=1,...,q$, we have
\begin{align}\label{s6}
i_{T_j^\alpha(\xi_j,t)}\left(w_j^\beta(\xi_j,t) \right)=i_{\sum_{\gamma=1}^n T_j^{\alpha\gamma}(\xi_j,t)\frac{\partial}{\partial \xi_j^\gamma} }\left( \sum_{\eta=1}^n w_j^{\beta\eta}(\xi_j,t)d\xi_j^\eta \right)= \sum_{\gamma=1}^n T_j^{\alpha\gamma}(\xi_j,t) w_j^{\beta \gamma}(\xi_j,t)=0
\end{align}
\end{remark}

Given a foliated complex analytic family $\left(\mathcal{M}, \Theta_\mathcal{F}, \mathcal{N}_\mathcal{F}^*,  B, \pi \right)$ with $\Theta_\mathcal{F}$ and $\mathcal{N}_\mathcal{F}^*$ locally free of deformations of $\left(M, \Theta_{\mathcal{F}_0} , \mathcal{N}_{\mathcal{F}_0}^* \right)=\pi^{-1}(0),0\in B$, the infinitesimal foliated deformations in the foliated complex analytic family  are encoded in $\mathbb{H}^1\left( M,\mathcal{F}_0^\bullet \right)$ which is the first cohomology group of the following complex
{\small{\begin{center}
$\mathcal{F}_0^\bullet: \begin{CD}
\cdots \\
@AF_3AA \\
\mathscr{H}om_{\mathcal{O}_{M}}\left(\bigwedge^3 \Theta_{\mathcal{F}_0},  \frac{\Theta_{M}}{\Theta_{\mathcal{F}_0}} \right) \bigoplus \mathscr{H}om_{\mathcal{O}_{M}}\left(\mathcal{N}_{\mathcal{F}_0}^*, \tilde{\mathcal{S}}^3 \right) \bigoplus \mathscr{H}om_{\mathcal{O}_{M}}\left(\mathcal{N}_{\mathcal{F}_0}^*,\bigwedge^2 \Theta_{\mathcal{F}_0}^* \right)       \\
@AF_2AA \\
\mathscr{H}om_{\mathcal{O}_{M}}\left(\bigwedge^2 \Theta_{\mathcal{F}_0}, \frac{\Theta_{M}}{\Theta_{\mathcal{F}_0} }\right) \bigoplus \mathscr{H}om_{\mathcal{O}_{M}}\left(\mathcal{N}_{\mathcal{F}_0}^*, \tilde{\mathcal{S}}^2 \right) \bigoplus \mathscr{H}om_{\mathcal{O}_{M}}\left(\mathcal{N}_{\mathcal{F}_0}^*, \Theta_{\mathcal{F}_0}^* \right)       \\
@AF_1AA \\
\mathscr{H}om_{\mathcal{O}_{M}}\left(\Theta_{\mathcal{F}_0}, \frac{\Theta_{M}}{\Theta_{\mathcal{F}_0}} \right) \bigoplus \mathscr{H}om_{\mathcal{O}_M}\left( \mathcal{N}_{\mathcal{F}_0}^*, \frac{\Omega_{M}^1}{\mathcal{N}_{\mathcal{F}_0}^* } \right)  \\
@AF_0 AA \\
\Theta_{M}
\end{CD}$
\end{center}}}
We will denote the $i$-th cohomology group by $\mathbb{H}^1\left( M , \mathcal{F}_0^\bullet \right)$. We can compute $\mathbb{H}^i\left(M_0, \mathcal{F}_0^\bullet \right)$  by the following \v Cech resolution of $\mathcal{F}_0^\bullet$. Here $\delta$ is the \v Cech map and $\mathcal{U}_0= \mathcal{U}\cap M=\left\{ U_j :=\mathcal{U}_j \cap M | j=1,2,... \right\}$ is an open covering of $M$:
{\tiny{\begin{equation}\label{os1}
\begin{CD}
\cdots \\
@AF_2AA \\
C^0\left( \mathcal{U}_0, \left(\bigwedge^2 \Theta_{\mathcal{F}_0}^*\otimes \frac{\Theta_{M}}{\Theta_{\mathcal{F}_0}}\right) \bigoplus \left(\left(\mathcal{N}_{\mathcal{F}_0}^* \right)^*\otimes \tilde{\mathcal{S}}^2\right) \bigoplus \left( \left( \mathcal{N}_{\mathcal{F}_0}^* \right)^*\otimes \Theta_{\mathcal{F}_0}^* \right)   \right) @>-\delta>> \cdots \\
@AF_1AA @AF_1AA \\
C^0\left(\mathcal{U}_0, \left(\Theta_{\mathcal{F}_0}^*\otimes \frac{\Theta_{M}}{\Theta_{\mathcal{F}_0}}\right) \bigoplus \left(\left(\mathcal{N}_{\mathcal{F}_0}^* \right)^*\otimes \frac{\Omega_{M}^1}{\mathcal{N}_{\mathcal{F}_0}^*} \right) \right) @>\delta>> C^1\left(\mathcal{U}_0, \left(\Theta_{\mathcal{F}_0}^*\otimes \frac{\Theta_{M}}{\Theta_{\mathcal{F}_0}}\right) \bigoplus \left(\left(\mathcal{N}_{\mathcal{F}_0}^* \right)^*\otimes \frac{\Omega_{M}^1}{\mathcal{N}_{\mathcal{F}_0}^*} \right) \right) @>-\delta>> \cdots \\
@AF_0AA @AF_0AA  @AF_0AA \\
C^0\left(\mathcal{U}_0, \Theta_{M} \right) @>-\delta >> C^1\left( \mathcal{U}_0, \Theta_{M} \right) @>>> C^2\left(\mathcal{U}_0, \Theta_{M} \right) @>\delta>> \cdots
\end{CD}
\end{equation}}}
Then we can define the foliated Kodaira-Spencer map at $0\in B$ in terms of both tangent and cotangent sheaves which encodes the information of infinitesimal foliated deformations of $\pi^{-1}(0)=\left(M, \Theta_{\mathcal{F}_0}, \mathcal{N}_{\mathcal{F}_0}^* \right)$ in the foliated  analytic family $\left( \mathcal{M}, \Theta_\mathcal{F}, \mathcal{N}_\mathcal{F}^*, B, \pi \right)$ (see Part I)

\begin{definition}[foliated Kodaira-Spencer map in terms of both tangent sheaf and cotangent sheaf]\label{s8}
Let $(\mathcal{M}, \Theta_\mathcal{F},  \mathcal{N}_\mathcal{F}^*, B, \pi)$ with $\Theta_\mathcal{F}$ and $\mathcal{N}_{\mathcal{F}}^*$ locally free be a foliated complex analytic family of deformations of $\left(M_t, \Theta_{\mathcal{F}_t}, \mathcal{N}_{\mathcal{F}_t}^*\right)=\pi^{-1}(t),t\in B$, where $B$ is a domain of $\mathbb{C}^m$. We keep the notations in \textnormal{Definition \ref{ta12}} and \textnormal{Definition \ref{n15}}. For a tangent vector $\frac{\partial}{\partial t}=\sum_{\lambda=1}^m c_\lambda \frac{\partial}{\partial t_\lambda}, c_\lambda \in \mathbb{C}$, of $B$, the foliated Kodaira-Spencer map in terms of both tangent and cotangent sheaf at $t$ is defined to be a $\mathbb{C}$-linear map
{\small{\begin{align*}
\varphi_0:T_0(B)&\to \mathbb{H}^1\left(M, \mathcal{F}_0^\bullet \right)\\
                \frac{\partial}{\partial t}&\mapsto \frac{\partial \left(M_t, \Theta_{\mathcal{F}_t},  \mathcal{N}_{\mathcal{F}_t}^* \right)}{\partial t}|_{t=0}:=\left(\frac{\partial M_t}{\partial t}, \frac{ \partial \Theta_{\mathcal{F}_t}}{\partial t} , \frac{\partial \mathcal{N}_{\mathcal{F}_t}^*}{\partial t} \right)|_{t=0} \\
                 &\,\,\,\,\,\,\,\,\,\,\in C^1\left(\mathcal{U}_0, \Theta_{M} \right) \bigoplus C^0\left(\mathcal{U}_0, \mathscr{H}om_{\mathcal{O}_{M}} \left( \Theta_{\mathcal{F}_0}, \frac{\Theta_{M}}{\Theta_{\mathcal{F}_0}} \right) \right) \bigoplus C^0\left( \mathcal{U}_0, \mathscr{H}om_{\mathcal{O}_{M}}\left(\mathcal{N}_{\mathcal{F}_0}^*, \frac{\Omega_{M}^1}{\mathcal{N}_{\mathcal{F}_0}^*} \right) \right)
\end{align*}}}
\end{definition}

In Part II we have shown that the complex of sheaves
{\small{\begin{equation} \label{os2}
\mathcal{E}_{\mathcal{F}_0}^\bullet: \begin{CD}
\cdots \\
@AF_3AA \\
\mathscr{H}om_{\mathcal{O}_{M}}\left(\bigwedge^3 \Theta_{\mathcal{F}_0},  \frac{\Theta_{M}}{\Theta_{\mathcal{F}_0}} \right) \bigoplus \mathscr{H}om_{\mathcal{O}_{M}}\left(\mathcal{N}_{\mathcal{F}_0}^*, \tilde{\mathcal{S}}^3 \right) \bigoplus \mathscr{H}om_{\mathcal{O}_{M}}\left(\mathcal{N}_{\mathcal{F}_0}^*,\bigwedge^2 \Theta_{\mathcal{F}_0}^* \right)       \\
@AF_2AA \\
\mathscr{H}om_{\mathcal{O}_{M}}\left(\bigwedge^2 \Theta_{\mathcal{F}_0}, \frac{\Theta_{M}}{\Theta_{\mathcal{F}_0} }\right) \bigoplus \mathscr{H}om_{\mathcal{O}_{M}}\left(\mathcal{N}_{\mathcal{F}_0}^*, \tilde{\mathcal{S}}^2 \right) \bigoplus \mathscr{H}om_{\mathcal{O}_{M}}\left(\mathcal{N}_{\mathcal{F}_0}^*, \Theta_{\mathcal{F}_0}^* \right)       \\
@AF_1'AA \\
\mathscr{H}om_{\mathcal{O}_{M}}\left(\Theta_{\mathcal{F}_0},\Theta_M\right) \bigoplus \mathscr{H}om_{\mathcal{O}_M}\left( \mathcal{N}_{\mathcal{F}_0}^*, \Omega_M^1 \right)  \\
@AF_0' AA \\
\mathcal{E}_{\mathcal{F}_0}
\end{CD}
\end{equation}}}
also controls simultaneous foliated deformations of $\left( M , \Theta_{\mathcal{F}_0}, \mathcal{N}_{\mathcal{F}_0}^*  \right)$ in terms of both tangent and cotangent sheaves. We can actually study obstructions of foliated deformations of $\left(M, \Theta_{\mathcal{F}_0}, \mathcal{N}_{\mathcal{F}_0}^* \right)$ in terms of both tangent and cotangent sheaves by using the complex of sheaves $\mathcal{F}_0^\bullet$.  However obstructions of deformations of the underlying locally free sheaf $\Theta_{\mathcal{F}_0}$ and locally free sheaf $\mathcal{N}_{\mathcal{F}_0}^*$ are implicit in $\mathbb{H}^2\left( M, \mathcal{F}_0^\bullet \right)$, i.e. in \v Cech resolution $(\ref{os1})$ of $\mathcal{F}_0^\bullet$. However the complex of sheaves $\mathcal{E}_{\mathcal{F}_0}^\bullet$ explicitly shows the infinitesimal deformations and obstructions of deformations of the underlying locally free sheaf $\Theta_{\mathcal{F}_0}$ and locally free sheaf $\mathcal{N}_{\mathcal{F}_0}^*$. Hence it is desirable to use $\mathcal{E}_{\mathcal{F}_0}^\bullet$ instead of $\mathcal{F}_0^\bullet$ to identify obstructions clearly. Let us identify the infinitesimal deformations in terms of the complex (\ref{os2}) of sheaves  $\mathcal{E}_{\mathcal{F}_0}^\bullet$.

\subsection{Identification of infinitesimal simultaneous foliated deformations by $\mathcal{E}_{\mathcal{F}_0}^\bullet$}\label{its1} \

Let $\left( \mathcal{M}, \Theta_\mathcal{F}, \mathcal{N}_\mathcal{F}^* ,  B, \pi \right)$ be a foliated analytic family with $\Theta_\mathcal{F}$ and $\mathcal{N}_\mathcal{F}^*$ locally free as in Definition $\ref{s1}$. We keep the notations in Remark $\ref{o4}$, $\ref{oc2}$ and $\ref{os3}$. Let $\mathcal{U}_t=\mathcal{U} \cap M_t =\left\{ U_j^t:=\mathcal{U}_j\cap M_t| j=1,2,... \right\}$ be an open covering of $M_t$. We recall the definition of $\mathcal{E}_{\mathcal{F}_t}$ on $M_t$ with respect to the open covering $\mathcal{U}_t=\left\{U_j^t\right\}$ of $M_t$ from Part II, and let us consider the \v Cech resolution of $\mathcal{E}_{\mathcal{F}_t}^\bullet$:

{\tiny{\begin{equation}\label{os18}
\begin{CD}
\cdots \\
@AF_2^tAA \\
C^0\left( \mathcal{U}_t, \left(\bigwedge^2 \Theta_{\mathcal{F}_t}^*\otimes \frac{\Theta_{M_t}}{\Theta_{\mathcal{F}_t}}\right) \bigoplus \left(\left(\mathcal{N}_{\mathcal{F}_t}^* \right)^*\otimes \tilde{\mathcal{S}}_t^2\right) \bigoplus \left( \left( \mathcal{N}_{\mathcal{F}_t}^* \right)^*\otimes \Theta_{\mathcal{F}_t}^* \right)   \right) @>-\delta>> \cdots \\
@AF_1'^tAA @AF_1'^tAA \\
C^0\left(\mathcal{U}_t, \left(\Theta_{\mathcal{F}_t}^*\otimes \Theta_{M_t}\right) \bigoplus \left(\left(\mathcal{N}_{\mathcal{F}_t}^* \right)^*\otimes \Omega_{M_t}^1 \right) \right) @>\delta>> C^1\left(\mathcal{U}_t, \left(\Theta_{\mathcal{F}_t}^*\otimes  \Theta_{M_t} \right) \bigoplus \left(\left(\mathcal{N}_{\mathcal{F}_t}^* \right)^*\otimes \Omega_{M_t}^1 \right) \right) @>-\delta>> \cdots \\
@AF_0'^tAA @AF_0'^tAA  @AF_0'^tAA \\
C^0\left(\mathcal{U}_t, \mathcal{E}_{\mathcal{F}_t} \right) @>-\delta >> C^1\left( \mathcal{U}_t, \mathcal{E}_{\mathcal{F}_t} \right) @>\delta>> C^2\left(\mathcal{U}_t, \mathcal{E}_{\mathcal{F}_t} \right) @>-\delta>> \cdots
\end{CD}
\end{equation}}}

We keep the notations in subsection \ref{o41} and subsection \ref{oc8}. Then by taking the derivative of $(\ref{s6})$ with respect to $t$ as in Part I, we have the infinitesimal simultaneous foliated deformation along $\frac{\partial}{\partial t}$ at $t$:
\begin{align}
\frac{\partial \left( M_t, \Theta_{\mathcal{F}_t}, \mathcal{N}_{\mathcal{F}_t}^* \right) }{\partial t}=\left( \left\{ \alpha_j(t)\right\}, \left\{ \beta_j(t)\right\}, \left\{ \left( \theta_{jk}(t), A_{jk}(t), B_{jk}(t)\right)  \right\} \right)\in \mathbb{H}^1 \left( M_t, \mathcal{E}_{\mathcal{F}_t}^\bullet   \right) \label{its2}
\end{align}

\subsection{Description of obstructed elements in \v Cech resolution of $\mathcal{E}_{\mathcal{F}_0}^\bullet$}\label{its6} \

\begin{definition}\label{its3}
We say that a compact foliated complex manifold $\left(M, \mathcal{F}_0 \right)$ with $\Theta_{\mathcal{F}_0}$ and $\mathcal{N}_{\mathcal{F}_0}^*$ locally free is unobstructed in simultaneous foliated deformations in terms of both tangent and cotangent sheaves if there is a foliated analytic family $\left( \mathcal{M}, \Theta_\mathcal{F}, \mathcal{N}_\mathcal{F}^*,  B, \pi \right)$ in terms of both tangent and cotangent sheaves of deformations of $\pi^{-1}(0)=\left( M, \Theta_{\mathcal{F}_0} , \mathcal{N}_{\mathcal{F}_0}^* \right), 0\in B$ such that the foliated Kodaira-Spencer map $\varphi_0:T_0B\to \mathbb{H}^1\left( M, \mathcal{E}_{\mathcal{F}_0}^\bullet \right)$ in terms of both tangent and cotangent sheaves is an isomorphism at $0 \in B$. Otherwise, we say that $\left( M, \mathcal{F}_0 \right)$ is obstructed in simultaneous foliated deformations in terms of both tangent and cotangent sheaves.
\end{definition}

\begin{remark}\label{its4}
We note that by theorem of existence for deformations of foliated complex analytic structures in terms of both tangent and cotangent sheaves $($\textnormal{Part I}$)$, if a compact foliated complex manifold $\left(M, \mathcal{F}_0 \right)$ with $\Theta_{\mathcal{F}_0}$ and $\mathcal{N}_{\mathcal{F}_0}^*$ locally free satisfies $\mathbb{H}^2\left(  M,  \mathcal{E}_{\mathcal{F}_0}^\bullet \right)=0$, then there is a foliated complex analytic family $\left( \mathcal{M}, \Theta_\mathcal{F}, \mathcal{N}_\mathcal{F}^,  B, \pi \right)$ such that $\pi^{-1}(0)=(M, \Theta_{\mathcal{F}_0}, \mathcal{N}_{\mathcal{F}_0}^*), 0\in B$ and the foliated Kodaira-Spencer map in terms of both tangent and cotangent sheaves is an isomorphism at $0\in B$, so that $(M, \mathcal{F}_0 )$ is unobstructed in simultaneous foliated deformations in terms of both tangent and cotangent sheaves.
\end{remark}

\begin{definition}\label{its5}
Let $(M, \mathcal{F}_0)$ be a compact foliated complex manifold with $\Theta_{\mathcal{F}_0}$ and $\mathcal{N}_{\mathcal{F}_0}^*$ locally free. We say that $a\in \mathbb{H}^1\left( M, \mathcal{E}_{\mathcal{F}_0}^\bullet \right)$ is an obstructed element in simultaneous foliated deformations of $(M, \mathcal{F}_0 )$ in terms of both tangent and cotangent sheaves if there is no foliated analytic family $\left( \mathcal{M}, \Theta_\mathcal{F}, \mathcal{N}_\mathcal{F}^*, B, \pi \right)$ in terms of both tangent and cotangent sheaves of deformations of $\left(M, \Theta_{\mathcal{F}_0} , \mathcal{N}_{\mathcal{F}_0}^* \right)= \pi^{-1}(0), 0\in B$ such that $a\in \varphi_0\left(T_0 M \right)$ where $\varphi_0: T_0M  \to \mathbb{H}^1\left( M, \mathcal{E}_{\mathcal{F}_0}^\bullet \right)$ is the foliated Kodaira-Spencer map at $0\in B$ of $\left(\mathcal{M}, \Theta_\mathcal{F}, \mathcal{N}_\mathcal{F}^* ,  B, \pi \right)$ in terms of both tangent and cotangent sheaves.
\end{definition}

Hence if there is an obstructed element $a\in \mathbb{H}^1\left( M, \mathcal{E}_{\mathcal{F}_0}^\bullet \right)$ in simultaneous foliated deformations in terms of both tangent and cotangent sheaves, $\left( M, \mathcal{F}_0 \right)$ is obstructed. We now describe a condition when $a\in \mathbb{H}^1\left( M, \mathcal{E}_{\mathcal{F}_0}^\bullet \right)$ is an obstructed element in simultaneous foliated deformations in terms of both tangent and cotangent sheaves. We keep the notations in subsection \ref{o42} and subsection \ref{oc9}. Let $\left(\mathcal{M}, \Theta_{\mathcal{F}}, \mathcal{N}_\mathcal{F}^*, B, \pi \right)$ with $0 \in B \subset \mathbb{C}$ and $\Theta_{\mathcal{F}}, \mathcal{N}_{\mathcal{F}}^*$ locally free be a foliated analytic family in terms of both tangent and cotangent sheaves such that $\pi^{-1}(0) =\left(M, \Theta_{\mathcal{F}_0}, \mathcal{N}_{\mathcal{F}_0}^* \right)$. Take a small disk $\Delta$ with centre $0$ such that $0\in \Delta \subset B$, and represent $\left(\mathcal{M}_\Delta, \Theta_{\mathcal{F}_\Delta}, \mathcal{N}_{\mathcal{F}_\Delta}^* \right)=\pi^{-1}\left(\Delta \right)$ in the form
\begin{align*}
\left( \mathcal{M}_\Delta, \Theta_{\mathcal{F}_\Delta}, \mathcal{N}_{\mathcal{F}_\Delta}^* \right)= \bigcup_{j=1}^l \left( U_j\times \Delta, \Theta_{\mathcal{F}_\Delta}, \mathcal{N}_{\mathcal{F}_\Delta}^* \right)
\end{align*}
as in subsection \ref{o42} and subsection \ref{oc9}. We only need to consider the condition $i_{T_j^\alpha(\xi_j,t)}\left(w_j^\beta(\xi_j,t) \right)=0, \alpha=1,...,p, \beta=1,...,q$. By taking the derivative with respect to $t$, we have
\begin{align*}
i_{\frac{\partial T_j^\alpha(\xi_j,t)}{\partial t}}\left(w_j^\beta(\xi_j,t) \right) + i_{T_j^\alpha(\xi_j,t)}\left(\frac{\partial w_j^\beta(\xi_j,t)}{\partial t} \right)=0
\end{align*}
By taking the derivative with respect to $t$ again, we obtain
\begin{align}
&i_{\frac{\partial^2 T_j^\alpha(\xi_j,t)}{\partial t^2}} \left(w_j^\beta(\xi_j,t) \right)  + i_{T_j^\alpha(\xi_j,t)}\left( \frac{\partial w_j^\beta(\xi_j,t)}{\partial t^2} \right) = -2 i_{\frac{\partial T_j^\alpha(\xi_j,t)}{\partial t}}\left( \frac{\partial w_j^\beta(\xi_j,t)}{\partial t} \right)\label{its8} \\
&\Longrightarrow  i_{\dot{\alpha_j}(t)(T_j^\alpha)}\left( w_j^\beta \right) + i_{T_j^\alpha}\left(\dot{\beta}_j(t)\left(w_j^\beta \right) \right) = 2 i_{\alpha_j(t)(T_j^\alpha)}\left( \beta_j(t)\left(w_j^\beta\right) \right) \label{os8}
\end{align}

We define $\alpha_j(t)*\beta_j(t)\in \Gamma\left( U_j^t, \mathscr{H}om_{\mathcal{O}_M}\left( \mathcal{N}_{\mathcal{F}_t}^*,  \Theta_{\mathcal{F}_t}^* \right) \right)$
\begin{align}\label{os4}
\alpha_j(t)*\beta_j(t) :\Gamma\left( U_j^t, \mathcal{N}_{\mathcal{F}_t}^* \right) &\to \Gamma\left( U_j^t, \Theta_{\mathcal{F}_t}^* \right) \\
                                                         w_j^\beta &\mapsto \left( T_j^\alpha \mapsto 2i_{\alpha_j(t)\left(T_j^\alpha\right)} \left(\beta_j(t)\left(w_j^\beta \right) \right) \right) \notag
\end{align}

By putting $t=0$, we have from $(\ref{o40})$ and $(\ref{oc34})$ and $(\ref{os8})$
{\small{\begin{align}
&D_1'\left(\dot{\alpha}_j(0) \right)= \left[\alpha_j(0), \alpha_j(0) \right]_{\Theta_{\mathcal{F}_0}} \label{os40}  \\
\dot{\alpha}_k(0)- \dot{\alpha}_j(0) + D_0'\left(\left(\dot{\theta}_{jk}(0), \dot{A}_{jk}(0) \right)  \right)& = -\alpha_j(0)\circ A_{jk}(0) - \alpha_k(0)\circ A_{jk}(0) +\mathcal{L}_{\theta_{jk}(0)}\circ \alpha_j(0) + \mathcal{L}_{\theta_{jk}(0)}\circ \alpha_k(0)\label{os15} \\
\left(\dot{\theta}_{ij}(0), \dot{A}_{ij}(0) \right) -& \left( \dot{\theta}_{ik} (0) , \dot{A}_{ik} (0) \right)+ \left( \dot{\theta}_{jk}(0), \dot{A}_{jk}(0) \right)  = \left[ \left( \theta_{ij} (0), A_{ij} (0)   \right) , \left( \theta_{jk} (0) , A_{jk}  (0) \right) \right] \label{os9} \\ 
&E_1'\left( \dot{\beta}_j(0)\right)= \left[ \beta_j(0), \beta_j(0)\right]_{\mathcal{N}_{\mathcal{F}_0}^*}  \label{os13} \\
\dot{\beta}_k(0)-\dot{\beta}_j(0)+ E_0'\left( \left( \dot{\theta}_{jk} (0) , \dot{B}_{jk} (0) \right) \right) & = - \beta_j(0)\circ B_{jk}(0) - \beta_k (0) \circ B_{jk}(0)+ \mathcal{L}_{\theta_{jk}(0)}\circ \beta_j(0)+ \mathcal{L}_{\theta_{jk}(0)}\circ \beta_k(0)  \label{os16} \\
\left( \dot{\theta}_{ij}(0), \dot{B}_{ij}(0)\right)- &\left( \dot{\theta}_{ik}(0), \dot{B}_{ik}(0) \right) + \left( \dot{\theta}_{jk}(0), \dot{B}_{jk}(0) \right)= \left[  \left( \theta_{jk} (0) , B_{jk}(0) \right) , \left( \theta_{jk}(0), B_{jk}(0) \right)   \right] \label{os12}  \\
i_{\dot{\alpha_j}(0)\left(T_{0j}^\alpha\right)} \left( w_{0j}^\beta \right) + i_{T_{0j}^\alpha}&\left(\dot{\beta}_j(0)\left(w_{0j}^\beta \right) \right) = 2 i_{\alpha_j(0)(T_{0j}^\alpha)}\left( \beta_j(0)\left(w_{0j}^\beta\right) \right) , \alpha=1,...,p, \beta=1,...,q \label{os11}
\end{align}}} 

We define 
\begin{align*}
\left[ \left( \alpha_j(0), \beta_j(0) \right), \left(\alpha_j(0) , \beta_j(0) \right)    \right]_{\mathcal{F}_0} &= \left( \left[ \alpha_j(0), \alpha_j(0)\right]_{\Theta_{\mathcal{F}_0}}, \left[ \beta_j(0), \beta_j(0) \right]_{\mathcal{N}_{\mathcal{F}_0}^*} , \alpha_j(0) * \beta_j(0) \right)  \\
 &\in \Gamma\left( U_j^0, \mathscr{H}om_{\mathcal{O}_{M}} \left( \bigwedge^2 \Theta_{\mathcal{F}_0}, \frac{\Theta_{M}}{\Theta_{\mathcal{F}_0}} \right) \bigoplus \mathscr{H}om_{\mathcal{O}_{M}}\left( \mathcal{N}_{\mathcal{F}_0}^*, \tilde{\mathcal{S}}^2 \right) \bigoplus \mathscr{H}om_{\mathcal{O}_{M}}\left( \mathcal{N}_{\mathcal{F}_0}^*,  \Theta_{\mathcal{F}_0}^* \right) \right)
\end{align*}
Then by putting $(\ref{os40}), (\ref{os13})$ and  $(\ref{os11})$ together, and putting $(\ref{os15}$ and $(\ref{os16})$ together and combining $(\ref{os9})$ with $(\ref{os12})$ (see Appendix \ref{apd3}), we can write $(\ref{os40})-(\ref{os11})$ as follows:
{\Small{\begin{align}
&G_1'\left(\dot{\alpha}_j(0) , \dot{\beta}_j(0) \right)= \left[\left(\alpha_j(0), \beta_j(0) \right), \left( \alpha_j(0) , \beta_j(0) \right) \right]_{\mathcal{F}_0}  \label{its12} \\
&\left(\dot{\alpha}_k(0), \dot{\beta}_k(0)\right)- \left(\dot{\alpha}_j(0),\dot{\beta}_j(0)\right) + G_0'\left(\left(\dot{\theta}_{jk}(0), \dot{A}_{jk}(0),\dot{B}_{jk}(0) \right)  \right) \notag \\
& = \left( -\alpha_j(0)\circ A_{jk}(0) - \alpha_k(0)\circ A_{jk}(0) +\mathcal{L}_{\theta_{jk}(0)}\circ \alpha_j(0) + \mathcal{L}_{\theta_{jk}(0)}\circ \alpha_k(0) , - \beta_j(0)\circ B_{jk}(0) - \beta_k (0) \circ B_{jk}(0)+ \mathcal{L}_{\theta_{jk}(0)}\circ \beta_j(0)+ \mathcal{L}_{\theta_{jk}(0)}\circ \beta_k(0) \right) \notag \\
&\left( \dot{\theta}_{ij}(0), \dot{A}_{ij}(0), \dot{B}_{ij}(0)\right)- \left( \dot{\theta}_{ik}(0), \dot{A}_{ik}(0), \dot{B}_{ik}(0) \right) + \left( \dot{\theta}_{jk}(0), \dot{A}_{jk}(0), \dot{B}_{jk}(0) \right)= \left[  \left( \theta_{ij} (0) , A_{ij}(0), B_{ij}(0) \right) , \left( \theta_{jk}(0), A_{jk}(0), B_{jk}(0) \right)   \right] \notag
\end{align}}}

\begin{lemma}\label{us24}
We keep the notations of \textnormal{Lemma \ref{o54}} and \textnormal{Lemma \ref{oc29}}.  For any $1$-cocycle
\begin{align*}
\left( \left\{ \alpha_j\right\}, \left\{ \beta_j \right\}, \left\{ \left( \theta_{jk}, A_{jk}, B_{jk} \right) \right\}  \right) \in  C^0\left( \mathcal{U}, \mathscr{H}om_{\mathcal{O}_M} \left( \Theta_{\mathcal{F}_0}, \Theta_M \right)\right) \bigoplus C^0\left( \mathcal{U}, \mathscr{H}om_{\mathcal{O}_M}\left( \mathcal{N}_{\mathcal{F}_0}^*, \Omega_M^1 \right) \right) \bigoplus C^1\left( \mathcal{U}, \mathcal{E}_{\mathcal{F}_0} \right)
\end{align*}
where $\alpha_j:\Gamma\left( U_j,\Theta_{\mathcal{F}_0}\right) \to \Gamma\left(U_j, \Theta_M \right), \beta_j:\Gamma\left( U_j, \mathcal{N}_{\mathcal{F}_0}^* \right)\to \Gamma\left( U_j, \Omega_M^1 \right), \theta_{jk}\in \Gamma\left( U_{jk}, \Theta_M\right), A_{jk}: \Gamma\left( U_{jk}, \Theta_{\mathcal{F}_0} \right) \to \Gamma\left( U_{jk}, \Theta_{\mathcal{F}_0}\right)$ and $B_{jk}:\Gamma\left( U_{jk}, \mathcal{N}_{\mathcal{F}_0}^* \right) \to \Gamma\left( U_{jk}, \mathcal{N}_{\mathcal{F}_0}^* \right)$ satisfying
\begin{enumerate}
\item $G_1'(\alpha_j,\beta_j)=0$, i.e. $\left[ T_{0j}^\alpha, \alpha_j\left(T_{0j}^\beta\right) \right] - \left[ T_{0j}^\beta, \alpha_j\left(T_{0j}^\alpha \right) \right] - \sum_{\gamma=1}^p g_{0j\alpha\beta}^\gamma \alpha_j\left(T_{0j}^\gamma\right)= Q_{j\alpha\beta}\in \Gamma\left(U_j, \Theta_{\mathcal{F}_0}\right)$ for $\alpha, \beta=1,...,p$, and $w_{0j}^1\wedge \cdots \wedge w_{0j}^q \wedge d \beta_j\left(w_{0j}^\alpha\right) + \sum_{\gamma=1}^q w_{0j}^1 \wedge \cdots \wedge \beta_j\left(w_{0j}^\gamma\right)\wedge \cdots \wedge w_{0j}^q \wedge dw_{0j}^\alpha,\alpha=1,...,q$, and $i_{T_{0i}^\alpha}\left(\beta_j\left(w_{0i}^\beta \right) \right)+ i_{\alpha_j\left(T_{0i}^\alpha \right)}\left( w_{0j}^\beta\right)=0,\alpha=1,...,p, \beta=1,...,q$.
\item $\left(\alpha_k, \beta_k\right)-\left(\alpha_j,\beta_j\right)+G_0'\left( \left(\theta_{jk}, A_{jk}, B_{jk} \right) \right)=0$, i.e. $\alpha_k(T_j^\alpha)- \alpha_j(T_j^\alpha)+ \mathcal{L}_{\theta_{jk}}(T_j^\alpha) - A_{jk}(T_j^\alpha)=0$ for $\alpha=1,...,p$, and $\beta_k\left(w_{0j}^\alpha\right) - \beta_j\left(w_{0j}^\alpha \right) +\mathcal{L}_{\theta_{jk}}\left(w_{0j}^\alpha \right) - B_{jk}\left(w_{0j}^\alpha \right)=0$ for $\alpha=1,...,q$
\item $\left( \theta_{ij}, A_{ij} , B_{ij}\right)  - \left(\theta_{ik}, A_{ik}, B_{ik} \right) + \left(\theta_{jk}, A_{jk}, B_{jk} \right)=0$, i.e. $\theta_{ij}-\theta_{ik}+ \theta_{jk}=0$ and $A_{ij}- A_{ik}+ A_{jk}^i=0$ and $B_{ij}- B_{ik}+ B_{jk}^i =0$.
\end{enumerate}
we define $R_j\in \Gamma\left(U_j, \mathscr{H}om_{\mathcal{O}_M}\left( \bigwedge^2 \Theta_{\mathcal{F}_0}, \frac{\Theta_M}{\Theta_{\mathcal{F}_0}} \right) \right), S_j \in \Gamma\left( U_j, \mathscr{H}om_{\mathcal{O}_M}\left( \mathcal{N}_{\mathcal{F}_0}^*, \tilde{\mathcal{S}}^2 \right) \right), F_{jk} \in \Gamma\left( U_{jk}, \mathscr{H}om_{\mathcal{O}_M}\left(\Theta_{\mathcal{F}_0}, \Theta_M \right) \right)$, $G_{jk}\in \Gamma\left(U_{jk}, \mathscr{H}om_{\mathcal{O}_M}\left( \mathcal{N}_{\mathcal{F}_0}^*, \Omega_M^1 \right) \right)$ as in \textnormal{Lemma \ref{o54}} and \textnormal{Lemma \ref{oc29}}. On the other hand, we define $P_j \in \Gamma\left( U_j, \mathscr{H}om_{\mathcal{O}_M}\left( \mathcal{N}_{\mathcal{F}_0}^*, \Theta_{\mathcal{F}_0}^* \right) \right)$ by
\begin{align*}
P_j: \Gamma\left( U_j, \mathcal{N}_{\mathcal{F}_0}^* \right)&\to \Gamma\left( U_j, \Theta_{\mathcal{F}_0}^* \right)\\
          w_{0j}^\beta &\mapsto \left( T_{0j}^\alpha \mapsto 2 i_{\alpha_j\left( T_{0j}^\alpha \right)} \left( \beta_j \left( w_{0j}^\beta \right) \right)   \right)
\end{align*}
and linearly extends to $\Gamma\left( U_j, \Theta_{\mathcal{F}_0} \right)$ and to $\Gamma\left( U_j, \mathcal{N}_{\mathcal{F}_0}^* \right)$. Then
{\small{\begin{align}\label{its39}
&\left(\left\{ R_j \right\}, \left\{ S_j \right\}, \left\{ P_j \right\} ,  \left\{ F_{jk} \right\}, \left\{ G_{jk} \right\}, \left\{ \left[ \left(\theta_{ij}, A_{ij}, B_{ij} \right), \left( \theta_{jk}, A_{jk}, B_{jk} \right) \right] \right\} \right)\\
& \in C^0\left(\mathcal{U}, \mathscr{H}om_{\mathcal{O}_M}\left( \bigwedge^2 \Theta_{\mathcal{F}_0}, \frac{\Theta_M}{\Theta_{\mathcal{F}_0}} \right) \right) \bigoplus C^0\left(\mathcal{U}, \mathscr{H}om_{\mathcal{O}_M}\left(\mathcal{N}_{\mathcal{F}_0}^*, \tilde{\mathcal{S}}^2 \right) \right) \bigoplus C^0 \left( \mathcal{U}, \mathscr{H}om_{\mathcal{O}_M}\left( \mathcal{N}_{\mathcal{F}_0}^*, \Theta_{\mathcal{F}_0}^* \right) \right) \notag \\
& \bigoplus C^1\left(\mathcal{U}, \mathscr{H}om_{\mathcal{O}_M}\left(\Theta_{\mathcal{F}_0}, \Theta_M \right) \right) \bigoplus C^1\left(\mathcal{U}, \mathscr{H}om_{\mathcal{O}_M}\left(\mathcal{N}_{\mathcal{F}_0}^*, \Omega_M^1 \right) \right) \bigoplus C^2\left(\mathcal{U}, \mathcal{E}_{\mathcal{F}_0} \right) \notag
\end{align}}}
defines a $2$-cocycle in \v Cech resolution $(\ref{os18})$ of $\mathcal{E}_{\mathcal{F}_0}^\bullet$.
\end{lemma}

\begin{remark}
We will deduce the following equality $(\ref{us23})$ below which will be used in the proof of \textnormal{Lemma \ref{us24}}. We note that
{\small{\begin{align*}
&-  i_{\alpha_j\left(T_{0j}^\beta \right)} i_{T_{0j}^\alpha}d\left( \beta_j\left(w_{0j}^\gamma \right) \right)   + i_{\alpha_j\left(T_{0j}^\alpha \right)} i_{T_{0j}^\beta} d\left( \beta_j\left(w_{0j}^\gamma \right)\right)+ i_{\alpha_j\left(T_{0j}^\alpha \right)}i_{\alpha_j\left( T_{0j}^\beta \right)} dw_{0j}^\gamma\\
&= -i_{\alpha_j\left(T_{0j}^\beta \right)} \mathcal{L}_{T_{0j}^\alpha}\left( \beta_j\left(w_{0j}^\gamma \right) \right) + i_{\alpha_j\left(T_{0j}^\beta \right)}d i_{T_{0j}^\alpha}\left(\beta_j\left( w_{0j}^\gamma \right) \right) + i_{\alpha_j\left( T_{0j}^\alpha\right)} \mathcal{L}_{T_{0j}^\beta}\left(\beta_j\left( w_{0j}^\gamma \right) \right) - i_{\alpha_j\left( T_{0j}^\alpha \right)}d i_{T_{0j}^\beta}\left(\beta_j\left( w_{0j}^\gamma \right) \right)\\
&+ i_{\alpha_j\left( T_{0j}^\alpha\right)} \mathcal{L}_{\alpha_j\left(T_{0j}^\beta \right)}\left( w_{0j}^\gamma \right) - i_{\alpha_j\left(T_{0j}^\alpha \right)} d i_{\alpha_j\left(T_{0j}^\beta \right)}\left(w_{0j}^\gamma \right)\\
&= -i_{\alpha_j\left(T_{0j}^\beta \right)} \mathcal{L}_{T_{0j}^\alpha}\left( \beta_j\left(w_{0j}^\gamma \right) \right) - i_{\alpha_j\left(T_{0j}^\beta \right)}d i_{\alpha_j\left(T_{0j}^\alpha\right)}\left(w_{0j}^\gamma \right) + i_{\alpha_j\left( T_{0j}^\alpha\right)} \mathcal{L}_{T_{0j}^\beta}\left(\beta_j\left( w_{0j}^\gamma \right) \right) +\cancel{ i_{\alpha_j\left( T_{0j}^\alpha \right)}d i_{\alpha_j\left(T_{0j}^\beta\right)}\left(w_{0j}^\gamma  \right)}\\
&+ i_{\alpha_j\left( T_{0j}^\alpha\right)} \mathcal{L}_{\alpha_j\left(T_{0j}^\beta \right)}\left( w_{0j}^\gamma \right) - \cancel{i_{\alpha_j\left(T_{0j}^\alpha \right)} d i_{\alpha_j\left(T_{0j}^\beta \right)}\left(w_{0j}^\gamma \right) }\\
&= -i_{\alpha_j\left(T_{0j}^\beta \right)} \mathcal{L}_{T_{0j}^\alpha}\left( \beta_j\left(w_{0j}^\gamma \right) \right) - i_{\alpha_j\left(T_{0j}^\beta \right)}\mathcal{L}_{\alpha_j\left( T_{0j}^\alpha \right)} \left(w_{0j}^\gamma \right) + \cancel{i_{\alpha_j\left( T_{0j}^\beta\right)} i_{\alpha_j\left( T_{0j}^\alpha \right)} d w_{0j}^\gamma} + i_{\alpha_j\left( T_{0j}^\alpha\right)} \mathcal{L}_{T_{0j}^\beta}\left(\beta_j\left( w_{0j}^\gamma \right) \right)  \\
&+i_{\alpha_j\left( T_{0j}^\alpha \right)} d i_{\alpha_j\left( T_{0j}^\beta \right)} \left( w_{0j}^\gamma \right) + \cancel{i_{\alpha_j\left( T_{0j}^\alpha\right)}i_{\alpha_j\left(T_{0j}^\beta \right)} d\left(w_{0j}^\gamma \right)}
\end{align*}}}
Hence we have
{\small{\begin{align}\label{us23}
&-  i_{\alpha_j\left(T_{0j}^\beta \right)} i_{T_{0j}^\alpha}d\left( \beta_j\left(w_{0j}^\gamma \right) \right)   + i_{\alpha_j\left(T_{0j}^\alpha \right)} i_{T_{0j}^\beta} d\left( \beta_j\left(w_{0j}^\gamma \right)\right)+ i_{\alpha_j\left(T_{0j}^\alpha \right)}i_{\alpha_j\left( T_{0j}^\beta \right)} dw_{0j}^\gamma\\
&= -i_{\alpha_j\left(T_{0j}^\beta \right)} \mathcal{L}_{T_{0j}^\alpha}\left( \beta_j\left(w_{0j}^\gamma \right) \right) - i_{\alpha_j\left(T_{0j}^\beta \right)}\mathcal{L}_{\alpha_j\left( T_{0j}^\alpha \right)} \left(w_{0j}^\gamma \right)  + i_{\alpha_j\left( T_{0j}^\alpha\right)} \mathcal{L}_{T_{0j}^\beta}\left(\beta_j\left( w_{0j}^\gamma \right) \right) + \mathcal{L}_{\alpha_j\left( T_{0j}^\alpha \right)} i_{\alpha_j\left( T_{0j}^\beta \right)} \left( w_{0j}^\gamma \right) \notag
\end{align}}}
\end{remark}

\begin{proof}[Proof of \textnormal{Lemma \ref{us24}}]
By Lemma \ref{o54} and Lemma \ref{oc29} it remains to show that $\beta_2\left(R_j \right) + \alpha_2\left(S_j \right) - \tilde{E}_1\left(P_j \right)=0$ and $P_j- P_k + \alpha_1\left( F_{jk} \right) + \beta_1\left( G_{jk} \right)=0$. We show that $\beta_2\left(R_j \right) + \alpha_2\left(S_j \right) - \tilde{E}_1\left(P_j \right)=0$. In fact, for $x\in U_j$, choose a neighborhood of $U_x\subset U_j-S$. Then it is sufficient to check locally on $U_x$. Then
{\small{\begin{align}\label{us26}
&\frac{1}{2}i_{R_j\left(T_{0j}^\alpha \wedge T_{0j}^\beta \right)}\left(w_{0j}^\gamma \right)=i_{\left[\alpha_j\left(T_{0j}^\alpha \right), \alpha_j\left(T_{0j}^\beta \right) \right]} \left(w_{0j}^\gamma \right) -  i_{\alpha_j\left(Q_{j\alpha\beta} \right)}\left(w_{0j}^\gamma \right) = i_{\left[\alpha_j\left(T_{0j}^\alpha \right), \alpha_j\left(T_{0j}^\beta \right) \right]} \left(w_{0j}^\gamma \right) + i_{Q_{j\alpha\beta}}\left(\beta_j\left(w_{0j}^\gamma \right) \right)\\
&= i_{\left[\alpha_j\left( T_{0j}^\alpha\right), \alpha_j\left( T_{0j}^\beta \right) \right]}\left(w_{0j}^\gamma \right) + i_{\left[ T_{0j}^\alpha, \alpha_j\left(T_{0j}^\beta\right) \right] - \left[ T_{0j}^\beta, \alpha_j \left( T_{0j}^\alpha \right) \right]- \sum_{c=1}^p g_{0j\alpha\beta}^c \alpha_j\left( T_{0j}^c \right) } \left(\beta_j\left(w_{0j}^\gamma \right)  \right) \notag\\
&=\mathcal{L}_{\alpha_j\left(T_{0j}^\alpha \right) } i_{\alpha_j\left( T_{0j}^\beta\right)} \left(w_{0j}^\gamma \right) - i_{\alpha_j\left(T_{0j}^\beta \right)} \mathcal{L}_{\alpha_j\left(T_{0j}^\alpha \right)}\left(  w_{0j}^\gamma \right) + \mathcal{L}_{T_{0j}^\alpha}i_{\alpha_j\left(T_{0j}^\beta \right)}\left(\beta_j\left( w_{0j}^\gamma \right) \right) - i_{\alpha_j\left( T_{0j}^\beta \right)} \mathcal{L}_{T_{0j}^\alpha}\left(\beta_j\left( w_{0j}^\gamma \right) \right) \notag \\
&- \mathcal{L}_{T_{0j}^\beta} i_{\alpha_j\left(T_{0j}^\alpha \right)}\left(\beta_j\left(w_{0j}^\gamma \right) \right) + i_{\alpha_j\left(T_{0j}^\alpha \right)} \mathcal{L}_{T_{0j}^\beta}\left( \beta_j\left(w_{0j}^\gamma \right) \right) - \sum_{c=1}^p g_{0j\alpha\beta}^c i_{ \alpha_j\left( T_{0j}^c \right) } \left(\beta_j\left(w_{0j}^\gamma \right)  \right) \notag
\end{align}}}

We recall that from Remark $\ref{oc81}$, we have
\begin{align}\label{its34}
\frac{1}{2}i_{T_{0j}^\alpha \wedge T_{0j}^\beta }\left( S_j\left(w_{0j}^\gamma \right) \right) = - w_{0j}^1 \wedge \cdots\wedge w_{0j}^q \cdot i_{T_{0j}^\alpha \wedge T_{0j}^\beta}\left( \sum_{\eta=1}^q b_{j_x}^{\gamma \eta} \wedge \beta_j\left( w_{0j}^\eta  \right)\right)
\end{align}
Then let us compute
{\small{\begin{align}\label{us27}
&i_{T_{0j}^\alpha\wedge T_{0j}^\beta}\left(\sum_{\eta=1}^q b_{j_x}^{\gamma \eta}\wedge \beta_j\left(w_{0j}^\eta \right)\right) = \sum_{\eta=1}^q i_{T_{0j}^\alpha}(b_{j_x}^{\gamma \eta}) i_{T_{0j}^\beta}\left(\beta_j\left(w_{0j}^\eta \right)\right) - \sum_{\eta=1}^q i_{T_{0j}^\beta}(b_{i_x}^{\gamma \eta}) i_{T_{0j}^\alpha}\left(\beta_j\left(w_{0j}^\eta \right) \right) \\
&=- \sum_{\eta=1}^q i_{T_{0j}^\alpha}\left(b_{j_x}^{\gamma \eta}\right)i_{\alpha_j\left(T_{0j}^\beta \right)}\left(w_{0j}^\eta\right) +  \sum_{\eta=1}^q i_{T_{0i}^\beta}\left(b_{j_x}^{\gamma \eta}\right)i_{\alpha_j\left(T_{0j}^\alpha \right)}\left(w_{0j}^\eta\right)=-i_{\alpha_j\left(T_{0j}^\beta \right)}\left(\sum_{\eta=1}^q i_{T_{0j}^\alpha}\left(b_{j_x}^{\gamma\eta}\right) w_{0j}^\eta \right) + i_{\alpha_j\left( T_{0j}^\alpha \right)}\left(\sum_{\eta=1}^q i_{T_{0j}^\beta}\left(b_{j_x}^{\gamma \eta} \right) w_{0j}^\eta \right)\notag \\
&=-i_{\alpha_j\left(T_{0j}^\beta \right)}i_{T_{0j}^\alpha}\left(d \beta_j\left(w_{0j}^\gamma \right) - \sum_{\eta=1}^q a_{j_x}^{\gamma\eta} \wedge \beta_j\left( w_{0j}^\eta \right) \right) +i_{\alpha_j\left(T_{0j}^\alpha \right)}i_{T_{0j}^\beta}\left(d \beta_j\left(w_{0j}^\gamma \right) - \sum_{\eta=1}^q a_{j_x}^{\gamma\eta} \wedge \beta_j\left( w_{0j}^\eta \right)\right) \notag \\
&=-  i_{\alpha_j\left(T_{0j}^\beta \right)} i_{T_{0j}^\alpha}d\left( \beta_j\left(w_{0j}^\gamma \right) \right) - \sum_{\eta=1}^q i_{\alpha_j\left(T_{0j}^\beta \right)}\left(a_{j_x}^{\gamma\eta} \right) i_{T_{0j}^\alpha}\left(\beta_j\left(w_{0j}^\eta \right) \right)  + \sum_{\eta=1}^q i_{T_{0j}^\alpha}\left(a_{j_x}^{\gamma\eta}\right) i_{\alpha_j\left( T_{0j}^\beta\right)}\left( \beta_j\left( w_{0j}^\eta \right)\right) \notag \\
&+ i_{\alpha_j\left(T_{0j}^\alpha \right)} i_{T_{0j}^\beta} d\left( \beta_j\left(w_{0j}^\gamma \right)\right) - \sum_{\eta=1}^q i_{T_{0j}^\beta}\left(a_{j_x}^{\gamma \eta}\right) i_{\alpha_j\left( T_{0j}^\alpha \right)}\left(\beta_j\left( w_{0j}^\eta \right) \right) + \sum_{\eta=1}^q i_{\alpha_j\left( T_{0j}^\alpha \right)}\left(a_{j_x}^{\gamma \eta} \right) i_{T_{0j}^\beta}\left( \beta_j\left( w_{0j}^\eta \right)\right) \notag \\
&=-  i_{\alpha_j\left(T_{0j}^\beta \right)} i_{T_{0j}^\alpha}d\left( \beta_j\left(w_{0j}^\gamma \right) \right) + \sum_{\eta=1}^q i_{\alpha_j\left(T_{0j}^\beta \right)}\left(a_{j_x}^{\gamma\eta} \right) i_{ \alpha_j\left(T_{0j}^\alpha\right)}\left(w_{0j}^\eta \right)  + \sum_{\eta=1}^q i_{T_{0j}^\alpha}\left(a_{j_x}^{\gamma\eta}\right) i_{\alpha_j\left( T_{0j}^\beta\right)}\left( \beta_j\left( w_{0j}^\eta \right)\right) \notag \\
&+ i_{\alpha_j\left(T_{0j}^\alpha \right)} i_{T_{0j}^\beta} d\left( \beta_j\left(w_{0j}^\gamma \right)\right) - \sum_{\eta=1}^q i_{T_{0j}^\beta}\left(a_{j_x}^{\gamma \eta}\right) i_{\alpha_j\left( T_{0j}^\alpha \right)}\left(\beta_j\left( w_{0j}^\eta \right) \right) - \sum_{\eta=1}^q i_{\alpha_j\left( T_{0j}^\alpha \right)}\left(a_{j_x}^{\gamma \eta} \right) i_{\alpha_j\left(T_{0j}^\beta\right)}\left(  w_{0j}^\eta \right) \notag \\
&=-  i_{\alpha_j\left(T_{0j}^\beta \right)} i_{T_{0j}^\alpha}d\left( \beta_j\left(w_{0j}^\gamma \right) \right)   + \sum_{\eta=1}^q i_{T_{0j}^\alpha}\left(a_{j_x}^{\gamma\eta}\right) i_{\alpha_j\left( T_{0j}^\beta\right)}\left( \beta_j\left( w_{0j}^\eta \right)\right)+ i_{\alpha_j\left(T_{0j}^\alpha \right)} i_{T_{0j}^\beta} d\left( \beta_j\left(w_{0j}^\gamma \right)\right) - \sum_{\eta=1}^q i_{T_{0j}^\beta}\left(a_{j_x}^{\gamma \eta}\right) i_{\alpha_j\left( T_{0j}^\alpha \right)}\left(\beta_j\left( w_{0j}^\eta \right) \right)  \notag \\
&+ i_{\alpha_j\left(T_{0j}^\alpha \right)}i_{\alpha_j\left( T_{0j}^\beta \right)} dw_{0j}^\gamma \notag
\end{align}}}

Let us compute
{\small{\begin{align}\label{us25}
&\frac{1}{2}\tilde{E}_1\left(P_j \right)\left(T_{0j}^\alpha\wedge T_{0j}^\beta \right)\left(  w_{0j}^\gamma \right)\\
&=\frac{1}{2}\left(\left[ T_{0j}^\alpha, P_j\left(T_{0j}^\beta \right)\left(w_{0j}^\gamma \right) \right] - P_j\left(T_{0j}^\beta \right)\left(\mathcal{L}_{T_{0j}^\alpha}\left(w_{0j}^\gamma \right) \right) - \left[ T_{0j}^\beta, P_j\left( T_{0j}^\alpha\right)\left( w_{0j}^\gamma\right)\right] + P_j\left( T_{0j}^\alpha \right)\left(\mathcal{L}_{T_{0j}^\beta} \left(w_{0j}^\gamma \right) \right) - P_j\left(\left[ T_{0j}^\alpha, T_{0j}^\beta \right] \right)\left( w_{0j}^\gamma \right) \right) \notag \\
&=\mathcal{L}_{T_{0j}^\alpha}i_{\alpha_j\left(T_{0j}^\beta \right)}\left( \beta_j\left(w_{0j}^\gamma \right)\right) - P_j\left(T_{0j}^\beta \right)\left(i_{T_{0j}^\alpha} dw_{0j}^\gamma \right) - \mathcal{L}_{T_{0j}^\beta}i_{\alpha_j\left( T_{0j}^\alpha \right)} \left( \beta_j\left(w_{0j}^\gamma \right) \right) + P_j\left(T_{0j}^\alpha \right)\left(i_{T_{0j}^\beta} dw_{0j}^\gamma  \right) - \sum_{c=1}^p g_{0j\alpha\beta}^c i_{\alpha_j\left( T_{0j}^c \right)} \left( \beta_j\left( w_{0j}^\gamma \right) \right) \notag \\
&=\mathcal{L}_{T_{0j}^\alpha}i_{\alpha_j\left(T_{0j}^\beta \right)}\left( \beta_j\left(w_{0j}^\gamma \right)\right) - P_j\left(T_{0j}^\beta \right)\left(\sum_{\eta=1}^q i_{T_{0j}^\alpha}\left( a_{j_x}^{\gamma \eta} \right) w_{0j}^\eta \right) - \mathcal{L}_{T_{0j}^\beta}i_{\alpha_j\left( T_{0j}^\alpha \right)} \left( \beta_j\left(w_{0j}^\gamma \right) \right) + P_j\left(T_{0j}^\alpha \right)\left(\sum_{\eta=1}^q i_{T_{0j}^\beta}\left( a_{j_x}^{\gamma \eta} \right) w_{0j}^\eta  \right) \notag \\
&- \sum_{c=1}^p g_{0j\alpha\beta}^c i_{\alpha_j\left( T_{0j}^c \right)} \left( \beta_j\left( w_{0j}^\gamma \right) \right) \notag \\
&=\mathcal{L}_{T_{0j}^\alpha}i_{\alpha_j\left(T_{0j}^\beta \right)}\left( \beta_j\left(w_{0j}^\gamma \right)\right) -\sum_{\eta=1}^q i_{T_{0j}^\alpha}\left( a_{j_x}^{\gamma \eta} \right) i_{\alpha_j\left( T_{0j}^\beta \right)} \left(\beta_j\left( w_{0j}^\eta \right) \right) - \mathcal{L}_{T_{0j}^\beta}i_{\alpha_j\left( T_{0j}^\alpha \right)} \left( \beta_j\left(w_{0j}^\gamma \right) \right) + \sum_{\eta=1}^q i_{T_{0j}^\beta}\left( a_{j_x}^{\gamma \eta} \right) i_{\alpha_j\left( T_{0j}^\alpha \right)}\left( \beta_j\left( w_{0j}^\eta \right) \right) \notag \\
&- \sum_{c=1}^p g_{0j\alpha\beta}^c i_{\alpha_j\left( T_{0j}^c \right)} \left( \beta_j\left( w_{0j}^\gamma \right) \right) \notag
\end{align}}}

Then  from $(\ref{us26}), (\ref{us27}),(\ref{us25}), (\ref{us23})$, we have
\begin{align}\label{its20}
&\frac{1}{2}i_{R_j}\left(T_{0j}^\alpha \wedge T_{0j}^\beta \right)\left(w_{0j}^\gamma \right) - i_{T_{0j}^\alpha\wedge T_{0j}^\beta}\left(\sum_{\eta=1}^q b_{j_x}^{\gamma \eta}\wedge \beta_j\left(w_{0j}^\eta \right)\right) - \frac{1}{2}\tilde{E}_1\left( P_j \right)\left(T_{0j}^\alpha \wedge T_{0j}^\beta \right)\left(w_{0j}^\gamma \right)\\
&=\mathcal{L}_{\alpha_j\left(T_{0j}^\alpha \right) } i_{\alpha_j\left( T_{0j}^\beta\right)} \left(w_{0j}^\gamma \right) - i_{\alpha_j\left(T_{0j}^\beta \right)} \mathcal{L}_{\alpha_j\left(T_{0j}^\alpha \right)}\left(  w_{0j}^\gamma \right) + \xcancel{\mathcal{L}_{T_{0j}^\alpha}i_{\alpha_j\left(T_{0j}^\beta \right)}\left(\beta_j\left( w_{0j}^\gamma \right) \right)} - i_{\alpha_j\left( T_{0j}^\beta \right)} \mathcal{L}_{T_{0j}^\alpha}\left(\beta_j\left( w_{0j}^\gamma \right) \right) \notag \\
&- \xcancel{\mathcal{L}_{T_{0j}^\beta} i_{\alpha_j\left(T_{0j}^\alpha \right)}\left(\beta_j\left(w_{0j}^\gamma \right) \right)} + i_{\alpha_j\left(T_{0j}^\alpha \right)} \mathcal{L}_{T_{0j}^\beta}\left( \beta_j\left(w_{0j}^\gamma \right) \right) - \xcancel{ \sum_{c=1}^p g_{0j\alpha\beta}^c i_{ \alpha_j\left( T_{0j}^c \right) } \left(\beta_j\left(w_{0j}^\gamma \right)  \right)} \notag \\
&+  i_{\alpha_j\left(T_{0j}^\beta \right)} i_{T_{0j}^\alpha}d\left( \beta_j\left(w_{0j}^\gamma \right) \right)   - \cancel{\sum_{\eta=1}^q i_{T_{0j}^\alpha}\left(a_{j_x}^{\gamma\eta}\right) i_{\alpha_j\left( T_{0j}^\beta\right)}\left( \beta_j\left( w_{0j}^\eta \right)\right)}- i_{\alpha_j\left(T_{0j}^\alpha \right)} i_{T_{0j}^\beta} d\left( \beta_j\left(w_{0j}^\gamma \right)\right) + \bcancel{ \sum_{\eta=1}^q i_{T_{0j}^\beta}\left(a_{j_x}^{\gamma \eta}\right) i_{\alpha_j\left( T_{0j}^\alpha \right)}\left(\beta_j\left( w_{0j}^\eta \right) \right)} \notag \\
&- i_{\alpha_j\left(T_{0j}^\alpha \right)}i_{\alpha_j\left( T_{0j}^\beta \right)} dw_{0j}^\gamma \notag \\
&-\xcancel{\mathcal{L}_{T_{0j}^\alpha}i_{\alpha_j\left(T_{0j}^\beta \right)}\left( \beta_j\left(w_{0j}^\gamma \right)\right) }+\cancel{\sum_{\eta=1}^q i_{T_{0j}^\alpha}\left( a_{j_x}^{\gamma \eta} \right) i_{\alpha_j\left( T_{0j}^\beta \right)} \left(\beta_j\left( w_{0j}^\eta \right) \right) }+ \xcancel{\mathcal{L}_{T_{0j}^\beta}i_{\alpha_j\left( T_{0j}^\alpha \right)} \left( \beta_j\left(w_{0j}^\gamma \right) \right) } - \bcancel{ \sum_{\eta=1}^q i_{T_{0j}^\beta}\left( a_{j_x}^{\gamma \eta} \right) i_{\alpha_j\left( T_{0j}^\alpha \right)}\left( \beta_j\left( w_{0j}^\eta \right) \right)} \notag \\
&+\xcancel{ \sum_{c=1}^p g_{0j\alpha\beta}^c i_{\alpha_j\left( T_{0j}^c \right)} \left( \beta_j\left( w_{0j}^\gamma \right) \right)} \notag \\
&=0 \notag
\end{align}
This implies $\beta_2\left(R_j \right) + \alpha_2\left(S_j \right) - \tilde{E}_1\left(P_j \right)=0$. Next we show that $P_j- P_k+ \alpha_1\left( F_{jk} \right) + \beta_1\left( G_{jk} \right)=0$. In fact, we note that
{\small{\begin{align}\label{os20}
F_{jk}\left( T_{0j}^\alpha \right)&= -\alpha_j\left( A_{jk}^\alpha \right) - \alpha_k \left(T_{0j}^\alpha \right) +\left[ \theta_{jk}, \alpha_j\left(T_{0j}^\alpha \right) \right] +\left[ \theta_{jk}, \alpha_k\left( T_{0j}^\alpha \right) \right] \\
&= -2 \alpha_j\left( A_{jk}^\alpha \right) + \sum_{\beta=1}^q A_{jk}^{\alpha\beta} \left( \mathcal{L}_{\theta_{jk}}  - A_{jk}    \right) \left( T_{0j}^\beta \right) +2 \left[  \theta_{jk}, \alpha_k\left( T_{0j}^\alpha \right)  \right] + \left[ \theta_{jk}, \left[ \theta_{jk}, T_{0j}^\alpha \right] \right]  - \left[ \theta_{jk}, A_{jk}^\alpha \right] \notag \\
&=-2 \alpha_j\left( A_{jk}^\alpha \right)  +2 \left[  \theta_{jk}, \alpha_k\left( T_{0j}^\alpha \right)  \right] + \left[ \theta_{jk}, \left[ \theta_{jk}, T_{0j}^\alpha \right] \right] + \Theta_{\mathcal{F}_0} \notag\\
G_{jk}\left(w_{0j}^\beta \right)&= - \beta_j\left( B_{jk}^\beta \right)- \beta_k\left( B_{jk}^\beta \right) + \mathcal{L}_{\theta_{jk}}\left(\beta_j\left(w_{0j}^\beta \right) \right) + \mathcal{L}_{\theta_{jk}}\left( \beta_k \left( w_{0j}^\beta \right) \right) \notag \\
&=-2 \beta_k\left(B_{jk}^\beta \right) - \sum_{\gamma=1}^q B_{jk}^{\beta \gamma}\left( \mathcal{L}_{\theta_{jk}}- B_{jk} \right)\left( w_{0j}^\gamma \right) +2 \mathcal{L}_{\theta_{jk}}\left( \beta_j\left( w_{0j}^\beta \right) \right) - \mathcal{L}_{\theta_{jk}}\mathcal{L}_{\theta_{jk}}\left( w_{0j}^\beta \right) + \mathcal{L}_{\theta_{jk}}\left( B_{jk}^\beta \right) \notag \\
&= -2 \beta_k\left(B_{jk}^\beta \right) +2 \mathcal{L}_{\theta_{jk}}\left(\beta_j\left( w_{0j}^\beta \right) \right) - \mathcal{L}_{\theta_{jk}}\mathcal{L}_{\theta_{jk}}\left( w_{0j}^\beta \right) + \mathcal{N}_{\mathcal{F}_0}^* \notag
\end{align}}}
and
{\small{\begin{align}\label{os21}
&2 \mathcal{L}_{\theta_{jk}}i_{\alpha_k\left(T_{0j}^\alpha\right)}\left(w_{0j}^\beta\right) +2\mathcal{L}_{\theta_{jk}}i_{T_{0j}^\alpha}\left( \beta_j\left(w_{0j}^\beta\right) \right) = -2 \mathcal{L}_{\theta_{jk}}i_{T_{0j}^\alpha}\left(\beta_k\left( w_{0j}^\beta\right) - \beta_j\left( w_{0j}^\beta \right) \right) = 2\mathcal{L}_{\theta_{jk}} i_{T_{0j}^\alpha}\mathcal{L}_{\theta_{jk}}\left(w_{0j}^\beta \right) \\
&= - 2\mathcal{L}_{\theta_{jk}}i_{\left[\theta_{jk}, T_{0j}^\alpha \right]}\left(w_{0j}^\beta \right) = -2 \mathcal{L}_{\left[ \theta_{jk}, \left[ \theta_{jk}, T_{0j}^\alpha \right] \right]}\left(w_{0j}^\beta \right)  - 2  i_{\left[\theta_{jk}, T_{0j}^\alpha \right]}\mathcal{L}_{\theta_{jk}}\left(w_{0j}^\beta \right) \notag \\
&= -2 \mathcal{L}_{\left[ \theta_{jk}, \left[ \theta_{jk}, T_{0j}^\alpha \right] \right]}\left(w_{0j}^\beta \right) -2 \mathcal{L}_{\theta_{jk}}i_{\theta_{jk}}\mathcal{L}_{\theta_{jk}}\left(w_{0j}^\beta \right) + 2 i_{T_{0j}^\alpha} \mathcal{L}_{\theta_{jk}}\mathcal{L}_{\theta_{jk}}\left(w_{0j}^\beta \right) =  - \mathcal{L}_{\left[ \theta_{jk}, \left[ \theta_{jk}, T_{0j}^\alpha \right] \right]}\left(w_{0j}^\beta \right)  +  i_{T_{0j}^\alpha} \mathcal{L}_{\theta_{jk}}\mathcal{L}_{\theta_{jk}}\left(w_{0j}^\beta \right) \notag
\end{align}}}
Then from $(\ref{os20})$ and $(\ref{os21})$, we have
{\small{\begin{align*}
&2i_{\alpha_j\left(T_{0j}^\alpha \right)} \left(\beta_j\left( w_{0j}^\beta \right) \right)=2 i_{\left(\alpha_k\left(T_{0j}^\alpha\right) +\left[\theta_{jk}, T_j^\alpha  \right] - A_{jk}^\alpha   \right)}\left(\beta_k\left(w_{0j}^\beta \right) +\mathcal{L}_{\theta_{jk}}\left(w_j^\beta \right)- B_{jk}^\beta \right)\\
&=2\sum_{\delta=1}^p \sum_{\eta=1}^q r_{0jk}^{\alpha \delta} h_{0jk}^{\beta \eta} i_{\alpha_k\left(T_{0k}^\delta  \right)} \beta_k \left(w_{0k}^\eta \right) +2i_{\alpha_k\left(T_{0j}^\alpha\right)}\mathcal{L}_{\theta_{jk}}\left(w_{0j}^\beta\right) -2 i_{\alpha_k\left(T_{0j}^\alpha\right)}\left( B_{jk}^\beta \right)+2i_{\left[\theta_{jk}, T_j^\alpha \right]}\left( \beta_j\left(w_{0j}^\beta\right) \right)-2 i_{A_{jk}^\alpha}\left(\beta_j\left(w_{0j}^\beta\right)\right)\\
&=2\sum_{\delta=1}^p \sum_{\eta=1}^q r_{0jk}^{\alpha \delta} h_{0jk}^{\beta \eta} i_{\alpha_k\left(T_{0k}^\delta  \right)} \beta_k \left(w_{0k}^\eta \right) - 2i_{\left[\theta_{jk}, \alpha_k\left(T_{0j}^\alpha\right) \right]}\left(w_{0j}^\beta\right) + 2 \mathcal{L}_{\theta_{jk}}i_{\alpha_k\left(T_{0j}^\alpha\right)}\left(w_{0j}^\beta\right) + 2i_{T_{0j}^\alpha}\left(\beta_k\left( B_{jk}^\beta\right) \right) \\
&+2\mathcal{L}_{\theta_{jk}}i_{T_{0j}^\alpha}\left( \beta_j\left(w_{0j}^\beta\right) \right) -2 i_{T_{0j}^\alpha}\mathcal{L}_{\theta_{jk}}\left(\beta_j\left(w_{0j}^\beta\right) \right)+2  i_{\alpha_j\left(A_{jk}^\alpha \right) }\left(w_{0j}^\beta \right) \\
&= 2\sum_{\delta=1}^p \sum_{\eta=1}^q r_{0jk}^{\alpha \delta} h_{0jk}^{\beta \eta} i_{\alpha_k\left(T_{0k}^\delta  \right)} \beta_k \left(w_{0k}^\eta \right)  - i_{T_{0j}^\alpha}\left( G_{jk}\left(w_{0j}^\beta \right) \right) - i_{F_{jk}\left(T_{0j}^\alpha \right)} \left( w_{0j}^\beta \right)
\end{align*}}}
This implies that $P_j-P_k + \alpha_1\left(G_{jk} \right) + \beta_1\left(F_{jk} \right)=0$. This completes the proof of Lemma \ref{us24}.

\end{proof}

From $(\ref{its12})$ and Lemma \ref{us24}, we obtain
\begin{theorem}
Suppose given a compact foliated complex manifold $\left(M, \Theta_{\mathcal{F}_0}, \mathcal{N}_{\mathcal{F}_0}^* \right)$ with $\Theta_{\mathcal{F}_0}$ and $\mathcal{N}_{\mathcal{F}_0}^*$ locally free, and $\left(\left\{ \alpha_j \right\}, \left\{ \beta_j \right\}, \left\{\left( \theta_{jk}, A_{jk}, B_{jk} \right) \right\} \right)\in \mathbb{H}^1\left( M, \mathcal{E}_{\mathcal{F}_0}^\bullet \right)$. In order that there may exist a foliated analytic family $\left( \mathcal{M}, \Theta_\mathcal{F}, \mathcal{N}_\mathcal{F}^*, B, \pi \right)$ in terms of both tangent and cotangent sheaves such that $\pi^{-1}(0)=\left( M, \Theta_{\mathcal{F}_0}, \mathcal{N}_{\mathcal{F}_0}^* \right)$, and that $\left( \frac{\partial\left( \mathcal{M}, \Theta_\mathcal{F}, \mathcal{N}_\mathcal{F}^* \right)}{\partial t} \right)|_{t=0}= \left( \left\{ \alpha_j \right\}, \left\{ \beta_j \right\}, \left\{ \left( \theta_{jk} , A_{jk}, B_{jk} \right)\right\} \right)$, it is necessary that
{\small{\begin{align*}
&\left(\left\{ R_j \right\}, \left\{ S_j \right\}, \left\{ P_j \right\} ,  \left\{ F_{jk} \right\}, \left\{ G_{jk} \right\}, \left\{ \left[ \left(\theta_{ij}, A_{ij}, B_{ij} \right), \left( \theta_{jk}, A_{jk}, B_{jk} \right) \right] \right\} \right)\\
& \in C^0\left(\mathcal{U}, \mathscr{H}om_{\mathcal{O}_M}\left( \bigwedge^2 \Theta_{\mathcal{F}_0}, \frac{\Theta_M}{\Theta_{\mathcal{F}_0}} \right) \right) \bigoplus C^0\left(\mathcal{U}, \mathscr{H}om_{\mathcal{O}_M}\left(\mathcal{N}_{\mathcal{F}_0}^*, \tilde{\mathcal{S}}^2 \right) \right) \bigoplus C^0 \left( \mathcal{U}, \mathscr{H}om_{\mathcal{O}_M}\left( \mathcal{N}_{\mathcal{F}_0}^*, \Theta_{\mathcal{F}_0}^* \right) \right)  \\
& \bigoplus C^1\left(\mathcal{U}, \mathscr{H}om_{\mathcal{O}_M}\left(\Theta_{\mathcal{F}_0}, \Theta_M \right) \right) \bigoplus C^1\left(\mathcal{U}, \mathscr{H}om_{\mathcal{O}_M}\left(\mathcal{N}_{\mathcal{F}_0}^*, \Omega_M^1 \right) \right) \bigoplus C^2\left(\mathcal{U}, \mathcal{E}_{\mathcal{F}_0} \right) 
\end{align*}}}
from $(\ref{its39})$ in \textnormal{Lemma \ref{us24}} defines the $0$ class in $\mathbb{H}^2\left( M, \mathcal{E}_{\mathcal{F}_0}^\bullet \right)$.
\end{theorem}
Given $\left(\left\{ \alpha_j \right\}, \left\{ \beta_j \right\}, \left\{\left( \theta_{jk}, A_{jk}, B_{jk} \right) \right\} \right)\in \mathbb{H}^1\left( M, \mathcal{E}_{\mathcal{F}_0}^\bullet \right)$, if the $2$-cocycle $\left(\left\{ R_j \right\}, \left\{ S_j \right\}, \left\{ P_j \right\} ,  \left\{ F_{jk} \right\}, \left\{ G_{jk} \right\}\right.$, 

$\left. \left\{ \left[ \left(\theta_{ij}, A_{ij}, B_{ij} \right), \left( \theta_{jk}, A_{jk}, B_{jk} \right) \right] \right\} \right)$ is not $0$ in $\mathbb{H}^2\left( M, \mathcal{E}_{\mathcal{F}_0}^\bullet \right)$, there is no foliated deformation $\left( M_t, \Theta_{\mathcal{F}_t}, \mathcal{N}_{\mathcal{F}_t}^*\right)$ in terms of both tangent and cotangent sheaves with $\left(M_0, \Theta_{\mathcal{F}_0}, \mathcal{N}_{\mathcal{F}_0}^* \right)=\left( M, \Theta_{\mathcal{F}_0}, \mathcal{N}_{\mathcal{F}_0}^* \right)$ and $\left( \frac{\partial\left( \mathcal{M}, \Theta_\mathcal{F}, \mathcal{N}_\mathcal{F}^* \right)}{\partial t} \right)|_{t=0}=\left( \left\{\alpha_j \right\}, \left\{\beta_j \right\}, \left\{\left( \theta_{jk} , A_{jk}, B_{jk} \right) \right\}\right)$

\section{Infinitesimal foliated deformations and obstructions in terms of both tangent and cotangent sheave by using Dolbeault resolution of $\mathcal{E}_{\mathcal{F}_0}^\bullet$}\label{itt4}

Let $\left(\mathcal{M}, \Theta_\mathcal{F}, \mathcal{N}_{\mathcal{F}}^*, B, \pi   \right)$ with $\Theta_\mathcal{F}$ and $\mathcal{N}_\mathcal{F}^*$ locally free be a foliated analytic family in terms of both tangent and cotangent sheaves of deformations of $\left(M, \Theta_{\mathcal{F}_0}, \mathcal{N}_{\mathcal{F}_0}^* \right)=\pi^{-1}(0)$. We keep the notations in subsection \ref{oi12}- \ref{it20} and subsection \ref{oci80}- \ref{oci81}. Then it remains to consider the condition
\begin{align}\label{it15}
i_{T_j^\alpha(\xi_j,t)}\left( w_j^\beta(\xi_j,t) \right) =0 ,\,\,\,\,\,\,\,\alpha=1,...,p, \beta=1,...,q
\end{align}
We note that $(\ref{it15})$ is equivalent to
\begin{align*}
0=\sum_{\eta=1}^n T_j^{\alpha\eta}\left(\xi_j(z,t),t\right) w_j^{\beta \eta}\left(\xi_j(z,t),t \right)= \sum_{\eta,\gamma=1}^n V_j^{\alpha \gamma}(z,t)\frac{\partial \xi_j^\eta}{\partial z^\gamma} w_j^{\beta \eta}\left(\xi_j(z,t),t \right) = \sum_{\gamma=1}^n V_j^{\alpha \gamma}(z,t) W_j^{\beta \gamma}(z,t),
\end{align*}
so that we obtain
\begin{align*}
i_{V_j^\alpha(z,t)}\left(W_j^\alpha(z,t) \right)=0,\,\,\,\,\,\alpha=1,...,p, \beta=1,...,q
\end{align*}
This implies that
\begin{align}\label{osi8}
i_{X_j^\alpha(t)}\left(Y_j^\beta (t) \right)=i_{\sum_{\delta=1}^p \left(c_j^{-1}\right)^{\alpha \delta} V_j^\delta} \left( \sum_{\eta=1}^q \left(g_j^{-1}\right)^{\beta \eta} W_j^\eta    \right) = \sum_{\delta=1}^p \sum_{\eta=1}^q \left(c_j^{-1}\right)^{\alpha \delta} \left(g_j^{-1}\right)^{\beta \eta} i_{V_j^\delta} \left(W_j^\eta\right)=0
\end{align}

Then from $(\ref{oi26})-(\ref{oi29})$ and $(\ref{oci6})-(\ref{oci7})$ and $(\ref{osi8})$, we have the integrability conditions:
\begin{align}
&\bar{\partial} \varphi - \frac{1}{2}\left[\varphi, \varphi \right]=0 \\
&\left[ X_j^\alpha, X_j^\beta\right]= \sum_{\gamma=1}^p K_{j\alpha\beta}^\gamma X_j^\gamma \\
&\bar{\partial} X_j^\alpha - \left[ \varphi, X_j^\alpha \right]+ \sum_{\eta=1}^p U_j^{\alpha \eta} X_j^\eta=0 \\
&\bar{\partial} U_j^{\alpha \delta}- \left[ \varphi, U_j^{\alpha \delta} \right]  + \sum_{\sigma=1}^p U_j^{\alpha \sigma} \wedge U_j^{\sigma \delta} \\
&Y_j^1\wedge \cdots \wedge Y_j^q \wedge \partial Y_j^\alpha =0\\
&\bar{\partial} Y_j^\alpha- \mathcal{L}_{\varphi}(Y_j^\alpha) + \sum_{\eta=1}^q Z_j^{\alpha \eta} Y_j^\eta =0  \\
&\bar{\partial} Z_j^{\alpha \delta}- \left[ \varphi, Z_j^{\alpha \delta} \right]  + \sum_{\sigma=1}^p Z_j^{\alpha \sigma} \wedge Z_j^{\sigma \delta} =0  \\
& i_{X_j^\alpha(t)}\left( Y_j^\beta(t) \right) = 0 \label{it16}
\end{align}

By taking the derivative of $(\ref{it16})$ with respect to $t$, we have
\begin{align}\label{osi3}
i_{\frac{\partial X_j^\alpha }{\partial t}}\left( Y_j^\beta \right) + i_{X_j^\alpha }\left( \frac{\partial Y_j^\beta}{\partial t}\right)=0,\,\,\,\,\,\alpha=1,...,p, \beta=1,...,q
\end{align}
By setting $t=0$, we obtain
\begin{align}\label{it17}
i_{\frac{\partial X_j}{\partial t}|_{t=0}}\left( w_{0j}^\beta \right) + i_{T_{0j}^\alpha} \left( \frac{\partial Y_j^\beta}{\partial t } |_{t=0} \right) = 0
\end{align}

We keep the notations in subsection \ref{it23} and subsection \ref{oci83}. Then $(\ref{oi42})-(\ref{oi46}), (\ref{it21})-(\ref{it22}), (\ref{it17})$ implies that 
{\Small{\begin{align*}
\left(\frac{\partial X}{\partial t}|_{t=0}, \frac{\partial Y}{\partial t}|_{t=0}, \left\{\left(\frac{\partial \varphi}{\partial t}|_{t=0}, \frac{\partial U_j}{\partial t}|_{t=0}, \frac{\partial Z_j}{\partial t}|_{t=0} \right) \right\} \right)\in A^{0,0}\left( M, \mathscr{H}om_{\mathcal{O}_M}\left( \Theta_{\mathcal{F}_0}, \Theta_M \right)\right) \bigoplus A^{0,0}\left( M, \mathscr{H}om_{\mathcal{O}_M}\left( \mathcal{N}_{\mathcal{F}_0}^*, \Omega_M^1 \right) \right) \bigoplus A^{0,1}\left( M, \mathcal{E}_{\mathcal{F}_0} \right) 
\end{align*}}}
defines a $1$-cocycle in the following bicomplex
{\Tiny{\begin{equation}\label{it25}
\begin{CD}
\cdots\\
@AAA\\
\frac{A^{0,0}\left( M, \left(\bigwedge^3 \Theta_{\mathcal{F}_0}\right) \otimes \Theta_M \right)}{A^{0,0}\left(M, \left(\bigwedge^3 \Theta_{\mathcal{F}_0} \right) \otimes \Theta_{\mathcal{F}_0} \right)}\bigoplus A^{0,0}\left(M-S, \left(\mathcal{N}_{\mathcal{F}_0}^* \right)^* \otimes \tilde{\mathcal{S}}^3 \right)\bigoplus A^{0,0}\left( M-S, \left( \mathcal{N}_{\mathcal{F}_0}^* \right)^*\otimes \bigwedge^2 \Theta_{\mathcal{F}_0} \right) @>-\bar{\partial}>> \cdots \\
@AAA @AAA  \\
\frac{A^{0,0}\left( M, \left(\bigwedge^2 \Theta_{\mathcal{F}_0}\right) \otimes \Theta_M \right)}{A^{0,0}\left(M, \left(\bigwedge^2 \Theta_{\mathcal{F}_0} \right) \otimes \Theta_{\mathcal{F}_0} \right)}\bigoplus A^{0,0}\left(M-S, \left(\mathcal{N}_{\mathcal{F}_0}^* \right)^* \otimes \tilde{\mathcal{S}}^2 \right)\bigoplus A^{0,0}\left( M, \left( \mathcal{N}_{\mathcal{F}_0}^* \right)^*\otimes \Theta_{\mathcal{F}_0} \right)@>\bar{\partial}>> \cdots \\
@AAA @AAA \\
A^{0,0}\left( M, \Theta_{\mathcal{F}_0}^*\otimes \Theta_M \right) \bigoplus A^{0,0}\left( M, \left( \mathcal{N}_{\mathcal{F}_0}^* \right)^* \otimes \Omega_M^1 \right) @> -\bar{\partial} >> A^{0,1}\left( M, \Theta_{\mathcal{F}_0}^*\otimes \Theta_M \right) \bigoplus A^{0,1}\left( M, \left( \mathcal{N}_{\mathcal{F}_0}^* \right)^* \otimes \Omega_M^1 \right) @>\bar{\partial}>>  \\
@AAA @AAA  \\
A^{0,0}\left( M, \mathcal{E}_{\mathcal{F}_0} \right) @>\bar{\partial}>> A^{0,1}\left( M, \mathcal{E}_{\mathcal{F}_0} \right) @>-\bar{\partial}>> 
\end{CD}
\end{equation}}}

We will denote the $i$-th cohomology group of the above double complex by $\textnormal{S}^i$. In the next section \ref{nii7} we will show that we have
\begin{align*}
\textnormal{S}^i \cong \mathbb{H}^i\left( M, \mathcal{E}_{\mathcal{F}_0}^\bullet \right),\,\,\,\,\,\textnormal{for}\,\,\,  i = 0,1,2
\end{align*}

Then we have the foliated Kodaira-Spencer map in terms of both tangent and cotangent sheaves
{\small{\begin{align}
\varphi_0: T_0 B \to & \mathbb{H}^1\left( M, \mathcal{E}_{\mathcal{F}_0}^\bullet \right)  \label{a17}\\
 \frac{\partial}{\partial t} \mapsto &\left( \frac{\partial X}{\partial t}|_{t=0}, \frac{\partial Y}{\partial t}|_{t=0}, \left\{\left(\frac{\partial \varphi}{\partial t}|_{t=0}, \frac{\partial U_j}{\partial t}|_{t=0}, \frac{\partial Z_j}{\partial t}|_{t=0} \right) \right\} \right) \notag \\
 &\in  A^{0,0}\left( M, \mathscr{H}om_{\mathcal{O}_M}\left( \Theta_{\mathcal{F}_0}, \Theta_M \right)\right) \bigoplus A^{0,0}\left( M, \mathscr{H}om_{\mathcal{O}_M}\left( \mathcal{N}_{\mathcal{F}_0}^*, \Omega_M^1 \right) \right) \bigoplus A^{0,1}\left( M, \mathcal{E}_{\mathcal{F}_0} \right) \notag
\end{align}}}

\subsection{Identification of obstructions by using Dolbeault resolution of $\mathcal{E}_{\mathcal{F}_0}^\bullet$ } \label{itt10} \

Now assume that the foliated analytic family $\left( \mathcal{M}, \Theta_\mathcal{F}, \mathcal{N}_\mathcal{F}^*, B, \pi \right)$ is a one-parameter family, i.e. $\Delta \subset \mathbb{C}$. We keep the notations in subsection \ref{it30} and subsection \ref{oci84}. Then by taking the derivative of $(\ref{osi3})$ with respect to $t$ again and setting $t=0$, we have
\begin{align*}
i_{\frac{\partial^2 X}{\partial t^2}|_{t=0}}\left(w_{0j}^\beta \right) + 2 i_{\frac{\partial X_j^\alpha}{\partial t}|_{t=0}}\left( \frac{\partial Y_j^\alpha}{\partial t}|_{t=0} \right) + i_{T_{0j}^\alpha}\left(  \frac{\partial^2 Y_j^\beta}{\partial t^2}|_{t=0} \right) = 0 
\end{align*}

Then we define an element $ \frac{\partial X_j}{\partial t}|_{t=0} * \frac{\partial Y_j}{\partial t}|_{t=0}\in \Gamma\left( U_j, \mathcal{A}^{0,0}\left(\mathscr{H}om_{\mathcal{O}_M}\left( \mathcal{N}_{\mathcal{F}_0}^*, \Theta_{\mathcal{F}_0}^* \right) \right) \right)$ 
\begin{align*}
 \frac{\partial X_j}{\partial t}|_{t=0} * \frac{\partial Y_j}{\partial t}|_{t=0}:\Gamma\left( U_j,  \mathcal{N}_{\mathcal{F}_0}^* \right) & \to \Gamma\left( U_j, \mathcal{A}^{0,0}\left(\Theta_{\mathcal{F}_0}^* \right) \right) \\
 w_{0j}^\beta  & \mapsto \left(  T_{0j}^\alpha \mapsto -2 i_{\frac{\partial X_j^\alpha}{\partial t}|_{t=0}} \left( \frac{\partial Y_j^\beta}{\partial t}|_{t=0} \right)    \right) 
\end{align*}
Then $\left\{ \frac{\partial X_j}{\partial t}|_{t=0} * \frac{\partial Y_j}{\partial t} |_{t=0} \right\}$ defines a global section $\frac{\partial X}{\partial t}|_{t=0} * \frac{\partial Y}{\partial t}|_{t=0}\in A^{0,0} \left( M, \mathscr{H}om_{\mathcal{O}_M}\left( \mathcal{N}_{\mathcal{F}_0}^*, \Theta_{\mathcal{F}_0}^* \right) \right)$ and
\begin{align} \label{it27}
i_{\frac{\partial^2 X}{\partial t^2}|_{t=0}}\left( w_{0j}^\beta \right)+ i_{T_{0j}^\alpha}\left( \frac{\partial^2 Y}{\partial t^2}|_{t=0} \right)=\frac{\partial X}{\partial t}|_{t=0} * \frac{\partial Y}{\partial t}|_{t=0}\left( w_{0j}^\beta, T_{0j}^\alpha \right)
\end{align}

\begin{lemma}\label{its28}
We keep the notations in \textnormal{Lemmas \ref{its30}} and \textnormal{\ref{oci85}}. For any $1$-cocycle $\left( \Phi=\left\{ \Phi_j\right\}, \Psi=\left\{ \Psi_j\right\}, \left\{ \psi, a_j, k_j \right\} \right)\in A^{0,0}\left( M, \mathscr{H}om_{\mathcal{O}_M}\left( \Theta_{\mathcal{F}_0}, \Theta_M \right)\right) \bigoplus A^{0,0}\left( M, \mathscr{H}om_{\mathcal{O}_M}\left( \mathcal{N}_{\mathcal{F}_0}^*, \Omega_M^1 \right) \right) \bigoplus A^{0,1}\left( M,  \mathcal{E}_{\mathcal{F}_0} \right)$ as in \textnormal{Lemma \ref{its30} and Lemma \ref{oci85}} satisfying
\begin{enumerate}
\item $\left[ T_{0j}^\alpha, \Phi_j\left( T_{0j}^\beta \right) \right]- \left[ T_{0j}^\beta, \Phi_j\left( T_{0j}^\alpha \right) \right]- \sum_{\gamma=1}^p g_{0j\alpha\beta}^\gamma \Phi_j\left( T_{0j}^\gamma \right) = K_{j\alpha\beta} $ for some $K_{j\alpha\beta}\in \Gamma\left( U_j, \mathcal{A}^{0,0}\left( \Theta_{\mathcal{F}_0} \right) \right),\alpha,\beta=1,...,p$.
\item $- \bar{\partial} \Phi_j^\alpha + \left[ \psi, T_{0j}^\alpha \right]- \sum_{\beta=1}^p a_j^{\alpha\beta} T_{0j}^\beta=0,\alpha=1,...,p$.
\item $\sum_{\beta=1}^q w_{0j}^1 \wedge \cdots \wedge \Psi_j^\beta \wedge \cdots \wedge  w_{0j}^q \wedge dw_{0j}^\alpha + w_{0j}^1 \wedge \cdots \wedge w_{0j}^q\wedge \partial \Psi_j^\alpha=0,\alpha=1,...,q$,
\item $- \bar{\partial} \Psi_j^\alpha + \mathcal{L}_\psi \left( w_{0j}^\alpha \right) - \sum_{\beta=1}^q k_j^{\alpha\beta} w_{0j}^\beta,\alpha=1,...,q$.
\item $- \bar{\partial} \psi=0, -\bar{\partial} a_j^{\alpha\beta}=0,\alpha,\beta=1,...,p$, and $- \bar{\partial} k_j^{\gamma \eta}=0, \gamma,\eta=1,...,q$.
\item $i_{\Phi_j^\alpha}\left( w_{0j}^\beta \right) + i_{T_{0j}^\alpha}\left( \Psi_j^\beta \right)=0,\alpha=1,...,p, \beta=1,...,q$.
\end{enumerate}
We define $\tilde{\rho}=\left\{\tilde{\rho}_i \right\}\in A^{0,0}\left( M, \mathscr{H}om_{\mathcal{O}_M}\left( \bigwedge^2 \Theta_{\mathcal{F}_0}, \Theta_M \right) \right)$, and define $\rho$ to be the image of $\frac{A^{0,0}\left( M, \mathscr{H}om_{\mathcal{O}_M}\left( \bigwedge^2 \Theta_{\mathcal{F}_0}, \Theta_M \right) \right)}{A^{0,0}\left( M, \mathscr{H}om_{\mathcal{O}_M}\left( \bigwedge^2 \Theta_{\mathcal{F}_0}, \Theta_{\mathcal{F}_0} \right) \right)}$, and define $\eta=\{\eta_i\}\in A^{0,0}\left( M, \mathscr{H}om_{\mathcal{O}_M}\left( \Theta_{\mathcal{F}_0}, \Theta_M \right)\right)$ as in \textnormal{Lemma \ref{its30}}. We define $l\in A^{0,0}\left( M-S, \mathscr{H}om_{\mathcal{O}_M}\left( \mathcal{N}_{\mathcal{F}_0}^*, \tilde{\mathcal{S}}^2 \right)\right)$ and $\kappa=\left\{ \kappa_i\right\}  \in A^{0,1}\left( M, \mathscr{H}om_{\mathcal{O}_M}\left( \mathcal{N}_{\mathcal{F}_0}^*, \Omega_M^1 \right) \right)$ as in \textnormal{Lemma \ref{oci85}}. On the other hand, we define an element $p_j\in \Gamma\left( U_j , \mathcal{A}^{0,1}\left( \mathscr{H}om_{\mathcal{O}_M}\left( \mathcal{N}_{\mathcal{F}_0}^*, \Theta_{\mathcal{F}_0}^* \right) \right)\right)$ by
\begin{align*}
p_j:\Gamma\left( U_j, \mathcal{N}_{\mathcal{F}_0}^* \right) &\to \Gamma\left( U_j, \mathcal{A}^{0,0}\left(\Theta_{\mathcal{F}_0}^* \right) \right)\\
      w_{0j}^\beta &\mapsto \left( T_{0j}^\alpha \mapsto  - 2 i_{\Phi_j^\alpha}\left(\Psi_j^\beta \right)\right)
\end{align*}
Then $p=\left\{ p_j \right\}$ defines a global section in $A^{0,0}\left( M, \mathscr{H}om_{\mathcal{O}_M}\left( \mathcal{N}_{\mathcal{F}_0}^*, \Theta_{\mathcal{F}_0}^* \right) \right)$, and we note that $\left[ \left( \psi, a_j, k_j \right), \left(\psi, a_j, k_j \right) \right]\in A^{0,2}\left( M, \mathcal{E}_{\mathcal{F}_0} \right)$ from \textnormal{Lemma \ref{its30} and Lemma \ref{oci85}}. Then
\begin{align}\label{it33}
&\left(\rho, l, \eta,  \kappa, p, \left\{ -\left[\left( \psi, a_j, k_j \right), \left( \psi, a_j, k_j \right) \right] \right\}\right)\\
&\in  \frac{A^{0,0}\left( M, \mathscr{H}om_{\mathcal{O}_M}\left( \bigwedge^2 \Theta_{\mathcal{F}_0}, \Theta_M \right) \right)}{A^{0,0}\left( M, \mathscr{H}om_{\mathcal{O}_M}\left( \bigwedge^2 \Theta_{\mathcal{F}_0}, \Theta_{\mathcal{F}_0} \right)\right)} \bigoplus A^{0,0}\left(M-S, \mathscr{H}om_{\mathcal{O}_M}\left( \mathcal{N}_{\mathcal{F}_0}^*, \tilde{\mathcal{S}}^2 \right) \right) \bigoplus A^{0,0}\left(M, \mathscr{H}om_{\mathcal{O}_M}\left( \mathcal{N}_{\mathcal{F}_0}^*, \Theta_{\mathcal{F}_0}^* \right) \right) \notag \\
&\bigoplus A^{0,1}\left( \mathscr{H}om_{\mathcal{O}_M}\left( \Theta_{\mathcal{F}_0}, \Theta_M \right) \right) \bigoplus A^{0,1}\left( M, \mathscr{H}om_{\mathcal{O}_M}\left( \mathcal{N}_{\mathcal{F}_0}^*, \Omega_M^1 \right) \right) \bigoplus A^{0,2}\left( M, \mathcal{E}_{\mathcal{F}_0} \right) \notag
\end{align}
defines a $2$-cocycle in the Dolbeault type bicomplex $(\ref{it25})$.
\end{lemma}

\begin{proof}
By Lemma \ref{its30} and Lemma \ref{oci85}, it remains to show that $\beta_2 \left(\rho\right) + \alpha_2\left( l\right) + \hat{E}_1'\left(p_j \right)=0$ and $\bar{\partial} p + \hat{\beta}_1\left(  \eta \right) + \hat{\alpha}_1\left(\kappa \right)=0$. In fact, as in the same way with $(\ref{us26}),(\ref{us27}),(\ref{us25}),(\ref{its20})$, we can show $\hat{\beta}_2 \left(\rho\right) + \hat{\alpha}_2\left( l\right) - \hat{E}_1'\left(p_j \right)=0$ on $M-S$. On the other hand, since $i_{\Phi_j^\alpha}\left(w_{0j}^\beta \right) + i_{T_{0j}^\alpha}\left(\Psi_j^\beta \right)=0$, we have
\begin{align*}
&\bar{\partial}\left(-2 i_{\Phi_j^\alpha}\left(\Psi_j^\beta \right) \right) =-2 i_{\bar{\partial} \Phi_j^\alpha}\left(\Psi_j^\beta \right)+2 i_{\Phi_j^\alpha}\left(\bar{\partial} \Psi_j^\beta  \right) = -2 i_{\left[ \psi, T_{0j}^\alpha \right] - \sum_{\delta=1}^p a_j^{\alpha \delta} \wedge T_{0j}^\delta } \left( \Phi_j^\beta   \right) + 2 i_{\Phi_j^\alpha}\left( \mathcal{L}_{\psi} \left(w_{0j}^\beta \right)- \sum_{\eta=1}^q k_j^{\beta \eta} \wedge w_{0j}^\eta \right)\\
&= -2 \mathcal{L}_\psi i_{T_{0j}^\alpha}\left(\Psi_j^\beta \right) - 2 i_{T_{0j}^\alpha} \mathcal{L}_\psi\left( \Psi_j^\beta \right) + 2\sum_{\delta=1}^p a_j^{\alpha \delta} i_{T_{0j}^\delta}\left( \Psi_j^\beta \right) + 2i_{\Phi_j^\alpha} \mathcal{L}_\psi\left(w_{0j}^\beta \right) +  2 \sum_{\eta=1}^q k_j^{\beta \eta} i_{\Phi_j^\alpha}\left(w_{0j}^\eta \right) \\
&= 2 \mathcal{L}_\psi i_{\Phi_j^\alpha}\left(w_{0j}^\beta \right) - 2 i_{T_{0j}^\alpha} \mathcal{L}_\psi \left( \Psi_j^\beta \right) - 2 \sum_{\delta=1}^p a_j^{\alpha \delta} i_{\Phi_j^\delta }\left( w_{0j}^\beta \right) + 2 i_{\Phi_j^\alpha} \mathcal{L}_\psi\left( w_{0j}^\beta \right) - 2 \sum_{\eta=1}^q k_j^{\beta \eta} i_{T_{0j}^\alpha}\left(\Psi_j^\eta \right)\\
&= 2 i_{\left[ \psi, \Phi_j^\alpha \right] - \sum_{\delta=1}^p a_j^{\alpha \delta} \wedge \Phi_j^\delta }\left( w_{0j}^\beta \right) - 2i_{T_{0j}^\alpha}\left(  \mathcal{L}_\psi \left( \Psi_j^\beta \right) - \sum_{\eta=1}^q k_j^{\beta \eta} \wedge \Psi_j^\eta   \right)
\end{align*}
This implies that $\bar{\partial} p + \hat{\beta}_1\left(  \eta \right) + \hat{\alpha}_1\left(\kappa \right)=0$ on $M-S$. This completes the proof of Lemma \ref{its28}.
\end{proof}

Then from $(\ref{it8})-(\ref{it9})$ and $(\ref{it31})-(\ref{it32})$ and $(\ref{it27})$ and Lemma \ref{its28}, we obtain
\begin{theorem}
Suppose given a compact foliated complex manifold $\left( M,\Theta_{ \mathcal{F}_0}, \mathcal{N}_{\mathcal{F}_0}^* \right)$ with $\Theta_{\mathcal{F}_0}$ and $\mathcal{N}_{\mathcal{F}_0}^*$ locally free, and $\left( \Phi, \Psi, \left\{\left(\psi, a_j, k_j\right) \right\} \right) \in \mathbb{H}^1\left( M, \mathcal{E}_{\mathcal{F}_0}^\bullet \right)$. In order that there may exist a foliated analytic family $\left(\mathcal{M}, \Theta_\mathcal{F}, \mathcal{N}_\mathcal{F}^*, B, \pi \right)$ in terms of both tangent and cotangent sheaves such that $\pi^{-1}(0)=\left( M, \Theta_{\mathcal{F}_0}, \mathcal{N}_{\mathcal{F}_0}^* \right)$ and that $\left( \frac{\partial\left(\mathcal{M}, \Theta_{\mathcal{F}}, \mathcal{N}_{\mathcal{F}}^* \right)}{\partial t} \right)|_{t=0}=\left( \Phi, \Psi, \left\{\left( \psi, a_j, k_j \right) \right\} \right)$, it is necessary that
\begin{align*}
&\left(\rho, l, \eta,  \kappa, p, \left\{ -\left[\left( \psi, a_j, k_j \right), \left( \psi, a_j, k_j \right) \right] \right\}\right)\\
&\in  \frac{A^{0,0}\left( M, \mathscr{H}om_{\mathcal{O}_M}\left( \bigwedge^2 \Theta_{\mathcal{F}_0}, \Theta_M \right) \right)}{A^{0,0}\left( M, \mathscr{H}om_{\mathcal{O}_M}\left( \bigwedge^2 \Theta_{\mathcal{F}_0}, \Theta_{\mathcal{F}_0} \right)\right)} \bigoplus A^{0,0}\left(M-S, \mathscr{H}om_{\mathcal{O}_M}\left( \mathcal{N}_{\mathcal{F}_0}^*, \tilde{\mathcal{S}}^2 \right) \right) \bigoplus A^{0,0}\left(M, \mathscr{H}om_{\mathcal{O}_M}\left( \mathcal{N}_{\mathcal{F}_0}^*, \Theta_{\mathcal{F}_0}^* \right) \right) \notag \\
&\bigoplus A^{0,1}\left( \mathscr{H}om_{\mathcal{O}_M}\left( \Theta_{\mathcal{F}_0}, \Theta_M \right) \right) \bigoplus A^{0,1}\left( M, \mathscr{H}om_{\mathcal{O}_M}\left( \mathcal{N}_{\mathcal{F}_0}^*, \Omega_M^1 \right) \right) \bigoplus A^{0,2}\left( M, \mathcal{E}_{\mathcal{F}_0} \right) \notag
\end{align*}
from $(\ref{it33})$ in \textnormal{Lemma \ref{its28}} defines the $0$ class in $\mathbb{H}^2\left( M, \mathcal{E}_{\mathcal{F}_0}^\bullet\right)$.
\end{theorem}

\section{Isomorphisms of infinitesimal foliated deformations and obstructions in terms of cotangent sheaves between \v Cech resolution and Dolbeault type bicomplex associated to $\mathcal{E}_{\mathcal{F}_0}^\bullet$}\label{nii7}

In this section we explicitly construct isomorphisms of first and second cohomology groups between \v Cech resolution of $\mathcal{E}_{\mathcal{F}_0}^\bullet$ and Dolbeault resolution of $\mathcal{E}_{\mathcal{F}_0}^\bullet$. We keep the  notations in section \ref{ti2} and \ref{ni2}. It remains to consider the condition: $i_{T_{0j}^\alpha}\left(w_{0j}^\beta  \right)=0$ for $\alpha=1,...,p, \beta=1,...,q$ on $U_j$.

\subsection{Isomorphisms of infinitesimal deformations between \v Cech resolution and Dolbeault type bicomplex associated to $\mathcal{E}_{\mathcal{F}_0}^\bullet$ } \label{a25} \

We keep the notations in subsection \ref{ti1} and \ref{ni1}. Then we define the correspondence of $1$-cocycle between \v Cech resolution and Dolbeault type bicomplex of $\mathcal{E}_{\mathcal{F}_0}^\bullet$ by 
\begin{align}\label{a18}
&\left( \left\{\alpha_j \right\},  \left\{ \beta_j\right\}, \left\{ \left( \theta_{jk}, A_{jk}, B_{jk} \right)      \right\}\right) \in C^0\left( \mathcal{U}, \mathscr{H}om_{\mathcal{O}_M}\left( \Theta_{\mathcal{F}_0}, \Theta_M \right)\right) \bigoplus  C^0\left( \mathcal{U}, \mathscr{H}om_{\mathcal{O}_M}\left( \mathcal{N}_{\mathcal{F}_0}^*, \Omega_M^1 \right) \right) \bigoplus C^1\left( \mathcal{U}, \mathcal{E}_{\mathcal{F}_0}  \right)\\
&\mapsto \left( \Phi, \Psi, \left\{ \left( \psi, \bar{\partial} A_j, \bar{\partial} B_j \right)  \right\} \right) \in \textnormal{S}^1 \notag
\end{align}
In fact, it remains to show that since $i_{T_{0j}^\alpha}\left( \beta_j\left(w_{0j}^\beta \right)\right) + i_{\alpha_j\left( T_{0j}^\alpha \right)}\left( w_{0j}^\beta \right)=0$, we have
\begin{align*}
i_{T_{0j}^\alpha}\left(\Psi_j^\beta  \right) + i_{\Phi_j}\left( w_{0j}^\beta \right) &= i_{T_{0j}^\alpha}\left( \mathcal{L}_{\xi_j}\left( w_{0j}^\alpha  \right) - B_j\left(w_{0j}^\alpha  \right) - \beta_j\left(w_{0j}^\alpha \right)\right)  + i_{\left[ \xi_j, T_{0j}^\alpha\right] - A_j\left( T_{0j}^\alpha \right) - \alpha_j\left( T_{0j}^\alpha \right)} \left( w_{0j}^\beta \right)\\
&= i_{T_{0j}^\alpha}\mathcal{L}_{\xi_j}\left(w_{0j}^\beta \right) + \mathcal{L}_{\xi_j} i_{T_{0j}^\alpha } \left(w_{0j}^\beta \right) - i_{ T_{0j}^\alpha} \mathcal{L}_{\xi_j}\left(w_{0j}^\beta \right) = 0
\end{align*}

Then by combining Proposition \ref{b18} and Proposition \ref{c1}, we obtain

\begin{proposition}
We recall that given a foliated analytic family $\left(\mathcal{M}, \Theta_{\mathcal{F}}, \mathcal{N}_\mathcal{F}^*, B, \pi \right)$ of deformations of $\left( M, \Theta_{\mathcal{F}_0}, \mathcal{N}_{\mathcal{F}_0}^* \right)= \pi^{-1}(0)$ as in \textnormal{Remark \ref{os3}}, the infinitesimal deformation at $t=0$ is given by 
\begin{align} \label{a19}
\left( \left\{ \alpha_j(0)\right\}, \left\{ \beta_j(0)\right\}, \left\{ \left( \theta_{jk}(0), A_{jk}(0), B_{jk}(0) \right) \right\} \right)
\end{align}
 from $(\ref{its2})$ in the first cohomology group of the \v Cech resolution of $\mathcal{E}_{\mathcal{F}_0}^\bullet$. On the other hand, the infinitesimal deformation at $t=0$ is also given by 
 \begin{align}\label{a20}
 \left( \frac{\partial X_j}{\partial t}|_{t=0}, \frac{\partial Y_j}{\partial t}|_{t=0},  \left\{\left( \frac{\partial \varphi}{\partial t}|_{t=0}, \frac{\partial U_j}{\partial t}|_{t=0}, \frac{\partial Z_j}{\partial t}|_{t=0} \right) \right\}\right)
 \end{align}
from $(\ref{a17})$ in the first cohomology group $\textnormal{S}^1$ of the Dolbeault type bicomplex associated to $\mathcal{E}_{\mathcal{F}_0}^\bullet$. Then $(\ref{a19})$ corresponds to $(\ref{a20})$ under the isomorphism $(\ref{a18})$.
\end{proposition}

\subsection{Isomorphism of obstructions between \v Cech resolution and Dolbeault type bicomplex associated to $\mathcal{E}_{\mathcal{F}_0}^\bullet$ }\label{a33} \ 

We construct an isomorphism between the second cohomology group $\mathbb{H}^2\left( M, \mathcal{E}_{\mathcal{F}_0}^\bullet \right)$ from the \v Cech resolution and the second cohomology group $\textnormal{S}^2$ from the Doubault type bicomplex associated to $\mathcal{E}_{\mathcal{F}_0}^\bullet$. We keep the notations in subsection \ref{a15} and subsection \ref{a16}. 

Let $\left(\left\{ C_i \right\}, \left\{ c_i \right\}, \left\{ J_{j} \right\}, \left\{ B_{jk}\right\}, \left\{ b_{jk}\right\},  \left\{\left( a_{ijk}, G_{ijk}, g_{ijk} \right) \right\} \right)$ be a $2$-cocycle, where $J_j\in C^0\left(\mathcal{U}, \mathscr{H}om_{\mathcal{O}_M}\left( \mathcal{N}_{\mathcal{F}_0}^*, \Theta_{\mathcal{F}_0}^* \right) \right)$. Then we note that
\begin{align*}
&J_j^{\alpha \beta} + i_{T_{0j}^\alpha }\left( f_j^\beta \right) + i_{F_j^\alpha}\left( w_{0j}^\beta \right) - \sum_{\delta=1}^p \sum_{\eta=1}^q r_{0jk}^{\alpha \delta} h_{0jk}^{\beta \eta} \left( J_k^{\delta \eta}  + i_{T_{0k}^\delta} \left(f_k^\eta \right)  + i_{F_k^\delta}\left( w_{0j}^\eta  \right)   \right) \\
&= J_j^{\alpha \beta}   - \sum_{\delta=1}^p \sum_{\eta=1}^q r_{0jk}^{\alpha \delta} h_{0jk}^{\beta \eta} J_k^{\delta \eta} + i_{T_{0j}^\alpha } \left( \sum_{\eta=1}^q h_{0jk}^{\beta \eta} f_k^\eta - \tilde{\mathcal{L}}_{d_{jk}}\left( w_{0j}^\beta\right) + e_{jk}^\beta + b_{jk}^\beta   \right) + i_{\sum_{\delta=1}^p r_{0jk}^{\alpha \delta} F_k^\delta  - \left[ d_{jk}, T_{0j}^\alpha  \right] + E_{jk}^\alpha + B_{jk}^\alpha}\left(w_{0j}^\beta  \right) \\
&- \sum_{\delta=1}^p \sum_{\eta=1}^q r_{0jk}^{\alpha \delta} h_{0jk}^{\beta \eta} \left(  i_{T_{0k}^\delta} \left(f_k^\eta \right)  + i_{F_k^\delta}\left( w_{0j}^\eta  \right)   \right)=0
\end{align*}

Let us define $I_j\in \Gamma\left( U_j, \mathcal{A}^{0,0}\left( \mathscr{H}om_{\mathcal{O}_M}\left( \mathcal{N}_{\mathcal{F}_0}^*, \Theta_{\mathcal{F}_0}^* \right) \right) \right)$ by
\begin{align*}
I_j: \Gamma\left( U_j, \mathcal{N}_{\mathcal{F}_0}^* \right)&\to \Gamma\left( U_j, \mathcal{A}^{0,0}\left(  \Theta_{\mathcal{F}_0}^* \right) \right) \\
 w_{0j}^\beta &\mapsto \left( T_{0j}^\alpha \mapsto - J_j^{\alpha\beta} - i_{T_{0j}^\alpha}\left( f_j^\beta \right) - i_{F_j^\alpha}\left( w_{0j}^\beta \right) \right)
\end{align*}
Then $I:=\left\{I_j \right\}$ defines a global section in $A^{0,0}\left( M, \mathscr{H}om_{\mathcal{O}_M}\left( \mathcal{N}_{\mathcal{F}_0}^*, \Theta_{\mathcal{F}_0}^* \right) \right)$. Then we define the correspondence of $2$-cocycles between \v Cech resolution of $\mathcal{E}_{\mathcal{F}_0}^\bullet$ and the Dolbeault type bicomplex associated to $\mathcal{E}_{\mathcal{F}_0}^\bullet$ by
\begin{align}\label{a21}
\left(\left\{ C_i \right\}, \left\{ c_i \right\}, \left\{ J_{j} \right\}, \left\{ B_{jk}\right\}, \left\{ b_{jk}\right\},  \left\{\left( a_{ijk}, G_{ijk}, g_{ijk} \right) \right\} \right) \mapsto \left( \overline{P}, R, I,  L, Q,  \left\{ \left( \bar{\partial} p_j, \bar{\partial} U_j, \bar{\partial} u_j \right) \right\} \right)
\end{align}

Conversely, let $\left(P, R, I, L,  Q, \left( \phi,  \left\{\left( M_j ,N_j  \right)\right\} \right) \right)$ be a $2$-cocycle, where $I \in A^{0,0}\left( M, \mathscr{H}om_{\mathcal{O}_M}\left( \mathcal{N}_{\mathcal{F}_0}^*, \Theta_{\mathcal{F}_0}^\bullet \right) \right)$. Then we note that since we have
\begin{align*}
0&=\bar{\partial}I_j\left( w_{0j}^\beta \right)\left( T_{0j}^\alpha\right) - i_{T_{0j}^\alpha}\left( Q_j^\beta \right) + i_{L_j^\alpha}\left(w_{0j}^\beta \right)\\
&=\bar{\partial} I_j^{\alpha\beta} - i_{T_{0j}^\alpha}\left( \bar{\partial} f_j^\beta - \mathcal{L}_{p_j}\left( w_{0j}^\beta \right) + u_j^\beta \right) + i_{\bar{\partial} F_j^\alpha - \left[ p_j, T_{0j}^\alpha\right] + U_j^\alpha } \left(w_{0j}^\beta \right) \\
&= \bar{\partial}\left( I_j^{\alpha\beta} + i_{T_{0j}^\alpha} \left( f_j^\beta \right) + i_{F_j^\alpha} \left( w_{0j}^\beta \right)   \right)
\end{align*}
Then we define $J_j \in \Gamma\left( U_j, \mathscr{H}om_{\mathcal{O}_M}\left( \mathcal{N}_{\mathcal{F}_0}^*,  \Theta_{\mathcal{F}_0}^* \right) \right)$ by
\begin{align*}
J_j: \Gamma\left( U_j, \mathcal{N}_{\mathcal{F}_0}^* \right) &\to \Gamma\left( U_j, \Theta_{\mathcal{F}_0}^* \right) \\
                       w_{0j}^\beta &\mapsto \left( T_{0j}^\alpha \mapsto - I_j^{\alpha\beta} - i_{T_{0j}^\alpha}\left( f_j^\beta \right)  - i_{F_j^\alpha} \left( w_{0j}^\beta \right) \right)
\end{align*}

Then we define the inverse correspondence by 
\begin{align*}
\left( \overline{P}, R,  I,  L, Q ,  \left( \phi, \left\{  \left(M_j , N_j \right)\right\} \right) \right)\mapsto \left(\left\{C_j\right\}, \left\{c_j \right\}, \left\{ J_j \right\}, \left\{ B_{jk}\right\}, \left\{ b_{jk}\right\}, \left\{\left( a_{ijk}, G_{ijk}, g_{ijk} \right) \right\} \right).
\end{align*} 

\begin{proposition}
We recall from \textnormal{Lemma \ref{us24}} that given a $1$-cocycle $\left( \left\{ \alpha_j\right\}, \left\{ \beta_j \right\}, \left\{\left( \theta_{jk}, A_{jk}, B_{jk} \right) \right\} \right)$ in $\mathbb{H}^1\left( M, \mathcal{E}_{\mathcal{F}_0}^\bullet \right)$, we have the associated obstruction class from $(\ref{its39})$
\begin{align}\label{a22}
\left(\left\{ R_j \right\}, \left\{ S_j \right\}, \left\{ P_j \right\} ,  \left\{ F_{jk} \right\}, \left\{ G_{jk} \right\}, \left\{ \left[ \left(\theta_{ij}, A_{ij}, B_{ij} \right), \left( \theta_{jk}, A_{jk}, B_{jk} \right) \right] \right\} \right) \in \mathbb{H}^2\left( M, \mathcal{E}_{\mathcal{F}_0}^\bullet \right)
\end{align}
from the \v Cech resolution of $\mathcal{E}_{\mathcal{F}_0}^\bullet$. Assume that the $1$-cocycle $\left(\left\{ \alpha_j \right\}, \left\{ \beta_j \right\}, \left\{ \left( \theta_{jk}, A_{jk}, B_{jk} \right) \right\} \right)$ from the \v Cech resolution of $\mathcal{E}_{\mathcal{F}_0}^\bullet$ corresponds to the $1$-cocycle $\left( \Phi, \Psi, \left\{\psi, a_j:=\bar{\partial} A_j, k_j:=\bar{\partial} B_j \right\}  \right)$ from the Dolbeault type bicomplex associated to $\mathcal{E}_{\mathcal{F}_0}^\bullet$ under the isomorphism $(\ref{a18})$ from $\mathbb{H}^1\left( M, \mathcal{E}_{\mathcal{F}_0}^\bullet \right)$ to $\textnormal{S}^1$. We recall from \textnormal{Lemma \ref{its28}} that we have the associated obstruction class from $(\ref{it33})$
\begin{align}\label{a23}
\left( \rho, l, p, \eta, \kappa,   \left\{ -\left[ \left(\psi, a_j, k_j \right), \left( \psi, a_j, k_j \right)       \right]    \right\} \right)\in \textnormal{S}^2
\end{align}
from the Dolbeault type bicomplex associated to $\mathcal{E}_{\Theta_{\mathcal{F}_0}}^\bullet$. Then $(\ref{a22})$ corresponds to $(\ref{a23})$ under the isomorphism $(\ref{a21})$ from $\mathbb{H}^2\left( M, \mathcal{E}_{\mathcal{F}_0}^\bullet \right)$ to $\textnormal{S}^2$.
\end{proposition}

\begin{proof}
Let $\left( \left\{ \alpha_j \right\}, \left\{  \beta_j \right\} , \left\{ \left( \theta_{jk}, A_{jk}, B_{jk} \right) \right\} \right)$ be a $1$-cocycle in Lemma \ref{us24}. Then $\left(\left\{ R_j \right\}, \left\{ S_j \right\}, \left\{ P_j \right\} ,  \left\{ F_{jk} \right\}, \left\{ G_{jk} \right\},\right.$ $\left. \left\{ \left[ \left(\theta_{ij}, A_{ij}, B_{ij} \right), \left( \theta_{jk}, A_{jk}, B_{jk} \right) \right] \right\} \right)$ defines a $2$-cocycle by Lemma \ref{us24}. We have the corresponding $1$-cocycle $\left( \Phi , \Psi, \left\{\psi, a_j=\bar{\partial} A_j, k_j=\bar{\partial} B_j \right\} \right)$ from the Dolbeault type bicomplex associated to $\mathcal{E}_{\mathcal{F}_0}^\bullet$ from $(\ref{a18})$. We keep the notations in the proof of Proposition \ref{a13} and the proof of Proposition \ref{a24}. Then it remains to show that $\left\{ P_j \right\}$ corresponds to $p$. In fact,
{\small{\begin{align}\label{iit2}
&P_j\left( T_{0j}^\alpha \right)\left( w_{0j}^\beta\right) + i_{F_j^\alpha }\left( w_{0j}^\beta\right) + i_{T_{0j}^\alpha}\left(f_j^\beta \right)\\
&= 2 i_{\alpha_j\left(  T_{0j}^\alpha \right)} \left(\beta_j\left(w_{0j}^\beta \right) \right)      +  i_{-2\alpha_j\left(A_j^\alpha \right)+ 2 \left[ \xi_j, \alpha_j\left(T_{0j}^\alpha \right)\right] +2 \left[ \xi_j, A_j^\alpha \right] - \left[ \xi_j ,  \left[ \xi_j,  T_{0j}^\alpha \right] \right] } \left(w_{0j}^\beta \right)  \notag   \\
&+ i_{T_{0j}^\alpha}\left( -2 \beta_j\left(B_j^\beta \right)  + 2\mathcal{L}_{\xi_j}\left(\beta_j\left(w_{0j}^\beta \right) \right) + 2\mathcal{L}_{\xi_j}\left( B_j^\beta \right) - \mathcal{L}_{\xi_j}\mathcal{L}_{\xi_j}\left(w_{0j}^\beta \right)  \right) \notag
\end{align}}}

On the other hand, by the following equality
{\small{\begin{align*}
&2 i_{\left[ \xi_j, T_{0j}^\alpha \right]}\mathcal{L}_{\xi_j}\left( w_{0j}^\beta \right) =- 2i_{\left[ \xi_j, \left[ \xi_j, T_{0j}^\alpha \right] \right]}\left(  w_{0j}^\beta \right) - 2 \mathcal{L}_{\xi_j}i_{T_{0j}^\alpha} \mathcal{L}_{\xi_j}\left( w_{0j}^\beta \right)  = -2 i_{\left[ \xi_j, \left[ \xi_j, T_{0j}^\alpha \right] \right]}\left(  w_{0j}^\beta \right)- 2 i_{\left[ \xi_j, T_{0j}^\alpha \right]}\mathcal{L}_{\xi_j}\left( w_{0j}^\beta \right) - 2 i_{T_{0j}^\alpha} \mathcal{L}_{\xi_j}\mathcal{L}_{\xi_j}\left( w_{0j}^\beta \right) \\
&\Longrightarrow 2i_{ \left[ \xi_j, T_{0j}^\alpha \right]}\mathcal{L}_{\xi_j}\left( w_{0j}^\beta \right) = - i_{\left[ \xi_j, \left[ \xi_j, T_{0j}^\alpha \right]\right]} \left( w_{0j}^\beta \right) - i_{T_{0j}^\alpha} \mathcal{L}_{\xi_j} \mathcal{L}_{\xi_j} \left( w_{0j}^\beta \right)
\end{align*}}}
we have
\begin{align}\label{iit1}
&2i_{\Phi_j^\alpha}\left(\Psi_j^\beta \right) =2i_{\left[ \xi_j, T_{0j}^\alpha \right] - A_j^\alpha - \alpha_j\left(T_{0j}^\alpha \right)} \left( \mathcal{L}_{\xi_j}\left(w_{0j}^\beta \right) - B_j^\beta - \beta_j\left(w_{0j}^\beta \right)\right)\\
&= 2i_{\left[ \xi_j, T_{0j}^\alpha \right]}\mathcal{L}_{\xi_j}\left( w_{0j}^\beta \right) - 2 i_{A_j^\alpha}\left(\mathcal{L}_{\xi_j}\left( w_{0j}^\beta \right) \right) - 2 i_{\alpha_j\left( T_{0j}^\alpha \right)}\mathcal{L}_{\xi_j}\left(w_{0j}^\beta \right) - 2i_{\left[ \xi_j,  T_{0j}^\alpha \right]}\left( B_j^\beta \right) + 2 i_{\alpha_j\left( T_{0j}^\alpha \right)}\left( B_j^\alpha \right) \notag \\
&- 2i_{\left[ \xi_j, T_{0j}^\alpha \right]}\left( \beta_j\left( w_{0j}^\beta \right) \right) + 2 i_{A_j^\alpha}\beta_j\left( w_{0j}^\beta \right) + 2 i_{\alpha_j\left( T_{0j}^\alpha \right)} \left( \beta_j\left( w_{0j}^\beta \right) \right) \notag \\
&=  - i_{\left[ \xi_j, \left[ \xi_j, T_{0j}^\alpha \right]\right]} \left( w_{0j}^\beta \right) - i_{T_{0j}^\alpha} \mathcal{L}_{\xi_j} \mathcal{L}_{\xi_j} \left( w_{0j}^\beta \right) +2 i_{\left[ \xi_j, A_j^\alpha \right]}\left( w_{0j}^\beta \right)  + 2 i_{\left[ \xi_j, \alpha_j\left( T_{0j}^\alpha \right)\right]} \left( w_{0j}^\beta \right) - \cancel{ 2 \mathcal{L}_{\xi_j} i_{\alpha_j\left( T_{0j}^\alpha \right)}\left( w_{0j}^\beta \right)}  \notag \\
&+2 i_{T_{0j}^\alpha}\mathcal{L}_{\xi_j}\left(B_j^\alpha \right) - 2i_{T_{0j}^\alpha}\left( \beta_j\left( B_j^\beta \right) \right) - \cancel{2 \mathcal{L}_{\xi_j} i_{T_{0j}^\alpha}\left( \beta_j \left(w_{0j}^\beta \right) \right)} + 2 i_{T_{0j}^\alpha} \mathcal{L}_{\xi_j}\left( \beta_j \left( w_{0j}^\beta \right) \right) - 2i_{\alpha_j\left( A_j^\alpha \right) } \left( w_{0j}^\beta \right) + 2 i_{\alpha_j\left( T_{0j}^\alpha \right) } \left( \beta_j\left( w_{0j}^\beta \right) \right) \notag
\end{align} 

Then from $(\ref{iit2})$ and $(\ref{iit1})$, $\left\{ P_j \right\} $ corresponds to $p$.

\end{proof}

\section{Infinitesimal foliated deformations and obstructions in terms of both tangent and cotangent sheaves by using \v Cech resolution of $\mathcal{E}_{\mathcal{F}_0}'^{\bullet}$ }

We recall Definition $\ref{s1}$ of a foliated complex analytic family $\left( \mathcal{M},\Theta_\mathcal{F}, \mathcal{N}_\mathcal{F}^*, B, \pi \right)$ in terms of both tangent and cotangent sheaves. We assume that $\Theta_\mathcal{F}$ and $\mathcal{N}_{\mathcal{F}}^*$ are locally free. We will replace the integrability condition $(6)$ $d_{\mathcal{M}/B}\left(\mathcal{N}_\mathcal{F}^* \right) \subset \mathcal{N}_\mathcal{F}^* \bigwedge \Omega_{\mathcal{M}/B}^1$ at every point of $\mathcal{M}-\bigcup_{t\in B} S_t$, where $S_t=\textnormal{Sing}\left( \mathcal{F}_t \right)$ for $t\in B$ in Definition \ref{s1} by the condition $(6)'$ $d_{\mathcal{M}/B} \left(\mathcal{N}_\mathcal{F}^*\right)$ vanishes on $\bigwedge^2 \Theta_\mathcal{F}$. We shall keep the notations in Remark \ref{o4} and Remark \ref{oc2} and Remark \ref{os3}. Then we only need to discuss the condition $(6)'$. The condition $(6)'$ implies that for $\alpha,\beta=1,...,p$ and $\gamma=1,...,q$, we have
\begin{align}\label{osi2}
i_{T_j^\alpha(\xi_j, t)\wedge T_j^\beta(\xi_j,t)}\left( dw_j^\gamma(\xi_j,t) \right) = 0
\end{align}

Given a foliated complex analytic family $\left(\mathcal{M}, \Theta_\mathcal{F}, \mathcal{N}_\mathcal{F}^*,  B, \pi \right)$ with $\Theta_\mathcal{F}$ and $\mathcal{N}_\mathcal{F}^*$ locally free of deformations of $\left(M, \Theta_{\mathcal{F}_0} , \mathcal{N}_{\mathcal{F}_0}^* \right)=\pi^{-1}(0),0\in B$, the infinitesimal foliated deformations in the foliated complex analytic family  are encoded in $\mathbb{H}^1\left( M,\mathcal{F}_0^\bullet \right)$ which is the first cohomology group of the following complex
{\small{\begin{center}
$\mathcal{F}_0'^\bullet: \begin{CD}
\cdots \\
@AF_3'AA \\
\mathscr{H}om_{\mathcal{O}_{M}}\left(\bigwedge^3 \Theta_{\mathcal{F}_0},  \frac{\Theta_{M}}{\Theta_{\mathcal{F}_0}} \right) \bigoplus \mathscr{H}om_{\mathcal{O}_{M}}\left(\mathcal{N}_{\mathcal{F}_0}^*, \bigwedge^3 \Theta_{\mathcal{F}_0}^* \right) \bigoplus \mathscr{H}om_{\mathcal{O}_{M}}\left(\mathcal{N}_{\mathcal{F}_0}^*,\bigwedge^2 \Theta_{\mathcal{F}_0}^* \right)       \\
@AF_2'AA \\
\mathscr{H}om_{\mathcal{O}_{M}}\left(\bigwedge^2 \Theta_{\mathcal{F}_0}, \frac{\Theta_{M}}{\Theta_{\mathcal{F}_0} }\right) \bigoplus \mathscr{H}om_{\mathcal{O}_{M}}\left(\mathcal{N}_{\mathcal{F}_0}^*, \bigwedge^2 \Theta_{\mathcal{F}_0}^* \right) \bigoplus \mathscr{H}om_{\mathcal{O}_{M}}\left(\mathcal{N}_{\mathcal{F}_0}^*, \Theta_{\mathcal{F}_0}^* \right)       \\
@AF_1'AA \\
\mathscr{H}om_{\mathcal{O}_{M}}\left(\Theta_{\mathcal{F}_0}, \frac{\Theta_{M}}{\Theta_{\mathcal{F}_0}} \right) \bigoplus \mathscr{H}om_{\mathcal{O}_M}\left( \mathcal{N}_{\mathcal{F}_0}^*, \frac{\Omega_{M}^1}{\mathcal{N}_{\mathcal{F}_0}^* } \right)  \\
@AF_0 AA \\
\Theta_{M}
\end{CD}$
\end{center}}}
We will denote the $i$-th cohomology group by $\mathbb{H}^1\left( M , \mathcal{F}_0'^\bullet \right)$. We can compute $\mathbb{H}^i\left(M_0, \mathcal{F}_0^\bullet \right)$  by the following \v Cech resolution of $\mathcal{F}_0^\bullet$. Here $\delta$ is the \v Cech map and $\mathcal{U}_0= \mathcal{U}\cap M=\left\{ U_j :=\mathcal{U}_j \cap M | j=1,2,... \right\}$ is an open covering of $M$:
{\tiny{\begin{equation}\label{osi4}
\begin{CD}
\cdots \\
@AF_2'AA \\
C^0\left( \mathcal{U}_0, \left(\bigwedge^2 \Theta_{\mathcal{F}_0}^*\otimes \frac{\Theta_{M}}{\Theta_{\mathcal{F}_0}}\right) \bigoplus \left(\left(\mathcal{N}_{\mathcal{F}_0}^* \right)^*\otimes \bigwedge^2 \Theta_{\mathcal{F}_0}^*\right) \bigoplus \left( \left( \mathcal{N}_{\mathcal{F}_0}^* \right)^*\otimes \Theta_{\mathcal{F}_0}^* \right)   \right) @>-\delta>> \cdots \\
@AF_1'AA @AF_1'AA \\
C^0\left(\mathcal{U}_0, \left(\Theta_{\mathcal{F}_0}^*\otimes \frac{\Theta_{M}}{\Theta_{\mathcal{F}_0}}\right) \bigoplus \left(\left(\mathcal{N}_{\mathcal{F}_0}^* \right)^*\otimes \frac{\Omega_{M}^1}{\mathcal{N}_{\mathcal{F}_0}^*} \right) \right) @>\delta>> C^1\left(\mathcal{U}_0, \left(\Theta_{\mathcal{F}_0}^*\otimes \frac{\Theta_{M}}{\Theta_{\mathcal{F}_0}}\right) \bigoplus \left(\left(\mathcal{N}_{\mathcal{F}_0}^* \right)^*\otimes \frac{\Omega_{M}^1}{\mathcal{N}_{\mathcal{F}_0}^*} \right) \right) @>-\delta>> \cdots \\
@AF_0AA @AF_0AA  @AF_0AA \\
C^0\left(\mathcal{U}_0, \Theta_{M} \right) @>-\delta >> C^1\left( \mathcal{U}_0, \Theta_{M} \right) @>>> C^2\left(\mathcal{U}_0, \Theta_{M} \right) 
\end{CD}
\end{equation}}}
Then we can define the foliated Kodaira-Spencer map at $0\in B$ in terms of both tangent and cotangent sheaves which encodes the information of infinitesimal foliated deformations of $\pi^{-1}(0)=\left(M, \Theta_{\mathcal{F}_0}, \mathcal{N}_{\mathcal{F}_0}^* \right)$ in the foliated  analytic family $\left( \mathcal{M}, \Theta_\mathcal{F}, \mathcal{N}_\mathcal{F}^*, B, \pi \right)$ (see Part I) as in Defintion \ref{s8}. In Part II we have shown that the complex of sheaves
{\small{\begin{equation} \label{osi1}
\mathcal{E}_{\mathcal{F}_0}'^\bullet: \begin{CD}
\cdots \\
@AF_3'AA \\
\mathscr{H}om_{\mathcal{O}_{M}}\left(\bigwedge^3 \Theta_{\mathcal{F}_0},  \frac{\Theta_{M}}{\Theta_{\mathcal{F}_0}} \right) \bigoplus \mathscr{H}om_{\mathcal{O}_{M}}\left(\mathcal{N}_{\mathcal{F}_0}^*, \tilde{\mathcal{S}}^3 \right) \bigoplus \mathscr{H}om_{\mathcal{O}_{M}}\left(\mathcal{N}_{\mathcal{F}_0}^*,\bigwedge^2 \Theta_{\mathcal{F}_0}^* \right)       \\
@AF_2'AA \\
\mathscr{H}om_{\mathcal{O}_{M}}\left(\bigwedge^2 \Theta_{\mathcal{F}_0}, \frac{\Theta_{M}}{\Theta_{\mathcal{F}_0} }\right) \bigoplus \mathscr{H}om_{\mathcal{O}_{M}}\left(\mathcal{N}_{\mathcal{F}_0}^*, \tilde{\mathcal{S}}^2 \right) \bigoplus \mathscr{H}om_{\mathcal{O}_{M}}\left(\mathcal{N}_{\mathcal{F}_0}^*, \Theta_{\mathcal{F}_0}^* \right)       \\
@AF_1''AA \\
\mathscr{H}om_{\mathcal{O}_{M}}\left(\Theta_{\mathcal{F}_0},\Theta_M\right) \bigoplus \mathscr{H}om_{\mathcal{O}_M}\left( \mathcal{N}_{\mathcal{F}_0}^*, \Omega_M^1 \right)  \\
@AF_0' AA \\
\mathcal{E}_{\mathcal{F}_0}
\end{CD}
\end{equation}}}
also controls simultaneous foliated deformations of $\left( M , \Theta_{\mathcal{F}_0}, \mathcal{N}_{\mathcal{F}_0}^*  \right)$ in terms of both tangent and cotangent sheaves.  Let us identify the infinitesimal deformations in terms of the complex (\ref{osi1}) of sheaves  $\mathcal{E}_{\mathcal{F}_0}'^\bullet$.

\section{Infinitesimal foliated deformations and obstructions in terms of both tangent and cotangent sheaves by using Dolbeault resolution of $\mathcal{E}_{\mathcal{F}_0}'^{\bullet}$}

\subsection{Identification of infinitesimal simultaneous foliated deformations by $\mathcal{E}_{\mathcal{F}_0}'^\bullet$} \

Let $\left( \mathcal{M}, \Theta_\mathcal{F}, \mathcal{N}_\mathcal{F}^*, B, \pi \right)$ be a foliated analytic family with $\Theta_{\mathcal{F}}$ and $\mathcal{N}_\mathcal{F}^*$ locally free as in Definition \ref{s1} with the integrability condtion $(6)'$ $d_{\mathcal{M}/B}\left(\mathcal{N}_\mathcal{F}^* \right)$ vanishes on $\bigwedge^2 \Theta_\mathcal{F}$. instead of $(6)$. We keep the notations in subsection \ref{its1}. Then as in $(\ref{its2})$, we have the infinitesimal simultaneous deformation along $\frac{\partial}{\partial t}$ at $t$:
\begin{align}\label{a30}
\frac{\partial \left( M_t, \Theta_{\mathcal{F}_t}, \mathcal{N}_{\mathcal{F}_t}^* \right) }{\partial t}=\left( \left\{ \alpha_j(t)\right\}, \left\{ \beta_j(t)\right\}, \left\{ \left( \theta_{jk}(t), A_{jk}(t), B_{jk}(t)\right)  \right\} \right)\in \mathbb{H}^1 \left( M_t, \mathcal{E}_{\mathcal{F}_t}'^\bullet   \right)
\end{align}
in the following \v Cech resolution of $\mathcal{E}_{\mathcal{F}_0}'^\bullet$:
{\tiny{\begin{equation}\label{its25}
\begin{CD}
\cdots \\
@AF_2'^tAA \\
C^0\left( \mathcal{U}_t, \left(\bigwedge^2 \Theta_{\mathcal{F}_t}^*\otimes \frac{\Theta_{M_t}}{\Theta_{\mathcal{F}_t}}\right) \bigoplus \left(\left(\mathcal{N}_{\mathcal{F}_t}^* \right)^*\otimes  \bigwedge^2 \Theta_{\mathcal{F}_t}^*\right) \bigoplus \left( \left( \mathcal{N}_{\mathcal{F}_t}^* \right)^*\otimes \Theta_{\mathcal{F}_t}^* \right)   \right) @>-\delta>> \cdots \\
@AF_1''^tAA @AF_1''^tAA \\
C^0\left(\mathcal{U}_t, \left(\Theta_{\mathcal{F}_t}^*\otimes \Theta_{M_t}\right) \bigoplus \left(\left(\mathcal{N}_{\mathcal{F}_t}^* \right)^*\otimes \Omega_{M_t}^1 \right) \right) @>\delta>> C^1\left(\mathcal{U}_t, \left(\Theta_{\mathcal{F}_t}^*\otimes  \Theta_{M_t} \right) \bigoplus \left(\left(\mathcal{N}_{\mathcal{F}_t}^* \right)^*\otimes \Omega_{M_t}^1 \right) \right) @>-\delta>> \cdots \\
@AF_0'^tAA @AF_0'^tAA  @AF_0'^tAA \\
C^0\left(\mathcal{U}_t, \mathcal{E}_{\mathcal{F}_t} \right) @>-\delta >> C^1\left( \mathcal{U}_t, \mathcal{E}_{\mathcal{F}_t} \right) @>\delta>> C^2\left(\mathcal{U}_t, \mathcal{E}_{\mathcal{F}_t} \right) 
\end{CD}
\end{equation}}}

\subsection{Description of obstructed elements in \v Cech resolution of $\mathcal{E}_{\mathcal{F}_0}'^\bullet$ } \

We have the exact same statements with Definition \ref{its3}, Remark \ref{its4}, Definition \ref{its5} in terms of the complex of sheaves $\mathcal{E}_{\mathcal{F}_0}'^\bullet$. We now describe a condition when $a\in \mathbb{H}^1\left( M, \mathcal{E}_{\mathcal{F}_0}'^\bullet \right)$ is an obstructed element in simultaneous foliated deformations in terms of both tangent and cotangent sheaves. We keep the notations in subsection \ref{its6}. Let $\left(\mathcal{M}, \Theta_\mathcal{F}, \mathcal{N}_\mathcal{F}^*, B, \pi \right)$ with $0\in B \subset \mathbb{C}$ and $\Theta_\mathcal{F}, \mathcal{N}_{\mathcal{F}}^*$ locally free be a foliated analytic family in terms of both tangent and cotangent sheaves such that $\pi^{-1}(0)=\left( M , \Theta_{\mathcal{F}_0}, \mathcal{N}_{\mathcal{F}_0}^* \right)$. Take a small disk $\Delta$ with centre $0$ such that $0\in \Delta \subset B$, and represent $\left(\mathcal{M}_\Delta, \Theta_{\mathcal{F}_\Delta}, \mathcal{N}_{\mathcal{F}_\Delta}^* \right)=\pi^{-1}\left(\Delta \right)$ in the form
\begin{align*}
\left( \mathcal{M}_\Delta, \Theta_{\mathcal{F}_\Delta}, \mathcal{N}_{\mathcal{F}_\Delta}^* \right)= \bigcup_{j=1}^l \left( U_j\times \Delta, \Theta_{\mathcal{F}_\Delta}, \mathcal{N}_{\mathcal{F}_\Delta}^* \right)
\end{align*}
as in subsection $(\ref{its6})$. We only need to consider the condition $i_{T_j^\alpha(\xi_j, t)\wedge T_j^\beta(\xi_j,t)}\left( dw_j^\gamma(\xi_j,t) \right)=0,\alpha,\beta=1,...,p,\gamma=1,...,q$. By taking the derivative with respect to $t$, we have
\begin{align}\label{its7}
i_{\frac{\partial T_j^\alpha(\xi_j,t)}{\partial t}\wedge T_j^\beta(\xi_j,t) }\left( dw_j^\gamma(\xi_j,t) \right) + i_{T_j^\alpha(\xi_j,t) \wedge \frac{\partial T_j^\beta(\xi_j,t)}{\partial t} } \left( dw_j^\gamma(\xi_j,t) \right) + i_{T_j^\alpha(\xi_j,t)\wedge T_j^\beta(\xi_j,t)}\left( d\left(\frac{\partial w_j^\gamma(\xi_j,t)}{\partial t} \right)\right) = 0
\end{align}

By taking the derivative of $(\ref{its7})$ with respect to $t$ again, we obtain
\begin{align} 
&i_{\frac{\partial^2 T_j^\alpha(\xi_j,t)}{\partial t^2} \wedge T_j^\beta(\xi_j,t)} \left(d w_j^\gamma(\xi_j,t) \right) + 2i_{\frac{\partial T_j^\alpha(\xi_j,t)}{\partial t} \wedge \frac{\partial T_j^\beta(\xi_j,t)}{\partial t}}\left( dw_j^\gamma(\xi_j,t) \right) + 2i_{\frac{\partial T_j^\alpha(\xi_j,t)}{\partial t} \wedge T_j^\beta(\xi_j,t)}\left( d\left(\frac{\partial w_j^\gamma(\xi_j,t)}{\partial t} \right) \right) \label{its10} \\
&+ i_{T_j^\alpha(\xi_j,t) \wedge \frac{\partial^2 T_j^\beta(\xi_j,t)}{\partial t^2} } \left( dw_j^\gamma(\xi_j,t) \right) + 2 i_{T_j^\alpha(\xi_j,t)\wedge \frac{\partial T_j^\beta(\xi_j,t)}{\partial t}}\left(  d\left( \frac{\partial w_j^\gamma(\xi_j,t)}{\partial t} \right) \right) + i_{T_j^\alpha(\xi_j,t) \wedge T_j^\beta(\xi_j,t) }\left(d\left( \frac{\partial^2 w_j^\gamma(\xi_j,t)}{\partial t^2} \right) \right) = 0 \notag
\end{align}

We recall from Part I that we set
\begin{align}\label{its9}
i_{T_j^\alpha(\xi_j,t)}\left( dw_j^\gamma(\xi_j,t)\right) = \sum_{\eta=1}^q a_{j\alpha}^{\gamma \eta}(\xi_j,t) w_j^\eta(\xi_j,t)\,\,\,\,\,\,\,\, \textnormal{for some $a_{j\alpha}^{\gamma \eta}(\xi_j,t)\in \Gamma\left( \mathcal{U}, \mathcal{O}_\mathcal{M}\right)$}
\end{align}

Then from $(\ref{its9})$, we see that $(\ref{its10})$ is equivalent to 
\begin{align}
&- \sum_{\eta=1}^q a_{j\beta}^{\gamma \eta}(\xi_j,t) i_{\frac{\partial^2 T_j^\alpha(\xi_j,t)}{\partial t}}\left( w_j^\eta(\xi_j,t) \right) + 2i_{\frac{\partial T_j^\alpha(\xi_j,t)}{\partial t} \wedge \frac{\partial T_j^\beta(\xi_j,t)}{\partial t}}\left( dw_j^\gamma(\xi_j,t) \right) + 2i_{\frac{\partial T_j^\alpha(\xi_j,t)}{\partial t} \wedge T_j^\beta(\xi_j,t)}\left( d\left(\frac{\partial w_j^\gamma(\xi_j,t)}{\partial t} \right) \right) \label{its11}\\
&+ \sum_{\eta=1}^q a_{j\alpha}^{\gamma \eta}(\xi_j,t) i_{\frac{\partial^2 T_j^\beta(\xi_j,t)}{\partial t^2} }\left( w_j^\gamma(\xi_j,t) \right) + 2 i_{T_j^\alpha(\xi_j,t)\wedge \frac{\partial T_j^\beta(\xi_j,t)}{\partial t}}\left(  d\left( \frac{\partial w_j^\gamma(\xi_j,t)}{\partial t} \right) \right) + i_{T_j^\alpha(\xi_j,t) \wedge T_j^\beta(\xi_j,t) }\left(d\left( \frac{\partial^2 w_j^\gamma(\xi_j,t)}{\partial t^2} \right) \right) = 0 \notag
\end{align}

From $(\ref{its8})$ $i_{\frac{\partial^2 T_j^\alpha(\xi_j,t)}{\partial t^2}} \left(w_j^\beta(\xi_j,t) \right) =  - i_{T_j^\alpha(\xi_j,t)}\left( \frac{\partial^2 w_j^\beta(\xi_j,t)}{\partial t^2} \right) -2 i_{\frac{\partial T_j^\alpha(\xi_j,t)}{\partial t}}\left( \frac{\partial w_j^\beta(\xi_j,t)}{\partial t} \right)$, we see that $(\ref{its11})$ is equivalent to
\begin{align*}
&\sum_{\eta=1}^q a_{j\beta}^{\gamma \eta}(\xi_j,t)i_{T_j^\alpha(\xi_j,t)}\left(\frac{\partial^2 w_j^\gamma(\xi_j,t)}{\partial t^2} \right) - \sum_{\eta=1}^q a_{j\alpha}^{\gamma \eta}(\xi_j,t) i_{T_j^\beta(\xi_j,t)}\left( \frac{\partial^2 w_j^\gamma(\xi_j,t)}{\partial t^2} \right) + i_{T_j^\alpha(\xi_j,t) \wedge T_j^\beta(\xi_j,t) }\left(d\left( \frac{\partial^2 w_j^\gamma(\xi_j,t)}{\partial t^2} \right) \right)\\
&= -2\sum_{\eta=1}^q a_{j\beta}^{\gamma \eta}(\xi_j,t) i_{\frac{\partial T_j^\alpha(\xi_j,t)}{\partial t}}\left( \frac{\partial w_j^\gamma(\xi_j,t)}{\partial t} \right) + 2\sum_{\eta=1}^q a_{j\alpha}^{\gamma \eta}(\xi_j,t) i_{\frac{\partial T_j^\beta(\xi_j,t)}{\partial t}}\left( \frac{\partial w_j^\gamma (\xi_j,t)}{\partial t} \right) - 2i_{\frac{\partial T_j^\alpha(\xi_j,t)}{\partial t} \wedge \frac{\partial T_j^\beta(\xi_j,t)}{\partial t}}\left( dw_j^\gamma(\xi_j,t) \right) \\
& - 2i_{\frac{\partial T_j^\alpha(\xi_j,t)}{\partial t} \wedge T_j^\beta(\xi_j,t)}\left( d\left(\frac{\partial w_j^\gamma(\xi_j,t)}{\partial t} \right) \right) -  2 i_{T_j^\alpha(\xi_j,t)\wedge \frac{\partial T_j^\beta(\xi_j,t)}{\partial t}}\left(  d\left( \frac{\partial w_j^\gamma(\xi_j,t)}{\partial t} \right) \right)
\end{align*}
which implies that
\begin{align*}
&\sum_{\eta=1}^q a_{j\beta}^{\gamma \eta} i_{T_j^\alpha}\left(\dot{\beta}_j(t)\left(w_j^\gamma \right) \right) - \sum_{\eta=1}^q a_{j\alpha}^{\gamma \eta} i_{T_j^\beta}\left( \dot{\beta}_j(t)\left(w_j^\gamma \right) \right) + i_{T_j^\alpha \wedge T_j^\beta}\left( d \dot{\beta}_j(t)\left( w_j^\gamma \right) \right)\\
&= 2 \sum_{\eta=1}^q a_{j\beta}^{\gamma \eta} i_{\alpha_j\left( T_j^\alpha \right)} \left(\beta_j(t)\left( w_j^\gamma \right) \right) -2 \sum_{\eta=1}^q a_{j\alpha}^{\gamma \eta} i_{\alpha_j(t)\left( T_j^\beta \right)} \left( \beta_j(t)\left( w_j^\gamma \right) \right) + 2 i_{\alpha_j(t)\left(T_j^\alpha \right) \wedge \alpha_j(t)\left( T_j^\beta \right)}\left( dw_j^\gamma \right)\\
&+2 i_{\alpha_j(t)\left( T_j^\alpha \right) \wedge T_j^\beta}\left( d\beta_j(t)\left( w_j^\gamma \right)\right) + 2 i_{T_j^\alpha \wedge \alpha_j(t) \left( T_j^\beta \right)}\left(d  \beta_j(t)\left(  w_j^\gamma \right) \right)
\end{align*}

We define $\left[ \left(\alpha_j(t), \beta_j(t) \right), \left( \alpha_j(t), \beta_j(t) \right)   \right]_{\mathcal{F}_0^\sharp}\in \Gamma\left( U_j^t, \mathscr{H}om_{\mathcal{O}_{M_t}} \left( \bigwedge^2 \Theta_{\mathcal{F}_t}, \frac{\Theta_M}{\Theta_{\mathcal{F}_t}}  \right) \right) \bigoplus \Gamma\left( U_j^t, \mathscr{H}om_{\mathcal{O}_{M_t}}\left( \mathcal{N}_{\mathcal{F}_t}^*, \bigwedge^2 \Theta_{\mathcal{F}_t}^* \right) \right)$ by
{\small{\begin{align}
\left[ \left(\alpha_j(t), \beta_j(t) \right), \left( \alpha_j(t), \beta_j(t) \right)   \right]_{\mathcal{F}_0^\sharp}:& \Gamma\left( U_j^t, \bigwedge^2 \Theta_{\mathcal{F}_t} \right)\bigoplus \Gamma\left( U_j^t, \mathcal{N}_{\mathcal{F}_t}^* \right) \to \Gamma\left( U_j^t, \frac{\Theta_{M_t}}{\Theta_{\mathcal{F}_t}} \right) \bigoplus \Gamma\left( U_j^t, \bigwedge^2 \Theta_{\mathcal{F}_t}^* \right) \\
\left(T_j^\alpha \wedge T_j^\beta, w_j^\gamma\right)\mapsto &  \overline{2\left[\alpha_j\left(T_j^\alpha \right), \alpha_j\left(  T_j^\beta\right) \right] + 2\alpha_j(t)\left(\sum_{\gamma=1}^p \frac{\partial g_{j\alpha\beta}^\gamma}{\partial t} T_j^\gamma \right)} \notag \\
&\bigoplus 2 \sum_{\eta=1}^q a_{j\beta}^{\gamma \eta} i_{\alpha_j\left( T_j^\alpha \right)} \left(\beta_j(t)\left( w_j^\gamma \right) \right) -2 \sum_{\eta=1}^q a_{j\alpha}^{\gamma \eta} i_{\alpha_j(t)\left( T_j^\beta \right)} \left( \beta_j(t)\left( w_j^\gamma \right) \right) + 2 i_{\alpha_j(t)\left(T_j^\alpha \right) \wedge \alpha_j(t)\left( T_j^\beta \right)}\left( dw_j^\gamma \right) \notag\\
&\,\,\,\,\,\,\,\,\,+2 i_{\alpha_j(t)\left( T_j^\alpha \right) \wedge T_j^\beta}\left( d\beta_j(t)\left( w_j^\gamma \right)\right) + 2 i_{T_j^\alpha \wedge \alpha_j(t) \left( T_j^\beta \right)}\left( d \beta_j(t)\left(  w_j^\gamma \right) \right) \notag
\end{align}}}
and linearly extends to $\Gamma\left( U_j^t,  \bigwedge^2 \Theta_{\mathcal{F}_t}\right)\bigoplus \Gamma\left( U_j^t, \mathcal{N}_{\mathcal{F}_t}^* \right)$.

By putting $t=0$, we define
{\small{\begin{align*}
\left[ \left(\alpha_j(0), \beta_j(0) \right)  , \left( \alpha_j(0), \beta_j(0) \right) \right]_{\mathcal{F}_0'}&=\left( \left[ \left(\alpha_j(0), \beta_j(0) \right)  , \left( \alpha_j(0), \beta_j(0) \right) \right]_{\mathcal{F}_0^\sharp}, \alpha_j(0)* \beta_j(0) \right) \\
&\in \Gamma\left( U_j^0, \mathscr{H}om_{\mathcal{O}_M}\left( \bigwedge^2 \Theta_{\mathcal{F}_0}, \frac{\Theta_M}{\Theta_{\mathcal{F}_0}} \right) \bigoplus \mathscr{H}om_{\mathcal{O}_M}\left(\mathcal{N}_{\mathcal{F}_0}^*, \bigwedge^2 \Theta_{\mathcal{F}_0}^* \right) \bigoplus  \mathscr{H}om_{\mathcal{O}_M}\left( \mathcal{N}_{\mathcal{F}_0}^*, \Theta_{\mathcal{F}_0}^* \right)\right)
\end{align*}}}
Then from $(\ref{its12})$, we have
{\Small{\begin{align}
&F_1''\left(\dot{\alpha}_j(0) , \dot{\beta}_j(0) \right)= \left[\left(\alpha_j(0), \beta_j(0) \right), \left( \alpha_j(0) , \beta_j(0) \right) \right]_{\mathcal{F}_0'} \label{itt2}  \\
&\left(\dot{\alpha}_k(0), \dot{\beta}_k(0)\right)- \left(\dot{\alpha}_j(0),\dot{\beta}_j(0)\right) + F_0'\left(\left(\dot{\theta}_{jk}(0), \dot{A}_{jk}(0),\dot{B}_{jk}(0) \right)  \right) \notag \\
& = \left( -\alpha_j(0)\circ A_{jk}(0) - \alpha_k(0)\circ A_{jk}(0) +\mathcal{L}_{\theta_{jk}(0)}\circ \alpha_j(0) + \mathcal{L}_{\theta_{jk}(0)}\circ \alpha_k(0) , - \beta_j(0)\circ B_{jk}(0) - \beta_k (0) \circ B_{jk}(0)+ \mathcal{L}_{\theta_{jk}(0)}\circ \beta_j(0)+ \mathcal{L}_{\theta_{jk}(0)}\circ \beta_k(0) \right) \notag \\
&\left( \dot{\theta}_{ij}(0), \dot{A}_{ij}(0), \dot{B}_{ij}(0)\right)- \left( \dot{\theta}_{ik}(0), \dot{A}_{ik}(0), \dot{B}_{ik}(0) \right) + \left( \dot{\theta}_{jk}(0), \dot{A}_{jk}(0), \dot{B}_{jk}(0) \right)= \left[  \left( \theta_{ij} (0) , A_{ij}(0), B_{ij}(0) \right) , \left( \theta_{jk}(0), A_{jk}(0), B_{jk}(0) \right)   \right] \notag
\end{align}}}

\begin{lemma}\label{itt1}
We keep the notations of \textnormal{Lemma \ref{us24}}. For any $1$-cocycle
\begin{align*}
\left(\left\{ \alpha_j \right\}, \left\{ \beta_j \right\}, \left\{\left( \theta_{jk}, A_{jk}, B_{jk} \right) \right\} \right)\in C^0\left( \mathcal{U}, \mathscr{H}om_{\mathcal{O}_M}\left( \Theta_{\mathcal{F}_0},  \Theta_M \right)\right) \bigoplus C^0\left( \mathcal{U}, \mathscr{H}om_{\mathcal{O}_M}\left( \mathcal{N}_{\mathcal{F}_0}^*, \Omega_M^1 \right) \right) \bigoplus C^1\left(  \mathcal{U}, \mathcal{E}_{\mathcal{F}_0}\right)
\end{align*}
where $\alpha_j:\Gamma\left(U_j, \Theta_{\mathcal{F}_0} \right)\to \Gamma\left( U_j, \Theta_M \right)$, $\beta_j:\Gamma\left(U_j, \mathcal{N}_{\mathcal{F}_0}^* \right) \to \Gamma\left( U_j, \Omega_M^1 \right), \theta_{jk}\in \Gamma\left( U_{jk}, \Theta_M \right)$, $A_{jk}:\Gamma\left(U_{jk}, \Theta_{\mathcal{F}_0} \right)\to \Gamma\left( U_{jk}, \Theta_{\mathcal{F}_0} \right)$ and $B_{jk}:\Gamma\left( U_{jk}, \mathcal{N}_{\mathcal{F}_0}^* \right) \to \Gamma\left( U_{jk}, \mathcal{N}_{\mathcal{F}_0}^* \right)$ satisfying
\begin{enumerate}
\item $F_1''\left(\alpha_j, \beta_j \right)=0$, i.e. $\left[ T_{0j}^\alpha, \alpha_j\left( T_{0j}^\beta \right)\right]- \left[ T_{0j}^\beta, \alpha_j\left( T_{0j}^\alpha \right)\right]- \sum_{\gamma=1}^p g_{0j\alpha\beta}^\gamma \alpha_j\left( T_{0j}^\gamma \right)=Q_{j\alpha\beta} \in \Gamma\left( U_j, \Theta_{\mathcal{F}_0} \right)$ for $\alpha,\beta=1,...,p$ and $\sum_{\eta=1}^q a_{0j\beta}^{\delta \eta} i_{T_{0j}^\alpha}\left( \beta_j\left( w_{0j}^\eta \right) \right)- \sum_{\eta=1}^q a_{0j\alpha}^{\gamma \eta} i_{T_{0j}^\beta}\left( \beta_j \left( w_{0j}^\eta \right) \right) + i_{T_{0j}^\alpha \wedge T_{0j}^\beta}\left(d \beta_j\left(w_{0j}^\gamma \right) \right)=0$ for $\alpha,\beta=1,...,p, \gamma=1,...,q$ and $i_{T_{0j}^\alpha}\left( \beta_j\left( w_{0j}^\beta \right)\right) + i_{\alpha_j\left( T_{0j}^\alpha \right)}\left( w_{0j}^\beta \right)=0$ for $\alpha=1,...,p,\beta=1,...,q$.
\item $\left(\alpha_k, \beta_k\right)-\left(\alpha_j,\beta_j\right)+G_0'\left( \left(\theta_{jk}, A_{jk}, B_{jk} \right) \right)=0$, i.e. $\alpha_k(T_j^\alpha)- \alpha_j(T_j^\alpha)+ \mathcal{L}_{\theta_{jk}}(T_j^\alpha) - A_{jk}(T_j^\alpha)=0$ for $\alpha=1,...,p$, and $\beta_k\left(w_{0j}^\alpha\right) - \beta_j\left(w_{0j}^\alpha \right) +\mathcal{L}_{\theta_{jk}}\left(w_{0j}^\alpha \right) - B_{jk}\left(w_{0j}^\alpha \right)=0$ for $\alpha=1,...,q$
\item $\left( \theta_{ij}, A_{ij} , B_{ij}\right)  - \left(\theta_{ik}, A_{ik}, B_{ik} \right) + \left(\theta_{jk}, A_{jk}, B_{jk} \right)=0$, i.e. $\theta_{ij}-\theta_{ik}+ \theta_{jk}=0$ and $A_{ij}- A_{ik}+ A_{jk}^i=0$ and $B_{ij}- B_{ik}+ B_{jk}^i =0$.
\end{enumerate}
we define $R_j\in \Gamma\left( U_j, \mathscr{H}om_{\mathcal{O}_M}\left( \bigwedge^2 \Theta_{\mathcal{F}_0} \right) \right), F_{jk}\in \Gamma\left( U_{jk}, \mathscr{H}om_{\mathcal{O}_M}\left( \Theta_{\mathcal{F}_0}, \Theta_M \right) \right), G_{jk}\in \Gamma\left( U_{jk}, \mathscr{H}om_{\mathcal{O}_M}\left(\mathcal{N}_{\mathcal{F}_0}^*, \Omega_M^1 \right) \right)$ and $P_j\in \Gamma\left( U_j, \mathscr{H}om_{\mathcal{O}_M}\left( \mathcal{N}_{\mathcal{F}_0}^*, \Theta_{\mathcal{F}_0}^* \right) \right)$ as in \textnormal{Lemma \ref{us24}}. We define $S_j'\in \Gamma\left( U_j, \mathscr{H}om_{\mathcal{O}_M}\left( \mathcal{N}_{\mathcal{F}_0}^*, \bigwedge^2 \Theta_{\mathcal{F}_0}^* \right) \right)$ by
\begin{align*}
S_j' : \Gamma\left( U_j, \mathcal{N}_{\mathcal{F}_0}^* \right) &\to \Gamma\left( U_j, \bigwedge^2 \Theta_{\mathcal{F}_0}^* \right) \\
              w_{0j}^\gamma &\mapsto \left( T_{0j}^\alpha \wedge T_{0j}^\beta \mapsto  2 \sum_{\eta=1}^q a_{0j\beta}^{\gamma \eta} i_{\alpha_j\left( T_{0j}^\alpha \right)}\left(\beta_j\left( w_{0j}^\eta \right) \right) - 2 \sum_{\eta=1}^q a_{0j\alpha}^{\gamma \eta} i_{\alpha_j\left( T_{0j}^\beta \right)} \left(\beta_j\left( w_{0j}^\eta \right) \right) + 2 i_{\alpha_j\left( T_{0j}^\alpha \right)\wedge \alpha_j\left( T_{0j}^\beta \right)}\left( d w_{0j}^\gamma \right)\right. \\
              &\left.\,\,\,\,\,\,\,\,\,\,\,\,\,\,\,\,\,\,\,\,\,\,\,\,\, \,\,\,\,\,\,\,\,\,\,\,\,\,\,\,\,\,\,\,\,\,\,\,\,\,\,+2 i_{\alpha_j\left( T_{0j}^\alpha \right)\wedge T_{0j}^\beta}\left( d\beta_j\left( w_{0j}^\gamma \right) \right) + 2 i_{T_{0j}^\alpha \wedge \alpha_j\left( T_{0j}^\beta \right)}\left( d \beta_j\left( w_{0j}^\gamma \right) \right) \right)
\end{align*}
and linearly extends to $\Gamma\left( U_j, \mathcal{N}_{\mathcal{F}_0}^* \otimes \bigwedge^2 \Theta_{\mathcal{F}_0} \right)$. Then
{\small{\begin{align} \label{itt3}
&\left(\left\{ R_j \right\}, \left\{ S_j' \right\}, \left\{ P_j \right\} ,  \left\{ F_{jk} \right\}, \left\{ G_{jk} \right\}, \left\{ \left[ \left(\theta_{ij}, A_{ij}, B_{ij} \right), \left( \theta_{jk}, A_{jk}, B_{jk} \right) \right] \right\} \right)\\
& \in C^0\left(\mathcal{U}, \mathscr{H}om_{\mathcal{O}_M}\left( \bigwedge^2 \Theta_{\mathcal{F}_0}, \frac{\Theta_M}{\Theta_{\mathcal{F}_0}} \right) \right) \bigoplus C^0\left(\mathcal{U}, \mathscr{H}om_{\mathcal{O}_M}\left(\mathcal{N}_{\mathcal{F}_0}^*, \bigwedge^2 \Theta_{\mathcal{F}_0}^* \right) \right) \bigoplus C^0 \left( \mathcal{U}, \mathscr{H}om_{\mathcal{O}_M}\left( \mathcal{N}_{\mathcal{F}_0}^*, \Theta_{\mathcal{F}_0}^* \right) \right) \notag\\
& \bigoplus C^1\left(\mathcal{U}, \mathscr{H}om_{\mathcal{O}_M}\left(\Theta_{\mathcal{F}_0}, \Theta_M \right) \right) \bigoplus C^1\left(\mathcal{U}, \mathscr{H}om_{\mathcal{O}_M}\left(\mathcal{N}_{\mathcal{F}_0}^*, \Omega_M^1 \right) \right) \bigoplus C^2\left(\mathcal{U}, \mathcal{E}_{\mathcal{F}_0} \right) \notag
\end{align}}}
defines a $2$-cocycle in \v Cech resolution $(\ref{its25})$ of $\mathcal{E}_{\mathcal{F}_0}'^\bullet$.
\end{lemma}

\begin{proof}
By Lemma \ref{us24}, it remains to show that $\tilde{E}_2\left(S_j' \right)=0$ and $\beta_2\circ R_j + S_j'- \tilde{E}_1\left( P_j\right)=0$ and $-\delta\left( \left\{ S_j'\right\} \right)+ \tilde{E}_1\left( \left\{ G_{jk}\right\} \right)=0$. In fact, we note that from $(\ref{us27})$, we see that
\begin{align}\label{its26}
\frac{1}{2}i_{T_{0j}^\alpha \wedge T_{0j}^\beta}\left(S_j\left( w_{0j}^\gamma \right) \right)= w_{0j}^1\wedge \cdots \wedge w_{0j}^q \cdot \frac{1}{2} S_j'\left(T_{0j}^\alpha, T_{0j}^\beta \right)\left( w_{0j}^\gamma \right)
\end{align}
Then from $(\ref{its20})$, we have
\begin{align}\label{its21}
&i_{R_j}\left(T_{0j}^\alpha \wedge T_{0j}^\beta \right)\left(w_{0j}^\gamma \right) + S_j'\left(T_{0j}^\alpha, T_{0j}^\beta \right)\left( w_{0j}^\gamma \right) - \tilde{E}_1 \left( P_j \right)\left(T_{0j}^\alpha \wedge T_{0j}^\beta \right)\left(w_{0j}^\gamma \right) = 0 \Longrightarrow \beta_2\circ R_j + S_j'- \tilde{E}_1\left( P_j\right)=0 
\end{align}
On the other hand, since $\beta_2\circ D_2\left(R_j \right) +E_2'\left( S_j' \right)- E_2' \left( \beta_1\circ R_j + S_j' - E_1'\left( P_j\right)  \right)=0$ and $D_2\left(R_j \right)=0$ from Lemma \ref{o54}, and $ \beta_2\circ R_j + S_j' - E_1'\left( P_j\right)=0$ from $(\ref{its21})$, we have $E_2'\left(S_j' \right)=0$. From $(\ref{its26})$ and $-\delta\left(\left\{ S_j \right\} \right)+ E_1'\left(\left\{ G_{jk}\right\} \right)=0$ in the proof of Lemma \ref{oc29}, we have $-\delta\left( \left\{ S_j'\right\} \right)+ \tilde{E}_1\left( \left\{ G_{jk}\right\} \right)=0$. This completes the proof of Lemma \ref{itt1}.

\end{proof}

From $(\ref{itt2})$ and Lemma \ref{itt1}, we obtain
\begin{theorem}
Suppose given a compact foliated complex manifold $\left(M, \Theta_{\mathcal{F}_0}, \mathcal{N}_{\mathcal{F}_0}^* \right)$ with $\Theta_{\mathcal{F}_0}$ and $\mathcal{N}_{\mathcal{F}_0}^*$ locally free, and $\left(\left\{ \alpha_j \right\}, \left\{ \beta_j \right\}, \left\{\left( \theta_{jk}, A_{jk}, B_{jk} \right) \right\} \right)\in \mathbb{H}^1\left( M, \mathcal{E}_{\mathcal{F}_0}'^\bullet \right)$. In order that there may exist a foliated analytic family $\left( \mathcal{M}, \Theta_\mathcal{F}, \mathcal{N}_\mathcal{F}^*, B, \pi \right)$ in terms of both tangent and cotangent sheaves such that $\pi^{-1}(0)=\left( M, \Theta_{\mathcal{F}_0}, \mathcal{N}_{\mathcal{F}_0}^* \right)$, and that $\left( \frac{\partial\left( \mathcal{M}, \Theta_\mathcal{F}, \mathcal{N}_\mathcal{F}^* \right)}{\partial t} \right)|_{t=0}= \left( \left\{ \alpha_j \right\}, \left\{ \beta_j \right\}, \left\{ \left( \theta_{jk} , A_{jk}, B_{jk} \right)\right\} \right)$, it is necessary that
{\small{\begin{align*}
&\left(\left\{ R_j \right\}, \left\{ S_j' \right\}, \left\{ P_j \right\} ,  \left\{ F_{jk} \right\}, \left\{ G_{jk} \right\}, \left\{ \left[ \left(\theta_{ij}, A_{ij}, B_{ij} \right), \left( \theta_{jk}, A_{jk}, B_{jk} \right) \right] \right\} \right)\\
& \in C^0\left(\mathcal{U}, \mathscr{H}om_{\mathcal{O}_M}\left( \bigwedge^2 \Theta_{\mathcal{F}_0}, \frac{\Theta_M}{\Theta_{\mathcal{F}_0}} \right) \right) \bigoplus C^0\left(\mathcal{U}, \mathscr{H}om_{\mathcal{O}_M}\left(\mathcal{N}_{\mathcal{F}_0}^*, \bigwedge^2 \Theta_{\mathcal{F}_0}^* \right) \right) \bigoplus C^0 \left( \mathcal{U}, \mathscr{H}om_{\mathcal{O}_M}\left( \mathcal{N}_{\mathcal{F}_0}^*, \Theta_{\mathcal{F}_0}^* \right) \right)  \\
& \bigoplus C^1\left(\mathcal{U}, \mathscr{H}om_{\mathcal{O}_M}\left(\Theta_{\mathcal{F}_0}, \Theta_M \right) \right) \bigoplus C^1\left(\mathcal{U}, \mathscr{H}om_{\mathcal{O}_M}\left(\mathcal{N}_{\mathcal{F}_0}^*, \Omega_M^1 \right) \right) \bigoplus C^2\left(\mathcal{U}, \mathcal{E}_{\mathcal{F}_0} \right) 
\end{align*}}}
from $(\ref{itt3})$ in \textnormal{Lemma \ref{itt1}} defines the $0$ class in $\mathbb{H}^2\left( M, \mathcal{E}_{\mathcal{F}_0}'^\bullet \right)$.
\end{theorem}
Given $\left(\left\{ \alpha_j \right\}, \left\{ \beta_j \right\}, \left\{\left( \theta_{jk}, A_{jk}, B_{jk} \right) \right\} \right)\in \mathbb{H}^1\left( M, \mathcal{E}_{\mathcal{F}_0}'^\bullet \right)$, if the $2$-cocycle $\left(\left\{ R_j \right\}, \left\{ S_j' \right\}, \left\{ P_j \right\} ,  \left\{ F_{jk} \right\}, \left\{ G_{jk} \right\}\right.$, $\left. \left\{ \left[ \left(\theta_{ij}, A_{ij}, B_{ij} \right), \left( \theta_{jk}, A_{jk}, B_{jk} \right) \right] \right\} \right)$ is not $0$ in $\mathbb{H}^2\left( M, \mathcal{E}_{\mathcal{F}_0}'^\bullet \right)$, there is no foliated deformation $\left( M_t, \Theta_{\mathcal{F}_t}, \mathcal{N}_{\mathcal{F}_t}^* \right)$ in terms of both tangent and cotangent sheaves with $\left(M_0, \Theta_{\mathcal{F}_0}, \mathcal{N}_{\mathcal{F}_0}^* \right)=\left( M, \Theta_{\mathcal{F}_0}, \mathcal{N}_{\mathcal{F}_0}^* \right)$ and $\left( \frac{\partial\left( \mathcal{M}, \Theta_\mathcal{F}, \mathcal{N}_\mathcal{F}^* \right)}{\partial t} \right)|_{t=0}=\left( \left\{\alpha_j \right\}, \left\{\beta_j \right\}, \left\{\left( \theta_{jk} , A_{jk}, B_{jk} \right) \right\}\right)$

\section{Infinitesimal foliated deformations and obstructions in terms of both tangent and cotangent sheaves by using Dolbeault resolution of $\mathcal{E}_{\mathcal{F}_0}'^{\bullet}$}

Let $\left(\mathcal{M}, \Theta_\mathcal{F}, \mathcal{N}_\mathcal{F}^*, \pi, B \right)$ with $\Theta_\mathcal{F}$ and $\mathcal{N}_\mathcal{F}^*$ locally free be a foliated analytic family in terms of both tangent and cotangent sheaves of deformations of $\left( M, \Theta_{\mathcal{F}_0}, \mathcal{N}_{\mathcal{F}_0}^* \right)=\pi^{-1}(0)$.  We keep the notations in section \ref{itt4}. Then it remains to consider the condition
\begin{align*}
i_{T_j^\alpha(\xi_j,t)\wedge T_j^\beta(\xi_j, t)}\left(d w_j^\gamma(\xi_j,t) \right)=0 
\end{align*}
which is equivalent to
\begin{align}\label{osi5}
0&=i_{\left(\sum_{\eta, \sigma=1}^n V_j^{\alpha \eta}\frac{\partial \xi_j^\sigma}{\partial z^\eta}\frac{\partial}{\partial \xi_j^\sigma} \right) \wedge\left( \sum_{c , d=1}^n V_j^{\beta c}\frac{\partial \xi_j^d}{\partial z^c} \frac{\partial}{\partial \xi_j^d} \right)}d\left(   \sum_{\delta=1}^n w_j^{\gamma \delta} d\xi_j^\delta \right)\\
&=  i_{\left(\sum_{\eta, \sigma=1}^n V_j^{\alpha \eta}\frac{\partial \xi_j^\sigma}{\partial z^\eta}\frac{\partial}{\partial \xi_j^\sigma} \right) \wedge\left( \sum_{c , d=1}^n V_j^{\beta c}\frac{\partial \xi_j^d}{\partial z^c} \frac{\partial}{\partial \xi_j^d} \right)} \left(  \sum_{a, \delta=1}^n \frac{\partial w_j^{\gamma \delta}}{\partial \xi_j^a} d\xi_j^a \wedge d\xi_j^\delta     \right) \notag \\
&= \sum_{\eta, \sigma , c, d=1}^n V_j^{\alpha \eta} V_j^{\beta c} \frac{\partial w_j^{\gamma d}}{\partial \xi_j^\sigma}\frac{\partial \xi_j^\sigma}{\partial z^\eta} \frac{\partial \xi_j^d}{\partial z^c} - \sum_{\eta, \sigma , c, d=1}^n V_j^{\alpha \eta} V_j^{\beta c} \frac{\partial w_j^{\gamma \sigma}}{\partial \xi_j^d}\frac{\partial \xi_j^\sigma}{\partial z^\eta} \frac{\partial \xi_j^d}{\partial z^c} \notag
\end{align}

On the other hand, we note that
{\small{\begin{align}\label{osi6}
&i_{ \left( \sum_{\eta=1}^n V_j^{\alpha \eta} \frac{\partial}{\partial z^\eta} \right) \wedge \left( \sum_{c=1}^n V_j^{\beta c}\frac{\partial}{\partial z^c} \right) } \left(\partial \left( \sum_{\delta=1}^n W_j^{\gamma \delta} dz^\delta    \right)    \right)= i_{ \left( \sum_{\eta=1}^n V_j^{\alpha \eta} \frac{\partial}{\partial z^\eta} \right) \wedge \left( \sum_{c=1}^n V_j^{\beta c}\frac{\partial}{\partial z^c} \right) } \left(\partial \left( \sum_{\delta, \sigma=1}^n w_j^{\gamma \sigma}  \frac{\partial \xi_j^\sigma }{\partial z^\delta} dz^\delta    \right)    \right)\\
&= i_{ \left( \sum_{\eta=1}^n V_j^{\alpha \eta} \frac{\partial}{\partial z^\eta} \right) \wedge \left( \sum_{c=1}^n V_j^{\beta c}\frac{\partial}{\partial z^c} \right)  } \left(  \sum_{\delta, \sigma, a, b=1}^n \frac{\partial w_j^{\gamma \sigma}}{\partial \xi_j^a}\frac{\partial \xi_j^a}{\partial z^b} \frac{\partial \xi_j^\sigma}{\partial z^\delta} dz^b \wedge dz^\delta            \right) +  i_{ \left( \sum_{\eta=1}^n V_j^{\alpha \eta} \frac{\partial}{\partial z^\eta} \right) \wedge \left( \sum_{c=1}^n V_j^{\beta c}\frac{\partial}{\partial z^c} \right)  } \left(\sum_{\delta, \sigma, b=1}^n w_j^{\gamma \sigma} \frac{\partial^2 \xi_j^\sigma}{\partial z^b \partial z^\delta} dz^b \wedge dz^\delta     \right) \notag \\
&= \sum_{\eta, c , \sigma , a=1}^n V_j^{\alpha\eta} V_j^{\beta c} \frac{\partial w_j^{\gamma \sigma}}{\partial \xi_j^a} \frac{\partial \xi_j^a}{\partial z^\eta}\frac{\partial \xi_j^\sigma}{\partial z^c} - \sum_{\eta, c , \sigma , a=1}^n V_j^{\alpha\eta} V_j^{\beta c} \frac{\partial w_j^{\gamma \sigma}}{\partial \xi_j^a} \frac{\partial \xi_j^a}{\partial z^c}\frac{\partial \xi_j^\sigma}{\partial z^\eta} = 0, \notag
\end{align}}}
Then $(\ref{osi5})$ and $(\ref{osi6})$ imply that
\begin{align*}
i_{V_j^\alpha(z,t) \wedge V_j^\beta(z,t) }\left( \partial W_j^\gamma (z,t) \right) = 0
\end{align*}
Then
\begin{align*}
&i_{X_j^\alpha(t)\wedge X_j^\beta(t)}\left( \partial Y_j^\gamma(t) \right) = i_{\left( \sum_{\xi=1}^p \left( c_j^{-1} \right)^{\alpha \xi}  V_j^\xi \right)\wedge \left(\sum_{\eta=1}^p \left( c_j^{-1} \right)^{\beta \eta} V_j^\eta \right) }\left( \partial \left( \sum_{\delta=1}^q \left( g^{-1} \right)^{ \gamma \delta} W_j^\delta \right)      \right)\\
&=\sum_{\xi, \eta=1}^p\sum_{\delta=1}^q  \left( c_j^{-1} \right)^{\alpha \xi} \left( c_j^{-1} \right)^{\beta \eta} \left(  g_j^{-1}\right)^{\gamma \delta} i_{V_j^\xi \wedge V_j^\eta}\left( \partial W_j^\delta \right) + \sum_{\xi, \eta, \delta=1}^q \left( c_j^{-1} \right)^{\alpha \xi} \left( c_j^{-1} \right)^{\beta \eta} \left( g_j^{-1} \right)^{\gamma \delta} i_{V_j^\xi \wedge V_j^\eta}\left( \partial  \left( g_j^{-1} \right)^{\gamma \delta}  \wedge W_j^\delta \right)\\
&=0
\end{align*}
Hence we have
\begin{align}\label{osi7}
i_{X_j^\alpha(t)\wedge X_j^\beta(t)}\left( \partial Y_j^\gamma(t) \right) =0
\end{align}

Then from $(\ref{oi26})-(\ref{oi29})$ and $(\ref{oci6})-(\ref{oci7})$ and $(\ref{osi8})$ and $(\ref{osi7})$, we have the integrability conditions:
\begin{align}
&\bar{\partial} \varphi - \frac{1}{2}\left[\varphi, \varphi \right]=0 \\
&\left[ X_j^\alpha, X_j^\beta\right]= \sum_{\gamma=1}^p K_{j\alpha\beta}^\gamma X_j^\gamma \\
&\bar{\partial} X_j^\alpha - \left[ \varphi, X_j^\alpha \right]+ \sum_{\eta=1}^p U_j^{\alpha \eta} X_j^\eta=0 \\
&\bar{\partial} U_j^{\alpha \delta}- \left[ \varphi, U_j^{\alpha \delta} \right]  + \sum_{\sigma=1}^p U_j^{\alpha \sigma} \wedge U_j^{\sigma \delta} \\
&i_{X_j^\alpha(t)\wedge X_j^\beta(t)}\left( \partial Y_j^\gamma(t) \right) =0 \label{osi9}\\
&\bar{\partial} Y_j^\alpha- \mathcal{L}_{\varphi}(Y_j^\alpha) + \sum_{\eta=1}^q Z_j^{\alpha \eta} Y_j^\eta =0  \\
&\bar{\partial} Z_j^{\alpha \delta}- \left[ \varphi, Z_j^{\alpha \delta} \right]  + \sum_{\sigma=1}^p Z_j^{\alpha \sigma} \wedge Z_j^{\sigma \delta} =0  \\
& i_{X_j^\alpha(t)}\left( Y_j^\beta(t) \right) = 0
\end{align}
By taking the derivative of $(\ref{osi9})$ with respect to $t$,
\begin{align}\label{its27}
i_{\frac{\partial X_j^\alpha}{\partial t} \wedge X_j^\beta}\left(\partial Y_j^\gamma \right) + i_{X_j^\alpha \wedge \frac{\partial X_j^\beta}{\partial t}}\left( \partial Y_j^\gamma \right) + i_{X_j^\alpha \wedge X_j^\beta}\left(\partial \left( \frac{\partial Y_j^\gamma}{\partial t} \right) \right) = 0
\end{align}

By setting $t=0$, we have
\begin{align}\label{itt5}
i_{\frac{\partial X_j^\alpha}{\partial t}|_{t=0} \wedge T_{0j}^\beta}\left( dw_{0j}^\gamma \right) + i_{T_{0j}^\alpha \wedge \frac{\partial X_j^\beta}{\partial t}|_{t=0}}\left( dw_{0j}^\gamma \right) + i_{T_{0i}^\alpha \wedge T_{0i}^\beta}\left( \partial \left( \frac{\partial Y_j^\gamma}{\partial t}|_{t=0} \right) \right) = 0
\end{align}
 Since we have $i_{\frac{\partial X_j}{\partial t}|_{t=0}}\left( w_{0j}^\beta \right) + i_{T_{0j}^\alpha}\left( \frac{\partial Y_j^\beta}{\partial t}|_{t=0} \right)=0$ from $(\ref{it17})$, we have from $(\ref{itt5})$
\begin{align}
&-i_{\frac{\partial X_j^\alpha}{\partial t}|_{t=0}}  i_{T_{0j}^\beta}\left( dw_{0j}^\gamma \right) + i_{ \frac{\partial X_j^\beta}{\partial t}|_{t=0}} i_{T_{0j}^\alpha}\left( dw_{0j}^\gamma \right) + i_{T_{0i}^\alpha \wedge T_{0i}^\beta}\left( \partial \left( \frac{\partial Y_j^\gamma}{\partial t}|_{t=0} \right) \right) =  0 \notag\\
&\iff - i_{\frac{\partial X_j^\alpha}{\partial t}|_{t=0} } \left( \sum_{\eta=1}^q a_{i\beta}^{\gamma \eta} w_{0j}^\eta    \right) + i_{\frac{\partial X_j^\beta}{\partial t}|_{t=0}}\left( \sum_{\eta=1}^q a_{i\alpha}^{\gamma \eta } w_{0j}^\eta \right)+ i_{T_{0i}^\alpha \wedge T_{0i}^\beta}\left( \partial \left( \frac{\partial Y_j^\gamma}{\partial t}|_{t=0} \right) \right)=0\notag\\
&\iff \sum_{\eta=1}^q a_{i\beta}^{\gamma \eta} i_{T_{0j}^\alpha}\left( \frac{\partial Y_j^\eta}{\partial t}|_{t=0}  \right) - \sum_{\eta=1}^q a_{i\alpha}^{\gamma \eta} i_{T_{0j}^\beta}\left( \frac{\partial Y_j^\eta}{\partial t}|_{t=0} \right) +  i_{T_{0i}^\alpha \wedge T_{0i}^\beta}\left( \partial \left( \frac{\partial Y_j^\gamma}{\partial t}|_{t=0} \right) \right) =0 \label{itt6}
\end{align}
Then $(\ref{oi42})-(\ref{oi46}), (\ref{itt7}),(\ref{it22})$, $(\ref{it17})$ and $(\ref{itt6})$ implies that
{\Small{\begin{align*}
\left(\frac{\partial X}{\partial t}|_{t=0}, \frac{\partial Y}{\partial t}|_{t=0}, \left\{\left(\frac{\partial \varphi}{\partial t}|_{t=0}, \frac{\partial U_j}{\partial t}|_{t=0}, \frac{\partial Z_j}{\partial t}|_{t=0} \right) \right\} \right)\in A^{0,0}\left( M, \mathscr{H}om_{\mathcal{O}_M}\left( \Theta_{\mathcal{F}_0}, \Theta_M \right)\right) \bigoplus A^{0,0}\left( M, \mathscr{H}om_{\mathcal{O}_M}\left( \mathcal{N}_{\mathcal{F}_0}^*, \Omega_M^1 \right) \right) \bigoplus A^{0,1}\left( M, \mathcal{E}_{\mathcal{F}_0} \right) 
\end{align*}}}
defines a $1$-cocycle in the following bicomplex
{\Tiny{\begin{equation}\label{itt8}
\begin{CD}
\cdots\\
@AAA\\
\frac{A^{0,0}\left( M, \left(\bigwedge^3 \Theta_{\mathcal{F}_0}\right) \otimes \Theta_M \right)}{A^{0,0}\left(M, \left(\bigwedge^3 \Theta_{\mathcal{F}_0} \right) \otimes \Theta_{\mathcal{F}_0} \right)}\bigoplus A^{0,0}\left(M, \left(\mathcal{N}_{\mathcal{F}_0}^* \right)^* \otimes \bigwedge^3 \Theta_{\mathcal{F}_0}^* \right)\bigoplus A^{0,0}\left( M, \left( \mathcal{N}_{\mathcal{F}_0}^* \right)^*\otimes \bigwedge^2 \Theta_{\mathcal{F}_0} \right) @>-\bar{\partial}>> \cdots \\
@AAA @AAA  \\
\frac{A^{0,0}\left( M, \left(\bigwedge^2 \Theta_{\mathcal{F}_0}\right) \otimes \Theta_M \right)}{A^{0,0}\left(M, \left(\bigwedge^2 \Theta_{\mathcal{F}_0} \right) \otimes \Theta_{\mathcal{F}_0} \right)}\bigoplus A^{0,0}\left(M, \left(\mathcal{N}_{\mathcal{F}_0}^* \right)^* \otimes \bigwedge^2 \Theta_{\mathcal{F}_0}^* \right)\bigoplus A^{0,0}\left( M, \left( \mathcal{N}_{\mathcal{F}_0}^* \right)^*\otimes \Theta_{\mathcal{F}_0} \right)@>\bar{\partial}>> \cdots \\
@AAA @AAA \\
A^{0,0}\left( M, \Theta_{\mathcal{F}_0}^*\otimes \Theta_M \right) \bigoplus A^{0,0}\left( M, \left( \mathcal{N}_{\mathcal{F}_0}^* \right)^* \otimes \Omega_M^1 \right) @> -\bar{\partial} >> A^{0,1}\left( M, \Theta_{\mathcal{F}_0}^*\otimes \Theta_M \right) \bigoplus A^{0,1}\left( M, \left( \mathcal{N}_{\mathcal{F}_0}^* \right)^* \otimes \Omega_M^1 \right) @>\bar{\partial}>>  \\
@AAA @AAA  \\
A^{0,0}\left( M, \mathcal{E}_{\mathcal{F}_0} \right) @>\bar{\partial}>> A^{0,1}\left( M, \mathcal{E}_{\mathcal{F}_0} \right) @>-\bar{\partial}>> 
\end{CD}
\end{equation}}}

We will denote the $i$-th cohomology group of the above double complex associated to the Dolbeault resolution of $\mathcal{E}_{\mathcal{F}_0}'^\bullet$ by $\textnormal{S}'^i$. Then we have $\textnormal{S}'^i\cong \mathbb{H}^i\left(M, \mathcal{E}_{\mathcal{F}_0}'^\bullet \right),i\geq 0$. In the next section \ref{nii9} we will explicitly describe the correspondence between the Dolbeault resolution and \v Cech resolution
\begin{align*}
\textnormal{S}'^i \cong \mathbb{H}^i\left( M, \mathcal{E}_{\mathcal{F}_0}'^\bullet \right),\,\,\,\,\,\textnormal{for}\,\,\,  i = 0,1,2
\end{align*}

Then we have the foliated Kodaira-Spencer map in terms of both tangent and cotangent sheaves
{\small{\begin{align}
\varphi_0: T_0 B \to & \mathbb{H}^1\left( M, \mathcal{E}_{\mathcal{F}_0}'^\bullet \right)  \label{a32}\\
 \frac{\partial}{\partial t} \mapsto &\left( \frac{\partial X}{\partial t}|_{t=0}, \frac{\partial Y}{\partial t}|_{t=0}, \left\{\left(\frac{\partial \varphi}{\partial t}|_{t=0}, \frac{\partial U_j}{\partial t}|_{t=0}, \frac{\partial Z_j}{\partial t}|_{t=0} \right) \right\} \right) \notag \\
 &\in  A^{0,0}\left( M, \mathscr{H}om_{\mathcal{O}_M}\left( \Theta_{\mathcal{F}_0}, \Theta_M \right)\right) \bigoplus A^{0,0}\left( M, \mathscr{H}om_{\mathcal{O}_M}\left( \mathcal{N}_{\mathcal{F}_0}^*, \Omega_M^1 \right) \right) \bigoplus A^{0,1}\left( M, \mathcal{E}_{\mathcal{F}_0} \right) \notag
\end{align}}}

\subsection{Identification of obstruction by using the Dolbeault type bicomplex associated to $\mathcal{E}_{\mathcal{F}_0}'^\bullet$}\

Now assume that the foliated analytic family $\left( \mathcal{M}, \Theta_{\mathcal{F}}, \mathcal{N}_\mathcal{F}^*, B, \pi \right)$ with $\Theta_{\mathcal{F}}, \mathcal{N}_\mathcal{F}^*$ locally free is an one-parameter family, i.e. $\Delta \subset \mathbb{C}$. We keep the notations in subsection \ref{itt10}. Then by taking the derivative of $(\ref{its27})$ with respect to $t$ again, we obtain
{\tiny{\begin{align*}
i_{\frac{\partial^2 X_j^\alpha}{\partial t} \wedge X_j^\beta}\left(\partial Y_j^\gamma \right) + 2 i_{\frac{\partial X_j^\alpha}{\partial t} \wedge \frac{\partial X_j^\beta}{\partial t}}\left(\partial Y_j^\gamma \right) + 2 i_{\frac{\partial X_j^\alpha }{\partial t}\wedge X_j^\beta}\left(\partial \left( \frac{\partial Y_j^\gamma}{\partial t} \right) \right) + i_{X_j^\alpha \wedge \frac{\partial^2 X_j^\beta}{ \partial t^2}}\left(\partial Y_j^\gamma \right) + 2 i_{X_j^\alpha \wedge \frac{\partial X_j^\beta}{\partial t}}\left( \partial \left( \frac{\partial Y_j^\gamma}{\partial t} \right) \right) + i_{X_j^\alpha \wedge X_j^\beta}\left( \partial \left(  \frac{\partial^2 Y_j^\gamma}{\partial t^2} \right) \right) = 0
\end{align*}}}
and setting $t=0$, we obtain
\begin{align*}
&i_{\frac{\partial^2 X_j^\alpha}{\partial t^2}|_{t=0} \wedge T_{0j}^\beta}\left( dw_{0j}^\gamma \right) + 2 i_{\frac{\partial X_j^\alpha}{\partial t}|_{t=0} \wedge \frac{\partial X_j^\beta}{\partial t}|_{t=0}}\left( dw_{0j}^\gamma \right) + 2 i_{\frac{\partial X_j^\alpha}{\partial t}|_{t=0} \wedge T_{0j}^\beta} \left( \partial \left( \frac{\partial Y_j^\gamma}{\partial t}|_{t=0} \right) \right) + 2 i_{T_{0j}^\alpha \wedge \frac{\partial X_j^\beta}{\partial t}|_{t=0}} \left( \partial \left( \frac{\partial Y_j^\gamma}{\partial t} |_{t=0} \right) \right)  \\
&+ i_{T_{0j}^\alpha \wedge \frac{\partial^2 X_j^\beta}{\partial t^2}|_{t=0}}\left(  dw_{0j}^\gamma \right) + i_{T_{0j}^\alpha \wedge T_{0j}^\beta}\left( \partial \left( \frac{\partial^2 Y_j^\gamma}{\partial t^2}|_{t=0} \right) \right) = 0
\end{align*}
which is equivalent to
\begin{align*}
&- \sum_{\eta=1}^q a_{0j\beta}^{\gamma \eta} i_{\frac{\partial^2 X_j^\alpha}{\partial t^2} |_{t=0}} \left(  w_{0j}^\eta \right) + 2 i_{\frac{\partial X_j^\alpha}{\partial t}|_{t=0} \wedge \frac{\partial X_j^\beta}{\partial t}|_{t=0}}\left( dw_{0j}^\gamma \right) + 2 i_{\frac{\partial X_j^\alpha}{\partial t}|_{t=0} \wedge T_{0j}^\beta} \left( \partial \left( \frac{\partial Y_j^\gamma}{\partial t}|_{t=0} \right) \right) + 2 i_{T_{0j}^\alpha \wedge \frac{\partial X_j^\beta}{\partial t}|_{t=0}} \left( \partial \left( \frac{\partial Y_j^\gamma}{\partial t} |_{t=0} \right) \right)\\
& + \sum_{\eta=1}^q a_{0j\alpha}^{\gamma \eta} i_{\frac{\partial^2 X_j^\beta}{\partial t^2}|_{t=0}}\left( w_{0j}^\eta \right) + i_{T_{0j}^\alpha \wedge T_{0j}^\beta}\left( \partial \left( \frac{\partial^2 Y_j^\gamma}{\partial t^2}|_{t=0} \right) \right) =0,
\end{align*}
so that we have
{\small{\begin{align}
&\sum_{\eta=1}^q a_{0j\beta}^{\gamma \eta} i_{T_{0j}^\alpha}\left( \frac{\partial^2 Y_j^\eta}{\partial t^2}|_{t=0} \right) - \sum_{\eta=1}^q a_{0j\alpha}^{\gamma \eta} i_{T_{0j}^\beta}\left( \frac{\partial Y_j^\eta}{\partial t}|_{t=0} \right) + i_{T_{0j}^\alpha \wedge T_{0j}^\beta}\left( \partial \left( \frac{\partial^2 Y_j^\gamma}{\partial t^2}|_{t=0}   \right) \right)\label{itt12}\\
&= - 2 \sum_{\eta=1}^q a_{0j\beta}^{\gamma \eta}i_{\frac{\partial X_j^\alpha}{\partial t}|_{t=0}}\left( \frac{\partial Y_j^\eta}{\partial t}|_{t=0} \right) - 2 i_{\frac{\partial X_j^\alpha}{\partial t}|_{t=0} \wedge \frac{\partial X_j^\beta}{\partial t}|_{t=0}}\left( dw_{0j}^\gamma \right) - 2 i_{\frac{\partial X_j^\alpha}{\partial t}|_{t=0} \wedge T_{0j}^\beta} \left( \partial \left( \frac{\partial Y_j^\gamma}{\partial t}|_{t=0} \right) \right) - 2 i_{T_{0j}^\alpha \wedge \frac{\partial X_j^\beta}{\partial t}|_{t=0}} \left( \partial \left( \frac{\partial Y_j^\gamma}{\partial t} |_{t=0} \right) \right) \notag\\
&+ 2 \sum_{\eta=1}^q a_{0j\alpha}^{\gamma \eta} i_{\frac{\partial X_j^\beta}{\partial t}|_{t=0}}\left( \frac{\partial Y_j^\eta}{\partial t}|_{t=0} \right) \notag
\end{align}}}

\begin{lemma}\label{its29}
We keep the notations in \textnormal{Lemma \ref{its28}}. For any $1$-cocycle $\left(\Phi=\left\{ \Phi_j\right\},  \Psi=\left\{\Psi_j \right\}, \left\{ \left( \psi, a_j,  k_j   \right) \right\} \right)\in A^{0,0}\left( M, \mathscr{H}om_{\mathcal{O}_M}\left( \Theta_{\mathcal{F}_0}, \Theta_M \right) \right) \bigoplus A^{0,0}\left( M, \mathscr{H}om_{\mathcal{O}_M}\left( \mathcal{N}_{\mathcal{F}_0}^*, \Omega_M^1 \right) \right) \bigoplus A^{0,1}\left( M, \mathcal{E}_{\mathcal{F}_0} \right)$ satisfying
\begin{enumerate}
\item $\left[ T_{0j}^\alpha, \Phi_j\left( T_{0j}^\beta \right) \right]- \left[ T_{0j}^\beta, \Phi_j\left( T_{0j}^\alpha \right) \right] - \sum_{\gamma=1}^p g_{0j\alpha\beta}^\gamma \Phi_j\left(T_{0j}^\gamma \right)= K_{j\alpha\beta}$ for some $K_{j\alpha\beta}\in \Gamma\left( U_j, \mathcal{A}^{0,0}\left( \Theta_{\mathcal{F}_0} \right) \right),\alpha,\beta=1,...,p$.
\item $-\bar{\partial}\Phi_j\left( T_{0j}^\alpha \right) + \left[ \psi, T_{0j}^\alpha \right] - \sum_{\beta=1}^p a_j^{\alpha\beta} T_{0j}^\beta = 0 ,\alpha=1,...,p$.
\item $\sum_{\eta=1}^q a_{0j\beta}^{\gamma \eta} i_{T_{0j}^\alpha}\left(\Psi_j^\eta \right) - \sum_{\eta=1}^q a_{0j\alpha}^{\gamma \eta} i_{T_{0j}^\beta}\left(\Psi_j^\eta \right)+ i_{T_{0j}^\alpha \wedge T_{0j}^\beta}\left(\partial \Psi_j^\gamma \right)=0, \alpha,\beta=1,...,p, \gamma=1,...q$.\\
\item $- \bar{\partial} \Psi_j^\alpha + \mathcal{L}_\psi\left( w_{0j}^\alpha \right) - \sum_{\beta=1}^q k_j^{\alpha\beta} w_{0j}^\beta, \alpha=1,...,q$.\\
\item $-\bar{\partial}\psi=0, -\bar{\partial} a_j^{\alpha\beta},\alpha,\beta=1,...,p$, and $-\bar{\partial} k_j^{\gamma \eta}=0,\gamma,\eta=1,...,q$.
\item $i_{\Phi_j^\alpha}\left( w_{0j}^\beta \right) + i_{T_{0j}^\alpha}\left(\Psi_j^\beta \ \right)=0,\alpha=1,...,p,  \beta=1,...,q$.
\end{enumerate}
We define $\tilde{\rho}=\left\{\tilde{\rho}_i\right\}\in A^{0,0}\left( M, \mathscr{H}om_{\mathcal{O}_M}\left( \bigwedge^2 \Theta_{\mathcal{F}_0}, \Theta_M \right) \right)$, and define $\rho$ to be the image of $\tilde{\rho}$ in $\frac{A^{0,0}\left( M, \mathscr{H}om_{\mathcal{O}_M}\left( \bigwedge^2 \Theta_{\mathcal{F}_0}, \Theta_M \right) \right)}{A^{0,0}\left( M, \mathscr{H}om_{\mathcal{O}_M}\left( \bigwedge^2 \Theta_{\mathcal{F}_0}, \Theta_{\mathcal{F}_0} \right) \right)}$ and define $\eta=\{\eta_i\}\in A^{0,0}\left( M,  \mathscr{H}om_{\mathcal{O}_M}\left(\Theta_{\mathcal{F}_0}, \Theta_M \right) \right)$ and define $\kappa=\left\{ \kappa_i\right\} \in A^{0,1}\left( M, \mathscr{H}om_{\mathcal{O}_M}\left( \mathcal{N}_{\mathcal{F}_0}^*, \Omega_M^1 \right) \right)$ and $p\in A^{0,0}\left( M, \mathscr{H}om_{\mathcal{O}_M}\left( \mathcal{N}_{\mathcal{F}_0}^*, \Theta_{\mathcal{F}_0}^* \right) \right)$ as in \textnormal{Lemma \ref{its28}}. On the other hand, we define an element $n_j\in \Gamma\left( U_j, \mathcal{A}^{0,0}\left(\mathscr{H}om_{\mathcal{O}_M}\left( \mathcal{N}_{\mathcal{F}_0}^*, \bigwedge^2 \Theta_{\mathcal{F}_0}^* \right)  \right) \right)$ by
\begin{align*}
&n_j:\Gamma\left( U_j, \mathcal{N}_{\mathcal{F}_0}^* \right) \to \Gamma\left( U_j, \mathcal{A}^{0,0}\left( \bigwedge^2 \Theta_{\mathcal{F}_0}^* \right) \right) \\
&w_{0j}^\gamma \mapsto \left( T_{0j}^\alpha \wedge T_{0j}^\beta \mapsto -2 \sum_{\eta=1}^q a_{0j\beta}^{\gamma \eta} i_{\Phi_j^\alpha}\left( \Psi_j^\eta \right) + 2 \sum_{\eta=1}^q a_{0j\alpha}^{\gamma \eta} i_{\Phi_j^\beta}\left(\Psi_j^\eta \right)-2 i_{\Phi_j^\alpha \wedge \Phi_j^\beta}\left( dw_{0j}^\gamma \right)- 2i_{\Phi_j^\alpha \wedge T_{0j}^\beta }\left(\partial \Psi_j^\gamma \right) -2 i_{T_{0j}^\alpha \wedge \Phi_j^\beta}\left( \partial \Psi_j^\gamma \right)\right)
\end{align*}
and linearly extends to $\Gamma\left(U_j, \mathcal{N}_{\mathcal{F}_0}^* \otimes \bigwedge^2 \Theta_{\mathcal{F}_0} \right)$. Then $n:=\{n_j\}$ defines a global section in $A^{0,0}\left( M, \mathscr{H}om_{\mathcal{O}_M}\left( \mathcal{N}_{\mathcal{F}_0}^*, \bigwedge^2 \Theta_{\mathcal{F}_0}^* \right) \right)$. Then
\begin{align}\label{itt11}
&\left( \rho, n, p,  \eta,  \kappa, \left\{- \left[ \left(\psi, a_j, k_j \right), \left( \psi, a_j, k_j \right) \right]     \right\}  \right) \\
&\in \frac{A^{0,0}\left( M, \mathscr{H}om_{\mathcal{O}_M}\left( \bigwedge^2 \Theta_{\mathcal{F}_0} , \Theta_M \right) \right)}{A^{0,0}\left( M, \mathscr{H}om_{\mathcal{O}_M}\left( \bigwedge^2 \Theta_{\mathcal{F}_0}, \Theta_{\mathcal{F}_0} \right) \right) } \bigoplus A^{0,0}\left( M, \mathscr{H}om_{\mathcal{O}_M}\left( \mathcal{N}_{\mathcal{F}_0}^*, \bigwedge^2 \Theta_{\mathcal{F}_0}^* \right) \right) \bigoplus A^{0,0}\left( M, \mathscr{H}om_{\mathcal{O}_M}\left( \mathcal{N}_{\mathcal{F}_0}^*, \Theta_{\mathcal{F}_0}^* \right) \right) \notag\\
&\bigoplus A^{0,1}\left( M, \mathscr{H}om_{\mathcal{O}_M}\left( \Theta_{\mathcal{F}_0}, \Theta_M \right) \right) \bigoplus A^{0,1}\left( M, \mathscr{H}om_{\mathcal{O}_M}\left( \mathcal{N}_{\mathcal{F}_0}^*, \Omega_M^1 \right) \right) \bigoplus A^{0,2}\left( M, \mathcal{E}_{\mathcal{F}_0} \right) \notag
\end{align}
defines a $2$-cocycle in the Dolbeault resolution of $\mathcal{E}_{\mathcal{F}_0}'^\bullet$.

\end{lemma}

\begin{remark}
We show that $\{n_j\}$ defines a global section in $A^{0,0}\left( M, \mathscr{H}om_{\mathcal{O}_M}\left( \mathcal{N}_{\mathcal{F}_0}^*, \bigwedge^2 \Theta_{\mathcal{F}_0}^* \right) \right)$. In fact,
{\small{\begin{align*}
&-2 \sum_{\eta=1}^q a_{0j\beta}^{\gamma \eta} i_{\Phi_j^\alpha}\left( \Psi_j^\eta \right) + 2 \sum_{\eta=1}^q a_{0j\alpha}^{\gamma \eta} i_{\Phi_j^\beta}\left(\Psi_j^\eta \right)-2 i_{\Phi_j^\alpha \wedge \Phi_j^\beta}\left( dw_{0j}^\gamma \right)- 2i_{\Phi_j^\alpha \wedge T_{0j}^\beta }\left(\partial \Psi_j^\gamma \right) -2 i_{T_{0j}^\alpha \wedge \Phi_j^\beta}\left( \partial \Psi_j^\gamma \right)\\
&= - 2\sum_{\eta, \delta =1}^q\sum_{\xi=1}^p a_{0j\beta}^{\gamma \eta} r_{0jk}^{\alpha \xi} h_{0jk}^{\eta \delta} i_{\Phi_k^\xi} \left( \Psi_k^\delta \right) + 2 \sum_{\eta , \delta =1}^q \sum_{\xi=1}^p a_{0j\alpha}^{\gamma \eta} r_{0jk}^{\beta \xi} h_{0jk}^{\eta \delta} i_{\Phi_k^\xi}\left(\Psi_k^\delta \right) - 2 \sum_{\xi, \sigma=1}^p \sum_{\delta=1}^q r_{0jk}^{\alpha \xi} r_{0jk}^{\beta \sigma} i_{\Phi_k^\xi \wedge \Phi_k^\sigma}\left( d h_{0jk}^{\gamma \delta} \wedge w_{0k}^\delta + h_{0jk}^{\gamma \delta} dw_{0k}^\delta \right)\\
&-2 \sum_{\xi, \sigma=1}^p \sum_{\delta=1}^q  r_{0jk}^{\alpha \xi} r_{0jk}^{\beta \sigma} i_{\Phi_k^\xi \wedge T_{0k}^\sigma}\left( d h_{0jk}^{\gamma \delta} \wedge \Phi_k^\delta  + h_{0jk}^{\gamma \delta} \partial \Phi_k^\delta \right) - 2 \sum_{\xi, \sigma=1}^p \sum_{\delta=1}^q r_{0jk}^{\alpha \xi} r_{0jk}^{\beta \sigma} i_{T_{0k}^\xi \wedge \Phi_k^\sigma} \left( d h_{0jk}^{\gamma \delta} \wedge \Phi_k^\delta + h_{0jk}^{\gamma \delta} \partial \Phi_k^\delta \right)\\
&=- \cancel{2\sum_{\xi=1}^p \sum_{\delta=1}^q r_{0jk}^{\alpha \xi} i_{T_{0j}^\beta} \left( d h_{0jk}^{\gamma \delta} \right) i_{\Phi_k^\xi}\left( \Psi_k^\delta  \right)} -2 \sum_{\xi, \sigma=1}^p \sum_{\eta, \delta=1}^q r_{0jk}^{\alpha \xi} r_{0jk}^{\beta \sigma} h_{0jk}^{\gamma\eta} a_{0k\sigma}^{\eta \delta} i_{\Phi_k^\xi}\left( \Psi_k^\delta \right) + \bcancel{2\sum_{\xi=1}^p \sum_{\delta=1}^q r_{0jk}^{\beta \xi} i_{T_{0j}^\alpha} \left( d h_{0jk}^{\gamma \delta} \right) i_{\Phi_k^\xi}\left( \Psi_k^\delta  \right)} \\
&+2 \sum_{\xi, \sigma=1}^p \sum_{\eta, \delta=1}^q r_{0jk}^{\beta \xi} r_{0jk}^{\alpha \sigma} h_{0jk}^{\gamma\eta} a_{0k\sigma}^{\eta \delta} i_{\Phi_k^\xi}\left( \Psi_k^\delta \right)- \xcancel{2 \sum_{\xi, \sigma=1}^p \sum_{\delta=1}^q r_{0jk}^{\alpha \xi} r_{0jk}^{\beta \sigma} i_{\Phi_k^\xi}\left(dh_{0jk}^{\gamma \delta} \right) i_{\Phi_k^\sigma}\left(w_{0k}^\delta \right) } + \cancel{2\sum_{\xi, \sigma=1}^p \sum_{\delta=1}^q r_{0jk}^{\alpha \xi} r_{0jk}^{\beta \sigma} i_{\Phi_k^\sigma}\left(dh_{0jk}^{\gamma \delta} \right) i_{\Phi_k^\xi}\left(w_{0k}^\delta \right)}^1 \\
&-2 \sum_{\xi, \sigma=1}^p \sum_{\delta=1}^q r_{0jk}^{\alpha \xi} r_{0jk}^{\beta \sigma} h_{0jk}^{\gamma \delta} i_{\Phi_k^\xi \wedge \Phi_k^\sigma}\left( d w_{0k}^\delta \right) - \xcancel{2 \sum_{\xi, \sigma=1}^p r_{0jk}^{\alpha \xi} r_{0jk}^{\beta \sigma} i_{\Phi_k^\xi}\left( d h_{0jk}^{\gamma \delta}  \right) i_{T_{0k}^\sigma}\left(\Phi_k^\delta \right)}+ 2 \cancel{\sum_{\xi, \sigma=1}^p r_{0jk}^{\alpha \xi} r_{0jk}^{\beta \sigma} i_{T_{0k}^\sigma}\left( d h_{0jk}^{\gamma \delta}  \right) i_{\Phi_k^\xi}\left(\Psi_k^\delta \right)}\\
&- 2\sum_{\xi, \sigma=1}^p \sum_{\delta=1}^q r_{0jk}^{\alpha \xi} r_{0jk}^{\beta \sigma} h_{0jk}^{\gamma \delta} i_{\Phi_k^\xi \wedge T_{0k}^\sigma}\left( \partial \Psi_k^{\delta}\right) -  \bcancel{2 \sum_{\xi, \sigma=1}^p \sum_{\delta=1}^q r_{0jk}^{\alpha \xi} r_{0jk}^{\beta \sigma} i_{T_{0k}^\xi}\left( dh_{0jk}^{\gamma \delta} \right) i_{\Phi_k^\sigma}\left( \Psi_k^\delta\right) }+  \cancel{ 2 \sum_{\xi, \sigma=1}^p \sum_{\delta=1}^q r_{0jk}^{\alpha \xi} r_{0jk}^{\beta \sigma} i_{\Phi_k^\delta}\left( dh_{0jk}^{\gamma \delta} \right) i_{T_{0j}^\xi}\left( \Psi_k^\delta\right)}^1 \\
&- 2 \sum_{\xi, \sigma=1}^p \sum_{\delta=1}^q r_{0jk}^{\alpha \xi} r_{0jk}^{\beta \sigma} h_{0jk}^{\gamma \delta} i_{T_{0k}^\xi \wedge \Phi_k^\sigma}\left(\partial \Psi_k^\delta \right)\\
&= \sum_{\xi, \sigma=1}^p \sum_{\delta=1}^q r_{0jk}^{\alpha \xi} r_{0jk}^{\beta \sigma} h_{0jk}^{\gamma \delta}\left(  -2 \sum_{\eta=1}^q a_{0k\sigma}^{\eta \delta} i_{\Phi_k^\xi} \left(\Psi_k^\delta \right) +2 \sum_{\eta=1}^q a_{0k\xi}^{\eta \delta} i_{\Phi_k^\sigma}\left( \Psi_k^\delta \right) -2 i_{\Phi_k^\xi \wedge \Phi_k^\sigma}\left(  dw_{0k}^\delta \right)  - 2 i_{\Phi_k^\xi \wedge T_{0k}^\sigma}\left(\partial \Psi_k^\delta \right) - 2 i_{T_{0k}^\xi \wedge \Phi_k^\sigma}\left( \partial \Psi_k^\delta \right) \right)
\end{align*}}}
\end{remark}

\begin{proof}[Proof of \textnormal{Lemma \ref{its29}}]
We keep the notations in the proof of Lemma \ref{its28}. Then it remains to show that $\hat{\tilde{E}}_2'\left( n \right)=0$ and $\bar{\partial} n + \hat{\tilde{E}}_1\left( \kappa\right)=0$ and $\hat{\beta}_2\left(\rho \right) + n - \hat{\tilde{E}}_1\left( p \right)=0$. In fact, as in the same way with $(\ref{us26}),(\ref{us25}),(\ref{its20})$, we have $\hat{\beta}_2\left(\rho \right) + n - \hat{\tilde{E}}_1\left( p \right)=0$ on $M$.  On the other hand, we note that
{\small{\begin{align*}
&-\frac{1}{2}\bar{\partial}n_j= \bar{\partial}\left(\sum_{\eta=1}^q a_{0j\beta}^{\gamma \eta}i_{\Phi_j^\alpha}\left( \Psi_j^\eta  \right) - \sum_{\eta=1}^q a_{0j\alpha}^{\gamma \eta}i_{\Phi_j^\beta}\left( \Psi_j^\eta \right) + i_{\Phi_j^\beta} i_{\Phi_j^\alpha}\left( dw_{0j}^\gamma \right) - i_{\Phi_j^\alpha}i_{ T_{0j}^\beta}\left(\partial \Psi_j^\gamma \right) + i_{ \Phi_j^\beta}i_{T_{0j}^\alpha}\left( \partial \Psi_j^\gamma \right) \right)\\
&= \sum_{\eta=1}^q a_{0j\beta}^{\gamma \eta} i_{\left[ \psi, T_{0j}^\alpha\right] - \sum_{\xi=1}^p a_j^{\alpha \xi} T_{0j}^\xi}\left(\Psi_j^\eta \right) - \sum_{\eta=1}^q a_{0j\beta}^{\gamma \eta} i_{\Phi_j^\alpha}\left( \mathcal{L}_\psi\left( w_{0j}^\eta \right) - \sum_{\sigma=1}^q k_j^{\eta\sigma} w_{0j}^\sigma \right)\\
& -\sum_{\eta=1}^q a_{0j\alpha}^{\gamma \eta} i_{\left[ \psi, T_{0j}^\beta \right] - \sum_{\xi=1}^p a_j^{\beta \xi} T_{0j}^\xi}\left(\Psi_j^\eta \right) + \sum_{\eta=1}^q a_{0j\alpha}^{\gamma \eta} i_{\Phi_j^\beta}\left( \mathcal{L}_\psi\left( w_{0j}^\eta \right) - \sum_{\sigma=1}^q k_j^{\eta\sigma} w_{0j}^\sigma \right) \\
&+ i_{\left[ \psi, T_{0j}^\beta\right]-\sum_{\xi=1}^p a_j^{\beta \xi} T_{0j}^\xi} i_{\Phi_j^\alpha} \left( dw_{0j}^\gamma \right) - i_{\Phi_j^\beta} i_{\left[ \psi, T_{0j}^\alpha \right]- \sum_{\xi=1}^p a_j^{\alpha \xi} T_{0j}^\xi}\left( dw_{0j}^\gamma \right) - i_{\left[ \psi, T_{0j}^\alpha \right] - \sum_{\xi=1}^p a_j^{\alpha \xi} T_{0j}^\xi} i_{T_{0j}^\beta}\left(\partial \Psi_j^\gamma \right)\\
& + i_{\Phi_j^\alpha} i_{T_{0j}^\beta}\left(\partial \left( \mathcal{L}_{\psi}\left( w_{0j}^\gamma \right) - \sum_{\sigma=1}^q k_j^{\gamma \sigma} w_{0j}^\sigma \right) \right) +   i_{\left[ \psi, T_{0j}^\beta \right] - \sum_{\xi=1}^p a_j^{\beta \xi} T_{0j}^\xi} i_{T_{0j}^\alpha}\left(\partial \Psi_j^\gamma \right) - i_{\Phi_j^\beta} i_{T_{0j}^\alpha}\left(\partial \left( \mathcal{L}_{\psi}\left( w_{0j}^\gamma \right) - \sum_{\sigma=1}^q k_j^{\gamma \sigma} w_{0j}^\sigma \right) \right)
\end{align*}}}
We continue the computation. Then
{\small{\begin{align*}
&= \sum_{\eta=1}^q a_{0j\beta}^{\gamma \eta} i_{\left[\psi, T_{0j}^\alpha \right]}\left(\Psi_j^\eta \right) + \cancel{\sum_{\eta=1}^q \sum_{\xi=1}^p a_{0j\beta}^{\gamma \eta} a_j^{\alpha \xi} i_{\Phi_j^\xi}\left( w_{0j}^\eta \right)} - \sum_{\eta=1}^q a_{0j\beta}^{\gamma \eta} i_{\Phi_j^\alpha} \mathcal{L}_{\psi}\left( w_{0j}^\eta \right) + \sum_{\eta,\sigma=1}^q a_{0j\beta}^{\gamma \eta}  k_j^{\eta \sigma} i_{T_{0j}^\alpha}\left( \Psi_j^\sigma \right)\\
&- \sum_{\eta=1}^q a_{0j\alpha}^{\gamma \eta} i_{\left[\psi, T_{0j}^\beta \right]}\left(\Psi_j^\eta \right) - \cancel{\sum_{\eta=1}^q \sum_{\xi=1}^p a_{0j\alpha}^{\gamma \eta} a_j^{\beta \xi} i_{\Phi_j^\xi}\left( w_{0j}^\eta \right)} + \sum_{\eta=1}^q a_{0j\alpha}^{\gamma \eta} i_{\Phi_j^\beta} \mathcal{L}_{\psi}\left( w_{0j}^\eta \right) - \sum_{\eta,\sigma=1}^q a_{0j\alpha}^{\gamma \eta}  k_j^{\eta \sigma} i_{T_{0j}^\beta}\left( \Psi_j^\sigma \right)\\
&+ i_{\left[ \psi, T_{0j}^\beta \right]}i_{\Phi_j^\alpha}\left( dw_{0j}^\gamma \right)  +\cancel{ \sum_{\xi=1}^p \sum_{\eta=1}^q a_j^{\beta \xi} a_{0j\xi}^{\gamma \eta} i_{\Phi_j^\alpha}\left( w_{0j}^\eta \right)} -  i_{\left[ \psi, T_{0j}^\alpha \right]}i_{\Phi_j^\beta}\left( dw_{0j}^\gamma \right)  -  \cancel{ \sum_{\xi=1}^p \sum_{\eta=1}^q a_j^{\alpha \xi} a_{0j\xi}^{\gamma \eta} i_{\Phi_j^\beta}\left( w_{0j}^\eta \right) } \\
&- i_{\left[ \psi, T_{0j}^\alpha \right]} i_{T_{0j}^\beta} \left(\partial \Psi_j^\gamma \right) + \cancel{ \sum_{\xi=1}^p a_j^{\alpha \xi} i_{T_{0j}^\beta \wedge T_{0j}^\xi}\left(\partial \Psi_j^\gamma \right) } + i_{\left[ \psi, T_{0j}^\beta \right]} i_{T_{0j}^\alpha} \left(\partial \Psi_j^\gamma \right) - \cancel{ \sum_{\xi=1}^p a_j^{\beta \xi} i_{T_{0j}^\alpha \wedge T_{0j}^\xi}\left(\partial \Psi_j^\gamma \right) }\\
&- i_{\Phi_j^\alpha}i_{T_{0j}^\beta}\mathcal{L}_\psi\left( dw_{0j}^\gamma \right) - \sum_{\sigma=1}^q i_{\Phi_j^\alpha} i_{T_{0j}^\beta} \left( \partial k_j^{\gamma\sigma} \wedge w_{0j}^\sigma - k_j^{\gamma \sigma} \wedge dw_{0j}^\sigma \right) + i_{\Phi_j^\beta}i_{T_{0j}^\alpha}\mathcal{L}_\psi\left( dw_{0j}^\gamma \right) + \sum_{\sigma=1}^q i_{\Phi_j^\beta} i_{T_{0j}^\alpha} \left( \partial k_j^{\gamma\sigma} \wedge w_{0j}^\sigma - k_j^{\gamma \sigma} \wedge dw_{0j}^\sigma \right)\\
&=  \sum_{\eta=1}^q a_{0j\beta}^{\gamma \eta}\mathcal{L}_\psi i_{T_{0j}^\alpha}\left(\Psi_j^\eta \right) + \sum_{\eta=1}^q a_{0j\beta}^{\gamma \eta} i_{T_{0j}^\alpha}\mathcal{L}_\psi\left( \Psi_j^\eta\right) - \sum_{\eta=1}^q a_{0j\beta}^{\gamma \eta} i_{\Phi_j^\alpha} \mathcal{L}_{\psi}\left( w_{0j}^\eta \right) + \sum_{\eta,\sigma=1}^q a_{0j\beta}^{\gamma \eta}  k_j^{\eta \sigma} i_{T_{0j}^\alpha}\left( \Psi_j^\sigma \right)\\
& - \sum_{\eta=1}^q a_{0j\alpha}^{\gamma \eta}\mathcal{L}_\psi i_{T_{0j}^\beta}\left(\Psi_j^\eta \right) - \sum_{\eta=1}^q a_{0j\alpha}^{\gamma \eta} i_{T_{0j}^\beta}\mathcal{L}_\psi\left( \Psi_j^\eta\right) + \sum_{\eta=1}^q a_{0j\alpha}^{\gamma \eta} i_{\Phi_j^\beta} \mathcal{L}_{\psi}\left( w_{0j}^\eta \right) - \sum_{\eta,\sigma=1}^q a_{0j\alpha}^{\gamma \eta}  k_j^{\eta \sigma} i_{T_{0j}^\beta}\left( \Psi_j^\sigma \right)\\
&+ \cancel{ i_{\left[ \psi, T_{0j}^\beta \right]} i_{\Phi_j^\alpha}\left( dw_{0j}^\gamma \right)} - \cancel{ i_{\left[ \psi , T_{0j}^\alpha \right]}i_{\Phi_j^\beta}\left( dw_{0j}^\gamma \right)} - i_{\left[ \psi, T_{0j}^\alpha \right]} i_{T_{0j}^\beta}\left(\partial \Psi_j^\gamma \right) + i_{\left[ \psi, T_{0j}^\beta \right]} i_{T_{0j}^\alpha}\left( \partial \Psi_j^\gamma \right) - \cancel{ i_{\Phi_j^\alpha} i_{\left[\psi, T_{0j}^\beta \right]}\left( dw_{0j}^\gamma \right) } + i_{\Phi_j^\alpha}\mathcal{L}_\psi\left(\sum_{\eta=1}^q a_{0j\beta}^{\gamma \eta} w_{0j}^\eta \right)\\
&+ \sum_{\sigma=1}^q  i_{T_{0j}^\beta} \left(\partial k_j^{\gamma \sigma} \right)i_{\Phi_j^\alpha}\left( w_{0j}^\sigma  \right) - \sum_{\sigma=1}^q i_{\Phi_j^\alpha}\left( k_j^{\gamma \sigma} \wedge \left(\sum_{\delta=1}^q a_{0j\beta}^{\sigma \delta} w_{0j}^\sigma  \right)\right) + \cancel{ i_{\Phi_j^\beta} i_{\left[\psi, T_{0j}^\alpha \right]}\left( dw_{0j}^\gamma \right) } - i_{\Phi_j^\beta}\mathcal{L}_\psi\left(\sum_{\eta=1}^q a_{0j\alpha}^{\gamma \eta} w_{0j}^\eta \right)\\
& -  \sum_{\sigma=1}^q  i_{T_{0j}^\alpha} \left(\partial k_j^{\gamma \sigma} \right)i_{\Phi_j^\beta}\left( w_{0j}^\sigma  \right) + \sum_{\sigma=1}^q i_{\Phi_j^\beta}\left( k_j^{\gamma \sigma} \wedge \left(\sum_{\delta=1}^q a_{0j\alpha}^{\sigma \delta} w_{0j}^\sigma  \right)\right)
\end{align*}}}
Then we have
\begin{align}\label{its38}
-\frac{1}{2}\bar{\partial} n_j&= \sum_{\eta=1}^q a_{0j\beta}^{\gamma \eta} \mathcal{L}_\psi i_{T_{0j}^\alpha}\left(\Psi_j^\eta \right) + \sum_{\eta=1}^q a_{0j\beta}^{\gamma \eta} i_{T_{0j}^\alpha} \mathcal{L}_\psi \left(  \Psi_j^\eta \right) - \sum_{\eta=1}^q \left[ \psi, a_{0j\beta}^{\gamma \eta} \right] i_{\Phi_j^\alpha}\left( w_{0j}^\eta \right)  + \sum_{\eta, \sigma=1}^q a_{0j\beta}^{\gamma \eta} k_j^{\eta \sigma} i_{T_{0j}^\alpha} \left( \Psi_j^\sigma \right) - i_{\left[ \psi, T_{0j}^\alpha\right]} i_{T_{0j}^\beta}\left(\partial \Psi_j^\gamma \right) \\
&-\sum_{\eta=1}^q a_{0j\alpha}^{\gamma \eta} \mathcal{L}_\psi i_{T_{0j}^\beta}\left(\Psi_j^\eta \right) - \sum_{\eta=1}^q a_{0j\alpha}^{\gamma \eta} i_{T_{0j}^\beta} \mathcal{L}_\psi \left(  \Psi_j^\eta \right) + \sum_{\eta=1}^q \left[ \psi, a_{0j\alpha}^{\gamma \eta} \right] i_{\Phi_j^\beta}\left( w_{0j}^\eta \right)  - \sum_{\eta, \sigma=1}^q a_{0j\alpha}^{\gamma \eta} k_j^{\eta \sigma} i_{T_{0j}^\beta} \left( \Psi_j^\sigma \right) + i_{\left[ \psi, T_{0j}^\beta\right]} i_{T_{0j}^\alpha}\left( \partial \Psi_j^\gamma \right) \notag \\
&- \sum_{\sigma=1}^q i_{T_{0j}^\beta} \left( \partial k_j^{\gamma \sigma} \right) i_{T_{0j}^\alpha}\left(\Psi_j^\sigma \right)  - \sum_{\sigma , \delta =1}^q k_j^{\gamma \sigma} a_{0j\beta}^{\sigma \delta}i_{T_{0j}^\alpha}\left( \Psi_j^\delta \right) + \sum_{\sigma=1}^q i_{T_{0j}^\alpha} \left( \partial k_j^{\gamma \sigma} \right) i_{T_{0j}^\beta}\left(\Psi_j^\sigma \right)  + \sum_{\sigma , \delta =1}^q k_j^{\gamma \sigma} a_{0j\alpha}^{\sigma \delta}i_{T_{0j}^\beta}\left( \Psi_j^\delta \right) \notag
\end{align}
On the other hand, we note that
{\small{\begin{align}\label{its36}
&\sum_{\eta=1}^q a_{0j\beta}^{\gamma \eta}i_{T_{0j}^\alpha}\left(  -\mathcal{L}_{\psi}\left(\Psi_j^\eta \right)   + \sum_{\sigma=1}^q k_j^{\eta \sigma} \wedge \Psi_j^\sigma \right) - \sum_{\eta=1}^q a_{0j\alpha}^{\gamma \eta} i_{T_{0j}^\beta}\left( - \mathcal{L}_{\psi}\left(\Psi_j^\eta \right) + \sum_{\sigma=1}^q k_j^{\eta \sigma} \wedge \Psi_j^\sigma  \right) + i_{T_{0j}^\alpha\wedge T_{0j}^\beta}\left(\partial \left(- \mathcal{L}_{\psi}  \left( \Psi_j^\gamma  \right) + \sum_{\sigma=1}^q k_j^{\gamma \sigma} \wedge \Psi_j^{\sigma}\right)\right)  \\
&= - \sum_{\eta=1}^q a_{0j\beta}^{\gamma \eta} i_{T_{0j}^\alpha}\mathcal{L}_\psi\left(\Psi_j^\eta \right) - \sum_{\eta,\sigma=1}^q a_{0j\beta}^{\gamma \eta} k_j^{\eta \sigma} i_{T_{0j}^\alpha}\left( \Psi_j^\sigma \right) + \sum_{\eta=1}^q a_{0j\alpha}^{\gamma \eta} i_{T_{0j}^\beta}\mathcal{L}_\psi\left(\Psi_j^\eta \right) + \sum_{\eta,\sigma=1}^q a_{0j\alpha}^{\gamma \eta} k_j^{\eta \sigma} i_{T_{0j}^\beta}\left( \Psi_j^\sigma \right) \notag  \\
&- i_{T_{0j}^\alpha}  i_{ T_{0j}^\beta}\left(\mathcal{L}_\psi \left(\partial \Psi_j^\gamma \right) \right) - \sum_{\sigma=1}^q i_{T_{0j}^\alpha}i_{ T_{0j}^\beta }\left( \partial k_j^{\gamma \sigma} \wedge \Psi_j^\sigma - k_j^{\gamma \sigma} \wedge \partial \Psi_j^\sigma \right) \notag
\end{align}}}
Let us consider the last part:
{\small{\begin{align}\label{its35}
&- i_{T_{0j}^\alpha}  i_{ T_{0j}^\beta}\left(\mathcal{L}_\psi \left(\partial \Psi_j^\gamma \right) \right) - \sum_{\sigma=1}^q i_{T_{0j}^\alpha}i_{T_{0j}^\beta}\left( \partial k_j^{\gamma \sigma} \wedge \Psi_j^\sigma - k_j^{\gamma \sigma} \wedge \partial \Psi_j^\sigma \right)\\
&= - i_{T_{0j}^\alpha} i_{\left[ \psi, T_{0j}^\beta \right]} \left(\partial \Psi_j^\gamma \right) + i_{T_{0j}^\alpha}\mathcal{L}_\psi i_{T_{0j}^\beta} \left(\partial \Psi_j^\gamma \right) + \sum_{\sigma=1}^q i_{T_{0j}^\beta} \left( \partial k_j^{\gamma \sigma}\right) i_{T_{0j}^\alpha}\left(\Psi_j^\sigma \right) - \sum_{\sigma=1}^q i_{T_{0j}^\alpha} \left( \partial k_j^{\gamma \sigma} \right) i_{T_{0j}^\beta}\left( \Psi_j^\sigma \right)  - \sum_{\sigma=1}^q k_j^{\gamma \sigma}  i_{T_{0j}^\alpha \wedge T_{0j}^\beta} \left(\partial \Psi_j^\beta \right) \notag \\
&= - i_{T_{0j}^\alpha} i_{\left[ \psi, T_{0j}^\beta \right]} \left(\partial \Psi_j^\gamma \right)  + i_{ \left[ \psi, T_{0j}^\alpha \right]} i_{T_{0j}^\beta}\left( \partial \Psi_j^\beta \right) + \mathcal{L}_\psi \left(i_{T_{0j}^\alpha \wedge T_{0j}^\beta}\left(\partial \Psi_j^\gamma \right) \right) + \sum_{\sigma=1}^q i_{T_{0j}^\beta} \left( \partial k_j^{\gamma \sigma}\right) i_{T_{0j}^\alpha}\left(\Psi_j^\sigma \right) - \sum_{\sigma=1}^q i_{T_{0j}^\alpha} \left( \partial k_j^{\gamma \sigma} \right) i_{T_{0j}^\beta}\left( \Psi_j^\sigma \right) \notag  \\
&- \sum_{\sigma=1}^q k_j^{\gamma \sigma}\left( \sum_{\delta=1}^q a_{0 j\alpha}^{\sigma \delta} i_{T_{0j}^\beta}\left(\Psi_j^\delta \right) - a_{0 j\beta}^{\sigma \delta} i_{T_{0j}^\alpha}\left( \Psi_j^\delta \right)   \right) \notag
\end{align}}}
Then $(\ref{its38}),(\ref{its36})$ and $(\ref{its35})$ imply that $\bar{\partial} n + \hat{\tilde{E}}_1\left(\kappa \right)=0$. Since $\hat{D}_2\left( \rho\right)=0$ and $\hat{\beta}_2\left( \rho \right)+ n - \hat{\tilde{E}}_1(p)=0$ and so $\alpha_3\circ\hat{D}_2\left(\rho \right) + \hat{E}_2'(n)=0$, we have $\hat{E}_2'(n)=0$. This completes the proof of Lemma \ref{its29}.
\end{proof}

Then from $(\ref{it8})-(\ref{it9})$ and $(\ref{it13}),(\ref{it32})$ and $(\ref{it27})$ and $(\ref{itt12})$ and Lemma \ref{its29}, we obtain
\begin{theorem}
Suppose given a compact foliated complex manifold $\left( M,\Theta_{ \mathcal{F}_0}, \mathcal{N}_{\mathcal{F}_0}^* \right)$ with $\Theta_{\mathcal{F}_0}$ and $\mathcal{N}_{\mathcal{F}_0}^*$ locally free, and $\left( \Phi, \Psi, \left\{\left(\psi, a_j, k_j\right) \right\} \right) \in \mathbb{H}^1\left( M, \mathcal{E}_{\mathcal{F}_0}'^\bullet \right)$. In order that there may exist a foliated analytic family $\left(\mathcal{M}, \Theta_\mathcal{F}, \mathcal{N}_\mathcal{F}^*, B, \pi \right)$ in terms of both tangent and cotangent sheaves such that $\pi^{-1}(0)=\left( M, \Theta_{\mathcal{F}_0}, \mathcal{N}_{\mathcal{F}_0}^* \right)$ and that $\left( \frac{\partial\left(\mathcal{M}, \Theta_{\mathcal{F}}, \mathcal{N}_{\mathcal{F}}^* \right)}{\partial t} \right)|_{t=0}=\left( \Phi, \Psi, \left\{\left( \psi, a_j, k_j \right) \right\} \right)$, it is necessary that
\begin{align*}
&\left( \rho, n, p,  \eta,  \kappa, \left\{- \left[ \left(\psi, a_j, k_j \right), \left( \psi, a_j, k_j \right) \right]     \right\}  \right) \\
&\in \frac{A^{0,0}\left( M, \mathscr{H}om_{\mathcal{O}_M}\left( \bigwedge^2 \Theta_{\mathcal{F}_0} , \Theta_M \right) \right)}{A^{0,0}\left( M, \mathscr{H}om_{\mathcal{O}_M}\left( \bigwedge^2 \Theta_{\mathcal{F}_0}, \Theta_{\mathcal{F}_0} \right) \right) } \bigoplus A^{0,0}\left( M, \mathscr{H}om_{\mathcal{O}_M}\left( \mathcal{N}_{\mathcal{F}_0}^*, \bigwedge^2 \Theta_{\mathcal{F}_0}^* \right) \right) \bigoplus A^{0,0}\left( M, \mathscr{H}om_{\mathcal{O}_M}\left( \mathcal{N}_{\mathcal{F}_0}^*, \Theta_{\mathcal{F}_0}^* \right) \right) \notag\\
&\bigoplus A^{0,1}\left( M, \mathscr{H}om_{\mathcal{O}_M}\left( \Theta_{\mathcal{F}_0}, \Theta_M \right) \right) \bigoplus A^{0,1}\left( M, \mathscr{H}om_{\mathcal{O}_M}\left( \mathcal{N}_{\mathcal{F}_0}^*, \Omega_M^1 \right) \right) \bigoplus A^{0,2}\left( M, \mathcal{E}_{\mathcal{F}_0} \right) \notag
\end{align*}
from $(\ref{itt11})$ in \textnormal{Lemma \ref{its29}} defines the $0$ class in $\mathbb{H}^2\left( M, \mathcal{E}_{\mathcal{F}_0}'^\bullet \right)$.
\end{theorem}

\section{Isomorphisms of infinitesimal foliated deformations and obstructions in terms of both tangent and cotangent sheaves between \v Cech resolution and Dolbeualt resolution of $\mathcal{E}_{\mathcal{F}_0}'^\bullet$ } \label{nii9}

\subsection{Isomorphism of infinitesimal simultaneous foliated deformations between \v Cech resolution and Dolbeault resolution of $\mathcal{E}_{\mathcal{F}_0}'$}\

We keep the notations in subsection \ref{a25}. Then we define the correspondence of $1$-cocycles between \v Cech resolution and Dolbeault resolution of $\mathcal{E}_{\mathcal{F}_0}'^\bullet$ by
{\small{\begin{align}\label{a31}
&\left( \left\{\alpha_j \right\},  \left\{ \beta_j\right\}, \left\{ \left( \theta_{jk}, A_{jk}, B_{jk} \right)      \right\}\right) \in C^0\left( \mathcal{U}, \mathscr{H}om_{\mathcal{O}_M}\left( \Theta_{\mathcal{F}_0}, \Theta_M \right)\right) \bigoplus  C^0\left( \mathcal{U}, \mathscr{H}om_{\mathcal{O}_M}\left( \mathcal{N}_{\mathcal{F}_0}^*, \Omega_M^1 \right) \right) \bigoplus C^1\left( \mathcal{U}, \mathcal{E}_{\mathcal{F}_0}  \right)\\
&\mapsto \left( \Phi, \Psi, \left\{ \left( \psi, a_j:=\bar{\partial} A_j, k_j:=\bar{\partial} B_j \right)  \right\} \right) \in A^{0,0}\left( M, \mathscr{H}om_{\mathcal{O}_M}\left( \Theta_{\mathcal{F}_0}, \Theta_M \right) \right) \bigoplus A^{0,0}\left( M, \mathscr{H}om_{\mathcal{O}_M}\left( \mathcal{N}_{\mathcal{F}_0}^*, \Omega_M^1 \right) \right) \bigoplus A^{0,1}\left( M, \mathcal{E}_{\mathcal{F}_0} \right) \notag
\end{align}}}

Then by combining Proposition \ref{b18} and Proposition \ref{c1}, we obtain

\begin{proposition}
We recall that given a foliated analytic family $\left(\mathcal{M}, \Theta_{\mathcal{F}}, \mathcal{N}_\mathcal{F}^*, B, \pi \right)$ of deformations of $\left( M, \Theta_{\mathcal{F}_0}, \mathcal{N}_{\mathcal{F}_0}^* \right)= \pi^{-1}(0)$ as in \textnormal{Remark \ref{os3}}, the infinitesimal deformation at $t=0$ is given by 
\begin{align} \label{a27}
\left( \left\{ \alpha_j(0)\right\}, \left\{ \beta_j(0)\right\}, \left\{ \left( \theta_{jk}(0), A_{jk}(0), B_{jk}(0) \right) \right\} \right)
\end{align}
 from $(\ref{a30})$ in the first cohomology group of the \v Cech resolution of $\mathcal{E}_{\mathcal{F}_0}^\bullet$. On the other hand, the infinitesimal deformation at $t=0$ is also given by 
 \begin{align}\label{a28}
 \left( \frac{\partial X_j}{\partial t}|_{t=0}, \frac{\partial Y_j}{\partial t}|_{t=0},  \left\{\left( \frac{\partial \varphi}{\partial t}|_{t=0}, \frac{\partial U_j}{\partial t}|_{t=0}, \frac{\partial Z_j}{\partial t}|_{t=0} \right) \right\}\right)
 \end{align}
from $(\ref{a32})$ in the first cohomology group of the Dolbeault resolution of $\mathcal{E}_{\mathcal{F}_0}^\bullet$. Then $(\ref{a27})$ corresponds to $(\ref{a28})$ under the isomorphism $(\ref{a31})$.
\end{proposition}

\subsection{Isomorphism of obstructions between \v Cech resolution and Dolbeault resolution of $\mathcal{E}_{\mathcal{F}_0}'$}\

We construct an isomorphism between the second cohomology group $\mathbb{H}^2\left(M, \mathcal{E}_{\mathcal{F}_0}'^\bullet \right)$ from the \v Cech resolution and the second cohomology group from the Dolbeault resolution of $\mathcal{E}_{\mathcal{F}_0}'^\bullet$. We keep the notations in subsection \ref{a15} and subsection \ref{a16} and subsection \ref{a33}. 

Let $\left(\left\{ C_i \right\}, \left\{ c_i' \right\}, \left\{ J_{j} \right\}, \left\{ B_{jk}\right\}, \left\{ b_{jk}\right\},  \left\{\left( a_{ijk}, G_{ijk}, g_{ijk} \right) \right\} \right)$, where $\left\{ c_i'\right\}\in C^0\left( \mathcal{U}, \mathscr{H}om_{\mathcal{O}_M}\left( \mathcal{N}_{\mathcal{F}_0}^*, \bigwedge^2 \Theta_{\mathcal{F}_0}^* \right) \right)$ such that $c_i'\left(T_{0i}^\alpha \wedge T_{0i}^\beta \right)\left( w_{0i}^\gamma\right)=c_{i\alpha\beta}'^\gamma\in \Gamma\left( U_i, \mathcal{O}_M \right)$. We recall from Part I that
\begin{align*}
\sum_{\eta=1}^q a_{0j\alpha}^{\gamma \eta} h_{0jk}^{\eta \delta}= i_{T_{0j}^\alpha}\left( d h_{0jk}^{\gamma \delta} \right) + \sum_{\sigma=1}^q \sum_{\beta=1}^p h_{0jk}^{\gamma \sigma} r_{0jk}^{\alpha\beta} a_{0k\beta}^{\sigma \delta}
\end{align*}

Then we note that 

{\tiny{\begin{align*}
&c_{j\alpha\beta}'^\gamma  +\sum_{\eta=1}^q a_{0j\beta}^{\gamma \eta} i_{T_{0j}^\alpha}\left(f_j^\eta \right) - \sum_{\eta=1}^q a_{0j\alpha}^{\gamma \eta} i_{T_{0j}^\beta}\left( f_j^\eta \right) + i_{T_{0j}^\alpha \wedge T_{0j}^\beta} \left(\partial f_j^\gamma \right)\\
&- \sum_{\delta, \eta=1}^p \sum_{\sigma=1}^q r_{0jk}^{\alpha \delta} r_{0jk}^{\beta \eta} h_{0jk}^{\gamma \sigma}\left( c_{k\delta \eta}'^\sigma + \sum_{\rho=1}^q a_{0k\eta}^{\sigma \rho} i_{T_{0k}^\delta}\left( f_k^\rho \right) - \sum_{\rho=1}^q a_{0k\delta}^{\sigma \rho} i_{T_{0k}^\eta} \left( f_k^\rho\right) + i_{T_{0k}^\delta\wedge T_{0k}^{\eta}} \left( {\partial} f_k^\sigma \right)\right)\\
&= c_{j\alpha\beta}'^\gamma - \sum_{\delta, \eta=1}^p\sum_{\sigma=1}^q r_{0jk}^{\alpha \delta} r_{0jk}^{\beta \eta} h_{0jk}^{\gamma \sigma} c_{k\delta\eta}'^\sigma  +\sum_{\eta=1}^q a_{0j\beta}^{\gamma \eta} i_{T_{0j}^\alpha}\left( \sum_{\delta=1}^q h_{0jk}^{\eta \delta} f_k^\delta - \tilde{\mathcal{L}}_{d_{jk}}\left( w_{0j}^\eta\right) + e_{jk}^\eta + b_{jk}^\eta   \right) - \sum_{\eta=1}^q a_{0j\alpha}^{\gamma \eta} i_{T_{0j}^\beta}\left(\sum_{\delta=1}^q h_{0jk}^{\eta \delta} f_k^\delta - \tilde{\mathcal{L}}_{d_{jk}}\left( w_{0j}^\eta\right) + e_{jk}^\eta + b_{jk}^\eta \right) \\
&+ i_{T_{0j}^\alpha \wedge T_{0j}^\beta} \left(\partial \left( \sum_{\delta=1}^q h_{0jk}^{\gamma \delta} f_k^\delta - \tilde{\mathcal{L}}_{d_{jk}}\left( w_{0j}^\gamma\right) + e_{jk}^\gamma + b_{jk}^\gamma  \right) \right)- \sum_{\delta, \eta=1}^p \sum_{\sigma=1}^q r_{0jk}^{\alpha \delta} r_{0jk}^{\beta \eta} h_{0jk}^{\gamma \sigma}\left(  \sum_{\rho=1}^q a_{0k\eta}^{\sigma \rho} i_{T_{0k}^\delta}\left( f_k^\rho \right) - \sum_{\rho=1}^q a_{0k\delta}^{\sigma \rho} i_{T_{0k}^\eta} \left( f_k^\rho\right) + i_{T_{0k}^\delta\wedge T_{0k}^{\eta}} \left( {\partial} f_k^\sigma \right)\right)\\
&=\sum_{\eta,\delta=1}^q a_{0j\beta}^{\gamma \eta} h_{0jk}^{\eta \delta} i_{T_{0j}^\alpha}\left( f_k^\delta \right) - \sum_{\eta,\delta=1}^q a_{0j\alpha}^{\gamma \eta} h_{0jk}^{\eta \delta} i_{T_{0j}^\beta}\left( f_k^\delta \right) + \sum_{\delta=1}^q i_{T_{0j}^\alpha \wedge T_{0j}^\beta}\left( d h_{0jk}^{\gamma \delta} \wedge f_k^\delta + h_{0jk}^{\gamma \delta} \partial f_k^\delta \right) \\
&- \sum_{\delta, \eta=1}^p \sum_{\sigma, \rho=1}^q r_{0jk}^{\alpha \delta} r_{0jk}^{\beta \eta} h_{0jk}^{\gamma \sigma} a_{0k\eta}^{\sigma \rho} i_{T_{0k}^\delta}\left( f_k^\rho \right) + \sum_{\delta, \eta=1}^p \sum_{\sigma, \rho=1}^q r_{0jk}^{\alpha \delta} r_{0jk}^{\beta \eta} h_{0jk}^{\gamma \sigma} a_{0k\delta}^{\sigma \rho} i_{T_{0k}^\eta}\left( f_k^\rho \right) - \sum_{\delta, \eta=1}^p \sum_{\sigma=1}^q r_{0jk}^{\alpha \delta} r_{0jk}^{\beta \eta} h_{0jk}^{\gamma \sigma} i_{T_{0k}^\delta \wedge T_{0k}^\eta}\left( \partial f_k^\sigma \right)\\
&=\sum_{\delta=1}^q \left( i_{T_{0j}^\beta}\left( d h_{0jk}^{\gamma \delta} \right) + \sum_{\sigma=1}^q \sum_{\eta=1}^p h_{0jk}^{\gamma \sigma} r_{0jk}^{\beta \eta} a_{0k\eta}^{\sigma \delta}\right) i_{T_{0j}^\alpha}\left( f_k^\delta \right) - \sum_{\eta,\delta=1}^q \left( i_{T_{0j}^\alpha}\left( d h_{0jk}^{\gamma \delta} \right) + \sum_{\sigma=1}^q \sum_{\beta=1}^p h_{0jk}^{\gamma \sigma} r_{0jk}^{\alpha\beta} a_{0k\beta}^{\sigma \delta} \right) i_{T_{0j}^\beta}\left( f_k^\delta \right)  \\
&+ \sum_{\delta=1}^q i_{T_{0j}^\alpha}\left( d h_{0jk}^{\gamma \delta} \right) i_{T_{0j}^\beta}\left( f_k^\delta\right) - \sum_{\delta=1}^q i_{T_{0j}^\beta}\left( d h_{0jk}^{\gamma \delta} \right) i_{T_{0j}^\alpha}\left( f_k^\delta \right)  - \sum_{\delta, \eta=1}^p \sum_{\sigma, \rho=1}^q r_{0jk}^{\alpha \delta} r_{0jk}^{\beta \eta} h_{0jk}^{\gamma \sigma} a_{0k\eta}^{\sigma \rho} i_{T_{0k}^\delta}\left( f_k^\rho \right) + \sum_{\delta, \eta=1}^p \sum_{\sigma, \rho=1}^q r_{0jk}^{\alpha \delta} r_{0jk}^{\beta \eta} h_{0jk}^{\gamma \sigma} a_{0k\delta}^{\sigma \rho} i_{T_{0k}^\eta}\left( f_k^\rho \right) = 0
\end{align*}}}

Let us define $R_j'\in \Gamma\left( U_j, \mathcal{A}^{0,0}\left( \mathscr{H}om_{\mathcal{O}_M}\left( \mathcal{N}_{\mathcal{F}_0}^*, \bigwedge^2 \Theta_{\mathcal{F}_0}^* \right) \right) \right)$ by
\begin{align*}
R_j':\Gamma\left( U_j, \mathcal{N}_{\mathcal{F}_0}^* \right)&\to \Gamma\left( U_j , \bigwedge^2 \Theta_{\mathcal{F}_0}^* \right)\\
 w_{0j}^\gamma &\mapsto \left( T_{0j}^\alpha \wedge T_{0j}^\beta \mapsto - c_{j\alpha\beta}'^\gamma  -\sum_{\eta=1}^q a_{0j\beta}^{\gamma \eta} i_{T_{0j}^\alpha}\left(f_j^\eta \right) + \sum_{\eta=1}^q a_{0j\alpha}^{\gamma \eta} i_{T_{0j}^\beta}\left( f_j^\eta \right) - i_{T_{0j}^\alpha \wedge T_{0j}^\beta} \left(\partial f_j^\gamma \right)    \right)
\end{align*}
Then $R':=\left\{ R_j'\right\}$ defines a global section in $A^{0,0}\left( M, \mathscr{H}om_{\mathcal{O}_M}\left( \mathcal{N}_{\mathcal{F}_0}^*, \bigwedge^2 \Theta_{\mathcal{F}_0}^* \right) \right)$. Then we define the correspondence of $2$-cocycles between \v Cech resolution and Dolbeault resolution of $\mathcal{E}_{\mathcal{F}_0}'^\bullet$ by
\begin{align}\label{a38}
\left(\left\{ C_i \right\}, \left\{ c_i' \right\}, \left\{ J_{j} \right\}, \left\{ B_{jk}\right\}, \left\{ b_{jk}\right\},  \left\{\left( a_{ijk}, G_{ijk}, g_{ijk} \right) \right\} \right) \mapsto \left( \overline{P}, R', I,  L, Q,  \left\{ \left( \bar{\partial} p_j, \bar{\partial} U_j, \bar{\partial} u_j \right) \right\} \right)
\end{align}

\begin{proposition}\label{a35}
We recall from \textnormal{Lemma \ref{itt1}} that given a $1$-cocycle $\left( \left\{ \alpha_j \right\}, \left\{ \beta_j \right\} , \left\{ \left( \theta_{jk}, A_{jk}, B_{jk} \right) \right\} \right)$ in $\mathbb{H}^1\left(  M, \mathcal{E}_{\mathcal{F}_0}'^\bullet \right)$ from the \v Cech resolution of $\mathcal{E}_{\mathcal{F}_0}'^\bullet$, we have the associated obstruction class from $(\ref{itt3})$
\begin{align}\label{a36}
\left(\left\{ R_j \right\}, \left\{ S_j' \right\}, \left\{ P_j \right\} ,  \left\{ F_{jk} \right\}, \left\{ G_{jk} \right\}, \left\{ \left[ \left(\theta_{ij}, A_{ij}, B_{ij} \right), \left( \theta_{jk}, A_{jk}, B_{jk} \right) \right] \right\} \right) \in \mathbb{H}^2\left( M, \mathcal{E}_{\mathcal{F}_0}'^\bullet \right)
\end{align}
from the \v Cech resolution of $\mathcal{E}_{\mathcal{F}_0}'^\bullet$. Assume that the $1$-cocycle $\left( \left\{ \alpha_j\right\}, \left\{ \beta_j \right\}, \left\{ \left( \theta_{jk}, A_{jk}, B_{jk} \right) \right\}\right)$ from the \v Cech resolution of $\mathcal{E}_{\mathcal{F}_0}'^\bullet$ corresponds to the $1$-cocycle $\left( \Phi, \Psi, \left\{\left(\psi, a_j:=\bar{\partial} A_j, k_j:=\bar{\partial} B_j \right) \right\} \right)$ from the Dolbeault resolution of $\mathcal{E}_{\mathcal{F}_0}'^\bullet$ under the isomorphism $(\ref{a31})$. We recall from \textnormal{Lemma \ref{its29}} that we have the associated obstruction class from $(\ref{itt11})$
\begin{align}\label{a37}
\left( \rho, n, p,  \eta,  \kappa, \left\{- \left[ \left(\psi, a_j, k_j \right), \left( \psi, a_j, k_j \right) \right]     \right\}  \right) \in \mathbb{H}^2 \left( M, \mathcal{E}_{\mathcal{F}_0}'^\bullet \right)
\end{align}
from the Dolbeault resolution of $\mathcal{E}_{\mathcal{F}_0}'^\bullet$. Then $(\ref{a36})$ corresponds to $(\ref{a37})$ under the isomorphism $(\ref{a38})$.
\end{proposition}

\begin{proof}
Let $\left( \left\{ \alpha_j \right\}, \left\{  \beta_j \right\} , \left\{ \left( \theta_{jk}, A_{jk}, B_{jk} \right) \right\} \right)$ be a $1$-cocycle in Lemma \ref{itt1}. Then $\left(\left\{ R_j \right\}, \left\{ S_j \right\}, \left\{ P_j \right\} ,  \left\{ F_{jk} \right\}, \left\{ G_{jk} \right\},\right.$ $\left. \left\{ \left[ \left(\theta_{ij}, A_{ij}, B_{ij} \right), \left( \theta_{jk}, A_{jk}, B_{jk} \right) \right] \right\} \right)$ defines a $2$-cocycle by Lemma \ref{itt1}. We have the corresponding $1$-cocycle $\left( \Phi , \Psi, \left\{\psi, a_j=\bar{\partial} A_j, k_j=\bar{\partial} B_j \right\} \right)$ from the Dolbeault resolution of $\mathcal{E}_{\mathcal{F}_0}^\bullet$ from $(\ref{a31})$. We keep the notations in the proof of Proposition \ref{a13} and the proof of Proposition \ref{a24} and the proof of Proposition \ref{a22}. Then it remains to show that $\left\{ S_j' \right\}$ corresponds to $n$. In fact, we note that
 {\small{\begin{align}
&S_j'\left(T_{0j}^\alpha \wedge T_{0j}^\beta \right)\left( w_{0j}^\gamma \right) +\sum_{\eta=1}^q a_{0j\beta}^{\gamma \eta} i_{T_{0j}^\alpha}\left(f_j^\gamma \right) - \sum_{\eta=1}^q a_{0j\alpha}^{\gamma \eta} i_{T_{0j}^\beta}\left( f_j^\gamma \right) + i_{T_{0j}^\alpha \wedge T_{0j}^\beta} \left(\partial f_j^\gamma \right)\label{a39}\\
&=2\sum_{\eta=1}^q a_{0j\beta}^{\gamma \eta} i_{\alpha_j\left( T_{0j}^\alpha \right)}\left( \beta_j\left(w_{0j}^\gamma \right)\right) -2 \sum_{\eta=1}^q a_{0j\alpha}^{\gamma \eta} i_{\alpha_j\left( T_{0j}^\beta \right)}\left(\beta_j\left( w_{0j}^\gamma \right) \right) + 2 i_{\alpha_j\left( T_{0j}^\alpha \right) \wedge \alpha_j\left( T_{0j}^\beta \right)}\left( d w_{0j}^\gamma \right) \notag \\
&+ 2 i_{\alpha_j\left(T_{0j}^\alpha \right)\wedge T_{0j}^\beta}\left( d\beta_j\left( w_{0j}^\gamma \right) \right) + 2 i_{T_{0j}^\alpha \wedge \alpha_j\left( T_{0j}^\beta \right)} \left( d\beta_j\left( w_{0j}^\gamma \right) \right) + \sum_{\eta=1}^q a_{0j\beta}^{\gamma \eta} i_{T_{0j}^\alpha}\left(f_j^\gamma \right) - \sum_{\eta=1}^q a_{0j\alpha}^{\gamma \eta} i_{T_{0j}^\beta}\left( f_j^\gamma \right) + i_{T_{0j}^\alpha \wedge T_{0j}^\beta} \left(\partial f_j^\gamma \right)\notag\\
&= 2 \sum_{\eta=1}^q a_{0j\beta}^{\gamma \eta} i_{\Phi_j^\alpha}\left( \Psi_j^\eta \right) - 2 \sum_{\eta=1}^q a_{0j\alpha}^{\gamma \eta} i_{\Phi_j^\beta}\left(\Psi_j^\eta \right) + \sum_{\eta=1}^q a_{0j\beta}^{\gamma \eta} i_{2 \alpha_j\left(A_j^\alpha \right) -2 \left[\xi_j, \alpha_j\left(T_{0j}^\alpha \right) \right] - 2\left[ \xi_j, A_j^\alpha \right] +\left[ \xi_j, \left[ \xi_j, T_{0j}^\alpha \right] \right] }\left(w_{0j}^\eta \right) \notag\\
&- \sum_{\eta=1}^q a_{0j\alpha}^{\gamma \eta} i_{ \alpha_j \left(A_j^\beta \right) -2 \left[\xi_j, \alpha_j\left(T_{0j}^\beta \right) \right] - 2\left[ \xi_j, A_j^\beta \right] +\left[ \xi_j, \left[ \xi_j, T_{0j}^\beta \right] \right] }\left(w_{0j}^\eta \right) + 2 i_{\alpha_j\left( T_{0j}^\alpha \right) \wedge \alpha_j\left( T_{0j}^\beta \right)}\left( d w_{0j}^\gamma \right) + 2 i_{\alpha_j\left(T_{0j}^\alpha \right)\wedge T_{0j}^\beta}\left( d\beta_j\left( w_{0j}^\gamma \right) \right) \notag\\
& + 2 i_{T_{0j}^\alpha \wedge \alpha_j\left( T_{0j}^\beta \right)} \left( d\beta_j\left( w_{0j}^\gamma \right) \right)  + i_{T_{0j}^\alpha \wedge T_{0j}^\beta }\left( \partial\left( -2 \beta_j\left( B_j^\gamma \right) + 2\mathcal{L}_{\xi_j}\left( \beta_j\left(w_{0j}^\gamma\right) \right) +2 \mathcal{L}_{\xi_j}\left( B_j^\gamma \right) -\mathcal{L}_{\xi_j}\mathcal{L}_{\xi_j}\left( w_{0j}^\gamma \right)  \right) \right) \notag\\
&= 2 \sum_{\eta=1}^q a_{0j\beta}^{\gamma \eta} i_{\Phi_j^\alpha}\left( \Psi_j^\eta \right) - 2 \sum_{\eta=1}^q a_{0j\alpha}^{\gamma \eta} i_{\Phi_j^\beta}\left(\Psi_j^\eta \right) + 2 i_{\alpha_j\left( A_j^\alpha \right)} i_{T_{0j}^\beta}\left(  dw_{0j}^\gamma \right) - 2 i_{\left[ \xi_j, \alpha_j\left( T_{0j}^\alpha \right) \right]} i_{T_{0j}^\beta}\left( dw_{0j}^\gamma \right) - 2 i_{\left[ \xi_j, A_j^\alpha \right]} i_{T_{0j}^\beta}\left( dw_{0j}^\gamma \right) \notag\\
& + i_{\left[ \xi_j, \left[ \xi_j, T_{0j}^\alpha\right] \right]} i_{T_{0j}^\beta}\left( dw_{0j}^\gamma  \right) - 2 i_{\alpha_j\left( A_j^\beta \right)} i_{T_{0j}^\alpha}\left(  dw_{0j}^\gamma \right) + 2 i_{\left[ \xi_j, \alpha_j\left( T_{0j}^\beta \right) \right]} i_{T_{0j}^\alpha}\left( dw_{0j}^\gamma \right) +2 i_{\left[ \xi_j, A_j^\beta \right]} i_{T_{0j}^\alpha}\left( dw_{0j}^\gamma \right)  - i_{\left[ \xi_j, \left[ \xi_j, T_{0j}^\beta\right] \right]} i_{T_{0j}^\alpha}\left( dw_{0j}^\gamma  \right)\notag\\
&+2 i_{\alpha_j\left( T_{0j}^\alpha \right) \wedge \alpha_j\left( T_{0j}^\beta \right)}\left( d w_{0j}^\gamma \right) + 2 i_{\alpha_j\left(T_{0j}^\alpha \right)\wedge T_{0j}^\beta}\left( d\beta_j\left( w_{0j}^\gamma \right) \right) + 2 i_{T_{0j}^\alpha \wedge \alpha_j\left( T_{0j}^\beta \right)}\left( d\beta_j\left( w_{0j}^\gamma \right) \right) \notag \\
&- 2 i_{T_{0j}^\alpha \wedge T_{0j}^\beta} \left(\partial \beta_j\left( B_j^\gamma \right) \right) +  2 i_{T_{0j}^\alpha \wedge T_{0j}^\beta}\left( \mathcal{L}_{\xi_j}\left( d\beta_j\left( w_{0j}^\gamma \right) \right) \right) + 2 i_{T_{0j}^\alpha \wedge T_{0j}^\beta } \left( \mathcal{L}_{\xi_j} \left( \partial B_j^\gamma \right) \right) - i_{T_{0j}^\alpha \wedge T_{0j}^\beta}\left( \mathcal{L}_{\xi_j} \mathcal{L}_{\xi_j} \left( d w_{0j}^\gamma \right) \right) \notag
\end{align}}}

On the other hand, we note that
{\small{\begin{align}
&2\sum_{\eta=1}^q a_{0j\beta}^{\gamma \eta} i_{\Phi_j^\alpha }\left( \Psi_j^\eta\right) -2 \sum_{\eta=1}^q a_{0j\alpha}^{\gamma \eta} i_{\Phi_j^\beta}\left(\Psi_j^\eta \right) + 2 i_{\Phi_j^\alpha \wedge \Phi_j^\beta }\left( d w_{0j}^\gamma \right)+ 2 i_{\Phi_j^\alpha \wedge T_{0j}^\beta}\left( \partial  \Psi_j^\gamma \right) + 2 i_{T_{0j}^\alpha \wedge \Phi_j^\beta } \left( \partial \Psi_j^\gamma \right)\label{a40}\\
&= 2 i_{\left(\left[ \xi_j, T_{0j}^\alpha \right]- A_j^\alpha - \alpha_j\left(T_{0j}^\alpha \right)\right)\wedge\left(\left[ \xi_j, T_{0j}^\beta\right]- A_j^\beta - \alpha\left(T_{0j}^\beta \right) \right)}\left( dw_{0j}^\gamma \right) + 2i_{\left(\left[\xi_j, T_{0j}^\alpha \right]- A_j^\alpha - \alpha_j\left(T_{0j}^\alpha \right)  \right) \wedge T_{0j}^\beta}\left(\partial\left( \mathcal{L}_{\xi_j}\left(w_{0j}^\gamma \right) - B_j^\gamma - \beta_j\left(w_{0j}^\gamma \right) \right) \right)\notag\\
&+ 2 i_{T_{0j}^\alpha \wedge \left( \left[ \xi_j, T_{0j}^\beta \right]- B_j^\beta - \beta_j\left( T_{0j}^\beta \right) \right)}\left(\partial \left( \mathcal{L}_{\xi_j}\left(w_{0j}^\gamma \right) - B_j^\gamma - \beta_j\left( w_{0j}^\gamma \right) \right) \right) + 2\sum_{\eta=1}^q a_{0j\beta}^{\gamma \eta} i_{\Phi_j^\alpha }\left( \Psi_j^\eta\right) -2 \sum_{\eta=1}^q a_{0j\alpha}^{\gamma \eta} i_{\Phi_j^\beta}\left(\Psi_j^\eta \right) \notag\\
&=  -2 i_{\left[ \xi_j, T_{0j}^\alpha\right]} i_{\left[ \xi_j , T_{0j}^\beta \right]}\left(dw_{0j}^\gamma \right) - 2 i_{A_j^\beta} i_{\left[ \xi_j, T_{0j}^\alpha \right]}\left( dw_{0j}^\gamma \right) -  2 i_{\alpha_j\left( T_{0j}^\beta \right)} i_{\left[ \xi_j,  T_{0j}^\alpha \right] } \left( dw_{0j}^\gamma \right)^* + 2 i_{A_j^\alpha} i_{\left[ \xi_j, T_{0j}^\beta \right]}\left( dw_{0j}^\gamma \right) -2 i_{A_j^\alpha} i_{\alpha_j\left( T_{0j}^\beta \right)}\left( dw_{0j}^\gamma \right) \notag \\
&+ 2i_{\alpha_j\left( T_{0j}^\alpha \right)} i_{\left[ \xi_j, T_{0j}^\beta \right]} \left( d w_{0j}^\gamma \right)^* + 2 i_{\alpha_j\left( T_{0j}^\alpha \right) \wedge \alpha_j\left( T_{0j}^\beta \right)} \left( dw_{0j}^\gamma \right) - 2 i_{\left[ \xi_j, T_{0j}^\alpha \right]} i_{T_{0j}^\beta} \mathcal{L}_{\xi_j}\left( dw_{0j}^\gamma \right) + 2 i_{A_j^\alpha} i_{T_{0j}^\beta} \mathcal{L}_{\xi_j}\left(d w_{0j}^\gamma \right) +2 i_{\alpha_j\left( T_{0j}^\alpha \right)} i_{T_{0j}^\beta} \mathcal{L}_{\xi_j}\left( dw_{0j}^\gamma \right)^* \notag\\
&- 2i_{T_{0j}^\beta} i_{\left[ \xi_j, T_{0j}^\alpha \right]}\left( \partial B_j^\gamma \right) + 2i_{T_{0j}^\beta} i_{\alpha_j\left( T_{0j}^\alpha \right)}\left(\partial B_j^\gamma \right) - 2 i_{T_{0j}^\beta} i_{\left[  \xi_j, T_{0j}^\alpha \right]}\left( d\beta_j\left(w_{0j}^\gamma \right) \right)+ 2i_{T_{0j}^\beta} i_{A_j^\alpha}\left( d\beta_j\left( w_{0j}^\gamma \right) \right) + 2i_{\alpha_j\left( T_{0j}^\alpha \right) \wedge T_{0j}^\beta }\left( d\beta_j \left( w_{0j}^\gamma \right) \right) \notag\\
&+ 2 i_{\left[ \xi_j, T_{0j}^\beta \right]} i_{T_{0j}^\alpha} \mathcal{L}_{\xi_j}\left( dw_{0j}^\gamma \right) - 2 i_{A_j^\beta} i_{T_{0j}^\alpha} \mathcal{L}_{\xi_j}\left(d w_{0j}^\gamma \right)  - 2 i_{\alpha_j\left( T_{0j}^\beta \right)} i_{T_{0j}^\alpha} \mathcal{L}_{\xi_j}\left( dw_{0j}^\gamma \right)^* + 2\sum_{\eta=1}^q a_{0j\beta}^{\gamma \eta} i_{\Phi_j^\alpha }\left( \Psi_j^\eta\right) -2 \sum_{\eta=1}^q a_{0j\alpha}^{\gamma \eta} i_{\Phi_j^\beta}\left(\Psi_j^\eta \right) \notag\\
&+ 2i_{T_{0j}^\alpha} i_{\left[ \xi_j, T_{0j}^\beta \right]}\left( \partial B_j^\gamma \right) - 2i_{T_{0j}^\alpha} i_{\alpha_j\left( T_{0j}^\beta \right)}\left(\partial B_j^\gamma \right) + 2 i_{T_{0j}^\alpha} i_{\left[  \xi_j, T_{0j}^\beta \right]}\left( d\beta_j\left(w_{0j}^\gamma \right) \right)- 2i_{T_{0j}^\alpha} i_{A_j^\beta}\left( d\beta_j\left( w_{0j}^\gamma \right) \right) - 2i_{\alpha_j\left( T_{0j}^\beta \right) \wedge T_{0j}^\alpha }\left( d\beta_j \left( w_{0j}^\gamma \right) \right) \notag
\end{align}}}

In view of $(\ref{a39})$ and $(\ref{a40})$, we have to show that
\begin{align}
& 2 i_{\alpha_j\left( A_j^\alpha \right)} i_{T_{0j}^\beta}\left(  dw_{0j}^\gamma \right) - 2 i_{\left[ \xi_j, \alpha_j\left( T_{0j}^\alpha \right) \right]} i_{T_{0j}^\beta}\left( dw_{0j}^\gamma \right) - 2 i_{\left[ \xi_j, A_j^\alpha \right]} i_{T_{0j}^\beta}\left( dw_{0j}^\gamma \right) \label{a41}\\
& + i_{\left[ \xi_j, \left[ \xi_j, T_{0j}^\alpha\right] \right]} i_{T_{0j}^\beta}\left( dw_{0j}^\gamma  \right) - 2 i_{\alpha_j\left( A_j^\beta \right)} i_{T_{0j}^\alpha}\left(  dw_{0j}^\gamma \right) + 2 i_{\left[ \xi_j, \alpha_j\left( T_{0j}^\beta \right) \right]} i_{T_{0j}^\alpha}\left( dw_{0j}^\gamma \right) +2 i_{\left[ \xi_j, A_j^\beta \right]} i_{T_{0j}^\alpha}\left( dw_{0j}^\gamma \right)  - i_{\left[ \xi_j, \left[ \xi_j, T_{0j}^\beta\right] \right]} i_{T_{0j}^\alpha}\left( dw_{0j}^\gamma  \right)\notag\\
&- 2 i_{T_{0j}^\alpha \wedge T_{0j}^\beta} \left(\partial \beta_j\left( B_j^\gamma \right) \right) +  2 i_{T_{0j}^\alpha \wedge T_{0j}^\beta}\left( \mathcal{L}_{\xi_j}\left( d\beta_j\left( w_{0j}^\gamma \right) \right) \right) + 2 i_{T_{0j}^\alpha \wedge T_{0j}^\beta } \left( \mathcal{L}_{\xi_j} \left( \partial B_j^\gamma \right) \right) - i_{T_{0j}^\alpha \wedge T_{0j}^\beta}\left( \mathcal{L}_{\xi_j} \mathcal{L}_{\xi_j} \left( d w_{0j}^\gamma \right) \right) \notag\\
&=  -2 i_{\left[ \xi_j, T_{0j}^\alpha\right]} i_{\left[ \xi_j , T_{0j}^\beta \right]}\left(dw_{0j}^\gamma \right) - 2 i_{A_j^\beta} i_{\left[ \xi_j, T_{0j}^\alpha \right]}\left( dw_{0j}^\gamma \right) -  2 i_{\alpha_j\left( T_{0j}^\beta \right)} i_{\left[ \xi_j,  T_{0j}^\alpha \right] } \left( dw_{0j}^\gamma \right)^* + 2 i_{A_j^\alpha} i_{\left[ \xi_j, T_{0j}^\beta \right]}\left( dw_{0j}^\gamma \right) -2 i_{A_j^\alpha} i_{\alpha_j\left( T_{0j}^\beta \right)}\left( dw_{0j}^\gamma \right) \notag\\
&+ 2i_{\alpha_j\left( T_{0j}^\alpha \right)} i_{\left[ \xi_j, T_{0j}^\beta \right]} \left( d w_{0j}^\gamma \right)^*  - 2 i_{\left[ \xi_j, T_{0j}^\alpha \right]} i_{T_{0j}^\beta} \mathcal{L}_{\xi_j}\left( dw_{0j}^\gamma \right) + 2 i_{A_j^\alpha} i_{T_{0j}^\beta} \mathcal{L}_{\xi_j}\left(d w_{0j}^\gamma \right) +2 i_{\alpha_j\left( T_{0j}^\alpha \right)} i_{T_{0j}^\beta} \mathcal{L}_{\xi_j}\left( dw_{0j}^\gamma \right)^* \notag \\
&- 2i_{T_{0j}^\beta} i_{\left[ \xi_j, T_{0j}^\alpha \right]}\left( \partial B_j^\gamma \right) + 2i_{T_{0j}^\beta} i_{\alpha_j\left( T_{0j}^\alpha \right)}\left(\partial B_j^\gamma \right) - 2 i_{T_{0j}^\beta} i_{\left[  \xi_j, T_{0j}^\alpha \right]}\left( d\beta_j\left(w_{0j}^\gamma \right) \right)^* + 2i_{T_{0j}^\beta} i_{A_j^\alpha}\left( d\beta_j\left( w_{0j}^\gamma \right) \right) \notag \\
&+ 2 i_{\left[ \xi_j, T_{0j}^\beta \right]} i_{T_{0j}^\alpha} \mathcal{L}_{\xi_j}\left( dw_{0j}^\gamma \right) - 2 i_{A_j^\beta} i_{T_{0j}^\alpha} \mathcal{L}_{\xi_j}\left(d w_{0j}^\gamma \right)  - 2 i_{\alpha_j\left( T_{0j}^\beta \right)} i_{T_{0j}^\alpha} \mathcal{L}_{\xi_j}\left( dw_{0j}^\gamma \right)^*  \notag\\
&+ 2i_{T_{0j}^\alpha} i_{\left[ \xi_j, T_{0j}^\beta \right]}\left( \partial B_j^\gamma \right) - 2i_{T_{0j}^\alpha} i_{\alpha_j\left( T_{0j}^\beta \right)}\left(\partial B_j^\gamma \right) + 2 i_{T_{0j}^\alpha} i_{\left[  \xi_j, T_{0j}^\beta \right]}\left( d\beta_j\left(w_{0j}^\gamma \right) \right)^* - 2i_{T_{0j}^\alpha} i_{A_j^\beta}\left( d\beta_j\left( w_{0j}^\gamma \right) \right) \notag
\end{align}

First we compute the parts $*$ in $(\ref{a41})$.
{\tiny{\begin{align}
&-2 i_{\alpha_j\left( T_{0j}^\beta \right)} i_{\left[ \xi_j, T_{0j}^\alpha \right]}\left( d w_{0j}^\gamma \right) - 2i_{\alpha_j\left( T_{0j}^\beta\right)} i_{T_{0j}^\alpha} \mathcal{L}_{\xi_j} \left( dw_{0j}^\gamma \right) + 2i_{\alpha_j\left( T_{0j}^\alpha \right)} i_{\left[ \xi_j, T_{0j}^\beta \right]}\left( dw_{0j}^\gamma \right) + 2 i_{\alpha_j\left( T_{0j}^\alpha \right)} i_{T_{0j}^\beta} \mathcal{L}_{\xi_j}\left( dw_{0j}^\gamma \right)\label{a42}\\
& -2 i_{T_{0j}^\beta} i_{\left[ \xi_j, T_{0j}^\alpha \right]} \left( d\beta_j\left( w_{0j}^\gamma \right) \right)+ 2i_{T_{0j}^\alpha} i_{\left[ \xi_j, T_{0j}^\beta \right]} \left( d \beta_j \left( w_{0j}^\gamma \right) \right)\notag\\
&= - 2i_{\alpha_j\left( T_{0j}^\beta\right)} \mathcal{L}_{\xi_j} i_{T_{0j}^\alpha}\left( dw_{0j}^\gamma \right) + 2 i_{\alpha_j\left( T_{0j}^\alpha \right)} \mathcal{L}_{\xi_j} i_{T_{0j}^\beta}\left( dw_{0j}^\gamma \right) -2 i_{T_{0j}^\beta} \mathcal{L}_{\xi_j} i_{T_{0j}^\alpha} \left( d\beta_j\left( w_{0j}^\gamma \right) \right) +2 i_{T_{0j}^\beta} i_{T_{0j}^\alpha}\left( \mathcal{L}_{\xi_j}\left( d \beta_j\left( w_{0j}^\gamma \right) \right) \right) + 2i_{T_{0j}^\alpha} i_{\left[ \xi_j, T_{0j}^\beta \right]} \left( d \beta_j \left( w_{0j}^\gamma \right) \right)\notag\\
&=2 i_{\left[ \xi_j, \alpha_j\left( T_{0j}^\beta \right)\right]} i_{T_{0j}^\alpha}\left( dw_{0j}^\gamma \right) - 2 \mathcal{L}_{\xi_j} i_{\alpha_j\left( T_{0j}^\beta\right)} i_{T_{0j}^\alpha} \left( dw_{0j}^\gamma \right) - 2i_{\left[ \xi_j,  \alpha_j \left( T_{0j}^\alpha \right) \right]} i_{T_{0j}^\beta} \left( dw_{0j}^\gamma \right) +2 \mathcal{L}_{\xi_j} i_{\alpha_j\left( T_{0j}^\alpha \right)} i_{T_{0j}^\beta} \left( dw_{0j}^\gamma \right) - 2 \mathcal{L}_{\xi_j} i_{T_{0j}^\beta} i_{T_{0j}^\alpha } \left( d\beta_j\left(w_{0j}^\gamma \right) \right) \notag\\
&+ 2 i_{T_{0j}^\alpha \wedge T_{0j}^\beta} \left( \mathcal{L}_{\xi_j}\left( d\beta_j\left( w_{0j}^\gamma \right) \right) \right) = 2 i_{\left[ \xi_j, \alpha_j\left( T_{0j}^\beta \right)\right]} i_{T_{0j}^\alpha}\left( dw_{0j}^\gamma \right) - 2i_{\left[ \xi_j,  \alpha_j \left( T_{0j}^\alpha \right) \right]} i_{T_{0j}^\beta} \left( dw_{0j}^\gamma \right) + 2 i_{T_{0j}^\alpha \wedge T_{0j}^\beta} \left( \mathcal{L}_{\xi_j}\left( d\beta_j\left( w_{0j}^\gamma \right) \right) \right) \notag
\end{align}}}

In view of $(\ref{a41})$ and $(\ref{a42})$, we have to show that
\begin{align}
& 2 i_{\alpha_j\left( A_j^\alpha \right)} i_{T_{0j}^\beta}\left(  dw_{0j}^\gamma \right)  - 2 i_{\left[ \xi_j, A_j^\alpha \right]} i_{T_{0j}^\beta}\left( dw_{0j}^\gamma \right) \label{a45} \\
& + i_{\left[ \xi_j, \left[ \xi_j, T_{0j}^\alpha\right] \right]} i_{T_{0j}^\beta}\left( dw_{0j}^\gamma  \right) - 2 i_{\alpha_j\left( A_j^\beta \right)} i_{T_{0j}^\alpha}\left(  dw_{0j}^\gamma \right)  +2 i_{\left[ \xi_j, A_j^\beta \right]} i_{T_{0j}^\alpha}\left( dw_{0j}^\gamma \right)  - i_{\left[ \xi_j, \left[ \xi_j, T_{0j}^\beta\right] \right]} i_{T_{0j}^\alpha}\left( dw_{0j}^\gamma  \right) \notag \\
&- 2 i_{T_{0j}^\alpha \wedge T_{0j}^\beta} \left(\partial \beta_j\left( B_j^\gamma \right) \right) + 2 i_{T_{0j}^\alpha \wedge T_{0j}^\beta } \left( \mathcal{L}_{\xi_j} \left( \partial B_j^\gamma \right) \right) - i_{T_{0j}^\alpha \wedge T_{0j}^\beta}\left( \mathcal{L}_{\xi_j} \mathcal{L}_{\xi_j} \left( d w_{0j}^\gamma \right) \right) \notag\\
&=  -2 i_{\left[ \xi_j, T_{0j}^\alpha\right]} i_{\left[ \xi_j , T_{0j}^\beta \right]}\left(dw_{0j}^\gamma \right) - 2 i_{A_j^\beta} i_{\left[ \xi_j, T_{0j}^\alpha \right]}\left( dw_{0j}^\gamma \right)^*  + 2 i_{A_j^\alpha} i_{\left[ \xi_j, T_{0j}^\beta \right]}\left( dw_{0j}^\gamma \right)^* -2 i_{A_j^\alpha} i_{\alpha_j\left( T_{0j}^\beta \right)}\left( dw_{0j}^\gamma \right)^\sharp + 2i_{A_j^\beta} i_{\alpha_j\left( T_{0j}^\alpha \right)} \left( d w_{0j}^\gamma \right)^\sharp \notag \\
&  - 2 i_{\left[ \xi_j, T_{0j}^\alpha \right]} i_{T_{0j}^\beta} \mathcal{L}_{\xi_j}\left( dw_{0j}^\gamma \right) + 2 i_{A_j^\alpha} i_{T_{0j}^\beta} \mathcal{L}_{\xi_j}\left(d w_{0j}^\gamma \right)^* - 2i_{T_{0j}^\beta} i_{\left[ \xi_j, T_{0j}^\alpha \right]}\left( \partial B_j^\gamma \right)^{\square} + 2i_{T_{0j}^\beta} i_{\alpha_j\left( T_{0j}^\alpha \right)}\left(\partial B_j^\gamma \right)^\circ  + 2i_{T_{0j}^\beta} i_{A_j^\alpha}\left( d\beta_j\left( w_{0j}^\gamma \right) \right)^\sharp \notag \\
&+ 2 i_{\left[ \xi_j, T_{0j}^\beta \right]} i_{T_{0j}^\alpha} \mathcal{L}_{\xi_j}\left( dw_{0j}^\gamma \right) - 2 i_{A_j^\beta} i_{T_{0j}^\alpha} \mathcal{L}_{\xi_j}\left(d w_{0j}^\gamma \right)^* + 2i_{T_{0j}^\alpha} i_{\left[ \xi_j, T_{0j}^\beta \right]}\left( \partial B_j^\gamma \right)^\square - 2i_{T_{0j}^\alpha} i_{\alpha_j\left( T_{0j}^\beta \right)}\left(\partial B_j^\gamma \right)^\circ - 2i_{T_{0j}^\alpha} i_{A_j^\beta}\left( d\beta_j\left( w_{0j}^\gamma \right) \right)^\sharp \notag
\end{align}

Let us compute the parts $*$ in $(\ref{a45})$
\begin{align}\label{a48}
&- 2 i_{A_j^\beta} i_{\left[ \xi_j, T_{0j}^\alpha \right]} \left( dw_{0j}^\gamma \right) + 2 i_{A_j^\alpha} i_{\left[ \xi_j, T_{0j}^\beta \right]}\left( dw_{0j}^\gamma \right) + 2 i_{A_j^\alpha} i_{T_{0j}^\beta} \mathcal{L}_{\xi_j}\left( dw_{0j}^\gamma \right)  - 2 i_{A_j^\beta} i_{T_{0j}^\alpha} \mathcal{L}_{\xi_j} \left( dw_{0j}^\gamma \right)\\
&= - 2 i_{A_j^\beta} \mathcal{L}_{\xi_j} i_{T_{0j}^\alpha} \left( dw_{0j}^\gamma \right) + 2 i_{A_j^\alpha} \mathcal{L}_{\xi_j} i_{T_{0j}^\beta} \left( d w_{0j}^\gamma \right)  = 2 i_{\left[ \xi_j, A_j^\beta \right]} i_{T_{0j}^\alpha} \left( dw_{0j}^\gamma \right) - 2 i_{\left[ \xi_j, A_j^\alpha \right]}  i_{T_{0j}^\beta}\left( d w_{0j}^\gamma \right) \notag
\end{align}

Let us compute the parts $\sharp$ in $(\ref{a45})$. Since $A_j^\alpha = \sum_{\eta=1}^p A_j^{\alpha \eta} T_{0j}^\eta$, we have
{\tiny{\begin{align}
&- 2 i_{A_j^\alpha} i_{\alpha_j\left(  T_{0j}^\beta \right)} \left( dw_{0j}^\gamma \right) +  2 i_{A_j^\beta} i_{\alpha_j\left( T_{0j}^\alpha \right)}\left( dw_{0j}^\gamma \right) + 2 i_{T_{0j}^\beta} i_{A_j^\alpha} \left( d\beta_j\left( w_{0j}^\gamma \right) \right) - 2 i_{T_{0j}^\alpha} i_{A_j^\beta} \left( d\beta_j\left( w_{0j}^\beta \right) \right) \label{a49}\\
&= 2\sum_{\eta=1}^p A_j^{\alpha \eta} i_{\alpha_j\left(T_{0j}^\beta \right)}i_{T_{0j}^\eta} \left(  dw_{0j}^\gamma \right) - 2 \sum_{\eta=1}^p A_j^{\beta \eta} i_{\alpha_j\left( T_{0j}^\alpha \right)} i_{T_{0j}^\eta} \left( dw_{0j}^\gamma \right) + 2 \sum_{\eta=1}^p A_j^{\alpha \eta} i_{T_{0j}^\eta \wedge T_{0j}^\beta}\left( d\beta_j\left( w_{0j}^\gamma \right) \right) - 2 \sum_{\eta=1}^p A_j^{\beta\eta} i_{T_{0j}^\eta} \wedge i_{T_{0j}^\alpha} \left(d \beta_j\left( w_{0j}^\gamma \right) \right) \notag \\
&= 2 \sum_{\eta=1}^p \sum_{\delta=1}^q A_j^{\alpha \eta} a_{0j\eta}^{\gamma \delta} i_{\alpha_j\left( T_{0j}^\beta \right)} \left( w_{0j}^\delta \right) + 2 \sum_{\eta=1}^p A_j^{\alpha \eta} i_{T_{0j}^\eta \wedge T_{0j}^\beta} \left( d\beta_j\left( w_{0j}^\gamma \right) \right) - 2 \sum_{\eta=1}^p\sum_{\delta=1}^q  A_j^{\beta \eta} a_{0j\eta}^{\gamma \delta} i_{\alpha_j\left( T_{0j}^\alpha \right)}\left( w_{0j}^\delta \right) - 2 \sum_{\eta=1}^p A_j^{\beta\eta} i_{T_{0j}^\eta} \wedge i_{T_{0j}^\alpha} \left(d \beta_j\left( w_{0j}^\gamma \right) \right) \notag\\
&= -2 \sum_{\eta=1}^p \sum_{\delta=1}^q A_j^{\alpha \eta} a_{0j\eta}^{\gamma \delta} i_{T_{0j}^\beta } \left(  \beta_j\left( w_{0j}^\delta \right) \right) + 2 \sum_{\eta=1}^p A_j^{\alpha \eta} i_{T_{0j}^\eta \wedge T_{0j}^\beta} \left( d\beta_j\left( w_{0j}^\gamma \right) \right) + 2 \sum_{\eta=1}^p\sum_{\delta=1}^q  A_j^{\beta \eta} a_{0j\eta}^{\gamma \delta} i_{ T_{0j}^\alpha }\left( \beta_j\left( w_{0j}^\delta \right) \right) - 2 \sum_{\eta=1}^p A_j^{\beta\eta} i_{T_{0j}^\eta} \wedge i_{T_{0j}^\alpha} \left(d \beta_j\left( w_{0j}^\gamma \right) \right) \notag
\end{align}}}

Let us compute that parts $\square$ in $(\ref{a45})$.
\begin{align*}
&-2 i_{T_{0j}^\beta} i_{\left[ \xi_j, T_{0j}^\alpha \right]} \left( \partial B_j^\gamma \right) + 2 i_{T_{0j}^\alpha} i_{\left[ \xi_j, T_{0j}^\beta \right]} \left( \partial B_j^\gamma \right)\\
&= - 2i_{T_{0j}^\beta}\mathcal{L}_{\xi_j} i_{T_{0j}^\alpha}\left(\partial B_j^\gamma \right) + 2 i_{T_{0j}^\beta} i_{T_{0j}^\alpha}\mathcal{L}_{\xi_j}\left( \partial B_j^\gamma \right) + 2i_{T_{0j}^\alpha}\mathcal{L}_{\xi_j} i_{T_{0j}^\beta}\left( \partial B_j^\gamma \right) - 2i_{T_{0j}^\alpha} i_{T_{0j}^\beta}\mathcal{L}_{\xi_j}\left( \partial B_j^\gamma \right) \\
&= 2 i_{\left[ \xi_j, T_{0j}^\beta \right]} i_{T_{0j}^\alpha}\left( \partial B_j^\gamma \right) + 4 i_{T_{0j}^\alpha \wedge T_{0j}^\beta}\left(\partial \left( \mathcal{L}_{\xi_j}\left( B_j^\gamma \right) \right) \right) - 2i_{\left[ \xi_j, T_{0j}^\alpha \right]} i_{T_{0j}^\beta}\left( \partial B_j^\gamma \right) 
\end{align*}
This implies that
\begin{align}\label{a47}
-2 i_{T_{0j}^\beta} i_{\left[ \xi_j, T_{0j}^\alpha \right]} \left( \partial B_j^\gamma \right) + 2 i_{T_{0j}^\alpha} i_{\left[ \xi_j, T_{0j}^\beta \right]} \left( \partial B_j^\gamma \right) = 2 i_{T_{0j}^\alpha \wedge T_{0j}^\beta}\left( \partial \left( \mathcal{L}_{\xi_j} \left( B_j^\gamma \right) \right) \right)
\end{align}

Let us compute the parts $\circ$ in $(\ref{a45})$.
{\tiny{\begin{align*}
& 2 i_{T_{0j}^\beta} i_{\alpha_j\left( T_{0j}^\alpha \right)} \left( \partial B_j^\gamma \right)  - 2 i_{T_{0j}^\alpha} i_{\alpha_j\left( T_{0j}^\beta \right)}\left( \partial B_j^\gamma \right)= - 2 i_{\alpha_j\left( T_{0j}^\alpha \right)} i_{T_{0j}^\beta}\left( \partial B_j^\gamma \right) + 2i_{\alpha_j\left( T_{0j}^\beta \right)}i_{T_{0j}^\alpha}\left( \partial B_j^\gamma \right)\\
&=- 2i_{\alpha_j\left( T_{0j}^\alpha \right)} \mathcal{L}_{T_{0j}^\beta}\left( B_j^\gamma \right) + 2i_{\alpha_j\left(  T_{0j}^\beta \right)} \mathcal{L}_{T_{0j}^\alpha}\left( B_j^\gamma \right) \\
&= 2 i_{\left[ T_{0j}^\beta, \alpha_j\left(T_{0j}^\alpha \right) \right]} \left(B_j^\gamma \right) - 2 \mathcal{L}_{T_{0j}^\beta} i_{\alpha_j\left(T_{0j}^\alpha \right)} \left(B_j^\gamma \right) - 2 i_{\left[ T_{0j}^\alpha, \alpha_j\left( T_{0j}^\beta \right) \right]}\left(B_j^\gamma \right)  +2 \mathcal{L}_{T_{0j}^\alpha} i_{\alpha_j\left(T_{0j}^\beta \right)}\left(B_j^\gamma \right)\\
&=2 i_{\left[ T_{0j}^\beta, \alpha_j\left(T_{0j}^\alpha \right) \right]} \left(B_j^\gamma \right) + 2 \mathcal{L}_{T_{0j}^\beta} i_{T_{0j}^\alpha } \left(\beta_j\left( B_j^\gamma \right) \right) -2 i_{\left[ T_{0j}^\alpha, \alpha_j\left( T_{0j}^\beta \right) \right]}\left( B_j^\gamma \right) -2 \mathcal{L}_{T_{0j}^\alpha} i_{T_{0j}^\beta}\left( \beta_j \left( B_j^\gamma \right) \right)\\
&= 2 i_{\left[ T_{0j}^\beta, \alpha_j\left(T_{0j}^\alpha \right) \right]} \left(B_j^\gamma \right)+ 2i_{\left[ T_{0j}^\beta, T_{0j}^\alpha \right]}\left(\beta_j\left( B_j^\gamma \right) \right) + 2 i_{T_{0j}^\alpha} \mathcal{L}_{T_{0j}^\beta} \left(  \beta_j\left( B_j^\gamma \right)\right) - 2 i_{\left[ T_{0j}^\alpha, \alpha_j\left( T_{0j}^\beta \right)\right]}\left(B_j^\gamma \right)- 2i_{\left[ T_{0j}^\alpha, T_{0j}^\beta \right]}\left(\beta_j\left( B_j^\gamma \right) \right) - i_{T_{0j}^\beta} \mathcal{L}_{T_{0j}^\alpha}\left(\beta_j\left(B_j^\gamma \right) \right)\\
&=2 i_{\left[ T_{0j}^\beta, \alpha_j\left( T_{0j}^\alpha \right)\right]}\left(B_j^\gamma \right) +2i_{\left[ T_{0j}^\beta, T_{0j}^\alpha \right]}\left(\beta_j\left( B_j^\gamma \right) \right) + 2 i_{T_{0j}^\alpha} \partial i_{T_{0j}^\beta} \left(  \beta_j\left( B_j^\gamma \right)\right) + 2 i_{T_{0j}^\alpha} i_{T_{0j}^\beta} \left( \partial \beta_j\left(  B_j^\gamma \right) \right) \\
&-2 i_{\left[ T_{0j}^\alpha, \alpha_j\left( T_{0j}^\beta \right)\right]}\left(B_j^\gamma \right) -2i_{\left[ T_{0j}^\alpha, T_{0j}^\beta \right]}\left(\beta_j\left( B_j^\gamma \right) \right) - 2 i_{T_{0j}^\beta} \partial i_{T_{0j}^\alpha} \left(  \beta_j\left( B_j^\gamma \right)\right) - 2 i_{T_{0j}^\beta} i_{T_{0j}^\alpha} \left( \partial \beta_j\left(  B_j^\gamma \right) \right)\\
&= 2 i_{\left[ T_{0j}^\beta, \alpha_j\left(T_{0j}^\alpha \right)\right]}\left(B_j^\gamma \right) + 4 i_{\alpha_j\left(\left[ T_{0j}^\alpha, T_{0j}^\beta \right] \right)}\left(B_j^\gamma \right) - 2 i_{T_{0j}^\alpha}\partial i_{\alpha_j\left( T_{0j}^\beta \right)}\left(B_j^\gamma \right)- 4 i_{T_{0j}^\alpha \wedge T_{0j}^\beta}\left(\partial \beta_j\left(B_j^\gamma \right) \right)- 2 i_{\left[ T_{0j}^\alpha, \alpha_j\left(T_{0j}^\beta \right) \right]}\left(B_j^\gamma \right) + 2i_{T_{0j}^\beta}\partial i_{\alpha_j\left(T_{0j}^\alpha \right)}\left( B_j^\gamma \right) \\
&= 2 i_{\left[ T_{0j}^\beta, \alpha_j\left(T_{0j}^\alpha \right)\right]}\left(B_j^\gamma \right) + 4 i_{\alpha_j\left(\left[ T_{0j}^\alpha, T_{0j}^\beta \right] \right)}\left(B_j^\gamma \right) -  2 i_{T_{0j}^\alpha} \mathcal{L}_{\alpha_j\left(T_{0j}^\beta \right)}\left(B_j^\gamma \right) + 2 i_{T_{0j}^\alpha}i_{\alpha_j\left(T_{0j}^\beta \right)}  \left(\partial B_j^\gamma \right) - 4 i_{T_{0j}^\alpha \wedge T_{0j}^\beta}\left(\partial \beta_j\left(B_j^\gamma \right) \right)- 2 i_{\left[ T_{0j}^\alpha, \alpha_j\left(T_{0j}^\beta \right) \right]}\left(B_j^\gamma \right) \\
&+ 2 i_{T_{0j}^\beta}\mathcal{L}_{\alpha_j\left( T_{0j}^\alpha \right)}\left( B_j^\gamma \right) - 2 i_{T_{0j}^\beta} i_{\alpha_j\left( T_{0j}^\alpha \right)} \left( \partial B_j^\gamma \right) \\
&=  \cancel{ 2 i_{\left[ T_{0j}^\beta, \alpha_j\left(T_{0j}^\alpha \right)\right]}\left(B_j^\gamma \right) }+ \cancel{4 i_{\alpha_j\left(\left[ T_{0j}^\alpha, T_{0j}^\beta \right] \right)}\left(B_j^\gamma \right)} + \cancel{2 i_{\left[ \alpha_j\left(  T_{0j}^\beta \right), T_{0j}^\alpha \right]}} + 2 i_{T_{0j}^\alpha} i_{\alpha_j\left( T_{0j}^\beta \right)}\left( \partial B_j^\gamma \right) - 4 i_{T_{0j}^\alpha \wedge T_{0j}^\beta}\left(\partial \beta_j\left(B_j^\gamma \right) \right)-  \cancel{ 2 i_{\left[ T_{0j}^\alpha, \alpha_j\left(T_{0j}^\beta \right) \right]}\left(B_j^\gamma \right)} \\
& - \cancel{2 i_{\left[ \alpha_j\left( T_{0j}^\alpha \right) , T_{0j}^\beta \right]} \left( B_j^\gamma \right)} - 2i_{T_{0j}^\beta} i_{\alpha_j\left( T_{0j}^\alpha \right)} \left( \partial B_j^\gamma \right)=  2 i_{T_{0j}^\alpha} i_{\alpha_j\left( T_{0j}^\beta \right)}\left( \partial B_j^\gamma \right) - 4 i_{T_{0j}^\alpha \wedge T_{0j}^\beta}\left(\partial \beta_j\left(B_j^\gamma \right) \right) - 2i_{T_{0j}^\beta} i_{\alpha_j\left( T_{0j}^\alpha \right)} \left( \partial B_j^\gamma \right)
\end{align*}}}
This implies that
\begin{align} \label{a46}
 2 i_{T_{0j}^\beta} i_{\alpha_j\left( T_{0j}^\alpha \right)} \left( \partial B_j^\gamma \right)  - 2 i_{T_{0j}^\alpha} i_{\alpha_j\left( T_{0j}^\beta \right)}\left( \partial B_j^\gamma \right) = - 2 i_{T_{0j}^\alpha \wedge T_{0j}^\beta}\left( \partial \beta_j\left( B_j^\gamma \right) \right)
\end{align}

In view of $(\ref{a45})$ and $(\ref{a48})$ and $(\ref{a49})$ and $(\ref{a47})$ and $(\ref{a46})$, we have to show that
\begin{align*}
&i_{\left[ \xi_j, \left[ \xi_j, T_{0j}^\alpha \right] \right]} i_{T_{0j}^\beta}\left( d w_{0j}^\gamma \right) - i_{\left[ \xi_j, \left[ \xi_j, T_{0j}^\beta \right] \right]} i_{T_{0j}^\alpha} \left( d w_{0j}^\gamma \right) - i_{T_{0j}^\alpha \wedge T_{0j}^\beta}\left(\mathcal{L}_{\xi_j} \mathcal{L}_{\xi_j}\left( dw_{0j}^\gamma \right) \right)\\
&= -2 i_{\left[ \xi_j, T_{0j}^\alpha \right]} i_{\left[ \xi_j,  T_{0j}^\beta \right]} \left( dw_{0j}^\gamma \right) - 2 i_{\left[ \xi_j, T_{0j}^\alpha \right]} i_{T_{0j}^\beta} \mathcal{L}_{\xi_j} \left( dw_{0j}^\gamma \right) + 2 i_{\left[ \xi_j, T_{0j}^\beta \right]} i_{T_{0j}^\alpha} \mathcal{L}_{\xi_j} \left( d w_{0j}^\gamma \right)
\end{align*}
In fact, we note that (here we omit $\left( dw_{0j}^\gamma \right)$).

\begin{align*}
&- i_{\left[ \xi_j, T_{0j}^\alpha\right]} i_{\left[ \xi_j, T_{0j}^\beta\right]} - i_{\left[ \xi_j, T_{0j}^\alpha \right]}i_{T_{0j}^\beta } \mathcal{L}_{\xi_j} + i_{\left[ \xi_j, T_{0j}^\beta \right]} i_{T_{0j}^\alpha} \mathcal{L}_{\xi_j}\\
&=- i_{\left[ \xi_j, T_{0j}^\alpha\right]} i_{\left[ \xi_j, T_{0j}^\beta\right]} + i_{\left[ \xi_j, T_{0j}^\alpha \right]} i_{\left[ \xi_j, T_{0j}^\beta \right]}  - i_{\left[ \xi_j, T_{0j}^\alpha \right]} \mathcal{L}_{\xi_j} i_{T_{0j}^\beta } - i_{\left[ \xi_j, T_{0j}^\beta \right]} i_{\left[ \xi_j, T_{0j}^\alpha \right]} + i_{\left[ \xi_j, T_{0j}^\beta\right]} \mathcal{L}_{\xi_j} i_{T_{0j}^\alpha } \\
&= i_{\left[ \xi_j, \left[ \xi_j, T_{0j}^\alpha \right] \right]} i_{T_{0j}^\beta} - \mathcal{L}_{\xi_j} i_{\left[ \xi_j, T_{0j}^\alpha \right]} i_{T_{0j}^\beta} - i_{\left[ \xi_j, T_{0j}^\beta\right]} i_{\left[ \xi_j, T_{0j}^\alpha \right]} - i_{\left[ \xi_j, \left[ \xi_j, T_{0j}^\beta \right] \right]} i_{T_{0j}^\alpha} + \mathcal{L}_{\xi_j} i_{\left[ \xi_j, T_{0j}^\beta \right]} i_{T_{0j}^\alpha}\\
&= i_{\left[ \xi_j, \left[ \xi_j, T_{0j}^\alpha \right] \right]} i_{T_{0j}^\beta}  - i_{\left[ \xi_j, T_{0j}^\beta\right]} i_{\left[ \xi_j, T_{0j}^\alpha \right]} - i_{\left[ \xi_j, \left[ \xi_j, T_{0j}^\beta \right] \right]} i_{T_{0j}^\alpha}  + \mathcal{L}_{\xi_j} i_{T_{0j}^\beta} \mathcal{L}_{\xi_j} i_{T_{0j}^\alpha} - \mathcal{L}_{\xi_j} i_{T_{0j}^\beta} i_{T_{0j}^\alpha} \mathcal{L}_{\xi_j} - \mathcal{L}_{\xi_j} i_{T_{0j}^\alpha}  \mathcal{L}_{\xi_j}  i_{T_{0j}^\beta} + \mathcal{L}_{\xi_j} i_{T_{0j}^\alpha} i_{T_{0j}^\beta} \mathcal{L}_{\xi_j}\\
&=  i_{\left[ \xi_j, \left[ \xi_j, T_{0j}^\alpha \right] \right]} i_{T_{0j}^\beta}  - i_{\left[ \xi_j, T_{0j}^\beta\right]} i_{\left[ \xi_j, T_{0j}^\alpha \right]} - i_{\left[ \xi_j, \left[ \xi_j, T_{0j}^\beta \right] \right]} i_{T_{0j}^\alpha}  - i_{\left[ \xi_j, T_{0j}^\beta \right]} i_{T_{0j}^\alpha} \mathcal{L}_{\xi_j} - i_{T_{0j}^\beta} \mathcal{L}_{\xi_j} i_{T_{0j}^\alpha} \mathcal{L}_{\xi_j} + i_{\left[ \xi_j, T_{0j}^\alpha \right]} i_{T_{0j}^\beta} \mathcal{L}_{\xi_j} + i_{T_{0j}^\alpha} \mathcal{L}_{\xi_j} i_{T_{0j}^\beta} \mathcal{L}_{\xi_j}\\
&+ \mathcal{L}_{\xi_j} i_{T_{0j}^\beta} \mathcal{L}_{\xi_j} i_{T_{0j}^\alpha} - \mathcal{L}_{\xi_j} i_{T_{0j}^\alpha} \mathcal{L}_{\xi_j} i_{T_{0j}^\beta}
\end{align*}
Hence it remains to show that (here we omit $\left( dw_{0j}^\gamma \right)$).
\begin{align*}
- i_{T_{0j}^\alpha \wedge T_{0j}^\beta} \left(\mathcal{L}_{\xi_j} \mathcal{L}_{\xi_j} \right)=\mathcal{L}_{\xi_j} i_{T_{0j}^\beta}\mathcal{L}_{\xi_j} i_{T_{0j}^\alpha} - \mathcal{L}_{\xi_j} i_{T_{0j}^\alpha} \mathcal{L}_{\xi_j} i_{T_{0j}^\beta} - i_{T_{0j}^\beta} \mathcal{L}_{\xi_j} i_{T_{0j}^\alpha} \mathcal{L}_{\xi_j} + i_{T_{0j}^\alpha} \mathcal{L}_{\xi_j} i_{T_{0j}^\beta} \mathcal{L}_{\xi_j}
\end{align*}
which follows from (here we omit $\left(dw_{0j}^\gamma\right)$)
\begin{align*}
&- i_{T_{0j}^\alpha \wedge T_{0j}^\beta} \left(\mathcal{L}_{\xi_j} \mathcal{L}_{\xi_j} \right)= i_{T_{0j}^\alpha} i_{T_{0j}^\beta} \mathcal{L}_{\xi_j} \mathcal{L}_{\xi_j}= - i_{T_{0j}^\alpha} i_{\left[\xi_j, T_{0j}^\beta \right]}\mathcal{L}_{\xi_j} + i_{T_{0j}^\alpha} \mathcal{L}_{\xi_j} i_{T_{0j}^\beta} \mathcal{L}_{\xi_j} = i_{\left[ \xi_j, T_{0j}^\beta \right]} i_{T_{0j}^\alpha} \mathcal{L}_{\xi_j} + i_{T_{0j}^\alpha} \mathcal{L}_{\xi_j} i_{T_{0j}^\beta} \mathcal{L}_{\xi_j}\\
&=\mathcal{L}_{\xi_j} i_{T_{0j}^\beta} i_{T_{0j}^\alpha} \mathcal{L}_{\xi_j} - i_{T_{0j}^\beta} \mathcal{L}_{\xi_j} i_{T_{0j}^\alpha} \mathcal{L}_{\xi_j} + i_{T_{0j}^\alpha} \mathcal{L}_{\xi_j} i_{T_{0j}^\alpha} \mathcal{L}_{\xi_j} = - \mathcal{L}_{\xi_j} i_{T_{0j}^\beta} i_{\left[ \xi_j, T_{0j}^\alpha \right]} + \mathcal{L}_{\xi_j} i_{T_{0j}^\beta}  \mathcal{L}_{\xi_j} i_{T_{0j}^\alpha} - i_{T_{0j}^\beta} \mathcal{L}_{\xi_j} i_{T_{0j}^\alpha} \mathcal{L}_{\xi_j} + i_{T_{0j}^\alpha} \mathcal{L}_{\xi_j} i_{T_{0j}^\alpha} \mathcal{L}_{\xi_j}\\
&=\mathcal{L}_{\xi_j} i_{\left[ \xi_j, T_{0j}^\alpha \right]} i_{T_{0j}^\beta} + \mathcal{L}_{\xi_j} i_{T_{0j}^\beta}  \mathcal{L}_{\xi_j} i_{T_{0j}^\alpha} - i_{T_{0j}^\beta} \mathcal{L}_{\xi_j} i_{T_{0j}^\alpha} \mathcal{L}_{\xi_j} + i_{T_{0j}^\alpha} \mathcal{L}_{\xi_j} i_{T_{0j}^\alpha} \mathcal{L}_{\xi_j}\\
&= -\mathcal{L}_{\xi_j} i_{T_{0j}^\alpha} \mathcal{L}_{\xi_j} i_{T_{0j}^\beta} + \mathcal{L}_{\xi_j} i_{T_{0j}^\beta}  \mathcal{L}_{\xi_j} i_{T_{0j}^\alpha} - i_{T_{0j}^\beta} \mathcal{L}_{\xi_j} i_{T_{0j}^\alpha} \mathcal{L}_{\xi_j} + i_{T_{0j}^\alpha} \mathcal{L}_{\xi_j} i_{T_{0j}^\alpha} \mathcal{L}_{\xi_j}
\end{align*}
since $i_{T_{0j}^\alpha \wedge T_{0j}^\beta}\left( dw_{0j}^\gamma \right)=0$. This completes the proof of Proposition \ref{a35}.
\end{proof}

\section{Infinitesimal unfoldings and obstructions by using \v Cech resolution of $\tilde{\mathcal{E}}_{\mathcal{N}_{\mathcal{F}_0}^*}^\bullet$}

\begin{definition} \label{u1}
Suppose that given a domain $B\subset \mathbb{C}^m$, there is a set $\left\{\left(M_t, \mathcal{N}_{\mathcal{F}_t}^*\right)|t\in  B \right\}$ of $n$-dimensional compact $($singularly$)$ foliated complex manifolds $\left(M_t, \mathcal{N}_{\mathcal{F}_t}^*\right)$ of codimension $q$, so that for each $t\in B$, we have an exact sequence
\begin{align}
0\to \mathcal{N}_{\mathcal{F}_t}^*\to \Omega_{M_t}^1\to \Omega_{M_t}^1/\mathcal{N}_{\mathcal{F}_t}^*\to 0
\end{align}
where $\Omega_{M_t}^1/\mathcal{N}_{\mathcal{F}_t}^*$ is a torsion-free $\mathcal{O}_{M_t}$-module. We say that $\left\{\left(M_t, \mathcal{N}_{\mathcal{F}_t}^*\right)|t\in B \right\}$ is a family of unfoldings or unfolded $($complex$)$ analytic family if there exist a complex manifold $\mathcal{M}$ and a holomorphic map $\pi:\mathcal{M}\to B$ such that
\begin{enumerate}
\item we have an exact sequence
\begin{align*}
0\to \tilde{\mathcal{N}}_\mathcal{F}^*\to \Omega_\mathcal{M}^1\to \Omega_\mathcal{M}^1/\tilde{\mathcal{N}}_{\mathcal{F}}^*\to 0
\end{align*}
such that $\tilde{\mathcal{N}}_{\mathcal{F}}^*$ is a coherent subsheaf of the cotangent sheaf $\Omega_\mathcal{M}^1$ on $\mathcal{M}$.
\item the rank of the Jacobian of $\pi$ is equal to $m$ at every point of $\mathcal{M}$.
\item from the exact sequence $0\to \pi^* \Omega_B^1\to \Omega_\mathcal{M}^1\xrightarrow{\alpha} \Omega_{\mathcal{M}/B}^1 \to 0$, let us denote the image of $\tilde{\mathcal{N}}_{\mathcal{F}}^*$ under $\alpha$ by $\mathcal{N}_{\mathcal{F}}^*$. Then $\alpha$ is isomorphic on $\tilde{\mathcal{N}}_{\mathcal{F}}^*$, i.e $ \tilde{\mathcal{N}}_{\mathcal{F}}^* \cong \alpha\left( \tilde{\mathcal{N}}_{\mathcal{F}}^*  \right)= \mathcal{\mathcal{N}}_{\mathcal{F}}^*$ and the exact sequence
\begin{align*}
0\to \mathcal{N}_{\mathcal{F}}^*\to \Omega_{\mathcal{M}/B}^1\to \Omega_{\mathcal{M}/B}^1/\mathcal{N}_{\mathcal{F}}^*\to 0
\end{align*}
defines a foliated $($complex$)$ analytic family of $\left\{\left(M_t, \mathcal{N}_{\mathcal{F}_t}^* \right)\right\}$ in terms of cotangent sheaves in the sense of  \textnormal{Definition \ref{n3}}, so that $\Omega_{\mathcal{M}/B}^1/\mathcal{N}_{\mathcal{F}}^*$ is flat over $B$.
\item $\tilde{\mathcal{N}}_{\mathcal{F}}^*$ satisfies the integrability condition. In other words,  let $\tilde{S}$ be the locus of points of $\mathcal{M}$ such that $\frac{\Omega_{\mathcal{M}}^1}{ \tilde{\mathcal{N}}_\mathcal{F}^*}$ is not locally free.  Then we have $d \tilde{\mathcal{N}}_{\mathcal{F}}^*\subset \tilde{\mathcal{N}}_{\mathcal{F}}^* \bigwedge \Omega_\mathcal{M}^1$ at every point of $\mathcal{M}- \tilde{S}$.  
\end{enumerate}
We will denote the unfolded analytic family by $(\mathcal{M}, \tilde{\mathcal{N}}_{\mathcal{F}}^*, B, \pi)$ and the associated foliated analytic family $\left(\mathcal{M}, \mathcal{N}_\mathcal{F}^*, B, \pi \right)$ in terms of cotangent sheaves.
\end{definition}

\begin{remark}
Let $\left( \mathcal{M}, \tilde{\mathcal{N}}_{\mathcal{F}}^*, B, \pi \right)$ with $\tilde{\mathcal{N}}_\mathcal{F}^*$ locally free be an unfolded analytic family. Let $\Delta$ be an open set of $B$. Then the restriction $\left( \mathcal{M}_\Delta =\pi^{-1}\left(\Delta \right) , \tilde{\mathcal{N}}_\mathcal{F}^*|_\Delta,  \Delta, \pi |_\Delta   \right)$ is also an unfolded complex analytic family. We will denote the family by $\left( \mathcal{M}_\Delta, \tilde{\mathcal{N}}_{\mathcal{F}_\Delta}^*, \Delta, \pi \right)$.
\end{remark}

\begin{remark}\label{ou1}
Given an unfolded analytic family $\left( \mathcal{M}, \tilde{\mathcal{N}}_\mathcal{F}^*, B, \pi \right)$ with $\tilde{\mathcal{N}}_{\mathcal{F}}^*$ locally free as in \textnormal{Definition \ref{u1}}, we can choose a system of local complex coordinates $\{\xi_i,...,\xi_j,...\}, \xi_j:p\to \xi_j(p)$, and coordinate polydisks $\mathcal{U}_j$ with respect to $\xi_j$, satisfying $(1)$ $\xi_j(p)=\left(\xi_j^1(p),...,\xi_j^n(p),t_1,...,t_m\right),\left(t_1,...,t_m\right)=\pi(p)$ and $(2)$ $\mathcal{U}=\left\{\mathcal{U}_j|j=1,2,...\right\}$ is locally finite. Then $\left\{p\to (\xi_j^1(p),...,\xi_j^n(p))|\mathcal{U}_j\cap M_t\ne \emptyset \right\}$ gives a system of local complex coordinates on $M_t$. In terms of these coordinates, $\pi$ is the projection given by $\pi: \left(\xi_j^1,...,\xi_j^n, t_1,...,t_m \right)\to \left(t_1,...,t_m\right)$. For $j,k$ with $\mathcal{U}_j\cap \mathcal{U}_k\ne \emptyset$, we denote that coordinate transformation from $\xi_k$ to $\xi_j$ by $f_{jk}:\left(\xi_k^1,...,\xi_k^n,t \right)\to \left(\xi_j^1,...,\xi_j^n,t \right)=f_{jk}\left(\xi_k^1,...,\xi_k^n,t \right)$. Thus $f_{jk}$ is given by $\xi_j^\alpha=f_{jk}^\alpha\left(\xi_k^1,...,\xi_k^n,t_1,...,t_m\right), \alpha=1,...,n$. On the other hand, we may assume that the locally free sheaf $\tilde{\mathcal{N}}_\mathcal{F}^*$ is trivialized in terms of the covering $\mathcal{U}$ so that $\Gamma\left(\mathcal{U}_j, \tilde{\mathcal{N} }_\mathcal{F}^* \right)\cong \bigoplus^q \Gamma\left(\mathcal{U}_j, \mathcal{O}_\mathcal{M} \right)$ which is generated by holomorphic $1$-forms $\omega_j^\alpha(\xi_j,t)=w_j^\alpha(\xi_j,t) + v_j^\alpha(\xi_j,t)\in \Gamma\left(\mathcal{U}_j,\tilde{\mathcal{N}}_\mathcal{F}^* \right),\alpha=1,...,q$ of the form
\begin{align}\label{u2}
\omega_j^\alpha(\xi_j,t):=w_j^\alpha(\xi_j,t)+v_j^\alpha(\xi_j,t)=\sum_{\beta=1}^n w_j^{\alpha\beta}(\xi_j,t_j)d\xi_j^\beta+\sum_{\eta=1}^m v_j^{\alpha\eta}(\xi_j,t_j) dt_j^\eta,\,\,\,\,\alpha=1,...,q
\end{align}
for some $w_j^{\alpha}(\xi_j,t)$ and $v_j^{\alpha\eta}(\xi_j,t)\in \Gamma\left( \mathcal{U}_j, \mathcal{O}_\mathcal{M} \right)$, and satisfies
\begin{align}\label{u8}
\omega_j^1(\xi_j,t)\wedge \cdots \wedge \omega_j^q(\xi_j,t)\wedge d\omega_j^\alpha(\xi_j,t)=0,\,\,\,\,\,\alpha=1,...,q\,\,\,\textnormal{on}\,\,\,\mathcal{U}_j.
\end{align}
For $j,k$ with $\mathcal{U}_j\cap \mathcal{U}_k\ne \emptyset$, we have
\begin{align}\label{ou2}
\omega_j^\alpha(\xi_j,t)= \sum_{\beta=1}^q h_{jk}^{\alpha\beta}(\xi_k,t) \omega_k^\beta(\xi_k,t) \iff \left( w_j^\alpha(\xi_j,t)+ v_j^\alpha(\xi_j,t) \right) = \sum_{\beta=1}^q h_{jk}^{\alpha\beta} \left( w_k^\beta(\xi_k,t)+ v_k^\beta(\xi_k,t) \right)
\end{align}
for some $h_{jk}^{\alpha\beta}(\xi_k,t)\in \Gamma\left( \mathcal{U}_j\cap \mathcal{U}_j, \mathcal{O}_\mathcal{M} \right)$,  and for $i,j,k$ with $\mathcal{U}_i\cap \mathcal{U}_j\cap \mathcal{U}_k\ne \emptyset$, we have
\begin{align}
h_{ik}^{\alpha \gamma}(\xi_k,t)= \sum_{\beta=1}^q h_{ij}^{\alpha\beta}\left(f_{jk}(\xi_k,t),t\right) h_{jk}^{\beta \gamma}(\xi_k,t)
\end{align}
More specifically $(\ref{ou2})$ implies that
\begin{align}
\sum_{\gamma,\eta=1}^n w_j^{\alpha\gamma}(f_{jk}(\xi_k,t),t)\frac{\partial f_{jk}^\gamma(\xi_k,t)}{\partial \xi_k^\eta} d\xi_k^\eta=\sum_{\beta=1}^q\sum_{\eta=1}^n &h_{jk}^{\alpha\beta}(\xi_k,t)w_k^{\beta\eta}(\xi_k,t)d\xi_k^\eta,\,\,\,\,\,\alpha=1,...,q\\
\sum_{\gamma=1}^n\sum_{\eta=1}^m w_j^{\alpha\gamma}\left(f_{jk}(\xi_k,t),t\right)\frac{\partial f_{jk}^\gamma(\xi_k,t)}{\partial t^\eta}dt^\eta + \sum_{\eta=1}^m v_j^{\alpha\eta}\left(f_{jk}(\xi_k,t),t \right)dt^\eta & = \sum_{\beta=1}^q \sum_{\eta=1}^m h_{jk}^{\alpha\beta}(\xi_k,t) v_k^{\beta\eta}(\xi_k,t)dt^\eta,\,\,\,\,\,\,\,\,\,\alpha=1,...,q. \label{u6}
\end{align}

\end{remark}

Given an unfolded complex analytic family $\left( \mathcal{M} , \tilde{\mathcal{N}}_\mathcal{F}^*, B, \pi \right)$ with $\tilde{\mathcal{N}}_\mathcal{F}^*$ locally free of unfoldings of $\left( M, \mathcal{N}_{\mathcal{F}_0}^* \right) = \pi^{-1}(0), 0 \in B$, the infinitesimal unfoldings in the unfolded complex analytic family are encoded in $\mathbb{H}^1\left( M, \tilde{\mathcal{N}}_{\mathcal{F}_0}^{*\bullet} \right)$ which is the first cohomology group of the following unfolded dual leaf complex
\begin{align}
\tilde{\mathcal{N}}_{\mathcal{F}_0}^{*\bullet}:  \Theta_M \to \mathscr{H}om_{\mathcal{O}_M}\left( \mathcal{N}_{\mathcal{F}_0}^*, \frac{\Omega_M^1}{\mathcal{N}_{\mathcal{F}_0}^*} \right) \bigoplus \mathscr{H}om_{\mathcal{O}_M}\left( \mathcal{N}_{\mathcal{F}_0}^*, \mathcal{O}_M \right) \to \mathscr{H}om_{\mathcal{O}_M}\left( \mathcal{N}_{\mathcal{F}_0}^*, \tilde{\mathcal{S}}^2 \right) \bigoplus \mathscr{H}om_{\mathcal{O}_M}\left( \mathcal{N}_{\mathcal{F}_0}^*, \tilde{\mathcal{S}}^1 \right) \to \cdots
\end{align} 
We will denote the $i$-th cohomology group by $\mathbb{H}^i\left( M, \tilde{\mathcal{N}}_{\mathcal{F}_0}^{*\bullet} \right)$. We can compute $\mathbb{H}^i\left( M, \tilde{\mathcal{N}}_{\mathcal{F}_0}^{*\bullet} \right)$ by the following \v Cech resolution of $\tilde{\mathcal{N}}_{\mathcal{F}_0}^{*\bullet}$. Here $\delta$ is the \v Cech map and $\mathcal{U}_0=\mathcal{U}\cap M = \left\{ U_j:=\mathcal{U}_j\cap M |j=1, 2, ... \right\}$ is an open covering of $M$:
{\Small{\begin{equation}\label{ou5}
\begin{CD}
\cdots\\
@AG_2AA \\
C^0\left(\mathcal{U}_0, \left( \left( \mathcal{N}_{\mathcal{F}_0}^* \right)^*\otimes  \tilde{\mathcal{S}}^2 \right) \bigoplus \left( \left( \mathcal{N}_{\mathcal{F}_0}^* \right)^* \otimes \tilde{\mathcal{S}}^1 \right) \right) @>-\delta>> \cdots \\
@AG_1AA @AG_1AA \\
C^0\left( \mathcal{U}_0,  \left(\left(\mathcal{N}_{\mathcal{F}_0}^*\right)^* \otimes \frac{\Omega_{M}^1}{\mathcal{N}_{\mathcal{F}_0}^*}\right) \bigoplus \left( \mathcal{N}_{\mathcal{F}_0}^* \right)^*\right) @>\delta>> C^1\left( \mathcal{U}_0, \left(\left(\mathcal{N}_{\mathcal{F}_0}^*\right)^* \otimes \frac{\Omega_{M}^1}{\mathcal{N}_{\mathcal{F}_0}^*}\right) \bigoplus \left( \mathcal{N}_{\mathcal{F}_0}^* \right)^*\right)  @>-\delta>> \cdots \\
@AG_0AA @AG_0AA @AG_0AA \\
C^0\left(\mathcal{U}_0, \Theta_{M} \right) @>-\delta>> C^1\left( \mathcal{U}_0, \Theta_{M} \right) @>\delta>> C^2\left( \mathcal{U}_0, \Theta_{M} \right) @>-\delta >> \cdots
\end{CD}
\end{equation}}}
Then we can define the unfolded Kodaira-Spencer map at $0\in B$ which encodes the information of infinitesimal unfoldings of $\pi^{-1}=\left( M, \mathcal{N}_{\mathcal{F}_0}^* \right)$ in the unfolded analytic family $\left( \mathcal{M}, \tilde{\mathcal{N}}_\mathcal{F}^*, B, \pi \right)$ (see Part III).

\begin{definition}[unfolded Kodaira-Spencer map]\label{u10}
Let $\left(\mathcal{M}, \tilde{\mathcal{N}}_\mathcal{F}^*, B, \pi\right)$ with $\tilde{\mathcal{N}}_{\mathcal{F}}^*$ locally free be an unfolded complex analytic family of unfoldings of $\left(M_t, \mathcal{N}_{\mathcal{F}_t}^*\right)=\pi^{-1}(t),t\in B$, where $B$ is a domain of $\mathbb{C}^m$. We keep the notations in \textnormal{Remark \ref{ou1}}. For a tangent vector $\frac{\partial}{\partial t}=\sum_{\eta=1}^m c_\eta \frac{\partial}{\partial t^\eta}, c_\eta \in \mathbb{C}$, of $B$,  the unfolded Kodaira-Spencer map at $t$ is defined to be a $\mathbb{C}$-linear map
{\Small{\begin{align*}
\tilde{\varphi}_t:T_t(B)&\to \mathbb{H}^1\left(M_t, \tilde{\mathcal{N}}_{\mathcal{F}_t}^{*\bullet} \right)\\
                \frac{\partial}{\partial t}&\mapsto \frac{\partial \left(M_t, \tilde{\mathcal{N}}_{\mathcal{F}_t}^* \right)}{\partial t}\\
                &:=\left(\left\{  \sum_{\alpha=1}^n \frac{\partial f_{jk}^\alpha(\xi_k,t)}{\partial t}\frac{\partial}{\partial \xi_j^\alpha} \right\}, \left\{ \left(w_j^\alpha(\xi_j,t)\mapsto -\overline{\frac{\partial w_j^\alpha(\xi_j,t)}{\partial t}}\right)_{\alpha=1,...,q} \right\} , \left\{ \left( w_j^\alpha(\xi_j,t) \mapsto - \sum_{\eta=1}^m c_\eta v_j^{\alpha \eta}(\xi_j,t) \right)_{\alpha=1,...,q}   \right\} \right) \\
                 &\in C^1\left(\mathcal{U}_t, \Theta_{M_t} \right) \bigoplus C^0\left( \mathcal{U}_t, \mathscr{H}om_{\mathcal{O}_{M_t}}\left(\mathcal{N}_{\mathcal{F}_t}^*, \frac{\Omega_{M_t}^1}{\mathcal{N}_{\mathcal{F}_t}^*} \right) \right) \bigoplus C^0\left( \mathcal{U}_t,  \mathscr{H}om_{\mathcal{O}_{M_t}} \left( \mathcal{N}_{\mathcal{F}_t}^*, \mathcal{O}_{M_t} \right) \right)
\end{align*}}}
Here $\frac{\partial w_j^\alpha(\xi_j,t)}{\partial t}|_{t=0}:=\sum_{\beta=1}^n \frac{\partial w_j^{\alpha\beta}(\xi_j,t)}{\partial t}|_{t=0}d\xi_j^\beta $ is the derivative of $w_j^\alpha(\xi_j,t)$ in $\left(\ref{u2} \right)$ with respect to $t$ at $0$, and $\overline{\frac{\partial w_j^\alpha(\xi_j,t)}{\partial t}}|_{t=0}$ is the image of $\frac{\partial w_j^\alpha(\xi_j,t)}{\partial t}|_{t=0}$ under $\Omega_M^1 \to \frac{\Omega_M^1}{\mathcal{N}_{\mathcal{F}_0}^*}$.
\end{definition}

In Part III we have shown that the complex of sheaves
{\small{\begin{align}\label{ou4}
\tilde{\mathcal{E}}_{\mathcal{N}_{\mathcal{F}_0}^*}^\bullet: \mathcal{E}_{\mathcal{N}_{\mathcal{F}_0}^*}\to \mathscr{H}om_{\mathcal{O}_M}\left( \mathcal{N}_{\mathcal{F}_0}^*,  \Omega_M^1 \right)  \bigoplus \mathscr{H}om_{\mathcal{O}_M}&\left( \mathcal{N}_{\mathcal{F}_0}^*, \mathcal{O}_M \right) \to  \mathscr{H}om_{\mathcal{O}_M}\left( \mathcal{N}_{\mathcal{F}_0}^* , \tilde{\mathcal{S}}^2 \right) \bigoplus \mathscr{H}om_{\mathcal{O}_M} \left( \mathcal{N}_{\mathcal{F}_0}^*, \tilde{\mathcal{S}}^1 \right) \\&\to \mathscr{H}om_{\mathcal{O}_M}\left( \mathcal{N}_{\mathcal{F}_0}^* , \tilde{\mathcal{S}}^3 \right) \bigoplus \mathscr{H}om_{\mathcal{O}_M}\left( \mathcal{N}_{\mathcal{F}_0}^*,  \tilde{\mathcal{S}}^2 \right) \to \cdots \notag
\end{align}}}
also controls unfoldings of $\left( M , \mathcal{N}_{\mathcal{F}_0}^* \right)$. We can actually study obstructions of unfoldings of $\left(M, \mathcal{N}_{\mathcal{F}_0}^* \right)$ by using the complex of sheaves $\tilde{\mathcal{N}}_{\mathcal{F}_0}^{*\bullet}$.  However obstructions of deformations of the underlying locally free sheaf $\mathcal{N}_{\mathcal{F}_0}^*$ are implicit in $\mathbb{H}^2\left( M, \tilde{\mathcal{N}}_{\mathcal{F}_0}^{*\bullet} \right)$, i.e. in \v Cech resolution $(\ref{ou5})$ of $\tilde{\mathcal{N}}_{\mathcal{F}_0}^{*\bullet}$. However the complex of sheaves $\tilde{\mathcal{E}}_{\mathcal{N}_{\mathcal{F}_0}^*}^\bullet$ explicitly shows the infinitesimal deformations and obstructions of deformations of the underlying locally free sheaf $\mathcal{N}_{\mathcal{F}_0}^*$. Hence it is desirable to use $\tilde{\mathcal{E}}_{\mathcal{N}_{\mathcal{F}_0}^*}^\bullet$ instead of $\tilde{\mathcal{N}}_{\mathcal{F}_0}^{*\bullet}$ to identify obstructions clearly. Let us identify the infinitesimal unfoldings in terms of the complex (\ref{oc4}) of sheaves  $\tilde{\mathcal{E}}_{\mathcal{N}_{\mathcal{F}_0}^*}^\bullet$.

\subsection{Identification of infinitesimal unfoldings by using \v Cech resolution of $\tilde{\mathcal{E}}_{\mathcal{N}_{\mathcal{F}_0}^*}^\bullet$} \

Let $\left( \mathcal{M}, \tilde{\mathcal{N}}_\mathcal{F}^*, B, \pi \right)$ be an unfolded analytic family with $\tilde{\mathcal{N}}_\mathcal{F}^*$ locally free as in Definition \ref{u1}. We keep the notations in \textnormal{Remark \ref{ou1}}, and we keep the notations in subsection \ref{oc8} for the underlying foliated analytic family $\left( \mathcal{M} , \mathcal{N}_\mathcal{F}^*, B, \pi \right)$ induced from $\left( \mathcal{M}, \tilde{\mathcal{N}}_\mathcal{F}^*, B, \pi \right)$. Let $\mathcal{U}_t= \mathcal{U}\cap M_t =\left\{ U_j^t := \mathcal{U}_j\cap M_t \right\}$ be an open covering of $M_t$. We recall that the definition of $\mathcal{E}_{\mathcal{N}_{\mathcal{F}_t}^*}$ on $M_t$ with respect to the open covering $\left\{ U_j^t \right\}$ from Part II, and let us consider the \v Cech resolution of $\tilde{\mathcal{E}}_{\mathcal{N}_{\mathcal{F}_0}^*}^\bullet$.

{\Small{\begin{equation}\label{ou3}
\begin{CD}
\cdots\\
@AG_2^t AA \\
C^0\left(\mathcal{U}_t, \left( \left( \mathcal{N}_{\mathcal{F}_t}^* \right)^*\otimes  \tilde{\mathcal{S}}_t^2 \right) \bigoplus \left( \left( \mathcal{N}_{\mathcal{F}_0}^* \right)^* \otimes \tilde{\mathcal{S}}^1 \right) \right) @>-\delta>> \cdots \\
@AG_1'^t AA @AG_1'^t AA \\
C^0\left( \mathcal{U}_t, \left(\left(\mathcal{N}_{\mathcal{F}_t}^*\right)^* \otimes \Omega_{M_t}^1 \right) \bigoplus \left( \mathcal{N}_{\mathcal{F}_t}^* \right)^*\right) @>\delta>> C^1\left( \mathcal{U}_t, \left(\left(\mathcal{N}_{\mathcal{F}_t}^*\right)^* \otimes  \Omega_{M_t}^1 \right) \bigoplus \left( \mathcal{N}_{\mathcal{F}_t}^* \right)^*\right)  @>-\delta>> \cdots \\
@AG_0'^tAA @AG_0'^tAA @AG_0'^tAA \\
C^0\left(\mathcal{U}_t, \mathcal{E}_{\mathcal{N}_{\mathcal{F}_0}^*} \right) @>-\delta>> C^1\left( \mathcal{U}_t, \mathcal{E}_{\mathcal{N}_{\mathcal{F}_t}^*} \right) @>\delta>> C^2\left( \mathcal{U}_t, \mathcal{E}_{\mathcal{N}_{\mathcal{F}_0}^*}  \right) @>-\delta >> \cdots
\end{CD}
\end{equation}}}

For $\frac{\partial}{\partial t}=\sum_{\eta=1}^n c_\eta \frac{\partial}{\partial t^\eta} \in T_t B$ with $c_\eta \in \mathbb{C}$, we define an element $\gamma_j(t)\in \Gamma\left( U_j^t, \mathscr{H}om_{\mathcal{O}_{M_t} } \left( \mathcal{N}_{\mathcal{F}_t}^*, \mathcal{O}_{M_t} \right) \right)$ in the following way:
\begin{align*}
\gamma_j(t):\Gamma\left( \mathcal{U}_j^t,  \mathcal{N}_{\mathcal{F}_t}^* \right)&\to \Gamma\left( U_j^t, \mathcal{O}_{M_t} \right)\\
w_j^\alpha(\xi_j,t) &\mapsto -i_{\frac{\partial}{\partial t}}\left(v_j^\alpha(\xi_j,t) \right)= -i_{\frac{\partial}{\partial t}}\left(\sum_{\eta=1}^m v_j^{\alpha\eta}(\xi_j,t_j)dt_j^\eta \right)= -\sum_{\eta=1}^m c_\eta v_j^{\alpha\eta}(\xi_j,t),\,\,\,\,\,\,\alpha=1,...,q
\end{align*}
and linearly extends to $\Gamma\left( U_{jk}^t, \mathcal{O}_{M_t} \right)$. Then we claim that $\left( \left\{ \beta_j(t) \right\}, \left\{ \gamma_j(t)\right\}, \left\{\left( \theta_{jk}(t), B_{jk}(t) \right\}\right)    \right)$ defines a $1$-cocycle in the above \v Cech resolution $(\ref{ou3})$ of $\tilde{\mathcal{E}}_{\mathcal{N}_{\mathcal{F}_0}^*}^\bullet$. From subsection $(\ref{oc8})$, it remains to show $G_1'^t\left(\beta_j(t), \gamma_j(t) \right)=0$ which follows from $(\ref{u8})$ (for the detail, see Part III), and $\gamma_k(t)\left( w_{j}^\alpha\right)- \gamma_j(t)\left( w_j^\alpha \right) +i_{\theta_{jk}(t)}\left( w_j^\alpha \right)=0,\alpha=1,...,q$ which follows from $(\ref{u6})$. Hence we have the infinitesimal unfolding along $\frac{\partial}{\partial t}$ at $t$:
\begin{align}\label{ou6}
\frac{\partial \left( M_t, \tilde{\mathcal{N}}_{\mathcal{F}_t}^* \right) }{\partial t} = \left( \left\{ \beta_j(t) \right\}, \left\{ \gamma_j(t) \right\}, \left\{ \left( \theta_{jk}(t), A_{jk}(t) \right) \right\} \right) \in \mathbb{H}^1\left( M_t, \tilde{\mathcal{E}}_{\mathcal{N}_{\mathcal{F}_t}^*}^\bullet \right)
\end{align}

\subsection{Description of obstructed elements in \v Cech resolution of $\tilde{\mathcal{E}}_{\mathcal{N}_{\mathcal{F}_0}^*}^\bullet$}

\begin{definition}
We say that a compact foliated complex manifold $\left(M, \mathcal{N}_{\mathcal{F}_0}^* \right)$ with $\mathcal{N}_{\mathcal{F}_0}^*$ locally free is unobstructed in unfoldings if there is an unfolded analytic family $\left( \mathcal{M}, \tilde{\mathcal{N}}_\mathcal{F}^*, B, \pi \right)$ of unfoldings of $\pi^{-1}(0)=\left( M, \mathcal{N}_{\mathcal{F}_0}^* \right), 0\in B$ such that the unfolded Kodaira-Spencer map $\tilde{\varphi}_0:T_0B\to \mathbb{H}^1\left( M, \tilde{\mathcal{E}}_{\mathcal{N}_{\mathcal{F}_0}^*}^\bullet \right)$ is an isomorphism at $0 \in B$. Otherwise, we say that $\left( M, \mathcal{N}_{\mathcal{F}_0}^* \right)$ is obstructed in unfoldings.
\end{definition}

\begin{remark}
Since we have not established the complete form of theorem of existence for unfoldings, and only established theorem of existence for unfoldings for one-parameter families $($see \textnormal{Part III}$)$, we cannot say that $\mathbb{H}^2\left( M, \tilde{\mathcal{E}}_{\mathcal{N}_{\mathcal{F}_0}^*}^\bullet \right)=0$ implies unobstructedness of unfoldings except for the case $\dim_\mathbb{C} \mathbb{H}^1\left( M, \tilde{\mathcal{E}}_{\mathcal{N}_{\mathcal{F}_0}^*}^\bullet \right)=1$. We leave the complete form of theorem of existence for unfoldings as a conjecture. However we still have the notion of obstructed elements in unfoldings as follows.
\end{remark}

\begin{definition}
Let $\left(M, \mathcal{N}_{\mathcal{F}_0}^*\right)$ be a compact foliated complex manifold with $\mathcal{N}_{\mathcal{F}_0}^*$ locally free. We say that $a\in \mathbb{H}^1\left( M, \tilde{\mathcal{E}}_{\mathcal{N}_{\mathcal{F}_0}^*}^\bullet \right)$ is an obstructed element in unfoldings of $\left(M, \mathcal{N}_{\mathcal{F}_0}^* \right)$  if there is no unfolded analytic family $\left( \mathcal{M}, \tilde{\mathcal{N}}_\mathcal{F}^*, B, \pi \right)$ of unfoldings of $\left(M, \mathcal{N}_{\mathcal{F}_0}^* \right)= \pi^{-1}(0), 0\in B$ such that $a\in \tilde{\varphi}_0\left(T_0 M \right)$ where $\tilde{\varphi}_0: T_0M  \to \mathbb{H}^1\left( M, \tilde{\mathcal{E}}_{\mathcal{N}_{\mathcal{F}_0}^*}^\bullet \right)$ is the unfolded Kodaira-Spencer map at $0\in B$ of $\left(\mathcal{M}, \tilde{\mathcal{N}}_\mathcal{F}^* , B, \pi \right)$.
\end{definition}

Hence if there is an obstructed element $a\in \mathbb{H}^1\left( M, \tilde{\mathcal{E}}_{\mathcal{N}_{\mathcal{F}_0}^*}^\bullet \right)$ in unfoldings, $\left( M, \mathcal{N}_{\mathcal{F}_0}^* \right)$ is obstructed in unfoldings. We now describe a condition when $a\in \mathbb{H}^1\left( M, \tilde{\mathcal{E}}_{\mathcal{N}_{\mathcal{F}_0}^*}^\bullet \right)$ is an obstructed element in unfoldings. We extend the arguments in \cite{Kod05} p.210-214 in the context of an unfolded analytic family, and we keep the notations in subsection \ref{oc9} for the underlying foliated analytic family $\left( \mathcal{M}, \mathcal{N}_\mathcal{F}^*, B, \pi \right)$ induced from $\left( \mathcal{M}, \tilde{\mathcal{N}}_\mathcal{F}^*, B, \pi \right)$.

Let $\left( \mathcal{M}, \tilde{\mathcal{N}}_\mathcal{F}^*, B, \pi \right)$ with $0\in B \subset \mathbb{C}$ and $\tilde{\mathcal{N}}_{\mathcal{F}}^*$ locally free be an unfolded analytic family such that $\pi^{-1}(0)=\left(M, \mathcal{N}_{\mathcal{F}_0}^* \right)$. Take a small disk $\Delta$ with centre $0$ such that $0\in \Delta \subset B$, and represent $\left(\mathcal{M}_\Delta, \tilde{\mathcal{N}}_{\mathcal{F}_\Delta}^* \right)=\pi^{-1}(\Delta)$ in the form
\begin{align*}
\left(\mathcal{M}_\Delta, \tilde{\mathcal{N}}_{\mathcal{F}_\Delta}^* \right) = \bigcup_{j=1}^l \left( U_j\times \Delta, \tilde{\mathcal{N}}_{\mathcal{F}_\Delta}^* \right)
\end{align*}
where each $U_j$ is a polydisk, and $\tilde{\mathcal{N}}_{\mathcal{F}_\Delta}^*$ on $U_j\times \Delta$ is generated by $\omega_j^\alpha(\xi_j,t):=w_j^\alpha(\xi_j,t)+v_j^\alpha(\xi_j,t)=\sum_{\beta=1}^n w_j^{\alpha\beta}(\xi_j,t_j)d\xi_j^\beta+ v_j^{\alpha}(\xi_j,t_j) dt_j ,\alpha=1,...,q$ such that $\omega_j^1(z_j,t)\wedge \cdots \wedge \omega_j^q(z_j,t)\wedge d\omega_j^\alpha(z_j,t)=0,\alpha=1,...,q$, and $(\xi_j,t)\in U_j\times \Delta$ and $(\xi_k,t)\in U_k\times \Delta$ are the same points on $\mathcal{M}_\Delta$ if $\xi_j^\alpha= f_{jk}^\alpha(\xi_k,t)$ for $\alpha=1,...,n$. Here each $f_{jk}^\alpha(\xi_k,t)=f_{jk}^\alpha(\xi_k^1,..., \xi_k^n, t)$ is a holomorphic function of $\xi_k^1,..., \xi_k^n,t$ defined on $\left(U_k\times \Delta \right) \cap \left( U_j\times \Delta \right)\ne \emptyset$. On the other hand, $\omega_j^\alpha(\xi_j,t)= \sum_{\beta=1}^q h_{jk}^{\alpha\beta}\left(\xi_j,t \right) \omega_k^\beta(\xi_k,t),\alpha=1,...,q$ on $\left(U_k\times \Delta \right) \cap \left( U_j\times \Delta \right)\ne \emptyset$, and $h_{ik}^{\alpha\gamma}\left(\xi_k,t \right)= \sum_{\beta=1}^p h_{ij}^{\alpha\beta}(\xi_j,t) h_{jk}^{\beta \gamma}(\xi_k,t)$ on $\left( U_i\times \Delta \right) \cap \left( U_j\times \Delta \right) \cap \left( U_k\times \Delta \right)\ne \emptyset$.

As in $(\ref{ou6})$, the infinitesimal unfolding $\frac{\partial \left( M_t , \tilde{\mathcal{N}}_{\mathcal{F}_t}^* \right)}{\partial t}\in \mathbb{H}^1\left( M_t, \tilde{\mathcal{E}}_{\mathcal{N}_{\mathcal{F}_t}^*}^\bullet \right)$ is the cohomology class of the $1$-cocycle $\left( \left\{ \beta_j (t)\right\}, \left\{ \gamma_j(t) \right\}, \left\{ \left(\theta_{jk}(t), B_{jk}(t) \right) \right\}  \right)\in C^0\left( \mathcal{U}_t, \mathscr{H}om_{\mathcal{O}_{M_t}} \left( \mathcal{N}_{\mathcal{F}_t}^*, \Omega_{M_t}^1 \right) \right)\bigoplus C^0\left(  \mathcal{U}_t, \mathscr{H}om_{\mathcal{O}_{M_t}}\left( \mathcal{N}_{\mathcal{F}_t}^*, \mathcal{O}_{M_t} \right) \right)$ $\bigoplus C^1\left(\mathcal{U}_t, \mathcal{E}_{\mathcal{N}_{\mathcal{F}_t}^*} \right)$, where $\mathcal{U}_t=\left\{U_j^t:= U_j\times t \right\}$ such that
\begin{align*}
\beta_j(t):\Gamma\left( U_j^t, \mathcal{N}_{\mathcal{F}_t}^* \right)&\to \Gamma\left( U_j^t, \Omega_{M_t}^1 \right) \\
 w_j^\alpha(\xi_j,t) &\mapsto  \beta_j^\alpha(\xi_j,t):=\sum_{\beta=1}^n \beta_j^{\alpha\beta}(\xi_j,t)d\xi_j^\beta ,\,\,\,\,\,\,\,\,\textnormal{where}\,\,\,\beta_j^{\alpha\beta}(\xi_j, t)= - \frac{\partial w_j^{\alpha\beta}(\xi_j,t)}{\partial t},\alpha=1,...,q
\end{align*}
and
\begin{align*}
\gamma_j(t): \Gamma\left( U_j^t, \mathcal{N}_{\mathcal{F}_t}^* \right) &\to \Gamma\left( U_j^t, \mathcal{O}_{M_t} \right)\\
 w_j^\alpha(\xi_j,t) &\mapsto v_j^\alpha\left(\xi_j,t \right) ,\,\,\,\,\,\,\,\alpha=1,...,q
\end{align*}
and the vector field
\begin{align*}
\theta_{jk}(t)=\sum_{\alpha=1}^n \theta_{jk}^\alpha(\xi_j,t)\frac{\partial}{\partial \xi_j^\alpha},\,\,\,\,\,\,\,\,\textnormal{where}\,\,\,\,\theta_{jk}^\alpha(\xi_j,t)=\frac{\partial f_{jk}(\xi_k,t)}{\partial t},\,\,\,\xi_k=f_{kj}(\xi_j,t)
\end{align*}
and
\begin{align*}
B_{jk}(t):\Gamma\left( U_j^t, \mathcal{N}_{\mathcal{F}_t}^* \right)&\to \Gamma\left( U_j^t, \mathcal{N}_{\mathcal{F}_t}^* \right) \\
     w_j^\alpha(\xi_j,t)&\mapsto \sum_{\gamma=1}^p B_{jk}^{\alpha\gamma}(\xi_j,t)w_j^\gamma(\xi_j,t),\,\,\textnormal{where}\,\,\, B_{jk}^{\alpha\gamma}(\xi_j,t)=\sum_{\beta=1}^q\frac{\partial h_{jk}^{\alpha\beta}(\xi_k,t)}{\partial t}h_{kj}^{\beta \gamma}(\xi_j,t),\xi_k=f_{kj}(\xi_j,t),\alpha=1,...,q
\end{align*}

On $\left(U_k\times \Delta \right)\cap \left(U_j\times \Delta \right) \ne \emptyset$, we have the equalities
\begin{align}
f_{ik}^\alpha(\xi_k,t)&=f_{ij}^\alpha\left(f_{jk}(\xi_k,t),t \right),\,\,\,\alpha=1,...,n \label{ou13}\\
\sum_{\gamma=1}^n w_j^{\alpha\gamma}(f_{jk}(\xi_k,t),t)\frac{\partial f_{jk}^\gamma(\xi_k,t)}{\partial \xi_k^\eta}&=\sum_{\beta=1}^qh_{jk}^{\alpha\beta}(\xi_k,t)w_k^{\beta\eta}(\xi_k,t),\,\,\,\alpha=1,...,q, \eta=1,...,n \label{ou10} \\
\sum_{\gamma=1}^nw_j^{\alpha\gamma}\left(f_{jk}(\xi_k,t),t\right)\frac{\partial f_{jk}^\gamma(\xi_k,t)}{\partial t} & +  v_j^{\alpha}\left(f_{jk}(\xi_k,t),t \right)  = \sum_{\beta=1}^q h_{jk}^{\alpha\beta}(\xi_k,t) v_k^{\beta }(\xi_k,t) \label{ou7}
\end{align}
and on $\left(U_i\cap \Delta \right)\cap \left(  U_j \times \Delta \right)\cap \left(U_k\times \Delta \right)\ne \emptyset$, we have the equalities
\begin{align}
h_{ik}^{\alpha\gamma}\left(\xi_k,t \right) &= \sum_{\beta=1}^q h_{ij}^{\alpha\beta}\left(f_{jk}(\xi_k,t),t \right) h_{jk}^{\beta \gamma}(\xi_k,t),\,\,\,\,\,\alpha,\gamma=1,...,q. \label{ou11}
\end{align}
and on $U_j\times \Delta$, we have the equalities 
\begin{align}
w_j^1(\xi_j,t)\wedge \cdots \wedge w_j^q(\xi_j,t)\wedge d_t w_j^\alpha(\xi_j,t)=0,\,\,\,\,\,\,\alpha=1,...,q \label{ou12}
\end{align}
where $d_t$ is the differential on $M_t$ (by abuse of notation we shall denote $d_t$ by $d$ in the following), and
\begin{align}\label{ou8}
\sum_{\gamma=1}^q (-1)^{q-\gamma} v_j^\gamma(\xi_j,t) w_j^1(\xi_j,t) \wedge \cdots \wedge &\widehat{w}_j^\gamma(\xi_j,t) \wedge \cdots \wedge w_j^q(\xi_j,t) \wedge dw_j^\alpha(\xi_j,t) + w_j^1(\xi_j,t)\wedge \cdots \wedge w_j^q(\xi_j,t) \wedge dv_j^\alpha(\xi_j,t) \\
&- w_j^1(\xi_j,t)\wedge \cdots \wedge w_j^q(\xi_j,t) \wedge \left(\frac{\partial w_j^\alpha(\xi_j,t)}{\partial t} \right) ,\,\,\,\,\,\alpha=1,...,q \notag
\end{align}

By differentiating both sides of $(\ref{ou13}),(\ref{ou10}), (\ref{ou7}), (\ref{ou11})$, $(\ref{ou12})$ in $t$, we have the equalities
{\small{\begin{align}
&\theta_{ik}^\alpha(\xi_i,t)=\theta_{ij}^\alpha(\xi_i,t) + \sum_{\beta=1}^n \frac{\partial \xi_i^\alpha}{\partial \xi_j^\beta}\theta_{jk}^\beta(\xi_j,t),\,\,\,\alpha=1,...,n. \label{ou14}\\
&\sum_{\gamma,\delta=1}^n \frac{\partial w_j^{\alpha\gamma}}{\partial \xi_j^\delta}\left(f_{jk}(\xi_k,t),t) \right)\frac{\partial f_{jk}^\delta}{\partial t}\frac{\partial f_{jk}^\gamma}{\partial \xi_k^\eta}+\sum_{\gamma=1}^n \frac{\partial w_j^{\alpha\gamma}}{\partial t}\left(f_{jk}(\xi_k,t),t \right)\frac{\partial f_{jk}^\gamma}{\partial \xi_k^\eta} +\sum_{\gamma=1}^n w_j^{\alpha\gamma}\left(f_{jk}(\xi_k,t),t \right)\frac{\partial}{\partial \xi_k^\eta}\left(\frac{\partial f_{jk}^\gamma}{\partial t} \right) \label{ou15}\\
&= \sum_{\beta=1}^q \frac{\partial h_{jk}^{\alpha\beta}}{\partial t}w_k^{\beta\eta} +\sum_{\beta=1}^q h_{jk}^{\alpha\beta}\frac{\partial w_k^{\beta\eta}}{\partial t} ,\,\,\,\,\,\,\,\,\alpha=1,...,q, \eta=1,...,n\notag
\end{align}}}
and for $\alpha,\gamma=1,...,q$,
{\small{\begin{align}
\frac{\partial h_{ik}^{\alpha \gamma}}{\partial t} = \sum_{\beta=1}^q \sum_{\eta=1}^n \frac{\partial h_{ij}^{\alpha\beta}}{\partial \xi_j^\eta}\left(f_{jk}(\xi_k, t) , t\right) \frac{\partial f_{jk}^{\eta}}{\partial t} h_{jk}^{\beta \gamma} + \sum_{\beta=1}^q \frac{\partial h_{ij}^{\alpha\beta}}{\partial t}\left(f_{jk}(\xi_k,t),t\right) h_{jk}^{\beta \gamma} + \sum_{\beta=1}^q h_{ij}^{\alpha\beta}\left(f_{jk}(\xi_k,t),t\right)\frac{\partial h_{jk}^{\beta \gamma}}{\partial t} \label{ou16}
\end{align}}}
\begin{align}
w_j^1\wedge \cdots \wedge w_j^q \wedge d\left(\frac{\partial w_j^\alpha}{\partial t} \right)+\sum_{\gamma=1}^q w_j^1\wedge \cdots \wedge \overbrace{\frac{\partial w_j^\gamma}{\partial t}}^{\gamma-\textnormal{th}} \wedge \cdots \wedge w_j^q\wedge dw_j^\alpha=0 ,\,\,\,\alpha=1,...,q\label{ou17}
\end{align}

Next we will take the derivative of $(\ref{ou14}), (\ref{ou17})$ in $t$ as functions of $\xi_j^1,..., \xi_j^n,t$ and $(\ref{ou15}), (\ref{ou16})$ in $t$ as functions of $\xi_k^1,..., \xi_k^n,t$ again, and take the derivative of $(\ref{ou7})$ as functions of $\xi_k^1,..., \xi_k^n, t$ and $(\ref{ou8})$ as functions of $\xi_j^1,..., \xi_j^n, t$.  We recall $(\ref{o53})$ in subsection \ref{o42}.

By taking the derivative of $(\ref{ou14})$ with respect to $t$, and putting $\dot{\theta}_{ij}(t)= \sum_{\alpha=1}^n \dot{\theta}_{ij}^\alpha (\xi_i,t)\frac{\partial}{\partial \xi_i^\alpha}$, where $\dot{\theta}_{ij}^\alpha(\xi_i,t)=\frac{\partial \theta_{ij}^\alpha(z_i,t)}{\partial t}$, we obtain (for the detail, see \cite{Kod05} p.211-213)
\begin{align}\label{ou18}
\dot{\theta}_{ij}(t) -\dot{\theta}_{ik}(t) + \dot{\theta}_{jk} (t) = \left[ \theta_{ij}(t), \theta_{jk}(t) \right]
\end{align}
Let us define $\dot{\beta}_j(t) \in \Gamma\left( U_j^t, \mathscr{H}om_{\mathcal{O}_{M_t}}\left( \mathcal{N}_{\mathcal{F}_t}^* , \Omega_{M_t}^1\right) \right)$ by 
\begin{align*}
\dot{\beta}_j(t):\Gamma\left(U_j^t, \mathcal{N}_{\mathcal{F}_t}^* \right) &\to \Gamma\left( U_j^t, \Omega_{M_t}^1 \right) \\
 w_j^\alpha &\mapsto \dot{\beta}_j^\alpha(t):= \dot{\beta}_j^\alpha(\xi_j,t)= \frac{\partial \beta_j^{\alpha}(\xi_j,t)}{\partial t}= \sum_{\beta=1}^n \frac{\partial \beta_j^{\alpha\beta}(\xi_j,t)}{\partial t} d\xi_j^\beta = - \sum_{\beta=1}^n \frac{\partial^2 w_j^{\alpha\beta}(\xi_j,t)}{\partial t^2}d\xi_j^\beta, \,\,\,\alpha=1,...,q
\end{align*}
and linearly extend to $\Gamma\left( U_j^t, \mathcal{N}_{\mathcal{F}_t}^* \right)$. Then by taking the derivative of $(\ref{ou17})$, we obtain 
\begin{align}\label{ou30}
E_1'\left( \dot{\beta}_j(t)\right)=\left[\beta_j(t), \beta_j(t) \right]_{\mathcal{N}_{\mathcal{F}_t}^*}
\end{align}
as in subsection \ref{oc9}. Let us define
\begin{align*}
\dot{B}_{jk}(t):\Gamma\left( U_{jk}^t,  \mathcal{N}_{\mathcal{F}_t}^* \right)\to &\Gamma\left( U_{jk}^t, \mathcal{N}_{\mathcal{F}_t}^* \right) \\
 w_j^\alpha \mapsto &\sum_{\gamma=1}^q \frac{\partial}{\partial t_j}\left( B_{jk}^{\alpha \gamma}(\xi_j,t)\right) w_j^\gamma = \sum_{\beta, \gamma=1}^q \frac{\partial}{\partial t_j}\left( \frac{\partial h_{jk}^{\alpha \beta}(\xi_k,t)}{\partial t} h_{kj}^{\beta \gamma}(\xi_j,t) \right) w_j^\gamma
\end{align*}
and linearly extends to $\Gamma\left( U_{jk}^t, \mathcal{N}_{\mathcal{F}_t}^* \right)$. Then $\left\{\left( \dot{\theta}_{jk}(t), \dot{B}_{jk}(t) \right) \right\} \in C^1\left(\mathcal{U}_t, \mathcal{E}_{ \mathcal{N}_{\mathcal{F}_t}^*} \right)$ as in subsection \ref{oc9}. By taking the derivative of $(\ref{ou16})$ with respect to $t_k$, we obtain as in subsection \ref{oc9}
\begin{align} \label{ou28}
\dot{B}_{ij} (t) - \dot{B}_{ik} (t) + \dot{B}_{jk}^i  (t) = \left[ \theta_{ij} (t) , B_{jk}^i(t)\right] - \left[ \theta_{jk}(t), B_{ij}(t)\right] + \left[ B_{ij} (t) ,   B_{jk}^i(t)\right]
\end{align}
and by taking the derivative of $(\ref{ou15})$ with respect to $t_k$, we have from Lemma \ref{oc24},
\begin{align}
&\dot{\beta}_k(t)-\dot{\beta}_j(t)+ E_0'\left( \left( \dot{\theta}_{jk} (t) , \dot{B}_{jk} (t) \right) \right) \label{ou24} \\
&= - \beta_j(t)\circ B_{jk}(t) - \beta_k (t) \circ B_{jk}(t)+ \mathcal{L}_{\theta_{jk}(t)}\circ \beta_j(t)+ \mathcal{L}_{\theta_{jk}(t)}\circ \beta_k(t) \notag
\end{align}

On the other hand, let us define $\dot{\gamma}_j(t) \in \Gamma\left( U_j^t,  \mathscr{H}om_{\mathcal{O}_{M_t} }\left( \mathcal{N}_{\mathcal{F}_t}^*, \mathcal{O}_{M_t} \right) \right)$ by
\begin{align*}
\dot{\gamma}_j(t) : \Gamma\left( U_j^t \mathcal{N}_{\mathcal{F}_t}^* \right) &\to \Gamma\left(  U_j^t, \mathcal{O}_{M_t} \right) \\
 w_j^\alpha &\mapsto \dot{\gamma}_j^\alpha(t):=\dot{\gamma}_j^\alpha(\xi_j, t)=\frac{\partial \gamma_j^\alpha(\xi_j,t)}{\partial t} = - \frac{\partial v_j^\alpha(\xi_j,t)}{\partial t}
\end{align*}
and linearly extend to $\Gamma\left( U_j^t,  \mathcal{N}_{\mathcal{F}_t}^* \right)$. Let us take the derivative of $(\ref{ou8})$ with respect to $t_j$, we have
{\small{\begin{align*}
&0= \sum_{\gamma=1}^q (-1)^{q-\gamma} \frac{\partial v_j^{\gamma}}{\partial t} w_j^1\wedge \cdots \wedge \widehat{w}_j^\gamma \wedge \cdots \wedge w_j^q \wedge d w_j^\alpha+\sum_{\gamma=1}^q\sum_{1\leq \beta <\gamma } (-1)^{q-\gamma} v_j^\gamma w_j^1\wedge \cdots \wedge \overbrace{\frac{\partial w_j^\beta }{\partial t} }^{\beta-\textnormal{th}} \wedge \cdots \wedge \hat{w}_j^\gamma \wedge \cdots \wedge w_j^q\wedge dw_j^\alpha \\
&+\sum_{\gamma=1}^q\sum_{\gamma<\beta \leq q} (-1)^{q-\gamma} v_j^\gamma w_j^1\wedge \cdots \wedge \hat{w}_j^\gamma \wedge \cdots \wedge \overbrace{ \frac{\partial w_j^\beta}{\partial t} }^{\beta - \textnormal{th}} \wedge \cdots \wedge w_j^q\wedge dw_j^\alpha  + \sum_{\gamma=1}^q (-1)^{q-\gamma} v_j^\gamma w_j^1\wedge \cdots \wedge \hat{w}_j^\gamma \wedge \cdots w_j^q \wedge d\left( \frac{\partial w_j^\alpha}{\partial t}\right) \\
& +\sum_{\gamma=1}^q w_j^1\wedge \cdots \wedge \overbrace{\frac{\partial w_j^\gamma}{\partial t}}^{\gamma-\textnormal{th}}\wedge \cdots \wedge w_j^q\wedge dv_j^\alpha + w_j^1\wedge \cdots \wedge w_j^q \wedge d\left( \frac{\partial v_j^\alpha}{\partial t} \right)\\
&-\sum_{\gamma-1}^q   w_j^1\wedge \cdots \wedge \overbrace{ \frac{\partial w_j^\gamma}{\partial t} }^{\gamma-\textnormal{th}} \wedge \cdots \wedge w_j^q\wedge \frac{\partial w_j^\alpha}{\partial t} - w_j^1\wedge \cdots \wedge w_j^q \wedge \frac{\partial^2 w_j^\alpha}{\partial t^2}
\end{align*}}}

Then we have
\begin{align} \label{ou20}
&\sum_{\gamma=1}^q(-1)^{q-\gamma} \dot{\gamma}_j^\gamma (t) w_j^1\wedge \cdots \wedge \hat{w}_j^\gamma \wedge \cdots \wedge w_j^q\wedge dw_j^\alpha +w_j^1\wedge \cdots \wedge w_j^q\wedge d\left( \dot{\gamma}_j^\alpha (t) \right) -w_j^1\wedge \cdots \wedge w_j^q \wedge \dot{\beta}_j^\alpha(t)\\
&=\sum_{\gamma=1}^q\sum_{1\leq \beta \ne \gamma \leq q } (-1)^{q-\gamma} \gamma_j^\gamma(t) w_j^1\wedge \cdots \wedge \overbrace{\beta_j^\beta(t) }^{\beta-\textnormal{th}} \wedge \cdots \wedge \hat{w}_j^\gamma \wedge \cdots \wedge w_j^q\wedge dw_j^\alpha \notag \\
&+\sum_{\gamma=1}^q (-1)^{q-\gamma} \gamma_j^\gamma(t) w_j^1\wedge \cdots \wedge \hat{w}_j^\gamma \wedge \cdots w_j^q \wedge d\left( \beta_j^\alpha(t)\right) + \sum_{\gamma=1}^q w_j^1\wedge \cdots \wedge \overbrace{\beta_j^\gamma(t)}^{\gamma-\textnormal{th}}\wedge \cdots \wedge w_j^q\wedge d\gamma_j^\alpha(t) \notag\\
& -  \sum_{\gamma=1}^q w_j^1\wedge \cdots \wedge  \overbrace{\beta_j^\gamma(t)}^{\gamma-\textnormal{th}} \wedge \cdots \wedge w_j^q\wedge \beta_j^\alpha(t) \notag
\end{align} 

Let us take the derivative of $(\ref{ou7})$ with respect to $t_k$:
\begin{align}\label{ou21}
\sum_{\gamma, \delta=1}^n \frac{\partial w_j^{\alpha\gamma}}{\partial \xi_j^\delta}\frac{\partial f_{jk}^\delta}{\partial t}\frac{\partial f_{jk}^\gamma}{\partial t}+\sum_{\gamma=1}^n \frac{\partial w_j^{\alpha\gamma}}{\partial t_j}\frac{\partial f_{jk}^\gamma}{\partial t} +\sum_{\gamma=1}^n w_j^{\alpha\gamma} \frac{\partial^2 f_{jk}^\gamma}{\partial t^2} +\sum_{\gamma=1}^n \frac{\partial v_j^{\alpha}}{\partial \xi_j^\gamma}\frac{\partial f_{jk}^\gamma}{\partial t} +\frac{\partial v_j^{\alpha}}{\partial t_j}=\sum_{\beta=1}^q\frac{\partial h_{jk}^{\alpha\beta}}{\partial t}v_k^{\beta}+\sum_{\beta=1}^q h_{jk}^{\alpha\beta} \frac{\partial v_k^{\beta}}{\partial t}
\end{align}

By noting $\dot{\theta}_{jk}^\gamma=\frac{\partial^2 f_{jk}^\gamma}{\partial t^2}- \sum_{\delta=1}^n \frac{\partial f_{jk}^{\delta}}{\partial t} \frac{\partial}{\partial \xi_j^\delta}\left(  \frac{\partial f_{jk}^\gamma}{\partial t} \right)$, $(\ref{ou21})$ is equivalent to 
{\Small{\begin{align}\label{ou22}
\sum_{\gamma, \delta=1}^n \frac{\partial w_j^{\alpha\gamma}}{\partial \xi_j^\delta}\frac{\partial f_{jk}^\delta}{\partial t}\frac{\partial f_{jk}^\gamma}{\partial t}+\sum_{\gamma=1}^n \frac{\partial w_j^{\alpha\gamma}}{\partial t_j}\frac{\partial f_{jk}^\gamma}{\partial t} +\sum_{\gamma=1}^n w_j^{\alpha\gamma} \left( \dot{\theta}_{jk}^\gamma + \sum_{\delta=1}^n \frac{\partial f_{jk}^\delta}{\partial t} \frac{\partial}{\partial \xi_j^\delta} \left(  \frac{\partial f_{jk}^\gamma}{\partial t} \right)   \right) +\sum_{\gamma=1}^n \frac{\partial v_j^{\alpha}}{\partial \xi_j^\gamma}\frac{\partial f_{jk}^\gamma}{\partial t} +\frac{\partial v_j^{\alpha}}{\partial t_j}=\sum_{\beta=1}^q\frac{\partial h_{jk}^{\alpha\beta}}{\partial t}v_k^{\beta}+\sum_{\beta=1}^q h_{jk}^{\alpha\beta} \frac{\partial v_k^{\beta}}{\partial t}
\end{align}}}
We note that
{\small{\begin{align}\label{ou23}
i_{\theta_{jk}}\mathcal{L}_{\theta_{jk}}\left(w_j^\alpha \right) & = i_{\sum_{\sigma=1}^n \frac{\partial f_{jk}^\sigma}{\partial t} \frac{\partial}{\partial \xi_j^\sigma}} \mathcal{L}_{\sum_{\delta=1}^n \frac{\partial f_{jk}^\delta}{\partial t} \frac{\partial }{\partial \xi_j^\delta}} \left(\sum_{\gamma=1}^n w_j^{\alpha \gamma} d\xi_j^\gamma \right)=  i_{\sum_{\sigma=1}^n \frac{\partial f_{jk}^\sigma}{\partial t} \frac{\partial}{\partial \xi_j^\sigma}} \left(  \sum_{\delta , \gamma=1}^n \frac{\partial f_{jk}^\delta}{\partial t}     \frac{\partial w_j^{\alpha \gamma}}{\partial \xi_j^\delta} d\xi_j^\gamma + w_j^{\alpha \gamma}\frac{\partial}{\partial \xi_j^\delta}\left( \frac{\partial f_{jk}^\gamma}{\partial t}\right) d\xi_j^\delta   \right)\\
&= \sum_{\delta, \gamma=1}^n \frac{\partial f_{jk}^\gamma}{\partial t} \frac{\partial f_{jk}^\delta}{\partial t} \frac{\partial w_j^{\alpha \gamma}}{\partial \xi_j^\delta} + \sum_{\gamma, \delta=1}^n w_j^{\alpha \gamma} \frac{\partial f_{jk}^\delta}{\partial t}\frac{\partial}{\partial \xi_j^\delta}\left(   \frac{\partial f_{jk}^\gamma}{\partial t}  \right)     \notag
\end{align}}}
The $(\ref{ou23})$ implies that $(\ref{ou22})$ is equivalent to
\begin{align}\label{ou25}
i_{\theta_{jk}}\mathcal{L}_{\theta_{jk}}(w_j^\alpha) + i_{\theta_{jk}}\left(\frac{\partial w_j^\alpha}{\partial t_j} \right) + i_{\dot{\theta}_{jk}}\left( w_j^\alpha \right) + \mathcal{L}_{\theta_{jk}}\left( v_j^\alpha \right) + \frac{\partial v_j^\alpha}{\partial t}= \sum_{\beta=1}^q \frac{\partial h_{jk}^{\alpha\beta}}{\partial t_k} v_k^\beta + \sum_{\beta=1}^q h_{jk}^{\alpha\beta} \frac{\partial v_k^\beta}{\partial t}
\end{align}

\begin{lemma}
We have the equality
\begin{align} \label{ou26}
&\dot{\gamma}_k(t)- \dot{\gamma}_j(t) + i_{\dot{\theta}_{jk}} \\
&=  -\frac{1}{2}\gamma_j(t)\circ B_{jk}(t) - \frac{1}{2}\gamma_k(t) \circ B_{jk}(t) + \frac{1}{2}\mathcal{L}_{\theta_{jk}(t)}\circ \gamma_j(t) + \frac{1}{2}\mathcal{L}_{\theta_{jk}(t)}\circ \gamma_k(t) + \frac{1}{2} i_{\theta_{jk}(t)}\circ \beta_j(t) + \frac{1}{2} i_{\theta_{jk}(t)} \circ \beta_k(t)  \notag
\end{align}
\end{lemma}

\begin{proof}
Since $\beta_k(t)- \beta_j(t) +\left(\mathcal{L}_{\theta_{jk}(t)}- B_{jk}(t) \right)=0$ and $\gamma_k(t)- \gamma_j(t)+i_{\theta_{jk}(t)}=0$, we have
{\small{\begin{align*}
 & \left( -\frac{1}{2}\gamma_j(t)\circ B_{jk}(t) - \frac{1}{2}\gamma_k(t) \circ B_{jk}(t) + \frac{1}{2}\mathcal{L}_{\theta_{jk}(t)}\circ \gamma_j(t) + \frac{1}{2}\mathcal{L}_{\theta_{jk}(t)}\circ \gamma_k(t) + \frac{1}{2} i_{\theta_{jk}(t)}\circ \beta_j(t) + \frac{1}{2} i_{\theta_{jk}(t)} \circ \beta_k(t)\right)\left(w_{j}^\alpha\right) \\
&= \left(- \gamma_k(t) \circ B_{jk}(t) - \frac{1}{2} i_{\theta_{jk}}\circ B_{jk}(t)  + \mathcal{L}_{\theta_{jk}} \circ \gamma_j(t) - \frac{1}{2}\mathcal{L}_{\theta_{jk}(t)} \circ i_{\theta_{jk}(t)} + i_{\theta_{jk}(t) } \circ \beta_j(t) - \frac{1}{2} i_{\theta_{jk}(t)} \circ \left( \mathcal{L}_{\theta_{jk}(t)} - B_{jk}(t) \right)\right)\left(w_j^\alpha \right) \\
&= \left( - \gamma_k(t) \circ B_{jk}(t) + \mathcal{L}_{\theta_{jk}(t)} \circ \gamma_j(t) - i_{\theta_{jk}(t)}\circ \mathcal{L}_{\theta_{jk}(t)} + i_{\theta_{jk}(t)} \circ \beta_j(t) \right)\left(w_j^\alpha \right) \\
&=\sum_{\beta=1}^q \frac{\partial h_{ij}^{\alpha\beta}}{\partial t} v_k^\beta - \mathcal{L}_{\theta_{jk}(t)} \left( v_j^\alpha \right) - i_{\theta_{jk}(t)}\mathcal{L}_{\theta_{jk}(t)} \left( w_j^\alpha \right) - i_{\theta_{jk}(t)}\left( \frac{\partial w_j^\alpha}{\partial t} \right)
\end{align*}}} 
Then $(\ref{ou25})$ implies $(\ref{ou26})$.

\end{proof}

From $(\ref{ou18}), (\ref{ou28}), (\ref{ou30}), (\ref{ou24}), (\ref{ou20})$ and $(\ref{ou26})$, substituting $t_0$, we obtain
\begin{align}\label{ou34}
\dot{\theta}_{ij}(0)- \dot{\theta}_{ik}(0) + \dot{\theta}_{jk}(0)= \left[ \theta_{ij}(0),\theta_{jk}(0) \right]
\end{align}
\begin{align} \label{ou31}
E_1'\left(\dot{\beta}_j(0) \right)=\left[\beta_j(t), \beta_j(0) \right]_{\mathcal{N}_{\mathcal{F}_0}^*}
\end{align}
\begin{align}\label{ou32}
\dot{\beta}_k(0)-\dot{\beta}_j(0)+ E_0'\left( \left( \dot{\theta}_{jk} (0) , \dot{B}_{jk} (0) \right) \right) = - \beta_j(0)\circ B_{jk}(0) - \beta_k (0) \circ B_{jk}(0)+ \mathcal{L}_{\theta_{jk}(0)}\circ \beta_j(0)+ \mathcal{L}_{\theta_{jk}(0)}\circ \beta_k(0)
\end{align}
\begin{align} \label{ou33}
\dot{B}_{ij} (0) - \dot{B}_{ik} (0) + \dot{B}_{jk}^i  (0) = \left[ \theta_{ij} (0) , B_{jk}^i(0)\right] - \left[ \theta_{jk}(0), B_{ij}(0)\right] + \left[ B_{ij} (0) ,   B_{jk}^i(0)\right]
\end{align}
\begin{align} \label{ou35}
&\sum_{\gamma=1}^q(-1)^{q-\gamma} \dot{\gamma}_j^\gamma (0) w_{0j}^1\wedge \cdots \wedge \hat{w}_{0j}^\gamma \wedge \cdots \wedge w_{0j}^q\wedge dw_{0j}^\alpha +w_{0j}^1\wedge \cdots \wedge w_{0j}^q\wedge d\left( \dot{\gamma}_j^\alpha (0) \right) -w_{0j}^1\wedge \cdots \wedge w_{0j}^q \wedge \dot{\beta}_j^\alpha(0)\\
&=\sum_{\gamma=1}^q\sum_{1\leq \beta \ne \gamma \leq q } (-1)^{q-\gamma} \gamma_j^\gamma (0) w_{0j}^1\wedge \cdots \wedge \overbrace{\beta_j^\beta(0) }^{\beta-\textnormal{th}} \wedge \cdots \wedge \hat{w}_{0j}^\gamma \wedge \cdots \wedge w_{0j}^q\wedge dw_{0j}^\alpha \notag \\
&+\sum_{\gamma=1}^q (-1)^{q-\gamma} \gamma_j^\gamma(0) w_{0j}^1\wedge \cdots \wedge \hat{w}_{0j}^\gamma \wedge \cdots w_{0j}^q \wedge d\left( \beta_j^\alpha(0)\right) + \sum_{\gamma=1}^q w_{0j}^1\wedge \cdots \wedge \overbrace{\beta_j^\gamma(0)}^{\gamma-\textnormal{th}}\wedge \cdots \wedge w_{0j}^q\wedge d\gamma_j^\alpha(0) \notag\\
& -  \sum_{\gamma=1}^q w_{0j}^1\wedge \cdots \wedge  \overbrace{\beta_j^\gamma(0)}^{\gamma-\textnormal{th}} \wedge \cdots \wedge w_{0j}^q\wedge \beta_j^\alpha(0) \notag
\end{align} 
\begin{align} \label{ou36}
&\dot{\gamma}_k(0)- \dot{\gamma}_j(0) + i_{\dot{\theta}_{jk}(0)} \\
&=  -\frac{1}{2}\gamma_j(0)\circ B_{jk}(0) - \frac{1}{2}\gamma_k(0) \circ B_{jk}(0) + \frac{1}{2}\mathcal{L}_{\theta_{jk}(0)}\circ \gamma_j(t) + \frac{1}{2}\mathcal{L}_{\theta_{jk}(0)}\circ \gamma_k(0) + \frac{1}{2} i_{\theta_{jk}(0)}\circ \beta_j(0) + \frac{1}{2} i_{\theta_{jk}(0)} \circ \beta_k(0)  \notag
\end{align}

Let us define $\left[ \left( \beta_j(0),\gamma_j(0)\right), \left( \beta_j(0), \gamma_j(0) \right) \right]_{ \tilde{\mathcal{N}}_{\mathcal{F}_0}^*} \in \Gamma\left( U_j^0, \mathscr{H}om_{\mathcal{O}_{M}}\left( \mathcal{N}_{\mathcal{F}_0}^*, \tilde{\mathcal{S}}^2\right) \bigoplus \mathscr{H}om_{\mathcal{O}_{M} } \left( \mathcal{N}_{\mathcal{F}_0}^*, \tilde{\mathcal{S}}^1 \right) \right)$ by
{\small{\begin{align*}
&\left[ \left( \beta_j(0),\gamma_j(0)\right), \left( \beta_j(0), \gamma_j(0) \right) \right]_{ \tilde{\mathcal{N}}_{\mathcal{F}_0}^*}:\Gamma\left( U_j^0, \mathcal{N}_{\mathcal{F}_0}^* \right) \to \Gamma\left( U_j^t, \tilde{\mathcal{S}}^2 \right) \bigoplus \Gamma\left( U_j^0  , \tilde{\mathcal{S}}^1 \right) \\
w_{0j}^\alpha \mapsto \left[ \beta_j(0), \beta_j(0)\right]_{\mathcal{N}_{\mathcal{F}_0}^*} &\oplus  \sum_{\gamma=1}^q\sum_{1\leq \beta \ne \gamma \leq q } (-1)^{q-\gamma} \gamma_j^\gamma (0) w_{0j}^1\wedge \cdots \wedge \overbrace{\beta_j^\beta(0) }^{\beta-\textnormal{th}} \wedge \cdots \wedge \hat{w}_{0j}^\gamma \wedge \cdots \wedge w_{0j}^q\wedge dw_{0j}^\alpha  \\
&+\sum_{\gamma=1}^q (-1)^{q-\gamma} \gamma_j^\gamma(0) w_{0j}^1\wedge \cdots \wedge \hat{w}_{0j}^\gamma \wedge \cdots w_{0j}^q \wedge d\left( \beta_j^\alpha(0)\right) + \sum_{\gamma=1}^q w_{0j}^1\wedge \cdots \wedge \overbrace{\beta_j^\gamma(0)}^{\gamma-\textnormal{th}}\wedge \cdots \wedge w_{0j}^q\wedge d\gamma_j^\alpha(0) \\
& -  \sum_{\gamma=1}^q w_{0j}^1\wedge \cdots \wedge  \overbrace{\beta_j^\gamma(0)}^{\gamma-\textnormal{th}} \wedge \cdots \wedge w_{0j}^q\wedge \beta_j^\alpha(0) 
\end{align*}}}
and linearly extends to $\Gamma\left( U_j^0, \mathcal{N}_{\mathcal{F}_0}^* \right)$. Then by combining $(\ref{ou34})$ and $(\ref{ou33})$ (see Appendix \ref{apd2}), we can write $(\ref{ou34})-(\ref{ou36})$ as follows:
\begin{align}\label{itt17}
&G_1' \left( \left( \dot{\beta}_j(0), \dot{\gamma}_j(0) \right) \right) =\left[ \left( \beta_j(0), \gamma_j(0) \right), \left( \beta_j(0), \gamma_j(0) \right)   \right]_{\tilde{\mathcal{N}}_{\mathcal{F}_0}^*}\\
&\left( \dot{\beta}_j(0) , \dot{\gamma}_j(0) \right)- \left( \dot{\beta}_j(0), \dot{\gamma}_j(0)  \right) + G_0'\left( \left( \dot{\theta}_{jk}(0), \dot{B}_{jk}(0) \right) \right)  \notag \\
&= \left( - \beta_j(0)\circ B_{jk}(0) - \beta_k (0) \circ B_{jk}(0)+ \mathcal{L}_{\theta_{jk}(0)}\circ \beta_j(0)+ \mathcal{L}_{\theta_{jk}(0)}\circ \beta_k(0), \right.   \notag \\
& \left. -\frac{1}{2}\gamma_j(0)\circ B_{jk}(0) - \frac{1}{2}\gamma_k(0) \circ B_{jk}(0) + \frac{1}{2}\mathcal{L}_{\theta_{jk}(0)}\circ \gamma_j(t) + \frac{1}{2}\mathcal{L}_{\theta_{jk}(0)}\circ \gamma_k(0) + \frac{1}{2} i_{\theta_{jk}(0)}\circ \beta_j(0) + \frac{1}{2} i_{\theta_{jk}(0)} \circ \beta_k(0) \right)  \notag \\
&\left( \dot{\theta}_{ij}(0), \dot{B}_{ij}(0)\right) - \left( \dot{\theta}_{ik}(0), \dot{B}_{ik}(0) \right) + \left( \dot{\theta}_{jk}(0), \dot{B}_{jk}(0) \right)= \left[  \left( \theta_{jk} (0) , B_{jk}(0) \right) , \left( \theta_{jk}(0), B_{jk}(0) \right)   \right] \notag
\end{align}

\begin{lemma}\label{ou50}
For any $1$-cocycle 
\begin{align*}
\left( \left\{\beta_j \right\} , \left\{ \gamma_j \right\}, \left\{\left( \theta_{jk}, B_{jk} \right) \right\} \right)\in C^0\left(\mathcal{U},  \mathscr{H}om_{\mathcal{O}_M}\left( \mathcal{N}_{\mathcal{F}_0}^*, \Omega_M^1 \right) \right) \bigoplus C^0\left( \mathcal{U}, \mathscr{H}om_{\mathcal{O}_M}\left( \mathcal{N}_{\mathcal{F}_0}^*, \mathcal{O}_M \right) \right)\bigoplus C^1\left(  \mathcal{U}, \mathcal{E}_{\mathcal{N}_{\mathcal{F}_0}^*} \right)
\end{align*}
 where $\beta_j:\Gamma\left( U_j, \mathcal{N}_{\mathcal{F}_0}^* \right)\to \Gamma\left( U_j, \Omega_M^1 \right) , \theta_{jk}\in\Gamma\left( U_{jk} , \Theta_M \right), \gamma_j : \Gamma\left( U_j, \mathcal{N}_{\mathcal{F}_0}^* \right) \to \Gamma\left( U_j, \mathcal{O}_M \right)$ and $B_{jk}:\Gamma\left(U_{jk}, \mathcal{N}_{\mathcal{F}_0}^* \right) \to \Gamma\left( U_{jk}, \mathcal{N}_{\mathcal{F}_0}^* \right)$ satisfying
\begin{enumerate}
\item $G_1'(\beta_j, \gamma_j)= 0$, i.e. $w_{0j}^1\wedge \cdots \wedge w_{0j}^q \wedge d \beta_j\left(w_{0j}^\alpha\right) + \sum_{\gamma=1}^q w_{0j}^1 \wedge \cdots \wedge \beta_j\left(w_{0j}^\gamma\right)\wedge \cdots \wedge w_{0j}^q \wedge dw_{0j}^\alpha$ and $\sum_{\gamma=1}^q (-1)^{q-\gamma} \gamma_j\left( w_{0j}^\gamma \right) w_{0j}^1 \wedge \cdots \wedge \widehat{w}_{0j}^\gamma \wedge \cdots \wedge w_{0j}^q \wedge dw_{0j}^\alpha + w_{0j}^1 \wedge \cdots \wedge w_{0j}^q \wedge d \gamma_j\left(w_{0j}^\alpha \right)- w_{0j}^1 \wedge \cdots \wedge w_{0j}^q \wedge \beta_j\left( w_{0j}^\alpha \right)=0$ for $\alpha=1,...,q$
\item $\left(\beta_k, \gamma_k \right) -  \left( \beta_j , \gamma_j \right) + G_0'\left( (\theta_{jk}, B_{jk})\right)=0$, i.e. $\beta_k\left(w_{0j}^\alpha\right) - \beta_j\left(w_{0j}^\alpha \right) +\mathcal{L}_{\theta_{jk}}\left(w_{0j}^\alpha \right) - B_{jk}\left(w_{0j}^\alpha \right)=0$ and $\gamma_k\left(w_{0j}^\alpha \right)- \gamma_j\left( w_{0j}^\alpha \right) + i_{\theta_{jk}}\left( w_{0j}^\alpha \right)=0$ for $\alpha=1,...,q$
\item $\left( \theta_{ij}, B_{ij} \right)- \left( \theta_{ik}, B_{ik} \right) + \left( \theta_{jk}, B_{jk} \right)=0$, i.e. $\theta_{ij}-\theta_{ik}+ \theta_{jk}=0$ and $B_{ij}- B_{ik} + B_{jk}^i=0$
\end{enumerate}
we define $S_j=\left[\beta_j, \beta_j \right]_{\mathcal{N}_{\mathcal{F}_0}^*} \in \Gamma\left( U_j, \mathscr{H}om_{\mathcal{O}_M}\left(  \mathcal{N}_{\mathcal{F}_0}^*, \tilde{\mathcal{S}}^2\right) \right)$ by
{\small{\begin{align*}
S_j:&\Gamma\left( U_j, \mathcal{N}_{\mathcal{F}_0}^* \right) \to \Gamma\left( U_j, \tilde{\mathcal{S}}^2\right) \\
w_{0j}^\alpha &\mapsto 2\sum_{\gamma=1}^q w_{0j}^1\wedge \cdots \wedge \beta_j\left(w_{0j}^\gamma\right)\wedge \cdots \wedge w_{0j}^q \wedge d \beta_j\left(w_{0j}^\alpha\right)+2\sum_{1\leq \gamma <\eta\leq q} w_{0j}^1\wedge \cdots \wedge \beta_j\left(w_{0j}^\gamma\right)\wedge \cdots  \wedge \beta_j\left(w_{0j}^\eta\right)\wedge \cdots \wedge w_{0j}^q\wedge dw_{0j}^\alpha
\end{align*}}}
and linearly extends to $\Gamma\left( U_j, \mathcal{N}_{\mathcal{F}_0}^* \right)$, and 
we define $Q_j\in \Gamma\left( U_j, \mathscr{H}om_{\mathcal{O}_M}\left( \mathcal{N}_{\mathcal{F}_0}^*,  \tilde{\mathcal{S}}^1 \right) \right)$ by
\begin{align*}
Q_j : \Gamma & \left( U_j, \mathcal{N}_{\mathcal{F}_0}^* \right)  \to \Gamma\left( U_j, \tilde{\mathcal{S}}^1 \right) \\
 w_{0j}^\alpha & \mapsto  \sum_{\gamma=1}^q\sum_{1\leq \beta \ne \gamma \leq q } (-1)^{q-\gamma} \gamma_j \left(w_{0j}^\gamma \right) w_{0j}^1\wedge \cdots \wedge \beta_j \left( w_{0j}^\beta \right)  \wedge \cdots \wedge \hat{w}_{0j}^\gamma \wedge \cdots \wedge w_{0j}^q\wedge dw_{0j}^\alpha  \\
&+\sum_{\gamma=1}^q (-1)^{q-\gamma} \gamma_j\left(w_{0j}^\gamma\right) w_{0j}^1\wedge \cdots \wedge \hat{w}_{0j}^\gamma \wedge \cdots w_{0j}^q \wedge d\left( \beta_j\left(w_{0j}^\alpha \right)\right) + \sum_{\gamma=1}^q w_{0j}^1\wedge \cdots \wedge \beta_j \left( w_{0j}^\gamma \right)\wedge \cdots \wedge w_{0j}^q\wedge d \left(\gamma_j\left(w_{0j}^\alpha \right) \right) \\
& -  \sum_{\gamma=1}^q w_{0j}^1\wedge \cdots \wedge  \beta_j \left(w_{0j}^\gamma \right)\wedge \cdots \wedge w_{0j}^q\wedge \beta_j\left( w_{0j}^\alpha \right)
\end{align*}
and linearly extends to $\Gamma\left( U_j, \mathcal{N}_{\mathcal{F}_0}^* \right)$. We define $G_{jk}\in \Gamma\left( U_{jk}, \mathscr{H}om_{\mathcal{O}_M} \left( \mathcal{N}_{\mathcal{F}_0}^*, \Omega_M^1 \right) \right)$ by 
\begin{align*}
G_{jk}: \Gamma\left( U_{jk}, \mathcal{N}_{\mathcal{F}_0}^* \right) &\to \Gamma\left(  U_{jk}, \Omega_M^1  \right) \\
w_{0j}^\alpha &\mapsto - \beta_j\left( B_{jk}\left(w_{0j}^\alpha \right) \right) - \beta_k\left( B_{jk}\left( w_{0j}^\alpha \right) \right) + \mathcal{L}_{\theta_{jk}}\left(\beta\left(w_{0j}^\alpha \right) \right)  +  \mathcal{L}_{\theta_{jk}}\left( \beta_k\left(w_{0j}^\alpha \right)\right)
\end{align*}
and linearly extends to $\Gamma\left( U_{jk}, \mathcal{N}_{\mathcal{F}_0}^* \right)$, and we define $C_{jk}\in \Gamma\left( U_{jk}, \mathscr{H}om_{\mathcal{O}_M}\left( \mathcal{N}_{\mathcal{F}_0}^*, \mathcal{O}_M   \right)  \right)$ by
\begin{align*}
&C_{jk}:\Gamma\left( U_{jk} , \mathcal{N}_{\mathcal{F}_0}^* \right)   \to \Gamma\left( U_{jk} , \mathcal{O}_M \right)\\
 &w_{0j}^\alpha   \mapsto - \frac{1}{2}\gamma_j\left( B_{jk}\left(w_{0j}^\alpha \right) \right) -  \frac{1}{2}\gamma_k \left( B_{jk}\left( w_{0j}^\alpha \right) \right) + \frac{1}{2}\left[ \theta_{jk}, \gamma_j\left(w_{0j}^\alpha \right) \right] + \frac{1}{2} \left[ \theta_{jk}, \gamma_k\left(w_{0j}^\alpha \right) \right] + \frac{1}{2} i_{\theta_{jk}}\left(\beta_j \left( w_{0j}^\alpha \right) \right) + \frac{1}{2} i_{\theta_{jk}}\left( \beta_k \left( w_{0j}^\alpha \right) \right)
\end{align*}
and linearly extends to $\Gamma\left( U_j, \mathcal{N}_{\mathcal{F}_0}^*  \right)$. Then
\begin{align} \label{itt19}
&\left( \left\{ S_j \right\}, \left\{ Q_j \right\}, \left\{ G_{jk} \right\}, \left\{ C_{jk} \right\} \left\{\left[\left( \theta_{ij}, B_{ij} \right) , \left(\theta_{jk}, B_{jk} \right) \right]  \right\} \right) \\
&\in C^0\left( \mathcal{U}, \mathscr{H}om_{\mathcal{O}_M}\left( \mathcal{N}_{\mathcal{F}_0}^*, \tilde{\mathcal{S}}^2 \right) \right) \bigoplus C^0\left( \mathcal{U}, \mathscr{H}om_{\mathcal{O}_M}\left( \mathcal{N}_{\mathcal{F}_0}^*, \tilde{\mathcal{S}}^1 \right) \right) \bigoplus C^1\left( \mathcal{U} , \mathscr{H}om_{\mathcal{O}_M}\left( \mathcal{N}_{\mathcal{F}_0}^*, \Omega_M^1 \right) \right) \notag \\
& \bigoplus C^1\left( \mathcal{U}, \mathscr{H}om_{\mathcal{O}_M}\left( \mathcal{N}_{\mathcal{F}_0}^*, \mathcal{O}_M \right) \right) \bigoplus C^2\left( \mathcal{U} , \mathcal{E}_{\mathcal{N}_{\mathcal{F}_0}^*}\right) \notag
\end{align}
defines a $2$-cocycle in \v Cech resolution of $(\ref{ou3})$ of $\mathcal{E}_{\tilde{\mathcal{N}}_{\mathcal{F}_0}^*}^\bullet$.
\end{lemma}

\begin{remark}
As in \textnormal{Remark \ref{oc80}}, we see that $\left\{G_{jk} \right\}\in C^1\left( \mathcal{U}, \mathscr{H}om_{\mathcal{O}_M}\left( \mathcal{N}_{\mathcal{F}_0}^*, \Omega_M^1 \right) \right)$. We show that $\left\{ C_{jk} \right\} \in C^1\left(\mathcal{U}, \mathscr{H}om_{\mathcal{O}_M} \left( \mathcal{N}_{\mathcal{F}_0}^*, \mathcal{O}_M \right) \right)$. In fact, from $(\ref{oc75})$, we have
{\Small{\begin{align*}
&C_{kj}\left(w_{0j}^\alpha \right)= \sum_{\beta=1}^q h_{0jk}^{\alpha\beta} C_{kj}\left(w_{0k}^\beta \right)\\
&=\sum_{\beta=1}^q h_{0jk}^{\alpha\beta}\left( - \frac{1}{2}\gamma_k\left( B_{kj}\left(w_{0k}^\beta \right) \right) -  \frac{1}{2}\gamma_j \left( B_{kj}\left( w_{0k}^\beta \right) \right) + \frac{1}{2}\left[ \theta_{kj}, \gamma_k\left(w_{0k}^\beta \right) \right] + \frac{1}{2} \left[ \theta_{kj}, \gamma_j\left(w_{0k}^\beta \right) \right] + \frac{1}{2} i_{\theta_{kj}}\left(\beta_k \left( w_{0k}^\beta \right) \right) + \frac{1}{2} i_{\theta_{kj}}\left( \beta_j \left( w_{0k}^\beta \right) \right) \right)\\
&= \frac{1}{2}\gamma_k\left( B_{jk}\left( w_{0j}^\alpha \right) \right) - \frac{1}{2}\sum_{\beta=1}^q \left[ \theta_{jk}, h_{0jk}^{\alpha\beta} \right] \gamma_k\left(w_{0k}^\beta \right) + \frac{1}{2} \gamma_j \left( B_{jk}\left( w_{0j}^\alpha \right) \right) - \frac{1}{2} \sum_{\beta=1}^q \left[ \theta_{jk}, h_{0jk}^{\alpha\beta} \right] \gamma_j\left( w_{0k}^\beta \right) - \frac{1}{2} \sum_{\beta=1}^q h_{0jk}^{\alpha\beta} \left[ \theta_{jk}, \gamma_k \left( w_{0k}^\beta  \right) \right] \\
& - \frac{1}{2} \sum_{\beta=1}^q h_{0jk}^{\alpha\beta} \left[ \theta_{jk}, \gamma_j\left( w_{0k}^\beta  \right) \right] - \frac{1}{2} i_{\theta_{jk} }\left( \beta_k\left( w_{0j}^\alpha \right) \right) - \frac{1}{2} i_{\theta_{jk}}\left( \beta_j \left( w_{0j}^\alpha \right) \right) \\
&= \frac{1}{2}\gamma_k \left( B_{jk}\left( w_{0j}^\alpha \right) \right) + \frac{1}{2} \gamma_j \left( B_{jk} \left( w_{0j}^\alpha \right) \right) - \frac{1}{2} \left[ \theta_{jk}, \gamma_k \left( w_{0j}^\alpha \right) \right] - \frac{1}{2} \left[ \theta_{jk}, \gamma_j \left( w_{0j}^\alpha \right) \right]  - \frac{1}{2} i_{\theta_{jk} }\left( \beta_k\left( w_{0j}^\alpha \right) \right) - \frac{1}{2} i_{\theta_{jk}}\left( \beta_j \left( w_{0j}^\alpha \right) \right) =  - C_{jk} \left( w_{0j}^\alpha \right)
\end{align*}}}
\end{remark}

\begin{remark}\label{os30}
As in \textnormal{Remark \ref{oc81}}, we see that $S_j$ is well-defined. We show that $Q_j$ is well-defined. In other words, $Q_j\left( w_{0j}^\alpha \right)$ actually lies in $\Gamma\left( U_j, \tilde{\mathcal{S}}^1 \right)$. In fact, for $x\in U_j-S$, there is a neighborhood $U_x\subset U_j-S$ such that $dw_{0j}^\alpha- \sum_{\delta=1}^q a_{j_x}^{\alpha\beta}\wedge w_{0j}^\delta=0$ on $U_x$ for some $a_{j_x}^{\alpha \delta} \in \Gamma\left(U_x, \Omega_M^1 \right),\alpha=1,...,q$. Then we have
{\small{\begin{align*}
w_{0j}^1\wedge \cdots \wedge w_{0j}^q \wedge \left(d\beta_j\left(w_{0j}^\alpha \right) - \sum_{\delta=1}^q a_{j_x}^{\alpha \delta} \wedge \beta_j\left(w_{0j}^\delta \right)\right)=0,\,\,\,\,\, w_{0j}^1\wedge \cdots \wedge w_{0j}^q  \wedge \left(d\gamma_j\left(w_{0j}^\alpha \right) - \sum_{\delta=1}^q \gamma_j\left(w_{0j}^\delta \right)  a_{j_x}^{\alpha \delta} - \beta_j\left(w_{0j}^\alpha \right) \right)=0
\end{align*}}}
so that we can write
\begin{align}
d\beta_j\left(w_{0j}^\alpha \right) &= \sum_{\delta=1}^q a_{j_x}^{\alpha\delta}\wedge \beta_j\left(w_{0j}^\delta \right) + \sum_{\eta=1}^q b_{j_x}^{\alpha \eta} \wedge w_{0j}^\eta \,\,\,\,\,\,\textnormal{on}\,\,\,U_x\,\,\,\,\textnormal{for some $b_{j_x}^{\alpha \eta}\in \Gamma\left( U_x, \Omega_M^1 \right)$} \label{os31}\\
d\gamma_j \left(w_{0j}^\alpha \right) &= \sum_{\delta=1}^q \gamma_j\left( w_{0j}^{\delta} \right)a_{j_x}^{\alpha\delta} + \beta_j\left( w_{0j}^\alpha  \right) + \sum_{\eta=1}^q c_{j_x}^{\alpha \eta} w_{0j}^\eta \,\,\,\,\,\,\textnormal{on}\,\,\,U_x\,\,\,\,\textnormal{for some $c_{j_x}^{\alpha \eta}\in \Gamma\left( U_x, \mathcal{O}_X \right)$} \label{os32}
\end{align}
Then we have
{\Small{\begin{align*}
&\sum_{\gamma=1}^q \sum_{\beta\ne \gamma} (-1)^{q-\gamma} \gamma_j\left(w_{0j}^\gamma\right) w_{0j}^1\wedge \cdots \wedge \beta_j\left(w_{0j}^\beta\right) \wedge \cdots \wedge \widehat{w}_{0j}^\gamma \wedge \cdots \wedge w_{0j}^q \wedge dw_{0j}^\alpha + \sum_{\gamma=1}^q (-1)^{q-\gamma} \gamma_j\left(w_{0j}^\gamma\right) w_j^1\wedge \cdots  \wedge \widehat{w}_{0j}^\gamma \wedge \cdots \wedge w_{0j}^q \wedge d\beta_j\left(w_{0j}^\alpha\right)\\
&+\sum_{\gamma=1}^q w_{0j}^1\wedge \cdots \wedge \beta_j\left(w_{0j}^\gamma\right)\wedge \cdots \wedge w_{0j}^q \wedge d\gamma_j\left(w_{0j}^\alpha\right)-\sum_{\gamma=1}^q w_{0j}^1 \wedge \cdots \wedge \beta_j\left(w_{0j}^\gamma\right) \wedge \cdots \wedge w_{0j}^q \wedge \beta_j \left(w_{0j}^\alpha\right)\\
&=\sum_{\gamma=1}^q \sum_{\beta \ne \gamma} (-1)^{q-r}\gamma_j\left(w_{0j}^\gamma\right) w_{0j}^1 \wedge \cdots \wedge \beta_j\left(w_{0j}^\beta\right) \wedge \cdots \wedge \widehat{w}_{0j}^\gamma \wedge \cdots \wedge w_{0j}^q \wedge \left( a_{j_x}^{\alpha \beta}\wedge w_{0j}^\beta + a_{j_x}^{\alpha\gamma}\wedge w_{0j}^\gamma      \right)\\
&+\sum_{\gamma=1}^q (-1)^{q-r}\gamma_j\left(w_{0j}^\gamma\right) w_{0j}^1 \wedge \cdots \wedge \widehat{w}_{0j}^\gamma \wedge \cdots \wedge w_{0j}^q \wedge \left(\sum_{\delta=1}^q a_{j_x}^{\alpha \delta}\wedge \beta_j\left(w_{0j}^\delta\right) + \sum_{\eta=1}^q b_{j_x}^{\alpha \eta} \wedge w_{0j}^\eta   \right)\\
&+ \sum_{\gamma=1}^q w_{0j}^1 \wedge \cdots \wedge \beta_j \left(w_{0j}^\gamma\right)\wedge \cdots \wedge w_{0j}^q \wedge\left( \sum_{\delta=1}^q  \gamma_j\left(w_{0j}^\delta\right) a_{j_x}^{\alpha \delta} +\sum_{\eta=1}^q c_{j_x}^{\alpha \eta} w_{0j}^\eta\right)\\
&= - \sum_{\gamma=1}^q \sum_{\beta \ne \gamma} (-1)^{q-\gamma}\gamma_j\left(w_{0j}^\gamma\right)w_{0j}^1 \wedge \cdots  \wedge \widehat{w}_{0j}^\gamma \wedge \cdots \wedge w_{0j}^q  \wedge \left( a_{j_x}^{\alpha\beta}\wedge \beta_j\left(w_{0j}^\beta\right)\right) - \sum_{\gamma=1}^q\sum_{\beta \ne \gamma} w_{0j}^1 \wedge \cdots \wedge \beta_j \left(w_{0j}^\beta\right)\wedge \cdots \wedge w_{0j}^q \wedge \left( \gamma_j\left(w_{0j}^\gamma\right) a_{j_x}^{\alpha \gamma} \right)\\
&+\sum_{\gamma=1}^q (-1)^{q-\gamma}\gamma_j\left(w_{0j}^\gamma\right) w_{0j}^1 \wedge\cdots \wedge \widehat{w}_{0j}^\gamma \wedge \cdots \wedge w_{0j}^q\wedge\left(\sum_{\delta=1}^q a_{j_x}^{\alpha\delta}\wedge \beta_j\left(w_{0j}^\delta\right) \right) - \sum_{\gamma=1}^q w_{0j}^1 \wedge \cdots \wedge w_{0j}^q \wedge\left( \gamma_j\left(w_{0j}^\gamma \right) b_{j_x}^{\alpha \gamma} \right)\\
&+ \sum_{\gamma=1}^q w_{0j}^1 \wedge \cdots \wedge \beta_j\left(w_{0j}^\gamma\right)\wedge \cdots \wedge w_{0j}^q \wedge\left( \sum_{\delta=1}^q  \gamma_j\left(w_{0j}^\delta\right) a_{j_x}^{\alpha \delta}+\sum_{\eta=1}^q c_{j_x}^{\alpha \eta} w_{0j}^\eta\right)
\end{align*}}}
We continue the computation. Then
{\Small{\begin{align*}
&= -\sum_{\gamma=1}^q  (-1)^{q-\gamma}\gamma_j\left(w_{0j}^\gamma\right)w_{0j}^1 \wedge \cdots  \wedge \widehat{w}_{0j}^\gamma \wedge \cdots \wedge w_{0j}^q  \wedge \left(\sum_{\beta=1}^q a_{j_x}^{\alpha\beta}\wedge \beta_j\left(w_{0j}^\beta\right) - a_{j_x}^{\alpha \gamma}\wedge \beta_j \left(w_{0j}^\gamma \right)\right) \\
& - \sum_{\gamma=1}^q w_{0j}^1 \wedge \cdots \wedge \beta_j\left(w_{0j}^\gamma\right)\wedge \cdots \wedge w_{0j}^q \wedge \left(\sum_{\beta=1}^q  \gamma_j\left(w_{0j}^\beta\right) a_{j_x}^{\alpha \beta} - \gamma_j\left(w_{0j}^\gamma\right) a_{j_x}^{\alpha \gamma}  \right)\\
&+\sum_{\gamma=1}^q (-1)^{q-\gamma}\gamma_j\left(w_{0j}^\gamma\right) w_{0j}^1 \wedge\cdots \wedge \widehat{w}_{0j}^\gamma \wedge \cdots \wedge w_{0j}^q\wedge\left(\sum_{\delta=1}^q a_{j_x}^{\alpha\delta}\wedge \beta_j \left(w_{0j}^\delta\right) \right) - \sum_{\gamma=1}^q w_{0j}^1 \wedge \cdots \wedge w_{0j}^q \wedge\left(\gamma_j\left(w_{0j}^\gamma\right) b_{j_x}^{\alpha \gamma}  \right)\\
&+ \sum_{\gamma=1}^q w_{0j}^1 \wedge \cdots \wedge \beta_j \left(w_{0j}^\gamma\right)\wedge \cdots \wedge w_{0j}^q \wedge\left( \sum_{\delta=1}^q  \gamma_j\left(w_{0j}^\delta\right) a_{j_x}^{\alpha \delta} +\sum_{\eta=1}^q c_{j_x}^{\alpha \eta} w_{0j}^\eta\right)\\
&=- w_{0j}^1 \wedge \cdots \wedge w_{0j}^q\wedge \left( \sum_{\gamma=1}^q  \gamma_j\left(w_{0j}^\gamma\right)  b_{j_x}^{\alpha \gamma} + \sum_{\gamma=1}^q c_{j_x}^{\alpha \gamma} \beta_j\left(w_{0j}^\gamma\right)    \right) \in \Gamma\left(U_x, \tilde{\mathcal{S}}^1\right)
\end{align*}}}

\end{remark}

\begin{remark}
We will deduce the equalities $(\ref{ou41}), (\ref{ou42}), (\ref{ou43}), (\ref{ou45}), (\ref{ou47})$ and $(\ref{ou48})$ below which will be used in the proof of \textnormal{Lemma \ref{ou50}}. Let $B_{jk}^\alpha:=B_{jk}\left(w_{0j}^\alpha \right) = \sum_{\beta=1}^q B_{jk}^{\alpha\beta} w_{0j}^\beta$ with $B_{jk}^{\alpha\beta} \in \Gamma\left( U_{jk}, \mathcal{O}_M \right)$. Then we have
{\small{\begin{align}\label{ou41}
&\sum_{\gamma=1}^q (-1)^{q-\gamma} \gamma_j\left(w_{0j}^\gamma \right) w_{0j}^1\wedge \cdots \wedge \widehat{w}_{0j}^\gamma \wedge \cdots \wedge w_{0j}^q\wedge dB_{jk}^\alpha +w_{0j}^1\wedge \cdots \wedge w_{0j}^q\wedge d\gamma_j\left(B_{jk}^\alpha\right)- w_{0j}^1\wedge \cdots \wedge w_{0j}^q\wedge \beta_j \left(B_{jk}^\alpha\right)\\
&=\sum_{\gamma,\eta=1}^q (-1)^{q-\gamma} \gamma_j\left(w_{0j}^\gamma\right) w_{0j}^1 \wedge \cdots \wedge \widehat{w}_{0j}^\gamma \wedge \cdots \wedge w_{0j}^q \wedge d\left( B_{jk}^{\alpha\eta} w_{0j}^\eta \right)+ \sum_{\eta=1}^q w_{0j}^1 \wedge \cdots \wedge w_{0j}^q \wedge d \left( B_{jk}^{\alpha\eta} \gamma_j\left(w_{0j}^\eta\right) \right) \notag \\
&-\sum_{\eta=1}^q B_{jk}^{\alpha\eta}  w_{0j}^1 \wedge \cdots \wedge w_{0j}^q \wedge \beta_j\left( w_{0j}^\eta\right) \notag \\
&=\sum_{\gamma,\eta=1}^q B_{jk}^{\alpha\eta} (-1)^{q-\gamma} \gamma_j\left(w_{0j}^\gamma\right) w_{0j}^1 \wedge \cdots \wedge \widehat{w}_{0j}^\gamma \wedge \cdots \wedge w_{0j}^q \wedge d w_{0j}^\eta+ \sum_{\eta=1}^q B_{jk}^{\alpha\eta} w_{0j}^1 \wedge \cdots \wedge w_{0j}^q \wedge d  \gamma_j\left(w_{0j}^\eta\right) \notag \\
&-\sum_{\eta=1}^q B_{jk}^{\alpha\eta}  w_{0j}^1 \wedge \cdots \wedge w_{0j}^q \wedge \beta_j\left( w_{0j}^\eta\right) -\sum_{\gamma=1}^q (-1)^{q-\gamma+q-\gamma} \gamma_j\left(w_{0j}^\gamma\right) w_{0j}^1 \wedge \cdots \wedge w_{0j}^q \wedge dB_{jk}^{\alpha\gamma} +\sum_{\eta=1}^q \gamma_j\left(w_{0j}^\eta\right) w_{0j}^1 \wedge \cdots \wedge w_{0j}^q \wedge dB_{jk}^{\alpha\eta}=0 \notag
\end{align}}}

We note that
{\small{\begin{align} \label{ou42}
&\sum_{\gamma=1}^q w_{0j}^1 \wedge \cdots \wedge \beta_j \left(w_{0j}^\gamma \right) \wedge \cdots \wedge w_{0j}^q \wedge B_{jk}^\alpha=\sum_{\gamma=1}^q w_{0j}^1\wedge \cdots \wedge \beta_j \left(w_{0j}^\gamma\right) \wedge \cdots \wedge w_{0j}^q \wedge\left(\sum_{\eta=1}^q B_{jk}^{\alpha\eta} w_{0j}^\eta \right)\\
&=-\sum_{\gamma=1}^q w_{0j}^1 \wedge \cdots \wedge w_{0j}^q \wedge B_{jk}^{\alpha\gamma} \beta_j \left(w_{0j}^\gamma\right)= -w_{0j}^1 \wedge \cdots \wedge w_{0j}^q \wedge \beta_j\left(B_{jk}^\alpha\right) \notag
\end{align}}}

We note that by taking $\mathcal{L}_{\theta_{jk}}$ on
\begin{align*}
\sum_{\gamma=1}^q (-1)^{q-\gamma} \gamma_k\left(w_{0j}^\gamma\right) w_{0j}^1\wedge \cdots \wedge \widehat{w}_{0j}^\gamma \wedge \cdots \wedge w_{0j}^q\wedge dw_{0j}^\alpha +w_{0j}^1\wedge \cdots \wedge w_{0j}^q\wedge d\gamma_k\left(w_{0j}^\alpha\right)- w_j^1\wedge \cdots \wedge w_{0j}^q\wedge \beta_k \left(w_{0j}^\alpha\right)=0
\end{align*}
we have
{\Small{\begin{align}\label{ou43}
&\sum_{\gamma=1}^q (-1)^{q-\gamma} \left[ \theta_{jk}, \gamma_k\left(w_{0j}^\gamma\right) \right] w_{0j}^1 \wedge \cdots \wedge \widehat{w}_{0j}^\gamma \wedge \cdots \wedge w_{0j}^q \wedge dw_j^\alpha +\sum_{\gamma=1}^q (-1)^{q-\gamma} \sum_{\beta\ne\gamma} \gamma_k\left(w_{0j}^\gamma\right) w_{0j}^1 \wedge \cdots \wedge \mathcal{L}_{\theta_{jk}}\left(w_{0j}^\beta\right) \wedge \cdots \wedge \widehat{w}_{0j}^\gamma \wedge\cdots \wedge w_{0j}^q\wedge dw_{0j}^\alpha\\
&+\sum_{\gamma=1}^q (-1)^{q-\gamma} \gamma_k\left(w_{0j}^\gamma\right) w_{0j}^1 \wedge \cdots \wedge \widehat{ w}_{0j}^\gamma \wedge \cdots \wedge w_{0j}^q \wedge d\mathcal{L}_{\theta_{jk}}\left(w_{0j}^\alpha\right)+\sum_{\gamma=1}^q w_{0j}^1 \wedge \cdots \wedge \mathcal{L}_{\theta_{jk}}\left( w_{0j}^\gamma  \right) \wedge \cdots \wedge w_{0j}^q \wedge d\gamma_k\left(w_{0j}^\alpha\right) \notag \\
&+w_{0j}^1 \wedge \cdots \wedge w_{0j}^q \wedge d\left[\theta_{jk}, \gamma_k\left(w_{0j}^\alpha\right)\right]-\sum_{\gamma=1}^q w_{0j}^1 \wedge \cdots \wedge L_{\theta_{jk}}\left(w_{0j}^\gamma\right)\wedge \cdots \wedge w_{0j}^q \wedge \beta_k\left(w_{0j}^\alpha\right) - w_{0j}^1 \wedge \cdots w_{0j}^q \wedge L_{\theta_{jk}}\left(\beta_k\left(w_{0j}^\alpha\right)\right)=0 \notag
\end{align}}}

We note that by taking $i_{\theta_{jk}}$ on  $\sum_{\beta=1}^q w_{0j}^1 \wedge \cdots \wedge \beta_k \left(w_{0j}^\beta \right)\wedge \cdots \wedge w_{0j}^q \wedge dw_{0j}^\alpha+w_{0j}^1\wedge \cdots \wedge w_{0j}^q \wedge d \beta_k\left( w_{0j}^\alpha\right)=0$,  we have
{\Small{\begin{align*}
0=&\sum_{\gamma=1}^q \sum_{\beta\ne \gamma} (-1)^{\gamma-1} i_{\theta_{jk}}\left(w_{0j}^\gamma\right) w_{0j}^1 \wedge \cdots \wedge\cdots \tilde{\beta}_k(w_j^\beta)\wedge \cdots \wedge  \hat{w}_j^\gamma \wedge \cdots \wedge w_j^q\wedge dw_j^\alpha+\sum_{\beta=1}^q(-1)^{\beta-1} w_j^1 \wedge \cdots \wedge i_{\theta_{jk}}(\tilde{\beta}_k(w_j^\beta))\wedge \cdots \wedge w_j^1 \wedge dw_j^\alpha\\
&+\sum_{\beta=1}^q(-1)^q w_j^1 \wedge \cdots \wedge \tilde{\beta}_k(w_j^\beta)\wedge \cdots \wedge w_j^q\wedge \left(i_{\theta_{jk}}dw_j^\alpha\right) +\sum_{\gamma=1}^q (-1)^{\gamma-1} w_j^1\wedge \cdots \wedge i_{\theta_{jk}}( w_j^\gamma ) \wedge \cdots \wedge \cdots \wedge w_j^q\wedge d \tilde{\beta}_k(w_j^\alpha)\\
&+ (-1)^q w_j^1 \wedge \cdots \wedge w_j^q \wedge \left(i_{\theta_{jk}}d\tilde{\beta}_k(w_j^\alpha) \right)\\
=&\sum_{\gamma=1}^q \sum_{\beta\ne \gamma} (-1)^{\gamma-1} i_{\theta_{jk}}(w_j^\gamma) w_j^1 \wedge \cdots \wedge\cdots \tilde{\beta}_k(w_j^\beta)\wedge \cdots \wedge  \hat{w}_j^\gamma \wedge \cdots \wedge w_j^q\wedge dw_j^\alpha+\sum_{\beta=1}^q(-1)^{\beta-1} w_j^1 \wedge \cdots \wedge i_{\theta_{jk}}(\tilde{\beta}_k(w_j^\beta))\wedge \cdots \wedge w_j^1 \wedge dw_j^\alpha\\
&+\sum_{\beta=1}^q(-1)^q w_j^1 \wedge \cdots \wedge \tilde{\beta}_k(w_j^\beta)\wedge \cdots \wedge w_j^q\wedge \left(L_{\theta_{jk}}(w_j^\alpha)\right)-\sum_{\beta=1}^q (-1)^q w_j^1 \wedge \cdots \wedge \tilde{\beta}_k(w_j^\beta)\wedge \cdots \wedge w_j^q \wedge d\left(i_{\theta_{jk}}(w_j^\alpha) \right) \\
&+\sum_{\gamma=1}^q (-1)^{\gamma-1} w_j^1\wedge \cdots \wedge i_{\theta_{jk}}( w_j^\gamma ) \wedge \cdots \wedge \cdots \wedge w_j^q\wedge d \tilde{\beta}_k(w_j^\alpha)\\
&+ (-1)^q w_j^1 \wedge \cdots \wedge w_j^q \wedge \left( L_{\theta_{jk}}(\tilde{\beta}_k(w_j^\alpha)) \right)-(-1)^q w_j^1 \wedge \cdots \wedge w_j^q \wedge d\left( i_{\theta_{jk}}(\tilde{\beta}_k(w_j^\alpha)) \right)
\end{align*}}}
We continue the computation. Then
{\small{\begin{align}\label{ou45}
=&\sum_{\gamma=1}^q \sum_{\beta\ne \gamma} (-1)^{q-\gamma} i_{\theta_{jk}}(w_j^\gamma) w_j^1 \wedge \cdots \wedge\cdots \tilde{\beta}_k(w_j^\beta)\wedge \cdots \wedge  \hat{w}_j^\gamma \wedge \cdots \wedge w_j^q\wedge dw_j^\alpha+\sum_{\beta=1}^q(-1)^{q-\beta} w_j^1 \wedge \cdots \wedge i_{\theta_{jk}}(\tilde{\beta}_k(w_j^\beta))\wedge \cdots \wedge w_j^1 \wedge dw_j^\alpha \\
&-\sum_{\beta=1}^q w_j^1 \wedge \cdots \wedge \tilde{\beta}_k(w_j^\beta)\wedge \cdots \wedge w_j^q\wedge \left(L_{\theta_{jk}}(w_j^\alpha)\right) +\sum_{\beta=1}^q w_j^1 \wedge \cdots \wedge \tilde{\beta}_k(w_j^\beta)\wedge \cdots \wedge w_j^q \wedge d\left(i_{\theta_{jk}}(w_j^\alpha) \right) \notag \\
&+\sum_{\gamma=1}^q (-1)^{q-\gamma} w_j^1\wedge \cdots \wedge i_{\theta_{jk}}( w_j^\gamma ) \wedge \cdots \wedge \cdots \wedge w_j^q\wedge d \tilde{\beta}_k(w_j^\alpha) \notag \\
&- w_j^1 \wedge \cdots \wedge w_j^q \wedge \left( L_{\theta_{jk}}(\tilde{\beta}_k(w_j^\alpha)) \right) +w_j^1 \wedge \cdots \wedge w_j^q \wedge d\left( i_{\theta_{jk}}(\tilde{\beta}_k\left(w_j^\alpha\right)) \right)=0 \notag
\end{align}}}

By taking $i_{\theta_{jk}}$ on $\sum_{\beta=1}^q w_{0j}^1 \wedge \cdots \wedge \mathcal{L}_{\theta_{jk}}\left(w_{0j}^\beta\right)\wedge \cdots \wedge w_j^q \wedge dw_j^\alpha+ w_j^1 \wedge \cdots \wedge w_j^q \wedge dL_{\theta_{jk}}(w_j^\alpha)=0$, we have $($we note that $L_{\theta_{jk}}i_{\theta_{jk}}=i_{\theta_{jk}}L_{\theta_{jk}}$$)$
\begin{align*}
&\sum_{\beta=1}^q \sum_{\beta \ne \gamma} (-1)^{\gamma-1} i_{\theta_{jk}}(w_j^\gamma) w_j^1 \wedge \cdots \wedge \hat{w}_j^\gamma  \wedge \cdots \wedge L_{\theta_{jk}}(w_j^\beta) \wedge \cdots \wedge w_j^q \wedge dw_j^\alpha\\
&+ \sum_{\beta=1}^q (-1)^{\beta-1} w_j^1 \wedge \cdots \wedge L_{\theta_{jk}}(i_{\theta_{jk}}(w_j^\beta))\wedge \cdots \wedge w_j^q \wedge dw_j^\alpha +\sum_{\beta=1}^q(-1)^q w_j^1 \wedge \cdots \wedge L_{\theta_{jk}}(w_j^\beta) \wedge \cdots \wedge w_j^q \wedge L_{\theta_{jk}}(w_j^\alpha)\\
&-\sum_{\beta=1}^q(-1)^q w_j^1 \wedge \cdots \wedge L_{\theta_{jk}}(w_j^\beta)\wedge \cdots \wedge w_j^q \wedge d \left( i_{\theta_{jk}}(w_j^\alpha) \right)\\
&+\sum_{\gamma=1}^q (-1)^{q-1} i_{\theta_{jk}}(w_j^\gamma) w_j^1 \wedge \cdots \wedge \hat{w}_j^\gamma \wedge \cdots \wedge w_j^q \wedge dL_{\theta_{jk}}(w_j^\alpha) +(-1)^q w_j^1\wedge \cdots \wedge w_j^q \wedge  \left(L_{\theta_{jk}}L_{\theta_{jk}}(w_j^\alpha)\right)\\
&-(-1)^q w_j^1 \wedge \cdots \wedge w_j^q \wedge dL_{\theta_{jk}}(i_{\theta_{jk}}(w_j^\alpha))
\end{align*}

\begin{align}\label{ou47}
&\sum_{\beta=1}^q \sum_{\beta \ne \gamma} (-1)^{q-\gamma} i_{\theta_{jk}}(w_j^\gamma) w_j^1 \wedge \cdots \wedge \hat{w}_j^\gamma  \wedge \cdots \wedge L_{\theta_{jk}}(w_j^\beta) \wedge \cdots \wedge w_j^q \wedge dw_j^\alpha\\
&+ \sum_{\beta=1}^q (-1)^{q-\beta} w_j^1 \wedge \cdots \wedge L_{\theta_{jk}}(i_{\theta_{jk}}(w_j^\beta))\wedge \cdots \wedge w_j^q \wedge dw_j^\alpha -\sum_{\beta=1}^q  w_j^1 \wedge \cdots \wedge L_{\theta_{jk}}(w_j^\beta) \wedge \cdots \wedge w_j^q \wedge L_{\theta_{jk}}(w_j^\alpha) \notag\\
&+\sum_{\beta=1}^q  w_j^1 \wedge \cdots \wedge L_{\theta_{jk}}(w_j^\beta)\wedge \cdots \wedge w_j^q \wedge d \left( i_{\theta_{jk}}(w_j^\alpha) \right) \notag\\
&+\sum_{\gamma=1}^q  i_{\theta_{jk}}(w_j^\gamma) w_j^1 \wedge \cdots \wedge \hat{w}_j^\gamma \wedge \cdots \wedge w_j^q \wedge dL_{\theta_{jk}}(w_j^\alpha) - w_j^1\wedge \cdots \wedge w_j^q \wedge  \left(L_{\theta_{jk}}L_{\theta_{jk}}(w_j^\alpha)\right) \notag \\
&+ w_j^1 \wedge \cdots \wedge w_j^q \wedge dL_{\theta_{jk}}(i_{\theta_{jk}}(w_j^\alpha)) \notag
\end{align}

We note that by setting $H_{0jk}= \det\left( h_{0ij}^{\alpha\beta} \right)$, we have
{\Small{\begin{align} \label{ou48}
&\sum_{\gamma=1}^q \sum_{\beta\ne \gamma} (-1)^{q-\gamma} \gamma_k\left(w_{0j}^\gamma\right) w_{0j}^1\wedge \cdots \wedge \beta_k\left(w_{0j}^\beta\right) \wedge \cdots \wedge \widehat{w}_{0j}^\gamma \wedge \cdots \wedge w_{0j}^q \wedge dw_{0j}^\alpha\\
& + \sum_{\gamma=1}^q (-1)^{q-\gamma} \gamma_k\left(w_{0j}^\gamma\right) w_{0j}^1\wedge \cdots  \wedge \widehat{w}_{0j}^\gamma \wedge \cdots \wedge w_{0j}^q \wedge d \beta_k\left(w_{0j}^\alpha\right) \notag\\
&+\sum_{\gamma=1}^q w_{0j}^1\wedge \cdots \wedge \beta_k \left(w_{0j}^\gamma\right)\wedge \cdots \wedge w_{0j}^q \wedge d\gamma_k\left(w_{0j}^\alpha \right)-\sum_{\gamma=1}^q w_{0j}^1 \wedge \cdots \wedge \beta_k \left(w_{0j}^\gamma\right) \wedge \cdots \wedge w_{0j}^q \wedge \beta_k \left(w_{0j}^\alpha\right) \notag\\
&=H_{0jk} \sum_{\gamma=1}^q \sum_{\beta\ne \gamma} (-1)^{q-\gamma} \gamma_k\left(w_{0j}^\gamma\right) w_{0j}^1\wedge \cdots \wedge \beta_k\left(w_{0k}^\beta\right) \wedge \cdots \wedge \widehat{w}_{0k}^\gamma \wedge \cdots \wedge w_{0k}^q \wedge d\left( \sum_{\eta=1}^q h_{0jk}^{\alpha\eta} w_{0k}^\eta \right) \notag \\
& + H_{0jk}\sum_{\gamma=1}^q (-1)^{q-\gamma} \gamma_k\left(w_{0k}^\gamma\right) w_{0k}^1\wedge \cdots  \wedge \widehat{w}_{0k}^\gamma \wedge \cdots \wedge w_{0k}^q \wedge d \left( \sum_{\beta=1}^q h_{0jk}^{\alpha\beta} \beta_k\left(w_{0k}^\beta\right)\right)\notag \\
&+H_{0jk}\sum_{\gamma=1}^q w_{0k}^1\wedge \cdots \wedge \beta_k \left(w_{0k}^\gamma\right)\wedge \cdots \wedge w_{0k}^q \wedge d\left(\sum_{\beta=1}^q h_{0jk}^{\alpha\beta}\gamma_k\left(w_{0k}^\beta \right) \right)-H_{0jk}\sum_{\gamma=1}^q w_{0j}^1 \wedge \cdots \wedge \beta_k \left(w_{0k}^\gamma\right) \wedge \cdots \wedge w_{0k}^q \wedge  \left(\sum_{\beta=1}^q h_{0jk}^{\alpha\beta} \beta_k \left(w_{0k}^\beta\right) \right) \notag \\
&= H_{0jk} \sum_{\beta=1}^q h_{0jk}^{\alpha\beta} Q_k\left(w_{0k}^\beta \right) -H_{0jk} \sum_{\gamma=1}^q \sum_{\beta\ne \gamma} (-1)^{q-\gamma} \gamma_k\left(w_{0j}^\gamma\right) w_{0j}^1\wedge \cdots \wedge \beta_k\left(w_{0k}^\beta\right) \wedge \cdots \wedge \widehat{w}_{0k}^\gamma \wedge \cdots \wedge w_{0k}^q \wedge \left(  w_{0k}^\beta \wedge d h_{jk}^{\alpha\beta} + w_{0k}^\gamma \wedge d h_{jk}^{\alpha \gamma} \right)\notag \\
&- H_{0jk} \sum_{\gamma, \beta=1}^q (-1)^{q-\gamma} \gamma_k\left(w_{0k}^\gamma \right) w_{0k}^1 \wedge \cdots \wedge \widehat{w}_{0k}^\gamma \wedge \cdots \wedge w_{0k}^q \wedge \beta_k\left(w_{0k}^\beta \right) \wedge d h_{0jk}^{\alpha\beta} + H_{0jk}\sum_{\gamma,\beta=1}^q \gamma_k\left( w_{0k}^\beta  \right) w_{0k}^1 \wedge \cdots \wedge \beta_k\left( w_{0k}^\gamma \right) \wedge \cdots \wedge w_{0k}^q \wedge d h_{0jk}^{\alpha\beta}  \notag\\
&=  H_{0jk} \sum_{\beta=1}^q h_{0jk}^{\alpha\beta} Q_k\left(w_{0k}^\beta \right) - H_{0jk} \sum_{\beta=1}^q (-1)^{q-\beta}\gamma_k\left(w_{0k}^\beta \right) w_{0k}^1 \wedge \cdots \wedge \widehat{w}_{0k}^\beta \wedge \cdots \wedge w_{0k}^q \wedge  \beta_k \left( w_{0k}^\beta \right) \wedge d h_{0jk}^{\alpha\beta} \notag \\
&+ H_{0jk} \sum_{\beta=1}^q \gamma_k\left(w_{0k}^\beta \right) w_{0k}^1 \wedge \cdots \wedge \beta_k\left( w_{0k}^\beta \right) \wedge \cdots \wedge w_{0k}^q \wedge d h_{0jk}^{\alpha\beta} =  H_{0jk}\sum_{\beta=1}^q h_{0jk}^{\alpha\beta} Q_k \left( w_{0k}^\beta \right) \notag
\end{align}}} 

\end{remark}

\begin{remark}\label{os33}
The following computation will be used in the proof of \textnormal{Lemma \ref{ou50}}. By applying $d$ on $(\ref{os32})$
\begin{align*}
d\gamma_j \left(w_{0j}^\alpha \right) &= \sum_{\delta=1}^q \gamma_j\left( w_{0j}^{\delta} \right)a_{j_x}^{\alpha\delta} + \beta_j\left( w_{0j}^\alpha  \right) + \sum_{\eta=1}^q c_{j_x}^{\alpha \eta} w_{0j}^\eta
\end{align*}
we have from $(\ref{os31}),(\ref{os32})$ and $(\ref{oc40})$
{\small{\begin{align*}
0=&\sum_{\delta=1}^q \gamma_j\left(w_{0j}^\delta\right)d a_{j_x}^{\alpha\delta}  -\sum_{\delta=1}^q a_{j_x}^{\alpha \delta} \wedge d\gamma_j\left(w_{0j}^\delta\right) + d \beta_j\left(w_{0j}^\alpha\right) + \sum_{\eta=1}^q dc_{j_x}^{\alpha \eta} \wedge w_{0j}^\eta + \sum_{\eta=1}^q c_{j_x}^{\alpha \eta} dw_{0j}^\eta\\
&= \sum_{\delta=1}^q \gamma_j\left(w_{0j}^\delta\right)d a_{j_x}^{\alpha \delta} - \sum_{\delta=1}^q a_{j_x}^{\alpha \delta}\wedge \left(\sum_{\eta=1}^q a_{j_x}^{\delta \eta} \gamma_j\left(w_{0j}^\eta\right)+ \beta_j\left(w_{0j}^\delta\right) + \sum_{\eta=1}^q c_{j_x}^{\delta \eta} w_{0j}^\eta  \right) \\
&+ \sum_{\delta=1}^q a_{j_x}^{\alpha \delta} \wedge \beta_j\left(w_{0j}^\delta\right) +\sum_{\gamma=1}^q b_{j_x}^{\alpha \gamma}\wedge w_{0j}^\gamma + \sum_{\eta=1}^q dc_{j_x}^{\alpha \eta} \wedge w_{0j}^\eta + \sum_{\eta, \delta=1}^q c_{j_x}^{\alpha \delta} a_{j_x}^{\delta \eta}\wedge w_{0j}^\eta\\
&= \sum_{\delta=1}^q \left( da_{j_x}^{\alpha \delta} - \sum_{\gamma=1}^q a_{j_x}^{\alpha \gamma}\wedge a_{j_x}^{\gamma \delta}    \right)\gamma_j\left(w_{0j}^\delta\right)  + \sum_{\eta=1}^q \left( dc_{j_x}^{\alpha \eta} + b_{j_x}^{\alpha \eta} -\sum_{\delta=1}^q c_{j_x}^{\delta \eta} a_{j_x}^{\alpha \delta}  + \sum_{\delta=1}^q c_{j_x}^{\alpha \delta} a_{j_x}^{\delta \eta}   \right) \wedge w_j^\eta \\
&=\sum_{\delta=1}^q \left( \sum_{\eta=1}^q g_{j_x}^{\alpha \delta \eta} \wedge w_{0j}^\eta \right) \gamma_j\left(w_{0j}^\delta\right)  + \sum_{\eta=1}^q \left( dc_{j_x}^{\alpha \eta} + b_{j_x}^{\alpha \eta} -\sum_{\delta=1}^q c_{j_x}^{\delta \eta} a_{j_x}^{\alpha \delta}  + \sum_{\delta=1}^q c_{j_x}^{\alpha \delta} a_{j_x}^{\delta \eta}   \right) \wedge w_{0j}^\eta  \\
&=\sum_{\eta=1}^q \left( \sum_{\delta=1}^q  g_{j_x}^{\alpha \delta \eta} \gamma_j\left(w_{0j}^\delta\right) + dc_{j_x}^{\alpha \eta} + b_{j_x}^{\alpha \eta} -\sum_{\delta=1}^q c_{j_x}^{\delta \eta} a_{j_x}^{\alpha \delta}  + \sum_{\delta=1}^q c_{j_x}^{\alpha \delta} a_{j_x}^{\delta \eta} \right) \wedge w_{0j}^\eta
\end{align*}}}
This implies that
\begin{align}\label{os34}
 dc_{j_x}^{\alpha \eta} + b_{j_x}^{\alpha \eta} -\sum_{\delta=1}^q c_{j_x}^{\delta \eta} a_{j_x}^{\alpha \delta}  + \sum_{\delta=1}^q c_{j_x}^{\alpha \delta} a_{j_x}^{\delta \eta} \equiv_{\mathcal{N}_{\mathcal{F}_0}^*} - \sum_{\delta=1}^q g_{j_x}^{\alpha \delta \eta} \gamma_j\left(w_{0j}^\delta \right)
\end{align}
\end{remark}

\begin{proof}[Proof of \textnormal{Lemma \ref{ou50}}]

By Lemma \ref{oc29}, it remains to show that $(-1)^2 S_j + E_1'\left( Q_j \right)= 0$ and $Q_j- Q_k - G_{jk} + \tilde{E}_0\left( C_{jk} \right)=0$ and $-\left(C_{ij}- C_{ik}+ C_{jk}\right) + i_{\left[ \theta_{ij}, \theta_{jk} \right]}=0$.

We show that $(-1)^2 S_j + E_1'\left( Q_j \right)= 0$. In fact, for $x\in U_j-S$, choose a neighborhood $U_x\subset U_j-S$ such that $dw_{0j}^\alpha- \sum_{\delta=1}^q a_{j_x}^{\alpha \delta} \wedge w_{0j}^\delta$ for some $a_{j_x}^{\alpha \delta}\in \Gamma\left( U_x, \Omega_M^1 \right)$. We can check $(-1)^2 S_j+ E_1'(Q_j)=0$ locally on $U_x$. We keep the notations in Remark \ref{oc81}, Remark \ref{oc43}, Remark \ref{os30} and Remark \ref{os33}, so that $S_j\left(w_{0j}^\alpha \right)=-2 w_{0j}^1 \wedge \cdots \wedge w_{0j}^q \wedge \left(\sum_{\gamma=1}^q b_{j_x}^{\alpha \gamma} \wedge \beta_j\left(w_{0j}^\gamma \right) \right)\in \Gamma\left(U_x, \tilde{\mathcal{S}}^2 \right)$ and $Q\left( w_{0j}^\alpha \right)= - w_{0j}^1 \wedge \cdots \wedge w_{0j}^q \wedge \left( \sum_{\gamma=1}^q \gamma_j\left( w_{0j}^\gamma \right) b_{j_x}^{\alpha \gamma} + \sum_{\gamma=1}^q c_{j_x}^{\alpha \gamma} \beta_j\left( w_{0j}^\gamma  \right)\right)\in \Gamma\left( U_x, \tilde{\mathcal{S}}^1 \right)$. Then from $(\ref{os31}), (\ref{os32})$, $(\ref{os34}), (\ref{oc42})$ and $(\ref{oc41})$, we have

{\small{\begin{align}\label{its49}
&d\left(  \sum_{\gamma=1}^q \gamma_j\left(w_{0j}^\gamma\right)  b_{j_x}^{\alpha \gamma} + \sum_{\gamma=1}^q c_{j_x}^{\alpha \gamma} \beta_j\left(w_{0j}^\gamma\right)   \right) - \sum_{\beta=1}^q a_{j_x}^{\alpha \beta} \wedge \left( \sum_{\gamma=1}^q \gamma_j\left(w_{0j}^\gamma\right)  b_{j_x}^{\beta \gamma} + \sum_{\gamma=1}^q c_{j_x}^{\beta \gamma} \beta_j \left(w_{0j}^\gamma\right)  \right)  + (-1)^2 2\sum_{\gamma=1}^q b_{j_x}^{\alpha \gamma} \wedge \beta_j \left(w_{0j}^\gamma\right) \\
&=  - \sum_{\gamma=1}^q b_{j_x}^{\alpha \gamma}\wedge d\gamma_j(w_j^\gamma) + \sum_{\gamma=1}^q \gamma_j(w_j^\gamma)db_{j_x}^{\alpha \gamma} + \sum_{\gamma=1}^q d c_{j_x}^{\alpha \gamma} \wedge \beta_j \left(w_{0j}^\gamma\right) + \sum_{\gamma=1}^q c_{j_x}^{\alpha \gamma} d\beta_j\left(w_{0j}^\gamma\right) \notag \\
&- \sum_{\beta,\gamma=1}^q \gamma_j\left(w_{0j}^\gamma\right) a_{j_x}^{\alpha \beta} \wedge b_{j_x}^{\beta\gamma} - \sum_{\beta,\gamma=1}^q c_{j_x}^{\beta\gamma} a_{j_x}^{\alpha\beta}\wedge \beta_j \left(w_{0j}^\gamma \right)+ 2\sum_{\gamma=1}^q b_{j_x}^{\alpha \gamma} \wedge \beta_j \left(w_{0j}^\gamma\right) \notag \\
&= - \sum_{\gamma=1}^q b_{j_x}^{\alpha \gamma} \wedge \left( \sum_{\delta=1}^q a_{j_x}^{\gamma \delta} \gamma_j\left(w_{0j}^\delta\right) + \beta_j \left(w_{0j}^\gamma\right) + \sum_{\eta=1}^q c_{j_x}^{\alpha \eta} w_{0j}^\eta  \right) + \sum_{\gamma=1}^q \gamma_j\left(w_{0j}^\gamma\right) db_{j_x}^{\alpha \gamma} + \sum_{\gamma=1}^q d c_{j_x}^{\alpha \gamma} \wedge \beta_j \left( w_{0j}^\gamma\right) \notag \\
& + \sum_{\gamma=1}^q c_{j_x}^{\alpha \gamma}\left(\sum_{\delta=1}^q a_{j_x}^{\gamma \delta} \wedge \beta_j \left(w_{0j}^\delta \right) + \sum_{\delta=1}^q b_{j_x}^{\gamma \delta} \wedge w_{0j}^\delta \right) - \sum_{\beta, \gamma=1}^q  \gamma_j\left(w_{0j}^\gamma\right) a_{j_x}^{\alpha \beta} \wedge b_{j_x}^{\beta \gamma} - \sum_{\beta, \gamma=1}^q c_{j_x}^{\beta \gamma} a_{j_x}^{\alpha \beta} \wedge \beta_j \left(w_{0j}^\gamma\right) + 2\sum_{\gamma=1}^q b_{j_x}^{\alpha \gamma} \wedge \beta_j \left(w_{0j}^\gamma \right) \notag \\
&\equiv_{\mathcal{N}_\mathcal{F}^*} \sum_{\delta=1}^q\left( - \sum_{\gamma=1}^q  b_{j_x}^{\alpha \gamma} \wedge a_{j_x}^{\gamma \delta}  + db_{j_x}^{\alpha \delta} - \sum_{\beta=1}^q a_{j_x}^{\alpha \beta}\wedge b_{j_x}^{\beta \delta} \right) \gamma_j\left(w_{0j}^\delta\right) + \sum_{\gamma=1}^q \left(  b_{j_x}^{\alpha \gamma} + d c_{j_x}^{\alpha \gamma} + \sum_{\delta=1}^q c_{j_x}^{\alpha \delta} a_{j_x}^{\delta \gamma} - \sum_{\beta=1}^q c_{j_x}^{\beta \gamma} a_{j_x}^{\alpha\beta}     \right)\wedge \beta_j \left(w_{0j}^\gamma \right) \notag \\
&\equiv_{\mathcal{N}_{\mathcal{F}}^*} \sum_{\delta=1}^q \left( \sum_{\eta=1}^q g_{j_x}^{\alpha \delta \eta} \wedge \beta_j \left(w_{0j}^\eta\right) \right)  \gamma_j\left(w_{0j}^\delta\right)  + \sum_{\gamma=1}^q \left(  - \sum_{\delta=1}^q g_{j_x}^{\alpha \delta \gamma}\gamma_j\left(w_{0j}^\delta\right)   \right) \wedge \beta_j \left(w_{0j}^\gamma\right) \equiv_{\mathcal{N}_{\mathcal{F}_0}^*} 0\notag
\end{align}}}
This implies that $(-1)^2 S_j+ E_1'(Q_j)=0$.

From (\ref{ou41}) and (\ref{ou42}), we have
{\Small{\begin{align}\label{itt16}
&Q_j\left( w_{0j}^\alpha \right)=\sum_{\gamma=1}^q \sum_{\beta\ne \gamma} (-1)^{q-\gamma} \gamma_j\left(w_{0j}^\gamma\right) w_{0j}^1\wedge \cdots \wedge \beta_j\left(w_{0j}^\beta\right) \wedge \cdots \wedge \widehat{w}_{0j}^\gamma \wedge \cdots \wedge w_{0j}^q \wedge dw_{0j}^\alpha\\
& + \sum_{\gamma=1}^q (-1)^{q-\gamma} \gamma_j\left(w_{0j}^\gamma\right) w_{0j}^1\wedge \cdots  \wedge \widehat{w}_{0j}^\gamma \wedge \cdots \wedge w_{0j}^q \wedge d \beta_j\left(w_{0j}^\alpha\right) \notag\\
&+\sum_{\gamma=1}^q w_{0j}^1\wedge \cdots \wedge \beta_j \left(w_{0j}^\gamma\right)\wedge \cdots \wedge w_{0j}^q \wedge d\gamma_j\left(w_{0j}^\alpha \right)-\sum_{\gamma=1}^q w_{0j}^1 \wedge \cdots \wedge \beta_j \left(w_{0j}^\gamma\right) \wedge \cdots \wedge w_{0j}^q \wedge \beta_j \left(w_{0j}^\alpha\right) \notag\\
&=\sum_{\gamma=1}^q (-1)^{q-\gamma} \sum_{\beta\ne \gamma} \left(\gamma_k\left(w_{0j}^\gamma\right)+i_{\theta_{jk}}\left(w_{0j}^\gamma\right)\right) w_{0j}^1 \wedge \cdots \wedge \left(\beta_k\left(w_{0j}^\beta\right) +\mathcal{L}_{\theta_{jk}}\left(w_{0j}^\beta\right)- B_{jk}^\beta \right) \wedge \cdots \wedge \widehat{w}_{0j}^\gamma \wedge \cdots \wedge w_{0j}^q \wedge dw_{0j}^\alpha \notag\\
&+\sum_{\gamma=1}^q (-1)^{q-\gamma} \left( \gamma_k\left(w_{0j}^\gamma\right)+i_{\theta_{jk}}\left(w_{0j}^\gamma\right) \right) w_{0j}^1 \wedge \cdots \wedge \widehat{w}_{0j}^\gamma \wedge \cdots \wedge w_{0j}^q \wedge d\left(\beta_k\left(w_{0j}^\alpha\right)+ \mathcal{L}_{\theta_{jk}}\left(w_{0j}^\alpha\right) -B_{jk}^\alpha \right) \notag\\
&+\sum_{\gamma=1}^q w_{0j}^1 \wedge \cdots \wedge \left(\beta_k\left(w_{0j}^\gamma\right)+\mathcal{L}_{\theta_{jk}}\left(w_{0j}^\gamma\right)-B_{jk}^\gamma \right) \wedge \cdots \wedge w_{0j}^q \wedge d\left( \gamma_k\left(w_{0j}^\alpha\right) +i_{\theta_{jk}}\left(w_{0j}^\alpha\right)\right) \notag\\
&-\sum_{\gamma=1}^q w_{0j}^1\wedge \cdots \wedge \left( \beta_k \left(w_{0j}^\gamma\right) +\mathcal{L}_{\theta_{jk}}\left(w_{0j}^\gamma\right)-B_{jk}^\gamma \right) \wedge \cdots \wedge w_{0j}^q \wedge \left( \beta_k\left(w_{0j}^\alpha\right) + \mathcal{L}_{\theta_{jk}}\left(w_{0j}^\alpha \right) - B_{jk}^\alpha  \right) \notag
\end{align}}}
We continue the computation. Then
{\tiny{\begin{align}
&=\sum_{\gamma=1}^q (-1)^{q-\gamma} \sum_{\beta\ne \gamma} \left(\gamma_k\left(w_{0j}^\gamma\right)+i_{\theta_{jk}}\left(w_{0j}^\gamma\right)\right) w_{0j}^1 \wedge \cdots \wedge \left(\beta_k\left(w_{0j}^\beta\right) +\mathcal{L}_{\theta_{jk}}\left(w_{0j}^\beta\right) \right) \wedge \cdots \wedge \widehat{w}_{0j}^\gamma \wedge \cdots \wedge w_{0j}^q \wedge dw_{0j}^\alpha \notag\\
&-\sum_{\gamma=1}^q (-1)^{q-\gamma} \sum_{\beta\ne \gamma} B_{jk}^{\beta \beta} \gamma_j\left(w_{0j}^\gamma\right) w_{0j}^1 \wedge \cdots \wedge w_{0j}^\beta \wedge \cdots \wedge \widehat{w}_{0j}^\gamma \wedge \cdots \wedge w_{0j}^q \wedge dw_{0j}^\alpha \notag \\
&-\sum_{\gamma=1}^q (-1)^{q-\gamma} \sum_{\beta\ne \gamma} B_{jk}^{\beta \gamma} \gamma_j\left(w_{0j}^\gamma\right) w_{0j}^1 \wedge \cdots \wedge \overbrace{ w_{0j}^\gamma}^{\beta-\textnormal{th}} \wedge \cdots \wedge \widehat{w}_{0j}^\gamma \wedge \cdots \wedge w_{0j}^q \wedge dw_{0j}^\alpha \notag \\
&+\sum_{\gamma=1}^q (-1)^{q-\gamma} \left( \gamma_k\left(w_{0j}^\gamma\right)+i_{\theta_{jk}}\left(w_{0j}^\gamma\right) \right) w_{0j}^1 \wedge \cdots \wedge \widehat{w}_{0j}^\gamma \wedge \cdots \wedge w_{0j}^q \wedge d\left(\beta_k\left(w_{0j}^\alpha\right)+ \mathcal{L}_{\theta_{jk}}\left(w_{0j}^\alpha\right) \right) \notag \\
&-\sum_{\gamma=1}^q (-1)^{q-\gamma} \gamma_j \left(w_{0j}^\gamma\right)  w_{0j}^1 \wedge \cdots \wedge \widehat{w}_{0j}^\gamma \wedge \cdots \wedge w_{0j}^q \wedge d B_{jk}^\alpha \notag \\
&+\sum_{\gamma=1}^q w_{0j}^1 \wedge \cdots \wedge \left(\beta_k\left(w_{0j}^\gamma\right)+\mathcal{L}_{\theta_{jk}}\left(w_{0j}^\gamma\right) \right) \wedge \cdots \wedge w_{0j}^q \wedge d\left( \gamma_k\left(w_{0j}^\alpha\right) +i_{\theta_{jk}}\left(w_{0j}^\alpha\right)\right) - \sum_{\beta=1}^q B_{jk}^{\beta \beta} w_{0j}^1 \wedge \cdots \wedge w_{0j}^\beta  \wedge \cdots \wedge w_{0j}^q \wedge d\gamma_j\left( w_{0j}^\alpha\right) \notag \\
&-\sum_{\gamma=1}^q w_{0j}^1\wedge \cdots \wedge \left( \beta_k \left(w_{0j}^\gamma\right) +\mathcal{L}_{\theta_{jk}}\left(w_{0j}^\gamma\right) \right) \wedge \cdots \wedge w_{0j}^q \wedge \left( \beta_k\left(w_{0j}^\alpha\right) + \mathcal{L}_{\theta_{jk}}\left(w_{0j}^\alpha \right)   \right) +\sum_{\beta=1}^q B_{jk}^{\beta \beta}w_{0j}^1\wedge \cdots \wedge w_{0j}^\beta  \wedge \cdots \wedge w_{0j}^q \wedge \beta_j \left( w_{0j}^\alpha \right) \notag \\
& + \sum_{\gamma=1}^q w_{0j}^1 \wedge \cdots \wedge \beta_j\left(w_{0j}^\gamma \right) \wedge \cdots \wedge w_{0j}^q \wedge B_{jk}^\alpha  \notag
\end{align}}}
We continue the computation. Then
{\tiny{\begin{align}
&=\sum_{\gamma=1}^q (-1)^{q-\gamma} \sum_{\beta\ne \gamma} \left(\gamma_k\left(w_{0j}^\gamma\right)+i_{\theta_{jk}}\left(w_{0j}^\gamma\right)\right) w_{0j}^1 \wedge \cdots \wedge \left(\beta_k\left(w_{0j}^\beta\right) +\mathcal{L}_{\theta_{jk}}\left(w_{0j}^\beta\right) \right) \wedge \cdots \wedge \widehat{w}_{0j}^\gamma \wedge \cdots \wedge w_{0j}^q \wedge dw_{0j}^\alpha \notag \\
&+ \sum_{\beta=1}^q (-1)^{q-\beta} B_{jk}^{\beta \beta}\gamma_j\left(w_{0j}^\beta \right) w_{0j}^1 \wedge \cdots \wedge \widehat{w}_{0j}^\beta \wedge \cdots \wedge w_{0j}^q \wedge dw_{0j}^\alpha + \sum_{\beta=1}^q \left( -1\right)^{q-\beta} \sum_{\beta \ne \gamma} B_{jk}^{\beta \gamma}\gamma_j\left( w_{0j}^\gamma \right) w_{0j}^1 \wedge \cdots \wedge \widehat{w}_{0j}^\beta \wedge \cdots \wedge w_{0j}^q \wedge  d w_{0j}^\alpha \notag \\
&+\sum_{\gamma=1}^q (-1)^{q-\gamma} \left( \gamma_k\left(w_{0j}^\gamma\right)+i_{\theta_{jk}}\left(w_{0j}^\gamma\right) \right) w_{0j}^1 \wedge \cdots \wedge \widehat{w}_{0j}^\gamma \wedge \cdots \wedge w_{0j}^q \wedge d\left(\beta_k\left(w_{0j}^\alpha\right)+ \mathcal{L}_{\theta_{jk}}\left(w_{0j}^\alpha\right) \right) \notag \\
& +w_{0j}^1\wedge \cdots \wedge w_{0j}^q\wedge d\gamma_j\left(B_{jk}^\alpha\right)- 2w_{0j}^1\wedge \cdots \wedge w_{0j}^q\wedge \beta_j \left(B_{jk}^\alpha\right) +\sum_{\gamma=1}^q w_{0j}^1 \wedge \cdots \wedge \left(\beta_k\left(w_{0j}^\gamma\right)+\mathcal{L}_{\theta_{jk}}\left(w_{0j}^\gamma\right) \right) \wedge \cdots \wedge w_{0j}^q \wedge d\left( \gamma_k\left(w_{0j}^\alpha\right) +i_{\theta_{jk}}\left(w_{0j}^\alpha\right)\right) \notag \\
&-\sum_{\gamma=1}^q w_{0j}^1\wedge \cdots \wedge \left( \beta_k \left(w_{0j}^\gamma\right) +\mathcal{L}_{\theta_{jk}}\left(w_{0j}^\gamma\right) \right) \wedge \cdots \wedge w_{0j}^q \wedge \left( \beta_k\left(w_{0j}^\alpha\right) + \mathcal{L}_{\theta_{jk}}\left(w_{0j}^\alpha \right)   \right)  \notag
\end{align}}}

Then  from $(\ref{ou43}), (\ref{ou45}), (\ref{ou47})$ and $(\ref{ou48})$, we have
{\Small{\begin{align*}
&Q_j\left(w_{0j}^\alpha \right)- H_{0jk} \sum_{\beta=1}^q h_{0jk}^{\alpha\beta} Q_k\left( w_{0k}^\beta \right)\\
&= \sum_{\beta=1}^q (-1)^{q-\beta} \gamma_j\left(B_{jk}^\beta \right) w_{0j}^1 \wedge \cdots \wedge \widehat{w}_{0j}^\beta \wedge \cdots \wedge w_{0j}^q \wedge dw_{0j}^\alpha + w_{0j}^1 \wedge \cdots \wedge w_{0j}^q \wedge d\gamma_j\left( B_{jk}^\alpha \right)- 2w_{0j}^1 \wedge \cdots \wedge w_{0j}^q \wedge \beta_j\left(B_{jk}^\alpha \right)\\
& -\sum_{\gamma=1}^q (-1)^{q-\gamma} \left[ \theta_{jk}, \gamma_k\left(w_{0j}^\gamma\right) \right] w_{0j}^1 \wedge \cdots \wedge \widehat{w}_{0j}^\gamma \wedge \cdots \wedge w_{0j}^q \wedge dw_j^\alpha - w_{0j}^1 \wedge \cdots \wedge w_{0j}^q \wedge d\left[\theta_{jk}, \gamma_k\left(w_{0j}^\alpha\right)\right] +  w_{0j}^1 \wedge \cdots w_{0j}^q \wedge \mathcal{L}_{\theta_{jk}}\left(\beta_k\left(w_{0j}^\alpha\right)\right)\\
&- \sum_{\beta=1}^{q-\beta} w_{0j}^1 \wedge \cdots \wedge i_{\theta_{jk}}\left(\beta_k\left(w_{0j}^\beta \right) \right) \wedge \cdots \wedge w_{0j}^q \wedge dw_{0j}^\alpha - w_{0j}^1 \wedge \cdots \wedge w_{0j}^q \wedge d\left( i_{\theta_{jk}}\left(w_{0j}^\alpha \right) \right) +  w_{0j}^1 \wedge \cdots \wedge w_{0j}^q \wedge \mathcal{L}_{\theta_{jk}}\left(\beta_k\left(w_{0j}^\alpha \right) \right)\\
&- \sum_{\beta=1}^q (-1)^{q-\beta} w_{0j}^1 \wedge \cdots \wedge \mathcal{L}_{\theta_{jk}}\left(i_{\theta_{jk}}\left( w_{0j}^\beta \right) \right) \wedge \cdots \wedge w_{0j}^q \wedge dw_{0j}^\alpha - w_{0j}^1 \wedge \cdots \wedge w_{0j}^q \wedge d\mathcal{L}_{\theta_{jk}}i_{\theta_{jk}}\left(w_{0j}^\alpha \right)  + w_{0j}^1\wedge \cdots \wedge w_{0j}^q \wedge \mathcal{L}_{\theta_{jk}} \mathcal{L}_{\theta_{jk}}\left(w_{0j}^\alpha \right)
\end{align*}}}

On the other hand, we note that 
\begin{align}
& - \frac{1}{2}\gamma_j\left( B_{jk}\left(w_{0j}^\alpha \right) \right) -  \frac{1}{2}\gamma_k \left( B_{jk}\left( w_{0j}^\alpha \right) \right) + \frac{1}{2}\left[ \theta_{jk}, \gamma_j\left(w_{0j}^\alpha \right) \right] + \frac{1}{2} \left[ \theta_{jk}, \gamma_k\left(w_{0j}^\alpha \right) \right] + \frac{1}{2} i_{\theta_{jk}}\left(\beta_j \left( w_{0j}^\alpha \right) \right) + \frac{1}{2} i_{\theta_{jk}}\left( \beta_k \left( w_{0j}^\alpha \right) \right)\label{itt15}\\
&=- \gamma_k\left( B_{jk}^\alpha \right) - \frac{1}{2} i_{\theta_{jk}}\left( B_{jk}^\alpha \right) + \left[ \theta_{jk}, \gamma_k\left(w_{0j}^\alpha \right) \right] +\frac{1}{2}\mathcal{L}_{\theta_{jk}}i_{\theta_{jk}}\left(w_{0j}^\alpha \right) + i_{\theta_{jk}}\left( \beta_k\left( w_{0j}^\alpha \right)\right) + \frac{1}{2}i_{\theta_{jk}}\mathcal{L}_{\theta_{jk}}\left(w_{0j}^\alpha \right)- \frac{1}{2} i_{\theta_{jk}}\left( B_{jk}^\alpha \right) \notag \\
&=-\gamma_k\left(B_{jk}^\alpha \right) - i_{\theta_{jk}}\left( B_{jk}^\alpha \right) + \left[ \theta_{jk}, \gamma_k\left( w_{0j}^\alpha  \right) \right] + i_{\theta_{jk}}\left( \beta_k\left( w_{0j}^\alpha \right)\right) + \mathcal{L}_{\theta_{jk}}i_{\theta_{jk}}\left( w_{0j}^\alpha \right) \notag
\end{align}

Then from $(\ref{itt16})$, and $(\ref{itt15})$, we have $Q_j- Q_k- G_{jk} + \tilde{E}_0\left(C_{jk} \right)=0$.

We show that $-\left(C_{ij}- C_{ik}+ C_{jk}\right) + i_{\left[ \theta_{ij}, \theta_{jk} \right]}=0$. We recall from $(\ref{oc69})$ that $B_{ij}^\alpha- B_{ik}^\alpha+ \sum_{\beta=1}^q h_{0ij}^{\alpha\beta} B_{jk}^\beta=- \sum_{\beta=1}^q \left[ \theta_{jk}, h_{0ij}^{\alpha\beta} \right] w_{0j}^\beta$. Then

\begin{align*}
&2C_{ij}\left( w_{0i}^\alpha \right) - 2C_{ik}\left( w_{0i}^\alpha \right) +2 \sum_{\beta=1}^q h_{0ij}^{\alpha\beta} C_{jk}\left( w_{0j}^\beta \right)\\
&=- \gamma_i\left(B_{ij}^\alpha\right)- \gamma_j \left(B_{ij}^\alpha\right) + \mathcal{L}_{\theta_{ij}}\left(\gamma_i\left(w_{0i}^\alpha\right)\right) +\mathcal{L}_{\theta_{ij}}\left(\gamma_j\left(w_{0i}^\alpha\right)\right) +i_{\theta_{ij}}\left(\beta_i\left(w_{0i}^\alpha\right)\right) + i_{\theta_{ij}}\left(\beta_j\left(w_{0i}^\alpha\right) \right)\\
&+ \gamma_i\left(B_{ik}^\alpha\right)+ \gamma_k \left(B_{ik}^\alpha\right) - \mathcal{L}_{\theta_{ik}}\left(\gamma_i\left(w_{0i}^\alpha\right)\right) -\mathcal{L}_{\theta_{ik}}\left(\gamma_k\left(w_{0i}^\alpha\right)\right) -i_{\theta_{ik}}\left(\beta_i\left(w_{0i}^\alpha\right)\right) - i_{\theta_{ik}}\left(\beta_k\left(w_{0i}^\alpha\right)\right)\\
&+\sum_{\beta=1}^q h_{0ij}^{\alpha\beta}\left( - \gamma_j\left(B_{jk}^\beta\right)- \gamma_k \left(B_{jk}^\beta\right) + \mathcal{L}_{\theta_{jk}}\left(\gamma_j\left(w_{0j}^\beta\right)\right) +\mathcal{L}_{\theta_{jk}}\left(\gamma_k\left(w_{0j}^\beta\right)\right) +i_{\theta_{jk}}\left(\beta_j\left(w_{0j}^\beta\right)\right) + i_{\theta_{jk}}\left(\beta_k\left(w_{0j}^\beta\right)\right) \right) \\
&= - \gamma_i\left(B_{ij}^\alpha\right)- \cancel{ \gamma_j \left(B_{ij}^\alpha\right)} + \mathcal{L}_{\theta_{ij}}\left(\gamma_i\left(w_{0i}^\alpha\right)\right) +\mathcal{L}_{\theta_{ij}}\left(\gamma_j\left(w_{0i}^\alpha\right)\right) +i_{\theta_{ij}}\left(\beta_i\left(w_{0i}^\alpha\right)\right) + i_{\theta_{ij}}\left(\beta_j\left(w_{0i}^\alpha\right) \right)\\
&+ \gamma_i\left(B_{ik}^\alpha\right)+ \xcancel{\gamma_k \left(B_{ik}^\alpha\right)} - \mathcal{L}_{\theta_{ik}}\left(\gamma_i\left(w_{0i}^\alpha\right)\right) -\mathcal{L}_{\theta_{ik}}\left(\gamma_k\left(w_{0i}^\alpha\right)\right) -i_{\theta_{ik}}\left(\beta_i\left(w_{0i}^\alpha\right)\right) - i_{\theta_{ik}}\left(\beta_k\left(w_{0i}^\alpha\right)\right)\\
&+ \cancel{\gamma_j\left(B_{ij}^\alpha\right) }- \gamma_j\left(B_{ik}^\alpha\right) + \cancel{\sum_{\beta=1}^q \left[ \theta_{jk}, h_{0ij}^{\alpha\beta} \right] \gamma_j\left(w_{0j}^\beta\right)} + \gamma_k\left(B_{ij}^\alpha\right) - \xcancel{\gamma_k\left(B_{ik}^\alpha\right)} + \bcancel{\sum_{\beta=1}^q  \left[ \theta_{jk}, h_{0ij}^{\alpha\beta}\right] \gamma_k\left(w_{0j}^\beta\right)}\\
&+ \mathcal{L}_{\theta_{jk}}\left(\gamma_j\left(w_{0i}^\alpha\right)\right)- \cancel{\sum_{\beta=1}^q \left[  \theta_{jk}, h_{0ij}^{\alpha\beta}\right] \gamma_j\left(w_{0j}^\beta\right)} + \mathcal{L}_{\theta_{jk}}\left(\gamma_k\left(w_{0i}^\alpha\right)\right)- \bcancel{\sum_{\beta=1}^q \left[ \theta_{jk}, h_{0ij}^{\alpha\beta} \right] \gamma_k\left(w_{0j}^\beta\right)} + i_{\theta_{jk}}\left(\beta_j\left(w_{0i}^\alpha\right)\right) + i_{\theta_{jk}}\left(\beta_k\left(w_{0i}^\alpha\right)\right)\\
&= - \gamma_i\left(B_{ij}^\alpha\right) + \mathcal{L}_{\theta_{ij}}\left(\gamma_i\left(w_{0i}^\alpha\right)\right) +\mathcal{L}_{\theta_{ij}}\left(\gamma_j\left(w_{0i}^\alpha\right)\right) +i_{\theta_{ij}}\left(\beta_i\left(w_{0i}^\alpha\right)\right) + i_{\theta_{ik}}\left(\beta_j\left(w_{0i}^\alpha\right) \right)\\
&+ \gamma_i\left(B_{ik}^\alpha\right) - \mathcal{L}_{\theta_{ik}}\left(\gamma_i\left(w_{0i}^\alpha\right)\right) -\mathcal{L}_{\theta_{ik}}\left(\gamma_k\left(w_{0i}^\alpha\right)\right) -i_{\theta_{ik}}\left(\beta_i\left(w_{0i}^\alpha\right)\right) - i_{\theta_{ij}}\left(\beta_k\left(w_{0i}^\alpha\right)\right)\\
&- \gamma_j\left(B_{ik}^\alpha\right) + \gamma_k\left(B_{ij}^\alpha\right)+ \mathcal{L}_{\theta_{jk}}\left(\gamma_j\left(w_{0i}^\alpha\right)\right) + \mathcal{L}_{\theta_{jk}}\left(\gamma_k\left(w_{0i}^\alpha\right)\right) \\
&=\gamma_i\left(B_{ik}^\alpha\right) - \gamma_j\left(B_{ik}^\alpha\right) + i_{\theta_{ij}}\left( \beta_i\left(w_{0i}^\alpha\right) - \beta_k\left(w_{0i}^\alpha\right) - \mathcal{L}_{\theta_{ik}}\left(w_{0i}^\alpha\right) + \mathcal{L}_{\theta_{ik}}\left(w_{0i}^\alpha \right) \right) - \mathcal{L}_{\theta_{jk}}\left(\gamma_i\left(w_{0i}^\alpha \right) -\gamma_k\left(w_{0i}^\alpha \right)\right)\\
&- \gamma_i\left(B_{ij}^\alpha\right) + \gamma_k\left(B_{ij}^\alpha\right) + i_{\theta_{ik}}\left(\beta_j\left(w_{0i}^\alpha\right) - \beta_i\left(w_{0i}^\alpha\right)  + \mathcal{L}_{\theta_{ij}}\left(w_{0i}^\alpha\right) - \mathcal{L}_{\theta_{ij}}\left(w_{0i}^\alpha \right)\right) +\mathcal{L}_{\theta_{ik}}\left( \gamma_j\left( w_{0i}^\alpha\right) - \gamma_k\left(w_{0j}^\alpha \right) \right)\\
&= i_{\theta_{ij}}\mathcal{L}_{\theta_{ik}}\left(w_{0i}^\alpha \right) - \mathcal{L}_{\theta_{jk}}i_{\theta_{ik}}\left(w_{0i}^\alpha \right) - i_{\theta_{ik}}\mathcal{L}_{\theta_{ij}}\left(w_{0i}^\alpha \right) + \mathcal{L}_{\theta_{ik}}i_{\theta_{jk}}\left(w_{0i}^\alpha \right)\\
&=\left( i_{\theta_{ij}}\mathcal{L}_{\theta_{ij}} + i_{\theta_{ij}}\mathcal{L}_{\theta_{jk}} - \mathcal{L}_{\theta_{jk}}i_{\theta_{ij}} - \mathcal{L}_{\theta_{jk}} i_{\theta_{jk}} - i_{\theta_{ij}}\mathcal{L}_{\theta_{ij}} - i_{\theta_{jk}} \mathcal{L}_{\theta_{ij}} + \mathcal{L}_{\theta_{ij}} i_{\theta_{jk}} + \mathcal{L}_{\theta_{jk}} i_{\theta_{jk}}    \right)\left(w_{0i}^\alpha \right) = 2i_{\left[\theta_{ij}, \theta_{jk} \right] } \left( w_{0i}^\alpha \right)
\end{align*}
This completes the proof of Lemma \ref{ou50}.
\end{proof}
Then from $(\ref{itt17})$ and Lemma \ref{ou50}, we obtain
\begin{theorem}
Suppose given a compact foliated complex manifold $\left( M, \mathcal{N}_{\mathcal{F}_0}^* \right)$ with $\mathcal{N}_{\mathcal{F}_0}^*$ locally free, and $\left( \left\{ \beta_j \right\}, \left\{\gamma_j\right\}, \left\{ \left( \theta_{jk}, B_{jk}\right) \right\} \right)\in \mathbb{H}^1\left( M, \tilde{\mathcal{E}}_{\mathcal{F}_0}^\bullet \right)$. In order that there may exist an unfolded analytic family $\left( \mathcal{M}, \tilde{\mathcal{N}}_\mathcal{F}^* , B, \pi \right)$ such that $\pi^{-1}(0)=\left( M, \mathcal{N}_{\mathcal{F}_0}^* \right)$ and $\left( \frac{\partial\left( \mathcal{M}, \tilde{\mathcal{N}}_{\mathcal{F}}^* \right)}{\partial t} \right)|_{t=0}= \left( \left\{ \beta_j \right\}, \left\{ \gamma_j \right\}, \left\{ \left( \theta_{jk}, B_{jk} \right) \right\} \right) $, it is necessary that
\begin{align*}
&\left( \left\{ S_j \right\}, \left\{ Q_j \right\}, \left\{ G_{jk} \right\}, \left\{ C_{jk} \right\} \left\{\left[\left( \theta_{ij}, B_{ij} \right) , \left(\theta_{jk}, B_{jk} \right) \right]  \right\} \right) \\
&\in C^0\left( \mathcal{U}, \mathscr{H}om_{\mathcal{O}_M}\left( \mathcal{N}_{\mathcal{F}_0}^*, \tilde{\mathcal{S}}^2 \right) \right) \bigoplus C^0\left( \mathcal{U}, \mathscr{H}om_{\mathcal{O}_M}\left( \mathcal{N}_{\mathcal{F}_0}^*, \tilde{\mathcal{S}}^1 \right) \right) \bigoplus C^1\left( \mathcal{U} , \mathscr{H}om_{\mathcal{O}_M}\left( \mathcal{N}_{\mathcal{F}_0}^*, \Omega_M^1 \right) \right) \notag \\
& \bigoplus C^1\left( \mathcal{U}, \mathscr{H}om_{\mathcal{O}_M}\left( \mathcal{N}_{\mathcal{F}_0}^*, \mathcal{O}_M \right) \right) \bigoplus C^2\left( \mathcal{U} , \mathcal{E}_{\mathcal{N}_{\mathcal{F}_0}^*}\right) \notag
\end{align*}
from $(\ref{itt19})$ in \textnormal{Lemma \ref{ou50}} defines the $0$ class in $\mathbb{H}^2\left( M, \tilde{\mathcal{E}}_{\mathcal{N}_{\mathcal{F}_0}^*}^\bullet \right)$.
\end{theorem}

\section{Infinitesimal unfoldings and obstructions by using Dolbeault resolution of $\tilde{\mathcal{E}}_{\mathcal{N}_{\mathcal{F}_0}^*}^\bullet$}

Let $\left( \mathcal{M}, \tilde{\mathcal{N}}_\mathcal{F}^*, B, \pi \right)$ with $\tilde{\mathcal{N}}_\mathcal{F}^*$ locally free be an unfolded analytic family of unfoldings of a compact foliated complex manifold $\left( M, \mathcal{N}_{\mathcal{F}_0}^* \right)=\pi^{-1}(0)$ with $\mathcal{N}_{\mathcal{F}_0}^*$ locally free, and $B\subset \mathbb{C}^m$ with the coordinate $(t_1,...,t_m)$. Assume that given an open covering $\{U_j\}$ of $M$ by coordinate neighborhoods, $\Gamma\left( U_j, \mathcal{N}_{\mathcal{F}_0}^* \right)$ is generated by $w_{0j}^1,..., w_{0j}^q\in \Gamma\left( U_j, \Omega_M^1 \right)$ with $w_{0j}^\alpha= \sum_{\beta=1}^q h_{0jk}^{\alpha\beta} w_{0k}^\beta$ on $U_j\cap U_k$ for some $h_{0jk}^{\alpha\beta} \in \Gamma\left( U_{jk}, \mathcal{O}_M \right),\alpha=1,...,q$ and $w_{0j}^1 \wedge \cdots \wedge w_{0j}^q \wedge dw_{0j}^\alpha=0,\alpha=1,...,q$. In this subsection, we show that in the unfolded analytic family $\left( \mathcal{M}, \tilde{\mathcal{N}}_\mathcal{F}^*, B, \pi \right)$, the unfoldings near $\left( M, \mathcal{N}_{\mathcal{F}_0}^* \right)$ are represented by a section $\left\{\left(\varphi(t), Z_j(t) \right) \right\} \in A^{0,1}\left( M, \mathcal{E}_{\mathcal{N}_{\mathcal{F}_0}^*} \right)$ where $\varphi(t)\in A^{0,1}\left( M, \Theta_M \right)$ and $Z_j(t)\in \Gamma\left( U_j, \mathcal{A}^{0,1}\left( \mathscr{H}om_{\mathcal{O}_M}\left( \mathcal{N}_{\mathcal{F}_0}^*, \mathcal{N}_{\mathcal{F}_0}^* \right) \right) \right)$ which is defined by $Z_j(t)\left( w_{0j}^\alpha \right) = \sum_{\beta=1}^q Z_j^{\alpha\beta} w_{0j}^\beta$ for some $Z_j^{\alpha\beta}(t) \in \Gamma\left( U_j, \mathcal{A}^{0,1} \right),\alpha = 1,...,q$, and $Y(t)= \left\{Y_j(t) \right\}\in A^{0,0}\left( M, \mathscr{H}om_{\mathcal{O}_M}\left( \mathcal{N}_{\mathcal{F}_0}^*, \Omega_M^1 \right) \right)$ with $Y_j(t) \in \Gamma\left( U_j, \mathcal{A}^{0,0}\left(\mathscr{H}om_{\mathcal{O}_M} \left( \mathcal{N}_{\mathcal{F}_0}^*, \Omega_M^1 \right) \right) \right)$ and for each $t_\gamma, \gamma=1,...,m$ we define $V^\gamma (t)=\left\{ V_j^\gamma (t) \right\}\in A^{0,0}\left(M, \mathscr{H}om_{\mathcal{O}_M}\left( \mathcal{N}_{\mathcal{F}_0}^*, \mathcal{O}_M \right) \right)$ with $V_j^\gamma(t)\in \Gamma\left(U_j, \mathcal{A}^{0,0}\left(\mathscr{H}om_{\mathcal{O}_M}\left( \mathcal{N}_{\mathcal{F}_0}^*, \mathcal{O}_M \right) \right) \right)$ satisfying $\varphi(0)=0$ and $Z_j(0)=0$, and $Y(0)=\textnormal{Id}$ and the integrability conditions : by setting $ Y_j(t)\left(w_{0j}^\alpha\right)=Y_j^\alpha(t)\in \Gamma\left(U_j, \mathcal{A}^{0,0}\left(\Omega_M^1 \right) \right)$ and $Z_j(t)\left( w_{0j}^\alpha \right)= Z_j^\alpha(t) \in \Gamma\left( U_j, \mathcal{A}^{0,0}\left( \Omega_M^1 \right) \right),\alpha=1,...,q$ and $V_j^\gamma (t)\left( w_{0j}^\alpha \right)= V_j^{\alpha \gamma}(t) \in \Gamma\left(U_j, \mathcal{O}_M \right)$ and $V_j^\alpha (t):=\sum_{\gamma=1}^m V_j^{ \alpha\gamma}(t) dt_\gamma$ 
\begin{align*}
&\bar{\partial} \varphi(t) - \frac{1}{2} \left[ \varphi(t), \varphi(t) \right]= 0\\
&\left(Y_j^1(t) + V_j^1(t) \right) \wedge \cdots \wedge  \left( Y_j^q(t) + V_j^q(t) \right) \wedge \partial \left( Y_j^\alpha(t) + V_j^\alpha(t)\right) =0,\,\,\,\alpha=1,...,q\\
&\bar{\partial} Y_j^\alpha(t) - \mathcal{L}_{\varphi(t)} \left( Y_j^\alpha(t) \right) + \sum_{\eta=1}^q Z_j^{\alpha \eta}(t) Y_j^\eta (t)=0,\,\,\,\alpha=1,...,q \\
& \bar{\partial} \left( V_j^{\alpha \gamma}(t) \right) - \left[ \varphi (t) , V_j^{\alpha\gamma} (t) \right]- i_{\frac{\partial \varphi(t)}{\partial t_\gamma}}\left(Y_j(t) \right) + \sum_{\eta=1}^q Z_j^{\alpha \eta} (t) V_j^{\eta \gamma} (t)=0, \,\,\,\,\,\,\,\,\alpha=1,...,q, \,\,\, \gamma=1,...,m \\
&\bar{\partial} Z_j^{\alpha \delta}(t)- \left[ \varphi(t), Z_j^{\alpha \delta}(t)\right] + \sum_{\sigma=1}^q Z_j^{\alpha \sigma}(t) \wedge Z_j^{\sigma \delta}(t) = 0 , \,\,\,\,\,\alpha, \delta=1,...,q
\end{align*}
where $t\in \Delta$ a sufficiently small polydisk in $B$ and $\partial$ in the second condition is the differential $\partial_{M\times \Delta}$ on $M\times \Delta$ with $d_{M\times \Delta}=\partial_{M\times \Delta} + \bar{\partial}_{M\times \Delta}$, not on $M$. To deduce the integrability condition, we extend Kodaira's approach (\cite{Kod05} section \S 5.3) in the context of an unfolded analytic family.

\subsection{Integrability conditions} \

Let $\left(\mathcal{M}, \tilde{\mathcal{N}}_\mathcal{F}^* , B, \pi \right)$ with $\tilde{\mathcal{N}}_\mathcal{F}^*$ locally free be an unfolded analytic family, where $B$ is a domain of $\mathbb{C}^m$ with the coordinate $(t_1,..., t_m)$ containing the origin $0$ as in Definition \ref{u1}. Define $|t|=\max_\lambda |t_\lambda|$ for $t=(t_1,...,t_m)\in \mathbb{C}^m$, and let $\Delta=\{t\in \mathbb{C}^m | |t|<r\}$ the polydisk  of radius $r>0$. If we take a sufficiently  small $\Delta \subset B$, then $\left(\mathcal{M}_\Delta, \mathcal{N}_{\mathcal{F}_\Delta}^* \right)=\pi^{-1}\left(\Delta\right)$ is represented in the form
\begin{align*}
\left(\mathcal{M}_\Delta, \tilde{\mathcal{N}}_{\mathcal{F}_\Delta}^* \right)=\bigcup_j \left(U_j\times \Delta, \tilde{\mathcal{N}}_\mathcal{F}^*|_{U_j\times \Delta}\right)
\end{align*}
We denote a point of $U_j$ by $\xi_j=\left(\xi_j^1,...,\xi_j^n\right)$ and $\Gamma\left(U_j\times \Delta ,  \tilde{\mathcal{N}}_\mathcal{F}^* \right)$ is generated by $\omega_j^1\left(\xi_j,t\right),...,\omega_j^q\left(\xi_j,t\right)$, where
\begin{align*}
\omega_j^\alpha\left(\xi_j,t\right)=\sum_{\beta=1}^n w_j^{\alpha\beta}\left(\xi_j, t\right) d\xi_j^\beta + \sum_{\gamma=1}^m v_j^{\alpha \gamma}\left( \xi_j,t \right) dt_\gamma
\end{align*}
for some $w_j^{\alpha\beta}(\xi_j,t)\in \Gamma \left( U_j \times \Delta , \mathcal{O}_\mathcal{M} \right)$ and some $v_j^{\alpha \gamma}(\xi_j,t) \in \Gamma\left( U_j\times \Delta, \mathcal{O}_{\mathcal{M}} \right)$, and satisfies
\begin{align}\label{oui1}
\omega_j^1(\xi_j,t)\wedge \cdots \wedge \omega_j^q(\xi_j,t)\wedge d\omega _j^\alpha(\xi_j,t)=0,\,\,\,\,\,\,\,\alpha=1,...,q,\,\,\,\,\,\textnormal{on}\,\,\,\,\, U_j \times \Delta.
\end{align}
For simplicity, we assume that $U_j=\left\{\xi_j\in \mathbb{C}^m | |\xi_j|<1\right\}$ where $|\xi_j|=\max_a |\xi_j^a|$. $(\xi_j,t)\in U_j\times \Delta$ and $(\xi_k, t)\in U_k\times \Delta$ are the same point on $\mathcal{M}_\Delta$ if $\xi_j^\alpha=f_{jk}^\alpha(\xi_k,t),\alpha=1,...,n$ where $f_{jk}(\xi_k,t)$ is a holomorphic map of $\xi_k^1,..., \xi_k^n, t_1,...,t_m$, defined on $(U_k\times \Delta) \cap (U_j\times \Delta)$, and we have the following relation
\begin{align}\label{oui9}
\omega_j^\alpha(\xi_j,t)=\sum_{\beta=1}^q h_{jk}^{\alpha\beta}(\xi_k,t) \omega_k^\beta(\xi_k,t),\,\,\,\,\,\alpha=1,...,q
\end{align}
for a holomorphic function $h_{jk}^{\alpha\beta}(\xi_k,t)$ on $(U_k\times \Delta) \cap ( U_j\times \Delta)$. More precisely
\begin{align}
\sum_{\gamma=1}^n w_j^{\alpha\gamma}(f_{jk}(\xi_k,t),t)\frac{\partial f_{jk}^\gamma(\xi_k,t)}{\partial z_k^\xi} &= \sum_{\beta=1}^q h_{jk}^{\alpha\beta}(\xi_k,t) w_k^{\beta\xi}(\xi_k,t)\\
\sum_{\gamma=1}^n\sum_{\eta=1}^m w_j^{\alpha\gamma}(f_{jk}(\xi_k,t),t)\frac{\partial f_{jk}^\gamma(\xi_k,t)}{\partial t_\eta}dt_\eta + \sum_{\eta=1}^m v_j^{\alpha\eta}(f_{jk}(\xi_k,t),t)dt_\eta & = \sum_{\beta=1}^q \sum_{\eta=1}^m h_{jk}^{\alpha\beta}(\xi_k,t) v_k^{\beta\eta}(\xi_k,t)dt_\eta,\,\,\,\,\,\,\,\,\,\alpha=1,...,q. \label{oui3}
\end{align} 
and on $\left(U_i\times \Delta \right)\cap \left(U_j \times \Delta \right) \cap \left( U_k \times \Delta \right)$, we have
\begin{align}
h_{ik}^{\alpha \gamma}(\xi_k,t)= \sum_{\beta=1}^q h_{ij}^{\alpha\beta}\left(f_{jk}(\xi_k,t),t\right) h_{jk}^{\beta \gamma}(\xi_k,t),\,\,\,\,\,\alpha,\gamma=1,...,q
\end{align}

On the other hand, as in the subsection \ref{oci80} (we keep the notations in subsection \ref{oci80}), we have a diffeomorphism $\Psi$ of $M\times \Delta$ onto $\mathcal{M}_\Delta$ with $\pi \circ \Psi$ is the projection $M\times \Delta\to \Delta$.  Let $z=\left(z^1,..., z^n \right)$ be local coordinates of $M=M_0$, and for $\Psi(z,t)\in U_j\times \Delta$, put $\Psi\left( z, t\right) = \left(\xi_j^1(z,t), ..., \xi_j^n(z,t), t_1,..., t_m \right)$ where $\xi_j^\alpha=\xi_j^\alpha(z,t),\alpha=1,...,n$ is a $C^\infty$ function. We note that since $(z^1,..., z^n)$ and $\left(\xi_j^1(z,0), ..., \xi_j^n(z,0)\right)$ are local complex coordinates on $M=M_0$, $\xi_j^\alpha(z,0)$ are holomorphic functions of $z^1,...,z^n,\alpha=1,...,n$, and if we take $\Delta$ sufficiently small, we have
\begin{align}\label{oui4}
\det\left( \frac{\partial \xi_j^\alpha(z,t)}{\partial z^\lambda} \right)_{\alpha, \lambda=1,...,n} \ne 0\,\,\,\,\,\, \, \textnormal{for}\,\,\,t\in \Delta
\end{align}
 Since $\tilde{\mathcal{N}}_{\mathcal{F}}^*\cong \mathcal{N}_\mathcal{F}^*$, as in subsection \ref{oci80}, we have invertible matrices $g_j(z,t):=\left(g_j^{\alpha\beta}(z,t) \right)$ where $g_j^{\alpha\beta}(z,t)$ is a $C^\infty$ function on $\mathscr{U}_j=\Psi^{-1}\left(U_j\times \Delta \right)$ such that
 \begin{align}\label{oui5}
h_{jk}^{\alpha\beta}\left(\xi_k(z,t),t\right)= \sum_{\delta, \eta=1}^q g_{j}^{\alpha \delta}(z,t) h_{0jk}^{\delta \eta}(z) \left(g_k^{-1}\right)^{\eta\beta}(z,t),\,\,\,\,\,\,\,\alpha,\beta=1,...,q
\end{align}
where $h_{0jk}^{\alpha\beta}(z):=h_{jk}^{\alpha\beta}(\xi_j(z,0),0)$ and $g_j^{-1}(z,t)=\left( \left(  g_j^{-1}\right)^{\alpha\beta}(z,t) \right)_{\alpha, \beta=1,.., q} $ is the inverse of $g_j(z,t)$.

\subsection{Identification of deformations of complex structures with $\varphi(t)\in A^{0,1}\left(M, \Theta_M \right)$}\

We keep the notations in subsection \ref{oi51}, so that we have $\varphi(t)\in A^{0,1}\left( M, \Theta_M \right)$ such that
\begin{align}\label{oui11}
\varphi(0)=0,\,\,\,\,\,\,\,\bar{\partial}\varphi(t)- \frac{1}{2}\left[ \varphi(t), \varphi(t) \right]=0
\end{align}

\subsection{Identification of deformations of complex structures and locally free sheaves with $\left\{ \left( \varphi(t), Z_j(t)\right) \right\} \in A^{0,1}\left( M, \mathcal{E}_{\mathcal{N}_{\mathcal{F}_0}^*} \right)$} \

As in the same way in subsection \ref{oi52}, by taking $Z_j^{\alpha \delta}(t):= \sum_{\eta=1}^q \left( g_j^{-1}\right)^{\alpha \eta}\left( \bar{\partial} g_j^{\eta \delta} - \left[ \varphi, g_j^{\eta \delta}\right] \right)$, let us define an element $Z_j (t)\in \Gamma\left(U_j, \mathcal{A}^{0,1}\left(\mathscr{H}om_{\mathcal{O}_M}\left( \mathcal{N}_{\mathcal{F}_0}^*, \mathcal{N}_{\mathcal{F}_0}^* \right) \right)    \right)$ by
\begin{align*}
Z_j(t):\Gamma\left( U_j, \mathcal{N}_{\mathcal{F}_0}^* \right)&\to \Gamma\left( U_j, \mathcal{A}^{0,0}\left( \mathcal{N}_{\mathcal{F}_0}^* \right) \right) \\
  w_{0j}^\alpha &\mapsto \sum_{\delta=1}^q Z_j^{\alpha \delta} w_{0j}^\delta:=\sum_{\eta, \delta=1}^q \left(g_j^{-1} \right)^{\alpha \eta} \left(\bar{\partial} g_j^{\eta \delta}- \left[\varphi, g_j^{\eta \delta} \right] \right) w_{0j}^\delta
\end{align*}
and linearly extends to $\Gamma\left(U_j, \mathcal{N}_{\mathcal{F}_0}^* \right)$. Then $Z_j\left( w_{0j}^\alpha \right) = Z_k\left( w_{0j}^\alpha \right) + \sum_{\beta=1}^q \left[ \varphi, h_{0jk}^{\alpha \beta} \right] w_{0j}^\beta$, so that $\left\{\left(\varphi(t), Z_j(t) \right)    \right\}\in A^{0,1}\left( M, \mathcal{E}_{\mathcal{N}_{\mathcal{F}_0}^*} \right)$ such that $Z_j(0)=0$ and
\begin{align} \label{oui6}
\bar{\partial} Z_j^{\alpha \eta}- \left[\varphi, Z_j^{\alpha \eta} \right] + \sum_{\sigma=1}^q Z_j^{\alpha \sigma} \wedge Z_j^{\sigma \eta} = 0
\end{align}

\subsection{Identification of deformations of foliated structures with $Y(t)\in A^{0,0}\left( M, \mathscr{H}om_{\mathcal{O}_M}\left( \mathcal{N}_{\mathcal{F}_0}^*, \Omega_M^1 \right) \right)$ } \

We keep the notations in subsection \ref{oci81}. Let us define $W_j^{\alpha\gamma}(z,t)= \sum_{\beta=1}^n w_j^{\alpha\beta}\left( \xi_j(z,t), t\right) \frac{\partial \xi_j^\beta(z,t)}{\partial z^\gamma},\alpha=1,...,q$ and set $W_j^\alpha(z,t)= \sum_{\gamma=1}^n W_j^{\alpha \gamma}(z,t) dz^\gamma$. Let us define
\begin{align*}
Y_j(t):\Gamma\left(U_j, \mathcal{N}_{\mathcal{F}_0}^*\right)&\to \Gamma \left(U_j, \mathcal{A}^{0,0}\left(\Omega_{M}^1\right)\right) \\
              w_{0j}^\alpha &\mapsto Y_j^\alpha :=\sum_{\delta=1}^q \left(g_j^{-1} \right)^{\alpha\delta} W_j^{\delta}
\end{align*}
and linearly extends to $\Gamma\left(U_j, \mathcal{N}_{\mathcal{F}_0}^* \right)$. Then  $Y= \left\{Y_j \right\}$ defines a global section in $A^{0,0}\left( M, \mathscr{H}om_{\mathcal{O}_M}\left(\mathcal{N}_{\mathcal{F}_0}^* , \Omega_M^1 \right)   \right)$, and satisfies 
\begin{align}
\bar{\partial} \left(Y_j\left(w_{0j}^\alpha\right)\right) - & \mathcal{L}_{\varphi}\left( Y_j\left(w_{0j}^\alpha\right) \right) + Y_j \left(Z_j\left( w_{0j}^\alpha \right) \right)=0 \iff \bar{\partial} Y_j^\alpha - \mathcal{L}_{ \varphi}\left( Y_j^\alpha \right) + \sum_{\beta=1}^q Z_j^{\alpha \beta} Y_j^\beta = 0 \label{oui7} \\
&Y_j^1(z,t)\wedge \cdots \wedge Y_j^q(z,t) \wedge  \partial Y_j^\alpha(z,t)= 0,\,\,\,\,\,\alpha=1,...,q \label{oui8}
\end{align}

\subsection{Identification of unfolded structures with $V^\gamma(t)\in A^{0,0}\left( M , \mathscr{H}om_{\mathcal{O}_M} \left( \mathcal{N}_{\mathcal{F}_0}^*, \mathcal{O}_M     \right) \right) , \gamma=1,...,m$}\

From $(\ref{oui9})$, we have
\begin{align*}
\sum_{\beta=1}^n w_j^{\alpha\beta}(\xi_j,t)d\xi_j^\beta+\sum_{\gamma=1}^m v_j^{\alpha\gamma}(\xi_j,t)dt_\gamma=\sum_{\eta=1}^q h_{jk}^{\alpha\eta}(\xi_k,t)\left(\sum_{\xi=1}^n w_k^{\eta\xi}(\xi_k,t)d\xi_k^\xi+\sum_{\gamma=1}^m v_k^{\eta \gamma}(\xi_k,t)dt_\gamma \right)
\end{align*}
This implies that by considering the coefficient of $dt_\gamma$ and from $(\ref{oui5})$,
\begin{align*}
&v_j^{\alpha\gamma}(\xi_j(z,t),t)+ \sum_{\beta=1}^n w_j^{\alpha\beta}(\xi_j(z,t),t)\frac{\partial \xi_j^\beta}{\partial t_\gamma}= \sum_{\eta=1}^q h_{jk}^{\alpha\eta}(\xi_k(z,t),t)\left( v_k^{ \eta\gamma} (\xi_k(z,t),t) + \sum_{\beta=1}^n w_k^{\eta\beta}(\xi_k(z,t),t)\frac{\partial \xi_k^\beta}{\partial t_\gamma}     \right)\\
&=\sum_{\delta, \xi, \eta=1}^q g_j^{\alpha \xi}(z,t) h_{0jk}^{\xi \delta}(z) \left(g_k^{-1}\right)^{\delta \eta} (z,t) \left(  v_k^{ \eta\gamma} (\xi_k(z,t),t) + \sum_{\beta=1}^n w_k^{\eta\beta}(\xi_k(z,t),t)\frac{\partial \xi_k^\beta}{\partial t_\gamma}       \right)
\end{align*}
By multiplying $\sum_{\alpha=1}^q \left(g_j^{-1}\right)^{\sigma \alpha}$,  we have
\begin{align*}
\sum_{\alpha=1}^q(g_j^{-1})^{\sigma \alpha }\left(v_j^{\alpha \gamma}(\xi_j(z,t),t) + \sum_{\beta=1}^q w_j^{\alpha\beta}(\xi_k(z,t),t)\frac{\partial \xi_j^\beta}{\partial t_\gamma}   \right)=\sum_{\delta, \eta=1}^q h_{0jk}^{\sigma \delta} \left(g_k^{-1}\right)^{\delta \eta}\left(  v_k^{ \eta\gamma} (\xi_k(z,t),t) + \sum_{\beta=1}^n w_k^{\eta\beta}(\xi_k(z,t),t)\frac{\partial \xi_k^\beta}{\partial t_\gamma}       \right)
\end{align*}

We set 
\begin{align}\label{c2}
Q_j^{\sigma\gamma}(z,t):= v_j^{\sigma \gamma}(\xi_j(z,t),t) + \sum_{\beta=1}^n w_j^{\sigma\beta}(\xi_j(z,t),t)\frac{\partial \xi_j^\beta}{\partial t_\gamma}
\end{align}
Then
\begin{align*}
V_j^\gamma(t):\Gamma\left(U_j, \mathcal{N}_{\mathcal{F}_0}^* \right) &\to \Gamma\left(U_j, \mathcal{O}_M\right)\\
                 w_{0j}^\alpha &\mapsto V_j^{\alpha \gamma}(t):=\sum_{\sigma=1}^q \left( g_j^{-1}\right)^{\alpha \sigma } Q_j^{\sigma \gamma}(t)=\sum_{\sigma=1}^q \left(g_j^{-1}\right)^{\alpha \sigma }\left(v_j^{\sigma \gamma}(\xi_j(z,t),t) + \sum_{\beta=1}^q w_j^{\sigma \beta}(\xi_k(z,t),t)\frac{\partial \xi_j^\beta}{\partial t_\gamma}   \right)
\end{align*}
defines a global section $V^\gamma=\{V_j^\gamma\}\in A^{0,0}\left( M, \mathscr{H}om_{\mathcal{O}_M}\left( \mathcal{N}_{\mathcal{F}_0}^*, \mathcal{O}_M \right)\right)$. We note that since $\bar{\partial} Q_j^{\sigma \gamma} -\left[ \varphi, Q_j^{\sigma \gamma} \right]-   i_{\frac{\partial \varphi}{\partial t}} \left( W_j^{\sigma} \right)=0$ from Part III, we have
{\small{\begin{align*}
&\bar{\partial}\left(V_j^{\alpha \gamma}\right) - \left[\varphi, V_j^{\alpha \gamma} \right] - i_{\frac{\partial\varphi}{\partial t_\gamma}}\left( Y_j^\alpha \right)\\
&= \sum_{\sigma=1}^q \bar{\partial} \left(g_j^{-1}\right)^{\alpha \sigma}\cdot \left(v_j^{\sigma \gamma}(\xi_j(z,t),t) + \sum_{\beta=1}^q w_j^{\sigma\beta}(\xi_k(z,t),t)\frac{\partial \xi_j^\beta}{\partial t_\gamma}   \right) - \sum_{\sigma=1}^q \left[\varphi, \left(g_j^{-1}\right)^{\alpha \sigma}   \right]\cdot \left( v_j^{\sigma \gamma}(\xi_j(z,t),t) + \sum_{\beta=1}^q w_j^{\sigma\beta}(\xi_k(z,t),t)\frac{\partial \xi_j^\beta}{\partial t_\gamma}  \right)\\
&=- \sum_{ \sigma, \delta, \eta =1}^q \left(g_j^{-1}\right)^{\alpha \delta}\left(\bar{\partial}  g_j^{\delta \eta} - \left[\varphi, g_j^{\delta \eta} \right]  \right) \left( g_j^{-1}\right)^{\eta \sigma} \cdot  \left( v_j^{\sigma \gamma}(\xi_j(z,t),t) + \sum_{\beta=1}^q w_j^{\sigma\beta}(\xi_k(z,t),t)\frac{\partial \xi_j^\beta}{\partial t_\gamma}  \right)\\
&=- \sum_{\delta, \sigma, \eta=1}^q \left( g_j^{-1}\right)^{\alpha \delta}\left( \bar{\partial} g_j^{\delta \eta}- \left[ \varphi, g_j^{\delta \eta} \right] \right) \cdot V_j^{\eta \gamma}
\end{align*}}}

Hence we have
\begin{align}\label{oui13}
\bar{\partial}\left(V_j^{\alpha \gamma}\right) - \left[\varphi, V_j^{\alpha \gamma} \right] - i_{\frac{\partial\varphi}{\partial t_\gamma}}\left( Y_j^\alpha \right) + \sum_{ \sigma=1}^q Z_j^{\alpha \sigma} \cdot V_j^{\sigma \gamma}=0 ,\,\,\,\,\,\,\,\,\alpha =1,...,q,\gamma=1,...,m
\end{align}

Since $\xi_j^\alpha(z,t),\alpha=1,...,n$ is holomorphic in $t_1,...,t_m$, we see that
\begin{align*}
&\omega_j^\alpha(\xi_j,t)= w_j^\alpha(\xi_j,t)+ v_j^\alpha(\xi_j,t)= W_j^\alpha(z,t) + \sum_{\gamma=1}^m Q_j^{\alpha \gamma}(z,t)dt_\gamma  + W_j'^\alpha(z,t) \\
&= \sum_{\beta,\eta=1}^n w_j^{\alpha \beta}\left( \xi_j(z,t),t  \right)\frac{\partial \xi_j^\beta}{\partial z^\eta} dz^\eta + \sum_{\gamma=1}^m \left( v_j^{\alpha \gamma}\left(\xi_j(z,t),t) \right) + \sum_{\beta=1}^n w_j^{\alpha\beta}\left(\xi_j(z,t),t) \right)\frac{\partial \xi_j^\beta}{\partial t_\gamma} \right)dt_\gamma + \sum_{\beta=1}^n w_j^{\alpha\beta} \left(\xi_j(z,t), t \right)\frac{\partial \xi_j^\beta}{\partial \bar{z}^\gamma} dz^\gamma
\end{align*}
By setting $Q_j^\alpha(z,t):=\sum_{\gamma=1}^m Q_j^{\alpha \gamma}(z,t)dt_\gamma$, we note that
\begin{align*}
0&=\omega_j^1(\xi_j,t)\wedge \cdots \wedge \omega_j^q(\xi_j,t) \wedge d \omega_j^\alpha(\xi_j,t)\\
 &= \left( W_j^1 + Q_j^1 + W_j'^1 \right) \wedge \cdots  \wedge \left(  W_j^q + Q_j^q + W_j'^q\right) \wedge d \left( W_j^\alpha + Q_j^\alpha + W_j'^\alpha \right)\\
 &=  \left( W_j^1 + Q_j^1 \right) \wedge \cdots  \wedge \left(  W_j^q + Q_j^q \right) \wedge \partial \left( W_j^\alpha + Q_j^\alpha  \right) + \cdots
\end{align*}
Here $\partial$ is the differential $\partial_{M\times \Delta}$ with $d= \partial_{M\times\Delta} + \bar{\partial}_{M\times \Delta}$ and $ \left( W_j^1 + Q_j^1 \right) \wedge \cdots  \wedge \left(  W_j^q + Q_j^q \right) \wedge \partial \left( W_j^\alpha + Q_j^\alpha  \right) $ only involves $dz^\xi, dt_\eta$ and no $d\bar{z}^\xi, d\bar{t}_\eta$, and other terms contains at least one $d\bar{z}^\xi$ and $d\bar{t}_\eta$. Hence we have
\begin{align}\label{oui10}
 \left( W_j^1 + Q_j^1 \right) \wedge \cdots  \wedge \left(  W_j^q + Q_j^q \right) \wedge \partial \left( W_j^\alpha + Q_j^\alpha  \right) = 0
\end{align}
We note that since $\partial \left(Y_j^\alpha + V_j^\alpha \right) = \sum_{\beta=1}^q \partial \left( g_j^{-1} \right)^{\alpha \beta} \wedge \left( W_j^\beta + Q_j^\beta \right) + \sum_{\beta=1}^q \left( g_j^{-1} \right)^{\alpha \beta }  \partial \left(W_j^\beta + Q_j^\beta \right)$, $(\ref{oui10})$ implies
\begin{align}  \label{oui12}
 \left( Y_j^1 + V_j^1 \right) \wedge \cdots  \wedge \left(  Y_j^q + V_j^q \right) \wedge \partial \left( Y_j^\alpha + V_j^\alpha  \right) = 0,\,\,\,\,\,\,\alpha=1,...,q
\end{align}

\subsection{Identification of Infinitesimal deformations by using Dolbeault resolution of $\tilde{\mathcal{E}}_{\mathcal{N}_{\mathcal{F}_0}^*}^\bullet$}\

In the previous subsection we showed from $(\ref{oui11}),(\ref{oui6}),(\ref{oui7}), (\ref{oui13})$ and $(\ref{oui12})$ that given an unfolded analytic family $\left(\mathcal{M}, \tilde{\mathcal{N}}_\mathcal{F}^*, B, \pi \right)$ of unfoldings of $\left(M, \mathcal{N}_{\mathcal{F}_0}^* \right)=\pi^{-1}(0)$, the deformations near $\left( M, \mathcal{N}_{\mathcal{F}_0}^* \right)$ are represented by a section $\left\{\left(\varphi(t), Z_j(t) \right) \right\} \in A^{0,1}\left( M, \mathcal{E}_{\mathcal{N}_{\mathcal{F}_0}^*} \right)$ where $\varphi(t)\in A^{0,1}\left( M, \Theta_M \right)$ and $Z_j(t)\in \Gamma\left( U_j, \mathcal{A}^{0,1}\left( \mathscr{H}om_{\mathcal{O}_M}\left( \mathcal{N}_{\mathcal{F}_0}^*, \mathcal{N}_{\mathcal{F}_0}^* \right) \right) \right)$ which is defined by $Z_j(t)\left( w_{0j}^\alpha \right) = \sum_{beta=1}^q Z_j^{\alpha\beta} w_{0j}^\beta$ for some $Z_j^{\alpha\beta}(t) \in \Gamma\left( U_j, \mathcal{A}^{0,1} \right),\alpha = 1,...,q$, and $Y(t)= \left\{Y_j(t) \right\}\in A^{0,0}\left( M, \mathscr{H}om_{\mathcal{O}_M}\left( \mathcal{N}_{\mathcal{F}_0}^*, \Omega_M^1 \right) \right)$ with $Y_j(t) \in \Gamma\left( U_j, \mathcal{A}^{0,0}\left(\mathscr{H}om_{\mathcal{O}_M} \left( \mathcal{N}_{\mathcal{F}_0}^*, \Omega_M^1 \right) \right) \right)$ and $V^\gamma (t)=\left\{ V_j^\gamma (t) \right\}\in A^{0,0}\left(M, \mathscr{H}om_{\mathcal{O}_M}\left( \mathcal{N}_{\mathcal{F}_0}^*, \mathcal{O}_M \right) \right)$ with $V_j^\gamma(t)\in \Gamma\left(U_j, \mathcal{A}^{0,0}\left(\mathscr{H}om_{\mathcal{O}_M}\left( \mathcal{N}_{\mathcal{F}_0}^*, \mathcal{O}_M \right) \right) \right),\gamma=1,...,m$ satisfying $\varphi(0)=0$ and $Z_j(0)=0$, and $Y(0)=\textnormal{Id}$ and the integrability conditions : by setting $ Y_j(t)\left(w_{0j}^\alpha\right)=Y_j^\alpha(t)\in \Gamma\left(U_j, \mathcal{A}^{0,0}\left(\Omega_M^1 \right) \right)$ and $Z_j(t)\left( w_{0j}^\alpha \right)= Z_j^\alpha(t) \in \Gamma\left( U_j, \mathcal{A}^{0,0}\left( \Omega_M^1 \right) \right),\alpha=1,...,q$ and $V_j^\gamma (t)\left( w_{0j}^\alpha \right)= V_j^{\alpha \gamma}(t) \in \Gamma\left(U_j, \mathcal{O}_M \right)$ and $V_j^\alpha (t):=\sum_{\gamma=1}^m V_j^{ \alpha\gamma}(t) dt_\gamma$ 
\begin{align}
&\bar{\partial} \varphi(t) - \frac{1}{2} \left[ \varphi(t), \varphi(t) \right]= 0 \label{oui14}\\
& Y_j^1(t) \wedge \cdots \wedge Y_j^q(t)\wedge \partial Y_j^\alpha(t) = 0 ,\,\,\,\,\,\alpha=1,...,q \label{oui16} \\
&\left(Y_j^1(t) + V_j^1(t) \right) \wedge \cdots \wedge  \left( Y_j^q(t) + V_j^q(t) \right) \wedge \partial \left( Y_j^\alpha(t) + V_j^\alpha(t)\right) =0,\,\,\,\alpha=1,...,q \label{oui17}\\
&\bar{\partial} Y_j^\alpha(t) - \mathcal{L}_{\varphi(t)} \left( Y_j^\alpha(t) \right) + \sum_{\eta=1}^q Z_j^{\alpha \eta}(t) Y_j^\eta (t)=0,\,\,\,\alpha=1,...,q \label{oui18} \\
& \bar{\partial} \left( V_j^{\alpha \gamma}(t) \right) - \left[ \varphi (t) , V_j^{\alpha\gamma} (t) \right]- i_{\frac{\partial \varphi(t)}{\partial t_\gamma}}\left(Y_j^\alpha (t) \right) + \sum_{\eta=1}^q Z_j^{\alpha \eta} (t) V_j^{\eta \gamma} (t)=0, \,\,\,\,\,\,\,\,\alpha=1,...,q, \gamma=1,...,m  \label{oui19} \\
&\bar{\partial} Z_j^{\alpha \delta}(t)- \left[ \varphi(t), Z_j^{\alpha \delta}(t)\right] + \sum_{\sigma=1}^q Z_j^{\alpha \sigma}(t) \wedge Z_j^{\sigma \delta}(t) = 0 , \,\,\,\,\,\alpha, \delta=1,...,q \label{oui15}
\end{align}
where $t\in \Delta$ a sufficiently small polydisk in $B$ and $\partial$ in the second condition is the differential $\partial_{M\times \Delta}$ on $M\times \Delta$ with $d_{M\times \Delta}=\partial_{M\times \Delta} + \bar{\partial}_{M\times \Delta}$, not on $M$. We identify infinitesimal deformations at $0$ by taking the derivative of $(\ref{oui14}),(\ref{oui16}),(\ref{oui18})$ and $(\ref{oui15})$ with respect to $t$ and setting $t=0$, and by considering the coefficient of $t_\gamma$ in $(\ref{oui17})$ and setting $t=0$, and by plugging $t=0$ in $(\ref{oui19})$. Let us apply $\frac{\partial}{\partial t}=\sum_{\gamma=1}^m c_\gamma \frac{\partial}{\partial t_\gamma }$ to $(\ref{oui14})-(\ref{oui15})$. In subsection \ref{oci83}, $(\ref{oui14}),(\ref{oui16}),(\ref{oui18})$ and $(\ref{oui15})$ have been studied (we keep the notations in subsection \ref{oci83}) and we have
\begin{align}
&\bar{\partial}\left(\frac{\partial \varphi}{\partial t}|_{t=0} \right)=0 \label{oui22}\\
&\sum_{\beta=1}^q w_{0j}^1\wedge \cdots \wedge \overbrace{\frac{\partial Y_j^\beta}{\partial t}|_{t=0} }^{\beta- \textnormal{th}} \wedge \cdots \wedge w_{0j}^q \wedge dw_{0j}^q + w_{0j}^1\wedge \cdots \wedge w_{0j}^q \wedge \partial \left( \frac{\partial Y_j^\alpha}{\partial t}|_{t=0} \right) = 0\\
&\bar{\partial}\left(\frac{\partial Y_j^\alpha}{\partial t}|_{t=0} \right) - \mathcal{L}_{\frac{\partial\varphi}{\partial t}|_{t=0}}\left( w_{0j}^\alpha \right)  + \sum_{\beta=1}^q \frac{\partial Z_j^{\alpha \eta}}{\partial t}|_{t=0} w_{0j}^\eta = 0\\
&\bar{\partial}\left( \frac{\partial Z_j^{\alpha \delta}}{\partial t}|_{t=0} \right) = 0 \label{oui23}
\end{align}
and $\left(\frac{\partial Y}{\partial t}|_{t=0}, \left\{ \left(\frac{\partial \varphi}{\partial t}|_{t=0}, \frac{\partial Z_j}{\partial t}|_{t=0} \right)\right\} \right) \in A^{0,0}\left( M, \mathscr{H}om_{\mathcal{O}_M}\left(\Theta_{\mathcal{F}_0}, \Omega_M^1 \right)\right)\bigoplus A^{0,1}\left(M, \mathcal{E}_{\mathcal{N}_{\mathcal{F}_0}^*} \right) $

Let us consider the condition $(\ref{oui17})$. By considering the coefficient of $t_\gamma$ in
\begin{align*}
\left(Y_j^1+ \sum_{\gamma=1}^r V_j^{1\gamma} dt_\gamma\right)\wedge \cdots \wedge \left( Y_j^q + \sum_{\gamma=1}^r V_j^{q\gamma} dt_\gamma \right) \wedge \partial_{M\times \Delta} \left( Y_j^\alpha + \sum_{\gamma=1}^r V_j^{\alpha \gamma} dt_\gamma  \right) = 0 , 
\end{align*}
we have
\begin{align}
 \sum_{\beta=1}^q (-1)^{q-\beta} V_j^{\beta \gamma} \cdot  Y_j^1\wedge \cdots \wedge \widehat{Y}_j^\beta \wedge \cdots \wedge Y_j^q \wedge \partial_z Y_j^\alpha + Y_j^1\wedge \cdots \wedge Y_j^q\wedge \partial_z V_j^{\alpha \gamma} - Y_j^1\wedge \cdots \wedge Y_j^q \wedge \frac{\partial Y_j^\alpha}{\partial t_\gamma}=0
\end{align}
Then by multiplying $\sum_{\gamma=1}^m c_\gamma$, we have
{\small{\begin{align}
&\sum_{\beta=1}^q (-1)^{q-\beta}\left( \sum_{\gamma=1}^r c_\gamma V_j^{\beta \gamma}     \right)\cdot Y_j^1 \wedge \cdots \wedge \widehat{Y}_j^\beta \wedge \cdots \wedge Y_j^q \wedge \partial_z Y_j^\alpha + Y_j^1\wedge \cdots \wedge Y_j^q \wedge \partial_z \left( \sum_{\gamma=1}^r c_\gamma V_j^{\alpha \gamma} \right) - Y_j^1 \wedge \cdots \wedge Y_j^q\wedge \left( \sum_{\gamma=1}^r c_\gamma \frac{\partial Y_j^\alpha}{\partial t_\gamma} \right)\notag \\
&=\sum_{\beta=1}^q (-1)^{q-\beta} i_{\frac{\partial}{\partial t}}\left( V_j^\beta \right)\cdot Y_j^1 \wedge \cdots \wedge \widehat{Y}_j^\beta \wedge \cdots \wedge Y_j^q \wedge d_z Y_j^\alpha + Y_j^1 \wedge \cdots \wedge Y_j^q \wedge \partial_z \left( i_{\frac{\partial}{\partial t}}\left( V_j^\alpha \right) \right)- Y_j^1 \wedge \cdots \wedge Y_j^q\wedge \frac{\partial Y_j^\alpha}{\partial t} =0\label{oui25}
\end{align}}}
By setting $t=0$, we have
{\small{\begin{align}\label{oui20}
\sum_{\beta=1}^q (-1)^{q-\beta} i_{\frac{\partial}{\partial t}}\left( V_j^\beta \right)|_{t=0} w_{0j}^1 \wedge \cdots \wedge \widehat{w}_{0j}^\beta \wedge \cdots \wedge w_{0j}^q\wedge dw_{0i}^\alpha + w_{0i}^1\wedge \cdots \wedge w_{0i}^q \wedge \partial_z \left( i_{\frac{\partial}{\partial t}}\left( V_j^\alpha \right)|_{t=0}  \right) - w_{0i}^1 \wedge \cdots \wedge w_{0i}^q \wedge \frac{\partial Y_j^\alpha}{\partial t}|_{t=0}=0
\end{align}}}
and by setting $t=0$ in $(\ref{oui19})$, we have
\begin{align*}
\bar{\partial}\left( V_j^{\alpha \gamma}(0) \right) - i_{\frac{\partial \varphi}{\partial t_\gamma }|_{t=0}}\left( w_{0j}^\alpha \right)=0,\,\,\,\,\,\,\gamma=1,...,q
\end{align*}
By multiplying $\sum_{\gamma=1}^m c_\gamma$, we have
\begin{align} \label{oui21}
\bar{\partial}\left( \sum_{\gamma=1}^m c_\gamma V_j^{\alpha \gamma}(0) \right) - i_{\sum_{\gamma=1}^m c_\gamma \frac{\partial \varphi}{\partial t_\gamma }|_{t=0}}\left( w_{0j}^\alpha \right)=0 \Longrightarrow  \bar{\partial}\left( i_{\frac{\partial}{\partial t}}\left(V_j^\alpha \right) |_{t=0}  \right) - i_{ \frac{\partial \varphi}{\partial t }|_{t=0}}\left( w_{0j}^\alpha \right)=0
\end{align}

Let us define $i_{\frac{\partial}{\partial t}}\left(V(t) \right)|_{t=0} \in A^{0,0}\left( M, \mathscr{H}om_{\mathcal{O}_M}\left( \mathcal{N}_{\mathcal{F}_0}^*, \mathcal{O}_M \right) \right)$ by
\begin{align*}
i_{\frac{\partial}{\partial t}}\left(V_j(t) \right)|_{t=0} : \Gamma\left( U_j, \mathcal{N}_{\mathcal{F}_0}^* \right) &\to \Gamma\left( U_j, \mathcal{A}^{0,0} \right) \\
                     w_{0j}^\alpha &\mapsto i_{\frac{\partial}{\partial t}}\left( V_j^\alpha \right)|_{t=0}
\end{align*}
and linearly extends to $\Gamma\left( U_j, \mathcal{N}_{\mathcal{F}_0}^* \right)$. Then $i_{\frac{\partial}{\partial t}}\left(V(t) \right)|_{t=0} : = \left\{ i_{\frac{\partial}{\partial t}}\left( V_j(t) \right)|_{t=0} \right\}$ defines a global section in $A^{0,0}\left( M, \mathscr{H}om_{\mathcal{O}_M} \left( \mathcal{N}_{\mathcal{F}_0}^*, \mathcal{O}_M \right) \right)$. Then $(\ref{oui22})-(\ref{oui23}), (\ref{oui20})$ and $(\ref{oui21})$ implies that 
{\small{\begin{align*}
\left( \frac{\partial Y}{\partial t}|_{t=0}, i_{\frac{\partial}{\partial t}}\left(V \right)|_{t=0}  , \left\{ \left( \frac{\partial \varphi}{\partial t}|_{t=0}, \frac{\partial Z_j}{\partial t}|_{t=0} \right) \right\}\right) \in A^{0,0}\left( M, \mathscr{H}om_{\mathcal{O}_M}\left( \mathcal{N}_{\mathcal{F}_0}^* \right) \right) \bigoplus A^{0,0}\left( M,  \mathscr{H}om_{\mathcal{O}_M}\left( \mathcal{N}_{\mathcal{F}_0}^*, \mathcal{O}_M \right) \right) \bigoplus A^{0,1}\left( M, \mathcal{E}_{\mathcal{N}_{\mathcal{F}_0}^*} \right)
\end{align*}}}
defines a $1$-cocycle in the following double complex
{\Tiny{\begin{equation} \label{oui24}
\begin{CD}
\cdots\\
@AAA \\
A^{0,0}\left( M-S, \left(\mathcal{N}_{\mathcal{F}_0}^* \right)^* \otimes \tilde{\mathcal{S}}^3 \right) \bigoplus A^{0,0}\left( M-S, \left( \mathcal{N}_{\mathcal{F}_0}^* \right)\otimes \tilde{\mathcal{S}}^2 \right)@>-(-1)^q\bar{\partial}>> A^{0,1}\left( M-S, \left(\mathcal{N}_{\mathcal{F}_0}^* \right)^* \otimes \tilde{\mathcal{S}}^3 \right) \bigoplus A^{0,1}\left( M-S, \left( \mathcal{N}_{\mathcal{F}_0}^* \right)\otimes \tilde{\mathcal{S}}^2 \right)@>(-1)^q\bar{\partial}>> \cdots \\
@AAA @AAA  \\
A^{0,0}\left( M-S, \left(\mathcal{N}_{\mathcal{F}_0}^* \right)^* \otimes \tilde{\mathcal{S}}^2 \right) \bigoplus A^{0,0}\left( M - S , \left( \mathcal{N}_{\mathcal{F}_0}^* \right)\otimes \tilde{\mathcal{S}}^1 \right) @>(-1)^q \bar{\partial}>>A^{0,1}\left( M-S, \left(\mathcal{N}_{\mathcal{F}_0}^* \right)^* \otimes \tilde{\mathcal{S}}^2 \right) \bigoplus A^{0,1}\left( M-S, \left( \mathcal{N}_{\mathcal{F}_0}^* \right)\otimes \tilde{\mathcal{S}}^1 \right) @>-(-1)^q\bar{\partial} >> \cdots \\
@AAA @AAA \\
A^{0,0}\left( M,  \left(\mathcal{N}_{\mathcal{F}_0}^*\right)^* \otimes \Omega_M^1 \right) \bigoplus A^{0,0}\left(M, \left(\mathcal{N}_{\mathcal{F}_0}^* \right)^* \right) @>-\bar{\partial}>>  A^{0,1}\left( M, \left( \mathcal{N}_{\mathcal{F}_0}^*\right)^* \otimes \Omega_M^1 \right) \bigoplus A^{0,1}\left(M,\left( \mathcal{N}_{\mathcal{F}_0}^*\right)^* \right) @>\bar{\partial}>> \cdots\\
@AAA @AAA @AAA \\
A^{0,0}\left(M, \mathcal{E}_{\mathcal{N}_{\mathcal{F}_0}^*}   \right) @>\bar{\partial}>> A^{0,1}\left( M, \mathcal{E}_{\mathcal{N}_{\mathcal{F}_0}^*} \right) @>-\bar{\partial}>> A^{0,2}\left( M, \mathcal{E}_{\mathcal{N}_{\mathcal{F}_0}^*} \right) @>\bar{\partial}>> \cdots
\end{CD}
\end{equation}}}

We will denote the $i$-th cohomology group of the above double complex $(\ref{oui24})$ by $\textnormal{C}^i$. In the next section \ref{ni3}, we will show that we have
\begin{align}\label{b3}
\textnormal{C}^i\cong \mathbb{H}^i\left( M, \tilde{\mathcal{E}}_{\mathcal{N}_{\mathcal{F}_0}^*}^\bullet  \right),\,\,\,\,\,i=0,1,2
\end{align}

Then we have the unfolded Kodaira-Spencer map
\begin{align} \label{b1}
\varphi_0: T_0 B \to &\mathbb{H}^1\left(  M, \tilde{\mathcal{E}}_{\mathcal{N}_{\mathcal{F}_0}^*}^\bullet  \right) \\
\frac{\partial}{\partial t} \mapsto &\left( \frac{\partial Y}{\partial t}|_{t=0}, i_{\frac{\partial}{\partial t}}\left(V \right)|_{t=0}  , \left\{ \left( \frac{\partial \varphi}{\partial t}|_{t=0}, \frac{\partial Z_j}{\partial t}|_{t=0} \right) \right\}\right) \notag \\
& \in A^{0,0}\left( M, \mathscr{H}om_{\mathcal{O}_M}\left( \mathcal{N}_{\mathcal{F}_0}^* , \Omega_M^1 \right) \right) \bigoplus A^{0,0}\left( M,  \mathscr{H}om_{\mathcal{O}_M}\left( \mathcal{N}_{\mathcal{F}_0}^*, \mathcal{O}_M \right) \right) \bigoplus A^{0,1}\left( M, \mathcal{E}_{\mathcal{N}_{\mathcal{F}_0}^*} \right) \notag
\end{align}

\subsection{Identification of obstructions by using Dolbeault resolution of $\tilde{\mathcal{E}}_{\mathcal{N}_{\mathcal{F}_0}^*}^\bullet$} \

We assume that the unfolded analytic family $\left( \mathcal{M}, \tilde{\mathcal{N}}_\mathcal{F}^*, B, \pi \right)$ in the previous subsection is a $1$-parameter family, i.e. $\Delta \subset \mathbb{C}$. Then as in subsection \ref{oci84}, we have $(\ref{it31})-(\ref{it32})$. It remains to take the derivative of $(\ref{oui25})$ and $(\ref{oui19})$ with respect to $t$ (here we are considering $1$-parameter family),
\begin{align}
&\sum_{\beta=1}^q (-1)^{q-\beta} i_{\frac{\partial}{\partial t}}\left( V_j^\beta \right)\cdot Y_j^1 \wedge \cdots \wedge \widehat{Y}_j^\beta \wedge \cdots \wedge Y_j^q \wedge \partial_z Y_j^\alpha + Y_j^1 \wedge \cdots \wedge Y_j^q \wedge \partial_z \left( i_{\frac{\partial}{\partial t}}\left( V_j^\alpha \right) \right)- Y_j^1 \wedge \cdots \wedge Y_j^q\wedge \frac{\partial Y_j^\alpha}{\partial t} =0 \label{oui26}\\
 & \bar{\partial} \left( i_{\frac{\partial}{\partial t}}\left(V_j^{\alpha}\right) \right) - \left[ \varphi  , i_{\frac{\partial}{\partial t}}\left( V_j^{\alpha}  \right) \right]- i_{\frac{\partial \varphi }{\partial t}}\left(Y_j^\alpha  \right) + \sum_{\eta=1}^q Z_j^{\alpha \eta} i_{\frac{\partial}{\partial t}}\left( V_j^{\eta} \right) =0, \,\,\,\,\,\,\,\,\alpha=1,...,q  \label{oui27}
\end{align}

Let us take the derivative of $(\ref{oui26})$ with respect to $t$ again:
\begin{align*}
&\sum_{\beta=1}^q (-1)^{q-\beta}\frac{\partial}{\partial t}\left( i_{\frac{\partial}{\partial t}}\left( V_j^\beta \right)   \right)\cdot  Y_j^1 \wedge \cdots \wedge \cdots \wedge \widehat{Y}_j^\beta \wedge \cdots \wedge Y_j^q \wedge \partial_z Y_j^\alpha\\
&+\sum_{\beta=1}^q \sum_{\eta\ne \beta}(-1)^{q-\beta}  i_{\frac{\partial}{\partial t}}\left( V_j^\beta \right)   \cdot  Y_j^1 \wedge \cdots \wedge \frac{\partial Y_j^\eta}{\partial t}\wedge \cdots \wedge \widehat{Y}_j^\beta \wedge \cdots \wedge Y_j^q \wedge \partial_z  Y_j^\alpha\\
&+ \sum_{\beta=1}^q (-1)^{q-\beta} i_{\frac{\partial}{\partial t}}\left( V_j^\beta \right)   \cdot  Y_j^1 \wedge \cdots \wedge \widehat{Y}_j^\beta \wedge \cdots \wedge Y_j^q \wedge \partial_z \left( \frac{\partial Y_j^\alpha}{\partial t} \right)\\
&+ \sum_{\beta=1}^q Y_j^1 \wedge \cdots \wedge \frac{\partial Y_j^\beta}{\partial t}\wedge \cdots \wedge Y_j^q \wedge \partial_z \left( i_{\frac{\partial}{\partial t}}\left( V_j^\alpha \right) \right) + Y_j^1 \wedge \cdots \wedge Y_j^q\wedge \partial_z \left( \frac{\partial}{\partial t}\left( i_{\frac{\partial}{\partial t}}\left( V_j^\alpha \right)   \right)    \right) \\
&- \sum_{\beta=1}^q Y_j^1 \wedge \cdots \wedge \frac{\partial Y_j^\beta}{\partial t} \wedge \cdots \wedge Y_j^q \wedge \frac{\partial Y_j^\alpha}{\partial t} - Y_j^1 \wedge \cdots \wedge Y_j^q \wedge \frac{\partial^2 Y_j^\alpha}{\partial t^2} = 0
\end{align*}

By setting $t=0$, we have
\begin{align}\label{itt21}
&\sum_{\beta=1}^q (-1)^{q-\beta}\frac{\partial}{\partial t}\left( i_{\frac{\partial}{\partial t}}\left( V_j^\beta \right)   \right)|_{t=0}\cdot  w_{0j}^1 \wedge \cdots \wedge \cdots \wedge \widehat{w}_{0j}^\beta \wedge \cdots \wedge w_{0j}^q \wedge d w_{0j}^\alpha + w_{0j}^1 \wedge \cdots \wedge w_{0j}^q\wedge \partial_z \left( \frac{\partial}{\partial t}\left( i_{\frac{\partial}{\partial t}}\left( V_j^\alpha \right) |_{t=0}  \right)    \right)\\
&- w_{0j}^1 \wedge \cdots \wedge w_{0j}^q \wedge \frac{\partial^2 Y_j^\alpha}{\partial t^2}|_{t=0} \notag \\
&=-\sum_{\beta=1}^q \sum_{\eta\ne \beta}(-1)^{q-\beta}  i_{\frac{\partial}{\partial t}}\left( V_j^\beta \right)|_{t=0}  \cdot  w_{0j}^1 \wedge \cdots \wedge \frac{\partial Y_j^\eta}{\partial t}|_{t=0} \wedge \cdots \wedge \widehat{w}_{0j}^\beta \wedge \cdots \wedge w_{0j}^q \wedge d w_{0j}^\alpha \notag \\
&- \sum_{\beta=1}^q (-1)^{q-\beta} i_{\frac{\partial}{\partial t}}\left( V_j^\beta \right) |_{t=0}   \cdot  w_{0j}^1 \wedge \cdots \wedge \widehat{w}_{0j}^\beta \wedge \cdots \wedge w_{0j}^q \wedge \partial_z \left( \frac{\partial Y_j^\alpha}{\partial t} |_{t=0 }\right) \notag \\
&- \sum_{\beta=1}^q w_{0j}^1 \wedge \cdots \wedge \frac{\partial Y_j^\beta}{\partial t}|_{t=0} \wedge \cdots \wedge w_{0j}^q \wedge \partial \left( i_{\frac{\partial}{\partial t}}\left( V_j^\alpha \right)|_{t=0} \right)  \notag \\
&+ \sum_{\beta=1}^q w_{0j}^1 \wedge \cdots \wedge \frac{\partial Y_j^\beta}{\partial t}|_{t=0} \wedge \cdots \wedge w_{0j}^q \wedge \frac{\partial Y_j^\alpha}{\partial t}|_{t=0}  \notag
\end{align}

By taking the derivative of $(\ref{oui27})$ with respect to $t$ and setting $t=0$, we have
\begin{align}\label{itt20}
-\bar{\partial}\left(  \frac{\partial}{\partial t }\left(i_{\frac{\partial}{\partial t}}\left( V_j^\alpha \right) \right)|_{t=0} \right) + i_{\frac{\partial^2 \varphi}{\partial t^2}|_{t=0}}\left( w_{0j}^\alpha \right)= - \left[ \frac{\partial \varphi}{\partial t}|_{t=0}, i_{\frac{\partial}{\partial t}}\left( V_j^\alpha \right)|_{t=0}  \right] - i_{\frac{\partial \varphi}{\partial t}|_{t=0}}\left( \frac{\partial Y_j^\alpha}{\partial t}|_{t=0} \right) + \sum_{\eta=1}^q \frac{\partial Z_j^{\alpha \eta}}{\partial t}|_{t=0} \cdot i_{\frac{\partial}{\partial t}}\left( V_j^\alpha \right)|_{t=0}
\end{align}

\begin{lemma}\label{oui28}
For any $1$-cocycle 
\begin{align*}
\left(\Psi=\{\Psi_j\}, \tau=\left\{ \tau_j\right\},  \left\{\left( \psi, k_j\right) \right\} \right)\in A^{0,0}\left(M, \mathscr{H}om_{\mathcal{O}_M}\left( \mathcal{N}_{\mathcal{F}_0}^*, \Omega_M^1 \right) \right) \bigoplus A^{0,0}\left(M, \mathscr{H}om_{\mathcal{O}_M}\left( \mathcal{N}_{\mathcal{F}_0}^* \mathcal{O}_M \right) \right) \bigoplus A^{0,1}\left( M, \mathcal{E}_{\mathcal{N}_{\mathcal{F}_0}^*} \right)
\end{align*}
where $\Psi_j\in \Gamma\left( U_j, \mathcal{A}^{0,0}\left(\mathscr{H}om_{\mathcal{O}_M}\left( \mathcal{N}_{\mathcal{F}_0}^*, \Omega_M^1 \right) \right) \right)$ and $\tau_j\in \Gamma\left( U_j, \mathscr{H}om_{\mathcal{O}_M}\left( \mathcal{N}_{\mathcal{F}_0}^*, \mathcal{O}_M \right) \right)$ and $\psi\in A^{0,1}\left( M, \Theta_M \right)$ and $k_j\in \Gamma\left( U_j, \mathscr{H}om_{\mathcal{O}_M}\left( \mathcal{N}_{\mathcal{F}_0}^*, \mathcal{N}_{\mathcal{F}_0}^* \right) \right) $ such that
\begin{align*}
\Psi_j: \Gamma\left( U_j, \mathcal{N}_{\mathcal{F}_0}^* \right) &\to \Gamma\left(U_j, \mathcal{A}^{0,0}\left( \Omega_M^1 \right) \right) \\
 w_{0j}^\alpha &\mapsto \Psi_j^\alpha
\end{align*}
with $\Psi_j=\Psi_k$ on $U_{jk}$, i.e. $\Psi_j^\alpha= \sum_{\beta=1}^q h_{0jk}^{\alpha\beta} \Psi_j^\beta$,  and 
\begin{align*}
\tau_j:\Gamma\left( U_j, \mathcal{N}_{\mathcal{F}_0}^* \right) &\to \Gamma\left( U_j, \mathcal{A}^{0,0} \right) \\
 w_{0j}^\alpha &\mapsto \tau_j^\alpha
\end{align*}
with $\tau_j= \tau_k$ on $U_{jk}$, and
\begin{align*}
k_j: \Gamma\left( U_j, \mathcal{N}_{\mathcal{F}_0}^* \right) &\to \Gamma\left( U_j, \mathcal{A}^{0,1}\left( \mathcal{N}_{\mathcal{F}_0}^* \right) \right) \\
w_{0j}^\alpha &\mapsto \sum_{\beta=1}^q k_j^{\alpha\beta} w_{0j}^\beta \,\,\,\,\,\,\,k_j^{\alpha\beta}\in \Gamma\left( U_j, \mathcal{A}^{0,1} \right)
\end{align*}
with $k_j\left( w_{0j}^\alpha \right)=k_k\left(w_{0j}^\alpha \right) + \sum_{\beta=1}^q \left[ \psi ,  h_{0jk}^{\alpha\beta} \right] w_{0k}^\beta , \alpha=1,...,q$ satisfying
\begin{enumerate}
\item $\sum_{\gamma=1}^q w_{0j}^1 \wedge \cdots \wedge \Psi_j^\gamma \wedge \cdots \wedge w_{0j}^q \wedge dw_{0j}^\alpha + w_{0j}^1 \wedge \cdots \wedge w_{0j}^q \wedge \partial \Psi_j^\alpha = 0 , \alpha=1,...,q$.
\item $\sum_{\gamma=1}^q (-1)^{q-\gamma} w_{0j}^1 \wedge \cdots \wedge \widehat{w}_{0j}^\gamma \wedge \cdots \wedge w_{0j}^q \wedge dw_{0j}^\alpha + w_{0j}^1 \wedge \cdots \wedge w_{0j}^q \wedge \partial \tau_j^\gamma - w_{0j}^1 \wedge \cdots \wedge w_{0j}^q \wedge \Psi_j^\alpha = 0 , \alpha=1,...,q$. 
\item $-\bar{\partial} \Psi_j\left( w_{0j}^\alpha \right) +\mathcal{L}_\psi\left( w_{0j}^\alpha \right)- \sum_{\beta=1}^q k_j^{\alpha\beta} w_{0j}^\beta, \alpha=1,...,q.$
\item $-\bar{\partial} \tau_j\left( w_{0j}^\alpha  \right) + i_\psi\left(w_{0j}^\alpha \right)=0, \alpha= 1 , ... , q. $
\item $-\bar{\partial} \psi=0$ and $-\bar{\partial} b_j^{\alpha\beta}=0, \alpha,\beta=1,..., q$.
\end{enumerate}
we define
\begin{align*}
l_j:\Gamma\left( U_j, \mathcal{N}_{\mathcal{F}_0}^* \right) &\to \Gamma\left( U_j, \mathcal{A}^{0,0}\left( \bigwedge^{q+2} \Omega_M^1 \otimes \left(\bigwedge^q \mathcal{N}_{\mathcal{F}_0}^* \right)^* \right) \right) \\
 w_{0j}^\alpha &\mapsto - 2\sum_{\gamma=1}^q w_{0j}^1 \wedge \cdots \wedge \Psi_j^\gamma \wedge \cdots \wedge w_{0j}^q \wedge \partial \Psi_j^\gamma - 2\sum_{1\leq \gamma < \eta \leq q} w_{0j}^1 \wedge \cdots \wedge \Psi_j^\gamma \wedge \cdots \wedge \Psi_j^\eta \wedge \cdots \wedge w_{0j}^q \wedge dw_{0j}^\alpha
\end{align*}
Then $\tilde{l}=\{l_j\}$ defines a global section in $A^{0,0}\left( M, \mathscr{H}om_{\mathcal{O}_M}\left(\mathcal{N}_{\mathcal{F}_0}^*,  \bigwedge^{q+2} \Omega_M^1 \otimes \left(\bigwedge^q \mathcal{N}_{\mathcal{F}_0}^* \right)^* \right)   \right)$ and the restriction $l$ of $\tilde{l}$ on $M-S$ defines a global section in $A^{0,0}\left( M-S, \mathscr{H}om_{\mathcal{O}_M}\left( \mathcal{N}_{\mathcal{F}_0}^*, \tilde{\mathcal{S}}^2 \right) \right)$, and define
\begin{align*}
\sigma_j :\Gamma&\left(U_j, \mathcal{N}_{\mathcal{F}_0}^* \right)\to \Gamma\left( U_j, \mathcal{A}^{0,0}\left( \bigwedge^{q+1} \Omega_M^1 \otimes \left(\bigwedge^q \mathcal{N}_{\mathcal{F}_0}^* \right)^* \right)  \right) \\
 w_{0j}^\alpha \mapsto &-\sum_{\gamma=1}^q\sum_{1\leq \beta \ne \gamma \leq q } (-1)^{q-\gamma} \tau_j^\gamma\cdot w_{0j}^1\wedge \cdots \wedge \Psi_j^\beta \wedge \cdots \wedge \hat{w}_{0j}^\gamma \wedge \cdots \wedge w_{0j}^q\wedge dw_{0j}^\alpha  \\
&-\sum_{\gamma=1}^q (-1)^{q-\gamma} \tau_j^\gamma \wedge \cdots \wedge \hat{w}_{0j}^\gamma \wedge \cdots w_{0j}^q \wedge \partial \left( \Psi_j^\alpha \right) - \sum_{\gamma=1}^q w_{0j}^1\wedge \cdots \wedge  \Psi_j^\gamma \wedge \cdots \wedge w_{0j}^q\wedge \partial \left(\tau_j^\alpha  \right) \\
& +  \sum_{\gamma=1}^q w_{0j}^1\wedge \cdots \wedge  \Psi_j^\gamma \wedge \cdots \wedge w_{0j}^q\wedge \Psi_j^\alpha
\end{align*}
Then $\tilde{\sigma}=\{\sigma_j\}$ defines a global section in $A^{0,0}\left( M , \mathscr{H}om_{\mathcal{O}_M}\left( \mathcal{N}_{\mathcal{F}_0}^*, \bigwedge^{q+1}  \Omega_M^1 \otimes \left( \mathcal{N}_{\mathcal{F}_0}^* \right)^* \right) \right)$ and the restriction $\sigma$ of $\tilde{\sigma}$ on $M-S$ defines a global section in $A^{0,0}\left(M-S, \mathscr{H}om_{\mathcal{O}_M}\left( \mathcal{N}_{\mathcal{F}_0}^*, \tilde{\mathcal{S}}^1 \right) \right)$, and we define $\kappa_i \in \Gamma\left( U_j, \mathcal{A}^{0,1}\left( \mathscr{H}om_{\mathcal{O}_M}\left( \mathcal{N}_{\mathcal{F}_0}^*, \Omega_M^1 \right) \right) \right)$ by
\begin{align*}
\kappa_i :  \Gamma\left( U_j , \mathcal{N}_{\mathcal{F}_0}^* \right) &\to  \Gamma\left( U_j, \mathcal{A}^{0,1} \left( \Omega_M^1 \right) \right) \\
 w_{0j}^\alpha &\mapsto - 2 \mathcal{L}_\psi \left( \Psi_j^\alpha \right) + 2 \sum_{\eta=1}^q k_j^{\alpha\eta} \wedge \Psi_j^\eta
\end{align*}
Then $\kappa=\{\kappa_i\}$ defines a global section in $A^{0,1}\left( M, \mathscr{H}om_{\mathcal{O}_M}\left( \mathcal{N}_{\mathcal{F}_0}^*, \Omega_M^1 \right) \right)$,  and define $\pi_j\in \Gamma\left( U_j, \mathcal{A}^{0,1}\left( \mathscr{H}om_{\mathcal{O}_M}\left( \mathcal{N}_{\mathcal{F}_0}^*, \mathcal{O}_M \right) \right) \right)$ by
\begin{align*}
\pi_j:\Gamma\left( U_j ,   \mathcal{N}_{\mathcal{F}_0}^*  \right) &\to \Gamma\left( U_j, \mathcal{A}^{0,1} \right) \\
w_{0j}^\alpha & \mapsto  \pi_j^\alpha:= -\left[ \psi, \tau_j^\alpha  \right] - i_{\psi}\left( \Psi_j^\alpha  \right) + \sum_{\eta=1}^q k_j^{\alpha \eta} \tau_j^\eta
\end{align*}
Then $\pi=\{\pi_j\}$ defines a global section in $A^{0,1}\left( M, \mathscr{H}om_{\mathcal{O}_M}\left( \mathcal{N}_{\mathcal{F}_0}^*, \mathcal{O}_M \right) \right)$ and we have $ - \left[\left(\psi, k_j \right), \left( \psi, k_j \right) \right] \in A^{0,2}\left(M, \mathcal{E}_{\mathcal{N}_{\mathcal{F}_0}^*} \right)$. Then
\begin{align}\label{itt22}
\left( l, \sigma, \kappa, \pi, \left\{ - \left[\left(\psi, k_j \right), \left(\psi, k_j \right) \right]   \right\} \right) \in &A^{0,0}\left( M-S, \mathscr{H}om_{\mathcal{O}_M}\left(\mathcal{N}_{\mathcal{F}_0}^* , \tilde{\mathcal{S}}^2 \right) \right) \bigoplus A^{0,0}\left( M-S , \mathscr{H}om_{\mathcal{O}_M}\left( \mathcal{N}_{\mathcal{F}_0}^*, \tilde{\mathcal{S}}^1 \right) \right) \\
&\bigoplus A^{0,1}\left( M, \mathscr{H}om_{\mathcal{O}_M}\left( \mathcal{N}_{\mathcal{F}_0}^*, \Omega_M^1 \right) \right) \bigoplus A^{0,1}\left( M, \mathscr{H}om_{\mathcal{O}_M}\left( \mathcal{N}_{\mathcal{F}_0}^*, \mathcal{O}_M \right) \right) \bigoplus A^{0,2}\left( M, \mathcal{E}_{\mathcal{N}_{\mathcal{F}_0}^*} \right) \notag
\end{align}
defines a $2$-cocycle in the double complex $(\ref{oui24})$.
\end{lemma}

\begin{remark}
As in the same way in \textnormal{Remark \ref{os30}}, we can show that $\sigma_j\left(w_{0j}^\alpha \right)$ on $U_j-S$ actually lies in $\Gamma\left( U_j-S, \mathcal{A}^{0,0}\left(\tilde{\mathcal{S}}^1 \right) \right)$, i.e. $\sigma_j|_{U_j-S} \in \Gamma\left( U_j-S, \mathcal{A}^{0,0}\left( \tilde{\mathcal{S}}^1 \right) \right)$. More precisely, for $x\in U_j-S$, choose a neighborhood $U_x \subset U_j-S$ such that $\partial \Psi_j^\alpha = \sum_{\delta=1}^q a_{j_x}^{\alpha \delta} \wedge \Psi_j^\delta + \sum_{\eta=1}^q b_{j_x}^{\alpha \eta } \wedge w_{0j}^\eta$ for some $b_{j_x}^{\alpha \eta} \in \Gamma\left( U_x, \mathcal{A}^{0,0}\left( \Omega_M^1 \right) \right)$, and $\partial \tau_j^\alpha = \sum_{\delta=1}^q \tau_j^\delta a_{j_x}^{\alpha \delta} + \Psi_j^\alpha + \sum_{\eta=1}^q c_{j_x}^{\alpha \eta} w_{0j}^\eta $ for some $c_{j_x}^{\alpha \eta}\in \Gamma\left(U_x, \mathcal{A}^{0,0} \right)$. Then we can write
\begin{align} \label{its48}
&-\sum_{\gamma=1}^q\sum_{1\leq \beta \ne \gamma \leq q } (-1)^{q-\gamma} \tau_j^\gamma\cdot w_{0j}^1\wedge \cdots \wedge \Psi_j^\beta \wedge \cdots \wedge \hat{w}_{0j}^\gamma \wedge \cdots \wedge w_{0j}^q\wedge dw_{0j}^\alpha  \\
&-\sum_{\gamma=1}^q (-1)^{q-\gamma} \tau_j^\gamma \wedge \cdots \wedge \hat{w}_{0j}^\gamma \wedge \cdots w_{0j}^q \wedge \partial_z \left( \Psi_j^\alpha \right) - \sum_{\gamma=1}^q w_{0j}^1\wedge \cdots \wedge  \Psi_j^\gamma \wedge \cdots \wedge w_{0j}^q\wedge \partial_z \left(\tau_j^\alpha  \right) \notag  \\
& +  \sum_{\gamma=1}^q w_{0j}^1\wedge \cdots \wedge  \Psi_j^\gamma \wedge \cdots \wedge w_{0j}^q\wedge \Psi_j^\alpha \notag \\
&= w_{0j}^1 \wedge \cdots \wedge w_{0j}^q \wedge \left( \sum_{\gamma=1}^q \tau_j^\gamma b_{j_x}^{\alpha \gamma} + \sum_{\gamma=1}^q c_{j_x}^{\alpha \gamma} \Psi_j^\gamma  \right) \in \Gamma\left( U_x, \mathcal{A}^{0,0}\left( \tilde{\mathcal{S}}^1 \right) \right) \notag
\end{align}
\end{remark}

\begin{remark}
We show that $\tilde{\sigma}$ defines a global section in $A^{0,0}\left( M, \mathscr{H}om_{\mathcal{O}_M}\left( \mathcal{N}_{\mathcal{F}_0}^*, \bigwedge^{q+1} \Omega_M^1 \otimes \mathcal{L}_0 \right) \right)$. In fact,
{\Small{\begin{align*}
&\tilde{\sigma}_j\left( w_{0j}^\alpha \right) \\
&=- \sum_{\gamma=1}^q \sum_{1\leq \beta \ne \gamma \leq q} (-1)^{q-\gamma} \tau_j^\gamma \cdot w_{0j}^1 \wedge \cdots \wedge \Psi_j^\beta \wedge \cdots \wedge \widehat{w}_{0j}^\gamma \wedge \cdots  \wedge w_{0j}^q \wedge d w_{0j}^\alpha\\
&-\sum_{\gamma=1}^q (-1)^{q-\gamma} \tau_j^\gamma w_{0j}^1 \wedge \cdots \wedge \hat{w}_{0j}^\gamma \wedge \cdots \wedge w_{0j}^q \wedge \partial \left( \Psi_j^\alpha \right) - \sum_{\gamma=1}^q w_{0j}^1\wedge \cdots \wedge  \Psi_j^\gamma \wedge \cdots \wedge w_{0j}^q\wedge \partial \left(\tau_j^\alpha  \right) +  \sum_{\gamma=1}^q w_{0j}^1\wedge \cdots \wedge  \Psi_j^\gamma \wedge \cdots \wedge w_{0j}^q\wedge \Psi_j^\alpha \\
&= -  H_{0jk}  \sum_{\gamma=1}^q \sum_{1\leq \eta \ne \gamma \leq q} (-1)^{q-\gamma} \tau_k^\gamma \cdot w_{0k}^1 \wedge \cdots \wedge \Psi_k^\eta \wedge \cdots \wedge \widehat{w}_{0k}^\gamma \wedge \cdots  \wedge w_{0k}^q \wedge d \left( \sum_{\beta=1}^q h_{0jk}^{\alpha\beta}w_{0k}^\beta \right)\\
&-H_{0jk}\sum_{\gamma=1}^q (-1)^{q-\gamma} \tau_k^\gamma w_{0k}^1 \wedge \cdots \wedge \hat{w}_{0k}^\gamma \wedge \cdots \wedge w_{0k}^q \wedge \partial \left( \sum_{\beta=1}^q h_{0ij}^{\alpha\beta} \Psi_k^\beta \right) -  H_{0jk}\sum_{\gamma=1}^q w_{0k}^1\wedge \cdots \wedge  \Psi_k^\gamma \wedge \cdots \wedge w_{0k}^q\wedge \partial \left( \sum_{\beta=1}^q h_{0ij}^{\alpha\beta} \tau_j^\beta  \right)\\
& + H_{0jk} \sum_{\gamma=1}^q w_{0k}^1\wedge \cdots \wedge  \Psi_k^\gamma \wedge \cdots \wedge w_{0k}^q\wedge \Psi_k^\alpha\\
&= H_{0jk}\sum_{\beta=1}^q  h_{0jk}^{\alpha\beta} \tilde{\sigma}_k\left( w_{0k}^\beta \right) -  H_{0jk}  \sum_{\gamma=1}^q \sum_{1\leq \eta \ne \gamma \leq q} (-1)^{q-\gamma} \tau_k^\gamma \cdot w_{0k}^1 \wedge \cdots \wedge \Psi_k^\eta \wedge \cdots \wedge \widehat{w}_{0k}^\gamma \wedge \cdots  \wedge w_{0k}^q \wedge  \left(  d h_{0jk}^{\alpha\eta} \wedge w_{0k}^\eta+ dh_{0jk}^{\alpha \gamma} \wedge w_{0k}^\gamma \right)\\
&- H_{0jk}\sum_{\gamma, \beta=1}^q (-1)^{q-\gamma} \tau_k^\gamma w_{0k}^1 \wedge \cdots \wedge \widehat{w}_{0k}^\gamma \wedge \cdots \wedge w_{0k}^q \wedge d h_{0jk}^{\alpha\beta} \wedge \Psi_k^\beta -  H_{0jk}\sum_{\gamma,\beta=1}^q w_{0k}^1\wedge \cdots \wedge  \Psi_k^\gamma \wedge \cdots \wedge w_{0k}^q\wedge d h_{0jk}^{\alpha\beta}\cdot \tau_k^\beta  \\
&= H_{0jk}\sum_{\beta=1}^q h_{0jk}^{\alpha\beta} \tilde{\sigma}_k\left( w_{0k}^\beta \right) - H_{0jk}\sum_{\beta=1}^q (-1)^{q-\beta} \tau_k^\beta w_{0k}^1 \wedge \cdots \wedge \widehat{w}_{0k}^\beta \wedge \cdots \wedge w_{0k}^q \wedge d h_{0jk}^{\alpha\beta} \wedge \Psi_k^\beta- H_{0jk} \sum_{\beta=1}^q w_{0k}^1\wedge \cdots \wedge \Psi_k^\beta \wedge \cdots \wedge w_{0k}^q \wedge d h_{0jk}^{\alpha\beta} \tau_k^\beta \\
&= H_{0jk} \sum_{\beta=1}^q h_{0jk}^{\alpha \beta} \tilde{\sigma}_k \left( w_{0k}^\beta \right)
\end{align*}}}

\end{remark}

\begin{remark}
We show that $\pi= \{\pi_j\}$ defines a global section in $A^{0,1}\left( M, \mathscr{H}om_{\mathcal{O}_M} \left( \mathcal{N}_{\mathcal{F}_0}^*, \mathcal{O}_M   \right) \right)$. In fact, since $\sum_{\beta=1}^q k_j^{\alpha\beta} h_{0jk}^{\beta \eta}=\sum_{\beta=1}^q h_{0jk}^{\alpha\beta} k_k^{\beta \eta} + \left[ \psi, h_{0jk}^{\alpha\eta} \right]$, we have
\begin{align*}
&-\left[ \psi, \tau_j^\alpha \right]- i_\psi\left( \Psi_j^\alpha \right) + \sum_{\eta=1}^q k_j^{\alpha\eta} \tau_j^\eta= - \sum_{\beta=1}^q \left[ \psi, h_{0jk}^{\alpha\beta} \tau_k^\beta \right] - \sum_{\beta=1}^q h_{0jk}^{\alpha\beta}  i_{\psi}\left( \Psi_k^\beta \right) + \sum_{\eta,\xi=1}^q  k_j^{\alpha \eta}h_{0jk}^{\eta \xi}  \tau_k^\xi\\
&= - \sum_{\beta=1}^q h_{0jk}^{\alpha\beta} \left[ \psi, \tau_k^\beta \right] - \sum_{\beta=1}^q \left[ \psi,  h_{0jk}^{\alpha\beta} \right] \tau_k^\beta - \sum_{\beta=1}^q h_{0jk}^{\alpha\beta} i_\psi \left( \Psi_k^\beta \right) + \sum_{\beta, \xi=1}^q h_{0jk}^{\alpha \beta} k_k^{\beta \xi} \tau_k^\xi + \sum_{\xi=1}^q \left[ \psi, h_{0jk}^{\alpha \xi} \right] \tau_k^\xi\\
&= \sum_{\beta=1}^q h_{0jk}^{\alpha\beta}\left( - \left[ \psi, \tau_j^\beta\right] - i_{\psi}\left(\Psi_k^\beta \right)  + \sum_{\xi=1}^q k_k^{\beta \xi} \tau_k^\xi  \right) 
\end{align*}

\end{remark}

\begin{remark}
The following equality $(\ref{oui29})$ will be used in the proof of  \textnormal{Lemma \ref{oui28}}. We note that we can write
\begin{align*}
\partial \tau_j^\alpha = \sum_{\delta=1}^q a_{j_x}^{\alpha \delta} \tau_j^\delta + \Psi_j^\alpha + \sum_{\eta=1}^q c_{j_x}^{\alpha \eta} w_{0j}^\eta  \,\,\,\,\,\,\,\textnormal{for some $c_{j_x}^{\alpha \eta}\in \Gamma\left( U_x, \mathcal{A}^{0,0}\right)$}
\end{align*}
Then we have
\begin{align*}
&-\partial \bar{\partial} \tau_j^\alpha= - \sum_{\delta=1}^q  a_{j_x}^{\alpha \delta} \wedge \bar{\partial} \tau_j^\delta + \bar{\partial} \Psi_j^\alpha + \sum_{\eta=1}^q \bar{\partial} c_{j_x}^{\alpha \eta} w_{0j}^\eta \\
&\Longrightarrow - \partial i_{\psi}\left(w_{0j}^\alpha \right) = -\sum_{\delta=1}^q a_{j_x}^{\alpha\delta} \wedge i_{\psi}\left( w_{0j}^\delta \right) + \mathcal{L}_{\psi}\left(w_{0j}^\alpha \right)  - \sum_{\beta=1}^q k_j^{\alpha \beta} w_{0j}^\beta  + \sum_{\eta=1}^q \bar{\partial} c_{j_x}^{\alpha \eta} w_{0j}^\eta \\
&\Longrightarrow - di_\psi\left(w_{0j}^\alpha \right) + \bar{\partial} i_\psi\left(w_{0j}^\alpha \right) = -\sum_{\delta=1}^q a_{j_x}^{\alpha\delta} \wedge i_{\psi}\left( w_{0j}^\delta \right) + \mathcal{L}_{\psi}\left(w_{0j}^\alpha \right)  - \sum_{\beta=1}^q k_j^{\alpha \beta} w_{0j}^\beta  + \sum_{\eta=1}^q \bar{\partial} c_{j_x}^{\alpha \eta} w_{0j}^\eta \\
&\Longrightarrow \mathcal{L}_\psi\left(w_{0j}^\alpha \right) - i_\psi\left(\sum_{\delta=1}^q a_{j_x}^{\alpha \delta} \wedge w_{0j}^\delta \right) = -\sum_{\delta=1}^q a_{j_x}^{\alpha\delta} \wedge i_{\psi}\left( w_{0j}^\delta \right) + \mathcal{L}_{\psi}\left(w_{0j}^\alpha \right)  - \sum_{\beta=1}^q k_j^{\alpha \beta} w_{0j}^\beta  + \sum_{\eta=1}^q \bar{\partial} c_{j_x}^{\alpha \eta} w_{0j}^\eta \\
&\Longrightarrow \sum_{\delta=1}^q \left(  i_\psi\left( a_{j_x}^{\alpha \delta} \right)  - k_j^{\alpha \delta} + \bar{\partial} c_{j_x}^{\alpha \delta} \right) \wedge w_{0j}^\delta = 0
\end{align*}
This implies that
\begin{align}\label{oui29}
i_\psi\left( a_{j_x}^{\alpha \delta} \right) - k_j^{\alpha \delta} + \bar{\partial} c_{j_x}^{\alpha\delta}\equiv_{\mathcal{N}_{\mathcal{F}_0}^*} 0 \,\,\,\,\,\,\,\,\,\textnormal{on}\,\,\, U_x
\end{align}

\end{remark}

\begin{proof}[Proof of \textnormal{Lemma \ref{oui28}}]
By Lemma \ref{oci85}, it remains to prove that $(-1)^q \sigma + \hat{E}_1'\left(\tau \right)- \kappa=0$ and $\bar{\partial}\pi_j\left(  w_{0j}^\alpha \right) - i_{\left[\psi, \psi\right]}\left( w_{0j}^\alpha \right)=0$, and $(-1)^2l_j+ \hat{E}_2'\left( \sigma_j \right)=0$. We show that $(-1)^q \sigma_j +\hat{E}_1'\left(\tau_j \right) - \kappa_j=0$. In fact,
\begin{align} \label{its51}
&\bar{\partial} \left(\sum_{\gamma=1}^q \tau_j^\gamma b_{j_x}^{\alpha \gamma} + \sum_{\gamma=1}^q c_{j_x}^{\alpha \gamma} \Psi_j^\gamma \right) = -\sum_{\gamma=1}^q   b_{j_x}^{\alpha \gamma} \wedge \bar{\partial} \tau_j^\gamma + \sum_{\beta=1}^q \bar{\partial} b_{j_x}^{\alpha \beta} \cdot \tau_j^\beta  + \sum_{\gamma=1}^q \bar{\partial} c_{j_x}^{\alpha \gamma} \wedge  \Psi_j^\gamma + \sum_{\gamma=1}^q c_{j_x}^{\alpha \gamma} \bar{\partial} \Psi_j^\gamma  \\
&\equiv_{\mathcal{N}_{\mathcal{F}_0}^*} - \sum_{\gamma=1}^q b_{j_x}^{\alpha \gamma} \wedge i_\psi\left( w_{0j}^\gamma \right) + \sum_{\beta=1}^q \left( \mathcal{L}_{\psi}\left(a_{j_x}^{\alpha\beta} \right)  + \partial k_j^{\alpha\beta} - \sum_{\gamma=1}^q k_j^{\alpha \gamma} \wedge a_{j_x}^{\gamma \beta} - \sum_{\gamma=1}^q a_{j_x}^{\alpha \gamma} \wedge k_j^{\gamma \beta} \right) \tau_j^\beta \notag \\
& + \sum_{\gamma=1}^q k_j^{\alpha \gamma} \wedge \Psi_j^\gamma - \sum_{\delta=1}^q i_\psi\left(a_{j_x}^{\alpha \delta} \right) \wedge \Psi_j^\delta + \sum_{\gamma=1}^q c_{j_x}^{\alpha \gamma} \mathcal{L}_\psi\left(w_{0j}^\gamma \right) \notag
\end{align}

From $(\ref{oui29})$, we have

{\small{\begin{align}\label{its52}
&-\partial\left( -\left[ \psi, \tau_j^\alpha  \right] - i_{\psi}\left( \Psi_j^\alpha  \right) + \sum_{\eta=1}^q k_j^{\alpha \eta} \tau_j^\eta    \right) + \sum_{\delta=1}^q a_{j_x}^{\alpha\delta}\wedge \left( -\left[ \psi, \tau_j^\delta\right] - i_\psi\left( \Psi_j^\delta\right) + \sum_{\eta=1}^q k_j^{\delta \eta} \tau_j^\eta  \right)\\
&= \partial \mathcal{L}_\psi\left(\tau_j^\alpha \right) + \partial i_\psi\left(\Psi_j^\alpha \right) - \sum_{\eta=1}^q \partial k_j^{\alpha \eta} \tau_j^\eta + \sum_{\eta=1}^q k_j^{\alpha \eta} \wedge \partial \tau_j^\eta - \sum_{\delta=1}^q a_{j_x}^{\alpha \delta} \wedge \left[\psi, \tau_j^\delta \right] - \sum_{\delta=1}^q a_{j_x}^{\alpha \delta} \wedge i_\psi\left(\Psi_j^\delta \right)  + \sum_{\delta, \eta=1}^q a_{j_x}^{\alpha \delta} \wedge k_j^{\delta \eta} \tau_j^\eta \notag \\
&=d\mathcal{L}_\psi\left(\tau_j^\alpha \right) -\bar{\partial} \mathcal{L}_\psi\left(\tau_j^\alpha \right) + di_\psi\left( \Psi_j^\alpha \right) - \bar{\partial} i_\psi\left(\Psi_j^\alpha \right) - \sum_{\eta=1}^q \partial k_j^{\alpha \eta} \tau_j^\eta  + \sum_{\eta=1}^q k_j^{\alpha \eta} \wedge \left( \sum_{\delta=1}^q a_{j_x}^{\eta \delta} \tau_j^\delta + \Psi_j^\eta + \sum_{\delta=1}^q c_{j_x}^{\eta \delta} w_{0j}^\delta \right) \notag \\
& - \sum_{\delta=1}^q a_{j_x}^{\alpha \delta} \wedge \left[\psi, \tau_j^\delta \right] - \sum_{\delta=1}^q a_{j_x}^{\alpha \delta} \wedge i_\psi\left(\Psi_j^\delta \right)  + \sum_{\delta, \eta=1}^q a_{j_x}^{\alpha \delta} \wedge k_j^{\delta \eta} \tau_j^\eta \notag \\
&= -\mathcal{L}_\psi\left(d \tau_j^\alpha \right) + \mathcal{L}_\psi\left(\bar{\partial} \tau_j^\alpha \right) + i_\psi\left(d\Psi_j^\alpha \right) -\mathcal{L}_\psi\left( \Psi_j^\alpha \right) - i_\psi\left(\bar{\partial} \Psi_j^\alpha \right) - \sum_{\eta=1}^q \partial k_j^{\alpha \eta} \tau_j^\eta + \sum_{\eta, \delta=1}^q k_j^{\alpha \eta} \wedge a_{j_x}^{\eta \delta} \tau_j^\delta + \sum_{\eta=1}^q k_j^{\alpha \eta} \wedge \Psi_j^\eta + \sum_{\eta, \delta=1}^q k_j^{\alpha \eta} \wedge c_{j_x}^{\eta \delta} w_{0j}^\delta \notag \\
& - \sum_{\delta=1}^q a_{j_x}^{\alpha \delta} \wedge \left[\psi, \tau_j^\delta \right] - \sum_{\delta=1}^q a_{j_x}^{\alpha \delta} \wedge i_\psi\left(\Psi_j^\delta \right)  + \sum_{\delta, \eta=1}^q a_{j_x}^{\alpha \delta} \wedge k_j^{\delta \eta} \tau_j^\eta \notag \\
&\equiv_{\mathcal{N}_{\mathcal{F}_0}^*} -\mathcal{L}_\psi\left(\partial \tau_j^\alpha \right) + i_\psi\left(\partial \Psi_j^\alpha \right) - \mathcal{L}_\psi\left( \Psi_j^\alpha \right) - \sum_{\eta=1}^q \partial k_j^{\alpha \eta} \tau_j^\eta + \sum_{\eta, \delta=1}^q k_j^{\alpha \eta} \wedge a_{j_x}^{\eta \delta} \tau_j^\delta + \sum_{\eta=1}^q k_j^{\alpha \eta} \wedge \Psi_j^\eta \notag  \\
& - \sum_{\delta=1}^q a_{j_x}^{\alpha \delta} \wedge \left[\psi, \tau_j^\delta \right] - \sum_{\delta=1}^q a_{j_x}^{\alpha \delta} \wedge i_\psi\left(\Psi_j^\delta \right)  + \sum_{\delta, \eta=1}^q a_{j_x}^{\alpha \delta} \wedge k_j^{\delta \eta} \tau_j^\eta  \notag \\
&=- \mathcal{L}_\psi\left( \sum_{\delta=1}^q a_{j_x}^{\alpha \delta} \tau_j^\delta + \Psi_j^\alpha + \sum_{\eta=1}^q c_{j_x}^{\alpha \eta} w_{0j}^\eta \right) + i_\psi\left( \sum_{\beta=1}^q a_{j_x}^{\alpha \beta} \wedge \Psi_j^\beta + \sum_{\beta=1}^q b_{j_x}^{\alpha \beta} \wedge w_{0j}^\beta \right) \notag \\
& - \mathcal{L}_\psi\left( \Psi_j^\alpha \right) - \sum_{\eta=1}^q \partial k_j^{\alpha \eta} \tau_j^\eta + \sum_{\eta, \delta=1}^q k_j^{\alpha \eta} \wedge a_{j_x}^{\eta \delta} \tau_j^\delta + \sum_{\eta=1}^q k_j^{\alpha \eta} \wedge \Psi_j^\eta  - \sum_{\delta=1}^q a_{j_x}^{\alpha \delta} \wedge \left[\psi, \tau_j^\delta \right] - \sum_{\delta=1}^q a_{j_x}^{\alpha \delta} \wedge i_\psi\left(\Psi_j^\delta \right)  + \sum_{\delta, \eta=1}^q a_{j_x}^{\alpha \delta} \wedge k_j^{\delta \eta} \tau_j^\eta  \notag \\
&\equiv_{\mathcal{N}_{\mathcal{F}_0}^*} - \sum_{\delta=1}^q \mathcal{L}_\psi\left( a_{j_x}^{\alpha \delta} \right)\tau_j^\delta - 2 \mathcal{L}_\psi\left(\Psi_j^\alpha \right) - \sum_{\eta=1}^q c_{j_x}^{\alpha \eta} \mathcal{L}_\psi\left( w_{0j}^\eta \right) + \sum_{\beta=1}^q i_\psi \left( a_{j_x}^{\alpha \beta} \right) \wedge \Psi_j^\beta + \sum_{\beta=1}^q b_{j_x}^{\alpha \beta} \wedge i_\psi \left( w_{0j}^\beta \right)  \notag \\
&- \sum_{\eta=1}^q \partial k_j^{\alpha \eta} \tau_j^\eta + \sum_{\eta, \delta=1}^q k_j^{\alpha \eta} \wedge a_{j_x}^{\eta \delta} \tau_j^\delta + \sum_{\eta=1}^q k_j^{\alpha \eta} \wedge \Psi_j^\eta    + \sum_{\delta, \eta=1}^q a_{j_x}^{\alpha \delta} \wedge k_j^{\delta \eta} \tau_j^\eta \notag
\end{align}}}
and 
\begin{align} \label{its50}
-\kappa_i\left( w_{0j}^\alpha \right)  = 2 \mathcal{L}_\psi\left(\Psi_j^\alpha \right) - 2\sum_{\eta=1}^q k_j^{\alpha \eta} \wedge \Psi_j^\eta
\end{align}
Then $(\ref{its48}), (\ref{its51}),(\ref{its52})$ and $(\ref{its50})$ imply that $(-1)^q \sigma_j + \hat{E}_1'(\tau_j)  - \kappa_j=0$. Next we show that $\bar{\partial} \pi_j \left(w_{0j}^\alpha \right)- i_{\left[ \psi, \psi \right]}\left(w_{0j}^\alpha\right)=0$. In fact, we note that
\begin{align*}
\bar{\partial} \pi_j\left(w_{0j}^\alpha \right) & =\bar{\partial}\left( -\left[ \psi, \tau_j^\alpha  \right] - i_{\psi}\left( \Psi_j^\alpha  \right) + \sum_{\eta=1}^q k_j^{\alpha \eta} \tau_j^\eta \right)  = \left[ \psi, \bar{\partial} \tau_j^\alpha \right] - i_\psi\left(\bar{\partial} \Psi_j^\alpha \right)  - \sum_{\eta=1}^q k_j^{\alpha \eta } \wedge \bar{\partial} \tau_j^\eta \\
&= \left[ \psi, i_\psi\left(w_{0j}^\alpha \right) \right] - i_\psi\left( \mathcal{L}_\psi\left( w_{0j}^\alpha \right) - \sum_{\beta=1}^q k_j^{\alpha \beta}\wedge  w_{0j}^\beta \right) - \sum_{\eta=1}^q k_j^{\alpha \eta} \wedge i_\psi\left( w_{0j}^\eta \right)\\
&= \mathcal{L}_\psi i_\psi\left(w_{0j}^\alpha \right)- i_\psi \mathcal{L}_\psi \left(w_{0j}^\alpha \right) =  i_{\left[\psi, \psi \right]}\left(w_{0j}^\alpha \right)
\end{align*}
From $(\ref{its48})$, as in the same way with $(\ref{its49})$, we can show that $(-1)^2l_j+\hat{E}'_1\left( \sigma_j \right)=0$ on $U_x$, and so $M-S$. This completes the proof of Lemma \ref{oui28}.
\end{proof}

Then from $(\ref{it31})-(\ref{it32})$ and $(\ref{itt21})$ and $(\ref{itt20})$, and Lemma \ref{oui28}, we obtain
\begin{theorem}
Suppose given a compact foliated complex manifold $\left( M, \mathcal{N}_{\mathcal{F}_0}^* \right)$ with $\mathcal{N}_{\mathcal{F}_0}^*$ locally free, and $\left( \Psi, \tau, \left\{\left( \psi, k_j \right) \right\} \right)\in \mathbb{H}^1\left( M, \tilde{\mathcal{N}}_{\mathcal{F}_0}^* \right)$. In order that there may exist an unfolded analytic family $\left( \mathcal{M}, \tilde{\mathcal{N}}_\mathcal{F}^*, B, \pi \right)$ such that $\pi^{-1}(0)=\left( M , \mathcal{N}_{\mathcal{F}_0}^* \right)$ and that $\left( \frac{\partial \left( \mathcal{M}, \tilde{\mathcal{N}}_\mathcal{F}^* \right)}{\partial t}\right)|_{t=0}=\left( \Psi, \tau, \left\{ \left( \psi, k_j \right) \right\} \right)$, it is necessary that
\begin{align*}
\left( l, \sigma, \kappa, \pi, \left\{ - \left[\left(\psi, k_j \right), \left(\psi, k_j \right) \right]   \right\} \right) \in &A^{0,0}\left( M-S, \mathscr{H}om_{\mathcal{O}_M}\left(\mathcal{N}_{\mathcal{F}_0}^* , \tilde{\mathcal{S}}^2 \right) \right) \bigoplus A^{0,0}\left( M-S , \mathscr{H}om_{\mathcal{O}_M}\left( \mathcal{N}_{\mathcal{F}_0}^*, \tilde{\mathcal{S}}^1 \right) \right) \\
&\bigoplus A^{0,1}\left( M, \mathscr{H}om_{\mathcal{O}_M}\left( \mathcal{N}_{\mathcal{F}_0}^*, \Omega_M^1 \right) \right) \bigoplus A^{0,1}\left( M, \mathscr{H}om_{\mathcal{O}_M}\left( \mathcal{N}_{\mathcal{F}_0}^*, \mathcal{O}_M \right) \right) \bigoplus A^{0,2}\left( M, \mathcal{E}_{\mathcal{N}_{\mathcal{F}_0}^*} \right) \notag
\end{align*}
from $(\ref{itt22})$ in \textnormal{Lemma \ref{oui28}} defines the $0$ class in $\mathbb{H}^2\left( M, \tilde{\mathcal{E}}_{\mathcal{N}_{\mathcal{F}_0}^*}^\bullet \right)$.
\end{theorem}

\section{Isomorphisms of infinitesimal unfoldings and obstructions between \v Cech resolution and Dolbeault type bicomplex associated to $\tilde{\mathcal{E}}_{\mathcal{N}_{\mathcal{F}_0}^*}^\bullet$} \label{ni3} \

In this section we explicitly construct isomorphisms of first and second cohomology groups between \v Cech resolution of $\tilde{\mathcal{E}}_{\mathcal{N}_{\mathcal{F}_0}^*}^\bullet$ and \v Dolbeault type bicomplex associated to $\tilde{\mathcal{E}}_{\mathcal{N}_{\mathcal{F}_0}^*}^\bullet$. Let $\mathcal{U}=\left\{U_j \right\}$ be an Stein open covering of $M$ by coordinate neighborhoods with coordinates $z_j=\left(z_j^1,..., z_j^n \right)$ on $U_j$, and $\mathcal{N}_{\mathcal{F}_0}^*$ is generated by $w_{0j}^1,..., w_{0j}^p \in \Gamma\left( U_j , \Omega_M^1 \right)$ with the relation $w_{0j}^\alpha= \sum_{\beta=1}^q h_{0jk}^{\alpha\beta} w_{0j}^\beta$ on $U_{jk}$ for some $h_{0jk}^{\alpha\beta}\in \Gamma\left( U_{jk}, \mathcal{O}_M \right)$ and satisfies $w_{0j}^1 \wedge \cdots \wedge w_{0j}^q \wedge dw_{0j}^\alpha=0,\alpha=1,...,q$.

\subsection{Isomorphisms of infinitesimal deformations between \v Cech resolution and Dolbeault type bicomplex associated to $\tilde{\mathcal{E}}_{\mathcal{N}_{\mathcal{F}_0}^*}^\bullet$ } \label{b4}\

 We will explicitly define an isomorphism between the first cohomology group $\mathbb{H}^1\left( M, \tilde{\mathcal{E}}_{\mathcal{N}_{\mathcal{F}_0}^*}^\bullet \right)$ from the \v Cech resolution of $\tilde{\mathcal{E}}_{\mathcal{N}_{\mathcal{F}_0}^*}^\bullet$ and the first cohomology group $\textnormal{C}^1$ from $(\ref{b3})$, and will show that the infinitesimal unfolding $\left( \left\{ \beta_j(0)\right\}, \left\{ \gamma_j(0)\right\},\left\{ \left( \theta_{jk}, B_{jk}\right)  \right\} \right)$ $(\ref{ou6})$ from \v Cech resolution of $\tilde{\mathcal{E}}_{\mathcal{N}_{\mathcal{F}_0}^*}^\bullet$ corresponds to $\left( \frac{\partial Y}{\partial t}|_{t=0}, i_{\frac{\partial}{\partial t}}\left(V \right)|_{t=0}  , \left\{ \left( \frac{\partial \varphi}{\partial t}|_{t=0}, \frac{\partial Z_j}{\partial t}|_{t=0} \right) \right\}\right)$ (\ref{b1}) under the isomorphism.

Let $\left( \left\{ \beta_j\right\}, \left\{\gamma_j \right\},  \left\{ \left( \theta_{jk}, B_{jk} \right)      \right\}\right) \in C^0\left( \mathcal{U}, \mathscr{H}om_{\mathcal{O}_M}\left( \mathcal{N}_{\mathcal{F}_0}^*, \Omega_M^1 \right) \right) \bigoplus C^0\left( \mathcal{U}, \mathscr{H}om_{\mathcal{O}_M}\left( \mathcal{N}_{\mathcal{F}_0}^*, \mathcal{O}_M \right) \right) \bigoplus C^1\left( \mathcal{U}, \mathcal{E}_{\mathcal{N}_{\mathcal{F}_0}^*} \right)$ be an $1$-cocycle in the \v Cech resolution of $\tilde{\mathcal{E}}_{\mathcal{N}_{\mathcal{F}_0}^*}^\bullet$ with respect to the covering $\mathcal{U}$, so that $\left(\theta_{ij}, B_{ij}\right)- \left(\theta_{ik}, B_{ik}\right)+ \left(\theta_{jk}, B_{jk}^i\right)=0$ and $\beta_k\left( w_{0j}^\alpha \right) - \beta_j\left(w_{0j}^\alpha \right)+ \mathcal{L}_{\theta_{jk}}\left(w_{0j}^\alpha \right)- B_{jk}\left( w_{0j}^\alpha \right)=0$, where $B_{jk}^\alpha=B_{jk}\left(w_{0j}^\alpha \right) = \sum_{\beta=1}^q B_{jk}^{\alpha \delta} w_{0j}^\delta$ with $B_{jk}^{\alpha\delta}\in \Gamma\left(U_{jk}, \mathcal{O}_M \right)$, and $\sum_{\beta=1}^q w_{0j}^1 \wedge \cdots \wedge \beta_j\left( w_{0j}^\beta\right)\wedge \cdots \wedge w_{0j}^q \wedge dw_{0j}^\alpha + w_{0j}^1 \wedge \cdots \wedge w_{0j}^q \wedge d\beta_j\left(w_{0j}^\alpha \right)=0, \alpha=1,...,q$, and $\gamma_k\left( w_{0j}^\alpha \right)-\gamma_j\left(w_{0j}^\alpha \right)+ i_{\theta_{jk}}\left( w_{0j}^\alpha \right)=0$ and $\sum_{\gamma=1}^q (-1)^{q-\gamma} \gamma_j\left( w_{0j}^\gamma \right)w_{0j}^1 \wedge \cdots \widehat{w}_{0j}^\gamma \wedge \cdots \wedge w_{0j}^q \wedge dw_{0j}^\alpha + w_{0j}^1 \wedge \cdots \wedge w_{0j}^q \wedge d \gamma_j\left(w_{0j}^\alpha \right) - w_{0j}^1 \wedge \cdots \wedge w_{0j}^q \wedge \beta_j\left( w_{0j}^\alpha \right)=0  $. Since $\delta\left( \left\{\theta_{jk}, B_{jk} \right\}\right)=0$, there exists $-\delta\left(\left\{ \left(\xi_j, B_j\right) \right\} \right)= \left\{ \left( \theta_{jk}, B_{jk}\right) \right\}$ where $\left\{\left( \xi_j, B_j\right) \right\}\in C^0\left(\mathcal{U}, \mathcal{A}^{0,0}\left( \mathcal{E}_{\mathcal{N}_{\mathcal{F}_0}^*} \right) \right)$ such that $\xi_j\in \Gamma\left( U_j, \mathcal{A}^{0,0}\left(\Theta_M \right) \right)$ and $B_j\in \Gamma\left( U_j, \mathcal{A}^{0,0}\left( \mathscr{H}om_{\mathcal{O}_M}\left(\mathcal{N}_{\mathcal{F}_0}^*, \mathcal{N}_{\mathcal{F}_0}^* \right)\right) \right)$ with $B_j\left(w_{0j}^\alpha \right) = \sum_{\beta=1}^q B_j^{\alpha \beta} w_{0j}^\beta$ satisfying $\xi_j-\xi_k= \theta_{jk}$ and $B_j-B_k^j= B_{jk}$ with $B_k^j\left(w_{0j}^\alpha \right)= B_k \left( w_{0j}^\alpha \right)+ \sum_{\beta=1}^q \left[ \xi_k, h_{0jk}^{\alpha\beta} \right] w_{0k}^\beta$. Then $\left\{\left( \psi =\bar{\partial} \xi_j, \bar{\partial} B_j \right) \right\}\in A^{0,0}\left(M, \mathcal{E}_{\mathcal{N}_{\mathcal{F}_0}^*} \right)$ where $\phi\in A^{0,1}\left( M, \Theta_M \right)$ and $\bar{\partial} B_j\in \Gamma\left(U_j, \mathcal{A}^{0,1}\left( \mathscr{H}om_{\mathcal{O}_M}\left( \mathcal{N}_{\mathcal{F}_0}^*, \mathcal{N}_{\mathcal{F}_0}^* \right) \right) \right)$. On the other hand, let us define $\Psi_j \in \Gamma\left(U_j, \mathcal{A}^{0,0}\left(\mathscr{H}om_{\mathcal{O}_M}\left( \Theta_{\mathcal{F}_0}, \Omega_M^1 \right) \right) \right)$ by
\begin{align*}
\Psi_j:\Gamma\left(U_j,  \mathcal{N}_{\mathcal{F}_0}^* \right) &\to \Gamma\left( U_j, \mathcal{A}^{0,0} \left( \Omega_M^1 \right) \right) \\
 w_{0j}^\alpha &\mapsto \mathcal{L}_{\xi_j}\left(w_{0j}^\alpha \right) - B_j\left( w_{0j}^\alpha \right) - \beta_j\left(w_{0j}^\alpha \right)
\end{align*}
Then as in subsection \ref{ni1}, $\Psi:=\left\{ \Psi_j \right\}\in A^{0,0}\left(M, \mathscr{H}om_{\mathcal{O}_M}\left(\mathcal{N}_{\mathcal{F}_0}^* , \Omega_M^1 \right) \right)$, and we define
\begin{align*}
\tau_j:\Gamma\left( U_j, \mathcal{N}_{\mathcal{F}_0}^* \right) &\to \Gamma\left( U_j, \mathcal{A}^{0,0} \right) \\
                    w_{0j}^\alpha & \mapsto i_{\xi_j}\left(w_{0j}^\alpha \right) - \gamma_j\left(w_{0j}^\alpha \right)
\end{align*}
Then we have
{\small{\begin{align*}
\tau_j\left( w_{0j}^\alpha \right) - \sum_{\beta=1}^q h_{0jk}^{\alpha\beta} \tau_k\left( w_{0k}^\beta \right) = i_{\xi_j}\left(w_{0j}^\alpha \right)- \gamma_j\left(w_{0j}^\alpha \right) - \sum_{\beta=1}^q h_{0jk}^{\alpha\beta} \left(i_{\xi_k}\left( w_{0k}^\beta  \right) - \gamma_k\left( w_{0k}^\beta \right) \right) = i_{\theta_{jk}}\left(w_{0j}^\alpha \right) - \gamma_j\left(w_{0j}^\alpha \right) + \gamma_k\left( w_{0j}^\alpha \right)=0
\end{align*}}}
This implies that $\tau:=\left\{ \tau_j \right\}\in A^{0,0}\left( M, \mathscr{H}om_{\mathcal{O}_M}\left( \mathcal{N}_{\mathcal{F}_0}^*, \mathcal{O}_M \right) \right)$. Then we claim that $\left( \Psi, \tau, \left(\psi, \bar{\partial} B_j \right) \right)$ defines a $1$-cocycle in the Dolbeault type bicomplex. In fact, from subsection \ref{ni1}, it is enough to check that
\begin{align*}
-\bar{\partial} \tau_j\left(w_{0j}^\alpha \right) + i_{\psi}\left(w_{0j}^\alpha \right)= - \partial\left( i_{\xi_j}\left(w_{0j}^\alpha \right) \right) + i_{\psi}\left( w_{0j}^\alpha \right)=0
\end{align*}
Then we define the correspondence of $1$-cocycle between \v Cech resolution and Dolbeault type bicomplex by 
\begin{align}\label{b2}
\left(\left\{\beta_j \right\},  \left\{\gamma_j \right\}, \left\{ \left( \theta_{jk}, B_{jk} \right) \right\} \right) \mapsto \left( \Psi ,\tau, \left\{ \left(\psi, \bar{\partial} B_j\right) \right\}\right).
\end{align}

\begin{proposition}\label{c3}
We recall that given an unfolded analytic family $\left( \mathcal{M}, \tilde{\mathcal{N}}_\mathcal{F}^*, B, \pi \right)$ of unfoldings of $\left( M, \mathcal{N}_{\mathcal{F}_0}^* \right)=\pi^{-1}(0)$ as in \textnormal{Remark \ref{ou1}}, the infinitesimal unfolding at $t=0$ is given by $\left( \left\{ \beta_j(0)\right\}, \left\{ \gamma_j(0)\right\},\left\{ \left( \theta_{jk}, B_{jk}\right)  \right\} \right)$ from $(\ref{ou6})$ in the first cohomology group of the \v Cech resolution of $\tilde{\mathcal{E}}_{\mathcal{N}_{\mathcal{F}_0}^*}^\bullet$. On the other hand, the infinitesimal deformation at $t=0$ is also given by $\left( \frac{\partial Y}{\partial t}|_{t=0}, i_{\frac{\partial}{\partial t}}\left(V \right)|_{t=0}  , \left\{ \left( \frac{\partial \varphi}{\partial t}|_{t=0}, \frac{\partial Z_j}{\partial t}|_{t=0} \right) \right\}\right) $ from $(\ref{b1})$ in the first cohomology group $\textnormal{C}^1$ of the Dolbeault type bicomplex associated to $\tilde{\mathcal{E}}_{\mathcal{N}_{\mathcal{F}_0}^*}^\bullet$. Under the isomorphism $(\ref{b2})$ from $\mathbb{H}^1\left(M, \tilde{\mathcal{E}}_{\mathcal{N}_{\mathcal{F}_0}^*}^\bullet \right)\to \textnormal{C}^1$, 
\begin{align*}
\left( \left\{ \beta_j(0)\right\}, \left\{ \gamma_j(0)\right\},\left\{ \left( \theta_{jk}, B_{jk}\right)  \right\} \right)\,\,\,\,\,\,\,\textnormal{corresponds to}\,\,\,\,\, \left( \frac{\partial Y}{\partial t}|_{t=0}, i_{\frac{\partial}{\partial t}}\left(V \right)|_{t=0}  , \left\{ \left( \frac{\partial \varphi}{\partial t}|_{t=0}, \frac{\partial Z_j}{\partial t}|_{t=0} \right) \right\}\right) 
\end{align*}
\end{proposition} 

\begin{proof}
We keep the notations in the proof of Proposition \ref{c1}. It remains to show that $\left\{ \gamma_j(0)\right\}$ corresponds to $i_{\frac{\partial}{\partial t}}\left( V \right)|_{t=0}$. In fact, by considering $(\ref{c2})$ and setting $t=0$, we have
\begin{align*}
Q_j^{\alpha \gamma}|_{t=0}= v_j^{\alpha \gamma}(\xi_j, 0)+ \sum_{\beta=1}^n w_{0j}^{\alpha \beta} \frac{\partial \xi_j^\beta}{\partial t_\gamma}|_{t=0}
\end{align*}
If we set $\frac{\partial}{\partial t}= \sum_{\gamma=1}^m c_\gamma \frac{\partial}{\partial t_\gamma}$, then
\begin{align*}
\sum_{\gamma=1}^m c_\gamma Q_j^{\alpha \gamma}|_{t=0} = \sum_{\gamma=1}^m c_\gamma v_j^{\alpha \gamma}(\xi_j,0)+ i_{\xi_j}\left( w_{0j}^\alpha \right) \Longrightarrow i_{\frac{\partial}{\partial t}}\left( V_j^\alpha  \right)|_{t=0} = i_{\xi_j}\left( w_{0j}^\alpha \right) - \left( - i_{\frac{\partial}{\partial t}}\left( v_j^\alpha \right)|_{t=0} \right)
\end{align*}
This completes the proof of Proposition \ref{c3}.
\end{proof}

\subsection{Isomorphism of obstructions between \v Cech resolution and Dolbeault type bicomplex associated to $\tilde{\mathcal{E}}_{\mathcal{N}_{\mathcal{F}_0}^*}^\bullet$} \

We will explicitly define an isomorphism between the second cohomology group $\mathbb{H}^2\left( M, \tilde{\mathcal{E}}_{\mathcal{N}_{\mathcal{F}_0}^*}^\bullet \right)$ from \v Cech resolution of $\tilde{\mathcal{E}}_{\mathcal{N}_{\mathcal{F}_0}^*}^\bullet$, and the second cohomology group $\textnormal{C}^2$ from the Dolbeault type bicomplex associated to $\tilde{\mathcal{E}}_{\mathcal{N}_{\mathcal{F}_0}^*}^\bullet$.

Let $\left(\left\{ c_i \right\}, \left\{ C_j \right\},  \left\{ b_{jk} \right\}, \left\{ B_{jk}\right\}, \left\{ a_{ijk}, g_{ijk} \right\} \right)\in C^0\left( \mathcal{U}, \mathscr{H}om_{\mathcal{O}_M}\left( \mathcal{N}_{\mathcal{F}_0}^* ,  \tilde{\mathcal{S}}^2 \right) \right) \bigoplus  C^0\left( \mathcal{U}, \mathscr{H}om_{\mathcal{O}_M}\left( \mathcal{N}_{\mathcal{F}_0}^*, \tilde{\mathcal{S}}^1 \right) \right) $ $\bigoplus $ $ C^1\left( \mathcal{U} , \mathscr{H}om_{\mathcal{O}_M}\left( \mathcal{N}_{\mathcal{F}_0}^*, \Omega_M^1 \right) \right) \bigoplus C^1\left( \mathcal{U}, \mathscr{H}om_{\mathcal{O}_M}\left( \mathcal{N}_{\mathcal{F}_0}^*, \mathcal{O}_M \right) \right) \bigoplus C^2\left( \mathcal{U}, \mathcal{E}_{\mathcal{N}_{\mathcal{F}_0}^*} \right) $ be a $2$-cocycle.

Since $-\delta\left( \left\{ \left( a_{ijk}, g_{ijk}\right) \right\} \right)=0$, we have $\left\{\left(d_{ij}, e_{ij} \right) \right\}\in C^1\left(\mathcal{U}, \mathcal{A}^{0,0}\left(\mathcal{E}_{\mathcal{N}_{\mathcal{F}_0}^*} \right) \right)$ such that $\delta\left( \left\{\left( d_{ij}, e_{ij} \right) \right\}\right)=\left\{ \left( a_{ijk}, g_{ijk}\right) \right\}$ where $d_{ij}\in \Gamma\left( U_{ij}, \mathcal{A}^{0,0}\left(\Theta_M \right) \right)$ and $e_{ij}\in \Gamma\left(U_{ij},  \mathcal{A}^{0,0}\left(\mathscr{H}om_{\mathcal{O}_M}\left( \mathcal{N}_{\mathcal{F}_0}^*, \mathcal{N}_{\mathcal{F}_0}^* \right) \right) \right)$ defined by
\begin{align*}
e_{ij}:\Gamma\left( U_{ij},  \mathcal{N}_{\mathcal{F}_0}^* \right) &\to \Gamma\left( U_{ij}, \mathcal{A}^{0,0}\left( \mathcal{N}_{\mathcal{F}_0}^* \right)\right)\\
 w_{0i}^\alpha &\mapsto e_{ij}\left(w_{0i}^\alpha \right)= e_{ij}^\alpha
\end{align*}
satisfying $d_{ij}-d_{ik}+ d_{jk}= a_{ijk}$ and $e_{ij}\left( w_{0i}^\alpha\right) - e_{ik}\left(w_{0i}^\alpha\right) + e_{jk}\left(w_{0i}^\alpha\right) + \sum_{\beta=1}^q \left[d_{jk}, h_{0ij}^{\alpha\beta} \right]w_{0j}^\beta= g_{ijk}\left(w_{0i}^\alpha \right)$. Since $\delta\left( \left\{\left(\bar{\partial} d_{jk} , \bar{\partial} e_{jk} \right) \right\}\right)=0$, there exists $\left\{\left( p_j, u_j \right) \right\} \in C^0\left( \mathcal{U}, \mathcal{A}^{0,1}\left(\mathcal{E}_{\mathcal{N}_{\mathcal{F}_0}^*}\right) \right)$ with $p_j \in \Gamma\left( U_j,  \mathcal{A}^{0,1}\left(\Theta_M \right) \right)$  and $u_j\in \Gamma\left( U_j, \mathcal{A}^{0,1}\left( \mathscr{H}om_{\mathcal{O}_M}\left( \mathcal{N}_{\mathcal{F}_0}^* , \mathcal{N}_{\mathcal{F}_0}^* \right) \right) \right)$ defined by
\begin{align*}
u_j:\Gamma\left( U_j, \mathcal{N}_{\mathcal{F}_0}^* \right) &\to \Gamma\left(U_j, \mathcal{A}^{0,1}\left(\mathcal{N}_{\mathcal{F}_0}^* \right) \right)\\
 w_{0j}^\alpha &\mapsto u_j^\alpha
\end{align*}
 such that $p_k- p_j=\bar{\partial} d_{jk}$, and $u_k\left( w_{0j}^\alpha \right)+ \sum_{\beta=1}^q \left[ p_k, h_{0jk}^{\alpha\beta} \right] w_{0k}^\beta - u_j\left( w_{0j}^\alpha \right)= \bar{\partial} e_{jk}^\alpha$. Then $\left\{\left(\bar{\partial} p_j, \bar{\partial} u_j \right)\right\}\in A^{0,2}\left( M, \mathcal{E}_{\mathcal{N}_{\mathcal{F}_0}^* } \right)$. On the second floor, as in subsection \ref{ni1}, there exists $f_j\in \Gamma\left( U_j, \mathcal{A}^{0,0} \left( \mathscr{H}om_{\mathcal{O}_M}\left( \mathcal{N}_{\mathcal{F}_0}^*, \Omega_M^1 \right) \right) \right)$ defined by
\begin{align*}
f_j:\Gamma\left( U_j, \mathcal{N}_{\mathcal{F}_0}^* \right)&\to \Gamma\left( U_j,  \mathcal{A}^{0,0}\left( \Omega_M^1 \right) \right)\\
 w_{0j}^\alpha &\mapsto f_j^\alpha
\end{align*}
such that $f_j\left( w_{0j}^\alpha \right)- f_k\left( w_{0j}^\alpha \right)= b_{jk}^\alpha-\mathcal{L}_{d_{jk}}\left(w_{0j}^\alpha \right)+ e_{jk}^\alpha$ and we define $Q_j\in \Gamma\left( U_j, \mathcal{A}^{0,1} \left( \mathscr{H}om_{\mathcal{O}_M}\left( \mathcal{N}_{\mathcal{F}_0}^*, \Omega_M^1 \right) \right) \right)$ by
\begin{align*}
Q_j:\Gamma\left( U_j, \mathcal{N}_{\mathcal{F}_0}^* \right) &\to \Gamma\left( U_j, \mathcal{A}^{0,1}\left(\Omega_M^1 \right) \right)\\
T_{0j}^\alpha &\mapsto  - \mathcal{L}_{p_j}\left( w_{0j}^\alpha \right) + u_j^\alpha + \bar{\partial} f_j^\alpha
\end{align*}
Then $Q=\{Q_j\} \in A^{0,1}\left( M, \mathscr{H}om_{\mathcal{O}_M}\left( \mathcal{N}_{\mathcal{F}_0}^*, \Omega_M^1 \right) \right)$ as in subsection \ref{a16}. On the other hand, we note that
\begin{align*}
i_{d_{ij}}\left( w_{0i}^\alpha\right) - i_{d_{ik}}\left( w_{0i}^\alpha \right) + i_{d_{jk}}\left(w_{0j}^\alpha \right) - \left(B_{ij}^\alpha- B_{ik}^\alpha + \sum_{\beta=1}^q h_{0ij}^{\alpha\beta} B_{jk}^\beta    \right) = i_{a_{ijk}}\left(w_{0i}^\alpha \right) - \left( B_{ij}^\alpha- B_{ik}^\alpha + \sum_{\beta=1}^q h_{0ij}^{\alpha\beta} B_{jk}^\beta \right) =0
\end{align*}

This implies that there exists $r_j\in \Gamma\left( U_j, \mathcal{A}^{0,0}\left(\mathscr{H}om_{\mathcal{O}_M}\left( \mathcal{N}_{\mathcal{F}_0}^*, \mathcal{O}_M \right) \right) \right)$ defined by
\begin{align*}
r_j:\Gamma\left( U_j, \mathcal{N}_{\mathcal{F}_0}^* \right) &\to \Gamma\left( U_j, \mathcal{A}^{0,0}\right) \\
w_{0j}^\alpha &\mapsto r_j^\alpha
\end{align*}
such that $r_j\left(w_{0j}^\alpha \right)- r_k\left( w_{0j}^\alpha \right)= B_{jk}^\alpha - i_{d_{jk}}\left( w_{0j}^\alpha \right)$. Then  $\bar{\partial} r_j\left(w_{0j}^\alpha \right)- \bar{\partial} r_k\left( w_{0j}^\alpha \right)=  - i_{\bar{\partial} d_{jk}}\left( w_{0j}^\alpha \right)$ and so
\begin{align*}
\sum_{\beta=1}^q h_{0jk}^{\alpha\beta}\left( i_{p_k} \left( w_{0k}^\beta \right) - \bar{\partial} r_k^\beta     \right) - \left( i_{p_j}\left(w_{0j}^\alpha \right) - \bar{\partial} r_j^\alpha \right) = i_{\bar{\partial} d_{jk}}\left(w_{0j}^\alpha \right) - i_{\bar{\partial} d_{jk}}\left( w_{0j}^\alpha \right) = 0
\end{align*}
Let us define $V_j\in \Gamma\left( U_j, \mathcal{A}^{0,1}\left(\mathscr{H}om_{\mathcal{O}_M}\left( \mathcal{N}_{\mathcal{F}_0}^*, \mathcal{O}_M \right) \right) \right)$ by
\begin{align*}
V_j:\Gamma\left( U_j, \mathcal{N}_{\mathcal{F}_0}^* \right) &\to \Gamma\left( U_j, \mathcal{A}^{0,1} \right) \\
        w_{0j}^\alpha &\mapsto - i_{p_j}\left(w_{0j}^\alpha \right) + \bar{\partial} r_j^\alpha
\end{align*}
Then we see that $V=\{V_j \}$ defines a global section in $A^{0,1}\left( M,  \mathscr{H}om_{\mathcal{O}_M}\left( \mathcal{N}_{\mathcal{F}_0}^*, \mathcal{O}_M \right) \right)$.

On the third floor, for $x\in U_{jk}- S$, choose a neighborhood $U_x\subset U_{jk}-S$ such that $dw_{0j}^\alpha - \sum_{\beta=1}^q a_{j_x}^{\alpha\beta} \wedge w_{0j}^\beta,\alpha=1,...q$ for some $a_{j_x}^{\alpha\beta}\in \Gamma\left( U_x, \Omega_M^1 \right)$. Then $\sum_{\beta=1}^q a_{j_x}^{\alpha\beta} h_{0jk}^{\beta \eta} -  d h_{0jk}^{\alpha \eta}- \sum_{\beta=1}^q h_{0jk}^{\alpha\beta} a_{k_x}^{\beta \eta}\in \Gamma\left(U_{x}, \mathcal{N}_{\mathcal{F}_0}^* \right)$. Let us define $R_j\in \Gamma\left( U_j- S,  \mathcal{A}^{0,0}\left(\mathscr{H}om_{\mathcal{O}_M}\left( \mathcal{N}_{\mathcal{F}_0}^*, \tilde{\mathcal{S}}^2 \right) \right) \right)$ by
\begin{align*}
R_j:\Gamma\left( U_j-S, \mathcal{N}_{\mathcal{F}_0}^*  \right) &\to \Gamma\left( U_j- S, \tilde{\mathcal{S}}^2 \right) \\
                                     w_{0j}^\alpha &\mapsto  -c_j^\alpha - \sum_{\beta=1}^q w_{0i}^1\wedge \cdots \wedge f_j^\beta \wedge \cdots \wedge w_{0j}^q \wedge dw_{0j}^\alpha - w_{0j}^1\wedge \cdots \wedge w_{0j}^q \wedge df_j^\alpha
\end{align*}
Then $R:=\{R_j\}$ defines a global section in $A^{0,0}\left( M-S, \mathscr{H}om_{\mathcal{O}_M}\left(  \mathcal{N}_{\mathcal{F}_0}^*,  \tilde{\mathcal{S}}^2 \right) \right)$ as in subsection \ref{a16}. On the other hand, we note that

{\small{\begin{align*}
&C_j^\alpha + w_{0j}^1\wedge \cdots \wedge w_{0j}^q\wedge \left( -f_j^\alpha + dr_j^\alpha - \sum_{\beta=1}^q a_{j_x}^{\alpha\beta} r_j^\beta    \right) -\sum_{\beta=1}^q h_{0jk}^{\alpha\beta} \det\left(h_{0jk}\right)\left( C_k^\beta  + w_{0k}^1\wedge \cdots \wedge w_{0k}^q\wedge \left(- f_k^\beta +dr_k^\beta - \sum_{\eta=1}^q a_{k_x}^{\beta \eta} r_k^\eta \right)      \right)\\
&= C_j^\alpha  - \sum_{\beta=1}^q h_{0jk}^{\alpha\beta} \det\left(h_{0jk}\right) C_k^\beta - \det\left(h_{0jk} \right)w_{0k}^1 \wedge \cdots \wedge w_{0k}^q \wedge \left( \sum_{\beta=1}^q h_{0jk}^{\alpha\beta} f_k^\beta + b_{jk}^\alpha - \mathcal{L}_{d_{jk}}\left(w_{0j}^\alpha \right) + e_{jk}^\alpha \right)\\
&+\det\left(h_{0jk}\right) w_{0k}^1 \wedge \cdots \wedge w_{0k}^q \wedge \left(  d\left(\sum_{\beta=1}^q h_{0jk}^{\alpha\beta} r_k^\beta -  i_{d_{jk}}\left(w_{0j}^\alpha\right) + B_{jk}^\alpha \right) - \sum_{\beta=1}^q a_{j_x}^{\alpha\beta} \wedge \left( \sum_{\eta=1}^q h_{0jk}^{\beta \eta} r_k^\eta  - i_{d_{jk}}\left(w_{0j}^\beta\right) + B_{jk}^\beta \right) \right)\\
& -\sum_{\beta=1}^q h_{0jk}^{\alpha\beta} \det(h_{0jk}) w_{0k}^1\wedge \cdots \wedge w_{0k}^q\wedge \left(-f_k^\beta + dr_k^\beta - \sum_{\eta=1}^q a_{k_x}^{\beta \eta} r_k^\eta \right)      \\
&= C_j^\alpha - \sum_{\beta=1}^q h_{0jk}^{\alpha\beta} \det\left( h_{0jk}\right) C_k^\beta+ w_{0j}^1\wedge \cdots \wedge w_{0j}^q \wedge \left( d B_{jk}^\alpha - \sum_{\beta=1}^q a_{j_x}^{\alpha\beta}  B_{jk}^\beta  - b_{jk}^\alpha  \right) = 0
\end{align*}}}
Let us define $G_j\in \Gamma\left(U_j-S, \mathcal{A}^{0,0}\left( \mathscr{H}om_{\mathcal{O}_M}\left( \mathcal{N}_{\mathcal{F}_0}^*, \tilde{\mathcal{S}}^1 \right) \right) \right)$ by
\begin{align*}
G_j:\Gamma\left( U_j- S, \mathcal{N}_{\mathcal{F}_0}^* \right) &\to \Gamma\left( U_j- S , \tilde{\mathcal{S}}^1   \right) \\ 
w_{0j}^\alpha &\mapsto w_{0j}^1\wedge \cdots \wedge w_{0j}^q \wedge f_j^\alpha - C_j^\alpha - \sum_{\beta=1}^q (-1)^{q-\beta} r_j^\beta w_{0j}^1 \wedge \cdots \wedge \widehat{w}_{0j}^\beta \wedge \cdots \wedge w_{0j}^q - w_{0j}^1 \wedge \cdots \wedge w_{0j}^q \wedge \partial r_j^\alpha
\end{align*}
Then $G:=\{G_j\}$ defines a global section in $A^{0,0}\left( M-S , \mathscr{H}om_{\mathcal{O}_M}\left( \mathcal{N}_{\mathcal{F}_0}^*, \tilde{\mathcal{S}}^1 \right) \right)$. Then we claim that
\begin{align*}
\left(R, G, Q, V, \left\{ \left(\bar{\partial} p_j, \bar{\partial} u_j \right) \right\} \right) \in A^{0,0}\left( M-S, \mathscr{H}om_{\mathcal{O}_M}\left( \mathcal{N}_{\mathcal{F}_0}^*,  \tilde{\mathcal{S}}^2 \right) \right) \bigoplus A^{0,0}\left( M-S,  \right)
\end{align*}
defines a $2$-cocycle. In fact, in view of subsection \ref{a16}, it is clear that $\hat{G}_2|_{M-S}\left( R_j, G_j \right)=0$ and $\bar{\partial}\left(\left(\bar{\partial} p_j, \bar{\partial} u_j \right) \right)=0$. On the other hand, on $U_x$, 
\begin{align*}
&(-1)^q\bar{\partial} G_j\left( w_{0j}^\alpha\right)=  - (-1)^q \bar{\partial}\left( w_{0j}^q \wedge \cdots  \wedge w_{0j}^q\wedge\left( - f_j^\alpha +\partial r_j^\alpha - \sum_{\beta=1}^q a_{j_x}^{\alpha\beta} \wedge r_j^\beta \right)\right)\\
&= w_{0j}^1\wedge \cdots \wedge w_{0j}^q\wedge \left( \bar{\partial} f_j^\alpha  + \bar{\partial}\left(\bar{\partial} r_j^\alpha \right)  - \sum_{\beta=1}^q a_{j_x}^{\alpha\beta} \wedge \bar{\partial} r_j^\beta \right) = - \hat{G}_i'\left( Q_j, V_j \right) \Longrightarrow (-1)^q\bar{\partial}\left(R_j, G_j \right)+ \hat{G}_1'\left( Q_j, V_j \right)=0
\end{align*}
and in view of subsection \ref{a16} we have
\begin{align*}
\bar{\partial} V_j\left( w_{0j}^\alpha \right)= \bar{\partial}\left(-i_{p_j}\left( w_{0j}^\alpha \right)+ \bar{\partial} r_j^\alpha \right)= - i_{\bar{\partial} p_j}\left(w_{0j}^\alpha \right) \Longrightarrow  \bar{\partial} V_j + \hat{G}_0'\left(\bar{\partial} p_j \right)= 0 \Longrightarrow \bar{\partial}\left( Q_j, V_j\right)+ \hat{G}_0\left( \bar{\partial} p_j, \bar{\partial} u_j \right)=0
\end{align*}

Then we define the correspondence of $2$-cocycle between \v Cech resolution and Dolbeault type bicomplex by
\begin{align}\label{b7}
\left(\left\{c_j \right\},\left\{C_j \right\},  \left\{ b_{jk}  \right\}, \left\{B_{jk}\right\},   \left\{ \left( a_{ijk}, g_{ijk} \right) \right\} \right)\mapsto \left(R, G, Q , V,   \left\{ \left(\bar{\partial} p_j, \bar{\partial} u_j  \right) \right\} \right)
\end{align}

Conversely, let $\left(R, G, Q, V,   \left( \phi, \left\{N_j \right\} \right) \right)$ be a $2$-cocycle in the Dolbeault type bicomplex associated to $\tilde{\mathcal{E}}_{\mathcal{N}_{\mathcal{F}_0}^*}^\bullet$. We keep the notations in subsection \ref{a16}. We note that since $\bar{\partial} V\left(w_{0j}^\alpha \right) + i_{\phi}\left( w_{0j}^\alpha \right)=0$, we have $\bar{\partial} \left(V\left( w_{0j}^\alpha\right) + i_{p_j}\left( w_{0j}^\alpha \right) \right)=0$. Then there exists $r_j\in \Gamma\left( U_j, \mathcal{A}^{0,0}\left( \mathcal{N}_{\mathcal{F}_0}^*, \mathcal{O}_M \right) \right)$ such that $\bar{\partial} r_j\left( w_{0j}^\alpha\right)= V_j^\alpha+ i_{p_j}\left( w_{0j}^\alpha \right)$. Then we note that
\begin{align*}
\bar{\partial} r_j\left( w_{0j}^\alpha \right)- \bar{\partial} r_k\left( w_{0j}^\alpha \right) = i_{p_j-p_k}\left( w_{0j}^\alpha \right) = i_{-\bar{\partial} d_{jk}}\left( w_{0j}^\alpha \right)
\end{align*}
Then we define
\begin{align*}
B_{jk}:\Gamma\left( U_{jk}, \mathcal{N}_{\mathcal{F}_0}^* \right) & \to \Gamma\left( U_{jk}, \mathcal{O}_M \right) \\
 w_{0j}^\alpha & \mapsto r_j\left( w_{0j}^\alpha \right)-r _k\left( w_{0j}^\alpha \right)+ i_{d_{jk}}\left( w_{0j}^\alpha \right)
\end{align*}
Then $\left\{ B_{jk}\right\} \in C^1\left(\mathcal{U}, \mathscr{H}om_{\mathcal{O}_M}\left( \mathcal{N}_{\mathcal{F}_0}^*, \mathcal{O}_M \right) \right)$. On the other hand, for $x\in U_j-S$, choose a neighborhood $U_x\subset U_j-S$ such that $dw_{0j}^\alpha - \sum_{\beta=1}^q a_{j_x}^{\alpha\beta } \wedge w_{0j}^\beta$ for some $a_{j_x}^{\alpha\beta}\in \Gamma\left( U_x, \Omega_M^1 \right),\alpha=1,...,q$. We can write $G\left( w_{0j}^\alpha \right)= w_{0j}^1\wedge \cdots \wedge w_{0j}^q \wedge G_{j_x}^\alpha$ for some $G_{j_x}\in \Gamma\left(U_x, \mathcal{A}^{0,0}\left(\Omega_M^1 \right) \right)$. Then since $w_{0j}^1\wedge \cdots \wedge w_{0j}^q \wedge \left(  \bar{\partial} G_{j_x}^\alpha - Q_j^\alpha - \partial V_j^\alpha + \sum_{\beta=1}^q a_{j_x}^{\alpha\beta}\cdot V_j^\beta \right)=0$. Then we have
\begin{align*}
0=&w_{0j}^1 \wedge \cdots \wedge w_{0j}^q \wedge \left( - \bar{\partial} f_j^\alpha + \mathcal{L}_{p_j}\left(w_{0j}^\alpha \right) - u_j^\alpha + \bar{\partial} G_{j_x}^\alpha - \partial \left( \bar{\partial} r_j^\alpha - i_{p_j}\left(w_{0j}^\alpha \right)  \right) + \sum_{\beta=1}^q a_{j_x}^{\alpha\beta} \wedge \left( \bar{\partial} r_j^\beta - i_{p_j}\left(w_{0j}^\beta \right) \right) \right) \\
&= w_{0j}^1 \wedge \cdots \wedge w_{0j}^q \wedge \left( - \bar{\partial} f_j^\alpha + \bar{\partial} G_{j_x}^\alpha +\bar{\partial} \left( \partial r_j^\alpha - \sum_{\beta=1}^q a_{j_x}^{\alpha\beta} \cdot r_j^\beta     \right) + \tilde{\mathcal{L}}_{p_j}\left( w_{0j}^\alpha \right)  + \partial i_{p_j}\left(w_{0j}^\alpha \right) - i_{p_j}\partial w_{0j}^\beta\right)  \\
&= w_{0j}^1 \wedge \cdots \wedge w_{0j}^q \wedge \left( - \bar{\partial} f_j^\alpha + \bar{\partial} G_{j_x}^\alpha +\bar{\partial} \left( \partial r_j^\alpha - \sum_{\beta=1}^q a_{j_x}^{\alpha\beta} \wedge f_j^\beta     \right)  \right)
\end{align*}

Then we have
{\small{\begin{align*}
C_j^\alpha&:= - G_j\left(w_{0j}^\alpha \right) + w_{0j}^1 \wedge \cdots \wedge w_{0j}^q \wedge f_j^\alpha -  \sum_{\beta=1}^q (-1)^{q-\beta} r_j^\beta w_{0j}^1 \wedge \cdots \wedge \widehat{w}_{0j}^\beta \wedge \cdots \wedge w_{0j}^q \wedge dw_{0j}^\alpha - w_{0j}^1 \wedge \cdots \wedge w_{0j}^q \wedge \partial r_j^\alpha  \in \Gamma\left( U_j- S, \tilde{\mathcal{S}}^1 \right)
\end{align*}}}
Hence $C_j^\alpha \in \Gamma\left( U_j, \tilde{\mathcal{S}}^1 \right)$. We define $C_j\in \Gamma\left( U_j, \mathcal{A}^{0,0}\left( \mathscr{H}om_{\mathcal{O}_M}\left( \mathcal{N}_{\mathcal{F}_0}^*, \tilde{\mathcal{S}}^1 \right) \right) \right)$ by
\begin{align*}
C_j:\Gamma\left( U_j, \mathcal{N}_{\mathcal{F}_0}^* \right) &\to \Gamma\left( U_j, \tilde{\mathcal{S}}^1 \right) \\
 w_{0j}^\alpha &\mapsto C_j^\alpha
\end{align*}

Then we define the inverse correspondence by $\left( R,  G, Q, V,  \left(\left\{ \phi, N_j \right\} \right) \right)\mapsto \left(\left\{c_j\right\}, \left\{C_j \right\}, \left\{ b_{jk}\right\}, \left\{ B_{jk} \right\}, \left\{\left( a_{ijk}, g_{ijk} \right) \right\} \right)$.

\begin{proposition} \label{b14}
We recall from \textnormal{Lemma \ref{ou50}} that given a $1$-cocycle $\left(\left\{ \beta_j\right\}, \left\{ \gamma_j\right\}, \left\{ \left( \theta_{jk}, B_{jk}\right) \right\} \right)$ in $\mathbb{H}^1\left( M, \tilde{\mathcal{E}}_{\mathcal{N}_{\mathcal{F}_0}^*}^\bullet \right)$, we have the associated obstruction class from $(\ref{itt19})$
\begin{align}\label{b5}
\left( \left\{ S_j \right\}, \left\{ Q_j \right\}, \left\{ G_{jk} \right\}, \left\{ C_{jk} \right\} \left\{\left[\left( \theta_{ij}, B_{ij} \right) , \left(\theta_{jk}, B_{jk} \right) \right]  \right\} \right) \in \mathbb{H}^2\left( M, \tilde{\mathcal{E}}_{\mathcal{N}_{\mathcal{F}_0}^*}^\bullet \right)
\end{align}
from \v Cech resolution of $\tilde{\mathcal{E}}_{\mathcal{N}_{\mathcal{F}_0^*}}^\bullet$. Assume that the $1$-cocycle $\left( \left\{ \beta_j \right\}, \left\{ \gamma_j\right\}, \left\{ \left( \theta_{jk}, B_{jk}\right) \right\} \right)$ from the \v Cech resolution of  $\tilde{\mathcal{E}}_{\mathcal{N}_{\mathcal{F}}^*}^\bullet$ corresponds to $\left( \Psi, \tau, \left\{\psi, k_j:=\bar{\partial} B_j \right\} \right)$ from the Dolbeault type bicomplex associated to $\tilde{\mathcal{E}}_{\mathcal{N}_{\mathcal{F}_0}^*}^\bullet$ under the isomorphism $(\ref{b2})$. We recall from \textnormal{Lemma \ref{oui28}} that we have the associated obstruction class from $(\ref{itt22})$
\begin{align}\label{b6}
\left( l, \sigma, \kappa, \pi, \left\{ - \left[\left(\psi, k_j \right), \left(\psi, k_j \right) \right]   \right\} \right) \in \textnormal{C}^2
\end{align}
from the Dolbeault type bicomplex associated to $\tilde{\mathcal{E}}_{\mathcal{N}_{\mathcal{F}_0}^*}^\bullet$. Then $(\ref{b5})$ corresponds to $(\ref{b6})$ under the isomorphism $(\ref{b7})$.
\end{proposition}

\begin{proof}

Let $\left(  \left\{ \beta_j \right\}, \left\{ \gamma_j \right\}, \left\{ \left( \theta_{jk}, B_{jk} \right) \right\} \right)$ be a $1$-cocycle in Lemma \ref{ou50}. Then $\left( \left\{ S_j \right\}, \left\{ Q_j \right\}, \left\{ G_{jk} \right\}, \left\{ C_{jk} \right\}, \right.$ $ \left. \left\{\left[\left( \theta_{ij}, B_{ij} \right) , \left(\theta_{jk}, B_{jk} \right) \right]  \right\} \right) $ defines a $2$-cocycle by Lemma \ref{ou50}. We have the corresponding $1$-cocycle $\left( \Psi, \tau, \left\{\left( \psi, k_j:= \bar{\partial} B_j \right) \right\} \right)$ from the Dolbeault type bicomplex associated to $\tilde{\mathcal{E}}_{\mathcal{N}_{\mathcal{F}_0}^*}^\bullet$. We keep the notations in the proof of Proposition \ref{a24}. Then it remains to show that $\left\{ Q_j \right\}$ corresponds to $\sigma$ and $\left\{ C_{jk}\right\}$ corresponds to $\pi$. In fact, since $\gamma_k\left(w_{0j}^\alpha \right)-\gamma_j \left( w_{0j}^\alpha \right) + i_{\theta_{jk}}\left( w_{0j}^\alpha \right)=0$ and $\beta_k\left(w_{0j}^\alpha \right)- \beta_j\left( w_{0j}^\alpha \right) +\mathcal{L}_{\theta_{jk}}\left( w_{0j}^\alpha \right) - B_{jk}^\alpha=0$, we note that $C_{jk}^\alpha$ is defined by
{\tiny{\begin{align*}
&-i_{\left[ \xi_j, \xi_k\right]}\left( w_{0j}^\alpha \right) - \frac{1}{2} \gamma_j\left( B_{jk}^\alpha \right) - \frac{1}{2}\gamma_k\left(B_{jk}^\alpha \right) + \frac{1}{2}\left[ \theta_{jk}, \gamma_j\left(w_{0j}^\alpha \right) \right] + \frac{1}{2}\left[ \theta_{jk}, \gamma_k\left(w_{0j}^\alpha \right) \right] + \frac{1}{2} i_{\theta_{jk}}\left(\beta_j\left( w_{0j}^\alpha \right) \right) + \frac{1}{2} i_{\theta_{jk}}\left(\beta_k\left(w_{0j}^\alpha \right) \right)\\
&=- i_{\left[ \xi_j, \xi_k\right]}\left( w_{0j}^\alpha \right) - \frac{1}{2} \sum_{\delta=1}^q \left(B_j^{\alpha \delta} - \sum_{\beta, \eta=1}^q h_{0jk}^{\alpha \beta} B_k^{\beta\eta} h_{0kj}^{\eta \delta} - \sum_{\gamma, \eta=1}^q \left[ \xi_k, h_{0jk}^{\alpha \gamma}\right] h_{0kj}^{\gamma \delta} \right)\gamma_j\left(  w_{0j}^\delta \right)\\
&- \frac{1}{2}\sum_{\delta=1}^q \left(B_j^{\alpha \delta} - \sum_{\beta, \eta=1}^q h_{0jk}^{\alpha \beta} B_k^{\beta\eta} h_{0kj}^{\eta \delta} - \sum_{\gamma, \eta=1}^q \left[ \xi_k, h_{0jk}^{\alpha \gamma}\right] h_{0kj}^{\gamma \delta} \right) \gamma_k\left(w_{0j}^\delta \right) + \frac{1}{2}\left[ \xi_j- \xi_k, \gamma_j\left(w_{0j}^\alpha \right) \right] + \frac{1}{2}\left[ \xi_j- \xi_k, \gamma_k\left(w_{0j}^\alpha \right) \right] \\
&+ \frac{1}{2} i_{\xi_j-\xi_k}\left(\beta_j\left(w_{0j}^\alpha \right) \right) + \frac{1}{2} i_{\xi_j- \xi_k}\left(\beta_k\left(w_{0j}^\alpha \right) \right)\\
&=-i_{\left[ \xi_j, \xi_k \right]}\left( w_{0j}^\alpha \right) - \frac{1}{2} \sum_{\delta=1}^q B_j^{\alpha \delta} \gamma_j\left( w_{0j}^\delta \right) + \frac{1}{2}\sum_{\delta, \beta, \eta=1}^q h_{0jk}^{\alpha\beta} B_k^{\beta \eta} h_{0kj}^{\eta \delta} \left( \gamma_k\left( w_{0j}^\alpha \right) + i_{\xi_j- \xi_k}\left(w_{0j}^\alpha \right) \right) + \frac{1}{2}\sum_{\delta, \gamma =1}^q \left[ \xi_k, h_{0jk}^{\alpha \gamma} \right] h_{0kj}^{\gamma \delta} \left(\gamma_k\left(w_{0j}^\delta \right) + i_{\xi_j-\xi_k}\left(w_{0j}^\delta \right) \right)\\
&- \frac{1}{2}\sum_{\delta=1}^q B_j^{\alpha \delta}\left(\gamma_j\left(w_{0j}^\delta \right)  - i_{\xi_j-\xi_k}\left(w_{0j}^\delta \right)\right) + \frac{1}{2}\sum_{\beta,\eta=1}^q h_{0jk}^{\alpha\beta}B_k^{\beta \eta} \gamma_k\left(w_{0k}^\eta \right) + \frac{1}{2}\sum_{\gamma, \eta=1}^q \left[ \xi_k, h_{0jk}^{\alpha \gamma} \right]  \gamma_k\left( w_{0k}^\gamma \right) + \frac{1}{2}\left[ \xi_j, \gamma_j\left(w_{0j}^\alpha \right) \right] \\
&-\frac{1}{2}\left[ \xi_k, \gamma_k\left(w_{0j}^\alpha \right)+ i_{\xi_j-\xi_k}\left(w_{0j}^\alpha \right)\right] + \frac{1}{2}\left[ \xi_j, \gamma_j\left(w_{0j}^\alpha \right) - i_{\xi_j-\xi_k}\left(w_{0j}^\alpha \right) \right]  - \frac{1}{2}\left[ \xi_k, \gamma_k\left(w_{0j}^\alpha \right) \right]+ \frac{1}{2} i_{\xi_j}\left(\beta_j\left(w_{0j}^\alpha \right) \right) \\
&- \frac{1}{2}i_{\xi_k}\left(\beta_k\left(w_{0j}^\alpha \right) + \mathcal{L}_{\xi_j-\xi_k}\left(w_{0j}^\alpha \right) - \sum_{\delta=1}^q\left( B_j^{\alpha \delta} - \sum_{\beta, \eta=1}^q h_{0jk}^{\alpha \beta} B_k^{\beta\eta} h_{0kj}^{\eta \delta} - \sum_{\gamma, \eta=1}^q \left[ \xi_k, h_{0jk}^{\alpha \gamma}\right] h_{0kj}^{\gamma \delta} \right)w_{0j}^\delta \right) -  \frac{1}{2}\sum_{\beta=1}^q h_{0jk}^{\alpha\beta}i_{\xi_k}\left(\beta_k\left( w_{0k}^\alpha \right) \right)\\
&+\frac{1}{2}i_{\xi_j}\left(\beta_j\left(w_{0j}^\alpha \right) - \mathcal{L}_{\xi_j-\xi_k}\left(w_{0j}^\alpha \right) + \sum_{\delta=1}^q \left( B_j^{\alpha \delta} - \sum_{\beta, \eta=1}^q h_{0jk}^{\alpha \beta} B_k^{\beta\eta} h_{0kj}^{\eta \delta} - \sum_{\gamma, \eta=1}^q \left[ \xi_k, h_{0jk}^{\alpha \gamma}\right] h_{0kj}^{\gamma \delta} \right)w_{0j}^\delta\right)
\end{align*}}}
We continue the computation. Then
{\tiny{\begin{align*}
&=- \frac{1}{2}i_{\left[ \xi_j, \xi_k\right]}\left( w_{0j}^\alpha \right) + \frac{1}{2} i_{\left[ \xi_k, \xi_j\right]}\left( w_{0j}^\alpha \right) - \sum_{\delta=1}^q B_j^{\alpha \delta}\gamma_j\left(w_{0j}^\delta \right) + \sum_{\beta, \eta=1}^q h_{0jk}^{\alpha \beta} B_k^{\beta \eta}\gamma_k\left(w_{0k}^\eta \right) + \frac{1}{2}\sum_{\delta, \beta, \eta=1}^q h_{0jk}^{\alpha \beta} B_k^{\beta \eta} i_{\xi_j}\left(w_{0k}^\eta  \right) -\frac{1}{2} \sum_{\beta, \eta=1}^q h_{0jk}^{\alpha \beta} B_k^{\beta \eta}i_{\xi_k}\left(w_{0k}^\eta \right) \\
& + \sum_{ \gamma=1}^q \left[ \xi_k, h_{0jk}^{\alpha \gamma} \right] \gamma_k\left( w_{0k}^\gamma  \right) + \frac{1}{2}\sum_{\gamma, \eta=1}^q \left[ \xi_k, h_{0jk}^{\alpha\gamma} \right] i_{\xi_j}\left(w_{0k}^\gamma \right) - \frac{1}{2}\sum_{\gamma, \eta=1}^q \left[ \xi_k, h_{0jk}^{\alpha\gamma} \right] i_{\xi_k}\left(w_{0k}^\gamma \right) + \frac{1}{2} \sum_{\delta=1}^q B_j^{\alpha \delta} i_{\xi_j}\left(w_{0j}^\delta \right) - \frac{1}{2} \sum_{\delta=1}^q B_j^{\alpha \delta} i_{\xi_k}\left(w_{0j}^\delta \right)\\
&+\left[ \xi_j, \gamma_j\left(w_{0j}^\alpha \right)\right] - \left[ \xi_k, \gamma_k\left( w_{0j}^\alpha \right) \right] - \frac{1}{2}\left[ \xi_k, i_{\xi_j}\left(w_{0j}^\alpha \right) \right] + \frac{1}{2}\left[ \xi_k, i_{\xi_k}\left(w_{0j}^\alpha \right) \right] - \frac{1}{2} \left[ \xi_j, i_{\xi_j}\left( w_{0j}^\alpha \right) \right] + \frac{1}{2}\left[ \xi_j, i_{\xi_k}\left( w_{0j}^\alpha \right) \right] + i_{\xi_j}\left(\beta_j\left(w_{0j}^\alpha \right) \right) \\
&- \sum_{\beta=1}^q h_{0jk}^{\alpha\beta} i_{\xi_k}\left(\beta_k\left(w_{0k}^\beta \right) \right) - \frac{1}{2}i_{\xi_k}\mathcal{L}_{\xi_j}\left(w_{0j}^\alpha \right) + \frac{1}{2}i_{\xi_k}\mathcal{L}_{\xi_k}\left(w_{0j}^\alpha \right) + \frac{1}{2}\sum_{\delta=1}^q B_j^{\alpha \delta} i_{\xi_k}\left(w_{0j}^\delta \right) - \frac{1}{2}\sum_{\beta, \eta=1}^q h_{0jk}^{\alpha\beta} B_k^{\beta \eta}i_{\xi_k}\left( w_{0k}^\eta\right) - \frac{1}{2}\sum_{\beta=1}^q\left[ \xi_k, h_{0jk}^{\alpha \beta} \right]i_{\xi_k}\left(w_{0k}^\beta \right) \\
&- \frac{1}{2} i_{\xi_j} \mathcal{L}_{\xi_j}\left( w_{0j}^\alpha \right) + \frac{1}{2} i_{\xi_j}\mathcal{L}_{\xi_k}\left(w_{0j}^\alpha \right) + \frac{1}{2} \sum_{\delta=1}^q B_j^{\alpha \delta} i_{\xi_j}\left(w_{0j}^\delta \right) - \frac{1}{2}\sum_{\beta, \eta=1}^q h_{0jk}^{\alpha\beta} B_k^{\beta \eta} i_{\xi_j}\left( w_{0k}^\eta\right) - \frac{1}{2} \sum_{\beta=1}^q \left[ \xi_k, h_{0jk}^{\alpha \beta}\right] i_{\xi_j}\left(w_{0k}^\beta \right)
\end{align*}}}
We continue the computation. Then
{\tiny{\begin{align*}
&= \cancel{-\frac{1}{2}\mathcal{L}_{\xi_j}i_{\xi_k}\left(w_{0j}^\alpha \right)} + \bcancel{\frac{1}{2}i_{\xi_k} \mathcal{L}_{\xi_j}\left(w_{0j}^\alpha \right)}+  \xcancel{\frac{1}{2}\mathcal{L}_{\xi_k} i_{\xi_j}\left(w_{0j}^\alpha \right)} - \cancel{ \frac{1}{2} i_{\xi_j} \mathcal{L}_{\xi_k} \left(w_{0j}^\alpha \right)} - \sum_{\delta=1}^q B_j^{\alpha \delta}\gamma_j\left(w_{0j}^\delta \right) + \sum_{\beta, \eta=1}^q h_{0jk}^{\alpha \beta} B_k^{\beta \eta}\gamma_k\left(w_{0k}^\eta \right) \\
& - \sum_{\beta, \eta=1}^q h_{0jk}^{\alpha \beta} B_k^{\beta \eta}i_{\xi_k}\left(w_{0k}^\eta \right)+ \sum_{ \gamma=1}^q \left[ \xi_k, h_{0jk}^{\alpha \gamma} \right] \gamma_k\left( w_{0k}^\gamma  \right) -\sum_{\gamma, \eta=1}^q \left[ \xi_k, h_{0jk}^{\alpha\gamma} \right] i_{\xi_k}\left(w_{0k}^\gamma \right) +  \sum_{\delta=1}^q B_j^{\alpha \delta} i_{\xi_j}\left(w_{0j}^\delta \right) \\
&+\left[ \xi_j, \gamma_j\left(w_{0j}^\alpha \right)\right] - \left[ \xi_k, \gamma_k\left( w_{0j}^\alpha \right) \right] - \xcancel{\frac{1}{2}\left[ \xi_k, i_{\xi_j}\left(w_{0j}^\alpha \right) \right]} + \frac{1}{2}\left[ \xi_k, i_{\xi_k}\left(w_{0j}^\alpha \right) \right] - \frac{1}{2} \left[ \xi_j, i_{\xi_j}\left( w_{0j}^\alpha \right) \right] + \cancel{ \frac{1}{2}\left[ \xi_j, i_{\xi_k}\left( w_{0j}^\alpha \right) \right] }+ i_{\xi_j}\left(\beta_j\left(w_{0j}^\alpha \right) \right) \\
&- \sum_{\beta=1}^q h_{0jk}^{\alpha\beta} i_{\xi_k}\left(\beta_k\left(w_{0k}^\beta \right) \right) - \bcancel{ \frac{1}{2}i_{\xi_k}\mathcal{L}_{\xi_j}\left(w_{0j}^\alpha \right) } + \frac{1}{2}i_{\xi_k}\mathcal{L}_{\xi_k}\left(w_{0j}^\alpha \right)  - \frac{1}{2} i_{\xi_j} \mathcal{L}_{\xi_j}\left( w_{0j}^\alpha \right) + \cancel{ \frac{1}{2} i_{\xi_j}\mathcal{L}_{\xi_k}\left(w_{0j}^\alpha \right)   }\\
&= - \sum_{\delta=1}^q B_j^{\alpha \delta}\gamma_j \left(w_{0j}^\delta \right) + \sum_{\beta,\eta=1}^q h_{0jk}^{\alpha \beta}B_k^{\beta \eta} \gamma_k\left(w_{0k}^\eta \right) - \sum_{\beta, \eta=1}^q h_{0jk}^{\alpha\beta} B_k^{\beta \eta} i_{\xi_k}\left(w_{0k}^\eta \right) - \sum_{\beta=1}^q h_{0jk}^{\alpha\beta} \left[ \xi_k, \gamma_k\left(w_{0k}^\beta \right) \right] + \sum_{\beta=1}^q h_{0jk}^{\alpha\beta} i_{\xi_k} \mathcal{L}_{\xi_k}\left(w_{0k}^\beta \right) \\
&+ \sum_{\delta=1}^q B_j^{\alpha \delta} i_{\xi_j}\left(w_{0j}^\delta \right) + \left[ \xi_j, \gamma_j\left(w_{0j}^\alpha \right) \right] - i_{\xi_j}\mathcal{L}_{\xi_j}\left( w_{0j}^\alpha \right) - \sum_{\beta=1}^q h_{0jk}^{\alpha\beta} i_{\xi_k}\left( \beta_k\left( w_{0k}^\beta \right)\right) + i_{\xi_j}\left(\beta_j\left(w_{0j}^\alpha \right) \right)\\
&= - \sum_{\delta=1}^q B_j^{\alpha \delta}\gamma_j \left(w_{0j}^\delta \right) + \sum_{\delta=1}^q B_j^{\alpha \delta} i_{\xi_j}\left(w_{0j}^\delta \right) + \left[ \xi_j, \gamma_j\left(w_{0j}^\alpha \right) \right] - i_{\xi_j}\mathcal{L}_{\xi_j}\left( w_{0j}^\alpha \right)  + i_{\xi_j}\left(\beta_j\left(w_{0j}^\alpha \right) \right)\\
 &- \sum_{\beta=1}^q h_{0jk}^{\alpha\beta}\left(  - \sum_{\eta=1}^q B_k^{\beta \eta} \gamma_k\left( w_{0k}^\eta \right)  + \sum_{\eta=1}^q B_k^{\beta \eta} i_{\xi_k}\left( w_{0k}^\eta \right) + \left[ \xi_k, \gamma_k\left( w_{0k}^\beta \right) \right]  - i_{\xi_k} \mathcal{L}_{\xi_k}\left( w_{0k}^\beta \right) + i_{\xi_k}\left( \beta_k\left( w_{0k}^\beta \right) \right)\right)
\end{align*}}}

Then we set
\begin{align*}
r_j^\alpha:&=  - \gamma_j\left( B_j^\alpha \right) + i_{\xi_j}\left(B_j^\alpha \right) + \left[\xi_j, \gamma_j\left( w_{0j}^\alpha \right) \right] - i_{\xi_j}\mathcal{L}_{\xi_j}\left( w_{0j}^\alpha \right) + i_{\xi_j}\left( \beta_j\left( w_{0j}^\alpha \right) \right)\\
 &= - \sum_{\delta=1}^q B_j^{\alpha \delta}\gamma_j \left(w_{0j}^\delta \right) + \sum_{\delta=1}^q B_j^{\alpha \delta} i_{\xi_j}\left(w_{0j}^\delta \right) + \left[ \xi_j, \gamma_j\left(w_{0j}^\alpha \right) \right] - i_{\xi_j}\mathcal{L}_{\xi_j}\left( w_{0j}^\alpha \right)  + i_{\xi_j}\left(\beta_j\left(w_{0j}^\alpha \right) \right)
\end{align*}
Then we note that
{\small{\begin{align*}
&- i_{p_j}\left(w_{0j}^\alpha \right) + \bar{\partial} r_j^\alpha\\
&= - i_{\left[ \psi, \xi_j \right]} \left( w_{0j}^\alpha \right) - \sum_{\delta=1}^q k_j^{\alpha \delta} \gamma_j \left( w_{0j}^\delta \right)+ \sum_{\delta=1}^q k_j^{\alpha \delta} i_{\xi_j}\left( w_{0j}^\delta \right) + i_{\psi}\left(B_j^\alpha \right) +\left[ \psi, \gamma_j\left(w_{0j}^\alpha \right) \right]- i_{\psi} \mathcal{L}_{\xi_j}\left(w_{0j}^\alpha \right) - i_{\xi_j}\mathcal{L}_{\psi}\left( w_{0j}^\alpha \right) + i_{\psi}\left( \beta_j\left( w_{0j}^\alpha \right) \right)\\
&= - \mathcal{L}_\psi i_{\xi_j}\left(w_{0j}^\alpha \right) + \cancel{ i_{\xi_j}\mathcal{L}_{\psi}\left(w_{0j}^\alpha \right) } + \sum_{\delta=1}^q k_j^{\alpha \delta}\left( i_{\xi_j}\left( w_{0j}^\delta \right) - \gamma_j\left( w_{0j}^\delta \right) \right) - i_{\psi}\left( \mathcal{L}_{\xi_j}\left( w_{0j}^\alpha \right) - B_j^\alpha - \beta_j\left( w_{0j}^\alpha \right) \right) +\left[ \psi, \gamma_j\left(w_{0j}^\alpha \right)\right]  -  \cancel{ i_{\xi_j} \mathcal{L}_{\psi}\left( w_{0j}^\alpha \right)} \\
&= - \left[ \psi, i_{\xi_j}\left( w_{0j}^\alpha \right) - \gamma_j\left( w_{0j}^\alpha \right) \right] + \sum_{\delta=1}^q k_j^{\alpha \delta} \tau_j^\delta - i_{\psi}\left(\Psi_j^\alpha \right) = - \left[ \psi, \tau_j^\alpha \right] - i_{\psi}\left( \Psi_j^\alpha \right) - \sum_{\delta=1}^q k_j^{\alpha \delta} \tau_j^\delta
\end{align*}}}
This implies that $\left\{ C_{jk}\right\}$ corresponds to $\pi$. On the other hand, we show that $\left\{Q_j \right\}$ corresponds to $\sigma$. In fact, we note that
{\small{\begin{align}
&Q_j^\alpha + \sum_{\beta=1}^q (-1)^{q-\beta} r_j^\beta w_{0j}^1 \wedge \cdots \wedge \hat{w}_{0j}^\beta \wedge \cdots \wedge w_{0j}^q \wedge dw_{0j}^\alpha + w_{0j}^1 \wedge \cdots \wedge w_{0j}^q\wedge \partial r_j^\alpha  - w_{0j}^1 \wedge \cdots \wedge w_{0j}^q \wedge f_j^\alpha \label{b8}\\
&=\sum_{\gamma=1}^q\sum_{1\leq \beta \ne \gamma \leq q } (-1)^{q-\gamma} \gamma_j \left(w_{0j}^\gamma \right) w_{0j}^1\wedge \cdots \wedge \beta_j \left( w_{0j}^\beta \right)  \wedge \cdots \wedge \hat{w}_{0j}^\gamma \wedge \cdots \wedge w_{0j}^q\wedge dw_{0j}^\alpha \notag \\
&+\sum_{\gamma=1}^q (-1)^{q-\gamma} \gamma_j\left(w_{0j}^\gamma\right) w_{0j}^1\wedge \cdots \wedge \hat{w}_{0j}^\gamma \wedge \cdots w_{0j}^q \wedge d\left( \beta_j\left(w_{0j}^\alpha \right)\right) + \sum_{\gamma=1}^q w_{0j}^1\wedge \cdots \wedge \beta_j \left( w_{0j}^\gamma \right)\wedge \cdots \wedge w_{0j}^q\wedge d \left(\gamma_j\left(w_{0j}^\alpha \right) \right) \notag \\
& -  \sum_{\gamma=1}^q w_{0j}^1\wedge \cdots \wedge  \beta_j \left(w_{0j}^\gamma \right)\wedge \cdots \wedge w_{0j}^q\wedge \beta_j\left( w_{0j}^\alpha \right) \notag \\
& + \sum_{\beta=1}^q (-1)^{q-\beta}\left(  - \gamma_j\left( B_j^\beta \right) + i_{\xi_j}\left(B_j^\beta \right) + \left[\xi_j, \gamma_j\left( w_{0j}^\beta \right) \right] - i_{\xi_j}\mathcal{L}_{\xi_j}\left( w_{0j}^\beta \right) + i_{\xi_j}\left( \beta_j\left( w_{0j}^\beta \right) \right) \right) w_{0j}^1 \wedge \cdots \wedge \hat{w}_{0j}^\beta \wedge \cdots \wedge w_{0j}^q\wedge dw_{0j}^\alpha \notag \\
&+ w_{0j}^1 \wedge \cdots \wedge w_{0j}^q \wedge \partial\left(  - \gamma_j\left( B_j^\alpha \right) + i_{\xi_j}\left(B_j^\alpha \right) + \left[\xi_j, \gamma_j\left( w_{0j}^\alpha \right) \right] - i_{\xi_j}\mathcal{L}_{\xi_j}\left( w_{0j}^\alpha \right) + i_{\xi_j}\left( \beta_j\left( w_{0j}^\alpha \right) \right) \right) \notag \\
&- w_{0j}^1 \wedge \cdots \wedge w_{0j}^q \wedge \left( - 2\beta_j\left( B_j^\alpha \right) + 2\mathcal{L}_{\xi_j}\left( \beta_j\left( w_{0j}^\alpha \right) \right) + 2 \mathcal{L}_{\xi_j}\left( B_j^\alpha \right) - \mathcal{L}_{\xi_j} \mathcal{L}_{\xi_j}\left( w_{0j}^\alpha \right) \right) \notag
\end{align}}}

On the other hand, in a similar way to $(\ref{ou41}), (\ref{ou42}), (\ref{ou43}),(\ref{ou45})$ and $(\ref{ou47})$,

{\Small{\begin{align}
&\sum_{\gamma=1}^q\sum_{1\leq \beta \ne \gamma \leq q } (-1)^{q-\gamma} \tau_j^\gamma\cdot w_{0j}^1\wedge \cdots \wedge \Psi_j^\beta \wedge \cdots \wedge \hat{w}_{0j}^\gamma \wedge \cdots \wedge w_{0j}^q\wedge dw_{0j}^\alpha \label{b9} \\
&+\sum_{\gamma=1}^q (-1)^{q-\gamma} \tau_j^\gamma \wedge \cdots \wedge \hat{w}_{0j}^\gamma \wedge \cdots w_{0j}^q \wedge \partial \left( \Psi_j^\alpha \right) + \sum_{\gamma=1}^q w_{0j}^1\wedge \cdots \wedge  \Psi_j^\gamma \wedge \cdots \wedge w_{0j}^q\wedge \partial \left(\tau_j^\alpha  \right) \notag \\
& -  \sum_{\gamma=1}^q w_{0j}^1\wedge \cdots \wedge  \Psi_j^\gamma \wedge \cdots \wedge w_{0j}^q\wedge \Psi_j^\alpha \notag\\
&=\sum_{\gamma=1}^q\sum_{1\leq \beta \ne \gamma \leq q } (-1)^{q-\gamma} \left(i_{\xi_j}\left( w_{0j}^\gamma \right) - \gamma_j\left(w_{0j}^\gamma \right)\right)\cdot w_{0j}^1\wedge \cdots \wedge \left( \mathcal{L}_{\xi_j}\left(  w_{0j}^\beta \right) - B_j^\beta - \beta_j\left( w_{0j}^\beta \right) \right) \wedge \cdots \wedge \hat{w}_{0j}^\gamma \wedge \cdots \wedge w_{0j}^q\wedge dw_{0j}^\alpha \notag  \\
&+\sum_{\gamma=1}^q (-1)^{q-\gamma} \left( i_{\xi_j}\left( w_{0j}^\gamma \right) - \gamma_j\left(w_{0j}^\gamma \right)\right) w_{0j}^1 \wedge \cdots \wedge \hat{w}_{0j}^\gamma \wedge \cdots w_{0j}^q \wedge \partial \left( \mathcal{L}_{\xi_j}\left(  w_{0j}^\alpha \right) - B_j^\alpha - \beta_j\left( w_{0j}^\alpha \right) \right) \notag \\
&+ \sum_{\gamma=1}^q w_{0j}^1\wedge \cdots \wedge  \left( \mathcal{L}_{\xi_j}\left(  w_{0j}^\gamma \right) - B_j^\gamma - \beta_j\left( w_{0j}^\gamma \right)\right) \wedge \cdots \wedge w_{0j}^q\wedge \partial \left( i_{\xi_j}\left( w_{0j}^\alpha \right) - \gamma_j\left(w_{0j}^\alpha \right)  \right) \notag \\
& -  \sum_{\gamma=1}^q w_{0j}^1\wedge \cdots \wedge  \left( \mathcal{L}_{\xi_j}\left(  w_{0j}^\gamma \right) - B_j^\gamma - \beta_j\left( w_{0j}^\gamma \right) \right) \wedge \cdots \wedge w_{0j}^q\wedge \left(\mathcal{L}_{\xi_j}\left(  w_{0j}^\alpha \right) - B_j^\alpha - \beta_j\left( w_{0j}^\alpha \right) \right)\notag
\end{align}}}
We continue the computation. Then
{\Small{\begin{align}
&=\sum_{\gamma=1}^q\sum_{1\leq \beta \ne \gamma \leq q } (-1)^{q-\gamma} \gamma_j \left(w_{0j}^\gamma \right) w_{0j}^1\wedge \cdots \wedge \beta_j \left( w_{0j}^\beta \right)  \wedge \cdots \wedge \hat{w}_{0j}^\gamma \wedge \cdots \wedge w_{0j}^q\wedge dw_{0j}^\alpha -  \sum_{\gamma=1}^q w_{0j}^1\wedge \cdots \wedge  \beta_j \left(w_{0j}^\gamma \right)\wedge \cdots \wedge w_{0j}^q\wedge \beta_j\left( w_{0j}^\alpha \right)  \notag \\
&+\sum_{\gamma=1}^q (-1)^{q-\gamma} \gamma_j\left(w_{0j}^\gamma\right) w_{0j}^1\wedge \cdots \wedge \hat{w}_{0j}^\gamma \wedge \cdots w_{0j}^q \wedge d\left( \beta_j\left(w_{0j}^\alpha \right)\right) + \sum_{\gamma=1}^q w_{0j}^1\wedge \cdots \wedge \beta_j \left( w_{0j}^\gamma \right)\wedge \cdots \wedge w_{0j}^q\wedge d \left(\gamma_j\left(w_{0j}^\alpha \right) \right)  \notag\\
& - w_{0j}^1 \wedge \cdots \wedge w_{0j}^q \wedge \partial \gamma_j\left( B_j^\alpha \right) + w_{0j}^1 \wedge \cdots \wedge w_{0j}^q \wedge \beta_j\left( B_j^\alpha \right) + w_{0j}^1 \wedge \cdots \wedge w_{0j}^q \wedge \beta_j\left(B_j^\alpha \right) \notag  \\
&+\sum_{\gamma=1}^q (-1)^{q-\gamma} \left[\xi_j, \gamma_j\left( w_{0j}^\gamma \right) \right] w_{0j}^1 \wedge \cdots \wedge \widehat{w}_{0j}^\gamma \wedge \cdots \wedge w_{0j}^q \wedge dw_{0j}^\alpha + w_{0j}^1 \wedge \cdots \wedge w_{0j}^q \wedge \partial\left(\left[ \xi_j, \gamma_j\left( w_{0j}^\alpha \right) \right] \right)  - w_{0j}^1 \wedge \cdots \wedge w_{0j}^q \wedge \mathcal{L}_{\xi_j}\left( \beta_j\left( w_{0j}^\alpha \right) \right) \notag\\
&+ \sum_{\beta=1}^q (-1)^{q-\beta} i_{\xi_j}\left(\beta_j\left( w_{0j}^\beta \right) \right)w_{0j}^1 \wedge \cdots \widehat{w}_{0j}^\beta \wedge \cdots \wedge w_{0j}^q \wedge dw_{0j}^\alpha - w_{0j}^1 \wedge \cdots \wedge w_{0j}^q \wedge \left( \mathcal{L}_{\xi_j}\left(\beta_j\left( w_{0j}^\alpha \right) \right) \right) + w_{0j}^1 \wedge \cdots \wedge w_{0j}^q \wedge \partial \left(i_{\xi_j}\left(\beta_j\left( w_{0j}^\alpha \right) \right) \right) \notag\\
&- \sum_{\beta=1}^q (-1)^{q-\beta} \mathcal{L}_{\xi_j}i_{\xi_j}\left(w_{0j}^\beta \right) w_{0j}^1 \wedge \cdots \wedge \hat{w}_{0j}^\beta \wedge \cdots \wedge w_{0j}^q \wedge dw_{0j}^\alpha + w_{0j}^1 \wedge \cdots \wedge w_{0j}^q \wedge \mathcal{L}_{\xi_j} \mathcal{L}_{\xi_j}\left( w_{0j}^\alpha \right) - w_{0j}^1 \wedge \cdots \wedge w_{0j}^q \wedge \partial \mathcal{L}_{\xi_j}i_{\xi_j}\left( w_{0j}^\alpha \right) \notag \\
&- \sum_{\gamma=1}^q \sum_{1\leq \beta \ne \gamma\leq q}(-1)^{q-\gamma} i_{\xi_j}\left( w_{0j}^\gamma \right) w_{0j}^1 \wedge \cdots \wedge \widehat{w}_{0j}^\gamma \wedge \cdots \wedge  B_j^\beta \wedge \cdots \wedge w_{0j}^q \wedge \cdots \wedge dw_{0j}^\alpha \notag \\
& + \sum_{\gamma=1}^q \sum_{1\leq \beta \ne \gamma\leq q}(-1)^{q-\gamma} \gamma_j \left( w_{0j}^\gamma \right) w_{0j}^1 \wedge \cdots \wedge \widehat{w}_{0j}^\gamma \wedge \cdots \wedge B_j^\beta \wedge \cdots \wedge w_{0j}^q \wedge \cdots \wedge dw_{0j}^\alpha \notag \\
&- \sum_{\gamma=1}^q (-1)^{q-\gamma} i_{\xi_j}\left( w_{0j}^\gamma \right) w_{0j}^1 \wedge \cdots \wedge \widehat{w}_{0j}^\gamma \wedge \cdots \wedge w_{0j}^q \wedge \partial B_j^\alpha + \sum_{\gamma=1}^q w_{0j}^1 \wedge \cdots \wedge \mathcal{L}_{\xi_j}\left( w_{0j}^\gamma \right) \wedge \cdots \wedge w_{0j}^q \wedge B_j^\alpha \notag \\
& + \sum_{\gamma=1}^q w_{0j}^1 \wedge \cdots \wedge B_j^\gamma \wedge \cdots \wedge w_{0j}^q \wedge \mathcal{L}_{\xi_j}\left( w_{0j}^\alpha \right) - \sum_{\gamma=1}^q w_{0j}^1 \wedge \cdots \wedge B_j^\gamma \wedge \cdots \wedge w_{0j}^q \wedge \partial i_{\xi_j}\left( w_{0j}^\alpha \right) + \sum_{\gamma=1}^q w_{0j}^1 \wedge \cdots \wedge B_j^\gamma \wedge \cdots \wedge w_{0j}^q \wedge \partial \gamma_j\left(w_{0j}^\alpha \right) \notag \\
& - \sum_{\gamma=1}^q w_{0j}^1 \wedge \cdots \wedge B_j^\gamma \wedge \cdots \wedge w_{0j}^q \wedge \beta_j\left(  w_{0j}^\alpha \right) \notag
\end{align}}}

In view of $(\ref{b8})$ and $(\ref{b9})$, we have to show that
{\tiny{\begin{align}
&\sum_{\beta=1}^q (-1)^{q-\beta}\left( -\gamma_j\left( B_j^\beta \right) + i_{\xi_j}\left(B_j^\beta \right) \right) w_{0j}^1 \wedge \cdots \wedge \widehat{w}_{0j}^\beta \wedge \cdots  \wedge w_{0j}^q \wedge dw_{0j}^\alpha + w_{0j}^1 \wedge \cdots \wedge w_{0j}^q \wedge \partial \left( i_{\xi_j} \left( B_j^\alpha  \right) \right) - 2 w_{0j}^1 \wedge \cdots \wedge w_{0j}^q \wedge \mathcal{L}_{\xi_j}\left(B_j^\alpha \right) \label{b13}\\
&= - \sum_{\gamma=1}^q \sum_{1\leq \beta \ne \gamma\leq q}(-1)^{q-\gamma} i_{\xi_j}\left( w_{0j}^\gamma \right) w_{0j}^1 \wedge \cdots \wedge \widehat{w}_{0j}^\gamma \wedge \cdots \wedge  B_j^\beta \wedge \cdots \wedge w_{0j}^q \wedge \cdots \wedge dw_{0j}^\alpha \notag \\
& + \sum_{\gamma=1}^q \sum_{1\leq \beta \ne \gamma\leq q}(-1)^{q-\gamma} \gamma_j \left( w_{0j}^\gamma \right) w_{0j}^1 \wedge \cdots \wedge \widehat{w}_{0j}^\gamma \wedge \cdots \wedge B_j^\beta \wedge \cdots \wedge w_{0j}^q \wedge \cdots \wedge dw_{0j}^\alpha \notag \\
&- \sum_{\gamma=1}^q (-1)^{q-\gamma} i_{\xi_j}\left( w_{0j}^\gamma \right) w_{0j}^1 \wedge \cdots \wedge \widehat{w}_{0j}^\gamma \wedge \cdots \wedge w_{0j}^q \wedge \partial B_j^\alpha + \sum_{\gamma=1}^q w_{0j}^1 \wedge \cdots \wedge \mathcal{L}_{\xi_j}\left( w_{0j}^\gamma \right) \wedge \cdots \wedge w_{0j}^q \wedge B_j^\alpha + \sum_{\gamma=1}^q w_{0j}^1 \wedge \cdots \wedge B_j^\gamma \wedge \cdots \wedge w_{0j}^q \wedge \mathcal{L}_{\xi_j}\left( w_{0j}^\alpha \right) \notag\\
&- \sum_{\gamma=1}^q w_{0j}^1 \wedge \cdots \wedge B_j^\gamma \wedge \cdots \wedge w_{0j}^q \wedge \partial i_{\xi_j}\left( w_{0j}^\alpha \right) + \sum_{\gamma=1}^q w_{0j}^1 \wedge \cdots \wedge B_j^\gamma \wedge \cdots \wedge w_{0j}^q \wedge \partial \gamma_j\left(w_{0j}^\alpha \right) - \sum_{\gamma=1}^q w_{0j}^1 \wedge \cdots \wedge B_j^\gamma \wedge \cdots \wedge w_{0j}^q \wedge \beta_j\left(  w_{0j}^\alpha \right) \notag
\end{align}}}

In fact, we note that
{\Small{\begin{align}
&\sum_{\gamma=1}^q w_{0j}^1 \wedge \cdots \wedge B_j^\gamma \wedge \cdots \wedge w_{0j}^q \wedge d w_{0j}^\alpha + w_{0j}^1 \wedge \cdots \wedge w_{0j}^q \wedge \partial B_j^\alpha = 0  \notag\\
&\Longrightarrow \sum_{\gamma=1}^q \sum_{1\leq \beta \ne \gamma \leq q} (-1)^{\beta-1} i_{\xi_j}\left( w_{0j}^\beta \right) w_{0j}^1 \wedge \cdots \wedge \widehat{w}_{0j}^\beta \wedge \cdots \wedge B_j^\gamma\wedge \cdots \wedge w_{0j}^q \wedge dw_{0j}^\alpha + \sum_{\gamma=1}^q (-1)^{\gamma -1} i_{\xi_j}\left( B_j^\gamma \right) w_{0j}^1 \wedge \cdots \wedge \widehat{w}_{0j}^\gamma \wedge \cdots \wedge w_{0j}^q \wedge dw_{0j}^\alpha \label{b10}\\ 
&+\sum_{\gamma=1}^q (-1)^q w_{0j}^1 \wedge \cdots \wedge B_j^\gamma \wedge \cdots \wedge w_{0j}^q \wedge \mathcal{L}_{\xi_j}\left( w_{0j}^\alpha \right) - \sum_{\gamma=1}^q (-1)^q w_{0j}^1 \wedge \cdots \wedge B_j^\gamma \wedge \cdots \wedge w_{0j}^q \wedge \partial \left(i_{\xi_j}\left(  w_{0j}^\alpha \right) \right) \notag\\
& + \sum_{\gamma=1}^q (-1)^{\gamma-1} i_{\xi_j}\left( w_{0j}^\gamma \right) w_{0j}^1 \wedge \cdots \wedge \widehat{w}_{0j}^\gamma \wedge \cdots \wedge w_{0j}^q \wedge \partial B_j^\alpha + (-1)^q w_{0j}^1 \wedge \cdots \wedge w_{0j}^q \wedge \mathcal{L}_{\xi_j}\left(B_j^\alpha \right) - (-1)^q w_{0j}^1 \wedge \cdots \wedge w_{0j}^q \wedge \partial \left( i_{\xi_j} \left(B_j^\alpha \right) \right)=0 \notag
\end{align}}}

and we note that
\begin{align}
&w_{0j}^1 \wedge \cdots \wedge w_{0j}^q \wedge B_j^\alpha =0  \notag \\
&\Longrightarrow \sum_{\gamma=1}^q w_{0j}^1 \wedge \cdots \wedge \mathcal{L}_{\xi_j} \left( w_{0j}^\gamma \right) \wedge \cdots \wedge w_{0j}^q \wedge B_j^\alpha + w_{0j}^1 \wedge \cdots \wedge w_{0j}^q \wedge \mathcal{L}_{\xi_j} \left( B_j^\alpha \right) = 0 \label{b11}
\end{align}

and we note that
{\tiny{\begin{align}
&\sum_{\beta=1}^q (-1)^{q-\beta} \gamma_j\left( B_j^\beta\right) w_{0j}^1 \wedge \cdots \wedge \widehat{w}_{0j}^\beta \wedge \cdots \wedge w_{0j}^q \wedge dw_{0j}^\alpha \label{b12}\\
&= \sum_{\beta=1}^q \sum_{\gamma=1}^q (-1)^{q-\beta} B_j^{\beta \gamma} \gamma_j\left(w_{0j}^\gamma \right) w_{0j}^1 \wedge \cdots \wedge \widehat{w}_{0j}^\beta \wedge \cdots \wedge w_{0j}^q \wedge d w_{0j}^\alpha \notag\\
&=\sum_{\beta \ne \gamma} (-1)^{q-\beta +\beta - \gamma+1} B_j^{\beta \gamma} \gamma_j\left( w_{0j}^\gamma \right) w_{0j}^1 \wedge \cdots \wedge \widehat{w}_{0j}^\gamma \wedge \cdots \wedge \overbrace{w_{0j}^\gamma}^{\beta-\textnormal{th}}\wedge \cdots \wedge w_{0j}^q \wedge dw_{0j}^\alpha + \sum_{\beta=1}^q (-1)^{q-\beta} B_j^{\beta \beta} \gamma_j\left( w_{0j}^\beta \right) w_{0j}^1 \wedge \cdots \wedge \widehat{w}_{0j}^\beta \wedge \cdots \wedge w_{0j}^q \wedge d w_{0j}^\alpha \notag\\
&=\sum_{\beta \ne \gamma} (-1)^{q-\gamma+1} \gamma_j\left(w_{0j}^\gamma \right) w_{0j}^1 \wedge \cdots \wedge \widehat{w}_{0j}^\gamma \wedge \cdots \wedge B_j^{\beta \gamma} w_{0j}^\gamma \wedge \cdots \wedge w_{0j}^q \wedge dw_{0j}^\alpha + \sum_{\beta=1}^q (-1)^{q-\beta} B_j^{\beta \beta} \gamma_j\left( w_{0j}^\beta \right) w_{0j}^1 \wedge \cdots \wedge \widehat{w}_{0j}^\beta \wedge \cdots \wedge w_{0j}^q \wedge d w_{0j}^\alpha \notag \\
&=\sum_{\beta \ne \gamma} (-1)^{q-\gamma+1} \gamma_j\left(w_{0j}^\gamma \right) w_{0j}^1 \wedge \cdots \wedge \widehat{w}_{0j}^\gamma \wedge \cdots \wedge \left( B_j^\beta - B_j^{\beta \beta} w_{0j}^\beta \right) \wedge \cdots \wedge w_{0j}^q \wedge dw_{0j}^\alpha + \sum_{\gamma=1}^q (-1)^{q-\gamma} B_j^{\gamma \gamma} \gamma_j\left( w_{0j}^\gamma \right) w_{0j}^1 \wedge \cdots \wedge \widehat{w}_{0j}^\gamma \wedge \cdots \wedge w_{0j}^q \wedge d w_{0j}^\alpha \notag\\
&=\sum_{\beta \ne \gamma} (-1)^{q-\gamma +1} \gamma_j\left( w_{0j}^\gamma \right) w_{0j}^1 \wedge \cdots \wedge \widehat{w}_{0j}^\gamma \wedge \cdots \wedge B_j^\beta \wedge \cdots \wedge w_{0j}^q \wedge dw_{0j}^\alpha + \sum_{\beta, \gamma=1}^q B_j^{\beta \beta} (-1)^{q-\gamma} \gamma_j\left( w_{0j}^\gamma \right) w_{0j}^1 \wedge \cdots \wedge \widehat{w}_{0j}^\gamma \wedge \cdots \wedge w_{0j}^q \wedge dw_{0j}^\alpha \notag \\
&= - \sum_{\beta \ne \gamma} (-1)^{q-\gamma} \gamma_j\left( w_{0j}^\gamma \right) w_{0j}^1 \wedge \cdots \wedge \widehat{w}_{0j}^\gamma \wedge \cdots \wedge B_j^\beta \wedge \cdots \wedge w_{0j}^q \wedge dw_{0j}^\alpha + \sum_{\beta=1}^q B_j^{\beta \beta}\left( w_{0j}^1 \wedge \cdots \wedge w_{0j}^q \wedge \left( \beta_j\left( w_{0j}^\alpha\right) - d\gamma_j\left( w_{0j}^\gamma \right) \right)    \right) \notag \\
&= - \sum_{\beta \ne \gamma} (-1)^{q-\gamma} \gamma_j\left( w_{0j}^\gamma \right) w_{0j}^1 \wedge \cdots \wedge \widehat{w}_{0j}^\gamma \wedge \cdots \wedge B_j^\beta \wedge \cdots  \wedge w_{0j}^q \wedge dw_{0j}^\alpha + \sum_{\beta=1}^q  w_{0j}^1 \wedge \cdots \wedge  B_j^\beta \wedge \cdots \wedge w_{0j}^q \wedge \left( \beta_j\left( w_{0j}^\alpha\right) - d\gamma_j\left( w_{0j}^\gamma \right) \right)     \notag
\end{align}}}

Then from $(\ref{b10})$ and $(\ref{b11})$ and $(\ref{b12})$, we have $(\ref{b13})$. Hence $\left\{ Q_j \right\}$ corresponds to $\sigma$. This completes the proof of Proposition \ref{b14}.

\end{proof}

\appendix

\section{Brackets and differentials on $\bigoplus^p A^{0,p}\left( M, \mathcal{E}_{\Theta_{\mathcal{F}_0}} \right), \bigoplus^p A^{0,p}\left( M, \mathcal{E}_{\mathcal{N}_{\mathcal{F}_0}^*}  \right)$ and $\bigoplus^p A^{0,p}\left(  M, \mathcal{E}_{\mathcal{F}_0}\right)$}

Let $\mathcal{F}$ be a vector bundle of rank $r$ on a compact complex manifold $M$. We recall that the associated Atiyah extension $\mathcal{E}_\mathcal{F}$ is defined in the following way: let $\mathcal{U}=\{U_i\}$ be an open covering of $M$ such that $\Gamma\left(U_i, \mathcal{F}\right)$ is generated by $w_i^1,...,w_i^r\in \Gamma\left(U_i, \mathcal{F}\right)$ and on $U_{ij}$, we have the relation $w_i^\alpha= \sum_{\beta=1}^q h_{ij}^{\alpha\beta} w_j^\beta,\alpha=1,...,r$. Then $\mathcal{E}_{\mathcal{F}}$ on $U_i$ is defined by $\mathcal{E}_\mathcal{F}|_{U_i}:= \Theta_M|_{U_i} \bigoplus \mathscr{H}om_{\mathcal{O}_M}\left( \mathcal{F}, \mathcal{F} \right)|_{U_i}$, and $\left(T_i, \phi_i \right)\in \Theta_M |_{U_i} \bigoplus \mathscr{H}om_{\mathcal{O}_M}\left( \mathcal{F}, \mathcal{F} \right)|_{U_i}$  and $\left(T_j, \phi_j\right)\in \Theta_M |_{U_j} \bigoplus \mathscr{H}om_{\mathcal{O}_M}\left(  \mathcal{F}, \mathcal{F} \right)|_{U_j}$ is identified on $U_i\cap U_j$ if and only if $T_i= T_j=T$ and $\phi_i\left(w_i^\alpha\right)= \phi_j\left(w_i^\alpha \right) + \sum_{\beta=1}^q \left[T, h_{ij}^{\alpha\beta} \right] w_j^\beta$. Then we have the exact sequence of sheaves
\begin{align*}
0\to \mathscr{H}om_{\mathcal{O}_M}\left(\mathcal{F}, \mathcal{F} \right) \to \mathcal{E}_{\mathcal{F}} \to \Theta_M \to 0
\end{align*}
Let us denote by $A^{0,p}(M, \mathcal{E}_\mathcal{F})$ the global sections of the sheaf $\mathcal{A}^{0,p}(\mathcal{E}_\mathcal{F})$ of $C^\infty$ $(0,p)$-forms with coefficients in $\mathcal{E}_\mathcal{F}$. Then we shall directly define the differential graded lie algebra structure on $\bigoplus^p A^{0,p}\left(M, \mathcal{E}_\mathcal{F} \right)$ by using the above definition of $\mathcal{E}_\mathcal{F}$  (which fits into the main body of the paper) without reference to a Hermitian connection and the associated curvature form (as in \cite{Hua95}).

We define the bracket 
\begin{align*}
\left[-,-\right]_{\mathcal{E}_{\mathcal{F}}}: A^{0,p}\left(M, \mathcal{E}_\mathcal{F}\right) \times A^{0,q}\left(M, \mathcal{E}_\mathcal{F}\right) \to A^{0,p+q} \left(M, \mathcal{E}_\mathcal{F}\right)
\end{align*}
which is locally defined by the following way: for $\left(\varphi_i, \eta_i \right), \left(\psi_i, \xi_i \right) \in \Gamma\left(U_i, \mathcal{A}^{0,p}\left( \Theta_M\right)\right)\bigoplus \Gamma\left(U_i, \mathcal{A}^{0, p} \left(\mathscr{H}om_{\mathcal{O}_M}\left(\mathcal{F}, \mathcal{F} \right)\right)\right)$
\begin{align*}
\left[ \left( \varphi_i, \eta_i \right), \left(\psi_i, \xi_i  \right) \right] =\left( \left[ \varphi_i, \psi_i \right], \left[ \eta_i, \xi_i \right] + \left[\varphi_i, \xi_i  \right]  - (-1)^{pq} \left[\psi_i, \eta_i \right]  \right) =\left( \left[ \varphi_i, \psi_i\right],   \eta_i\circ \xi_i - (-1)^{pq}\xi_i \circ \eta_i + \left[\varphi_i, \xi_i  \right]  - (-1)^{pq} \left[\psi_i, \eta_i \right]     \right)
\end{align*}
where $\eta_i \left(w_i^\alpha \right)=\sum_{ \beta=1}^q A_i^{\alpha\beta} w_i^\beta$ with $A_i^{\alpha\beta}\in \Gamma\left(U_i, \mathcal{A}^{0,p}\right)$, and $\xi_i\left(w_i^\alpha\right)=\sum_{\beta=1}^q B_i^{\alpha\beta}w_i^\beta $ with $B_i^{\alpha\beta}\in \Gamma\left(U_i, \mathcal{A}^{0,q}\right)$
\begin{align*}
\eta_i\circ \xi_i \left(w_i^\alpha\right)=\eta_i\left( \sum_{\beta=1}^r B_i^{\alpha\beta} w_i^\beta    \right) & = (-1)^{pq}\sum_{\beta=1}^r B_i^{\alpha\beta} \wedge \eta_i \left(w_i^\beta \right)  = (-1)^{pq} \sum_{\beta, \gamma=1}^r B_i^{\alpha\beta} \wedge A_i^{\beta \gamma} w_i^\gamma\\
\xi_i\circ \eta_i(w_i^\alpha)=\xi_i\left( \sum_{\beta=1}^r A_i^{\alpha\beta} w_i^\beta    \right) & = (-1)^{pq}\sum_{\beta=1}^r A_i^{\alpha\beta} \wedge \eta_i \left(w_i^\beta \right)   = (-1)^{pq} \sum_{\beta, \gamma=1}^r A_i^{\alpha\beta} \wedge B_i^{\beta \gamma} w_i^\gamma\\
\left[  \varphi_i, \xi_i   \right]\left(w_i^\alpha\right) -(-1)^{pq}\left[\psi_i, \xi_i  \right]\left(w_i^\alpha\right) & = \sum_{\beta=1}^r \left[ \varphi_i,  B_i^{\alpha\beta}  \right]w_i^\beta - (-1)^{pq}\sum_{\beta=1}^r \left[\psi_i, A_i^{\alpha\beta} \right] w_i^\beta
\end{align*}
and linearly extends to $\Gamma\left( U_i, \mathcal{F} \right)$. Then we have
\begin{align*}
\left[\eta_i, \xi_i \right]\left(w_i^\alpha\right)&:= \left(\eta_i \circ \xi_i - (-1)^{pq}\xi_i\circ \eta_i \right)\left( w_{i}^\alpha \right) =  (-1)^{pq} \sum_{\beta, \gamma=1}^r B_i^{\alpha\beta} \wedge A_i^{\beta \gamma} w_i^\gamma -  \sum_{\beta, \gamma=1}^r A_i^{\alpha\beta} \wedge B_i^{\beta \gamma} w_i^\gamma  \\
                   &= - (-1)^{pq}\left( (-1)^{pq} \sum_{\beta, \gamma=1}^r A_i^{\alpha\beta} \wedge B_i^{\beta \gamma} w_i^\gamma - \sum_{\beta, \gamma=1}^r  B_i^{\alpha\beta} \wedge A_i^{\beta \gamma} w_i^\gamma       \right) = -(-1)^{pq} \left[\xi_i, \eta_i   \right] \left(w_i^\alpha \right)
\end{align*}
and we have
\begin{align*}
&\left[  \varphi_i, \xi_i   \right](w_i^\alpha) -(-1)^{pq}\left[\psi_i, \xi_i  \right]\left(w_i^\alpha\right)= \sum_{\beta=1}^r \left[ \varphi_i,  B_i^{\alpha\beta}  \right]w_i^\beta - (-1)^{pq}\sum_{\beta=1}^r \left[\psi_i, A_i^{\alpha\beta} \right] w_i^\beta\\
&=-(-1)^{pq}\left( \sum_{\beta=1}^r  \left[\psi_i, A_i^{\alpha\beta} \right] w_i^\beta - (-1)^{pq} \sum_{\beta=1}^r \left[\varphi_i, B_i^{\alpha\beta} \right] w_i^\beta    \right)= -(-1)^{pq}\left( \left[  \xi, \varphi_i   \right](w_i^\alpha) -(-1)^{pq}\left[\xi_i, \varphi_i  \right] \right)
\end{align*}

We show that it is well-defined. In fact, if $\left(\varphi,\eta_i \right)= \left(\varphi, \eta_j \right)$ and $\left(\psi , \xi_i \right)= \left(\psi, \xi_j \right)$, in other words,
\begin{align*}
\eta_i \left(w_i^\alpha \right)= \eta_j\left(w_i^\alpha\right) + \sum_{\beta=1}^r \left[\varphi, h_{ij}^{\alpha\beta} \right]w_j^\beta & \iff \sum_{\gamma=1}^r A_i^{\alpha\gamma} w_i^\gamma = \sum_{\beta, \gamma=1}^r h_{ij}^{\alpha\beta} A_j^{\beta \gamma} w_j^\gamma + \sum_{\beta=1}^r \left[\varphi, h_{ij}^{\alpha\beta} \right] w_j^\beta\\
\xi_i \left(w_i^\alpha \right)= \xi_j \left(w_i^\alpha \right) + \sum_{\beta=1}^r \left[\psi, h_{ij}^{\alpha\beta} \right]w_j^\beta & \iff \sum_{\gamma=1}^r B_i^{\alpha \gamma} w_i^\gamma= \sum_{\beta, \gamma=1}^r h_{ij}^{\alpha \beta} B_j^{\beta \gamma} w_j^\gamma + \sum_{\beta=1}^r \left[ \psi, h_{ij}^{\alpha\beta} \right] w_j^\beta,
\end{align*}
then we have

{\small{\begin{align}\label{nn1}
&\left[ \eta_i, \xi_i \right]\left(w_i^\alpha\right) = (-1)^{pq} \sum_{\beta, \gamma=1}^r B_i^{\alpha\beta} \wedge A_i^{\beta \gamma} w_i^\gamma -  \sum_{\beta, \gamma=1}^r A_i^{\alpha\beta} \wedge B_i^{\beta \gamma} w_i^\gamma  = (-1)^{pq} \sum_{\beta, \gamma, \xi=1}^r B_i^{\alpha \beta} \wedge A_i^{\beta \gamma} h_{ij}^{\gamma \xi} w_j^\xi- \sum_{\beta, \gamma, \xi=1}^r A_i^{\alpha \beta} \wedge B_i^{\beta \gamma} h_{ij}^{\gamma \xi} w_j^\xi \\
 &= (-1)^{pq}\sum_{\beta,  \xi=1}^r B_i^{\alpha \beta} \wedge \left( \sum_{\gamma=1}^r h_{ij}^{\beta \gamma} A_j^{\gamma  \xi} + \left[\varphi, h_{ij}^{\beta \xi} \right] \right)w_j^\xi- \sum_{\beta,  \xi=1}^r A_i^{\alpha \beta} \wedge \left( \sum_{\gamma=1}^r h_{ij}^{\beta \gamma} B_j^{\gamma  \xi} + \left[\psi, h_{ij}^{\beta \xi} \right] \right)w_j^\xi \notag \\
 &=(-1)^{pq}\sum_{\beta, \xi,\gamma =1}^r B_i^{\alpha\beta} h_{ij}^{\beta \gamma} \wedge A_j^{\gamma \xi} w_j^\xi + (-1)^{pq} \sum_{\beta, \xi=1}^r B_i^{\alpha\beta} \wedge \left[ \varphi, h_{ij}^{\beta \xi}\right] w_j^\xi - \sum_{\beta, \xi,\gamma =1}^r A_i^{\alpha\beta} h_{ij}^{\beta \gamma} \wedge B_j^{\gamma \xi} w_j^\xi - \sum_{\beta, \xi=1}^r A_i^{\alpha\beta} \wedge \left[ \psi , h_{ij}^{\beta \xi}\right] w_j^\xi \notag \\
 &=(-1)^{pq}\sum_{\beta, \gamma , \xi=1}^r h_{ij}^{\alpha \beta} B_j^{\beta \gamma}\wedge A_j^{\gamma \xi} w_j^\xi +(-1)^{pq}\sum_{\gamma, \xi=1}^r \left[ \psi, h_{ij}^{\alpha \gamma}  \right]\wedge A_j^{\gamma \xi} w_j^\xi+ (-1)^{pq} \sum_{\beta, \xi=1}^r B_i^{\alpha\beta} \wedge \left[ \varphi, h_{ij}^{\beta \xi}\right] w_j^\xi \notag \\
 &-\sum_{\beta, \gamma , \xi=1}^r h_{ij}^{\alpha \beta} A_j^{\beta \gamma}\wedge B_j^{\gamma \xi} w_j^\xi - \sum_{\gamma, \xi=1}^r \left[ \varphi, h_{ij}^{\alpha \gamma}  \right]\wedge B_j^{\gamma \xi} w_j^\xi- \sum_{\beta, \xi=1}^r A_i^{\alpha\beta} \wedge \left[ \psi ,  h_{ij}^{\beta \xi}\right] w_j^\xi \notag \\
 &=\left[\eta_j, \xi_j \right]\left(w_i^\alpha\right) +(-1)^{pq}\sum_{\gamma, \xi=1}^r \left[ \psi, h_{ij}^{\alpha \gamma}  \right]\wedge A_j^{\gamma \xi} w_j^\xi+ (-1)^{pq} \sum_{\beta, \xi=1}^r B_i^{\alpha\beta} \wedge \left[ \varphi, h_{ij}^{\beta \xi}\right] w_j^\xi  - \sum_{\gamma, \xi=1}^r \left[ \varphi, h_{ij}^{\alpha \gamma}  \right]\wedge B_j^{\gamma \xi} w_j^\xi- \sum_{\beta, \xi=1}^r A_i^{\alpha\beta} \wedge \left[ \psi ,  h_{ij}^{\beta \xi}\right] w_j^\xi \notag
\end{align}}}
On the other hand, we have
{\small{\begin{align} \label{nn2}
&\left[  \varphi , \xi_j  \right]\left(w_i^\alpha\right) -(-1)^{pq}\left[\psi, \eta_j  \right]\left(w_i^\alpha \right)= \sum_{\beta=1}^r h_{ij}^{\alpha\beta} \left[\varphi, \xi_j \right] \left(w_j^\beta\right) - (-1)^{pq} \sum_{\beta=1}^r h_{ij}^{\alpha\beta} \left[\psi, \eta_j \right]\left(w_j^\beta\right) \\
& = \sum_{\beta,\gamma=1}^r h_{ij}^{\alpha\beta}\left[ \varphi,  B_j^{\beta \gamma}  \right]w_j^\gamma - (-1)^{pq}\sum_{\beta,\gamma=1}^r h_{ij}^{\alpha \beta}\left[\psi, A_i^{\beta\gamma} \right] w_j^\gamma \notag  \\
&=\sum_{ \gamma=1}^r \left[\varphi, \sum_{\beta=1}^r B_i^{\alpha \beta} h_{ij}^{\beta \gamma} -\left[\psi, h_{ij}^{\alpha\gamma} \right]     \right] w_j^\gamma - \sum_{\beta, \gamma=1}^r [\varphi, h_{ij}^{\alpha \beta}] \wedge  B_j^{\beta \gamma} w_j^\gamma - (-1)^{pq} \sum_{ \gamma=1}^r \left[ \psi,  \sum_{\beta=1}^r A_i^{\alpha \beta} h_{ij}^{\beta \gamma}- \left[ \varphi, h_{ij}^{\alpha \gamma} \right]  \right] w_j^\gamma \notag \\
&+ (-1)^{pq}\sum_{\beta, \gamma=1}^r \left[ \psi, h_{ij}^{\alpha\beta} \right] \wedge A_i^{\beta \gamma} w_j^\gamma  \notag \\
&=\sum_{\gamma, \beta=1}^r(-1)^{pq} B_i^{\alpha \beta} \wedge \left[\varphi, h_{ij}^{\beta \gamma} \right] w_j^\gamma  + \left[\varphi,  B_i^{\alpha \beta}   \right]w_i^\beta - \sum_{\gamma=1}^r \left[   \varphi, \left[\psi, h_{ij}^{\alpha \gamma} \right] \right] w_j^\gamma - \sum_{\beta, \gamma=1}^r \left[ \varphi, h_{ij}^{\alpha\beta} \right] \wedge B_j^{\beta \gamma} w_j^\gamma - \sum_{\beta, \gamma=1}^r A_i^{\alpha \gamma} \wedge \left[ \psi, h_{ij}^{\beta \gamma}\right] w_j^\gamma \notag \\
& - (-1)^{pq} \sum_{\gamma, \beta=1}^r \left[ \psi, A_i^{\alpha \beta} \right] w_i^\beta + (-1)^{pq}\sum_{\gamma=1}^r \left[ \psi, \left[ \varphi, h_{ij}^{\alpha \gamma} \right] \right] w_j^\gamma + (-1)^{pq} \sum_{\beta, \gamma=1}^r \left[ \psi, h_{ij}^{\alpha\beta} \right] \wedge A_i^{\beta \gamma} w_j^\gamma \notag \\
&=\left[ \varphi, \xi_i\right]\left( w_i^\alpha\right)-(-1)^{pq} \left[ \psi, \eta_i \right]\left( w_i^\alpha \right) - \sum_{\gamma=1}^r \left[ \left[ \varphi, \psi \right], h_{ij}^{\alpha \gamma} \right] w_j^\gamma + \sum_{\gamma, \beta=1}^r (-1)^{pq} B_i^{\alpha \beta} \wedge \left[ \varphi, h_{ij}^{\beta \gamma} \right] w_{j}^\gamma - \sum_{\beta ,\gamma=1}^r \left[ \varphi, h_{ij}^{\alpha\beta}\right] \wedge B_j^{\beta \gamma} w_j^\gamma  \notag \\
&- \sum_{\beta, \gamma=1}^r A_i^{\alpha \gamma} \wedge \left[ \psi, h_{ij}^{\beta \gamma} \right] w_j^\gamma + (-1)^{pq} \sum_{\beta, \gamma=1}^r \left[ \psi, h_{ij}^{\alpha\beta} \right] \wedge A_i^{\beta \gamma} w_j^\gamma \notag
\end{align}}}

Then from $\textnormal{(\ref{nn1})}$ and $\textnormal{(\ref{nn2})}$ we see that
\begin{align*}
\left[ \eta_i, \xi_i \right]\left( w_{i}^\alpha \right) +\left[ \varphi, \xi_i \right]\left( w_{0}^\alpha \right)  - (-1)^{pq} \left[ \psi, \eta_i \right] \left( w_{i}^\alpha \right)= \left[ \eta_j, \xi_j \right]\left(w_{0i}^\alpha \right) + \left[ \varphi, \xi_j\right] \left( w_{i}^\alpha \right) - (-1)^{pq} \left[ \psi, \eta_j \right]\left( w_{i}^\alpha \right) + \sum_{\gamma=1}^r \left[ \left[ \varphi, \psi \right], h_{ij}^{\alpha\gamma} \right] w_j^\gamma
\end{align*}
Hence the bracket $\left[ -, - \right]_{\mathcal{E}_\mathcal{F}}$ is well-defined. On the other hand, we define the differential $\bar{\partial}$ locally in the following way: for $\left( \varphi_i, \eta_i \right)\in \Gamma\left( U_i, \mathcal{A}^{0,p}\left( \mathcal{E}_\mathcal{F} \right) \right)$ with $\eta_i\left( w_i^\alpha\right)= \sum_{\beta=1}^r A_i^{\alpha\beta} w_i^\beta $ for some $A_i^{\alpha\beta}\in \Gamma\left(U_i,  \mathcal{A}^{0,p}\right)$,
\begin{align*}
\bar{\partial}(\varphi_i, \eta_i)=\left(\bar{\partial}\varphi_i, \bar{\partial}\eta_i \right)
\end{align*}
where
\begin{align*}
\left(\bar{\partial} \eta_i \right)(w_i^\alpha)= \sum_{\beta=1}^r \bar{\partial} A_i^{\alpha \beta} w_i^\beta
\end{align*}
It is well-defined since for $\left(\varphi, \eta_i\right)=\left( \varphi, \eta_j\right)$, i.e. $\sum_{\gamma=1}^r A_i^{\alpha \gamma} w_i^\gamma= \sum_{\beta, \gamma=1}^r h_{ij}^{\alpha\beta} A_j^{\beta \gamma} w_j^\gamma + \sum_{\beta=1}^r \left[ \varphi, h_{ij}^{\alpha\beta} \right] w_j^\beta$, we have
\begin{align*}
 \sum_{\gamma=1}^r \bar{\partial}A_i^{\alpha\gamma} w_i^\gamma = \sum_{\beta, \gamma=1}^r h_{ij}^{\alpha\beta} \bar{\partial}A_j^{\beta \gamma} w_j^\gamma + \sum_{\beta=1}^r \left[\bar{\partial}\varphi, h_{ij}^{\alpha\beta} \right] w_j^\beta
\end{align*}

\subsection{Bracket and differential on $\bigoplus^p A^{0,p}\left( M, \mathcal{E}_{\Theta_{\mathcal{F}_0}}\right)$} \label{apd1} \

Let $\left( M, \Theta_{\mathcal{F}_0} \right)$ be a compact foliated complex manifold with $\Theta_{\mathcal{F}_0}$ locally free. Let $\mathcal{U}=\left\{ U_i \right\}$ be a Stein open covering of $M$ such that $\Gamma\left( U_j, \Theta_{\mathcal{F}_0}\right)$ is generated by $T_{0j}^1,..., T_{0j}^p \in \Gamma\left( U_j , \Theta_M \right) $ with the relation $\sum_{\beta=1}^p r_{0ij}^{\alpha\beta} T_{0j}^\beta$. We recall from Part II that $\mathcal{E}_{\Theta_{\mathcal{F}_0}}$ is defined by $\mathcal{E}_{\Theta_{\mathcal{F}_0}}|_{U_i}:=\Theta_M|_{U_i} \bigoplus \mathscr{H}om_{\mathcal{O}_M}\left(\Theta_{\mathcal{F}_0}, \Theta_{\mathcal{F}_0} \right)|_{U_i}$, and $\left( T_i, \phi_i \right) \in \Theta_M|_{U_i} \bigoplus \mathscr{H}om_{\mathcal{O}_M}\left( \Theta_{\mathcal{F}_0}, \Theta_{\mathcal{F}_0} \right)|_{U_i}$ and $\left(T_j, \phi_j \right)\in \Theta_M|_{U_j}\bigoplus \mathscr{H}om_{\mathcal{O}_M}\left( \Theta_{\mathcal{F}_0}, \Theta_{\mathcal{F}_0} \right)|_{U_j}$ are identified on $U_i\cap U_j$ if and only if $T_i=T_j=T$ and $\phi_i\left( T_{0i}^\alpha \right)= \phi_j\left(T_{0i}^\alpha \right) + \sum_{\beta=1}^p \left[ T, r_{0ij}^{\alpha\beta} \right] T_{0j}^\beta$.

We define the bracket on $\bigoplus^p A^{0,p}\left( M, \mathcal{E}_{\Theta_{\mathcal{F}_0}} \right)$ by
\begin{align*}
\left[ - , - \right]_{\Theta_{\mathcal{F}_0}}: A^{0,p}\left( M, \mathcal{E}_{\Theta_{\mathcal{F}_0}} \right) \times A^{0,q}\left( M, \mathcal{E}_{\Theta_{\mathcal{F}_0}}\right) \to A^{0,p+q}\left( M, \mathcal{E}_{\Theta_{\mathcal{F}_0}} \right)
\end{align*}
which is locally defined in the following way: for $\left(\varphi_i, \eta_i \right), \left(\psi_i, \xi_i \right) \in \Gamma\left(U_i, \mathcal{A}^{0,p}\left( \mathcal{E}_{\Theta_{\mathcal{F}_0}} \right) \right)$
\begin{align*}
\left[ \left( \varphi_i, \eta_i \right), \left(\psi_i, \xi_i  \right) \right] =\left( \left[ \varphi_i, \psi_i \right], \left[ \eta_i, \xi_i \right] + \left[\varphi_i, \xi_i  \right]  - (-1)^{pq} \left[\psi_i, \eta_i \right]  \right) 
\end{align*}
and we define the differential $\bar{\partial}: A^{0,p}\left(M, \mathcal{E}_{\Theta_{\mathcal{F}_0}} \right) \to A^{0,p+1}\left( M, \mathcal{E}_{ \Theta_{\mathcal{F}_0}} \right)$ locally by
\begin{align*}
\bar{\partial}(\varphi_i, \eta_i)=\left(\bar{\partial}\varphi_i, \bar{\partial}\eta_i \right)
\end{align*}

\subsection{Bracket and differential on $\bigoplus^p A^{0,p}\left( M, \mathcal{E}_{\mathcal{N}_{\mathcal{F}_0}^*} \right)$} \label{apd2} \

Let $\left( M, \mathcal{N}_{\mathcal{F}_0}^* \right)$ be a compact foliated complex manifold with $\mathcal{N}_{\mathcal{F}_0}^*$ locally free. Let $\mathcal{U}=\left\{ U_i \right\}$ be a Stein open covering of $M$ such that $\Gamma\left( U_j, \mathcal{N}_{\mathcal{F}_0}^*\right)$ is generated by $w_{0j}^1,..., w_{0j}^q \in \Gamma\left( U_j , \Omega_M^1 \right) $ with the relation $\sum_{\beta=1}^q h_{0ij}^{\alpha\beta} w_{0j}^\beta$. We recall from Part II that $\mathcal{E}_{\mathcal{N}_{\mathcal{F}_0}^*}$ is defined by $\mathcal{E}_{\mathcal{N}_{\mathcal{F}_0}^*}|_{U_i}:=\Theta_M|_{U_i} \bigoplus \mathscr{H}om_{\mathcal{O}_M}\left(\mathcal{N}_{\mathcal{F}_0}^*, \mathcal{N}_{\mathcal{F}_0}^* \right)|_{U_i}$, and $\left( T_i, \phi_i' \right) \in \Theta_M|_{U_i} \bigoplus \mathscr{H}om_{\mathcal{O}_M}\left( \mathcal{N}_{\mathcal{F}_0}^*, \mathcal{N}_{\mathcal{F}_0}^* \right)|_{U_i}$ and $\left(T_j, \phi_j' \right)\in \Theta_M|_{U_j}\bigoplus \mathscr{H}om_{\mathcal{O}_M}\left( \mathcal{N}_{\mathcal{F}_0}^*, \mathcal{N}_{\mathcal{F}_0}^* \right)|_{U_j}$ are identified on $U_i\cap U_j$ if and only if $T_i=T_j=T$ and $\phi_i'\left( w_{0i}^\alpha \right)= \phi_j'\left(w_{0i}^\alpha \right) + \sum_{\beta=1}^q \left[ T, h_{0ij}^{\alpha\beta} \right] w_{0j}^\beta$.

We define the bracket on $\bigoplus^p A^{0,p}\left( M, \mathcal{E}_{\mathcal{N}_{\mathcal{F}_0}^*} \right)$ by
\begin{align*}
\left[ - , - \right]_{\mathcal{N}_{\mathcal{F}_0}^*}: A^{0,p}\left( M, \mathcal{E}_{\mathcal{N}_{\mathcal{F}_0}^*} \right) \times A^{0,q}\left( M, \mathcal{E}_{\mathcal{N}_{\mathcal{F}_0}^*}\right) \to A^{0,p+q}\left( M, \mathcal{E}_{\mathcal{N}_{\mathcal{F}_0}^*} \right)
\end{align*}
which is locally defined in the following way: for $\left(\varphi_i, \eta_i' \right), \left(\psi_i, \xi_i' \right) \in \Gamma\left(U_i, \mathcal{A}^{0,p}\left( \mathcal{E}_{\Theta_{\mathcal{F}_0}} \right) \right)$
\begin{align*}
\left[ \left( \varphi_i, \eta_i' \right), \left(\psi_i, \xi_i' \right) \right] =\left( \left[ \varphi_i, \psi_i \right], \left[ \eta_i', \xi_i' \right] + \left[\varphi_i, \xi_i'  \right]  - (-1)^{pq} \left[\psi_i, \eta_i' \right]  \right) 
\end{align*}
and we define the differential $\bar{\partial}: A^{0,p}\left(M, \mathcal{E}_{\Theta_{\mathcal{F}_0}} \right) \to A^{0,p+1}\left( M, \mathcal{E}_{ \Theta_{\mathcal{F}_0}} \right)$ locally by
\begin{align*}
\bar{\partial}(\varphi_i, \eta_i')=\left(\bar{\partial}\varphi_i, \bar{\partial}\eta_i' \right)
\end{align*}

\subsection{Bracket and differential on $\bigoplus^p A^{0,p}\left( M, \mathcal{E}_{\mathcal{F}_0}\right)$}\label{apd3} \

Let $\left( M, \Theta_{\mathcal{F}_0} , \mathcal{N}_{\mathcal{F}_0}^* \right)$ be a compact foliated complex manifold with $\Theta_{\mathcal{F}_0}$ and $\mathcal{N}_{\mathcal{F}_0}^*$ locally free. Let $\mathcal{U}=\left\{ U_i \right\}$ be a Stein open covering of $M$ such that $\Gamma\left( U_j, \Theta_{\mathcal{F}_0}\right)$ is generated by $T_{0j}^1,..., T_{0j}^p \in \Gamma\left( U_j , \Theta_M \right) $ with the relation $\sum_{\beta=1}^p r_{0ij}^{\alpha\beta} T_{0j}^\beta$, and on the other hand, $\Gamma\left( U_j, \mathcal{N}_{\mathcal{F}_0}^*\right)$ is generated by $w_{0j}^1,..., w_{0j}^q \in \Gamma\left( U_j , \Omega_M^1 \right) $ with the relation $\sum_{\beta=1}^q h_{0ij}^{\alpha\beta} w_{0j}^\beta$. We recall from Part II that $\mathcal{E}_{\mathcal{F}_0}$ is defined locally by $\mathcal{E}_{\mathcal{F}_0}|_{U_i} :=\Theta_M|_{U_i} \bigoplus \mathscr{H}om_{\mathcal{O}_M}\left(\Theta_{\mathcal{F}_0}, \Theta_{\mathcal{F}_0} \right)|_{U_i}$ $ \bigoplus \mathscr{H}om_{\mathcal{O}_M}\left( \mathcal{N}_{\mathcal{F}_0}^*, \mathcal{N}_{\mathcal{F}_0}^*  \right)|_{U_i}$, and $\left( T_i, \phi_i , \phi_i' \right) \in \Theta_M|_{U_i} \bigoplus \mathscr{H}om_{\mathcal{O}_M}\left( \Theta_{\mathcal{F}_0}, \Theta_{\mathcal{F}_0} \right)|_{U_i} \bigoplus \mathscr{H}om_{\mathcal{O}_M}\left( \mathcal{N}_{\mathcal{F}_0}^* , \mathcal{N}_{\mathcal{F}_0}^* \right)|_{U_i}$ and $\left(T_j, \phi_j , \phi_j' \right)\in \Theta_M|_{U_j}\bigoplus \mathscr{H}om_{\mathcal{O}_M}\left( \Theta_{\mathcal{F}_0}, \Theta_{\mathcal{F}_0} \right)|_{U_j} \bigoplus \mathscr{H}om_{\mathcal{O}_M}\left( \mathcal{N}_{\mathcal{F}_0}^*, \mathcal{N}_{\mathcal{F}_0}^* \right)|_{U_j}$ are identified on $U_i\cap U_j$ if and only if $T_i=T_j=T$ and $\phi_i\left( T_{0i}^\alpha \right)= \phi_j\left(T_{0i}^\alpha \right) + \sum_{\beta=1}^p \left[ T, r_{0ij}^{\alpha\beta} \right] T_{0j}^\beta$ and $\phi_i'\left( w_{0i}^\alpha \right)= \phi_j'\left( w_{0i}^\alpha \right) + \sum_{\beta=1}^q \left[ T, h_{0ij}^{\alpha\beta} \right] w_{0j}^\beta$.

We define the bracket on $\bigoplus^p A^{0,p}\left( M, \mathcal{E}_{\mathcal{F}_0} \right)$ by
\begin{align*}
\left[ - , - \right]_{\mathcal{F}_0}: A^{0,p}\left( M, \mathcal{E}_{\mathcal{F}_0} \right) \times A^{0,q}\left( M, \mathcal{E}_{\mathcal{F}_0}\right) \to A^{0,p+q}\left( M, \mathcal{E}_{\mathcal{F}_0} \right)
\end{align*}
which is locally defined in the following way: for $\left(\varphi_i, \eta_i \right), \left(\psi_i, \xi_i \right) \in \Gamma\left(U_i, \mathcal{A}^{0,p}\left( \mathcal{E}_{\mathcal{F}_0} \right) \right)$
\begin{align*}
\left[ \left( \varphi_i, \eta_i , \eta_i' \right), \left(\psi_i, \xi_i , \xi_i' \right) \right] =\left( \left[ \varphi_i, \psi_i \right], \left[ \eta_i, \xi_i \right] + \left[\varphi_i, \xi_i  \right]  - (-1)^{pq} \left[\psi_i, \eta_i \right] , \left[ \eta_i', \xi_i' \right] + \left[\varphi_i, \xi_i'  \right]  - (-1)^{pq} \left[\psi_i, \eta_i' \right]  \right) 
\end{align*}
and we define the differential $\bar{\partial}: A^{0,p}\left(M, \mathcal{E}_{\mathcal{F}_0} \right) \to A^{0,p+1}\left( M, \mathcal{E}_{ \mathcal{F}_0 } \right)$ locally by
\begin{align*}
\bar{\partial}(\varphi_i, \eta_i  , \eta_i ')=\left(\bar{\partial}\varphi_i, \bar{\partial}\eta_i , \bar{\partial}  \eta_i'  \right)
\end{align*}

\bibliographystyle{amsalpha}
\bibliography{References-Rev9}

\end{document}